%% file: main.tex
\documentclass[11pt,a4paper]{amsart}
\usepackage[T1]{fontenc}
\usepackage[margin=1.48in]{geometry}
\usepackage{graphicx}

\graphicspath{ {../../images/} }

\usepackage{amsthm, amsmath, amsfonts, amssymb, mathtools}
\usepackage{bbold}
\usepackage{booktabs}
\usepackage[dvipsnames]{xcolor}
\usepackage{tikz}
\usepackage{subcaption}
\usepackage{algorithm}
\usepackage{algpseudocode}
\usepackage{longtable}
\usepackage{makecell}
\usepackage{bbm}
\usepackage{pgfplots}
\pgfplotsset{compat=1.18}

\usepackage{xspace}
\usepackage{csquotes}
\usepackage[colorlinks=true, pdfstartview=FitV, linkcolor=Blue, citecolor=OliveGreen, urlcolor=BrickRed, breaklinks=false]{hyperref}
\hypersetup{hypertexnames=false}

\DeclarePairedDelimiter\abs{\lvert}{\rvert}%

\usetikzlibrary{calc}
\usetikzlibrary{arrows.meta}
\tikzset{
  axis/.style    = {thick, -{Stealth[length=2mm, width=1.3mm]}},
  grid/.style    = {black!30, thin},
  polygon/.style = {},
  newton/.style  = {black!20},
  curve/.style   = {line width=2pt},
}

\newcommand\1{\mathbb{1}}
\newcommand{\CC}{\mathbb{C}}
\newcommand{\FF}{\mathbb{F}}
\newcommand{\RR}{\mathbb{R}}
\newcommand{\PP}{\mathbb{P}}
\newcommand{\ZZ}{\mathbb{Z}}
\newcommand{\realized}[1]{{\boldmath $#1$}}
\newcommand{\openmark}[1]{$#1$\textsuperscript{*}}
\DeclareMathOperator\conv{conv} %
\DeclareMathOperator\interior{int}  %
\DeclareMathOperator{\parity}{parity}

\newcommand{\CP}{\CC\PP}
\newcommand{\RP}{\RR\PP}
\newcommand{\Sph}{\mathbb{S}}

\newcommand{\cC}{\mathcal{C}}

\newcommand{\cH}{\mathcal{H}}

\newcommand{\cS}{\mathcal{S}}
\newcommand{\frS}{\mathfrak{S}}
\newcommand{\cT}{\mathcal{T}}
\newcommand{\cZ}{\mathcal{Z}}
\newcommand{\GF}[1]{\mathbb{F}_{#1}}
\newcommand{\polymake}{\texttt{po\-ly\-ma\-ke}\xspace}

\newcommand{\Mathematica}{\texttt{Mathematica}\xspace}
\newcommand\smallSetOf[2]{\{#1 \mid #2\}}

\DeclareMathOperator{\Aut}{Aut} %
\newcommand\scheme[1]{\langle #1 \rangle}
\newcommand\Sym[1]{\mathrm{Sym}_{#1}} %
\newcommand{\highlight}[2]{\smallskip\begin{center}\begin{minipage}{.95\textwidth}\itshape #1\normalfont\quad\small (\mbox{#2})\end{minipage}\end{center}\medskip}

\newcommand{\degsixA}{bat}  %
\newcommand{\degsixAtypes}{53}
\newcommand{\degsixB}{moth}  %
\newcommand{\degsixBtypes}{44}

\newcommand{\degsevenA}{cen}  %
\newcommand{\degsevenAtypes}{115}
\newcommand{\degsevenB}{spl}  %
\newcommand{\degsevenBtypes}{107}
\newcommand{\degsevenC}{fra}  %
\newcommand{\degsevenCtypes}{103}
\newcommand{\degsevenD}{hon}  %
\newcommand{\degsevenDtypes}{47}

\newtheorem{theorem}{Theorem}
\newtheorem{corollary}[theorem]{Corollary}
\newtheorem{lemma}[theorem]{Lemma}
\newtheorem{proposition}[theorem]{Proposition}
\theoremstyle{definition}

\newtheorem{example}[theorem]{Example}
\newtheorem{remark}[theorem]{Remark}
\newtheorem{question}[theorem]{Question}

\usepackage[colorinlistoftodos, bordercolor=orange, backgroundcolor=orange!20, linecolor=orange, textsize=scriptsize]{todonotes}
\makeatletter
\renewcommand{\todo}[2][]{\tikzexternaldisable\@todo[#1]{#2}\tikzexternalenable}
\makeatother

\definecolor{myred}{RGB}{219, 86, 42}
\definecolor{myblue}{RGB}{8, 46, 79}

\definecolor{mycolor5}{RGB}{123, 64, 168}   %
\definecolor{mycolor4}{RGB}{52, 88, 164}      %
\definecolor{mycolor6}{RGB}{78, 140, 74}     %
\definecolor{mycolor3}{RGB}{224, 170, 38}   %
\definecolor{mycolor2}{RGB}{214, 110, 36}   %
\definecolor{mycolor1}{RGB}{192, 64, 64}       %

\usepackage[backend=biber,style=numeric,doi=false,isbn=false,url=false,giveninits=true,maxbibnames=50]{biblatex}

\title{Limits of combinatorial patchworking}
\author[Geiselmann et al.]{Zoe Geiselmann}
\author[]{Michael Joswig}
\author[]{Lars Kastner}
\author[]{Konrad Mundinger}
\author[]{Sebastian Pokutta}
\author[]{Christoph Spiegel}
\author[]{Marcel Wack}
\author[]{Max Zimmer}

\address[Z.~Geiselmann, M.~Joswig, L.~Kastner, M.~Wack]{Technische Universität Berlin, Chair of Discrete Mathe\-ma\-tics/Geo\-me\-try}
\email{\{geiselmann,joswig,kastner,wack\}@math.tu-berlin.de}

\address[K.~Mundinger, S.~Pokutta, C.~Spiegel, M.~Zimmer]{Zuse Institute Berlin, AI in Society, Science, and Technology, Takustraße 7, 14195 Berlin, Germany}
\email{\{mundinger,pokutta,spiegel,zimmer\}@zib.de}

\begin{document}

\begin{abstract}
  It is shown that there are real plane algebraic curves of degree eight that cannot be realized as T-curves, i.e., via combinatorial patchworking.
  In fact, this holds for several real schemes (i.e., ambient isotopy types) with the maximal number of real components, called $M$-curves.
  On the other hand, each nonempty real scheme of lower degree, maximal or not, arises as a T-curve.
  By constructing one patchwork of the dilated triangle $d\cdot\Delta_2$ for each nonempty real scheme of degree $d\leq 7$, we provide an explicit method for constructing polynomials realizing these real schemes.
  This resolves a question of Itenberg and Viro (1996).
\end{abstract}

\maketitle

\section{Introduction}

The modern topological classification of real plane projective algebraic curves begins with Harnack~\cite{1876:Harnack:VieltheiligkeitEbenen} in 1876, who showed that the number of real components of a real curve of degree $d$ is bounded by
\begin{equation}\label{eq:M}
  M = \tfrac{1}{2}(d-1)(d-2)+1 \enspace .
\end{equation}
Via an explicit construction he also proved that this bound is tight for all $d\geq 1$.
Since then these curves have been called \emph{maximal} or \emph{$M$-curves}.
Hilbert~\cite{1891:Hilbert:ReellenZuge} continued this work in 1891 by constructing more $M$-curves for $d\in\{6,7,8\}$.

The precise meaning of such a topological classification admits several interpretations; we focus on the \emph{ambient isotopy class} of a curve in the real projective plane, called its \emph{real scheme}.
The full classification based on this notion is not difficult for $d\leq 5$ and was already known to Harnack and Hilbert.
In 1900, Hilbert posed his celebrated 16th problem~\cite{1900:Hilbert:Mathematical}: the topological classification of real plane algebraic curves of arbitrary degree.
The classification for $d=6$ was completed by Gudkov in 1969 \cite{1969:Gudkov:ArrangementOvals, 1974:Gudkov:ProjectiveAlgebraic}, and Viro~\cite{1984:Viro:GluingAlgebraic} resolved degree seven, building on work of Rokhlin~\cite{1978:Rokhlin:TopologicalAlgebraic}, Nikulin~\cite{1980:Nikulin:SymmetricBilinear}, and others.
The case $d=8$ is almost settled, but not entirely; see Orevkov \cite{2002:Orevkov:CurvesIsotopy} for details.
For $d \geq 9$ our knowledge is rather fragmented; e.g., see \cite{2010:Orevkov:FormulasNests}.

Classifying real schemes for fixed degree $d$ requires two things: constructing sufficiently many curves, and excluding other possibilities.
For the latter, B\'ezout's theorem suffices for small degrees; higher degrees require further tools such as the Gudkov--Rokhlin congruence~\cite{1971:Arnold:ArrangementOvals, 1972:Rokhlin:GudkovsHypothesis}.
For the former, many ad hoc constructions exist; Shustin's $M$-curve of degree eight~\cite{1987:Shustin:CurveDegree}, for instance, relies on careful analysis of singularities of type $Z_{11}$. %
Much of the previous work was systematized by Viro~\cite{1986:Viro:TopologyAlgebraic}, who introduced patchworking as a general combinatorial method to produce many curves.
Nowadays, patchworking may be seen as a version of tropical geometry over the reals~\cite{2008:Viro:SixteenthHilbert}; for introductions to tropical geometry, see \cite{2009:ItenbergEtAl:AlgebraicGeometry,2015:MaclaganSturmfels:TropicalGeometry,2021:Joswig:EssentialsCombinatorics}.

The simplest variant of patchworking is \emph{combinatorial patchworking}: a unimodular regular triangulation of the lattice point set $A = \{(i,j) \in \ZZ^2 : i,j \geq 0,\, i+j \leq d\}$ together with a sign distribution $\sigma : A \to \GF{2}$ determines a real plane algebraic curve of degree $d$.
The resulting curves are known as \emph{T-curves} (for \enquote{triangle}).
Itenberg \cite[\S7]{1995:Itenberg:RagsdaleCurves} showed that, for sufficiently large degree, there are nonempty real schemes that cannot be realized by a T-curve of the same degree.
While realizability by T-curves seems to be known for $d \leq 7$ (see Appendix~\ref{app:known}), we are not aware of a prior explicit complete list of such realizations; for individual T-curves of degrees seven and eight see De Loera and Wicklin \cite{1998:DeLoeraWicklin:ConvexityPatchworking}.
In \cite[p.22]{1996:ItenbergViro:PatchworkingRagsdale} Itenberg and Viro explicitly asked whether every (nonempty) real scheme with $d=7$ or $d=8$ is realizable as a T-curve.
The purpose of this article is to answer this question for both degrees, with opposite outcomes.
\highlight{In degree $d=8$ there are $M$-curves not realizable as T-curves.}{Corollary~\ref{cor:3_19_impossible}}
We were informed by Alo\"{\i}s Demory that he had obtained the above result independently; it was announced in a talk held on April 3, 2026 in the S\'eminaire G\'eom\'etrie \'enum\'erative at the Institut de Math\'ematiques, Jussieu, Paris.

The situation concerning realizability is very different for lower degrees.
In degree $d\leq 7$, every nonempty real scheme of degree $d$ is realizable by a T-curve.
In particular, based on Viro's classification \cite{1984:Viro:GluingAlgebraic}, we show:
\highlight{Four triangulations cover all 121 real schemes in degree $d=7$.}{Theorem~\ref{thm:covering-7}}
Note that one triangulation can realize many nonempty real schemes through different sign distributions.
We also show that two triangulations suffice for the 55 nonempty real schemes of degree six (Theorem~\ref{thm:covering-6}), classified by Gudkov \cite{1969:Gudkov:ArrangementOvals, 1974:Gudkov:ProjectiveAlgebraic}.
For each degree $d \leq 5$, a single triangulation is enough (Proposition~\ref{thm:covering-5}).
In this way, the distinction between $d\leq 7$ and $d \geq 8$ delineates a limit of combinatorial patchworking.

Our 55 T-curves of degree six may be compared to the 64 curves in \cite{2019:KaihnsaEtAl:SixtyCurves}.
Kaihnsa et al.\ give explicit polynomials with integer coefficients, avoiding tropical methods.
The sampling results in \cite[Table 2]{2019:KaihnsaEtAl:SixtyCurves} are found by computing the real scheme of a given polynomial.
But this step is costly as it relies on cylindrical algebraic decomposition \cite{1975:Collins:QuantifierAlgebraic}.
One advantage of our approach is that the topological type of a patchworked curve is easy to compute \cite[Algorithm 1.4.E]{2006:Viro:PatchworkingAlgebraic}; see Section~\ref{sec:computation}.
In fact, for modest degree it is even feasible to accomplish by hand.
Thus, proving our main results does not require a computer, although the relevant patchworks were found by the large computer search described in Section~\ref{sec:realizations}.
The accompanying software and data are described in Appendix~\ref{app:software}.
For applications of plane curves, both classical and tropical representations can be useful.

If few triangulations suffice, which real schemes does each one produce?
Gayet and Welschinger~\cite[\S4.2]{2011:GayetWelschinger:ExponentialRarefaction} showed that curves close to the Harnack bound are measure-theoretically rare, and asked \cite[p.94]{2011:GayetWelschinger:ExponentialRarefaction} about the expected number of connected components.
We investigate the combinatorial analog: for a fixed triangulation, what is the distribution of real schemes as the signs vary?
\highlight{For each of the triangulations above, we exhaustively determine the frequency with which each nonempty real scheme arises over all sign distributions.}{Appendix~\ref{app:statistics}}

\subsection*{Outline}
Section~\ref{sec:viro} gives a brief overview of combinatorial patchworking, including a notion of equivalence of sign distributions that identifies patchworks agreeing up to reflection, considerably reducing the search space.
Section~\ref{sec:haas} revisits Haas' zone decompositions \cite{1997:Haas:AlgebraicCombinatorial}; these are used to exhibit degree-eight $M$-curves which do not arise as T-curves.
In this context we also verify that all maximal T-curves of degree eight satisfy Ragsdale's conjecture (Corollary~\ref{cor:ragsdale}).
Section~\ref{sec:realizations} describes our computational search and gives explicit T-curve realizations degree by degree: for all nonempty real schemes of degree at most six (Proposition~\ref{thm:covering-5} and Theorem~\ref{thm:covering-6}), for all real schemes of degree seven (Theorem~\ref{thm:covering-7}), and partially in degree eight (Section~\ref{sec:deg-8}).
Section~\ref{sec:outlook} concludes with open questions.
Appendix~\ref{app:statistics} provides statistics on the real schemes arising from a fixed triangulation as the signs vary, Appendix~\ref{app:software} describes available software and data and Appendix~\ref{app:known} sketches known constructions of T-curves of degree $d \leq 7$.

\subsection*{Acknowledgment}
We are much indebted to several colleagues for their valuable comments on preliminary versions of this article.
Erwan Brugall\'e provided a sketch for constructing degree seven curves as T-curves; see Appendix~\ref{app:known}.
Ilia Itenberg and an anonymous referee directed our attention to the work of Haas \cite{1997:Haas:AlgebraicCombinatorial}; this led to a near-complete rewrite and several new results.
With Bernd Sturmfels we discussed discriminants and the connection with \cite{2019:KaihnsaEtAl:SixtyCurves}.
This work was funded by the Deutsche Forschungsgemeinschaft (DFG, German Research Foundation) under Germany's Excellence Strategy -- \enquote{The Berlin Mathematics Research Center MATH$^+$} (EXC-2046/1, EXC-2046/2, pro\-ject ID 390685689), \enquote{Symbolic Tools in Mathematics and their Application} (TRR 195, project ID 286237555), \enquote{Mathematical Modelling, Simulation and Optimization Using the Example of Gas Networks} (SFB/TRR 154, project ID 239904186), and \enquote{Mathematical Research Data Initiative (MaRDI)} (project ID 460135501), as well as by the German Federal Ministry of Research, Technology and Space (Research Campus MODAL, grant numbers 05M14ZAM and 05M20ZBM) and VDI/VDE Innovation + Technik GmbH (grant number 16IS23025B).

\section{Viro's patchworking}
\label{sec:viro}
In this section we will fix our notation and recall the combinatorial framework for patchworked curves.
A \emph{real plane algebraic curve} of degree $d$ is the zero set $V_\RR(f) = \{f=0\} \subset \RP^2$ of a homogeneous polynomial $f \in \RR[x,y,z]$ of degree $d$, and its \emph{complexification} is denoted $V_\CC(f) = \{f=0\} \subset \CP^2$.
Throughout, all real curves are assumed to be \emph{smooth}, meaning that $V_\CC(f)$ is nonsingular.

\subsection{Real schemes and related classifications}
\label{subsec:real-schemes}
The connected components of $V_\RR(f)$ are called \emph{loops}; due to our smoothness assumption, each component is an embedding of the circle $\Sph^1$ into $\RP^2$.
A loop is called an \emph{oval} if it separates $\RP^2$ into a disk and a Möbius strip, and a \emph{pseudo-line} otherwise.
When $d$ is even, every loop is an oval; when $d$ is odd, there is exactly one pseudo-line.
By Harnack's theorem~\cite{1876:Harnack:VieltheiligkeitEbenen}, the number of loops is at most $M = \tfrac{1}{2}(d-1)(d-2)+1$.
The \emph{depth} of an oval is the number of other ovals surrounding it. An oval is \emph{even} (resp.\ \emph{odd}) if it lies in the interior of an even (resp.\ odd) number of other ovals; we write $p$ and $n$ for the number of even and odd ovals, respectively. 
For an $M$-curve of even degree $d=2k$, the Gudkov--Rokhlin congruence \cite{1971:Arnold:ArrangementOvals,1972:Rokhlin:GudkovsHypothesis} states that
\begin{equation}\label{eq:gudkov}
  p-n\equiv k^2\pmod 8 \enspace .
\end{equation}
The \emph{real scheme} of a curve is the ambient isotopy class of $V_\RR(f) \subset \RP^2$, or equivalently, the topological type of the pair $(\RP^2, V_\RR(f))$.
Concretely, the real scheme records the number of connected components and their nesting structure.
In the literature on Hilbert's 16th problem, several terms are used interchangeably for this notion: Viro~\cite{1986:Viro:TopologyAlgebraic} uses both \emph{isotopy type} and \emph{real scheme}; Itenberg--Viro~\cite{1995:Itenberg:RagsdaleCurves, 1996:ItenbergViro:PatchworkingRagsdale} speak of the \emph{topological type} of the pair.
We adopt \emph{real scheme} throughout to avoid confusion with the finer notion of \emph{rigid isotopy type} discussed below.

We use the \emph{Rokhlin--Viro notation} to encode real schemes:
$\scheme{1\scheme{X}}$ denotes an oval containing scheme $\scheme{X}$, and $\scheme{X \sqcup Y}$ two separable schemes.
The base case depends on parity: for even $d$, the empty scheme is denoted $\scheme{0}$; for odd $d$, a lone pseudo-line is denoted $\scheme{J}$.
The shortcut $\scheme{k\scheme{X}}$ means $k$ disjoint ovals each containing $\scheme{X}$, and $\scheme{k} = \scheme{k\scheme{0}}$.
This notation views a real scheme as a rooted tree describing the nesting structure of ovals.
For instance, the unique real scheme for $d=1$ is $\scheme{J}$; for $d=2$ we have $\scheme{0}$ and $\scheme{1}$; and for $d=3$ the real schemes are $\scheme{J}$ and $\scheme{J \sqcup 1}$.

There are finer equivalence relations on curves than ambient isotopy, which defines the real schemes considered here.
A \emph{rigid isotopy} between two smooth degree-$d$ curves is a continuous path of smooth degree-$d$ curves connecting them, i.e., a path in the complement of the discriminant; cf.\ \cite[\S11.5]{1994:GKZ:DiscriminantsResultants}.
The \emph{rigid isotopy type} of a curve is its equivalence class under this relation~\cite{1986:Viro:TopologyAlgebraic, 2019:KaihnsaEtAl:SixtyCurves}.
For $d = 6$, there are $56$ real schemes but $64$ rigid isotopy types (including the empty scheme)~\cite[Theorem 3.4.3]{1980:Nikulin:SymmetricBilinear}; the difference arises from eight real schemes that admit both \emph{dividing} (type~I) and \emph{non-dividing} (type~II) realizations, which are not rigidly isotopic~\cite{1978:Rokhlin:TopologicalAlgebraic, 1980:Nikulin:SymmetricBilinear}; for $d = 7$, the rigid isotopy classification remains open.
Another refinement is the \emph{complex scheme}, which for dividing curves additionally records the complex orientations of ovals~\cite{1978:Rokhlin:TopologicalAlgebraic}.
These orientations are invariant under rigid isotopy, so rigidly isotopic curves have the same complex scheme.

\subsection{Combinatorial patchworking}
\label{subsec:patchworking}

We denote the \emph{standard triangle} by $\Delta_2 = \conv\{(0,0),(1,0),(0,1)\}$ and its \emph{dilation} by
\begin{equation}
  d\cdot\Delta_2=\conv\{(0,0),(d,0),(0,d)\} \enspace .
\end{equation}
Note that the lattice point set $A$ from the introduction equals $d\cdot\Delta_2\cap\ZZ^2$.
Now let $\cT$ be a triangulation of $A$, i.e., $\cT$ is a triangulation of $d\cdot\Delta_2$ such that the vertices form a subset of $A$.
In our examples we will usually describe $\cT$ in terms of its edges.

We take four copies of $\cT$, which we call $\cT_{++}, \cT_{+-}, \cT_{-+}, \cT_{--}$, to build a triangulation of the real projective plane $\RP^2$.
As the vertex set we take the \emph{diamond} $A^\diamond=\smallSetOf{(i,j)\in\ZZ^2}{|i|+|j|\leq d}$.
Observe that the set $A$ is precisely the intersection of $A^\diamond$ with the positive quadrant.
We let $\cT_{++}$ be $\cT$.
Further, $\cT_{-+}$ is the reflection of $\cT$ across the $y$-axis, and $\cT_{+-}$ the reflection across the $x$-axis.
Finally, $\cT_{--}$ is the reflection of $\cT_{-+}$ across the $x$-axis or, equivalently, the reflection of $\cT_{+-}$ across the $y$-axis.
Now the union $\cT^\diamond = \cT_{++} \cup \cT_{+-} \cup \cT_{-+} \cup  \cT_{--}$ forms a triangulation of $A^\diamond$.
By identifying the boundary points $(i,d-i)$ with $-(i,d-i)$ and $(i,i-d)$ with $-(i,i-d)$, for $0 \leq i \leq d$, we obtain a cell decomposition $\cS=\cS(\cT)$ of $\RP^2$.
Although the maximal cells are triangles, the cell complex $\cS$ is not necessarily a simplicial complex.
This happens if and only if $\cT$ has an edge connecting one of the corner vertices of $d\cdot\Delta_2$, namely, $(0,0)$, $(0,d)$, or $(d,0)$, with some other boundary vertex.
Further, we consider a sign function $\sigma:A\to\GF{2}=\ZZ/2\ZZ$, which is extended to the entire diamond $A^\diamond$ by the rule
\begin{equation}\label{eq:signs}
  \sigma(i,-j) = \sigma(i,j) + j \,,\quad \sigma(-i,j) = \sigma(i,j) + i \,,\quad \sigma(-i,-j) = \sigma(i,j) + i + j  \enspace ,
\end{equation}
where $(i,j)\in A$ and the addition is taken modulo two.
These signs allow us to define the \emph{patchworked curve} of the \emph{patchwork} $(\cT,\sigma)$.
To this end, we consider the dual graph $\Gamma^*(\cS)$ of the polygonal surface $\cS$. Its nodes are the maximal cells of $\cS$, and two nodes are joined by a (dual) edge if the corresponding triangles share an edge in $\cS$.
Now the patchworked curve $\cC(\cT,\sigma)$ is the subgraph of $\Gamma^*(\cS)$ comprising the edges dual to edges of $\cS$ whose endpoints have distinct signs.
The dual graph $\Gamma^*(\cS)$, and thus its subgraph $\cC(\cT,\sigma)$, admit an embedding into the first barycentric subdivision $\beta\cS$.

The triangulation $\cT$ is \emph{unimodular} (or \enquote{primitive}) if each triangle has Euclidean area $\tfrac{1}{2}$.
Equivalently, each point in $A$ occurs as a vertex of $\cT$.
The triangulation $\cT$ is \emph{regular} (or \enquote{coherent} or \enquote{convex}) if it is induced by a height function such that $\cT$ arises by projecting the lower convex hull of the lifted points; for details, see \cite[\S7.1.C]{1994:GKZ:DiscriminantsResultants}, \cite[\S2.2.3]{2010:DeLoeraRambauSantos:TriangulationsAlgorithms}, or \cite[\S1.2]{2021:Joswig:EssentialsCombinatorics}.
The following key result holds more generally for hypersurfaces in arbitrary dimension, even without the unimodularity assumptions.
It was obtained by Viro \cite{1984:Viro:GluingAlgebraic}; for a proof in the context of toric varieties and discriminants, see Gel{'}fand, Kapranov, and Zelevinsky \cite[Theorem 5.6]{1994:GKZ:DiscriminantsResultants}; see also Deng, Rojas, and Telek \cite{2026:DengEtAl:VirosPatchworking}.

\begin{figure}[t]
  \newcommand{\scalefactor}{0.6}
  \begin{tabular}{llll}
  \begin{tikzpicture}[scale=\scalefactor]\input{tikz/deg2-111.tikz}\end{tikzpicture} &
  \begin{tikzpicture}[scale=\scalefactor]\input{tikz/deg2-011.tikz}\end{tikzpicture} &
  \begin{tikzpicture}[scale=\scalefactor]\input{tikz/deg2-101.tikz}\end{tikzpicture} &
  \begin{tikzpicture}[scale=\scalefactor]\input{tikz/deg2-001.tikz}\end{tikzpicture} \\[1.2ex]
  \begin{tikzpicture}[scale=\scalefactor]\input{tikz/deg2-110.tikz}\end{tikzpicture} &
  \begin{tikzpicture}[scale=\scalefactor]\input{tikz/deg2-010.tikz}\end{tikzpicture} &
  \begin{tikzpicture}[scale=\scalefactor]\input{tikz/deg2-100.tikz}\end{tikzpicture} &
  \begin{tikzpicture}[scale=\scalefactor]\input{tikz/deg2-000.tikz}\end{tikzpicture}
  \end{tabular}
  \caption{Eight T-curves of degree two arising from one unimodular triangulation of $2\cdot\Delta_2$.
    In each case the real scheme is $\scheme{1}$, and the interior of the unique oval is drawn shaded.}
  \label{fig:deg2-eight-curves}
\end{figure}

\begin{theorem}[Viro's Combinatorial Patchworking Theorem]\label{thm:patchworking}
  Let $\cT$ be a regular and unimodular triangulation of $A=d\cdot\Delta_2\cap\ZZ^2$ with lifting function $\omega:A\to\ZZ$, and let $\sigma:A\to\GF{2}$ be a sign distribution.
  Then there is a real number $t_0>0$ such that for all $t\in(0,t_0]$ the real projective plane curve $V_\RR(f_t)\subset\RP^2$, where
  \begin{equation}
    f_t(x,y,z) \ = \ \sum_{(i,j)\in A} (-1)^{\sigma(i,j)} t^{\omega(i,j)} x^iy^jz^{d-i-j} \enspace ,
  \end{equation}
  is ambient isotopic to the patchworked curve $\cC(\cT,\sigma)$.
\end{theorem}
We call the patchworked curve $\cC(\cT,\sigma)$ of a unimodular triangulation $\cT$ and a sign distribution $\sigma$ a \emph{T-curve}; see Figure~\ref{fig:deg2-eight-curves} for examples.
If $\cT$ is regular, Theorem~\ref{thm:patchworking} realizes the T-curve as a real algebraic curve $V_\RR(f)$ whose complex curve $V_\CC(f)$ is smooth; unless stated otherwise, our T-curves are regular, and we flag the occasional \emph{nonregular} ones explicitly.

\begin{remark}\label{rem:lex-ordering}
  Whenever we write sign distributions as vectors, we refer to the lexicographic ordering of $A$ from low to high, i.e., the points are ordered $(0,0)$, $(0,1)$, $\ldots$, $(0,d)$, $(1,0)$, $(1,1)$, $\ldots$, $(1,d-1)$, $(2,0)$, $\ldots$, $(d,0)$.
  This convention explains our notation in Tables~\ref{tab:deg-6} and~\ref{tab:deg-7} below.
\end{remark}

\subsection{Combinatorial topology}
\label{subsec:combinatorial-topology}

Since a patchworked curve is determined by the combinatorial data $(\cT,\sigma)$, the above topological information is also encoded combinatorially. Consider a not necessarily regular subdivision $\cT$ of $A$, a sign function $\sigma$, the induced patchworked curve $\cC=\cC(\cT, \sigma)$, and the cell decomposition $\cS=\cS(\cT)$ of $\RP^2$.
We can view both $\cS$ and its barycentric subdivision $\beta\cS$ as \emph{combinatorial surfaces} in the sense of Armstrong \cite[p.154]{1983:Armstrong:BasicTopology}.
A curve formed from edges of a combinatorial surface is called \emph{polygonal}. In this way, we may view any closed loop in $\cS$, as well as every connected component of the patchworked curve, as a polygonal curve in $\RP^2$.
Of course, the topological distinction between ovals and nonseparating curves applies to the polygonal case.

\begin{proposition}\label{prop:disjoint-curves}
  Consider two simple closed polygonal curves in the surface $\cS$ which are disjoint.
  Then at most one of them is nonseparating.
  In particular, every T-curve, whether regular or not, has at most one nonseparating loop.
\end{proposition}

\begin{proof}
  The Euler characteristic of $\RP^2$ equals one.
  Surgery along a simple closed polygonal curve increases the Euler characteristic; see \cite[Theorem 7.11]{1983:Armstrong:BasicTopology}.
  Consequently, after surgery along one nonseparating polygonal curve, we arrive at the sphere $\Sph^2$, with Euler characteristic two, which has no nonseparating closed polygonal curves by the Jordan curve theorem.
\end{proof}

\begin{remark}\label{rem:non-regular}
  Itenberg shows that the Harnack bound \eqref{eq:M} also applies to nonregular T-curves \cite[\S7]{1995:Itenberg:RagsdaleCurves}.
  Renaudineau and Shaw generalized the Harnack bound for curves to regular patchworked hypersurfaces \cite{2023:RenaudineauShaw:BettiHypersurfaces}.
  Brugall\'e, L\'opez de Medrano, and Rau further extended this result to arbitrary nonregular T-hypersurfaces \cite{2024:BrugalleEtAl:CombinatorialPatchworking}.
\end{remark}

\subsection{Regions and nesting}
\label{subsec:regions-nesting}

The $1$-skeleton of the triangulation $\cT^\diamond$ of $A^\diamond$ is a graph $\cT^\diamond_{\leq1}$, and the sign distribution induces a bicoloring on it.
This bicoloring is not proper, meaning that there are monocolored edges.
The antipodal identification $\equiv$ at the boundary of $A^\diamond$ induces the graph $\cS_{\leq 1}=\cT^\diamond_{\leq 1}/\equiv$, a subgraph of the $1$-skeleton of the cell decomposition $\cS$.
For odd $d$, the coloring is not well defined on $\cS_{\leq 1}$, because antipodal boundary points receive different signs by \eqref{eq:signs}.
However, the signs of the two endpoints of an identified boundary edge flip simultaneously.
Therefore, each edge of $\cS_{\leq 1}$ is unambiguously monocolored or bicolored.
For even $d$ the coloring is well defined. 
Removing all the bicolored edges, i.e., those intersected by components of $\cC$, leaves connected subgraphs of $\cS_{\leq 1}$ which we call the \emph{regions} of $\cS$ induced by $\sigma$.
For even $d$, and every region consists of vertices with the same color.
In particular, this implies that for even $d$, every component of $\cC$ is separating and hence an oval.

Let $R$ be a region. If $R$ lies inside an oval $O$, where $O$ is the innermost oval with this property, then we say that $O$ \emph{corresponds} to $R$. If there is no oval corresponding to $R$, we call $R$ the \emph{root region}.
When $d$ is odd, the unique nonseparating loop of the T-curve lies in the \enquote{closure} of the root region.

We say that two regions $R\neq S$ are \emph{neighbors} if there is an edge in $\cS_{\leq1}$ connecting a point in $R$ and a point in $S$.
Further, we call a region $S$ \emph{nested} in a region $R$, and write $R>S$, if either $R$ is the root region or $R$ is proper and corresponds to an oval whose interior contains the oval corresponding to $S$.
If additionally $R$ and $S$ are neighbors, we say $R$ is the \emph{parent} of $S$, and $S$ is a \emph{child} of $R$.
In particular, every proper region has exactly one parent, and neighboring regions are always comparable under $>$.
The \emph{depth} of a proper region is the depth of the corresponding oval, and the depth of the root region is $-1$. In particular, the depth of a child is always one more than the depth of its parent.
This partial ordering of the regions is exactly the Rokhlin--Viro notation encoded as a rooted tree, and the root region is the maximal element of that partial order.
If two patchworked curves $\cC(\cT,\sigma)$ and $\cC(\cT',\sigma')$ are isotopic, and if additionally there is some vertex in $A^\diamond/{\equiv}$ which for both curves lies in the respective root region, then we call $\cC(\cT,\sigma)$ and $\cC(\cT',\sigma')$ \emph{root isotopic}.
Recall that the two cellular surfaces $\cS(\cT)$ and $\cS(\cT')$ share the same vertex set $A^\diamond/{\equiv}$.

\subsection{Equivalence of sign distributions}
\label{sec:signs}
It is easy to see that for any triangulation there are always distinct sign distributions giving the same patchworked curve.
This comes from underlying symmetries, which we explore next.
While the results of this section are known, we include a brief sketch for the sake of completeness.

A \emph{unimodular transformation} of $\RR^2$ is an affine-linear transformation of the plane which leaves the integer lattice $\ZZ^2$ invariant.
For a finite set $X\subset\ZZ^2$ of lattice points, the \emph{group of symmetries} $\Aut(X)$ is the group of unimodular transformations fixing $X$.
The unimodular affine maps are precisely those that are area-preserving.
By construction, $\Aut(X)$ is always finite.
Again, we fix $d\geq 1$ and set $A = d \cdot \Delta_2\cap\ZZ^2$.
\begin{lemma}\label{lem:dihedral}
  The group $\Aut(A)$ is isomorphic to the symmetric group $\Sym{3}$ permuting the three vertices of the triangle $d\cdot\Delta_2$.
  The group $\Aut(A^\diamond)$ is the dihedral group of order eight. It contains the four maps
  \begin{equation}\label{eq:sign-reflections}
    s_{ij}\colon (x,y)\longmapsto \bigl((-1)^{i}x,\,(-1)^{j}y\bigr)\,, \qquad i,j\in\GF{2}\,,
  \end{equation}
  which form a Klein four-group; here $s_{10}$ and $s_{01}$ are the reflections in the $y$- and $x$-axes, $s_{11}$ is their composition, and $s_{00}$ the identity. Together with the diagonal reflection $t\colon(x,y)\mapsto(y,x)$, they generate $\Aut(A^\diamond)$.
\end{lemma}

For $\Sigma$ a subdivision of a finite set $X\subset\ZZ^2$, we define $\Aut(X,\Sigma)$ as the subgroup of $\Aut(X)$ which leaves $\Sigma$ invariant.
As before, we now consider a (unimodular) triangulation $\cT$ of $A$ and the induced triangulation $\cT^\diamond$ of $A^\diamond$.
\begin{lemma}\label{lem:symmetries-diamond}
  The group $\Aut(A^\diamond,\cT^\diamond)$ equals $\Aut(A^\diamond)=\langle s_{10},t\rangle$ if $\cT$ is symmetric with respect to the line $x=y$, and the Klein four-group $\{s_{ij}: i,j\in\GF{2}\}$ otherwise.
\end{lemma}
We call two sign distributions $\sigma,\sigma':A\to\GF{2}$ \emph{equivalent} if $\sigma'(u)=\epsilon+\sigma(s_{ij}(u))$ for all $u\in A$ and some $\epsilon,i,j\in\GF{2}$; as before we use the same symbols for the extensions $\sigma,\sigma':A^\diamond\to\GF{2}$.
For instance, if $s_{ij}$ is the identity (that is, $i=j=0$) and $\epsilon=1$, then all signs are flipped.
When we want to enumerate all patchworked curves coming from a fixed triangulation, it suffices to pick one sign distribution per equivalence class of signs.
This observation is the content of the next result.
\begin{proposition}
  If $\sigma,\sigma':A\to\GF{2}$ are equivalent sign distributions, then for any unimodular triangulation $\cT$ of $A$ the patchworked curves $\cC(\cT,\sigma)$ and $\cC(\cT,\sigma')$ are linearly isomorphic and thus, in particular, ambient isotopic.
\end{proposition}

\begin{proposition}\label{prop:canonical-signs}
  Each equivalence class of sign distributions has exactly eight elements.
\end{proposition}

The number of equivalence classes in degree $d$ equals
\begin{equation}
  2^{\tbinom{d+2}{2}-3} \enspace.
\end{equation}
For $2 \leq d \leq 8$ these values are listed in the last column of Table~\ref{tab:search-space}.

We pick $\sigma(0)=\sigma(e_1)=\sigma(e_2)=1$ for the \emph{canonical representative} of an equivalence class of sign distributions.
Figure~\ref{fig:deg2-eight-curves} shows those $2^{6-3}=8$ canonical representatives for one fixed triangulation of $2\cdot\Delta_2$.
The resulting patchworked curves share the real scheme $\langle 1\rangle$.
Taking the symmetries of $\cT$ into account, there may be further sign distributions that can be eliminated in an enumeration.

\section{Haas' zone decompositions}
\label{sec:haas}

In his PhD thesis \cite{1997:Haas:AlgebraicCombinatorial} Haas developed a powerful combinatorial framework for analyzing maximal T-curves.
Our non-realizability results crucially rely on Haas' ideas.
However, we cast his ideas in a slightly different language and begin with a discussion of a particular class of sign distributions and Harnack splits.

\subsection{Harnack signs} 
The \emph{parity} of a point $u=(x,y)\in A^\diamond$ is the tuple
\[
   \parity(u)=(x\bmod{2}, y\bmod{2})\in\FF_2^2 \enspace.
\] 
A point with parity $(0,0)$ is called \emph{even}; a point with any other parity is called \emph{odd}.
For arbitrary $d\geq 1$, let $A=d\cdot\Delta_2\cap\ZZ^2$. Following Itenberg \cite[\S4]{1995:Itenberg:RagsdaleCurves}, the \emph{Harnack sign distribution} $\eta:A\to\GF{2}$ is given by
\begin{equation}\label{eq:harnack-signs}
	\eta(u) = (x+1)(y+1) \,, \quad \text{for } u=(x,y)\in A \enspace .
\end{equation}
That is, $\eta$ assigns $1$ to each even point, and all odd points receive $0$.
Observe that the signs in the other three quadrants of $A^\diamond$ are similar.
For instance, a point in the upper left quadrant receives the sign $0$ if and only if its parity is $(0,1)$.
So the sign distribution on the vertices of $A^\diamond$ in the upper left quadrant, up to translation in the $y$-direction, looks like the inverse sign distribution $\1+\eta$, which is equivalent to~$\eta$.
We use the term \emph{Harnack sign distribution} also for any sign function in the equivalence class of $\eta$.

\begin{example}\label{exmp:honeycomb}
  The affine Coxeter arrangement of type A$_2$ is the infinite affine line arrangement in $\RR^2$ given by
  \begin{equation}
    \label{eq:arrangement:A}
    H_{ij}^\ell \ = \ \bigl\{ (z_1, z_2) \in\RR^2 \bigm| z_i-z_j =\ell \bigr\} \enspace \text{ for } \ell\in\ZZ \enspace ,
  \end{equation}
  where $0\leq i < j \leq 2$ and $z_0=0$; we have three parallel classes of lines: $H_{01}^\ell$, $H_{02}^\ell$, and $H_{12}^\ell$.
  The chambers of this line arrangement are translations of the standard triangle $\Delta_2$, and so we obtain a unimodular triangulation of the entire plane $\RR^2$.
  Restricting to $d\cdot\Delta_2$ we get a unimodular triangulation of that scaled triangle, called the \emph{honeycomb triangulation} of degree $d$, denoted $\cH_d$.
  The honeycomb triangulation is ubiquitous in geometric combinatorics and algebraic geometry; see, e.g., \cite[Remark 10.33]{2008:AbramenkoBrown:BuildingsTheory}.
  It is known to be regular; a lifting function with minimal integral values has been determined in \cite[Remark 15]{2024:CasabellaEtAl:WronskiHoneycomb}.
  For example, equipping $\cH_6$ with the Harnack sign distribution yields a maximal T-curve of degree six and type $\scheme{9\sqcup 1\scheme{1}}$; see Figure~\ref{fig:harnack-curve-6}.
\end{example}  

More generally, any patchworked curve whose sign distribution is Harnack is an $M$-curve.
Itenberg gives a proof for even degree \cite[\S4]{1995:Itenberg:RagsdaleCurves} that applies to odd degrees with minor adjustments.

\begin{proposition}\label{prop:special-harnack}
	Any unimodular triangulation $\cT$ of $A$ together with the Harnack sign distribution $\eta:A\to\GF{2}$ yields an $M$-curve of degree $d$ with real scheme
	\begin{align*}
		\scheme{\tfrac{3}{2}(k^2-k)\sqcup 1\scheme{\tfrac{1}{2}(k-1)(k-2)}} \,, \quad & \text{if } d=2k \text{ is even}\,,\\
		\scheme{J\sqcup \tfrac{1}{2}(d-1)(d-2)} \,, \quad  & \text{if $d$ is odd}\enspace .
	\end{align*}
	In particular, we have $p=\tfrac{3}{2}(k^2-k)+1$ and $n=\tfrac{1}{2}(k-1)(k-2)$ for $d$ even, and $(p,n)=(\tfrac{1}{2}(d-1)(d-2),0)$ if $d$ is odd.
	It is not necessary for $\cT$ to be regular.
\end{proposition}

For $d=6$, we get $p=10$, $n=1$, and the real scheme $\scheme{9\sqcup1\scheme{1}}$, as we saw in Example~\ref{exmp:honeycomb}; for $d=7$, we have $(p,n)=(15,0)$ and $\scheme{J\sqcup 15}$; for $d=8$, this becomes $(p,n)=(19,3)$ and $\scheme{18\sqcup1\scheme{3}}$.
Such curves were first constructed by Harnack \cite{1876:Harnack:VieltheiligkeitEbenen}.
Today the term \enquote{Harnack curve} usually refers to Mikhalkin's generalization~\cite{2000:Mikhalkin:AlgebraicAmoebas}.
The T-curves described here correspond to the original construction of Harnack; for simplicity we call them \emph{special Harnack curves}.

\begin{figure}[t]
	\newcommand{\scalefactor}{0.4}
	\centering
		\begin{tikzpicture}[scale=\scalefactor]
			\input{tikz/harnack_honeycomb_6.tikz}
		\end{tikzpicture}
		\caption{Degree six special Harnack curve, supported by the honeycomb triangulation.}
		\label{fig:harnack-curve-6}
\end{figure}

\subsection{Harnack splits}
Let $\alpha,\beta\in\FF_2^2$ be two distinct parity vectors. 
An \emph{$\{\alpha,\beta\}$-split} is a path in $\cT_{\leq 1}$ with vertices in $A$ connecting precisely two boundary points of $d\cdot\Delta_2$ and consisting of one or two primitive edges such that its vertices have alternating parities $\alpha$ and $\beta$.
A \emph{Harnack split} is an $\{\alpha,\beta\}$-split for some choice of $\alpha$ and $\beta$. 
It is called \emph{simple} if it has exactly one edge and \emph{double} if it has exactly two edges.
The endpoints of a Harnack split are the \emph{base points}; for a double Harnack split, the interior vertex is the \emph{apex}.
An $\{\alpha,\beta\}$-split is \emph{even} if $\alpha$ or $\beta$ is $(0,0)$; otherwise, it is \emph{odd}.
Examples of Harnack splits appear in Figure~\ref{fig:splits_induce-zon_dec_2}.

Every Harnack split $S$, simple or double, defines the decomposition $\cZ(S)=\{Z^+,Z^-\}$ into two \emph{zones} with $Z^+\cup Z^-= d\cdot\Delta_2$ and $Z^+\cap Z^-= S$.
Here $Z^+$ is the \emph{zone surrounded by the split $S$}, which is defined as follows.
The sets $Z^+\cap\{(0,0),(d,0),(0,d)\}$ and $Z^-\cap\{(0,0),(d,0),(0,d)\}$ are the sets of vertices of $d\cdot \Delta_2$ lying in the respective zones.
The zone containing fewer of these vertices is $Z^+$.
If $S$ itself contains such a vertex, it can happen that $Z^+$ and $Z^-$ both contain two vertices; in that case, the designation of $Z^+$ can be chosen arbitrarily.
The other zone, $Z^-$, is the \emph{root zone} of $S$.

When we now consider the copies of $d\cdot\Delta_2$ in the four quadrants, each zone of the split $S$ is copied, too.
In this way, each zone $Z^\pm$ defines a \emph{zone} of the quotient $\RP^2$ which is denoted by $Z^{\pm\diamond}$. 
Note that $Z^{+\diamond}$ may be disconnected, in which case there are precisely two connected components.

A set of Harnack splits is called \emph{weakly compatible} if there is a triangulation containing all of their line segments as edges.
A weakly compatible set is \emph{compatible} if no two Harnack splits in the collection share an edge.

\begin{remark}\label{rem:split}
  Simple Harnack splits are \enquote{splits} in the sense of \cite{2008:HerrmannJoswig:SplittingPolytopes,2009:Herrmann:FacetsPolytope}; see also \cite[\S5.3.3]{2010:DeLoeraRambauSantos:TriangulationsAlgorithms}.
  That is, they are (necessarily coarsest regular) subdivisions of the point configuration $A=d\cdot\Delta_2\cap\ZZ^2$ with exactly two maximal cells, i.e., the zones.
  Similarly, the notions of (weak) compatibility stem from that context.
\end{remark}

The \enquote{zones} of Haas \cite[\S~2.5]{1997:Haas:AlgebraicCombinatorial} are slightly more general than our zones, because of the parity condition that we require for Harnack splits.
In Haas' notation \cite[Definition 11.2.0.11]{1997:Haas:AlgebraicCombinatorial} our root zone is called \enquote{special zone}.
Haas' \enquote{boundary-nodes} and \enquote{inner-nodes} \cite[Definitions 8.3.0.2 and 8.3.0.3]{1997:Haas:AlgebraicCombinatorial} of a zone are, respectively, our base points and the apex of a double split.
Our reason for deviating slightly from various notions in \cite{1997:Haas:AlgebraicCombinatorial} is to highlight the connection to known concepts in polyhedral geometry; cf.\ Remark~\ref{rem:split}.

\subsection{Zone decompositions}

Next we relate our setup to the original framework of Haas \cite{1997:Haas:AlgebraicCombinatorial}.
Let $\frS=\{S_1,S_2,\dots,S_\ell\}$ be a collection of compatible Harnack splits of $A$.
A \emph{zone} of $\frS$ is the nonempty intersection $\bigcap_{i}Z_i^{\pm}$ of zones, where $\cZ(S_i)=\{Z_i^+,Z_i^-\}$.
The collection of all zones of $\frS$ is the \emph{zone decomposition} $\cZ(\frS)$.
We call such a zone decomposition \emph{valid} for some patchwork $(\cT,\sigma)$ if the triangulation $\cT$ refines $\cZ(\frS)$ and $\sigma$ induces a Harnack sign distribution on each zone such that the Harnack sign distributions for neighboring zones are distinct.

An \emph{H-zone decomposition}, meaning a zone decomposition in the sense of Haas, is a cell decomposition of $d\cdot\Delta_2$ in which every edge is a primitive line segment in the lattice $A=d\cdot\Delta_2\cap\ZZ^2$ and has an endpoint on the boundary of $d\cdot \Delta_2$; cf.\ \cite[Definition 6.3.1.1]{1997:Haas:AlgebraicCombinatorial}.
Now, consider an H-zone decomposition $\cZ_H$ and an interior lattice point $p\in A$ that is a vertex of $\cZ_H$.
Let $E(p)$ be the set of edges of zones in $\cZ_H$ that have $p$ as one of their endpoints.
The star of $p$ in the cell decomposition $\cZ_H$ is called a \enquote{cycle of zones} in \cite[Definition 7.1.0.5]{1997:Haas:AlgebraicCombinatorial}; an edge in $E(p)$ is called an \enquote{edge of the cycle}.
Further, let $E(p,\alpha)$ be the set of edges in $E(p)$ whose other endpoints have parity $\alpha\in\FF_2^2$. 
Then $\cZ_H$ is \emph{odd-cycle-free} if for any interior lattice point $p$ and any parity $\alpha$ the cardinality of $E(p,\alpha)$ is even. 

\begin{proposition}\label{prop:algo_zones_are_haas_zones}
  For a compatible collection $\frS$ of Harnack splits, the zone decomposition $\cZ(\frS)$ is an odd-cycle-free H-zone decomposition.
  Conversely, for any odd-cycle-free H-zone decomposition $\cZ_H$, there is a collection of compatible simple and double Harnack splits $\frS$ such that $\cZ_H=\cZ(\frS)$.
\end{proposition}

\begin{proof}
	Let $\frS$ be a collection of compatible simple and double Harnack splits and $\cZ(\frS)$ its zone decomposition. Observe that any zone decomposition is also an H-zone decomposition. 
	An interior lattice point $p\in A$ is either in the interior of some zone or the apex of one or more double splits. 
	In the latter case, there are two edges per double split whose remaining endpoints have the same parity $\alpha\neq \parity(p)$, implying that $\abs{E(p,\alpha)}$ is even.
	Thus, $\cZ(\frS)$ is an odd-cycle-free H-zone decomposition.
	
	Consider now an odd-cycle-free H-zone decomposition $\cZ_H$.
	Then for any interior lattice point $p\in A$ which is a vertex of $\cZ_H$, there is an even number of edges whose remaining endpoints have a given parity $\alpha\neq\parity(p)=:\beta$. 
	Any two of these edges make up an $\{\alpha,\beta\}$-split. 
	On the other hand, every edge of $\cZ_H$ that is not incident to any interior point is a primitive line segment in $A$ and hence a simple split.
	We can thus partition the set of edges in $\cZ_H$ into a collection of compatible simple and double Harnack splits $\frS$ such that $\cZ_H=\cZ(\frS)$. 
\end{proof}

Proposition~\ref{prop:algo_zones_are_haas_zones} allows us to derive the following characterization of maximal T-curves from Haas \cite[Theorem 7.3.0.10]{1997:Haas:AlgebraicCombinatorial}.

\begin{theorem}[Haas]\label{thm:haas}
	Let $\cC(\cT, \sigma)$ be any T-curve. 
	There is a collection $\frS$ of compatible (simple or double) Harnack splits such that the zone decomposition $\cZ(\frS)$ is valid for the patchwork $(\cT, \sigma)$ if and only if the T-curve $\cC(\cT, \sigma)$ is maximal.
\end{theorem}

In this case, we say that $\frS$ \emph{induces} the maximal T-curve $\cC(\cT, \sigma)$. Note that Theorem~\ref{thm:haas} makes no assumption on the regularity of the triangulation $\cT$.
This formulation of Haas' theorem is similar to those in \cite[Theorem 3.13]{2015:BrugalleEtAl:TropicalGeometry} and \cite[Theorem 4.6]{2017:BertrandEtAl:HaasRevisited}.

\subsection{Surgical twists}
We now show that any zone decomposition is valid for some patchwork.
Recall from~\eqref{eq:sign-reflections} the maps $s_{ij}$; a sign function $\sigma$ is equivalent to $\epsilon+\sigma\circ s_{ij}$ for any $\epsilon, i, j\in\FF_2$.
Consider a curve $\cC(\cT, \sigma)$ and an $\{\alpha,\beta\}$-split $S\subset\cT$ with $\alpha=(\alpha_1,\alpha_2)$ and $\beta=(\beta_1,\beta_2)$.
Let $Z^+$ be the zone surrounded by $S$. 
The \emph{twisted sign function} $\sigma'$ is defined by 
\[\FF_2 \ni \sigma'(u)=\begin{cases}
	\sigma(u)\,, &\textrm{ if } u\notin Z^+\,,\\
	\epsilon+\sigma(s_{ij}(u))\,, &\textrm{ if } u\in Z^+\,,
\end{cases}  \]
where $\epsilon=\alpha_1\beta_2+\alpha_2\beta_1$, $i=\alpha_2+\beta_2$, $j=\alpha_1+\beta_1$, and $u\in A$.
Exchanging $\cC(\cT, \sigma)$ for $\cC(\cT, \sigma')$ is a \emph{surgical twist}. The name is inspired by the notion of a \enquote{twist} used in \cite[Definition 5.1.0.4]{1997:Haas:AlgebraicCombinatorial} and adopted, for example, in \cite[Definitions 3.2 and 3.3]{2015:BrugalleEtAl:TropicalGeometry} and \cite[Definition 4.1]{2017:BertrandEtAl:HaasRevisited}. It is also reminiscent of the notion of \enquote{surgery} in \cite[Example 5.8]{1994:GKZ:DiscriminantsResultants}: one cuts off the zone surrounded by the split, reflects it across one or two axes, and pastes it back in; see Figure~\ref{fig:splits_induce-zon_dec_2}. The next lemma shows that this gluing always works.

\begin{lemma}\label{lem:agree-on-split}
The sign functions $\sigma$ and $\epsilon+\sigma\circ s_{ij}$ agree on $S$. No other equivalent sign function $\epsilon'+\sigma\circ s_{i'j'}$, with $(\epsilon',i',j')$ distinct from both $(0,0,0)$ and $(\epsilon,i,j)$, agrees with $\sigma$ on $S$.
\end{lemma}

\begin{proof}
By the extension rule~\eqref{eq:signs} one has $\sigma(s_{ab}(u))=\sigma(u)+ax+by$ for every $u=(x,y)\in A^\diamond$ and $a,b\in\GF{2}$.
Hence $\epsilon'+\sigma\circ s_{i'j'}$ agrees with $\sigma$ at $u$ if and only if $\epsilon'+i'x+j'y=0$.
Every vertex of $S$ has parity $\alpha$ or $\beta$, and $S$ contains at least one of each, so agreement on all of $S$ amounts to the two $\GF{2}$-linear equations
\[\epsilon'+i'\alpha_1+j'\alpha_2=0 \qquad\text{and}\qquad \epsilon'+i'\beta_1+j'\beta_2=0\enspace .\]
Since $\alpha\neq\beta$, the coefficient vectors $(1,\alpha_1,\alpha_2)$ and $(1,\beta_1,\beta_2)$ are distinct and nonzero, hence linearly independent. Thus, the system has exactly $2^{3-2}=2$ solutions.
Both $(0,0,0)$ and $(\epsilon,i,j)=(\alpha_1\beta_2+\alpha_2\beta_1,\,\alpha_2+\beta_2,\,\alpha_1+\beta_1)$ solve it, and they are distinct because $\alpha\neq\beta$ forces $(i,j)=(\alpha_2+\beta_2,\alpha_1+\beta_1)\neq(0,0)$; so these are the only two.
In particular, $\sigma$ and $\epsilon+\sigma\circ s_{ij}$ agree on $S$, while every other equivalent sign function differs from $\sigma$ somewhere on~$S$.
\end{proof}

Applying a surgical twist to the root zone instead of the zone surrounded by the split yields the same result up to reflection or rotation of $\RP^2$.
The triple $(\epsilon,i,j)$ encoding a surgical twist is Haas' \enquote{total parity} of the split.

\begin{figure}[t]
	\newcommand{\scalefactor}{0.4}
	\centering
	\begin{subfigure}{0.48\textwidth}
		\centering
		\begin{tikzpicture}[scale=\scalefactor]
			\input{tikz/collection_of_splits_1.tikz}
		\end{tikzpicture}
		\caption{A zone decomposition induced by three compatible Harnack splits}
		\label{fig:collection_of_splits_1}
	\end{subfigure}
	\hfill
	\begin{subfigure}{0.48\textwidth}
		\centering
		\begin{tikzpicture}[scale=\scalefactor]
			\input{tikz/collection_of_splits_2.tikz}
		\end{tikzpicture}
		\caption{A zone decomposition induced by four compatible Harnack splits}
		\label{fig:collection_of_splits_2}
	\end{subfigure}	
	\caption{Two zone decompositions connected by a surgical twist. The surgical twist corresponds to the densely dashed double $\{(1,0),(0,0)\}$-split and affects the gray area. 
    Two other double splits share an apex in such a way that their respective zones intersect.}
	\label{fig:splits_induce-zon_dec_2}
\end{figure}

\begin{lemma}
Two compatible surgical twists commute. That is, for a sign function $\sigma$ and two triples $(\epsilon, i, j)$ and $(\epsilon', i', j')$ in $\FF_2^3$, we have
\[\epsilon'+(\epsilon+\sigma\circ s_{ij})\circ s_{i'j'}=\epsilon+(\epsilon'+\sigma\circ s_{i'j'})\circ s_{ij}\enspace .\]
In particular, the surgical twists corresponding to two compatible Harnack splits in $\cT$ can be applied to $\cC(\cT, \sigma)$ in any order, resulting in the same patchworked curve.
\label{lem:twists_commute}
\end{lemma}

\begin{proof}
Let $u=(x,y)\in A$. Then we have 
\begin{align*}
	(\epsilon'+(\epsilon+\sigma\circ s_{ij})\circ s_{i'j'})(u) = 
	&\epsilon'+\epsilon+(\sigma\circ (s_{ij}\circ s_{i'j'}))(u)\\
	= & \epsilon'+\epsilon+\sigma((-1)^{i+i'}x, (-1)^{j+j'}y) \\
	= & \epsilon+\epsilon'+\sigma((-1)^{i'+i}x, (-1)^{j'+j}y) \\
	= & \epsilon+\epsilon'+(\sigma\circ (s_{i'j'}\circ s_{ij}))(u)\\
	=& (\epsilon+(\epsilon'+\sigma\circ s_{i'j'})\circ s_{ij})(u) \enspace .
\end{align*}
\end{proof}

\begin{lemma}
	Any zone decomposition $\cZ(\frS)$ is valid for some patchwork $(\cT, \sigma)$.
	\label{lem:splits_induce-zon_dec}
\end{lemma}

\begin{proof}
Fix an ordering of $\frS$, say $S_1,S_2, \dots, S_n$.
Let $\cT$ be a triangulation refining $\frS$, and consider the special Harnack curve $\cC(\cT, \eta)$.
Set $\sigma_0=\eta$, and, for each $i$, let $\cC(\cT, \sigma_i)$ be the curve obtained from $\cC(\cT, \sigma_{i-1})$ by the surgical twist induced by $S_i$. Set $\sigma=\sigma_n$; by Lemma~\ref{lem:twists_commute}, it is independent of the chosen ordering.
Recall that $\cZ(\frS)$ is the zone decomposition consisting of all intersections $Z=Z_1\cap\dots\cap Z_n$, where $Z_i$ is either of the zones induced by $S_i$. 
Since a surgical twist preserves the original sign function up to equivalence inside its corresponding zones, $\sigma$ restricted to a zone of $\cZ(\frS)$ coincides with some Harnack sign distribution.
On the other hand, two neighboring zones of $\cZ(\frS)$ are separated by exactly one split of $\frS$, whose surgical twist changes the sign function precisely at the points $u=(x,y)$ with $\epsilon+ix+jy=1$; see the proof of Lemma~\ref{lem:agree-on-split}.
Since $(i,j)\neq(0,0)$ and each zone contains a triangle of $\cT$, every zone contains such a point, so $\sigma$ restricts to different Harnack sign distributions on neighboring zones.
Thus, $\cZ(\frS)$ is valid for the patchwork $(\cT, \sigma)$.
\end{proof}

\subsection{The pair \texorpdfstring{$(p,n)$}{(p,n)} depends only on odd splits}

Together, Haas' theorem (Theorem~\ref{thm:haas}) and Lemma~\ref{lem:splits_induce-zon_dec} tell us that a maximal T-curve $\cC$ is defined by a collection $\frS$ of Harnack splits applied to a special Harnack curve $\cC(\cT, \eta)$; we say that the maximal T-curve $\cC(\frS)\coloneq\cC$ is \emph{induced} by $\frS$.
The notation $\cC(\frS)$ is independent of $\cT$, as every unimodular triangulation that refines $\frS$ will do. 

Haas also shows \cite[Theorem 10.6.0.5]{1997:Haas:AlgebraicCombinatorial} that we can ignore the even splits in $\frS$ when counting the number of even and odd ovals.

\begin{theorem}[Haas]\label{thm:haas_even}
	Consider a maximal T-curve of even degree $d$,  $\cC$ induced by a collection $\frS$ of Harnack splits. Then the maximal T-curve $\cC'$ induced by the collection obtained from $\frS$ by removing all even splits has the same numbers of even and odd ovals as $\cC$.
\end{theorem}

Indeed, an even $\{\alpha,(0,0)\}$-split has twist parameter $\epsilon=\alpha_1\cdot 0+\alpha_2\cdot 0=0$, so its surgical twist is a pure reflection; reflections preserve the parity of every lattice point and the nesting of the ovals, hence leave $(p,n)$ unchanged.

When investigating which pairs $(p,n)$ are possible for a maximal T-curve of fixed degree $d$, we need only enumerate all odd splits in $d\cdot \Delta_2$ and determine the effect of their induced surgical twists on a patchworked curve.
We now specialize to degree eight.
For any maximal T-curve $\cC$ of degree eight, we have $p+n=M=22$, so $(p,n)=(19-m,3+m)$ for some integer $m$. We say that $m$ is the \emph{effect} of $\frS$ on $(p,n)$, where $\frS$ is the collection of Harnack splits that induces $\cC$.
If $\frS$ contains only one split $S$, we say that $m$ is the \emph{effect} of $S$.

\begin{lemma}
For $d=8$ there are, up to equivalence, only seven odd simple or double Harnack splits whose surgical twists, when applied to a special Harnack curve, change the pair $(p,n)$.
\label{lem:odd_splits}
\end{lemma}

\begin{figure}
	\newcommand{\scalefactor}{0.35}
	\centering	
\begin{tikzpicture}[scale=\scalefactor]
	\begin{scope}[shift={(-9,0)}]
		\clip (0,9) -- (9,0) -- (0,-9) -- (-9,0) -- cycle;
		\input{tikz/simple_split_4.tikz}
	\end{scope}
	
	\begin{scope}[shift={(0,9)}]
		\input{tikz/simple_split_6.tikz}
	\end{scope}
	
	\begin{scope}[shift={(9,0)}]
		\input{tikz/simple_split_12.tikz}
	\end{scope}
	
	\begin{scope}[shift={(-9,-18)}]
		\input{tikz/double_split_4_var.tikz}
	\end{scope}
	
	\begin{scope}[shift={(0,-9)}]
		\input{tikz/double_split_6.tikz}
	\end{scope}
	
	\begin{scope}[shift={(9,-18)}]
		\input{tikz/double_split_8.tikz}
	\end{scope}
	
	\begin{scope}[shift={(0,-27)}]
		\input{tikz/double_split_4_narrow.tikz}
	\end{scope}
\end{tikzpicture}
	\caption{A patchwork for each surgical twist induced by an odd split from Lemma~\ref{lem:odd_splits} when applied to a special Harnack curve of degree eight with a suitable triangulation.}
	\label{fig:odd_splits}
\end{figure}

\begin{proof}
Consider the surgical twist from a special Harnack curve $\cC=\cC(\cT, \eta)$ to $\cC(\cT, \eta')$ induced by an $\{\alpha,\beta\}$-split $S$; here the edges of $S$ occur in $\cT$. 
We may assume that $\{\alpha,\beta\}=\{(0,1),(1,0)\}$ and, for a double split, that the apex has parity $(0,1)$ while the bases have parity $(1,0)$.
Otherwise, we apply a suitable projective transformation.

Thus, $\eta'$ agrees with the sign function $1+\eta\circ s_{11}$ in the zone $Z^+$ surrounded by the split $S$.
In other words, when executing the surgical twist in each of the four copies of $Z^+$, the sign function $\eta$ is exchanged with the inverted sign function of the copy diagonally opposite.
Therefore, within the zone $Z^{+\diamond}$ of $\RP^2$, all even ovals are moved to odd ovals in another quadrant and vice versa.
See Figure~\ref{fig:odd_splits}.

As $\eta$ is the Harnack sign distribution, the number of even (odd) ovals in $Z^{+\diamond}$ is exactly the number of odd (even) interior lattice points in $Z^+$.
Thus, for $\cC(\cT, \eta')$, we have $(p,n)=(19-m, 3+m)$ with
\[m=\abs{\interior(Z^+)\cap\ZZ^2}-2\abs{\interior(Z^+)\cap(2\ZZ)^2}\,, \]
where $\interior(Z^+)$ denotes the interior of $Z^+$. 

Let $S$ first be a simple split. By our assumptions, $S$ has vertices $(0,y)$ and $(x,0)$, where $x,y$ are odd and coprime. We may also assume that $x>y$; otherwise, consider the split between $(0,x)$ and $(y,0)$.
We have $y>1$, as otherwise there are no interior lattice points in $Z^+$. We are left with three options:
\[(x,y)=(5,3)\,,\quad (x,y)=(7,3)\,,\quad \text{and}\quad (x,y)=(7,5)\enspace .\]

Let $S$ now be a double split.
By our assumption, $S$ has vertices $(x_1,0)$, $(x_2,y)$, and $(x_3,0)$, with $x_1,x_3,y$ odd and $x_2$ even.
As before, we have $y>1$.
Three options for $x_1$, three options for $x_3$, and three options for $(x_2,y)$ remain:
\[x_1\in\{1,3,5\},\qquad x_3\in\{3,5,7\},\qquad \text{and}\qquad (x_2,y)\in\{(2,3),(2,5),(4,3)\}.\]
After eliminating the combinations that result in line segments with relative interior lattice points, we are left with four viable options:
\begin{multline*}
  (x_1,x_2,x_3,y)=(3,2,5,5)\,,\quad (3,2,7,3)\,,\quad (1,2,7,3)\,,\quad \text{and}\quad (1,2,5,5)\enspace.
\end{multline*}
There are seven options in total, displayed in Figure~\ref{fig:odd_splits}.
\end{proof}

\begin{remark}\label{rem:effect_of_odd_splits}
For all but one of the possible odd splits in the proof of Lemma~\ref{lem:odd_splits}, we get
\[m=\abs{\interior(Z^+)\cap\ZZ^2}-2\abs{\interior(Z^+)\cap(2\ZZ)^2}=4\,,\]
the remaining one is the simple split with $(x,y)=(7,5)$ and there we have 
\[m=\abs{\interior(Z^+)\cap\ZZ^2}-2\abs{\interior(Z^+)\cap(2\ZZ)^2}=8\enspace .\]
\end{remark}

\begin{remark}\label{rem:zones_nested_or_disjoint}
Suppose two compatible odd simple or double Harnack splits have zones $Z_1$ and $Z_2$ with $Z_2\subset Z_1$ and individual effects $m_1$ and $m_2$ on $(p,n)$. Applying both surgical twists has total effect $m_1-m_2$ because the ovals in $Z_2$ change their parity twice.

If two odd double Harnack splits share an apex such that their zones intersect non-trivially, both contain an edge whose endpoint lies on the same edge of $8\cdot\Delta_2$, so they have the same parity. We may thus interpret these four edges as a pair of splits with nested or disjoint zones instead. The patchworked curve is unaffected.
\end{remark}

\begin{theorem}\label{thm:3_19_impossible}
  The pair $(p,n)$ for a maximal T-curve of degree $d=8$ attains one of the following four values:
  \[(19,3)\,,\ (15,7)\,,\ (11,11)\,,\ (7,15) \enspace.\]
\end{theorem}

\begin{proof}
Let $\cC$ be a degree eight maximal T-curve.
By Theorems~\ref{thm:haas} and~\ref{thm:haas_even}, we may assume that $\cC=\cC(\frS)$ for a collection $\frS$ of odd Harnack splits that have a positive effect on $(p,n)$.
Lemma~\ref{lem:odd_splits} classifies those splits.
By Remark~\ref{rem:zones_nested_or_disjoint}, we can choose the splits such that the corresponding zones are nested or disjoint.
Thus, by Remark~\ref{rem:effect_of_odd_splits}, the total effect of $\frS$ is divisible by four.

The effect is maximal when the zones are disjoint. 
Note that the zone surrounded by any of the splits in Lemma~\ref{lem:odd_splits} contains the point $(4,0)$ on the $x$-axis.
The zone surrounded by the only split whose effect on $(p,n)$ is 8 rather than 4 also contains the point $(0,4)$ on the $y$-axis; cf.\ Remark~\ref{rem:effect_of_odd_splits}.
Thus, for each contribution of at most $4$ to the effect on $(p,n)$, some zone contains the midpoint of a boundary edge of $8\cdot\Delta_2$.
This property is preserved under projective transformations of the splits.
Since we considered disjoint zones, no boundary-edge midpoint is contained in more than one of them.
Thus, the effect is bounded by $12$ because $8\cdot\Delta_2$ has only three boundary edges.
We conclude that the effect of $\frS$ on $(p,n)$ lies in $\{0,4,8,12\}$.
\end{proof}

Note that the proof of Theorem~\ref{thm:3_19_impossible} is purely combinatorial and does not rely on the Gudkov--Rokhlin congruence \eqref{eq:gudkov}.

\begin{corollary}\label{cor:3_19_impossible}
  For $d=8$ there are $M$-curves that do not arise as T-curves.
\end{corollary}
\begin{proof}
   Viro \cite{1980:Viro:CurvesRagsdale} and Shustin \cite{1988:Shustin:CurvesDegree} constructed $M$-curves of degree eight with $(p,n)=(3,19)$.
  Theorem~\ref{thm:3_19_impossible} implies that none of these are realizable as T-curves.
\end{proof}

Table~\ref{tab:deg-8} reproduces Orevkov \cite[Table 1]{2002:Orevkov:CurvesIsotopy}. We make two comments about the table.
The Wiman scheme \cite{1923:Wiman:ReellenZuge} is misprinted there as $\langle 17\sqcup 3\langle 1\rangle\rangle$ (which would amount to $23$ ovals); we correct it to $\langle 16\sqcup 3\langle 1\rangle\rangle$.
The scheme $\langle 17\sqcup 1\langle 2\sqcup 1\langle 1\rangle\rangle\rangle$ is attributed to Viro \cite{1989:Viro:AlgebraicTopology} by Orevkov but to Hilbert \cite{1891:Hilbert:ReellenZuge} in Viro's own $1980$ table.
Concerning credits for the other curves, see \cite{2002:Orevkov:CurvesIsotopy}.
\begin{remark}
  An $(M{-}1)$-curve with real scheme $\scheme{2\sqcup1\scheme{18}}$, i.e., with $(p,n)=(3,18)$, can nevertheless be realized as a T-curve. For example, it is induced by the zone decomposition into three zones meeting at the vertex $(2,2)$ such that two of those zones meet at each of the vertices $(0,3)$, $(3,0)$, and $(3,5)$.
\end{remark}

Ragsdale \cite{1906:Ragsdale:AlgebraicCurves} conjectured that the numbers $p$ and $n$ of even and odd ovals of a real curve of even degree $d$ satisfy
\begin{equation}\label{eq:ragsdale}
  p\leq(3d^2-6d+8)/8 \quad \text{and} \quad  n\leq(3d^2-6d)/8 \enspace .
\end{equation}
For $d=8$ these inequalities read $p\leq 19$ and $n\leq 18$.
The Gudkov--Rokhlin congruence \eqref{eq:gudkov} implies that the only degree-eight $M$-scheme exceeding these bounds is $(3,19)$.
Viro~\cite[\S4]{1980:Viro:CurvesRagsdale} found a real algebraic $M$-curve with those parameters, thus providing a degree-eight counterexample to Ragsdale's conjecture.
However, by Theorem~\ref{thm:3_19_impossible}, that curve is not a T-curve.
\begin{corollary}\label{cor:ragsdale}
  Ragsdale's conjecture \eqref{eq:ragsdale} holds for maximal T-curves of degrees $d\leq 8$.
\end{corollary}

\begin{table}[t]
\centering\footnotesize
\caption{Real schemes of pseudo-holomorphic $M$-curves of degree eight, reproducing Orevkov's Table~1 \cite{2002:Orevkov:CurvesIsotopy}.
  The $38$ real schemes that we can realize as T-curves (see Section~\ref{sec:deg-8} below) are set in \textbf{bold}.
  The schemes in the last column are not realizable as T-curves by Theorem~\ref{thm:3_19_impossible}.
  Six schemes are not known to correspond to a real plane algebraic curve; these are marked with an asterisk.}
\label{tab:deg-8}
\smallskip
\setlength{\tabcolsep}{3pt}
\resizebox{\textwidth}{!}{%
\renewcommand*{\multicitedelim}{\addcomma}%
\medmuskip=1.5mu\relax
\renewcommand{\arraystretch}{1.15}%
\begin{tabular}{lllll}
\toprule
$(p,n)=(19,3)$ & $(p,n)=(15,7)$ & $(p,n)=(11,11)$ & $(p,n)=(7,15)$ & $(p,n)=(3,19)$ \\
\midrule
\realized{\scheme{18\sqcup 1\scheme{3}} } & \realized{\scheme{14\sqcup 1\scheme{7}} } & \realized{\scheme{10\sqcup 1\scheme{11}} } & \realized{\scheme{6\sqcup 1\scheme{15}} } & $\scheme{2\sqcup 1\scheme{19}} $ \\
\realized{\scheme{17\sqcup 1\scheme{1} \sqcup 1\scheme{2}} } & \realized{\scheme{13\sqcup 1\scheme{1} \sqcup 1\scheme{6}} } & \realized{\scheme{9\sqcup 1\scheme{1} \sqcup 1\scheme{10}} } & \realized{\scheme{5\sqcup 1\scheme{1} \sqcup 1\scheme{14}} } & \openmark{\scheme{1\sqcup 1\scheme{1} \sqcup 1\scheme{18}} } \\
\realized{\scheme{16\sqcup 3\scheme{1}} } & \realized{\scheme{13\sqcup 1\scheme{2} \sqcup 1\scheme{5}} } & \realized{\scheme{9\sqcup 1\scheme{2} \sqcup 1\scheme{9}} } & \realized{\scheme{5\sqcup 1\scheme{2} \sqcup 1\scheme{13}} } & $\scheme{1\sqcup 1\scheme{2} \sqcup 1\scheme{17}} $ \\
$\scheme{1\sqcup 1\scheme{2\sqcup 1\scheme{17}}} $ & \realized{\scheme{13\sqcup 1\scheme{3} \sqcup 1\scheme{4}} } & \realized{\scheme{9\sqcup 1\scheme{3} \sqcup 1\scheme{8}} } & \realized{\scheme{5\sqcup 1\scheme{3} \sqcup 1\scheme{12}} } & \openmark{\scheme{1\sqcup 1\scheme{4} \sqcup 1\scheme{15}} } \\
$\scheme{2\sqcup 1\scheme{2\sqcup 1\scheme{16}}} $ & \realized{\scheme{12\sqcup 2\scheme{1} \sqcup 1\scheme{5}} } & \realized{\scheme{9\sqcup 1\scheme{4} \sqcup 1\scheme{7}} } & \realized{\scheme{5\sqcup 1\scheme{4} \sqcup 1\scheme{11}} } & $\scheme{1\sqcup 1\scheme{5} \sqcup 1\scheme{14}} $ \\
$\scheme{3\sqcup 1\scheme{2\sqcup 1\scheme{15}}} $ & \realized{\scheme{12\sqcup 1\scheme{1} \sqcup 2\scheme{3}} } & \realized{\scheme{9\sqcup 1\scheme{5} \sqcup 1\scheme{6}} } & \realized{\scheme{5\sqcup 1\scheme{5} \sqcup 1\scheme{10}} } & \openmark{\scheme{1\sqcup 1\scheme{7} \sqcup 1\scheme{12}} } \\
\openmark{\scheme{4\sqcup 1\scheme{2\sqcup 1\scheme{14}}} } & $\scheme{1\sqcup 1\scheme{6\sqcup 1\scheme{13}}} $ & \realized{\scheme{8\sqcup 2\scheme{1} \sqcup 1\scheme{9}} } & $\scheme{5\sqcup 1\scheme{6} \sqcup 1\scheme{9}} $ & $\scheme{1\sqcup 1\scheme{8} \sqcup 1\scheme{11}} $ \\
$\scheme{5\sqcup 1\scheme{2\sqcup 1\scheme{13}}} $ & $\scheme{2\sqcup 1\scheme{6\sqcup 1\scheme{12}}} $ & \realized{\scheme{8\sqcup 1\scheme{1} \sqcup 1\scheme{3} \sqcup 1\scheme{7}} } & \realized{\scheme{5\sqcup 1\scheme{7} \sqcup 1\scheme{8}} } & \openmark{\scheme{1\sqcup 1\scheme{9} \sqcup 1\scheme{10}} } \\
$\scheme{6\sqcup 1\scheme{2\sqcup 1\scheme{12}}} $ & $\scheme{3\sqcup 1\scheme{6\sqcup 1\scheme{11}}} $ & \realized{\scheme{8\sqcup 1\scheme{1} \sqcup 2\scheme{5}} } & $\scheme{4\sqcup 2\scheme{1} \sqcup 1\scheme{13}} $ & $\scheme{2\scheme{1} \sqcup 1\scheme{17}} $ \\
$\scheme{7\sqcup 1\scheme{2\sqcup 1\scheme{11}}} $ & $\scheme{4\sqcup 1\scheme{6\sqcup 1\scheme{10}}} $ & \realized{\scheme{8\sqcup 2\scheme{3} \sqcup 1\scheme{5}} } & $\scheme{4\sqcup 1\scheme{1} \sqcup 1\scheme{3} \sqcup 1\scheme{11}} $ & $\scheme{1\scheme{1} \sqcup 1\scheme{7} \sqcup 1\scheme{11}} $ \\
$\scheme{8\sqcup 1\scheme{2\sqcup 1\scheme{10}}} $ & $\scheme{5\sqcup 1\scheme{6\sqcup 1\scheme{9}}} $ & $\scheme{1\sqcup 1\scheme{10\sqcup 1\scheme{9}}} $ & $\scheme{4\sqcup 1\scheme{1} \sqcup 1\scheme{5} \sqcup 1\scheme{9}} $ & $\scheme{1\scheme{5} \sqcup 2\scheme{7}} $ \\
$\scheme{9\sqcup 1\scheme{2\sqcup 1\scheme{9}}} $ & \realized{\scheme{6\sqcup 1\scheme{6\sqcup 1\scheme{8}}} } & $\scheme{2\sqcup 1\scheme{10\sqcup 1\scheme{8}}} $ & $\scheme{4\sqcup 1\scheme{1} \sqcup 2\scheme{7}} $ & $\scheme{1\sqcup 1\scheme{18\sqcup 1\scheme{1}}} $ \\
\realized{\scheme{10\sqcup 1\scheme{2\sqcup 1\scheme{8}}} } & $\scheme{7\sqcup 1\scheme{6\sqcup 1\scheme{7}}} $ & $\scheme{3\sqcup 1\scheme{10\sqcup 1\scheme{7}}} $ & $\scheme{4\sqcup 1\scheme{3} \sqcup 1\scheme{5} \sqcup 1\scheme{7}} $ &  \\
$\scheme{11\sqcup 1\scheme{2\sqcup 1\scheme{7}}} $ & $\scheme{8\sqcup 1\scheme{6\sqcup 1\scheme{6}}} $ & $\scheme{4\sqcup 1\scheme{10\sqcup 1\scheme{6}}} $ & $\scheme{4\sqcup 3\scheme{5}} $ &  \\
$\scheme{12\sqcup 1\scheme{2\sqcup 1\scheme{6}}} $ & \realized{\scheme{9\sqcup 1\scheme{6\sqcup 1\scheme{5}}} } & \realized{\scheme{5\sqcup 1\scheme{10\sqcup 1\scheme{5}}} } & $\scheme{1\sqcup 1\scheme{14\sqcup 1\scheme{5}}} $ &  \\
$\scheme{13\sqcup 1\scheme{2\sqcup 1\scheme{5}}} $ & \realized{\scheme{10\sqcup 1\scheme{6\sqcup 1\scheme{4}}} } & \realized{\scheme{6\sqcup 1\scheme{10\sqcup 1\scheme{4}}} } & $\scheme{2\sqcup 1\scheme{14\sqcup 1\scheme{4}}} $ &  \\
\openmark{\scheme{14\sqcup 1\scheme{2\sqcup 1\scheme{4}}} } & \realized{\scheme{11\sqcup 1\scheme{6\sqcup 1\scheme{3}}} } & \realized{\scheme{7\sqcup 1\scheme{10\sqcup 1\scheme{3}}} } & $\scheme{3\sqcup 1\scheme{14\sqcup 1\scheme{3}}} $ &  \\
$\scheme{15\sqcup 1\scheme{2\sqcup 1\scheme{3}}} $ & $\scheme{12\sqcup 1\scheme{6\sqcup 1\scheme{2}}} $ & $\scheme{8\sqcup 1\scheme{10\sqcup 1\scheme{2}}} $ & $\scheme{4\sqcup 1\scheme{14\sqcup 1\scheme{2}}} $ &  \\
$\scheme{16\sqcup 1\scheme{2\sqcup 1\scheme{2}}} $ & \realized{\scheme{13\sqcup 1\scheme{6\sqcup 1\scheme{1}}} } & \realized{\scheme{9\sqcup 1\scheme{10\sqcup 1\scheme{1}}} } & \realized{\scheme{5\sqcup 1\scheme{14\sqcup 1\scheme{1}}} } &  \\
\realized{\scheme{17\sqcup 1\scheme{2\sqcup 1\scheme{1}}} } &  &  &  &  \\
\bottomrule
\end{tabular}%
}
\end{table}

\section{Computing and realizing real schemes}
\label{sec:realizations}
\label{sec:computation}

The computations in this section rest on an algorithm, developed in the companion paper \cite{2026:GeiselmannEtAl:IsotopyComputation}, that computes the real scheme of a patchworked curve $\cC(\cT,\sigma)$ directly from the combinatorial data in near-quadratic time in the degree.
Using the regions and their nesting from Section~\ref{subsec:regions-nesting}, it groups the vertices of $\cT^\diamond$ into same-sign regions by a union-find pass \cite{1975:Tarjan:LinearAlgorithm}, identifies antipodal boundary points to pass from $\RR^2$ to $\RP^2$, and returns the real scheme as the canonical form of the resulting rooted tree of regions.
Each pair $(\cT,\sigma)$ is classified independently and with only fixed-size, thread-local state, so the computation is \enquote{embarrassingly parallel}; that is, the parallelization does not require any communication between workers.
This computation maps naturally onto a GPU: for example, a single NVIDIA~A100 classifies on the order of $10^8$ pairs per second in degree eight.
We use this classification in two complementary ways.

The first is a general search over pairs $(\cT,\sigma)$ that makes no assumption about the resulting curve.
Its basic operation is an exhaustive sweep of all sign distributions for a fixed triangulation; by the equivalence of Section~\ref{sec:signs} it suffices to take one representative per class.
Such a sweep takes seconds in degree six, minutes in degree seven, and hours in degree eight, and as a byproduct yields the distribution of real schemes over all sign distributions of that triangulation (Appendix~\ref{app:statistics}).
The number of triangulations grows quickly with the degree (Table~\ref{tab:search-space}), so the search proceeds differently in each case.
In degree six, it is simply exhaustive: a sweep over all symmetric triangulations already realizes every nonempty scheme.
In degree seven, we sweep a small, hand-curated set of triangulations, chosen for the variety of schemes they produce.
In degree eight, there are far too many triangulations to sweep, so we first generate promising ones by simulated annealing and tabu search, and then iterate in \enquote{ping-pong} fashion: sweep all symmetric triangulations against a pool of witness sign distributions, sweep all sign distributions over the most productive triangulations, and enlarge both pools by edge flips and weight-space interpolation.
Together, these procedures produce every nonempty real scheme through degree seven and the bulk of those in degree eight.
Restricting the search to symmetric triangulations is motivated by Brugall\'e's results on algebraic curves of degree seven \cite[Proposition 5.36]{Brugalle:Diss}.
\begin{table}[t]
\centering
\caption{Number of regular unimodular triangulations of $d\cdot\Delta_2$ up to the action of $\Sym{3}$, with the number of these that are symmetric under the reflection $(x,y)\mapsto(y,x)$ in parentheses; $M$ is Harnack's bound, and $|A|$ is the number of lattice points. A dash marks counts too large to enumerate. The final column lists the number of sign distributions per triangulation up to equivalence, $2^{|A|-3}$ (Section~\ref{sec:signs}).}
\label{tab:search-space}
\smallskip
\begin{tabular*}{\linewidth}{@{\extracolsep{\fill}}rrrrr@{}}
\toprule
$d$ & $|A|$ & $M$ & triangulations (symmetric) & sign distributions \\
\midrule
2 & 6  & 1  & $2$ \,($2$) & $8$ \\
3 & 10 & 2  & $18$ \,($7$) & $128$ \\
4 & 15 & 4  & $1{,}278$ \,($74$) & $4{,}096$ \\
5 & 21 & 7  & $561{,}885$ \,($1{,}194$) & $262{,}144$ \\
6 & 28 & 11 & $1{,}198{,}202{,}590$ \,($62{,}960$) & $33{,}554{,}432$ \\
7 & 36 & 16 & {--} \,($4{,}729{,}700$) & $8{,}589{,}934{,}592$ \\
8 & 45 & 22 & {--} \,($1{,}199{,}795{,}773$) & $4{,}398{,}046{,}511{,}104$ \\
\bottomrule
\end{tabular*}
\end{table}

The second method applies only to $M$-curves and builds on Haas' zone decompositions (Section~\ref{sec:haas}).
By Theorem~\ref{thm:haas}, maximal T-curve corresponds to a compatible collection of simple and double Harnack splits, so the $M$-curves can be enumerated by growing such a collection one split at a time. This process forms a directed acyclic graph whose root, the empty collection, represents the Harnack curve $\cC(\cT,\eta)$.
The full graph is intractable even in moderate degree, so in degree eight we unfold it only partially: to depth four from the root and around the split collections of $M$-curves already found by the general search.
Together, these limited unfolds realize eight $M$-curves that the general search misses.

\subsection{Curves of degree at most six}
\label{sec:deg-6}
The classification of real schemes of degree $d\leq 5$ was already accomplished by Harnack and Hilbert.
We summarize their results in terms of T-curves; see also \cite[Corollary 1.3]{1998:DeLoeraWicklin:ConvexityPatchworking}.

\begin{proposition}\label{thm:covering-5}
  All nonempty real schemes of curves of degree $d\leq 5$ can be generated as T-curves from the honeycomb triangulation of degree $d$.
\end{proposition}
\begin{proof}
  For each degree $d\leq 5$, there is a unique real scheme of an $M$-curve.
  The Harnack sign distribution $\eta:A\to\GF{2}$ and the honeycomb triangulation $\cH_d$ yield the special Harnack curve $\cC(\cH_d,\eta)$.
  By Proposition~\ref{prop:special-harnack}, this is an $M$-curve.

  The cases $d\leq 3$ are straightforward and thus omitted.
  For $d=4$, the real scheme of $\cC(\cH_4,\eta)$ is $\scheme{4}$.
  By flipping the signs of the points $(1,3)$ and $(3,1)$ in $\eta$, we obtain T-curves with real schemes $\scheme{3}$ and $\scheme{2}$.
  The constant sign distribution $\1$ yields the T-curve $\cC(\cH_4,\1)$ with real scheme $\scheme{1\scheme{1}}$.
  Flipping the sign of $(3,1)$ in $\1$ gives the final nonempty real scheme $\scheme{1}$.

  For $d=5$, the real scheme of $\cC(\cH_5,\eta)$ is $\scheme{J\sqcup 6}$.
  By flipping the signs of the points $(1,4)$, $(2,3)$, and $(3,2)$ in $\eta$, we obtain T-curves with real schemes $\scheme{J\sqcup 5}$, $\scheme{J\sqcup 4}$, and $\scheme{J\sqcup 3}$.
  The constant sign distribution $\1$ yields the T-curve $\cC(\cH_5,\1)$ with real scheme $\scheme{J\sqcup 1\scheme{1}}$.
  Flipping $(1,4)$ in $\1$ produces $\scheme{J\sqcup 1}$; additionally flipping $(2,3)$ produces $\scheme{J}$.
  Flipping $(1,4)$ and $(4,1)$ in $\1$ gives the final nonempty real scheme $\scheme{J\sqcup 2}$.
\end{proof}

\begingroup
\newcommand{\scalefactor}{0.8}
\input{tikz/deg6_triangulations.tex}
\endgroup

We now revisit the known classification of the real schemes in degree six, which is the first nontrivial case.
The 56 real schemes have been classified by Gudkov \cite{1969:Gudkov:ArrangementOvals, 1974:Gudkov:ProjectiveAlgebraic}.
The classification of the rigid isotopy types was obtained by Rokhlin \cite{1978:Rokhlin:TopologicalAlgebraic} and Nikulin \cite{1980:Nikulin:SymmetricBilinear}.
We copy the result from Viro \cite{1986:Viro:TopologyAlgebraic}, who uses the notation that we introduced earlier.
\begin{theorem}\label{thm:classification-6}
  There exist curves of degree six with the following real schemes:
  \begin{enumerate}
  \item $\langle \alpha \sqcup 1\langle\beta\rangle \rangle$ where $\alpha + \beta \leq 10$, $0 \leq \alpha \leq 9$, $1 \leq \beta \leq 9$, and
  \begin{align*}\alpha-\beta\equiv\begin{cases} 0\bmod{8}\,, &\text{ if } \alpha+\beta=10\,,\\
  		 \pm 1\bmod{8}\,, &\text{ if } \alpha+\beta=9\,;
  \end{cases}\end{align*}
  \item $\langle \alpha \rangle$ with $0 \leq \alpha \leq 10\,$;
  \item $\langle 1\langle 1\langle 1\rangle\rangle\rangle\,$.
  \end{enumerate}
  Any curve of degree six has one of these 56 real schemes.
\end{theorem}

The 56 real schemes include the empty scheme $\scheme{0}$ and three schemes of $M$-curves; here $M=11$.
The types $\langle 1 \sqcup 1\langle 9\rangle \rangle$ with $(p,n)=(2,9)$ and $\langle 9 \sqcup 1 \langle 1 \rangle\rangle$ with $(p,n)=(10,1)$ occur in Hilbert's article \cite[pp.~118--119]{1891:Hilbert:ReellenZuge}.
Gudkov found the third $M$-curve, of type $\langle 5 \sqcup 1\langle 5\rangle \rangle$, with $(p,n)=(6,5)$~\cite{1974:Gudkov:ProjectiveAlgebraic}.

While it is known that all nonempty real schemes of degree six arise as T-curves, we are not aware of an explicit list.
Viro sketches the construction of 53 out of the 56 types in~\cite[\S~4]{2008:Viro:SixteenthHilbert}; see also \cite[\S3.2]{1984:Viro:GluingAlgebraic}.
Note that a patchwork can never realize the empty type.
Our version of that classification shows that only two triangulations suffice to construct all relevant patchworks.
We say that a triangulation $\cT$ \emph{supports} a real scheme $\scheme{X}$ if there is a sign distribution $\sigma$ such that the real scheme of $\cC(\cT,\sigma)$ is $\scheme{X}$.

\begin{theorem}\label{thm:covering-6}
  All 55 nonempty real schemes of degree six can be generated as T-curves from the two triangulations of $6\cdot\Delta_2$ shown in Figure~\ref{fig:deg6-triangulations} by varying the sign distributions.
\end{theorem}

\begin{proof}
  First, we verify that the two triangulations are regular.
  To this end, we provide explicit lifting functions in Figure~\ref{fig:deg6-triangulations}.
  Checking whether a given lifting function $\omega:A\to\RR$ fits a given triangulation $\cT$ of $A$ requires verifying the strict folding conditions from \cite[Proposition 5.2.6]{2010:DeLoeraRambauSantos:TriangulationsAlgorithms}.
  Although finding a lifting function requires solving a linear program, the folding conditions for a given $\omega:A\to\RR$ are easy to check by hand, as follows.
  Each interior edge $e$ of $\cT$ with vertices $v$ and $x$ lies in exactly two triangles, whose third vertices we denote by $u$ and $w$.
  Since both triangles are unimodular, we have $u+w=v+x+k(x-v)$ for some $k\in\ZZ$.
  The edge $e$ belongs to the regular triangulation of $A$ induced by $\omega$ if and only if $\omega(u)+\omega(w)>(1-k)\,\omega(v)+(1+k)\,\omega(x)$.
  Checking all interior edges of $\cT$ suffices.
  For instance, the edge $(1,0) \sim (0,1)$ occurs in the triangulation {\degsixA} because here $k=0$ and $4+1>2+2$; see Figure~\ref{fig:bat}.
  For the edge $(0,1) \sim (2,2)$ of {\degsixA} we have $k=1$, and the condition reads $\omega(1,1)+\omega(3,3)>2\,\omega(2,2)$, that is, $1+0>0$.
  
  For each nonempty real scheme, Table~\ref{tab:deg-6} lists a sign distribution and one of the two triangulations in Figure~\ref{fig:deg6-triangulations}.
  The sign vectors are in lexicographic order; see Remark~\ref{rem:lex-ordering}.
  The T-curve $\cC(\cT,\sigma)$ can also be constructed by hand from $\cT$ and $\sigma:A\to\GF{2}$.
  Since this may be tedious, we provide each patchwork as a \polymake file, which can be loaded into our patchworking tool to obtain a visualization.
  That tool also recognizes the real scheme.
  Our algorithm is a fast version of Viro's sketch \cite[Algorithm 1.4.E]{2006:Viro:PatchworkingAlgebraic}; it is summarized at the beginning of this section and treated in detail in the companion paper \cite{2026:GeiselmannEtAl:IsotopyComputation}.
  For information on further software, see Appendix~\ref{app:software}.
\end{proof}

\begin{remark}\label{rem:coefficients}
  Obtaining explicit polynomials defining those 55 curves requires combining Theorems~\ref{thm:patchworking} and~\ref{thm:covering-6}.
  This entails finding a valid threshold $t_0$ for each patchwork $(\cT,\sigma)$ via bisection.
  Our approach should be compared to Kaihnsa et al.\ \cite{2019:KaihnsaEtAl:SixtyCurves}, who gave explicit polynomials for all real schemes of degree six.
  Most likely, the polynomials from \cite{2019:KaihnsaEtAl:SixtyCurves} have smaller coefficients than those obtained from Theorems~\ref{thm:patchworking} and~\ref{thm:covering-6}.
  However, it is highly nontrivial and impossible to do by hand to derive the real scheme from a given polynomial.
  Kaihnsa and co-authors employ quantifier elimination in \Mathematica.
  One advantage of our approach is that a single patchwork of moderate degree can be verified directly.
\end{remark}

\input{tables/proof_deg6.tex}

\clearpage

\begingroup
\newcommand{\scalefactor}{0.75}
\input{tikz/deg7_triangulations_A.tex}
\endgroup

\subsection{Curves of degree seven}
\label{sec:deg-7}

Viro's classification \cite{1986:Viro:TopologyAlgebraic} of smooth plane real curves of degree seven is an analog of Theorem~\ref{thm:classification-6}.
In this case, $M=16$.
\begin{theorem}[{Viro \cite[Theorem~3.3.A]{1984:Viro:GluingAlgebraic}}]\label{thm:classification-7}
  There exist curves of degree seven with the following real schemes:
  \begin{enumerate}
  \item $\langle J \sqcup \alpha \sqcup 1\langle\beta\rangle \rangle$ with $\alpha + \beta \leq 14$, $0 \leq \alpha \leq 13$, $1 \leq \beta \leq 13\,$;
  \item $\langle J \sqcup \alpha \rangle$ with $0 \leq \alpha \leq 15\,$;
  \item $\langle J \sqcup 1\langle 1\langle 1\rangle\rangle\rangle\,$.
  \end{enumerate}
  Any curve of degree seven has one of these 121 real schemes.
\end{theorem}

Hilbert constructs four $M$-curves of degree seven in \cite{1891:Hilbert:ReellenZuge}:
$\langle J \sqcup 12 \sqcup 1 \langle 2\rangle \rangle$,
$\langle J \sqcup 2 \sqcup 1 \langle 12\rangle \rangle$,
$\langle J \sqcup 11 \sqcup 1 \langle 3\rangle \rangle$, and
$\langle J \sqcup 1 \sqcup 1 \langle 13\rangle \rangle$; see also \cite[p.~58]{1974:Gudkov:ProjectiveAlgebraic}.
Viro sketches patchworks for all 121 types \cite[\S3.3]{1984:Viro:GluingAlgebraic}, but not as T-curves.
De Loera and Wicklin \cite{1998:DeLoeraWicklin:ConvexityPatchworking} list 24 T-curves of degree seven.
Although we need four triangulations here, rather than two as in degree six, that number is still very low.

\begingroup
\newcommand{\scalefactor}{0.75}
\input{tikz/deg7_triangulations_B.tex}
\endgroup

\begin{theorem}\label{thm:covering-7}
  All of the 121 real schemes of curves of degree seven can be generated as T-curves from the four triangulations of $7\cdot\Delta_2$ shown in Figures~\ref{fig:deg7-A-triangulations} and \ref{fig:deg7-B-triangulations} by varying the sign distributions.
\end{theorem}

\begin{proof}
  The proof is essentially the same as the one for Theorem~\ref{thm:covering-6}.
  For degree seven, the pairs of triangulations and sign distributions are given in Table~\ref{tab:deg-7}.
\end{proof}

The triangulation {\degsevenD} in Figure~\ref{fig:deg7-D} is the honeycomb triangulation $\cH_7$.
The lifting function stems from \cite[Remark 15]{2024:CasabellaEtAl:WronskiHoneycomb}; cf.\ Example~\ref{exmp:honeycomb}.
The name {\degsevenA} in Figure~\ref{fig:deg7-A} refers to the fact that the long edges of this triangulation appear to radiate centrifugally from one triangle. The remaining two triangulations appear to be a \enquote{split} and a \enquote{frayed} version of it; hence, we call them {\degsevenB} (Figure~\ref{fig:deg7-B}) and {\degsevenC} (Figure~\ref{fig:deg7-C}), respectively.

\input{tables/proof_deg7.tex}

\clearpage

\subsection{Curves of degree eight}
\label{sec:deg-8}
Unlike in lower degrees, the real schemes of degree eight are not classified, and outside the maximal case they appear not to have been studied systematically.
Making no assumption on the number of ovals, we realize $2{,}367$ distinct nonempty real schemes as T-curves through a combination of the approaches outlined above, with representatives for every number of ovals from $1$ to $22$; as the search is not exhaustive, this is a lower bound.

By contrast, the maximal case with $M=22$ ovals is well studied. Orevkov \cite{2002:Orevkov:CurvesIsotopy} determined that the 89 real schemes in Table~\ref{tab:deg-8} are realizable by flexible (pseudo-holomorphic) $M$-curves of degree eight.
Whether these are also realizable by genuine algebraic curves is a finer question: $83$ are, while for the remaining six it is unknown.
We show that at least $38$ can be realized as T-curves, each certified by an explicit patchwork: $5$ of the $20$ with $(p,n)=(19,3)$, $11$ of $19$ with $(15,7)$, $14$ of $19$ with $(11,11)$, and $8$ of $19$ with $(7,15)$.
In the opposite direction, Theorem~\ref{thm:3_19_impossible} shows that none of the $12$ schemes with $(p,n)=(3,19)$ can be realized as a T-curve, even though eight of them are known to be realizable by algebraic curves.
This leaves 39 $M$-curve schemes, among them the two algebraically open ones with $(p,n)=(19,3)$, for which T-curve realizability is open.

\section{Concluding remarks}
\label{sec:outlook}
It would be interesting to know which real schemes exactly can be realized as T-curves.
For degree at most seven, we have the complete picture: every nonempty real scheme is realizable by a T-curve.
Degree eight is the first degree for which this fails: the maximal schemes with $(p,n)=(3,19)$ are not T-curves (Theorem~\ref{thm:3_19_impossible}). In this degree, the picture is only partial: we realize $38$ of the $89$ maximal schemes, which leaves the following question open. See also Itenberg and Viro \cite[p.22]{1996:ItenbergViro:PatchworkingRagsdale}.
\begin{question}
  Which real schemes of degree eight are realizable as T-curves? For $39$ of the $89$ maximal schemes this is currently undecided.
\end{question}
The following construction was communicated to us by Ilia Itenberg.
By combining an argument from \cite[\S7]{1995:Itenberg:RagsdaleCurves} with a method described by Wiman \cite[p.~227]{1923:Wiman:ReellenZuge}, one can construct a real scheme of degree 36 that cannot be realized as a T-curve of the same degree.
Independently, we would like to know how T-curves are related to regular triangulations.
\begin{question}
  How many regular unimodular triangulations of $d\cdot\Delta_2$ are necessary to support every real scheme realizable by a T-curve of degree $d$?
\end{question}
By Proposition~\ref{thm:covering-5} and Theorems~\ref{thm:covering-6} and~\ref{thm:covering-7}, the number of such triangulations equals one for $d\leq 5$, at most two for $d=6$, and at most four for $d=7$.
We also see that, for all degrees $d\leq 7$, unimodular triangulations that are symmetric about the line $x=y$ suffice; this corresponds to Brugall\'e's results \cite[\S5]{Brugalle:Diss}.

\appendix

\section{Statistics of sign distributions}
\label{app:statistics}

In Tables~\ref{tab:deg-6} and~\ref{tab:deg-7} we displayed one patchwork per nonempty real scheme in degrees six (from two triangulations of $6\cdot\Delta_2$) and seven (from four triangulations of $7\cdot\Delta_2$).
This section reports a complete analysis of those six triangulations by considering all equivalence classes of sign distributions determined in Section~\ref{sec:signs}.
For each nonempty real scheme and each triangulation, we list the number of equivalence classes that produce that scheme.

These data are the condensed output of a large computational experiment.
In contrast to the results in the main body of this article, the reported frequencies are intrinsically computational and impossible to verify by hand.
Since even a computational verification is challenging, we refrain from formally stating the various outcomes as mathematical theorems.

\begin{figure}[bt]
  \centering
  \begin{subfigure}[b]{0.48\textwidth}
    \centering
    \resizebox{\textwidth}{!}{\input{tikz/gudkov_53_3_density.tex}}
    \caption{\degsixA}
  \end{subfigure}
  \hfill
  \begin{subfigure}[b]{0.48\textwidth}
    \centering
    \resizebox{\textwidth}{!}{\input{tikz/symmetric-14822_density.tex}}
    \caption{\degsixB}
  \end{subfigure}
  \caption{Distributions of the number of ovals for two triangulations in degree six.
    The maximum eleven marks the $M$-curves.}
  \label{fig:densities-deg-6}
\end{figure}

\subsection{Degree six}

For degree six, we have two triangulations, which look quite similar at first sight; cf.\ Figure~\ref{fig:deg6-triangulations}.
Figure~\ref{fig:densities-deg-6} shows the distributions of the number of loops.
The averages 4.41 (for {\degsixA}) and 3.70 (for {\degsixB}) provide computational data relevant to a question raised by Gayet and Welschinger about the expected number of loops of an arbitrary real plane curve \cite[p.94]{2011:GayetWelschinger:ExponentialRarefaction}.

The triangulation {\degsixA} can realize all the nonempty types except two: $\scheme{10}$ and $\scheme{1\scheme{1\scheme{1}}}$.
The triangulation {\degsixB} also supports this scheme.
The $(M{-}1)$-curve $\scheme{10}$ is supported by neither $\cH_6$ nor {\degsixA}.
This real scheme is also realized by a \enquote{Robinson curve} in \cite[\S2]{2019:KaihnsaEtAl:SixtyCurves}.
The triangulation $\cH_6$ realizes 37 of the 55 nonempty types, and an exhaustive check confirms that no single symmetric triangulation of $6\cdot\Delta_2$ combined with $\cH_6$ covers all 55 nonempty types.
We conjecture that there is no regular unimodular triangulation of $6\cdot\Delta_2$ that supports all 55 nonempty real schemes, but we do not have a proof.

\input{tables/degree6_cover.tex}

\clearpage

\begin{figure}[th]
  \centering
  \begin{subfigure}[b]{0.48\textwidth}
    \centering
    \resizebox{\textwidth}{!}{\input{tikz/Jv1v1_13__centrifugal_density.tex}}
    \caption{\degsevenA}
  \end{subfigure}
  \hfill
  \begin{subfigure}[b]{0.48\textwidth}
    \centering
    \resizebox{\textwidth}{!}{\input{tikz/Jv2v1_12__split_centrifugal_density.tex}}
    \caption{\degsevenB}
  \end{subfigure}

  \vspace{1em}

  \begin{subfigure}[b]{0.48\textwidth}
    \centering
     \resizebox{\textwidth}{!}{\input{tikz/Jv3v1_11__frayed_centrifugal_density.tex}}
    \caption{\degsevenC}
  \end{subfigure}
  \hfill
  \begin{subfigure}[b]{0.48\textwidth}
    \centering
    \resizebox{\textwidth}{!}{\input{tikz/honeycomb-7_density.tex}}
    \caption{\degsevenD}
  \end{subfigure}
  \caption{Distributions of the number of ovals for four triangulations of degree seven; these counts ignore the pseudo-line $J$. Thus, the maximum 15 marks the $M$-curves.}
  \label{fig:densities-deg-7}
\end{figure}

\subsection{Degree seven}

The statistics for $d=7$ are more involved than for $d=6$ in several ways.
First, there are four triangulations to consider instead of two.
Second, the triangulation that covers most of the types, {\degsevenA}, misses six real schemes, whereas {\degsixA} misses only two nonempty real schemes in degree six.
These six real schemes, $\scheme{ J \sqcup 10 \sqcup 1 \scheme{ 4 }}$, $\scheme{ J \sqcup 2 \sqcup 1 \scheme{ 12 }}$, $\scheme{ J \sqcup 3 \sqcup 1 \scheme{ 11 }}$, $\scheme{ J \sqcup 6 \sqcup 1 \scheme{ 8 }}$, $\scheme{ J \sqcup 1 \scheme{ 13 }}$, and $\scheme{ J \sqcup 1 \scheme{ 11 }}$, seem unlikely to be supported by single triangulation.
Third, while all the distributions in Figure~\ref{fig:densities-deg-7} look similar, their mean values and standard deviations vary visibly.
We do not know whether three or fewer triangulations of $7\cdot\Delta_2$ suffice to support all 121 real schemes.

\input{tables/degree7_cover.tex}

\section{Further software and data}
\label{app:software}

Known software systems for patchworking include the Combinatorial Patchworking Tool of El-Hilany, Rau, and Renaudineau \cite{2017:ElHilanyEtAl:CombinatorialPatchworking}, \texttt{Viro.Sage} by de Wolff, Kwaakwah, and O'Neill \cite{2021:DeWolffEtAl:ViroSage}, and \polymake \cite{2000:GawrilowJoswig:PolymakePolytopes,2020:JoswigVater:HyperfacesPatchworking}.

For this article, we developed an improved version of \cite{2017:ElHilanyEtAl:CombinatorialPatchworking}.
It allows one to interactively construct, modify, and analyze T-curves.
The tool can also load and save \polymake files, e.g., for further analysis along the lines of \cite{2020:JoswigVater:HyperfacesPatchworking}.
We provide one dataset for each nonempty real scheme of degree at most seven. These data form an electronic version of Tables~\ref{tab:deg-6} and~\ref{tab:deg-7}, extended to all $d\in\{2,3,\dots,7\}$.\footnote{\url{https://github.com/dmg-lab/CombinatorialPatchworking}}
For $d=8$, the repository additionally contains a compressed archive with one patchwork certificate for each of the $2{,}367$ nonempty real schemes from Section~\ref{sec:deg-8}.\footnote{\url{https://github.com/dmg-lab/CombinatorialPatchworking/blob/main/deg8.pcoms.txz}}

\section{Known constructions}
\label{app:known}
\newcommand{\newton}{\Delta}

The realizability of degree six and seven curves by T-curves seems to be known.
Since we found the relevant literature somewhat hard to navigate, for the benefit of the reader, here we sketch a construction, privately communicated to us by Erwan Brugall\'e; see also \cite{Brugalle:Diss}.

Viro's general patchworking theorem \cite{1984:Viro:GluingAlgebraic} stitches together bivariate polynomials.
We give a rough overview, omitting details and ignoring generalizations to more variables; cf.\ \cite[\S11.5.C]{1994:GKZ:DiscriminantsResultants}.
One starts with a list of polynomials $a_1,\ldots, a_n$ in two variables $x,y$, and their Newton polygons $\newton(a_i)$ in $\RR^2$, such that:
\begin{enumerate}
\item the polygons $\newton(a_1),\ldots,\newton(a_n)$ are the maximal cells of a regular subdivision of a convex polygon $\Delta$, induced by a height function $\nu:\Delta\to \RR$;
\item the polynomials $a_1, \dots, a_n$  are totally non-degenerate, in particular, the curves $V_{\CC}(a_i)$ are smooth;
\item the truncations $a_i^{\newton(a_i)\cap\newton(a_j)}$ and $a_j^{\newton(a_i)\cap\newton(a_j)}$ to the common intersection $\newton(a_i)\cap\newton(a_j)$ agree for all $i\neq j$; here the \emph{truncation} $a_i^{\Gamma}$ comprises the terms of $a_i$ whose exponent vectors lie in the set $\Gamma$;
\end{enumerate}
The second condition says that the coefficients of the pieces patch together: there is a polynomial $a=\sum_w a_w x^{w_1}y^{w_2}$ with $\Delta=\newton(a)$ and $a^{\newton(a_i)}=a_i$.
We obtain the family of polynomials
\begin{equation}\label{eq:patchworked}
  b_t(x,y)\ :=\ \sum_{w\in\Delta\cap\ZZ^2}a_w\cdot t^{\nu(w)}\cdot x^{w_1} y^{w_2}
\end{equation}
parameterized by $t$.
Each of the polynomials above comes with a \enquote{chart}.
Fully explaining that notion and further technical details is beyond our current scope.
\begin{theorem}[{Viro \cite{1984:Viro:GluingAlgebraic}}]\label{thm:generalpatchworking}
  For $t>0$ small enough the chart of $b_t$ is the gluing of the charts of the $a_i$, such that the affine curve $V_\RR(b_t)$ has controlled topology.
\end{theorem}

To make up for the lacking details, we reproduce a simple example.
\begin{example}[{\cite[Example~1.5.A]{2006:Viro:PatchworkingAlgebraic}}]
  For the two polynomials
  \[
    f(x,y)\ :=\ 8x^3 - x^2 + 4y^2\ \mbox{ and }\ g(x,y)\ :=\ 4y^2 - x^2 + 1
  \]
  the function $\nu$ attaining $0$ on $(0,0),(2,0),(0,2)$ and $1$ on $(3,0)$ is a height function on the four vertices of the subdivision of the triangle $\Delta=\conv\{(0,0),(3,0),(0,2)\}$ with the two maximal cells $\Delta(f)$ and $\Delta(g)$, the shaded regions in Figure~\ref{fig:viro_pw}.
  The patchworked family from \eqref{eq:patchworked} and Theorem~\ref{thm:generalpatchworking} is now
  \[
    b_t(x,y) \ = \ 8tx^3 - x^2 + 4y^2 + 1 \enspace .
  \]

  \begin{figure}[th]
    \begin{tabular}{ccc}
      Chart of $f$ & Chart of $g$ & Glued patchwork\\
\begin{tikzpicture}[scale=.7]
  \fill[newton] (3,0) -- (2,0) -- (0,2) -- cycle;
  \draw[grid] (-3.2, -2.2) grid (3.2, 2.2);
  \draw[axis] (-3.2,0) -- (3.5,0);
  \draw[axis] (0,-2.2) -- (0,2.4);
  \draw[polygon] (0,2) -- (3,0) -- (0,-2) -- (-3,0) -- cycle;
  \draw[polygon] (0,2) -- (2,0) -- (0,-2) -- (-2,0) -- cycle;
  \draw[curve] (1,1) .. controls (1.55,0.72) and (2.45,0.32) .. (2.45,0)
               .. controls (2.45,-0.32) and (1.55,-0.72) .. (1,-1);
  \draw[curve] (-1,1)  to[bend left=12] (-1.2,37/30);
  \draw[curve] (-1,-1) to[bend right=12]  (-1.2,-37/30);
\end{tikzpicture}
       &
\begin{tikzpicture}[scale=.7]
  \fill[newton] (0,0) -- (2,0) -- (0,2) -- cycle;
  \draw[grid] (-2.2, -2.2) grid (2.2, 2.2);
  \draw[axis] (-2.1,0) -- (2.4,0);
  \draw[axis] (0,-2.1) -- (0,2.4);
  \draw[polygon] (0,2) -- (2,0) -- (0,-2) -- (-2,0) -- cycle;
  \draw[curve] (1,1)  arc (126.87:233.13:1.25);
  \draw[curve] (-1,1) arc (53.13:-53.13:1.25);
\end{tikzpicture}
       &
\begin{tikzpicture}[scale=.7]
  \fill[newton] (0,0) -- (3,0) -- (0,2) -- cycle;
  \draw[grid] (-3.2, -2.2) grid (3.2, 2.2);
  \draw[axis] (-3.1,0) -- (3.5,0);
  \draw[axis] (0,-2.1) -- (0,2.4);
  \draw[polygon] (0,2) -- (3,0) -- (0,-2) -- (-3,0) -- cycle;
  \draw[polygon] (0,2) -- (2,0) -- (0,-2) -- (-2,0) -- cycle;
  \draw[curve] (1,1)  arc (126.87:233.13:1.25);
  \draw[curve] (-1,1) arc (53.13:-53.13:1.25);
  \draw[curve] (1,1) .. controls (1.55,0.72) and (2.45,0.32) .. (2.45,0)
               .. controls (2.45,-0.32) and (1.55,-0.72) .. (1,-1);
  \draw[curve] (-1,1)  to[bend left=12] (-1.2,37/30);
  \draw[curve] (-1,-1) to[bend right=12]  (-1.2,-37/30);
\end{tikzpicture}
    \end{tabular}
    
    Reproduction of Figure~12, 14 and 15 from \cite{2006:Viro:PatchworkingAlgebraic}.
    \caption{Illustration for general patchwork}
    \label{fig:viro_pw}
  \end{figure}
\end{example}

The charts can be thought of as diagrams of the vanishing loci $V_\RR(f)$, $V_\RR(g)$ and $V_\RR(b_t)$ mapped to the respective Newton polygons; this idea is captured in Figure~\ref{fig:viro_pw}.
Next we describe how Theorems~\ref{thm:patchworking} and \ref{thm:generalpatchworking} are related.
For a combinatorial patchwork, we are given a triangulation $\cT$ and a sign function $\sigma$.
Now for the $i$th triangle with vertices $u,v,w$ the polynomial
\[
  a_i(x,y)\ :=\ (-1)^{\sigma(u)}\cdot x^{u_1} y^{u_2}
  \ +\  (-1)^{\sigma(v)}\cdot x^{v_1} y^{v_2}
  \ +\  (-1)^{\sigma(w)}\cdot x^{w_1} y^{w_2}
\]
will be inserted in the general patchworking theorem.
The polytope $\Delta$ then becomes $d\cdot\Delta_2$ and the function $\nu$ induces the regular triangulation $\cT$.
In this way the combinatorial patchworks arise as a special case of Theorem~\ref{thm:generalpatchworking}.

Now T-curves can be obtained through Theorem~\ref{thm:generalpatchworking} by starting with some subdivision $\Delta:= \bigcup_i\newton(a_i)$ resulting in a singular curve and refining that subdivision into a unimodular triangulation.
We sketch the idea for degree $d=6$.
Let
\[
  s(x,y) := (y-\alpha x^2)(y-\beta x^2)(y-\gamma x^2)\, , \quad \text{for } \alpha,\beta,\gamma\in\RR \text{ pairwise distinct.}
\]
The affine curve $V_\RR(s)\subset\RR^2$ is the union of three real parabolas sharing an axis: pairwise tangent at the origin (common tangent line $y=0$) and pairwise tangent at the point at infinity $[0:1:0]$.
The curve $V_\RR(s)$ has exactly two singular points of type $J_{10}^-$ in Arnold's classification \cite[\S15.1]{Arnold:singularities}.
The Newton polygon $\Delta(s)$ is the line segment $[(0,3),(6,0)]$, which splits $6\cdot\Delta_2$ into two triangles.
In \cite{1984:Viro:GluingAlgebraic,2008:Viro:SixteenthHilbert} Viro shows that the two points can be smoothed independently and simultaneously, realizing $54$ of the $55$ nonempty real schemes of degree six; the missing type $\scheme{10}$ is constructed separately.

In degree $d=7$, Korchagin \cite{1988:Korchagin:IsotopySingular} classifies real schemes of smoothings of the singularity type $Z_{15}$.
A suitable variation of the above degree six strategy yields all real schemes.
The realizability of all smooth degree 7 curves by T-curves is implicitly used in work of Brugall\'e and L\'opez de Medrano \cite{2012:Brugalle-Lopez}, who proved that regular T-curves satisfying a mild assumption on the lifting function are maximally inflected.

\printbibliography

\end{document}

%% file: tikz/deg2-111.tikz
\coordinate (A0) at (0, 0);
\coordinate (A1) at (0, 1);
\coordinate (A2) at (0, 2);
\coordinate (A3) at (1, 0);
\coordinate (A4) at (1, 1);
\coordinate (A5) at (2, 0);
\coordinate (A6) at (0, -1);
\coordinate (A7) at (0, -2);
\coordinate (A8) at (1, -1);
\coordinate (A9) at (-1, -1);
\coordinate (A10) at (-1, 0);
\coordinate (A11) at (-1, 1);
\coordinate (A12) at (-2, 0);
%\definecolor{color1}{rgb}{0.70907765858261,0.451838304829136,0.686808616435137}
\colorlet{color1}{mycolor4}
\fill[color1!50] ($(A0)! 0.5!(A10)$) -- ($(A1)! 0.5!(A10)$) -- (A10) -- cycle;
\fill[color1!50] ($(A0)! 0.5!(A6)$) -- ($(A3)! 0.5!(A6)$) -- (A6) -- cycle;
\fill[color1!50] (A6) -- ($(A6)! 0.5!(A0)$) -- ($(A10)! 0.5!(A0)$) -- (A10) -- cycle;
\fill[color1!50] ($(A1)! 0.5!(A11)$) -- ($(A2)! 0.5!(A11)$) -- (A11) -- cycle;
\fill[color1!50] ($(A1)! 0.5!(A10)$) -- ($(A12)! 0.5!(A10)$) -- (A10) -- cycle;
\fill[color1!50] ($(A1)! 0.5!(A11)$) -- ($(A12)! 0.5!(A11)$) -- (A11) -- cycle;
\fill[color1!50] ($(A3)! 0.5!(A6)$) -- ($(A5)! 0.5!(A6)$) -- (A6) -- cycle;
\fill[color1!50] (A6) -- ($(A6)! 0.5!(A5)$) -- ($(A8)! 0.5!(A5)$) -- (A8) -- cycle;
\fill[color1!50] (A6) -- ($(A6)! 0.5!(A7)$) -- ($(A8)! 0.5!(A7)$) -- (A8) -- cycle;
\fill[color1!50] ($(A7)! 0.5!(A6)$) -- ($(A9)! 0.5!(A6)$) -- (A6) -- cycle;
\fill[color1!50] ($(A9)! 0.5!(A6)$) -- ($(A12)! 0.5!(A6)$) -- (A6) -- cycle;
\fill[color1!50] (A6) -- ($(A6)! 0.5!(A12)$) -- ($(A10)! 0.5!(A12)$) -- (A10) -- cycle;
\draw[color1, thick] ($(A0)! 0.5!(A10)$) -- ($(A1)! 0.5!(A10)$);
\draw[color1, thick] ($(A0)! 0.5!(A6)$) -- ($(A3)! 0.5!(A6)$);
\draw[color1, thick] ($(A6)! 0.5!(A0)$) -- ($(A10)! 0.5!(A0)$);
\draw[color1, thick] ($(A1)! 0.5!(A11)$) -- ($(A2)! 0.5!(A11)$);
\draw[color1, thick] ($(A1)! 0.5!(A10)$) -- ($(A12)! 0.5!(A10)$);
\draw[color1, thick] ($(A1)! 0.5!(A11)$) -- ($(A12)! 0.5!(A11)$);
\draw[color1, thick] ($(A3)! 0.5!(A6)$) -- ($(A5)! 0.5!(A6)$);
\draw[color1, thick] ($(A6)! 0.5!(A5)$) -- ($(A8)! 0.5!(A5)$);
\draw[color1, thick] ($(A6)! 0.5!(A7)$) -- ($(A8)! 0.5!(A7)$);
\draw[color1, thick] ($(A7)! 0.5!(A6)$) -- ($(A9)! 0.5!(A6)$);
\draw[color1, thick] ($(A9)! 0.5!(A6)$) -- ($(A12)! 0.5!(A6)$);
\draw[color1, thick] ($(A6)! 0.5!(A12)$) -- ($(A10)! 0.5!(A12)$);
\foreach \a/\b/\c in {0/1/3,0/1/10,0/3/6,0/6/10,1/2/4,1/2/11,1/3/5,1/4/5,1/10/12,1/11/12,3/5/6,5/6/8,6/7/8,6/7/9,6/9/12,6/10/12}{
  \draw[black!30] (A\a) -- (A\b) -- (A\c) -- cycle;
}
\foreach \i in {0,1,2,3,4,5,7,9,12}{
  \fill[myred] (A\i) circle (3pt);
}
\foreach \i in {6,8,10,11}{
  \fill[myblue] (A\i) circle (3pt);
}
\foreach \i in {0,1,...,12}{
%  \node[anchor=north] at (A\i) {\tiny{$\i$}};
}

%% file: tikz/deg2-011.tikz
\coordinate (A0) at (0, 0);
\coordinate (A1) at (0, 1);
\coordinate (A2) at (0, 2);
\coordinate (A3) at (1, 0);
\coordinate (A4) at (1, 1);
\coordinate (A5) at (2, 0);
\coordinate (A6) at (0, -1);
\coordinate (A7) at (0, -2);
\coordinate (A8) at (1, -1);
\coordinate (A9) at (-1, -1);
\coordinate (A10) at (-1, 0);
\coordinate (A11) at (-1, 1);
\coordinate (A12) at (-2, 0);
%\definecolor{color1}{rgb}{0.739873893682572,0.702714891006087,0.763619429969696}
\colorlet{color1}{mycolor4}
\fill[color1!50] ($(A0)! 0.5!(A10)$) -- ($(A1)! 0.5!(A10)$) -- (A10) -- cycle;
\fill[color1!50] ($(A0)! 0.5!(A6)$) -- ($(A3)! 0.5!(A6)$) -- (A6) -- cycle;
\fill[color1!50] (A6) -- ($(A6)! 0.5!(A0)$) -- ($(A10)! 0.5!(A0)$) -- (A10) -- cycle;
\fill[color1!50] ($(A1)! 0.5!(A2)$) -- ($(A4)! 0.5!(A2)$) -- (A2) -- cycle;
\fill[color1!50] (A2) -- ($(A2)! 0.5!(A1)$) -- ($(A11)! 0.5!(A1)$) -- (A11) -- cycle;
\fill[color1!50] ($(A1)! 0.5!(A10)$) -- ($(A12)! 0.5!(A10)$) -- (A10) -- cycle;
\fill[color1!50] ($(A1)! 0.5!(A11)$) -- ($(A12)! 0.5!(A11)$) -- (A11) -- cycle;
\fill[color1!50] ($(A3)! 0.5!(A6)$) -- ($(A5)! 0.5!(A6)$) -- (A6) -- cycle;
\fill[color1!50] (A6) -- ($(A6)! 0.5!(A5)$) -- ($(A8)! 0.5!(A5)$) -- (A8) -- cycle;
\fill[color1!50] (A6) -- (A7) -- (A8) -- cycle;
\fill[color1!50] (A6) -- ($(A6)! 0.5!(A9)$) -- ($(A7)! 0.5!(A9)$) -- (A7) -- cycle;
\fill[color1!50] ($(A9)! 0.5!(A6)$) -- ($(A12)! 0.5!(A6)$) -- (A6) -- cycle;
\fill[color1!50] (A6) -- ($(A6)! 0.5!(A12)$) -- ($(A10)! 0.5!(A12)$) -- (A10) -- cycle;
\draw[color1, thick] ($(A0)! 0.5!(A10)$) -- ($(A1)! 0.5!(A10)$);
\draw[color1, thick] ($(A0)! 0.5!(A6)$) -- ($(A3)! 0.5!(A6)$);
\draw[color1, thick] ($(A6)! 0.5!(A0)$) -- ($(A10)! 0.5!(A0)$);
\draw[color1, thick] ($(A1)! 0.5!(A2)$) -- ($(A4)! 0.5!(A2)$);
\draw[color1, thick] ($(A2)! 0.5!(A1)$) -- ($(A11)! 0.5!(A1)$);
\draw[color1, thick] ($(A1)! 0.5!(A10)$) -- ($(A12)! 0.5!(A10)$);
\draw[color1, thick] ($(A1)! 0.5!(A11)$) -- ($(A12)! 0.5!(A11)$);
\draw[color1, thick] ($(A3)! 0.5!(A6)$) -- ($(A5)! 0.5!(A6)$);
\draw[color1, thick] ($(A6)! 0.5!(A5)$) -- ($(A8)! 0.5!(A5)$);
\draw[color1, thick] ($(A6)! 0.5!(A9)$) -- ($(A7)! 0.5!(A9)$);
\draw[color1, thick] ($(A9)! 0.5!(A6)$) -- ($(A12)! 0.5!(A6)$);
\draw[color1, thick] ($(A6)! 0.5!(A12)$) -- ($(A10)! 0.5!(A12)$);
\foreach \a/\b/\c in {0/1/3,0/1/10,0/3/6,0/6/10,1/2/4,1/2/11,1/3/5,1/4/5,1/10/12,1/11/12,3/5/6,5/6/8,6/7/8,6/7/9,6/9/12,6/10/12}{
  \draw[black!30] (A\a) -- (A\b) -- (A\c) -- cycle;
}
\foreach \i in {0,1,3,4,5,9,12}{
  \fill[myred] (A\i) circle (3pt);
}
\foreach \i in {2,6,7,8,10,11}{
  \fill[myblue] (A\i) circle (3pt);
}
\foreach \i in {0,1,...,12}{
%  \node[anchor=north] at (A\i) {\tiny{$\i$}};
}

%% file: tikz/deg2-101.tikz
\coordinate (A0) at (0, 0);
\coordinate (A1) at (0, 1);
\coordinate (A2) at (0, 2);
\coordinate (A3) at (1, 0);
\coordinate (A4) at (1, 1);
\coordinate (A5) at (2, 0);
\coordinate (A6) at (0, -1);
\coordinate (A7) at (0, -2);
\coordinate (A8) at (1, -1);
\coordinate (A9) at (-1, -1);
\coordinate (A10) at (-1, 0);
\coordinate (A11) at (-1, 1);
\coordinate (A12) at (-2, 0);
%\definecolor{color1}{rgb}{0.798033942313851,0.748473264788535,0.686213767980413}
\colorlet{color1}{mycolor4}
\fill[color1!50] ($(A0)! 0.5!(A10)$) -- ($(A1)! 0.5!(A10)$) -- (A10) -- cycle;
\fill[color1!50] ($(A0)! 0.5!(A6)$) -- ($(A3)! 0.5!(A6)$) -- (A6) -- cycle;
\fill[color1!50] (A6) -- ($(A6)! 0.5!(A0)$) -- ($(A10)! 0.5!(A0)$) -- (A10) -- cycle;
\fill[color1!50] ($(A1)! 0.5!(A4)$) -- ($(A2)! 0.5!(A4)$) -- (A4) -- cycle;
\fill[color1!50] ($(A1)! 0.5!(A4)$) -- ($(A5)! 0.5!(A4)$) -- (A4) -- cycle;
\fill[color1!50] ($(A1)! 0.5!(A10)$) -- ($(A12)! 0.5!(A10)$) -- (A10) -- cycle;
\fill[color1!50] ($(A3)! 0.5!(A6)$) -- ($(A5)! 0.5!(A6)$) -- (A6) -- cycle;
\fill[color1!50] ($(A5)! 0.5!(A6)$) -- ($(A8)! 0.5!(A6)$) -- (A6) -- cycle;
\fill[color1!50] ($(A7)! 0.5!(A6)$) -- ($(A8)! 0.5!(A6)$) -- (A6) -- cycle;
\fill[color1!50] (A6) -- ($(A6)! 0.5!(A7)$) -- ($(A9)! 0.5!(A7)$) -- (A9) -- cycle;
\fill[color1!50] (A6) -- ($(A6)! 0.5!(A12)$) -- ($(A9)! 0.5!(A12)$) -- (A9) -- cycle;
\fill[color1!50] (A6) -- ($(A6)! 0.5!(A12)$) -- ($(A10)! 0.5!(A12)$) -- (A10) -- cycle;
\draw[color1, thick] ($(A0)! 0.5!(A10)$) -- ($(A1)! 0.5!(A10)$);
\draw[color1, thick] ($(A0)! 0.5!(A6)$) -- ($(A3)! 0.5!(A6)$);
\draw[color1, thick] ($(A6)! 0.5!(A0)$) -- ($(A10)! 0.5!(A0)$);
\draw[color1, thick] ($(A1)! 0.5!(A4)$) -- ($(A2)! 0.5!(A4)$);
\draw[color1, thick] ($(A1)! 0.5!(A4)$) -- ($(A5)! 0.5!(A4)$);
\draw[color1, thick] ($(A1)! 0.5!(A10)$) -- ($(A12)! 0.5!(A10)$);
\draw[color1, thick] ($(A3)! 0.5!(A6)$) -- ($(A5)! 0.5!(A6)$);
\draw[color1, thick] ($(A5)! 0.5!(A6)$) -- ($(A8)! 0.5!(A6)$);
\draw[color1, thick] ($(A7)! 0.5!(A6)$) -- ($(A8)! 0.5!(A6)$);
\draw[color1, thick] ($(A6)! 0.5!(A7)$) -- ($(A9)! 0.5!(A7)$);
\draw[color1, thick] ($(A6)! 0.5!(A12)$) -- ($(A9)! 0.5!(A12)$);
\draw[color1, thick] ($(A6)! 0.5!(A12)$) -- ($(A10)! 0.5!(A12)$);
\foreach \a/\b/\c in {0/1/3,0/1/10,0/3/6,0/6/10,1/2/4,1/2/11,1/3/5,1/4/5,1/10/12,1/11/12,3/5/6,5/6/8,6/7/8,6/7/9,6/9/12,6/10/12}{
  \draw[black!30] (A\a) -- (A\b) -- (A\c) -- cycle;
}
\foreach \i in {0,1,2,3,5,7,8,11,12}{
  \fill[myred] (A\i) circle (3pt);
}
\foreach \i in {4,6,9,10}{
  \fill[myblue] (A\i) circle (3pt);
}
\foreach \i in {0,1,...,12}{
%  \node[anchor=north] at (A\i) {\tiny{$\i$}};
}

%% file: tikz/deg2-001.tikz
\coordinate (A0) at (0, 0);
\coordinate (A1) at (0, 1);
\coordinate (A2) at (0, 2);
\coordinate (A3) at (1, 0);
\coordinate (A4) at (1, 1);
\coordinate (A5) at (2, 0);
\coordinate (A6) at (0, -1);
\coordinate (A7) at (0, -2);
\coordinate (A8) at (1, -1);
\coordinate (A9) at (-1, -1);
\coordinate (A10) at (-1, 0);
\coordinate (A11) at (-1, 1);
\coordinate (A12) at (-2, 0);
%\definecolor{color1}{rgb}{0.634778965494326,0.872095711001444,0.729207131452429}
\colorlet{color1}{mycolor4}
\fill[color1!50] ($(A0)! 0.5!(A10)$) -- ($(A1)! 0.5!(A10)$) -- (A10) -- cycle;
\fill[color1!50] ($(A0)! 0.5!(A6)$) -- ($(A3)! 0.5!(A6)$) -- (A6) -- cycle;
\fill[color1!50] (A6) -- ($(A6)! 0.5!(A0)$) -- ($(A10)! 0.5!(A0)$) -- (A10) -- cycle;
\fill[color1!50] (A2) -- ($(A2)! 0.5!(A1)$) -- ($(A4)! 0.5!(A1)$) -- (A4) -- cycle;
\fill[color1!50] ($(A1)! 0.5!(A2)$) -- ($(A11)! 0.5!(A2)$) -- (A2) -- cycle;
\fill[color1!50] ($(A1)! 0.5!(A4)$) -- ($(A5)! 0.5!(A4)$) -- (A4) -- cycle;
\fill[color1!50] ($(A1)! 0.5!(A10)$) -- ($(A12)! 0.5!(A10)$) -- (A10) -- cycle;
\fill[color1!50] ($(A3)! 0.5!(A6)$) -- ($(A5)! 0.5!(A6)$) -- (A6) -- cycle;
\fill[color1!50] ($(A5)! 0.5!(A6)$) -- ($(A8)! 0.5!(A6)$) -- (A6) -- cycle;
\fill[color1!50] (A6) -- ($(A6)! 0.5!(A8)$) -- ($(A7)! 0.5!(A8)$) -- (A7) -- cycle;
\fill[color1!50] (A6) -- (A7) -- (A9) -- cycle;
\fill[color1!50] (A6) -- ($(A6)! 0.5!(A12)$) -- ($(A9)! 0.5!(A12)$) -- (A9) -- cycle;
\fill[color1!50] (A6) -- ($(A6)! 0.5!(A12)$) -- ($(A10)! 0.5!(A12)$) -- (A10) -- cycle;
\draw[color1, thick] ($(A0)! 0.5!(A10)$) -- ($(A1)! 0.5!(A10)$);
\draw[color1, thick] ($(A0)! 0.5!(A6)$) -- ($(A3)! 0.5!(A6)$);
\draw[color1, thick] ($(A6)! 0.5!(A0)$) -- ($(A10)! 0.5!(A0)$);
\draw[color1, thick] ($(A2)! 0.5!(A1)$) -- ($(A4)! 0.5!(A1)$);
\draw[color1, thick] ($(A1)! 0.5!(A2)$) -- ($(A11)! 0.5!(A2)$);
\draw[color1, thick] ($(A1)! 0.5!(A4)$) -- ($(A5)! 0.5!(A4)$);
\draw[color1, thick] ($(A1)! 0.5!(A10)$) -- ($(A12)! 0.5!(A10)$);
\draw[color1, thick] ($(A3)! 0.5!(A6)$) -- ($(A5)! 0.5!(A6)$);
\draw[color1, thick] ($(A5)! 0.5!(A6)$) -- ($(A8)! 0.5!(A6)$);
\draw[color1, thick] ($(A6)! 0.5!(A8)$) -- ($(A7)! 0.5!(A8)$);
\draw[color1, thick] ($(A6)! 0.5!(A12)$) -- ($(A9)! 0.5!(A12)$);
\draw[color1, thick] ($(A6)! 0.5!(A12)$) -- ($(A10)! 0.5!(A12)$);
\foreach \a/\b/\c in {0/1/3,0/1/10,0/3/6,0/6/10,1/2/4,1/2/11,1/3/5,1/4/5,1/10/12,1/11/12,3/5/6,5/6/8,6/7/8,6/7/9,6/9/12,6/10/12}{
  \draw[black!30] (A\a) -- (A\b) -- (A\c) -- cycle;
}
\foreach \i in {0,1,3,5,8,11,12}{
  \fill[myred] (A\i) circle (3pt);
}
\foreach \i in {2,4,6,7,9,10}{
  \fill[myblue] (A\i) circle (3pt);
}
\foreach \i in {0,1,...,12}{
%  \node[anchor=north] at (A\i) {\tiny{$\i$}};
}

%% file: tikz/deg2-110.tikz
\coordinate (A0) at (0, 0);
\coordinate (A1) at (0, 1);
\coordinate (A2) at (0, 2);
\coordinate (A3) at (1, 0);
\coordinate (A4) at (1, 1);
\coordinate (A5) at (2, 0);
\coordinate (A6) at (0, -1);
\coordinate (A7) at (0, -2);
\coordinate (A8) at (1, -1);
\coordinate (A9) at (-1, -1);
\coordinate (A10) at (-1, 0);
\coordinate (A11) at (-1, 1);
\coordinate (A12) at (-2, 0);
%\definecolor{color1}{rgb}{0.816460005774183,0.669287465165293,0.795344288523899}
\colorlet{color1}{mycolor4}
\fill[color1!50] (A0) -- (A1) -- (A3) -- cycle;
\fill[color1!50] (A0) -- ($(A0)! 0.5!(A10)$) -- ($(A1)! 0.5!(A10)$) -- (A1) -- cycle;
\fill[color1!50] (A0) -- ($(A0)! 0.5!(A6)$) -- ($(A3)! 0.5!(A6)$) -- (A3) -- cycle;
\fill[color1!50] ($(A6)! 0.5!(A0)$) -- ($(A10)! 0.5!(A0)$) -- (A0) -- cycle;
\fill[color1!50] (A1) -- (A2) -- (A4) -- cycle;
\fill[color1!50] (A1) -- ($(A1)! 0.5!(A11)$) -- ($(A2)! 0.5!(A11)$) -- (A2) -- cycle;
\fill[color1!50] (A1) -- ($(A1)! 0.5!(A5)$) -- ($(A3)! 0.5!(A5)$) -- (A3) -- cycle;
\fill[color1!50] (A1) -- ($(A1)! 0.5!(A5)$) -- ($(A4)! 0.5!(A5)$) -- (A4) -- cycle;
\fill[color1!50] ($(A10)! 0.5!(A1)$) -- ($(A12)! 0.5!(A1)$) -- (A1) -- cycle;
\fill[color1!50] ($(A11)! 0.5!(A1)$) -- ($(A12)! 0.5!(A1)$) -- (A1) -- cycle;
\fill[color1!50] ($(A5)! 0.5!(A3)$) -- ($(A6)! 0.5!(A3)$) -- (A3) -- cycle;
\fill[color1!50] ($(A6)! 0.5!(A7)$) -- ($(A8)! 0.5!(A7)$) -- (A7) -- cycle;
\fill[color1!50] (A7) -- ($(A7)! 0.5!(A6)$) -- ($(A9)! 0.5!(A6)$) -- (A9) -- cycle;
\fill[color1!50] ($(A6)! 0.5!(A9)$) -- ($(A12)! 0.5!(A9)$) -- (A9) -- cycle;
\draw[color1, thick] ($(A0)! 0.5!(A10)$) -- ($(A1)! 0.5!(A10)$);
\draw[color1, thick] ($(A0)! 0.5!(A6)$) -- ($(A3)! 0.5!(A6)$);
\draw[color1, thick] ($(A6)! 0.5!(A0)$) -- ($(A10)! 0.5!(A0)$);
\draw[color1, thick] ($(A1)! 0.5!(A11)$) -- ($(A2)! 0.5!(A11)$);
\draw[color1, thick] ($(A1)! 0.5!(A5)$) -- ($(A3)! 0.5!(A5)$);
\draw[color1, thick] ($(A1)! 0.5!(A5)$) -- ($(A4)! 0.5!(A5)$);
\draw[color1, thick] ($(A10)! 0.5!(A1)$) -- ($(A12)! 0.5!(A1)$);
\draw[color1, thick] ($(A11)! 0.5!(A1)$) -- ($(A12)! 0.5!(A1)$);
\draw[color1, thick] ($(A5)! 0.5!(A3)$) -- ($(A6)! 0.5!(A3)$);
\draw[color1, thick] ($(A6)! 0.5!(A7)$) -- ($(A8)! 0.5!(A7)$);
\draw[color1, thick] ($(A7)! 0.5!(A6)$) -- ($(A9)! 0.5!(A6)$);
\draw[color1, thick] ($(A6)! 0.5!(A9)$) -- ($(A12)! 0.5!(A9)$);
\foreach \a/\b/\c in {0/1/3,0/1/10,0/3/6,0/6/10,1/2/4,1/2/11,1/3/5,1/4/5,1/10/12,1/11/12,3/5/6,5/6/8,6/7/8,6/7/9,6/9/12,6/10/12}{
  \draw[black!30] (A\a) -- (A\b) -- (A\c) -- cycle;
}
\foreach \i in {0,1,2,3,4,7,9}{
  \fill[myred] (A\i) circle (3pt);
}
\foreach \i in {5,6,8,10,11,12}{
  \fill[myblue] (A\i) circle (3pt);
}
\foreach \i in {0,1,...,12}{
%  \node[anchor=north] at (A\i) {\tiny{$\i$}};
}

%% file: tikz/deg2-010.tikz
\coordinate (A0) at (0, 0);
\coordinate (A1) at (0, 1);
\coordinate (A2) at (0, 2);
\coordinate (A3) at (1, 0);
\coordinate (A4) at (1, 1);
\coordinate (A5) at (2, 0);
\coordinate (A6) at (0, -1);
\coordinate (A7) at (0, -2);
\coordinate (A8) at (1, -1);
\coordinate (A9) at (-1, -1);
\coordinate (A10) at (-1, 0);
\coordinate (A11) at (-1, 1);
\coordinate (A12) at (-2, 0);
%\definecolor{color1}{rgb}{0.779235277740695,0.568443577831131,0.847682146429805}
\colorlet{color1}{mycolor4}
\fill[color1!50] (A0) -- (A1) -- (A3) -- cycle;
\fill[color1!50] (A0) -- ($(A0)! 0.5!(A10)$) -- ($(A1)! 0.5!(A10)$) -- (A1) -- cycle;
\fill[color1!50] (A0) -- ($(A0)! 0.5!(A6)$) -- ($(A3)! 0.5!(A6)$) -- (A3) -- cycle;
\fill[color1!50] ($(A6)! 0.5!(A0)$) -- ($(A10)! 0.5!(A0)$) -- (A0) -- cycle;
\fill[color1!50] (A1) -- ($(A1)! 0.5!(A2)$) -- ($(A4)! 0.5!(A2)$) -- (A4) -- cycle;
\fill[color1!50] ($(A2)! 0.5!(A1)$) -- ($(A11)! 0.5!(A1)$) -- (A1) -- cycle;
\fill[color1!50] (A1) -- ($(A1)! 0.5!(A5)$) -- ($(A3)! 0.5!(A5)$) -- (A3) -- cycle;
\fill[color1!50] (A1) -- ($(A1)! 0.5!(A5)$) -- ($(A4)! 0.5!(A5)$) -- (A4) -- cycle;
\fill[color1!50] ($(A10)! 0.5!(A1)$) -- ($(A12)! 0.5!(A1)$) -- (A1) -- cycle;
\fill[color1!50] ($(A11)! 0.5!(A1)$) -- ($(A12)! 0.5!(A1)$) -- (A1) -- cycle;
\fill[color1!50] ($(A5)! 0.5!(A3)$) -- ($(A6)! 0.5!(A3)$) -- (A3) -- cycle;
\fill[color1!50] ($(A6)! 0.5!(A9)$) -- ($(A7)! 0.5!(A9)$) -- (A9) -- cycle;
\fill[color1!50] ($(A6)! 0.5!(A9)$) -- ($(A12)! 0.5!(A9)$) -- (A9) -- cycle;
\draw[color1, thick] ($(A0)! 0.5!(A10)$) -- ($(A1)! 0.5!(A10)$);
\draw[color1, thick] ($(A0)! 0.5!(A6)$) -- ($(A3)! 0.5!(A6)$);
\draw[color1, thick] ($(A6)! 0.5!(A0)$) -- ($(A10)! 0.5!(A0)$);
\draw[color1, thick] ($(A1)! 0.5!(A2)$) -- ($(A4)! 0.5!(A2)$);
\draw[color1, thick] ($(A2)! 0.5!(A1)$) -- ($(A11)! 0.5!(A1)$);
\draw[color1, thick] ($(A1)! 0.5!(A5)$) -- ($(A3)! 0.5!(A5)$);
\draw[color1, thick] ($(A1)! 0.5!(A5)$) -- ($(A4)! 0.5!(A5)$);
\draw[color1, thick] ($(A10)! 0.5!(A1)$) -- ($(A12)! 0.5!(A1)$);
\draw[color1, thick] ($(A11)! 0.5!(A1)$) -- ($(A12)! 0.5!(A1)$);
\draw[color1, thick] ($(A5)! 0.5!(A3)$) -- ($(A6)! 0.5!(A3)$);
\draw[color1, thick] ($(A6)! 0.5!(A9)$) -- ($(A7)! 0.5!(A9)$);
\draw[color1, thick] ($(A6)! 0.5!(A9)$) -- ($(A12)! 0.5!(A9)$);
\foreach \a/\b/\c in {0/1/3,0/1/10,0/3/6,0/6/10,1/2/4,1/2/11,1/3/5,1/4/5,1/10/12,1/11/12,3/5/6,5/6/8,6/7/8,6/7/9,6/9/12,6/10/12}{
  \draw[black!30] (A\a) -- (A\b) -- (A\c) -- cycle;
}
\foreach \i in {0,1,3,4,9}{
  \fill[myred] (A\i) circle (3pt);
}
\foreach \i in {2,5,6,7,8,10,11,12}{
  \fill[myblue] (A\i) circle (3pt);
}
\foreach \i in {0,1,...,12}{
%  \node[anchor=north] at (A\i) {\tiny{$\i$}};
}

%% file: tikz/deg2-100.tikz
\coordinate (A0) at (0, 0);
\coordinate (A1) at (0, 1);
\coordinate (A2) at (0, 2);
\coordinate (A3) at (1, 0);
\coordinate (A4) at (1, 1);
\coordinate (A5) at (2, 0);
\coordinate (A6) at (0, -1);
\coordinate (A7) at (0, -2);
\coordinate (A8) at (1, -1);
\coordinate (A9) at (-1, -1);
\coordinate (A10) at (-1, 0);
\coordinate (A11) at (-1, 1);
\coordinate (A12) at (-2, 0);
%\definecolor{color1}{rgb}{0.697363765107675,0.587382114312838,0.762524557931145}
\colorlet{color1}{mycolor4}
\fill[color1!50] (A0) -- (A1) -- (A3) -- cycle;
\fill[color1!50] (A0) -- ($(A0)! 0.5!(A10)$) -- ($(A1)! 0.5!(A10)$) -- (A1) -- cycle;
\fill[color1!50] (A0) -- ($(A0)! 0.5!(A6)$) -- ($(A3)! 0.5!(A6)$) -- (A3) -- cycle;
\fill[color1!50] ($(A6)! 0.5!(A0)$) -- ($(A10)! 0.5!(A0)$) -- (A0) -- cycle;
\fill[color1!50] (A1) -- ($(A1)! 0.5!(A4)$) -- ($(A2)! 0.5!(A4)$) -- (A2) -- cycle;
\fill[color1!50] (A1) -- (A2) -- (A11) -- cycle;
\fill[color1!50] (A1) -- ($(A1)! 0.5!(A5)$) -- ($(A3)! 0.5!(A5)$) -- (A3) -- cycle;
\fill[color1!50] ($(A4)! 0.5!(A1)$) -- ($(A5)! 0.5!(A1)$) -- (A1) -- cycle;
\fill[color1!50] ($(A10)! 0.5!(A1)$) -- ($(A12)! 0.5!(A1)$) -- (A1) -- cycle;
\fill[color1!50] (A1) -- ($(A1)! 0.5!(A12)$) -- ($(A11)! 0.5!(A12)$) -- (A11) -- cycle;
\fill[color1!50] ($(A5)! 0.5!(A3)$) -- ($(A6)! 0.5!(A3)$) -- (A3) -- cycle;
\fill[color1!50] ($(A5)! 0.5!(A8)$) -- ($(A6)! 0.5!(A8)$) -- (A8) -- cycle;
\fill[color1!50] (A7) -- ($(A7)! 0.5!(A6)$) -- ($(A8)! 0.5!(A6)$) -- (A8) -- cycle;
\fill[color1!50] ($(A6)! 0.5!(A7)$) -- ($(A9)! 0.5!(A7)$) -- (A7) -- cycle;
\draw[color1, thick] ($(A0)! 0.5!(A10)$) -- ($(A1)! 0.5!(A10)$);
\draw[color1, thick] ($(A0)! 0.5!(A6)$) -- ($(A3)! 0.5!(A6)$);
\draw[color1, thick] ($(A6)! 0.5!(A0)$) -- ($(A10)! 0.5!(A0)$);
\draw[color1, thick] ($(A1)! 0.5!(A4)$) -- ($(A2)! 0.5!(A4)$);
\draw[color1, thick] ($(A1)! 0.5!(A5)$) -- ($(A3)! 0.5!(A5)$);
\draw[color1, thick] ($(A4)! 0.5!(A1)$) -- ($(A5)! 0.5!(A1)$);
\draw[color1, thick] ($(A10)! 0.5!(A1)$) -- ($(A12)! 0.5!(A1)$);
\draw[color1, thick] ($(A1)! 0.5!(A12)$) -- ($(A11)! 0.5!(A12)$);
\draw[color1, thick] ($(A5)! 0.5!(A3)$) -- ($(A6)! 0.5!(A3)$);
\draw[color1, thick] ($(A5)! 0.5!(A8)$) -- ($(A6)! 0.5!(A8)$);
\draw[color1, thick] ($(A7)! 0.5!(A6)$) -- ($(A8)! 0.5!(A6)$);
\draw[color1, thick] ($(A6)! 0.5!(A7)$) -- ($(A9)! 0.5!(A7)$);
\foreach \a/\b/\c in {0/1/3,0/1/10,0/3/6,0/6/10,1/2/4,1/2/11,1/3/5,1/4/5,1/10/12,1/11/12,3/5/6,5/6/8,6/7/8,6/7/9,6/9/12,6/10/12}{
  \draw[black!30] (A\a) -- (A\b) -- (A\c) -- cycle;
}
\foreach \i in {0,1,2,3,7,8,11}{
  \fill[myred] (A\i) circle (3pt);
}
\foreach \i in {4,5,6,9,10,12}{
  \fill[myblue] (A\i) circle (3pt);
}
\foreach \i in {0,1,...,12}{
%  \node[anchor=north] at (A\i) {\tiny{$\i$}};
}

%% file: tikz/deg2-000.tikz
\coordinate (A0) at (0, 0);
\coordinate (A1) at (0, 1);
\coordinate (A2) at (0, 2);
\coordinate (A3) at (1, 0);
\coordinate (A4) at (1, 1);
\coordinate (A5) at (2, 0);
\coordinate (A6) at (0, -1);
\coordinate (A7) at (0, -2);
\coordinate (A8) at (1, -1);
\coordinate (A9) at (-1, -1);
\coordinate (A10) at (-1, 0);
\coordinate (A11) at (-1, 1);
\coordinate (A12) at (-2, 0);
%\definecolor{color1}{rgb}{0.655585796009128,0.520126453511578,0.784413843115343}
\colorlet{color1}{mycolor4}
\fill[color1!50] (A0) -- (A1) -- (A3) -- cycle;
\fill[color1!50] (A0) -- ($(A0)! 0.5!(A10)$) -- ($(A1)! 0.5!(A10)$) -- (A1) -- cycle;
\fill[color1!50] (A0) -- ($(A0)! 0.5!(A6)$) -- ($(A3)! 0.5!(A6)$) -- (A3) -- cycle;
\fill[color1!50] ($(A6)! 0.5!(A0)$) -- ($(A10)! 0.5!(A0)$) -- (A0) -- cycle;
\fill[color1!50] ($(A2)! 0.5!(A1)$) -- ($(A4)! 0.5!(A1)$) -- (A1) -- cycle;
\fill[color1!50] (A1) -- ($(A1)! 0.5!(A2)$) -- ($(A11)! 0.5!(A2)$) -- (A11) -- cycle;
\fill[color1!50] (A1) -- ($(A1)! 0.5!(A5)$) -- ($(A3)! 0.5!(A5)$) -- (A3) -- cycle;
\fill[color1!50] ($(A4)! 0.5!(A1)$) -- ($(A5)! 0.5!(A1)$) -- (A1) -- cycle;
\fill[color1!50] ($(A10)! 0.5!(A1)$) -- ($(A12)! 0.5!(A1)$) -- (A1) -- cycle;
\fill[color1!50] (A1) -- ($(A1)! 0.5!(A12)$) -- ($(A11)! 0.5!(A12)$) -- (A11) -- cycle;
\fill[color1!50] ($(A5)! 0.5!(A3)$) -- ($(A6)! 0.5!(A3)$) -- (A3) -- cycle;
\fill[color1!50] ($(A5)! 0.5!(A8)$) -- ($(A6)! 0.5!(A8)$) -- (A8) -- cycle;
\fill[color1!50] ($(A6)! 0.5!(A8)$) -- ($(A7)! 0.5!(A8)$) -- (A8) -- cycle;
\draw[color1, thick] ($(A0)! 0.5!(A10)$) -- ($(A1)! 0.5!(A10)$);
\draw[color1, thick] ($(A0)! 0.5!(A6)$) -- ($(A3)! 0.5!(A6)$);
\draw[color1, thick] ($(A6)! 0.5!(A0)$) -- ($(A10)! 0.5!(A0)$);
\draw[color1, thick] ($(A2)! 0.5!(A1)$) -- ($(A4)! 0.5!(A1)$);
\draw[color1, thick] ($(A1)! 0.5!(A2)$) -- ($(A11)! 0.5!(A2)$);
\draw[color1, thick] ($(A1)! 0.5!(A5)$) -- ($(A3)! 0.5!(A5)$);
\draw[color1, thick] ($(A4)! 0.5!(A1)$) -- ($(A5)! 0.5!(A1)$);
\draw[color1, thick] ($(A10)! 0.5!(A1)$) -- ($(A12)! 0.5!(A1)$);
\draw[color1, thick] ($(A1)! 0.5!(A12)$) -- ($(A11)! 0.5!(A12)$);
\draw[color1, thick] ($(A5)! 0.5!(A3)$) -- ($(A6)! 0.5!(A3)$);
\draw[color1, thick] ($(A5)! 0.5!(A8)$) -- ($(A6)! 0.5!(A8)$);
\draw[color1, thick] ($(A6)! 0.5!(A8)$) -- ($(A7)! 0.5!(A8)$);
\foreach \a/\b/\c in {0/1/3,0/1/10,0/3/6,0/6/10,1/2/4,1/2/11,1/3/5,1/4/5,1/10/12,1/11/12,3/5/6,5/6/8,6/7/8,6/7/9,6/9/12,6/10/12}{
  \draw[black!30] (A\a) -- (A\b) -- (A\c) -- cycle;
}
\foreach \i in {0,1,3,8,11}{
  \fill[myred] (A\i) circle (3pt);
}
\foreach \i in {2,4,5,6,7,9,10,12}{
  \fill[myblue] (A\i) circle (3pt);
}
\foreach \i in {0,1,...,12}{
%  \node[anchor=north] at (A\i) {\tiny{$\i$}};
}

%% file: tikz/harnack_honeycomb_6.tikz
\coordinate (A0) at (0, 0);
\coordinate (A1) at (0, 1);
\coordinate (A2) at (0, 2);
\coordinate (A3) at (0, 3);
\coordinate (A4) at (0, 4);
\coordinate (A5) at (0, 5);
\coordinate (A6) at (0, 6);
\coordinate (A7) at (1, 0);
\coordinate (A8) at (1, 1);
\coordinate (A9) at (1, 2);
\coordinate (A10) at (1, 3);
\coordinate (A11) at (1, 4);
\coordinate (A12) at (1, 5);
\coordinate (A13) at (2, 0);
\coordinate (A14) at (2, 1);
\coordinate (A15) at (2, 2);
\coordinate (A16) at (2, 3);
\coordinate (A17) at (2, 4);
\coordinate (A18) at (3, 0);
\coordinate (A19) at (3, 1);
\coordinate (A20) at (3, 2);
\coordinate (A21) at (3, 3);
\coordinate (A22) at (4, 0);
\coordinate (A23) at (4, 1);
\coordinate (A24) at (4, 2);
\coordinate (A25) at (5, 0);
\coordinate (A26) at (5, 1);
\coordinate (A27) at (6, 0);
\coordinate (A28) at (0, -1);
\coordinate (A29) at (0, -2);
\coordinate (A30) at (0, -3);
\coordinate (A31) at (0, -4);
\coordinate (A32) at (0, -5);
\coordinate (A33) at (0, -6);
\coordinate (A34) at (1, -1);
\coordinate (A35) at (1, -2);
\coordinate (A36) at (1, -3);
\coordinate (A37) at (1, -4);
\coordinate (A38) at (1, -5);
\coordinate (A39) at (2, -1);
\coordinate (A40) at (2, -2);
\coordinate (A41) at (2, -3);
\coordinate (A42) at (2, -4);
\coordinate (A43) at (3, -1);
\coordinate (A44) at (3, -2);
\coordinate (A45) at (3, -3);
\coordinate (A46) at (4, -1);
\coordinate (A47) at (4, -2);
\coordinate (A48) at (5, -1);
\coordinate (A49) at (-1, -5);
\coordinate (A50) at (-1, -4);
\coordinate (A51) at (-1, -3);
\coordinate (A52) at (-1, -2);
\coordinate (A53) at (-1, -1);
\coordinate (A54) at (-1, 0);
\coordinate (A55) at (-1, 1);
\coordinate (A56) at (-1, 2);
\coordinate (A57) at (-1, 3);
\coordinate (A58) at (-1, 4);
\coordinate (A59) at (-1, 5);
\coordinate (A60) at (-2, -4);
\coordinate (A61) at (-2, -3);
\coordinate (A62) at (-2, -2);
\coordinate (A63) at (-2, -1);
\coordinate (A64) at (-2, 0);
\coordinate (A65) at (-2, 1);
\coordinate (A66) at (-2, 2);
\coordinate (A67) at (-2, 3);
\coordinate (A68) at (-2, 4);
\coordinate (A69) at (-3, -3);
\coordinate (A70) at (-3, -2);
\coordinate (A71) at (-3, -1);
\coordinate (A72) at (-3, 0);
\coordinate (A73) at (-3, 1);
\coordinate (A74) at (-3, 2);
\coordinate (A75) at (-3, 3);
\coordinate (A76) at (-4, -2);
\coordinate (A77) at (-4, -1);
\coordinate (A78) at (-4, 0);
\coordinate (A79) at (-4, 1);
\coordinate (A80) at (-4, 2);
\coordinate (A81) at (-5, -1);
\coordinate (A82) at (-5, 0);
\coordinate (A83) at (-5, 1);
\coordinate (A84) at (-6, 0);
\colorlet{color1}{mycolor4}
\fill[color1!50] (A1) -- ($(A1)! 0.5!(A0)$) -- ($(A7)! 0.5!(A0)$) -- (A7) -- cycle;
\fill[color1!50] ($(A0)! 0.5!(A1)$) -- ($(A54)! 0.5!(A1)$) -- (A1) -- cycle;
\fill[color1!50] ($(A0)! 0.5!(A7)$) -- ($(A28)! 0.5!(A7)$) -- (A7) -- cycle;
\fill[color1!50] (A1) -- ($(A1)! 0.5!(A2)$) -- ($(A8)! 0.5!(A2)$) -- (A8) -- cycle;
\fill[color1!50] ($(A2)! 0.5!(A1)$) -- ($(A55)! 0.5!(A1)$) -- (A1) -- cycle;
\fill[color1!50] (A1) -- (A7) -- (A8) -- cycle;
\fill[color1!50] ($(A54)! 0.5!(A1)$) -- ($(A55)! 0.5!(A1)$) -- (A1) -- cycle;
\fill[color1!50] (A3) -- ($(A3)! 0.5!(A2)$) -- ($(A9)! 0.5!(A2)$) -- (A9) -- cycle;
\fill[color1!50] ($(A2)! 0.5!(A3)$) -- ($(A56)! 0.5!(A3)$) -- (A3) -- cycle;
\fill[color1!50] (A8) -- ($(A8)! 0.5!(A2)$) -- ($(A9)! 0.5!(A2)$) -- (A9) -- cycle;
\fill[color1!50] (A3) -- ($(A3)! 0.5!(A4)$) -- ($(A10)! 0.5!(A4)$) -- (A10) -- cycle;
\fill[color1!50] ($(A4)! 0.5!(A3)$) -- ($(A57)! 0.5!(A3)$) -- (A3) -- cycle;
\fill[color1!50] (A3) -- (A9) -- (A10) -- cycle;
\fill[color1!50] ($(A56)! 0.5!(A3)$) -- ($(A57)! 0.5!(A3)$) -- (A3) -- cycle;
\fill[color1!50] (A5) -- ($(A5)! 0.5!(A4)$) -- ($(A11)! 0.5!(A4)$) -- (A11) -- cycle;
\fill[color1!50] ($(A4)! 0.5!(A5)$) -- ($(A58)! 0.5!(A5)$) -- (A5) -- cycle;
\fill[color1!50] (A10) -- ($(A10)! 0.5!(A4)$) -- ($(A11)! 0.5!(A4)$) -- (A11) -- cycle;
\fill[color1!50] (A5) -- ($(A5)! 0.5!(A6)$) -- ($(A12)! 0.5!(A6)$) -- (A12) -- cycle;
\fill[color1!50] ($(A6)! 0.5!(A5)$) -- ($(A59)! 0.5!(A5)$) -- (A5) -- cycle;
\fill[color1!50] (A5) -- (A11) -- (A12) -- cycle;
\fill[color1!50] ($(A58)! 0.5!(A5)$) -- ($(A59)! 0.5!(A5)$) -- (A5) -- cycle;
\fill[color1!50] (A7) -- ($(A7)! 0.5!(A13)$) -- ($(A8)! 0.5!(A13)$) -- (A8) -- cycle;
\fill[color1!50] ($(A13)! 0.5!(A7)$) -- ($(A34)! 0.5!(A7)$) -- (A7) -- cycle;
\fill[color1!50] ($(A28)! 0.5!(A7)$) -- ($(A34)! 0.5!(A7)$) -- (A7) -- cycle;
\fill[color1!50] (A8) -- (A9) -- (A14) -- cycle;
\fill[color1!50] (A8) -- ($(A8)! 0.5!(A13)$) -- ($(A14)! 0.5!(A13)$) -- (A14) -- cycle;
\fill[color1!50] (A9) -- ($(A9)! 0.5!(A15)$) -- ($(A10)! 0.5!(A15)$) -- (A10) -- cycle;
\fill[color1!50] (A9) -- ($(A9)! 0.5!(A15)$) -- ($(A14)! 0.5!(A15)$) -- (A14) -- cycle;
\fill[color1!50] (A10) -- (A11) -- (A16) -- cycle;
\fill[color1!50] (A10) -- ($(A10)! 0.5!(A15)$) -- ($(A16)! 0.5!(A15)$) -- (A16) -- cycle;
\fill[color1!50] (A11) -- ($(A11)! 0.5!(A17)$) -- ($(A12)! 0.5!(A17)$) -- (A12) -- cycle;
\fill[color1!50] (A11) -- ($(A11)! 0.5!(A17)$) -- ($(A16)! 0.5!(A17)$) -- (A16) -- cycle;
\fill[color1!50] (A14) -- ($(A14)! 0.5!(A13)$) -- ($(A18)! 0.5!(A13)$) -- (A18) -- cycle;
\fill[color1!50] ($(A13)! 0.5!(A18)$) -- ($(A39)! 0.5!(A18)$) -- (A18) -- cycle;
\fill[color1!50] (A14) -- ($(A14)! 0.5!(A15)$) -- ($(A19)! 0.5!(A15)$) -- (A19) -- cycle;
\fill[color1!50] (A14) -- (A18) -- (A19) -- cycle;
\fill[color1!50] (A16) -- ($(A16)! 0.5!(A15)$) -- ($(A20)! 0.5!(A15)$) -- (A20) -- cycle;
\fill[color1!50] (A19) -- ($(A19)! 0.5!(A15)$) -- ($(A20)! 0.5!(A15)$) -- (A20) -- cycle;
\fill[color1!50] (A16) -- ($(A16)! 0.5!(A17)$) -- ($(A21)! 0.5!(A17)$) -- (A21) -- cycle;
\fill[color1!50] (A16) -- (A20) -- (A21) -- cycle;
\fill[color1!50] (A18) -- ($(A18)! 0.5!(A22)$) -- ($(A19)! 0.5!(A22)$) -- (A19) -- cycle;
\fill[color1!50] ($(A22)! 0.5!(A18)$) -- ($(A43)! 0.5!(A18)$) -- (A18) -- cycle;
\fill[color1!50] ($(A39)! 0.5!(A18)$) -- ($(A43)! 0.5!(A18)$) -- (A18) -- cycle;
\fill[color1!50] (A19) -- (A20) -- (A23) -- cycle;
\fill[color1!50] (A19) -- ($(A19)! 0.5!(A22)$) -- ($(A23)! 0.5!(A22)$) -- (A23) -- cycle;
\fill[color1!50] (A20) -- ($(A20)! 0.5!(A24)$) -- ($(A21)! 0.5!(A24)$) -- (A21) -- cycle;
\fill[color1!50] (A20) -- ($(A20)! 0.5!(A24)$) -- ($(A23)! 0.5!(A24)$) -- (A23) -- cycle;
\fill[color1!50] (A23) -- ($(A23)! 0.5!(A22)$) -- ($(A25)! 0.5!(A22)$) -- (A25) -- cycle;
\fill[color1!50] ($(A22)! 0.5!(A25)$) -- ($(A46)! 0.5!(A25)$) -- (A25) -- cycle;
\fill[color1!50] (A23) -- ($(A23)! 0.5!(A24)$) -- ($(A26)! 0.5!(A24)$) -- (A26) -- cycle;
\fill[color1!50] (A23) -- (A25) -- (A26) -- cycle;
\fill[color1!50] (A25) -- ($(A25)! 0.5!(A27)$) -- ($(A26)! 0.5!(A27)$) -- (A26) -- cycle;
\fill[color1!50] ($(A27)! 0.5!(A25)$) -- ($(A48)! 0.5!(A25)$) -- (A25) -- cycle;
\fill[color1!50] ($(A46)! 0.5!(A25)$) -- ($(A48)! 0.5!(A25)$) -- (A25) -- cycle;
\fill[color1!50] ($(A32)! 0.5!(A49)$) -- ($(A33)! 0.5!(A49)$) -- (A49) -- cycle;
\fill[color1!50] ($(A32)! 0.5!(A49)$) -- ($(A50)! 0.5!(A49)$) -- (A49) -- cycle;
\fill[color1!50] ($(A50)! 0.5!(A49)$) -- ($(A60)! 0.5!(A49)$) -- (A49) -- cycle;
\fill[color1!50] ($(A60)! 0.5!(A69)$) -- ($(A61)! 0.5!(A69)$) -- (A69) -- cycle;
\fill[color1!50] ($(A61)! 0.5!(A69)$) -- ($(A70)! 0.5!(A69)$) -- (A69) -- cycle;
\fill[color1!50] ($(A70)! 0.5!(A69)$) -- ($(A76)! 0.5!(A69)$) -- (A69) -- cycle;
\fill[color1!50] ($(A76)! 0.5!(A81)$) -- ($(A77)! 0.5!(A81)$) -- (A81) -- cycle;
\fill[color1!50] ($(A77)! 0.5!(A81)$) -- ($(A82)! 0.5!(A81)$) -- (A81) -- cycle;
\fill[color1!50] ($(A82)! 0.5!(A81)$) -- ($(A84)! 0.5!(A81)$) -- (A81) -- cycle;
\colorlet{color2}{mycolor3}
\fill[color2!50] ($(A9)! 0.5!(A15)$) -- ($(A10)! 0.5!(A15)$) -- (A15) -- cycle;
\fill[color2!50] ($(A9)! 0.5!(A15)$) -- ($(A14)! 0.5!(A15)$) -- (A15) -- cycle;
\fill[color2!50] ($(A10)! 0.5!(A15)$) -- ($(A16)! 0.5!(A15)$) -- (A15) -- cycle;
\fill[color2!50] ($(A14)! 0.5!(A15)$) -- ($(A19)! 0.5!(A15)$) -- (A15) -- cycle;
\fill[color2!50] ($(A16)! 0.5!(A15)$) -- ($(A20)! 0.5!(A15)$) -- (A15) -- cycle;
\fill[color2!50] ($(A19)! 0.5!(A15)$) -- ($(A20)! 0.5!(A15)$) -- (A15) -- cycle;
\colorlet{color3}{mycolor3}
\fill[color3!50] ($(A29)! 0.5!(A35)$) -- ($(A30)! 0.5!(A35)$) -- (A35) -- cycle;
\fill[color3!50] ($(A29)! 0.5!(A35)$) -- ($(A34)! 0.5!(A35)$) -- (A35) -- cycle;
\fill[color3!50] ($(A30)! 0.5!(A35)$) -- ($(A36)! 0.5!(A35)$) -- (A35) -- cycle;
\fill[color3!50] ($(A34)! 0.5!(A35)$) -- ($(A39)! 0.5!(A35)$) -- (A35) -- cycle;
\fill[color3!50] ($(A36)! 0.5!(A35)$) -- ($(A40)! 0.5!(A35)$) -- (A35) -- cycle;
\fill[color3!50] ($(A39)! 0.5!(A35)$) -- ($(A40)! 0.5!(A35)$) -- (A35) -- cycle;
\colorlet{color4}{mycolor3}
\fill[color4!50] ($(A31)! 0.5!(A37)$) -- ($(A32)! 0.5!(A37)$) -- (A37) -- cycle;
\fill[color4!50] ($(A31)! 0.5!(A37)$) -- ($(A36)! 0.5!(A37)$) -- (A37) -- cycle;
\fill[color4!50] ($(A32)! 0.5!(A37)$) -- ($(A38)! 0.5!(A37)$) -- (A37) -- cycle;
\fill[color4!50] ($(A36)! 0.5!(A37)$) -- ($(A41)! 0.5!(A37)$) -- (A37) -- cycle;
\fill[color4!50] ($(A38)! 0.5!(A37)$) -- ($(A42)! 0.5!(A37)$) -- (A37) -- cycle;
\fill[color4!50] ($(A41)! 0.5!(A37)$) -- ($(A42)! 0.5!(A37)$) -- (A37) -- cycle;
\colorlet{color5}{mycolor3}
\fill[color5!50] ($(A40)! 0.5!(A44)$) -- ($(A41)! 0.5!(A44)$) -- (A44) -- cycle;
\fill[color5!50] ($(A40)! 0.5!(A44)$) -- ($(A43)! 0.5!(A44)$) -- (A44) -- cycle;
\fill[color5!50] ($(A41)! 0.5!(A44)$) -- ($(A45)! 0.5!(A44)$) -- (A44) -- cycle;
\fill[color5!50] ($(A43)! 0.5!(A44)$) -- ($(A46)! 0.5!(A44)$) -- (A44) -- cycle;
\fill[color5!50] ($(A45)! 0.5!(A44)$) -- ($(A47)! 0.5!(A44)$) -- (A44) -- cycle;
\fill[color5!50] ($(A46)! 0.5!(A44)$) -- ($(A47)! 0.5!(A44)$) -- (A44) -- cycle;
\colorlet{color6}{mycolor3}
\fill[color6!50] ($(A30)! 0.5!(A51)$) -- ($(A31)! 0.5!(A51)$) -- (A51) -- cycle;
\fill[color6!50] ($(A30)! 0.5!(A51)$) -- ($(A52)! 0.5!(A51)$) -- (A51) -- cycle;
\fill[color6!50] ($(A31)! 0.5!(A51)$) -- ($(A50)! 0.5!(A51)$) -- (A51) -- cycle;
\fill[color6!50] ($(A50)! 0.5!(A51)$) -- ($(A61)! 0.5!(A51)$) -- (A51) -- cycle;
\fill[color6!50] ($(A52)! 0.5!(A51)$) -- ($(A62)! 0.5!(A51)$) -- (A51) -- cycle;
\fill[color6!50] ($(A61)! 0.5!(A51)$) -- ($(A62)! 0.5!(A51)$) -- (A51) -- cycle;
\colorlet{color7}{mycolor3}
\fill[color7!50] ($(A28)! 0.5!(A53)$) -- ($(A29)! 0.5!(A53)$) -- (A53) -- cycle;
\fill[color7!50] ($(A28)! 0.5!(A53)$) -- ($(A54)! 0.5!(A53)$) -- (A53) -- cycle;
\fill[color7!50] ($(A29)! 0.5!(A53)$) -- ($(A52)! 0.5!(A53)$) -- (A53) -- cycle;
\fill[color7!50] ($(A52)! 0.5!(A53)$) -- ($(A63)! 0.5!(A53)$) -- (A53) -- cycle;
\fill[color7!50] ($(A54)! 0.5!(A53)$) -- ($(A64)! 0.5!(A53)$) -- (A53) -- cycle;
\fill[color7!50] ($(A63)! 0.5!(A53)$) -- ($(A64)! 0.5!(A53)$) -- (A53) -- cycle;
\colorlet{color8}{mycolor3}
\fill[color8!50] ($(A55)! 0.5!(A65)$) -- ($(A56)! 0.5!(A65)$) -- (A65) -- cycle;
\fill[color8!50] ($(A55)! 0.5!(A65)$) -- ($(A64)! 0.5!(A65)$) -- (A65) -- cycle;
\fill[color8!50] ($(A56)! 0.5!(A65)$) -- ($(A66)! 0.5!(A65)$) -- (A65) -- cycle;
\fill[color8!50] ($(A64)! 0.5!(A65)$) -- ($(A72)! 0.5!(A65)$) -- (A65) -- cycle;
\fill[color8!50] ($(A66)! 0.5!(A65)$) -- ($(A73)! 0.5!(A65)$) -- (A65) -- cycle;
\fill[color8!50] ($(A72)! 0.5!(A65)$) -- ($(A73)! 0.5!(A65)$) -- (A65) -- cycle;
\colorlet{color9}{mycolor3}
\fill[color9!50] ($(A57)! 0.5!(A67)$) -- ($(A58)! 0.5!(A67)$) -- (A67) -- cycle;
\fill[color9!50] ($(A57)! 0.5!(A67)$) -- ($(A66)! 0.5!(A67)$) -- (A67) -- cycle;
\fill[color9!50] ($(A58)! 0.5!(A67)$) -- ($(A68)! 0.5!(A67)$) -- (A67) -- cycle;
\fill[color9!50] ($(A66)! 0.5!(A67)$) -- ($(A74)! 0.5!(A67)$) -- (A67) -- cycle;
\fill[color9!50] ($(A68)! 0.5!(A67)$) -- ($(A75)! 0.5!(A67)$) -- (A67) -- cycle;
\fill[color9!50] ($(A74)! 0.5!(A67)$) -- ($(A75)! 0.5!(A67)$) -- (A67) -- cycle;
\colorlet{color10}{mycolor3}
\fill[color10!50] ($(A62)! 0.5!(A71)$) -- ($(A63)! 0.5!(A71)$) -- (A71) -- cycle;
\fill[color10!50] ($(A62)! 0.5!(A71)$) -- ($(A70)! 0.5!(A71)$) -- (A71) -- cycle;
\fill[color10!50] ($(A63)! 0.5!(A71)$) -- ($(A72)! 0.5!(A71)$) -- (A71) -- cycle;
\fill[color10!50] ($(A70)! 0.5!(A71)$) -- ($(A77)! 0.5!(A71)$) -- (A71) -- cycle;
\fill[color10!50] ($(A72)! 0.5!(A71)$) -- ($(A78)! 0.5!(A71)$) -- (A71) -- cycle;
\fill[color10!50] ($(A77)! 0.5!(A71)$) -- ($(A78)! 0.5!(A71)$) -- (A71) -- cycle;
\colorlet{color11}{mycolor3}
\fill[color11!50] ($(A73)! 0.5!(A79)$) -- ($(A74)! 0.5!(A79)$) -- (A79) -- cycle;
\fill[color11!50] ($(A73)! 0.5!(A79)$) -- ($(A78)! 0.5!(A79)$) -- (A79) -- cycle;
\fill[color11!50] ($(A74)! 0.5!(A79)$) -- ($(A80)! 0.5!(A79)$) -- (A79) -- cycle;
\fill[color11!50] ($(A78)! 0.5!(A79)$) -- ($(A82)! 0.5!(A79)$) -- (A79) -- cycle;
\fill[color11!50] ($(A80)! 0.5!(A79)$) -- ($(A83)! 0.5!(A79)$) -- (A79) -- cycle;
\fill[color11!50] ($(A82)! 0.5!(A79)$) -- ($(A83)! 0.5!(A79)$) -- (A79) -- cycle;
\draw[color11, thick] ($(A73)! 0.5!(A79)$) -- ($(A74)! 0.5!(A79)$);
\draw[color11, thick] ($(A73)! 0.5!(A79)$) -- ($(A78)! 0.5!(A79)$);
\draw[color11, thick] ($(A74)! 0.5!(A79)$) -- ($(A80)! 0.5!(A79)$);
\draw[color11, thick] ($(A78)! 0.5!(A79)$) -- ($(A82)! 0.5!(A79)$);
\draw[color11, thick] ($(A80)! 0.5!(A79)$) -- ($(A83)! 0.5!(A79)$);
\draw[color11, thick] ($(A82)! 0.5!(A79)$) -- ($(A83)! 0.5!(A79)$);
\draw[color10, thick] ($(A62)! 0.5!(A71)$) -- ($(A63)! 0.5!(A71)$);
\draw[color10, thick] ($(A62)! 0.5!(A71)$) -- ($(A70)! 0.5!(A71)$);
\draw[color10, thick] ($(A63)! 0.5!(A71)$) -- ($(A72)! 0.5!(A71)$);
\draw[color10, thick] ($(A70)! 0.5!(A71)$) -- ($(A77)! 0.5!(A71)$);
\draw[color10, thick] ($(A72)! 0.5!(A71)$) -- ($(A78)! 0.5!(A71)$);
\draw[color10, thick] ($(A77)! 0.5!(A71)$) -- ($(A78)! 0.5!(A71)$);
\draw[color9, thick] ($(A57)! 0.5!(A67)$) -- ($(A58)! 0.5!(A67)$);
\draw[color9, thick] ($(A57)! 0.5!(A67)$) -- ($(A66)! 0.5!(A67)$);
\draw[color9, thick] ($(A58)! 0.5!(A67)$) -- ($(A68)! 0.5!(A67)$);
\draw[color9, thick] ($(A66)! 0.5!(A67)$) -- ($(A74)! 0.5!(A67)$);
\draw[color9, thick] ($(A68)! 0.5!(A67)$) -- ($(A75)! 0.5!(A67)$);
\draw[color9, thick] ($(A74)! 0.5!(A67)$) -- ($(A75)! 0.5!(A67)$);
\draw[color8, thick] ($(A55)! 0.5!(A65)$) -- ($(A56)! 0.5!(A65)$);
\draw[color8, thick] ($(A55)! 0.5!(A65)$) -- ($(A64)! 0.5!(A65)$);
\draw[color8, thick] ($(A56)! 0.5!(A65)$) -- ($(A66)! 0.5!(A65)$);
\draw[color8, thick] ($(A64)! 0.5!(A65)$) -- ($(A72)! 0.5!(A65)$);
\draw[color8, thick] ($(A66)! 0.5!(A65)$) -- ($(A73)! 0.5!(A65)$);
\draw[color8, thick] ($(A72)! 0.5!(A65)$) -- ($(A73)! 0.5!(A65)$);
\draw[color7, thick] ($(A28)! 0.5!(A53)$) -- ($(A29)! 0.5!(A53)$);
\draw[color7, thick] ($(A28)! 0.5!(A53)$) -- ($(A54)! 0.5!(A53)$);
\draw[color7, thick] ($(A29)! 0.5!(A53)$) -- ($(A52)! 0.5!(A53)$);
\draw[color7, thick] ($(A52)! 0.5!(A53)$) -- ($(A63)! 0.5!(A53)$);
\draw[color7, thick] ($(A54)! 0.5!(A53)$) -- ($(A64)! 0.5!(A53)$);
\draw[color7, thick] ($(A63)! 0.5!(A53)$) -- ($(A64)! 0.5!(A53)$);
\draw[color6, thick] ($(A30)! 0.5!(A51)$) -- ($(A31)! 0.5!(A51)$);
\draw[color6, thick] ($(A30)! 0.5!(A51)$) -- ($(A52)! 0.5!(A51)$);
\draw[color6, thick] ($(A31)! 0.5!(A51)$) -- ($(A50)! 0.5!(A51)$);
\draw[color6, thick] ($(A50)! 0.5!(A51)$) -- ($(A61)! 0.5!(A51)$);
\draw[color6, thick] ($(A52)! 0.5!(A51)$) -- ($(A62)! 0.5!(A51)$);
\draw[color6, thick] ($(A61)! 0.5!(A51)$) -- ($(A62)! 0.5!(A51)$);
\draw[color5, thick] ($(A40)! 0.5!(A44)$) -- ($(A41)! 0.5!(A44)$);
\draw[color5, thick] ($(A40)! 0.5!(A44)$) -- ($(A43)! 0.5!(A44)$);
\draw[color5, thick] ($(A41)! 0.5!(A44)$) -- ($(A45)! 0.5!(A44)$);
\draw[color5, thick] ($(A43)! 0.5!(A44)$) -- ($(A46)! 0.5!(A44)$);
\draw[color5, thick] ($(A45)! 0.5!(A44)$) -- ($(A47)! 0.5!(A44)$);
\draw[color5, thick] ($(A46)! 0.5!(A44)$) -- ($(A47)! 0.5!(A44)$);
\draw[color4, thick] ($(A31)! 0.5!(A37)$) -- ($(A32)! 0.5!(A37)$);
\draw[color4, thick] ($(A31)! 0.5!(A37)$) -- ($(A36)! 0.5!(A37)$);
\draw[color4, thick] ($(A32)! 0.5!(A37)$) -- ($(A38)! 0.5!(A37)$);
\draw[color4, thick] ($(A36)! 0.5!(A37)$) -- ($(A41)! 0.5!(A37)$);
\draw[color4, thick] ($(A38)! 0.5!(A37)$) -- ($(A42)! 0.5!(A37)$);
\draw[color4, thick] ($(A41)! 0.5!(A37)$) -- ($(A42)! 0.5!(A37)$);
\draw[color3, thick] ($(A29)! 0.5!(A35)$) -- ($(A30)! 0.5!(A35)$);
\draw[color3, thick] ($(A29)! 0.5!(A35)$) -- ($(A34)! 0.5!(A35)$);
\draw[color3, thick] ($(A30)! 0.5!(A35)$) -- ($(A36)! 0.5!(A35)$);
\draw[color3, thick] ($(A34)! 0.5!(A35)$) -- ($(A39)! 0.5!(A35)$);
\draw[color3, thick] ($(A36)! 0.5!(A35)$) -- ($(A40)! 0.5!(A35)$);
\draw[color3, thick] ($(A39)! 0.5!(A35)$) -- ($(A40)! 0.5!(A35)$);
\draw[color1, thick] ($(A1)! 0.5!(A0)$) -- ($(A7)! 0.5!(A0)$);
\draw[color1, thick] ($(A0)! 0.5!(A1)$) -- ($(A54)! 0.5!(A1)$);
\draw[color1, thick] ($(A0)! 0.5!(A7)$) -- ($(A28)! 0.5!(A7)$);
\draw[color1, thick] ($(A1)! 0.5!(A2)$) -- ($(A8)! 0.5!(A2)$);
\draw[color1, thick] ($(A2)! 0.5!(A1)$) -- ($(A55)! 0.5!(A1)$);
\draw[color1, thick] ($(A54)! 0.5!(A1)$) -- ($(A55)! 0.5!(A1)$);
\draw[color1, thick] ($(A3)! 0.5!(A2)$) -- ($(A9)! 0.5!(A2)$);
\draw[color1, thick] ($(A2)! 0.5!(A3)$) -- ($(A56)! 0.5!(A3)$);
\draw[color1, thick] ($(A8)! 0.5!(A2)$) -- ($(A9)! 0.5!(A2)$);
\draw[color1, thick] ($(A3)! 0.5!(A4)$) -- ($(A10)! 0.5!(A4)$);
\draw[color1, thick] ($(A4)! 0.5!(A3)$) -- ($(A57)! 0.5!(A3)$);
\draw[color1, thick] ($(A56)! 0.5!(A3)$) -- ($(A57)! 0.5!(A3)$);
\draw[color1, thick] ($(A5)! 0.5!(A4)$) -- ($(A11)! 0.5!(A4)$);
\draw[color1, thick] ($(A4)! 0.5!(A5)$) -- ($(A58)! 0.5!(A5)$);
\draw[color1, thick] ($(A10)! 0.5!(A4)$) -- ($(A11)! 0.5!(A4)$);
\draw[color1, thick] ($(A5)! 0.5!(A6)$) -- ($(A12)! 0.5!(A6)$);
\draw[color1, thick] ($(A6)! 0.5!(A5)$) -- ($(A59)! 0.5!(A5)$);
\draw[color1, thick] ($(A58)! 0.5!(A5)$) -- ($(A59)! 0.5!(A5)$);
\draw[color1, thick] ($(A7)! 0.5!(A13)$) -- ($(A8)! 0.5!(A13)$);
\draw[color1, thick] ($(A13)! 0.5!(A7)$) -- ($(A34)! 0.5!(A7)$);
\draw[color1, thick] ($(A28)! 0.5!(A7)$) -- ($(A34)! 0.5!(A7)$);
\draw[color1, thick] ($(A8)! 0.5!(A13)$) -- ($(A14)! 0.5!(A13)$);
\draw[color1, thick] ($(A9)! 0.5!(A15)$) -- ($(A10)! 0.5!(A15)$);
\draw[color1, thick] ($(A9)! 0.5!(A15)$) -- ($(A14)! 0.5!(A15)$);
\draw[color1, thick] ($(A10)! 0.5!(A15)$) -- ($(A16)! 0.5!(A15)$);
\draw[color1, thick] ($(A11)! 0.5!(A17)$) -- ($(A12)! 0.5!(A17)$);
\draw[color1, thick] ($(A11)! 0.5!(A17)$) -- ($(A16)! 0.5!(A17)$);
\draw[color1, thick] ($(A14)! 0.5!(A13)$) -- ($(A18)! 0.5!(A13)$);
\draw[color1, thick] ($(A13)! 0.5!(A18)$) -- ($(A39)! 0.5!(A18)$);
\draw[color1, thick] ($(A14)! 0.5!(A15)$) -- ($(A19)! 0.5!(A15)$);
\draw[color1, thick] ($(A16)! 0.5!(A15)$) -- ($(A20)! 0.5!(A15)$);
\draw[color1, thick] ($(A19)! 0.5!(A15)$) -- ($(A20)! 0.5!(A15)$);
\draw[color1, thick] ($(A16)! 0.5!(A17)$) -- ($(A21)! 0.5!(A17)$);
\draw[color1, thick] ($(A18)! 0.5!(A22)$) -- ($(A19)! 0.5!(A22)$);
\draw[color1, thick] ($(A22)! 0.5!(A18)$) -- ($(A43)! 0.5!(A18)$);
\draw[color1, thick] ($(A39)! 0.5!(A18)$) -- ($(A43)! 0.5!(A18)$);
\draw[color1, thick] ($(A19)! 0.5!(A22)$) -- ($(A23)! 0.5!(A22)$);
\draw[color1, thick] ($(A20)! 0.5!(A24)$) -- ($(A21)! 0.5!(A24)$);
\draw[color1, thick] ($(A20)! 0.5!(A24)$) -- ($(A23)! 0.5!(A24)$);
\draw[color1, thick] ($(A23)! 0.5!(A22)$) -- ($(A25)! 0.5!(A22)$);
\draw[color1, thick] ($(A22)! 0.5!(A25)$) -- ($(A46)! 0.5!(A25)$);
\draw[color1, thick] ($(A23)! 0.5!(A24)$) -- ($(A26)! 0.5!(A24)$);
\draw[color1, thick] ($(A25)! 0.5!(A27)$) -- ($(A26)! 0.5!(A27)$);
\draw[color1, thick] ($(A27)! 0.5!(A25)$) -- ($(A48)! 0.5!(A25)$);
\draw[color1, thick] ($(A46)! 0.5!(A25)$) -- ($(A48)! 0.5!(A25)$);
\draw[color1, thick] ($(A32)! 0.5!(A49)$) -- ($(A33)! 0.5!(A49)$);
\draw[color1, thick] ($(A32)! 0.5!(A49)$) -- ($(A50)! 0.5!(A49)$);
\draw[color1, thick] ($(A50)! 0.5!(A49)$) -- ($(A60)! 0.5!(A49)$);
\draw[color1, thick] ($(A60)! 0.5!(A69)$) -- ($(A61)! 0.5!(A69)$);
\draw[color1, thick] ($(A61)! 0.5!(A69)$) -- ($(A70)! 0.5!(A69)$);
\draw[color1, thick] ($(A70)! 0.5!(A69)$) -- ($(A76)! 0.5!(A69)$);
\draw[color1, thick] ($(A76)! 0.5!(A81)$) -- ($(A77)! 0.5!(A81)$);
\draw[color1, thick] ($(A77)! 0.5!(A81)$) -- ($(A82)! 0.5!(A81)$);
\draw[color1, thick] ($(A82)! 0.5!(A81)$) -- ($(A84)! 0.5!(A81)$);
\draw[color2, thick] ($(A9)! 0.5!(A15)$) -- ($(A10)! 0.5!(A15)$);
\draw[color2, thick] ($(A9)! 0.5!(A15)$) -- ($(A14)! 0.5!(A15)$);
\draw[color2, thick] ($(A10)! 0.5!(A15)$) -- ($(A16)! 0.5!(A15)$);
\draw[color2, thick] ($(A14)! 0.5!(A15)$) -- ($(A19)! 0.5!(A15)$);
\draw[color2, thick] ($(A16)! 0.5!(A15)$) -- ($(A20)! 0.5!(A15)$);
\draw[color2, thick] ($(A19)! 0.5!(A15)$) -- ($(A20)! 0.5!(A15)$);
\foreach \a/\b/\c in {0/1/7,0/1/54,0/7/28,0/28/54,1/2/8,1/2/55,1/7/8,1/54/55,2/3/9,2/3/56,2/8/9,2/55/56,3/4/10,3/4/57,3/9/10,3/56/57,4/5/11,4/5/58,4/10/11,4/57/58,5/6/12,5/6/59,5/11/12,5/58/59,7/8/13,7/13/34,7/28/34,8/9/14,8/13/14,9/10/15,9/14/15,10/11/16,10/15/16,11/12/17,11/16/17,13/14/18,13/18/39,13/34/39,14/15/19,14/18/19,15/16/20,15/19/20,16/17/21,16/20/21,18/19/22,18/22/43,18/39/43,19/20/23,19/22/23,20/21/24,20/23/24,22/23/25,22/25/46,22/43/46,23/24/26,23/25/26,25/26/27,25/27/48,25/46/48,28/29/34,28/29/53,28/53/54,29/30/35,29/30/52,29/34/35,29/52/53,30/31/36,30/31/51,30/35/36,30/51/52,31/32/37,31/32/50,31/36/37,31/50/51,32/33/38,32/33/49,32/37/38,32/49/50,34/35/39,35/36/40,35/39/40,36/37/41,36/40/41,37/38/42,37/41/42,39/40/43,40/41/44,40/43/44,41/42/45,41/44/45,43/44/46,44/45/47,44/46/47,46/47/48,49/50/60,50/51/61,50/60/61,51/52/62,51/61/62,52/53/63,52/62/63,53/54/64,53/63/64,54/55/64,55/56/65,55/64/65,56/57/66,56/65/66,57/58/67,57/66/67,58/59/68,58/67/68,60/61/69,61/62/70,61/69/70,62/63/71,62/70/71,63/64/72,63/71/72,64/65/72,65/66/73,65/72/73,66/67/74,66/73/74,67/68/75,67/74/75,69/70/76,70/71/77,70/76/77,71/72/78,71/77/78,72/73/78,73/74/79,73/78/79,74/75/80,74/79/80,76/77/81,77/78/82,77/81/82,78/79/82,79/80/83,79/82/83,81/82/84,82/83/84}{
  \draw[black!30] (A\a) -- (A\b) -- (A\c) -- cycle;
}
\foreach \i in {1,3,5,7,8,9,10,11,12,14,16,18,19,20,21,23,25,26,35,37,44,49,51,53,65,67,69,71,79,81}{
  \fill[myred] (A\i) circle (3pt);
}
\foreach \i in {0,2,4,6,13,15,17,22,24,27,28,29,30,31,32,33,34,36,38,39,40,41,42,43,45,46,47,48,50,52,54,55,56,57,58,59,60,61,62,63,64,66,68,70,72,73,74,75,76,77,78,80,82,83,84}{
  \fill[myblue] (A\i) circle (3pt);
}
\iffalse
\foreach \i in {0,1,...,84}{
  \node[anchor=north] at (A\i) {\tiny{$\i$}};
}
\fi

%% file: tikz/collection_of_splits_1.tikz
\coordinate (A0) at (0, 0);
\coordinate (A1) at (0, 1);
\coordinate (A2) at (0, 2);
\coordinate (A3) at (0, 3);
\coordinate (A4) at (0, 4);
\coordinate (A5) at (0, 5);
\coordinate (A6) at (0, 6);
\coordinate (A7) at (1, 0);
\coordinate (A8) at (1, 1);
\coordinate (A9) at (1, 2);
\coordinate (A10) at (1, 3);
\coordinate (A11) at (1, 4);
\coordinate (A12) at (1, 5);
\coordinate (A13) at (2, 0);
\coordinate (A14) at (2, 1);
\coordinate (A15) at (2, 2);
\coordinate (A16) at (2, 3);
\coordinate (A17) at (2, 4);
\coordinate (A18) at (3, 0);
\coordinate (A19) at (3, 1);
\coordinate (A20) at (3, 2);
\coordinate (A21) at (3, 3);
\coordinate (A22) at (4, 0);
\coordinate (A23) at (4, 1);
\coordinate (A24) at (4, 2);
\coordinate (A25) at (5, 0);
\coordinate (A26) at (5, 1);
\coordinate (A27) at (6, 0);
\coordinate (A28) at (0, -1);
\coordinate (A29) at (0, -2);
\coordinate (A30) at (0, -3);
\coordinate (A31) at (0, -4);
\coordinate (A32) at (0, -5);
\coordinate (A33) at (0, -6);
\coordinate (A34) at (1, -1);
\coordinate (A35) at (1, -2);
\coordinate (A36) at (1, -3);
\coordinate (A37) at (1, -4);
\coordinate (A38) at (1, -5);
\coordinate (A39) at (2, -1);
\coordinate (A40) at (2, -2);
\coordinate (A41) at (2, -3);
\coordinate (A42) at (2, -4);
\coordinate (A43) at (3, -1);
\coordinate (A44) at (3, -2);
\coordinate (A45) at (3, -3);
\coordinate (A46) at (4, -1);
\coordinate (A47) at (4, -2);
\coordinate (A48) at (5, -1);
\coordinate (A49) at (-1, -5);
\coordinate (A50) at (-1, -4);
\coordinate (A51) at (-1, -3);
\coordinate (A52) at (-1, -2);
\coordinate (A53) at (-1, -1);
\coordinate (A54) at (-1, 0);
\coordinate (A55) at (-1, 1);
\coordinate (A56) at (-1, 2);
\coordinate (A57) at (-1, 3);
\coordinate (A58) at (-1, 4);
\coordinate (A59) at (-1, 5);
\coordinate (A60) at (-2, -4);
\coordinate (A61) at (-2, -3);
\coordinate (A62) at (-2, -2);
\coordinate (A63) at (-2, -1);
\coordinate (A64) at (-2, 0);
\coordinate (A65) at (-2, 1);
\coordinate (A66) at (-2, 2);
\coordinate (A67) at (-2, 3);
\coordinate (A68) at (-2, 4);
\coordinate (A69) at (-3, -3);
\coordinate (A70) at (-3, -2);
\coordinate (A71) at (-3, -1);
\coordinate (A72) at (-3, 0);
\coordinate (A73) at (-3, 1);
\coordinate (A74) at (-3, 2);
\coordinate (A75) at (-3, 3);
\coordinate (A76) at (-4, -2);
\coordinate (A77) at (-4, -1);
\coordinate (A78) at (-4, 0);
\coordinate (A79) at (-4, 1);
\coordinate (A80) at (-4, 2);
\coordinate (A81) at (-5, -1);
\coordinate (A82) at (-5, 0);
\coordinate (A83) at (-5, 1);
\coordinate (A84) at (-6, 0);
\colorlet{color1}{mycolor4}
\fill[color1!50] (A1) -- ($(A1)! 0.5!(A0)$) -- ($(A7)! 0.5!(A0)$) -- (A7) -- cycle;
\fill[color1!50] ($(A0)! 0.5!(A1)$) -- ($(A54)! 0.5!(A1)$) -- (A1) -- cycle;
\fill[color1!50] ($(A0)! 0.5!(A7)$) -- ($(A28)! 0.5!(A7)$) -- (A7) -- cycle;
\fill[color1!50] ($(A2)! 0.5!(A1)$) -- ($(A9)! 0.5!(A1)$) -- (A1) -- cycle;
\fill[color1!50] (A1) -- ($(A1)! 0.5!(A2)$) -- ($(A56)! 0.5!(A2)$) -- (A56) -- cycle;
\fill[color1!50] (A1) -- (A7) -- (A8) -- cycle;
\fill[color1!50] (A1) -- ($(A1)! 0.5!(A15)$) -- ($(A8)! 0.5!(A15)$) -- (A8) -- cycle;
\fill[color1!50] (A1) -- ($(A1)! 0.5!(A9)$) -- ($(A21)! 0.5!(A9)$) -- (A21) -- cycle;
\fill[color1!50] (A1) -- ($(A1)! 0.5!(A15)$) -- ($(A21)! 0.5!(A15)$) -- (A21) -- cycle;
\fill[color1!50] ($(A54)! 0.5!(A1)$) -- ($(A55)! 0.5!(A1)$) -- (A1) -- cycle;
\fill[color1!50] ($(A55)! 0.5!(A1)$) -- ($(A66)! 0.5!(A1)$) -- (A1) -- cycle;
\fill[color1!50] (A1) -- ($(A1)! 0.5!(A75)$) -- ($(A56)! 0.5!(A75)$) -- (A56) -- cycle;
\fill[color1!50] ($(A66)! 0.5!(A1)$) -- ($(A75)! 0.5!(A1)$) -- (A1) -- cycle;
\fill[color1!50] ($(A2)! 0.5!(A56)$) -- ($(A3)! 0.5!(A56)$) -- (A56) -- cycle;
\fill[color1!50] ($(A3)! 0.5!(A56)$) -- ($(A4)! 0.5!(A56)$) -- (A56) -- cycle;
\fill[color1!50] ($(A4)! 0.5!(A56)$) -- ($(A5)! 0.5!(A56)$) -- (A56) -- cycle;
\fill[color1!50] ($(A5)! 0.5!(A6)$) -- ($(A12)! 0.5!(A6)$) -- (A6) -- cycle;
\fill[color1!50] (A6) -- ($(A6)! 0.5!(A5)$) -- ($(A59)! 0.5!(A5)$) -- (A59) -- cycle;
\fill[color1!50] (A56) -- ($(A56)! 0.5!(A5)$) -- ($(A57)! 0.5!(A5)$) -- (A57) -- cycle;
\fill[color1!50] (A57) -- ($(A57)! 0.5!(A5)$) -- ($(A58)! 0.5!(A5)$) -- (A58) -- cycle;
\fill[color1!50] (A58) -- ($(A58)! 0.5!(A5)$) -- ($(A59)! 0.5!(A5)$) -- (A59) -- cycle;
\fill[color1!50] (A7) -- ($(A7)! 0.5!(A15)$) -- ($(A8)! 0.5!(A15)$) -- (A8) -- cycle;
\fill[color1!50] (A7) -- ($(A7)! 0.5!(A13)$) -- ($(A14)! 0.5!(A13)$) -- (A14) -- cycle;
\fill[color1!50] ($(A13)! 0.5!(A7)$) -- ($(A39)! 0.5!(A7)$) -- (A7) -- cycle;
\fill[color1!50] (A7) -- ($(A7)! 0.5!(A15)$) -- ($(A14)! 0.5!(A15)$) -- (A14) -- cycle;
\fill[color1!50] ($(A28)! 0.5!(A7)$) -- ($(A34)! 0.5!(A7)$) -- (A7) -- cycle;
\fill[color1!50] ($(A34)! 0.5!(A7)$) -- ($(A40)! 0.5!(A7)$) -- (A7) -- cycle;
\fill[color1!50] ($(A39)! 0.5!(A7)$) -- ($(A40)! 0.5!(A7)$) -- (A7) -- cycle;
\fill[color1!50] ($(A9)! 0.5!(A17)$) -- ($(A10)! 0.5!(A17)$) -- (A17) -- cycle;
\fill[color1!50] (A16) -- ($(A16)! 0.5!(A9)$) -- ($(A17)! 0.5!(A9)$) -- (A17) -- cycle;
\fill[color1!50] (A16) -- ($(A16)! 0.5!(A9)$) -- ($(A21)! 0.5!(A9)$) -- (A21) -- cycle;
\fill[color1!50] ($(A10)! 0.5!(A17)$) -- ($(A11)! 0.5!(A17)$) -- (A17) -- cycle;
\fill[color1!50] ($(A11)! 0.5!(A17)$) -- ($(A12)! 0.5!(A17)$) -- (A17) -- cycle;
\fill[color1!50] (A14) -- ($(A14)! 0.5!(A13)$) -- ($(A18)! 0.5!(A13)$) -- (A18) -- cycle;
\fill[color1!50] ($(A13)! 0.5!(A18)$) -- ($(A39)! 0.5!(A18)$) -- (A18) -- cycle;
\fill[color1!50] (A14) -- ($(A14)! 0.5!(A15)$) -- ($(A19)! 0.5!(A15)$) -- (A19) -- cycle;
\fill[color1!50] (A14) -- (A18) -- (A19) -- cycle;
\fill[color1!50] (A19) -- ($(A19)! 0.5!(A15)$) -- ($(A25)! 0.5!(A15)$) -- (A25) -- cycle;
\fill[color1!50] (A20) -- ($(A20)! 0.5!(A15)$) -- ($(A21)! 0.5!(A15)$) -- (A21) -- cycle;
\fill[color1!50] (A20) -- ($(A20)! 0.5!(A15)$) -- ($(A23)! 0.5!(A15)$) -- (A23) -- cycle;
\fill[color1!50] (A23) -- ($(A23)! 0.5!(A15)$) -- ($(A25)! 0.5!(A15)$) -- (A25) -- cycle;
\fill[color1!50] (A16) -- (A17) -- (A21) -- cycle;
\fill[color1!50] (A18) -- ($(A18)! 0.5!(A22)$) -- ($(A19)! 0.5!(A22)$) -- (A19) -- cycle;
\fill[color1!50] ($(A22)! 0.5!(A18)$) -- ($(A43)! 0.5!(A18)$) -- (A18) -- cycle;
\fill[color1!50] ($(A39)! 0.5!(A18)$) -- ($(A43)! 0.5!(A18)$) -- (A18) -- cycle;
\fill[color1!50] (A19) -- ($(A19)! 0.5!(A22)$) -- ($(A25)! 0.5!(A22)$) -- (A25) -- cycle;
\fill[color1!50] (A20) -- ($(A20)! 0.5!(A24)$) -- ($(A21)! 0.5!(A24)$) -- (A21) -- cycle;
\fill[color1!50] (A20) -- ($(A20)! 0.5!(A24)$) -- ($(A23)! 0.5!(A24)$) -- (A23) -- cycle;
\fill[color1!50] ($(A22)! 0.5!(A25)$) -- ($(A43)! 0.5!(A25)$) -- (A25) -- cycle;
\fill[color1!50] (A23) -- ($(A23)! 0.5!(A24)$) -- ($(A26)! 0.5!(A24)$) -- (A26) -- cycle;
\fill[color1!50] (A23) -- (A25) -- (A26) -- cycle;
\fill[color1!50] (A25) -- ($(A25)! 0.5!(A27)$) -- ($(A26)! 0.5!(A27)$) -- (A26) -- cycle;
\fill[color1!50] ($(A27)! 0.5!(A25)$) -- ($(A48)! 0.5!(A25)$) -- (A25) -- cycle;
\fill[color1!50] ($(A40)! 0.5!(A25)$) -- ($(A43)! 0.5!(A25)$) -- (A25) -- cycle;
\fill[color1!50] ($(A40)! 0.5!(A25)$) -- ($(A46)! 0.5!(A25)$) -- (A25) -- cycle;
\fill[color1!50] ($(A46)! 0.5!(A25)$) -- ($(A48)! 0.5!(A25)$) -- (A25) -- cycle;
\fill[color1!50] ($(A28)! 0.5!(A52)$) -- ($(A29)! 0.5!(A52)$) -- (A52) -- cycle;
\fill[color1!50] (A52) -- ($(A52)! 0.5!(A28)$) -- ($(A69)! 0.5!(A28)$) -- (A69) -- cycle;
\fill[color1!50] ($(A28)! 0.5!(A69)$) -- ($(A62)! 0.5!(A69)$) -- (A69) -- cycle;
\fill[color1!50] ($(A29)! 0.5!(A30)$) -- ($(A35)! 0.5!(A30)$) -- (A30) -- cycle;
\fill[color1!50] (A30) -- ($(A30)! 0.5!(A29)$) -- ($(A52)! 0.5!(A29)$) -- (A52) -- cycle;
\fill[color1!50] ($(A31)! 0.5!(A30)$) -- ($(A35)! 0.5!(A30)$) -- (A30) -- cycle;
\fill[color1!50] (A30) -- ($(A30)! 0.5!(A31)$) -- ($(A52)! 0.5!(A31)$) -- (A52) -- cycle;
\fill[color1!50] ($(A31)! 0.5!(A32)$) -- ($(A35)! 0.5!(A32)$) -- (A32) -- cycle;
\fill[color1!50] (A32) -- ($(A32)! 0.5!(A31)$) -- ($(A52)! 0.5!(A31)$) -- (A52) -- cycle;
\fill[color1!50] (A32) -- (A33) -- (A38) -- cycle;
\fill[color1!50] (A32) -- ($(A32)! 0.5!(A49)$) -- ($(A33)! 0.5!(A49)$) -- (A33) -- cycle;
\fill[color1!50] (A32) -- ($(A32)! 0.5!(A35)$) -- ($(A36)! 0.5!(A35)$) -- (A36) -- cycle;
\fill[color1!50] (A32) -- ($(A32)! 0.5!(A37)$) -- ($(A36)! 0.5!(A37)$) -- (A36) -- cycle;
\fill[color1!50] (A32) -- ($(A32)! 0.5!(A37)$) -- ($(A38)! 0.5!(A37)$) -- (A38) -- cycle;
\fill[color1!50] (A32) -- ($(A32)! 0.5!(A49)$) -- ($(A50)! 0.5!(A49)$) -- (A50) -- cycle;
\fill[color1!50] (A32) -- ($(A32)! 0.5!(A51)$) -- ($(A50)! 0.5!(A51)$) -- (A50) -- cycle;
\fill[color1!50] (A32) -- ($(A32)! 0.5!(A51)$) -- ($(A52)! 0.5!(A51)$) -- (A52) -- cycle;
\fill[color1!50] (A36) -- ($(A36)! 0.5!(A35)$) -- ($(A42)! 0.5!(A35)$) -- (A42) -- cycle;
\fill[color1!50] ($(A35)! 0.5!(A42)$) -- ($(A41)! 0.5!(A42)$) -- (A42) -- cycle;
\fill[color1!50] (A36) -- ($(A36)! 0.5!(A37)$) -- ($(A42)! 0.5!(A37)$) -- (A42) -- cycle;
\fill[color1!50] (A38) -- ($(A38)! 0.5!(A37)$) -- ($(A42)! 0.5!(A37)$) -- (A42) -- cycle;
\fill[color1!50] ($(A41)! 0.5!(A42)$) -- ($(A45)! 0.5!(A42)$) -- (A42) -- cycle;
\fill[color1!50] (A50) -- ($(A50)! 0.5!(A49)$) -- ($(A60)! 0.5!(A49)$) -- (A60) -- cycle;
\fill[color1!50] (A50) -- ($(A50)! 0.5!(A51)$) -- ($(A60)! 0.5!(A51)$) -- (A60) -- cycle;
\fill[color1!50] (A52) -- ($(A52)! 0.5!(A51)$) -- ($(A60)! 0.5!(A51)$) -- (A60) -- cycle;
\fill[color1!50] (A52) -- ($(A52)! 0.5!(A61)$) -- ($(A60)! 0.5!(A61)$) -- (A60) -- cycle;
\fill[color1!50] (A52) -- ($(A52)! 0.5!(A61)$) -- ($(A69)! 0.5!(A61)$) -- (A69) -- cycle;
\fill[color1!50] (A56) -- (A57) -- (A68) -- cycle;
\fill[color1!50] (A56) -- (A67) -- (A68) -- cycle;
\fill[color1!50] (A56) -- ($(A56)! 0.5!(A75)$) -- ($(A67)! 0.5!(A75)$) -- (A67) -- cycle;
\fill[color1!50] (A57) -- (A58) -- (A68) -- cycle;
\fill[color1!50] (A58) -- (A59) -- (A68) -- cycle;
\fill[color1!50] (A60) -- ($(A60)! 0.5!(A61)$) -- ($(A69)! 0.5!(A61)$) -- (A69) -- cycle;
\fill[color1!50] ($(A62)! 0.5!(A69)$) -- ($(A70)! 0.5!(A69)$) -- (A69) -- cycle;
\fill[color1!50] (A67) -- ($(A67)! 0.5!(A75)$) -- ($(A68)! 0.5!(A75)$) -- (A68) -- cycle;
\fill[color1!50] ($(A70)! 0.5!(A69)$) -- ($(A76)! 0.5!(A69)$) -- (A69) -- cycle;
\fill[color1!50] ($(A76)! 0.5!(A81)$) -- ($(A77)! 0.5!(A81)$) -- (A81) -- cycle;
\fill[color1!50] ($(A77)! 0.5!(A81)$) -- ($(A82)! 0.5!(A81)$) -- (A81) -- cycle;
\fill[color1!50] ($(A82)! 0.5!(A81)$) -- ($(A84)! 0.5!(A81)$) -- (A81) -- cycle;
\colorlet{color2}{mycolor1}
\fill[color2!50] (A2) -- ($(A2)! 0.5!(A1)$) -- ($(A9)! 0.5!(A1)$) -- (A9) -- cycle;
\fill[color2!50] ($(A1)! 0.5!(A2)$) -- ($(A56)! 0.5!(A2)$) -- (A2) -- cycle;
\fill[color2!50] ($(A1)! 0.5!(A9)$) -- ($(A21)! 0.5!(A9)$) -- (A9) -- cycle;
\fill[color2!50] (A2) -- (A3) -- (A9) -- cycle;
\fill[color2!50] (A2) -- ($(A2)! 0.5!(A56)$) -- ($(A3)! 0.5!(A56)$) -- (A3) -- cycle;
\fill[color2!50] (A3) -- (A4) -- (A9) -- cycle;
\fill[color2!50] (A3) -- ($(A3)! 0.5!(A56)$) -- ($(A4)! 0.5!(A56)$) -- (A4) -- cycle;
\fill[color2!50] (A4) -- (A5) -- (A9) -- cycle;
\fill[color2!50] (A4) -- ($(A4)! 0.5!(A56)$) -- ($(A5)! 0.5!(A56)$) -- (A5) -- cycle;
\fill[color2!50] (A5) -- ($(A5)! 0.5!(A6)$) -- ($(A12)! 0.5!(A6)$) -- (A12) -- cycle;
\fill[color2!50] ($(A6)! 0.5!(A5)$) -- ($(A59)! 0.5!(A5)$) -- (A5) -- cycle;
\fill[color2!50] (A5) -- (A9) -- (A10) -- cycle;
\fill[color2!50] (A5) -- (A10) -- (A11) -- cycle;
\fill[color2!50] (A5) -- (A11) -- (A12) -- cycle;
\fill[color2!50] ($(A56)! 0.5!(A5)$) -- ($(A57)! 0.5!(A5)$) -- (A5) -- cycle;
\fill[color2!50] ($(A57)! 0.5!(A5)$) -- ($(A58)! 0.5!(A5)$) -- (A5) -- cycle;
\fill[color2!50] ($(A58)! 0.5!(A5)$) -- ($(A59)! 0.5!(A5)$) -- (A5) -- cycle;
\fill[color2!50] (A9) -- ($(A9)! 0.5!(A17)$) -- ($(A10)! 0.5!(A17)$) -- (A10) -- cycle;
\fill[color2!50] ($(A16)! 0.5!(A9)$) -- ($(A17)! 0.5!(A9)$) -- (A9) -- cycle;
\fill[color2!50] ($(A16)! 0.5!(A9)$) -- ($(A21)! 0.5!(A9)$) -- (A9) -- cycle;
\fill[color2!50] (A10) -- ($(A10)! 0.5!(A17)$) -- ($(A11)! 0.5!(A17)$) -- (A11) -- cycle;
\fill[color2!50] (A11) -- ($(A11)! 0.5!(A17)$) -- ($(A12)! 0.5!(A17)$) -- (A12) -- cycle;
\fill[color2!50] ($(A32)! 0.5!(A49)$) -- ($(A33)! 0.5!(A49)$) -- (A49) -- cycle;
\fill[color2!50] ($(A32)! 0.5!(A49)$) -- ($(A50)! 0.5!(A49)$) -- (A49) -- cycle;
\fill[color2!50] ($(A50)! 0.5!(A49)$) -- ($(A60)! 0.5!(A49)$) -- (A49) -- cycle;
\colorlet{color3}{mycolor3}
\fill[color3!50] ($(A1)! 0.5!(A15)$) -- ($(A8)! 0.5!(A15)$) -- (A15) -- cycle;
\fill[color3!50] ($(A1)! 0.5!(A15)$) -- ($(A21)! 0.5!(A15)$) -- (A15) -- cycle;
\fill[color3!50] ($(A7)! 0.5!(A15)$) -- ($(A8)! 0.5!(A15)$) -- (A15) -- cycle;
\fill[color3!50] ($(A7)! 0.5!(A15)$) -- ($(A14)! 0.5!(A15)$) -- (A15) -- cycle;
\fill[color3!50] ($(A14)! 0.5!(A15)$) -- ($(A19)! 0.5!(A15)$) -- (A15) -- cycle;
\fill[color3!50] ($(A19)! 0.5!(A15)$) -- ($(A25)! 0.5!(A15)$) -- (A15) -- cycle;
\fill[color3!50] ($(A20)! 0.5!(A15)$) -- ($(A21)! 0.5!(A15)$) -- (A15) -- cycle;
\fill[color3!50] ($(A20)! 0.5!(A15)$) -- ($(A23)! 0.5!(A15)$) -- (A15) -- cycle;
\fill[color3!50] ($(A23)! 0.5!(A15)$) -- ($(A25)! 0.5!(A15)$) -- (A15) -- cycle;
\colorlet{color4}{mycolor3}
\fill[color4!50] ($(A32)! 0.5!(A37)$) -- ($(A36)! 0.5!(A37)$) -- (A37) -- cycle;
\fill[color4!50] ($(A32)! 0.5!(A37)$) -- ($(A38)! 0.5!(A37)$) -- (A37) -- cycle;
\fill[color4!50] ($(A36)! 0.5!(A37)$) -- ($(A42)! 0.5!(A37)$) -- (A37) -- cycle;
\fill[color4!50] ($(A38)! 0.5!(A37)$) -- ($(A42)! 0.5!(A37)$) -- (A37) -- cycle;
\colorlet{color5}{mycolor3}
\fill[color5!50] ($(A40)! 0.5!(A44)$) -- ($(A45)! 0.5!(A44)$) -- (A44) -- cycle;
\fill[color5!50] ($(A40)! 0.5!(A44)$) -- ($(A46)! 0.5!(A44)$) -- (A44) -- cycle;
\fill[color5!50] ($(A45)! 0.5!(A44)$) -- ($(A47)! 0.5!(A44)$) -- (A44) -- cycle;
\fill[color5!50] ($(A46)! 0.5!(A44)$) -- ($(A47)! 0.5!(A44)$) -- (A44) -- cycle;
\colorlet{color6}{mycolor3}
\fill[color6!50] ($(A32)! 0.5!(A51)$) -- ($(A50)! 0.5!(A51)$) -- (A51) -- cycle;
\fill[color6!50] ($(A32)! 0.5!(A51)$) -- ($(A52)! 0.5!(A51)$) -- (A51) -- cycle;
\fill[color6!50] ($(A50)! 0.5!(A51)$) -- ($(A60)! 0.5!(A51)$) -- (A51) -- cycle;
\fill[color6!50] ($(A52)! 0.5!(A51)$) -- ($(A60)! 0.5!(A51)$) -- (A51) -- cycle;
\colorlet{color7}{mycolor3}
\fill[color7!50] ($(A28)! 0.5!(A53)$) -- ($(A54)! 0.5!(A53)$) -- (A53) -- cycle;
\fill[color7!50] ($(A28)! 0.5!(A53)$) -- ($(A62)! 0.5!(A53)$) -- (A53) -- cycle;
\fill[color7!50] ($(A54)! 0.5!(A53)$) -- ($(A62)! 0.5!(A53)$) -- (A53) -- cycle;
\colorlet{color8}{mycolor1}
\fill[color8!50] ($(A52)! 0.5!(A61)$) -- ($(A60)! 0.5!(A61)$) -- (A61) -- cycle;
\fill[color8!50] ($(A52)! 0.5!(A61)$) -- ($(A69)! 0.5!(A61)$) -- (A61) -- cycle;
\fill[color8!50] ($(A60)! 0.5!(A61)$) -- ($(A69)! 0.5!(A61)$) -- (A61) -- cycle;
\colorlet{color9}{mycolor3}
\fill[color9!50] ($(A54)! 0.5!(A65)$) -- ($(A64)! 0.5!(A65)$) -- (A65) -- cycle;
\fill[color9!50] ($(A54)! 0.5!(A65)$) -- ($(A66)! 0.5!(A65)$) -- (A65) -- cycle;
\fill[color9!50] ($(A64)! 0.5!(A65)$) -- ($(A72)! 0.5!(A65)$) -- (A65) -- cycle;
\fill[color9!50] ($(A66)! 0.5!(A65)$) -- ($(A73)! 0.5!(A65)$) -- (A65) -- cycle;
\fill[color9!50] ($(A72)! 0.5!(A65)$) -- ($(A73)! 0.5!(A65)$) -- (A65) -- cycle;
\colorlet{color10}{mycolor3}
\fill[color10!50] ($(A62)! 0.5!(A71)$) -- ($(A63)! 0.5!(A71)$) -- (A71) -- cycle;
\fill[color10!50] ($(A62)! 0.5!(A71)$) -- ($(A82)! 0.5!(A71)$) -- (A71) -- cycle;
\fill[color10!50] ($(A63)! 0.5!(A71)$) -- ($(A72)! 0.5!(A71)$) -- (A71) -- cycle;
\fill[color10!50] ($(A72)! 0.5!(A71)$) -- ($(A78)! 0.5!(A71)$) -- (A71) -- cycle;
\fill[color10!50] ($(A78)! 0.5!(A71)$) -- ($(A82)! 0.5!(A71)$) -- (A71) -- cycle;
\colorlet{color11}{mycolor3}
\fill[color11!50] ($(A66)! 0.5!(A79)$) -- ($(A74)! 0.5!(A79)$) -- (A79) -- cycle;
\fill[color11!50] ($(A66)! 0.5!(A79)$) -- ($(A82)! 0.5!(A79)$) -- (A79) -- cycle;
\fill[color11!50] ($(A74)! 0.5!(A79)$) -- ($(A80)! 0.5!(A79)$) -- (A79) -- cycle;
\fill[color11!50] ($(A80)! 0.5!(A79)$) -- ($(A83)! 0.5!(A79)$) -- (A79) -- cycle;
\fill[color11!50] ($(A82)! 0.5!(A79)$) -- ($(A83)! 0.5!(A79)$) -- (A79) -- cycle;
\draw[color11, thick] ($(A66)! 0.5!(A79)$) -- ($(A74)! 0.5!(A79)$);
\draw[color11, thick] ($(A66)! 0.5!(A79)$) -- ($(A82)! 0.5!(A79)$);
\draw[color11, thick] ($(A74)! 0.5!(A79)$) -- ($(A80)! 0.5!(A79)$);
\draw[color11, thick] ($(A80)! 0.5!(A79)$) -- ($(A83)! 0.5!(A79)$);
\draw[color11, thick] ($(A82)! 0.5!(A79)$) -- ($(A83)! 0.5!(A79)$);
\draw[color10, thick] ($(A62)! 0.5!(A71)$) -- ($(A63)! 0.5!(A71)$);
\draw[color10, thick] ($(A62)! 0.5!(A71)$) -- ($(A82)! 0.5!(A71)$);
\draw[color10, thick] ($(A63)! 0.5!(A71)$) -- ($(A72)! 0.5!(A71)$);
\draw[color10, thick] ($(A72)! 0.5!(A71)$) -- ($(A78)! 0.5!(A71)$);
\draw[color10, thick] ($(A78)! 0.5!(A71)$) -- ($(A82)! 0.5!(A71)$);
\draw[color9, thick] ($(A54)! 0.5!(A65)$) -- ($(A64)! 0.5!(A65)$);
\draw[color9, thick] ($(A54)! 0.5!(A65)$) -- ($(A66)! 0.5!(A65)$);
\draw[color9, thick] ($(A64)! 0.5!(A65)$) -- ($(A72)! 0.5!(A65)$);
\draw[color9, thick] ($(A66)! 0.5!(A65)$) -- ($(A73)! 0.5!(A65)$);
\draw[color9, thick] ($(A72)! 0.5!(A65)$) -- ($(A73)! 0.5!(A65)$);
\draw[color7, thick] ($(A28)! 0.5!(A53)$) -- ($(A54)! 0.5!(A53)$);
\draw[color7, thick] ($(A28)! 0.5!(A53)$) -- ($(A62)! 0.5!(A53)$);
\draw[color7, thick] ($(A54)! 0.5!(A53)$) -- ($(A62)! 0.5!(A53)$);
\draw[color5, thick] ($(A40)! 0.5!(A44)$) -- ($(A45)! 0.5!(A44)$);
\draw[color5, thick] ($(A40)! 0.5!(A44)$) -- ($(A46)! 0.5!(A44)$);
\draw[color5, thick] ($(A45)! 0.5!(A44)$) -- ($(A47)! 0.5!(A44)$);
\draw[color5, thick] ($(A46)! 0.5!(A44)$) -- ($(A47)! 0.5!(A44)$);
\draw[color1, thick] ($(A1)! 0.5!(A0)$) -- ($(A7)! 0.5!(A0)$);
\draw[color1, thick] ($(A0)! 0.5!(A1)$) -- ($(A54)! 0.5!(A1)$);
\draw[color1, thick] ($(A0)! 0.5!(A7)$) -- ($(A28)! 0.5!(A7)$);
\draw[color1, thick] ($(A2)! 0.5!(A1)$) -- ($(A9)! 0.5!(A1)$);
\draw[color1, thick] ($(A1)! 0.5!(A2)$) -- ($(A56)! 0.5!(A2)$);
\draw[color1, thick] ($(A1)! 0.5!(A15)$) -- ($(A8)! 0.5!(A15)$);
\draw[color1, thick] ($(A1)! 0.5!(A9)$) -- ($(A21)! 0.5!(A9)$);
\draw[color1, thick] ($(A1)! 0.5!(A15)$) -- ($(A21)! 0.5!(A15)$);
\draw[color1, thick] ($(A54)! 0.5!(A1)$) -- ($(A55)! 0.5!(A1)$);
\draw[color1, thick] ($(A55)! 0.5!(A1)$) -- ($(A66)! 0.5!(A1)$);
\draw[color1, thick] ($(A1)! 0.5!(A75)$) -- ($(A56)! 0.5!(A75)$);
\draw[color1, thick] ($(A66)! 0.5!(A1)$) -- ($(A75)! 0.5!(A1)$);
\draw[color1, thick] ($(A2)! 0.5!(A56)$) -- ($(A3)! 0.5!(A56)$);
\draw[color1, thick] ($(A3)! 0.5!(A56)$) -- ($(A4)! 0.5!(A56)$);
\draw[color1, thick] ($(A4)! 0.5!(A56)$) -- ($(A5)! 0.5!(A56)$);
\draw[color1, thick] ($(A5)! 0.5!(A6)$) -- ($(A12)! 0.5!(A6)$);
\draw[color1, thick] ($(A6)! 0.5!(A5)$) -- ($(A59)! 0.5!(A5)$);
\draw[color1, thick] ($(A56)! 0.5!(A5)$) -- ($(A57)! 0.5!(A5)$);
\draw[color1, thick] ($(A57)! 0.5!(A5)$) -- ($(A58)! 0.5!(A5)$);
\draw[color1, thick] ($(A58)! 0.5!(A5)$) -- ($(A59)! 0.5!(A5)$);
\draw[color1, thick] ($(A7)! 0.5!(A15)$) -- ($(A8)! 0.5!(A15)$);
\draw[color1, thick] ($(A7)! 0.5!(A13)$) -- ($(A14)! 0.5!(A13)$);
\draw[color1, thick] ($(A13)! 0.5!(A7)$) -- ($(A39)! 0.5!(A7)$);
\draw[color1, thick] ($(A7)! 0.5!(A15)$) -- ($(A14)! 0.5!(A15)$);
\draw[color1, thick] ($(A28)! 0.5!(A7)$) -- ($(A34)! 0.5!(A7)$);
\draw[color1, thick] ($(A34)! 0.5!(A7)$) -- ($(A40)! 0.5!(A7)$);
\draw[color1, thick] ($(A39)! 0.5!(A7)$) -- ($(A40)! 0.5!(A7)$);
\draw[color1, thick] ($(A9)! 0.5!(A17)$) -- ($(A10)! 0.5!(A17)$);
\draw[color1, thick] ($(A16)! 0.5!(A9)$) -- ($(A17)! 0.5!(A9)$);
\draw[color1, thick] ($(A16)! 0.5!(A9)$) -- ($(A21)! 0.5!(A9)$);
\draw[color1, thick] ($(A10)! 0.5!(A17)$) -- ($(A11)! 0.5!(A17)$);
\draw[color1, thick] ($(A11)! 0.5!(A17)$) -- ($(A12)! 0.5!(A17)$);
\draw[color1, thick] ($(A14)! 0.5!(A13)$) -- ($(A18)! 0.5!(A13)$);
\draw[color1, thick] ($(A13)! 0.5!(A18)$) -- ($(A39)! 0.5!(A18)$);
\draw[color1, thick] ($(A14)! 0.5!(A15)$) -- ($(A19)! 0.5!(A15)$);
\draw[color1, thick] ($(A19)! 0.5!(A15)$) -- ($(A25)! 0.5!(A15)$);
\draw[color1, thick] ($(A20)! 0.5!(A15)$) -- ($(A21)! 0.5!(A15)$);
\draw[color1, thick] ($(A20)! 0.5!(A15)$) -- ($(A23)! 0.5!(A15)$);
\draw[color1, thick] ($(A23)! 0.5!(A15)$) -- ($(A25)! 0.5!(A15)$);
\draw[color1, thick] ($(A18)! 0.5!(A22)$) -- ($(A19)! 0.5!(A22)$);
\draw[color1, thick] ($(A22)! 0.5!(A18)$) -- ($(A43)! 0.5!(A18)$);
\draw[color1, thick] ($(A39)! 0.5!(A18)$) -- ($(A43)! 0.5!(A18)$);
\draw[color1, thick] ($(A19)! 0.5!(A22)$) -- ($(A25)! 0.5!(A22)$);
\draw[color1, thick] ($(A20)! 0.5!(A24)$) -- ($(A21)! 0.5!(A24)$);
\draw[color1, thick] ($(A20)! 0.5!(A24)$) -- ($(A23)! 0.5!(A24)$);
\draw[color1, thick] ($(A22)! 0.5!(A25)$) -- ($(A43)! 0.5!(A25)$);
\draw[color1, thick] ($(A23)! 0.5!(A24)$) -- ($(A26)! 0.5!(A24)$);
\draw[color1, thick] ($(A25)! 0.5!(A27)$) -- ($(A26)! 0.5!(A27)$);
\draw[color1, thick] ($(A27)! 0.5!(A25)$) -- ($(A48)! 0.5!(A25)$);
\draw[color1, thick] ($(A40)! 0.5!(A25)$) -- ($(A43)! 0.5!(A25)$);
\draw[color1, thick] ($(A40)! 0.5!(A25)$) -- ($(A46)! 0.5!(A25)$);
\draw[color1, thick] ($(A46)! 0.5!(A25)$) -- ($(A48)! 0.5!(A25)$);
\draw[color1, thick] ($(A28)! 0.5!(A52)$) -- ($(A29)! 0.5!(A52)$);
\draw[color1, thick] ($(A52)! 0.5!(A28)$) -- ($(A69)! 0.5!(A28)$);
\draw[color1, thick] ($(A28)! 0.5!(A69)$) -- ($(A62)! 0.5!(A69)$);
\draw[color1, thick] ($(A29)! 0.5!(A30)$) -- ($(A35)! 0.5!(A30)$);
\draw[color1, thick] ($(A30)! 0.5!(A29)$) -- ($(A52)! 0.5!(A29)$);
\draw[color1, thick] ($(A31)! 0.5!(A30)$) -- ($(A35)! 0.5!(A30)$);
\draw[color1, thick] ($(A30)! 0.5!(A31)$) -- ($(A52)! 0.5!(A31)$);
\draw[color1, thick] ($(A31)! 0.5!(A32)$) -- ($(A35)! 0.5!(A32)$);
\draw[color1, thick] ($(A32)! 0.5!(A31)$) -- ($(A52)! 0.5!(A31)$);
\draw[color1, thick] ($(A32)! 0.5!(A49)$) -- ($(A33)! 0.5!(A49)$);
\draw[color1, thick] ($(A32)! 0.5!(A35)$) -- ($(A36)! 0.5!(A35)$);
\draw[color1, thick] ($(A32)! 0.5!(A37)$) -- ($(A36)! 0.5!(A37)$);
\draw[color1, thick] ($(A32)! 0.5!(A37)$) -- ($(A38)! 0.5!(A37)$);
\draw[color1, thick] ($(A32)! 0.5!(A49)$) -- ($(A50)! 0.5!(A49)$);
\draw[color1, thick] ($(A32)! 0.5!(A51)$) -- ($(A50)! 0.5!(A51)$);
\draw[color1, thick] ($(A32)! 0.5!(A51)$) -- ($(A52)! 0.5!(A51)$);
\draw[color1, thick] ($(A36)! 0.5!(A35)$) -- ($(A42)! 0.5!(A35)$);
\draw[color1, thick] ($(A35)! 0.5!(A42)$) -- ($(A41)! 0.5!(A42)$);
\draw[color1, thick] ($(A36)! 0.5!(A37)$) -- ($(A42)! 0.5!(A37)$);
\draw[color1, thick] ($(A38)! 0.5!(A37)$) -- ($(A42)! 0.5!(A37)$);
\draw[color1, thick] ($(A41)! 0.5!(A42)$) -- ($(A45)! 0.5!(A42)$);
\draw[color1, thick] ($(A50)! 0.5!(A49)$) -- ($(A60)! 0.5!(A49)$);
\draw[color1, thick] ($(A50)! 0.5!(A51)$) -- ($(A60)! 0.5!(A51)$);
\draw[color1, thick] ($(A52)! 0.5!(A51)$) -- ($(A60)! 0.5!(A51)$);
\draw[color1, thick] ($(A52)! 0.5!(A61)$) -- ($(A60)! 0.5!(A61)$);
\draw[color1, thick] ($(A52)! 0.5!(A61)$) -- ($(A69)! 0.5!(A61)$);
\draw[color1, thick] ($(A56)! 0.5!(A75)$) -- ($(A67)! 0.5!(A75)$);
\draw[color1, thick] ($(A60)! 0.5!(A61)$) -- ($(A69)! 0.5!(A61)$);
\draw[color1, thick] ($(A62)! 0.5!(A69)$) -- ($(A70)! 0.5!(A69)$);
\draw[color1, thick] ($(A67)! 0.5!(A75)$) -- ($(A68)! 0.5!(A75)$);
\draw[color1, thick] ($(A70)! 0.5!(A69)$) -- ($(A76)! 0.5!(A69)$);
\draw[color1, thick] ($(A76)! 0.5!(A81)$) -- ($(A77)! 0.5!(A81)$);
\draw[color1, thick] ($(A77)! 0.5!(A81)$) -- ($(A82)! 0.5!(A81)$);
\draw[color1, thick] ($(A82)! 0.5!(A81)$) -- ($(A84)! 0.5!(A81)$);
\draw[color8, thick] ($(A52)! 0.5!(A61)$) -- ($(A60)! 0.5!(A61)$);
\draw[color8, thick] ($(A52)! 0.5!(A61)$) -- ($(A69)! 0.5!(A61)$);
\draw[color8, thick] ($(A60)! 0.5!(A61)$) -- ($(A69)! 0.5!(A61)$);
\draw[color6, thick] ($(A32)! 0.5!(A51)$) -- ($(A50)! 0.5!(A51)$);
\draw[color6, thick] ($(A32)! 0.5!(A51)$) -- ($(A52)! 0.5!(A51)$);
\draw[color6, thick] ($(A50)! 0.5!(A51)$) -- ($(A60)! 0.5!(A51)$);
\draw[color6, thick] ($(A52)! 0.5!(A51)$) -- ($(A60)! 0.5!(A51)$);
\draw[color4, thick] ($(A32)! 0.5!(A37)$) -- ($(A36)! 0.5!(A37)$);
\draw[color4, thick] ($(A32)! 0.5!(A37)$) -- ($(A38)! 0.5!(A37)$);
\draw[color4, thick] ($(A36)! 0.5!(A37)$) -- ($(A42)! 0.5!(A37)$);
\draw[color4, thick] ($(A38)! 0.5!(A37)$) -- ($(A42)! 0.5!(A37)$);
\draw[color3, thick] ($(A1)! 0.5!(A15)$) -- ($(A8)! 0.5!(A15)$);
\draw[color3, thick] ($(A1)! 0.5!(A15)$) -- ($(A21)! 0.5!(A15)$);
\draw[color3, thick] ($(A7)! 0.5!(A15)$) -- ($(A8)! 0.5!(A15)$);
\draw[color3, thick] ($(A7)! 0.5!(A15)$) -- ($(A14)! 0.5!(A15)$);
\draw[color3, thick] ($(A14)! 0.5!(A15)$) -- ($(A19)! 0.5!(A15)$);
\draw[color3, thick] ($(A19)! 0.5!(A15)$) -- ($(A25)! 0.5!(A15)$);
\draw[color3, thick] ($(A20)! 0.5!(A15)$) -- ($(A21)! 0.5!(A15)$);
\draw[color3, thick] ($(A20)! 0.5!(A15)$) -- ($(A23)! 0.5!(A15)$);
\draw[color3, thick] ($(A23)! 0.5!(A15)$) -- ($(A25)! 0.5!(A15)$);
\draw[color2, thick] ($(A2)! 0.5!(A1)$) -- ($(A9)! 0.5!(A1)$);
\draw[color2, thick] ($(A1)! 0.5!(A2)$) -- ($(A56)! 0.5!(A2)$);
\draw[color2, thick] ($(A1)! 0.5!(A9)$) -- ($(A21)! 0.5!(A9)$);
\draw[color2, thick] ($(A2)! 0.5!(A56)$) -- ($(A3)! 0.5!(A56)$);
\draw[color2, thick] ($(A3)! 0.5!(A56)$) -- ($(A4)! 0.5!(A56)$);
\draw[color2, thick] ($(A4)! 0.5!(A56)$) -- ($(A5)! 0.5!(A56)$);
\draw[color2, thick] ($(A5)! 0.5!(A6)$) -- ($(A12)! 0.5!(A6)$);
\draw[color2, thick] ($(A6)! 0.5!(A5)$) -- ($(A59)! 0.5!(A5)$);
\draw[color2, thick] ($(A56)! 0.5!(A5)$) -- ($(A57)! 0.5!(A5)$);
\draw[color2, thick] ($(A57)! 0.5!(A5)$) -- ($(A58)! 0.5!(A5)$);
\draw[color2, thick] ($(A58)! 0.5!(A5)$) -- ($(A59)! 0.5!(A5)$);
\draw[color2, thick] ($(A9)! 0.5!(A17)$) -- ($(A10)! 0.5!(A17)$);
\draw[color2, thick] ($(A16)! 0.5!(A9)$) -- ($(A17)! 0.5!(A9)$);
\draw[color2, thick] ($(A16)! 0.5!(A9)$) -- ($(A21)! 0.5!(A9)$);
\draw[color2, thick] ($(A10)! 0.5!(A17)$) -- ($(A11)! 0.5!(A17)$);
\draw[color2, thick] ($(A11)! 0.5!(A17)$) -- ($(A12)! 0.5!(A17)$);
\draw[color2, thick] ($(A32)! 0.5!(A49)$) -- ($(A33)! 0.5!(A49)$);
\draw[color2, thick] ($(A32)! 0.5!(A49)$) -- ($(A50)! 0.5!(A49)$);
\draw[color2, thick] ($(A50)! 0.5!(A49)$) -- ($(A60)! 0.5!(A49)$);

\fill[black, opacity=.2] (A7) -- (A15) -- (A25) -- (A40);
\fill[black, opacity=.2] (A54) -- (A62) -- (A82)  -- (A66);

\foreach \a/\b/\c in {0/1/7,0/1/54,0/7/28,0/28/54,1/2/9,1/2/56,1/7/8,1/8/15,1/9/21,1/15/21,1/54/55,1/55/66,1/56/75,1/66/75,2/3/9,2/3/56,3/4/9,3/4/56,4/5/9,4/5/56,5/6/12,5/6/59,5/9/10,5/10/11,5/11/12,5/56/57,5/57/58,5/58/59,7/8/15,7/13/14,7/13/39,7/14/15,7/28/34,7/34/40,7/39/40,9/10/17,9/16/17,9/16/21,10/11/17,11/12/17,13/14/18,13/18/39,14/15/19,14/18/19,15/19/25,15/20/21,15/20/23,15/23/25,16/17/21,18/19/22,18/22/43,18/39/43,19/22/25,20/21/24,20/23/24,22/25/43,23/24/26,23/25/26,25/26/27,25/27/48,25/40/43,25/40/46,25/46/48,28/29/35,28/29/52,28/34/40,28/35/45,28/40/45,28/52/69,28/53/54,28/53/62,28/62/69,29/30/35,29/30/52,30/31/35,30/31/52,31/32/35,31/32/52,32/33/38,32/33/49,32/35/36,32/36/37,32/37/38,32/49/50,32/50/51,32/51/52,35/36/42,35/41/42,35/41/45,36/37/42,37/38/42,39/40/43,40/44/45,40/44/46,41/42/45,44/45/47,44/46/47,46/47/48,49/50/60,50/51/60,51/52/60,52/60/61,52/61/69,53/54/62,54/55/66,54/62/63,54/63/64,54/64/65,54/65/66,56/57/68,56/67/68,56/67/75,57/58/68,58/59/68,60/61/69,62/63/71,62/69/70,62/70/77,62/71/82,62/77/82,63/64/72,63/71/72,64/65/72,65/66/73,65/72/73,66/73/82,66/74/75,66/74/79,66/79/82,67/68/75,69/70/76,70/76/77,71/72/78,71/78/82,72/73/78,73/78/82,74/75/80,74/79/80,76/77/81,77/81/82,79/80/83,79/82/83,81/82/84,82/83/84}{
  \draw[black!30] (A\a) -- (A\b) -- (A\c) -- cycle;
}

\foreach \a/\b in {28/45,28/35,35/32,29/35,35/42,28/69,28/52,52/32,29/52,52/60,1/75,1/56,56/5,2/56,56/68}{
	\draw[black, thick] (A\a) -- (A\b);
}
\foreach \a/\b in {1/9,9/5}{
	\draw[black, loosely dashed, thick] (A\a) -- (A\b);
}
\foreach \a/\b in {2/9,9/17}{
	\draw[black, densely dotted, thick] (A\a) -- (A\b);
}

\foreach \a/\b in {1/21}{
	\draw[black, densely dashdotted, thick] (A\a) -- (A\b);
}

\foreach \i in {1,6,7,8,14,16,17,18,19,20,21,23,25,26,30,32,33,36,38,42,44,50,52,53,56,57,58,59,60,65,67,68,69,71,79,81}{
  \fill[myred] (A\i) circle (3pt);
}
\foreach \i in {0,2,3,4,5,9,10,11,12,13,15,22,24,27,28,29,31,34,35,37,39,40,41,43,45,46,47,48,49,51,54,55,61,62,63,64,66,70,72,73,74,75,76,77,78,80,82,83,84}{
  \fill[myblue] (A\i) circle (3pt);
}

\foreach \i in {-6,-4,-2,0,2,4,6}{
	\foreach \j in {-6,-4,-2,0,2,4,6}{
		\node[] at (\i,\j) {\input{tikz/even_even.tikz}};
	}
	\foreach \j in {-5,-3,-1,1,3,5}{
		\node[] at (\i,\j) {\input{tikz/even_odd.tikz}};
	}
}
\foreach \i in {-5,-3,-1,1,3,5}{
	\foreach \j in {-6,-4,-2,0,2,4,6}{
		\node[rotate=90] at (\i,\j) {\input{tikz/even_odd.tikz}};
	}
	\foreach \j in {-5,-3,-1,1,3,5}{
		\node[] at (\i,\j) {\input{tikz/odd_odd.tikz}};
	}
}

\iffalse
\foreach \i in {0,1,...,84}{
  \node[anchor=north] at (A\i) {\tiny{$\i$}};
}
\fi

%% file: tikz/collection_of_splits_2.tikz
\coordinate (A0) at (0, 0);
\coordinate (A1) at (0, 1);
\coordinate (A2) at (0, 2);
\coordinate (A3) at (0, 3);
\coordinate (A4) at (0, 4);
\coordinate (A5) at (0, 5);
\coordinate (A6) at (0, 6);
\coordinate (A7) at (1, 0);
\coordinate (A8) at (1, 1);
\coordinate (A9) at (1, 2);
\coordinate (A10) at (1, 3);
\coordinate (A11) at (1, 4);
\coordinate (A12) at (1, 5);
\coordinate (A13) at (2, 0);
\coordinate (A14) at (2, 1);
\coordinate (A15) at (2, 2);
\coordinate (A16) at (2, 3);
\coordinate (A17) at (2, 4);
\coordinate (A18) at (3, 0);
\coordinate (A19) at (3, 1);
\coordinate (A20) at (3, 2);
\coordinate (A21) at (3, 3);
\coordinate (A22) at (4, 0);
\coordinate (A23) at (4, 1);
\coordinate (A24) at (4, 2);
\coordinate (A25) at (5, 0);
\coordinate (A26) at (5, 1);
\coordinate (A27) at (6, 0);
\coordinate (A28) at (0, -1);
\coordinate (A29) at (0, -2);
\coordinate (A30) at (0, -3);
\coordinate (A31) at (0, -4);
\coordinate (A32) at (0, -5);
\coordinate (A33) at (0, -6);
\coordinate (A34) at (1, -1);
\coordinate (A35) at (1, -2);
\coordinate (A36) at (1, -3);
\coordinate (A37) at (1, -4);
\coordinate (A38) at (1, -5);
\coordinate (A39) at (2, -1);
\coordinate (A40) at (2, -2);
\coordinate (A41) at (2, -3);
\coordinate (A42) at (2, -4);
\coordinate (A43) at (3, -1);
\coordinate (A44) at (3, -2);
\coordinate (A45) at (3, -3);
\coordinate (A46) at (4, -1);
\coordinate (A47) at (4, -2);
\coordinate (A48) at (5, -1);
\coordinate (A49) at (-1, -5);
\coordinate (A50) at (-1, -4);
\coordinate (A51) at (-1, -3);
\coordinate (A52) at (-1, -2);
\coordinate (A53) at (-1, -1);
\coordinate (A54) at (-1, 0);
\coordinate (A55) at (-1, 1);
\coordinate (A56) at (-1, 2);
\coordinate (A57) at (-1, 3);
\coordinate (A58) at (-1, 4);
\coordinate (A59) at (-1, 5);
\coordinate (A60) at (-2, -4);
\coordinate (A61) at (-2, -3);
\coordinate (A62) at (-2, -2);
\coordinate (A63) at (-2, -1);
\coordinate (A64) at (-2, 0);
\coordinate (A65) at (-2, 1);
\coordinate (A66) at (-2, 2);
\coordinate (A67) at (-2, 3);
\coordinate (A68) at (-2, 4);
\coordinate (A69) at (-3, -3);
\coordinate (A70) at (-3, -2);
\coordinate (A71) at (-3, -1);
\coordinate (A72) at (-3, 0);
\coordinate (A73) at (-3, 1);
\coordinate (A74) at (-3, 2);
\coordinate (A75) at (-3, 3);
\coordinate (A76) at (-4, -2);
\coordinate (A77) at (-4, -1);
\coordinate (A78) at (-4, 0);
\coordinate (A79) at (-4, 1);
\coordinate (A80) at (-4, 2);
\coordinate (A81) at (-5, -1);
\coordinate (A82) at (-5, 0);
\coordinate (A83) at (-5, 1);
\coordinate (A84) at (-6, 0);
\colorlet{color1}{mycolor4}
\fill[color1!50] (A1) -- ($(A1)! 0.5!(A0)$) -- ($(A7)! 0.5!(A0)$) -- (A7) -- cycle;
\fill[color1!50] ($(A0)! 0.5!(A1)$) -- ($(A54)! 0.5!(A1)$) -- (A1) -- cycle;
\fill[color1!50] ($(A0)! 0.5!(A7)$) -- ($(A28)! 0.5!(A7)$) -- (A7) -- cycle;
\fill[color1!50] ($(A2)! 0.5!(A1)$) -- ($(A9)! 0.5!(A1)$) -- (A1) -- cycle;
\fill[color1!50] (A1) -- ($(A1)! 0.5!(A2)$) -- ($(A56)! 0.5!(A2)$) -- (A56) -- cycle;
\fill[color1!50] (A1) -- (A7) -- (A8) -- cycle;
\fill[color1!50] (A1) -- ($(A1)! 0.5!(A15)$) -- ($(A8)! 0.5!(A15)$) -- (A8) -- cycle;
\fill[color1!50] (A1) -- ($(A1)! 0.5!(A9)$) -- ($(A21)! 0.5!(A9)$) -- (A21) -- cycle;
\fill[color1!50] (A1) -- ($(A1)! 0.5!(A15)$) -- ($(A21)! 0.5!(A15)$) -- (A21) -- cycle;
\fill[color1!50] ($(A54)! 0.5!(A1)$) -- ($(A55)! 0.5!(A1)$) -- (A1) -- cycle;
\fill[color1!50] ($(A55)! 0.5!(A1)$) -- ($(A66)! 0.5!(A1)$) -- (A1) -- cycle;
\fill[color1!50] (A1) -- ($(A1)! 0.5!(A75)$) -- ($(A56)! 0.5!(A75)$) -- (A56) -- cycle;
\fill[color1!50] ($(A66)! 0.5!(A1)$) -- ($(A75)! 0.5!(A1)$) -- (A1) -- cycle;
\fill[color1!50] ($(A2)! 0.5!(A56)$) -- ($(A3)! 0.5!(A56)$) -- (A56) -- cycle;
\fill[color1!50] ($(A3)! 0.5!(A56)$) -- ($(A4)! 0.5!(A56)$) -- (A56) -- cycle;
\fill[color1!50] ($(A4)! 0.5!(A56)$) -- ($(A5)! 0.5!(A56)$) -- (A56) -- cycle;
\fill[color1!50] ($(A5)! 0.5!(A6)$) -- ($(A12)! 0.5!(A6)$) -- (A6) -- cycle;
\fill[color1!50] (A6) -- ($(A6)! 0.5!(A5)$) -- ($(A59)! 0.5!(A5)$) -- (A59) -- cycle;
\fill[color1!50] (A56) -- ($(A56)! 0.5!(A5)$) -- ($(A57)! 0.5!(A5)$) -- (A57) -- cycle;
\fill[color1!50] (A57) -- ($(A57)! 0.5!(A5)$) -- ($(A58)! 0.5!(A5)$) -- (A58) -- cycle;
\fill[color1!50] (A58) -- ($(A58)! 0.5!(A5)$) -- ($(A59)! 0.5!(A5)$) -- (A59) -- cycle;
\fill[color1!50] (A7) -- ($(A7)! 0.5!(A15)$) -- ($(A8)! 0.5!(A15)$) -- (A8) -- cycle;
\fill[color1!50] ($(A13)! 0.5!(A7)$) -- ($(A14)! 0.5!(A7)$) -- (A7) -- cycle;
\fill[color1!50] (A7) -- ($(A7)! 0.5!(A13)$) -- ($(A39)! 0.5!(A13)$) -- (A39) -- cycle;
\fill[color1!50] ($(A14)! 0.5!(A7)$) -- ($(A15)! 0.5!(A7)$) -- (A7) -- cycle;
\fill[color1!50] ($(A28)! 0.5!(A7)$) -- ($(A34)! 0.5!(A7)$) -- (A7) -- cycle;
\fill[color1!50] ($(A34)! 0.5!(A7)$) -- ($(A40)! 0.5!(A7)$) -- (A7) -- cycle;
\fill[color1!50] (A7) -- ($(A7)! 0.5!(A40)$) -- ($(A39)! 0.5!(A40)$) -- (A39) -- cycle;
\fill[color1!50] ($(A9)! 0.5!(A17)$) -- ($(A10)! 0.5!(A17)$) -- (A17) -- cycle;
\fill[color1!50] (A16) -- ($(A16)! 0.5!(A9)$) -- ($(A17)! 0.5!(A9)$) -- (A17) -- cycle;
\fill[color1!50] (A16) -- ($(A16)! 0.5!(A9)$) -- ($(A21)! 0.5!(A9)$) -- (A21) -- cycle;
\fill[color1!50] ($(A10)! 0.5!(A17)$) -- ($(A11)! 0.5!(A17)$) -- (A17) -- cycle;
\fill[color1!50] ($(A11)! 0.5!(A17)$) -- ($(A12)! 0.5!(A17)$) -- (A17) -- cycle;
\fill[color1!50] ($(A13)! 0.5!(A18)$) -- ($(A14)! 0.5!(A18)$) -- (A18) -- cycle;
\fill[color1!50] (A18) -- ($(A18)! 0.5!(A13)$) -- ($(A39)! 0.5!(A13)$) -- (A39) -- cycle;
\fill[color1!50] ($(A14)! 0.5!(A18)$) -- ($(A19)! 0.5!(A18)$) -- (A18) -- cycle;
\fill[color1!50] ($(A15)! 0.5!(A25)$) -- ($(A19)! 0.5!(A25)$) -- (A25) -- cycle;
\fill[color1!50] (A20) -- ($(A20)! 0.5!(A15)$) -- ($(A21)! 0.5!(A15)$) -- (A21) -- cycle;
\fill[color1!50] (A20) -- ($(A20)! 0.5!(A15)$) -- ($(A23)! 0.5!(A15)$) -- (A23) -- cycle;
\fill[color1!50] (A23) -- ($(A23)! 0.5!(A15)$) -- ($(A25)! 0.5!(A15)$) -- (A25) -- cycle;
\fill[color1!50] (A16) -- (A17) -- (A21) -- cycle;
\fill[color1!50] ($(A19)! 0.5!(A18)$) -- ($(A22)! 0.5!(A18)$) -- (A18) -- cycle;
\fill[color1!50] (A18) -- ($(A18)! 0.5!(A22)$) -- ($(A43)! 0.5!(A22)$) -- (A43) -- cycle;
\fill[color1!50] (A18) -- (A39) -- (A43) -- cycle;
\fill[color1!50] ($(A19)! 0.5!(A25)$) -- ($(A22)! 0.5!(A25)$) -- (A25) -- cycle;
\fill[color1!50] (A20) -- ($(A20)! 0.5!(A24)$) -- ($(A21)! 0.5!(A24)$) -- (A21) -- cycle;
\fill[color1!50] (A20) -- ($(A20)! 0.5!(A24)$) -- ($(A23)! 0.5!(A24)$) -- (A23) -- cycle;
\fill[color1!50] (A25) -- ($(A25)! 0.5!(A22)$) -- ($(A43)! 0.5!(A22)$) -- (A43) -- cycle;
\fill[color1!50] (A23) -- ($(A23)! 0.5!(A24)$) -- ($(A26)! 0.5!(A24)$) -- (A26) -- cycle;
\fill[color1!50] (A23) -- (A25) -- (A26) -- cycle;
\fill[color1!50] (A25) -- ($(A25)! 0.5!(A27)$) -- ($(A26)! 0.5!(A27)$) -- (A26) -- cycle;
\fill[color1!50] ($(A27)! 0.5!(A25)$) -- ($(A48)! 0.5!(A25)$) -- (A25) -- cycle;
\fill[color1!50] (A25) -- ($(A25)! 0.5!(A40)$) -- ($(A43)! 0.5!(A40)$) -- (A43) -- cycle;
\fill[color1!50] ($(A40)! 0.5!(A25)$) -- ($(A46)! 0.5!(A25)$) -- (A25) -- cycle;
\fill[color1!50] ($(A46)! 0.5!(A25)$) -- ($(A48)! 0.5!(A25)$) -- (A25) -- cycle;
\fill[color1!50] ($(A28)! 0.5!(A52)$) -- ($(A29)! 0.5!(A52)$) -- (A52) -- cycle;
\fill[color1!50] (A52) -- ($(A52)! 0.5!(A28)$) -- ($(A69)! 0.5!(A28)$) -- (A69) -- cycle;
\fill[color1!50] ($(A28)! 0.5!(A69)$) -- ($(A62)! 0.5!(A69)$) -- (A69) -- cycle;
\fill[color1!50] ($(A29)! 0.5!(A30)$) -- ($(A35)! 0.5!(A30)$) -- (A30) -- cycle;
\fill[color1!50] (A30) -- ($(A30)! 0.5!(A29)$) -- ($(A52)! 0.5!(A29)$) -- (A52) -- cycle;
\fill[color1!50] ($(A31)! 0.5!(A30)$) -- ($(A35)! 0.5!(A30)$) -- (A30) -- cycle;
\fill[color1!50] (A30) -- ($(A30)! 0.5!(A31)$) -- ($(A52)! 0.5!(A31)$) -- (A52) -- cycle;
\fill[color1!50] ($(A31)! 0.5!(A32)$) -- ($(A35)! 0.5!(A32)$) -- (A32) -- cycle;
\fill[color1!50] (A32) -- ($(A32)! 0.5!(A31)$) -- ($(A52)! 0.5!(A31)$) -- (A52) -- cycle;
\fill[color1!50] (A32) -- (A33) -- (A38) -- cycle;
\fill[color1!50] (A32) -- ($(A32)! 0.5!(A49)$) -- ($(A33)! 0.5!(A49)$) -- (A33) -- cycle;
\fill[color1!50] (A32) -- ($(A32)! 0.5!(A35)$) -- ($(A36)! 0.5!(A35)$) -- (A36) -- cycle;
\fill[color1!50] (A32) -- ($(A32)! 0.5!(A37)$) -- ($(A36)! 0.5!(A37)$) -- (A36) -- cycle;
\fill[color1!50] (A32) -- ($(A32)! 0.5!(A37)$) -- ($(A38)! 0.5!(A37)$) -- (A38) -- cycle;
\fill[color1!50] (A32) -- ($(A32)! 0.5!(A49)$) -- ($(A50)! 0.5!(A49)$) -- (A50) -- cycle;
\fill[color1!50] (A32) -- ($(A32)! 0.5!(A51)$) -- ($(A50)! 0.5!(A51)$) -- (A50) -- cycle;
\fill[color1!50] (A32) -- ($(A32)! 0.5!(A51)$) -- ($(A52)! 0.5!(A51)$) -- (A52) -- cycle;
\fill[color1!50] (A36) -- ($(A36)! 0.5!(A35)$) -- ($(A42)! 0.5!(A35)$) -- (A42) -- cycle;
\fill[color1!50] ($(A35)! 0.5!(A42)$) -- ($(A41)! 0.5!(A42)$) -- (A42) -- cycle;
\fill[color1!50] (A36) -- ($(A36)! 0.5!(A37)$) -- ($(A42)! 0.5!(A37)$) -- (A42) -- cycle;
\fill[color1!50] (A38) -- ($(A38)! 0.5!(A37)$) -- ($(A42)! 0.5!(A37)$) -- (A42) -- cycle;
\fill[color1!50] (A39) -- ($(A39)! 0.5!(A40)$) -- ($(A43)! 0.5!(A40)$) -- (A43) -- cycle;
\fill[color1!50] ($(A41)! 0.5!(A42)$) -- ($(A45)! 0.5!(A42)$) -- (A42) -- cycle;
\fill[color1!50] (A50) -- ($(A50)! 0.5!(A49)$) -- ($(A60)! 0.5!(A49)$) -- (A60) -- cycle;
\fill[color1!50] (A50) -- ($(A50)! 0.5!(A51)$) -- ($(A60)! 0.5!(A51)$) -- (A60) -- cycle;
\fill[color1!50] (A52) -- ($(A52)! 0.5!(A51)$) -- ($(A60)! 0.5!(A51)$) -- (A60) -- cycle;
\fill[color1!50] (A52) -- ($(A52)! 0.5!(A61)$) -- ($(A60)! 0.5!(A61)$) -- (A60) -- cycle;
\fill[color1!50] (A52) -- ($(A52)! 0.5!(A61)$) -- ($(A69)! 0.5!(A61)$) -- (A69) -- cycle;
\fill[color1!50] (A56) -- (A57) -- (A68) -- cycle;
\fill[color1!50] (A56) -- (A67) -- (A68) -- cycle;
\fill[color1!50] (A56) -- ($(A56)! 0.5!(A75)$) -- ($(A67)! 0.5!(A75)$) -- (A67) -- cycle;
\fill[color1!50] (A57) -- (A58) -- (A68) -- cycle;
\fill[color1!50] (A58) -- (A59) -- (A68) -- cycle;
\fill[color1!50] (A60) -- ($(A60)! 0.5!(A61)$) -- ($(A69)! 0.5!(A61)$) -- (A69) -- cycle;
\fill[color1!50] ($(A62)! 0.5!(A69)$) -- ($(A70)! 0.5!(A69)$) -- (A69) -- cycle;
\fill[color1!50] (A67) -- ($(A67)! 0.5!(A75)$) -- ($(A68)! 0.5!(A75)$) -- (A68) -- cycle;
\fill[color1!50] ($(A70)! 0.5!(A69)$) -- ($(A76)! 0.5!(A69)$) -- (A69) -- cycle;
\fill[color1!50] ($(A76)! 0.5!(A81)$) -- ($(A77)! 0.5!(A81)$) -- (A81) -- cycle;
\fill[color1!50] ($(A77)! 0.5!(A81)$) -- ($(A82)! 0.5!(A81)$) -- (A81) -- cycle;
\fill[color1!50] ($(A82)! 0.5!(A81)$) -- ($(A84)! 0.5!(A81)$) -- (A81) -- cycle;
\colorlet{color2}{mycolor1}
\fill[color2!50] (A2) -- ($(A2)! 0.5!(A1)$) -- ($(A9)! 0.5!(A1)$) -- (A9) -- cycle;
\fill[color2!50] ($(A1)! 0.5!(A2)$) -- ($(A56)! 0.5!(A2)$) -- (A2) -- cycle;
\fill[color2!50] ($(A1)! 0.5!(A9)$) -- ($(A21)! 0.5!(A9)$) -- (A9) -- cycle;
\fill[color2!50] (A2) -- (A3) -- (A9) -- cycle;
\fill[color2!50] (A2) -- ($(A2)! 0.5!(A56)$) -- ($(A3)! 0.5!(A56)$) -- (A3) -- cycle;
\fill[color2!50] (A3) -- (A4) -- (A9) -- cycle;
\fill[color2!50] (A3) -- ($(A3)! 0.5!(A56)$) -- ($(A4)! 0.5!(A56)$) -- (A4) -- cycle;
\fill[color2!50] (A4) -- (A5) -- (A9) -- cycle;
\fill[color2!50] (A4) -- ($(A4)! 0.5!(A56)$) -- ($(A5)! 0.5!(A56)$) -- (A5) -- cycle;
\fill[color2!50] (A5) -- ($(A5)! 0.5!(A6)$) -- ($(A12)! 0.5!(A6)$) -- (A12) -- cycle;
\fill[color2!50] ($(A6)! 0.5!(A5)$) -- ($(A59)! 0.5!(A5)$) -- (A5) -- cycle;
\fill[color2!50] (A5) -- (A9) -- (A10) -- cycle;
\fill[color2!50] (A5) -- (A10) -- (A11) -- cycle;
\fill[color2!50] (A5) -- (A11) -- (A12) -- cycle;
\fill[color2!50] ($(A56)! 0.5!(A5)$) -- ($(A57)! 0.5!(A5)$) -- (A5) -- cycle;
\fill[color2!50] ($(A57)! 0.5!(A5)$) -- ($(A58)! 0.5!(A5)$) -- (A5) -- cycle;
\fill[color2!50] ($(A58)! 0.5!(A5)$) -- ($(A59)! 0.5!(A5)$) -- (A5) -- cycle;
\fill[color2!50] (A9) -- ($(A9)! 0.5!(A17)$) -- ($(A10)! 0.5!(A17)$) -- (A10) -- cycle;
\fill[color2!50] ($(A16)! 0.5!(A9)$) -- ($(A17)! 0.5!(A9)$) -- (A9) -- cycle;
\fill[color2!50] ($(A16)! 0.5!(A9)$) -- ($(A21)! 0.5!(A9)$) -- (A9) -- cycle;
\fill[color2!50] (A10) -- ($(A10)! 0.5!(A17)$) -- ($(A11)! 0.5!(A17)$) -- (A11) -- cycle;
\fill[color2!50] (A11) -- ($(A11)! 0.5!(A17)$) -- ($(A12)! 0.5!(A17)$) -- (A12) -- cycle;
\fill[color2!50] ($(A32)! 0.5!(A49)$) -- ($(A33)! 0.5!(A49)$) -- (A49) -- cycle;
\fill[color2!50] ($(A32)! 0.5!(A49)$) -- ($(A50)! 0.5!(A49)$) -- (A49) -- cycle;
\fill[color2!50] ($(A50)! 0.5!(A49)$) -- ($(A60)! 0.5!(A49)$) -- (A49) -- cycle;
\colorlet{color3}{mycolor3}
\fill[color3!50] ($(A1)! 0.5!(A15)$) -- ($(A8)! 0.5!(A15)$) -- (A15) -- cycle;
\fill[color3!50] ($(A1)! 0.5!(A15)$) -- ($(A21)! 0.5!(A15)$) -- (A15) -- cycle;
\fill[color3!50] ($(A7)! 0.5!(A15)$) -- ($(A8)! 0.5!(A15)$) -- (A15) -- cycle;
\fill[color3!50] (A13) -- ($(A13)! 0.5!(A7)$) -- ($(A14)! 0.5!(A7)$) -- (A14) -- cycle;
\fill[color3!50] ($(A7)! 0.5!(A13)$) -- ($(A39)! 0.5!(A13)$) -- (A13) -- cycle;
\fill[color3!50] (A14) -- ($(A14)! 0.5!(A7)$) -- ($(A15)! 0.5!(A7)$) -- (A15) -- cycle;
\fill[color3!50] (A13) -- ($(A13)! 0.5!(A18)$) -- ($(A14)! 0.5!(A18)$) -- (A14) -- cycle;
\fill[color3!50] ($(A18)! 0.5!(A13)$) -- ($(A39)! 0.5!(A13)$) -- (A13) -- cycle;
\fill[color3!50] (A14) -- (A15) -- (A19) -- cycle;
\fill[color3!50] (A14) -- ($(A14)! 0.5!(A18)$) -- ($(A19)! 0.5!(A18)$) -- (A19) -- cycle;
\fill[color3!50] (A15) -- ($(A15)! 0.5!(A25)$) -- ($(A19)! 0.5!(A25)$) -- (A19) -- cycle;
\fill[color3!50] ($(A20)! 0.5!(A15)$) -- ($(A21)! 0.5!(A15)$) -- (A15) -- cycle;
\fill[color3!50] ($(A20)! 0.5!(A15)$) -- ($(A23)! 0.5!(A15)$) -- (A15) -- cycle;
\fill[color3!50] ($(A23)! 0.5!(A15)$) -- ($(A25)! 0.5!(A15)$) -- (A15) -- cycle;
\fill[color3!50] (A19) -- ($(A19)! 0.5!(A18)$) -- ($(A22)! 0.5!(A18)$) -- (A22) -- cycle;
\fill[color3!50] ($(A18)! 0.5!(A22)$) -- ($(A43)! 0.5!(A22)$) -- (A22) -- cycle;
\fill[color3!50] (A19) -- ($(A19)! 0.5!(A25)$) -- ($(A22)! 0.5!(A25)$) -- (A22) -- cycle;
\fill[color3!50] ($(A25)! 0.5!(A22)$) -- ($(A43)! 0.5!(A22)$) -- (A22) -- cycle;
\colorlet{color4}{mycolor3}
\fill[color4!50] ($(A32)! 0.5!(A37)$) -- ($(A36)! 0.5!(A37)$) -- (A37) -- cycle;
\fill[color4!50] ($(A32)! 0.5!(A37)$) -- ($(A38)! 0.5!(A37)$) -- (A37) -- cycle;
\fill[color4!50] ($(A36)! 0.5!(A37)$) -- ($(A42)! 0.5!(A37)$) -- (A37) -- cycle;
\fill[color4!50] ($(A38)! 0.5!(A37)$) -- ($(A42)! 0.5!(A37)$) -- (A37) -- cycle;
\colorlet{color5}{mycolor3}
\fill[color5!50] ($(A40)! 0.5!(A44)$) -- ($(A45)! 0.5!(A44)$) -- (A44) -- cycle;
\fill[color5!50] ($(A40)! 0.5!(A44)$) -- ($(A46)! 0.5!(A44)$) -- (A44) -- cycle;
\fill[color5!50] ($(A45)! 0.5!(A44)$) -- ($(A47)! 0.5!(A44)$) -- (A44) -- cycle;
\fill[color5!50] ($(A46)! 0.5!(A44)$) -- ($(A47)! 0.5!(A44)$) -- (A44) -- cycle;
\colorlet{color6}{mycolor3}
\fill[color6!50] ($(A32)! 0.5!(A51)$) -- ($(A50)! 0.5!(A51)$) -- (A51) -- cycle;
\fill[color6!50] ($(A32)! 0.5!(A51)$) -- ($(A52)! 0.5!(A51)$) -- (A51) -- cycle;
\fill[color6!50] ($(A50)! 0.5!(A51)$) -- ($(A60)! 0.5!(A51)$) -- (A51) -- cycle;
\fill[color6!50] ($(A52)! 0.5!(A51)$) -- ($(A60)! 0.5!(A51)$) -- (A51) -- cycle;
\colorlet{color7}{mycolor3}
\fill[color7!50] ($(A28)! 0.5!(A53)$) -- ($(A54)! 0.5!(A53)$) -- (A53) -- cycle;
\fill[color7!50] ($(A28)! 0.5!(A53)$) -- ($(A62)! 0.5!(A53)$) -- (A53) -- cycle;
\fill[color7!50] ($(A54)! 0.5!(A53)$) -- ($(A62)! 0.5!(A53)$) -- (A53) -- cycle;
\colorlet{color8}{mycolor1}
\fill[color8!50] ($(A52)! 0.5!(A61)$) -- ($(A60)! 0.5!(A61)$) -- (A61) -- cycle;
\fill[color8!50] ($(A52)! 0.5!(A61)$) -- ($(A69)! 0.5!(A61)$) -- (A61) -- cycle;
\fill[color8!50] ($(A60)! 0.5!(A61)$) -- ($(A69)! 0.5!(A61)$) -- (A61) -- cycle;
\colorlet{color9}{mycolor2}
\fill[color9!50] ($(A54)! 0.5!(A63)$) -- ($(A62)! 0.5!(A63)$) -- (A63) -- cycle;
\fill[color9!50] ($(A54)! 0.5!(A63)$) -- ($(A64)! 0.5!(A63)$) -- (A63) -- cycle;
\fill[color9!50] ($(A62)! 0.5!(A63)$) -- ($(A71)! 0.5!(A63)$) -- (A63) -- cycle;
\fill[color9!50] ($(A64)! 0.5!(A63)$) -- ($(A72)! 0.5!(A63)$) -- (A63) -- cycle;
\fill[color9!50] ($(A71)! 0.5!(A63)$) -- ($(A72)! 0.5!(A63)$) -- (A63) -- cycle;
\colorlet{color10}{mycolor2}
\fill[color10!50] ($(A65)! 0.5!(A73)$) -- ($(A66)! 0.5!(A73)$) -- (A73) -- cycle;
\fill[color10!50] ($(A65)! 0.5!(A73)$) -- ($(A72)! 0.5!(A73)$) -- (A73) -- cycle;
\fill[color10!50] ($(A66)! 0.5!(A73)$) -- ($(A82)! 0.5!(A73)$) -- (A73) -- cycle;
\fill[color10!50] ($(A72)! 0.5!(A73)$) -- ($(A78)! 0.5!(A73)$) -- (A73) -- cycle;
\fill[color10!50] ($(A78)! 0.5!(A73)$) -- ($(A82)! 0.5!(A73)$) -- (A73) -- cycle;
\colorlet{color11}{mycolor3}
\fill[color11!50] ($(A66)! 0.5!(A79)$) -- ($(A74)! 0.5!(A79)$) -- (A79) -- cycle;
\fill[color11!50] ($(A66)! 0.5!(A79)$) -- ($(A82)! 0.5!(A79)$) -- (A79) -- cycle;
\fill[color11!50] ($(A74)! 0.5!(A79)$) -- ($(A80)! 0.5!(A79)$) -- (A79) -- cycle;
\fill[color11!50] ($(A80)! 0.5!(A79)$) -- ($(A83)! 0.5!(A79)$) -- (A79) -- cycle;
\fill[color11!50] ($(A82)! 0.5!(A79)$) -- ($(A83)! 0.5!(A79)$) -- (A79) -- cycle;
\draw[color11, thick] ($(A66)! 0.5!(A79)$) -- ($(A74)! 0.5!(A79)$);
\draw[color11, thick] ($(A66)! 0.5!(A79)$) -- ($(A82)! 0.5!(A79)$);
\draw[color11, thick] ($(A74)! 0.5!(A79)$) -- ($(A80)! 0.5!(A79)$);
\draw[color11, thick] ($(A80)! 0.5!(A79)$) -- ($(A83)! 0.5!(A79)$);
\draw[color11, thick] ($(A82)! 0.5!(A79)$) -- ($(A83)! 0.5!(A79)$);
\draw[color10, thick] ($(A65)! 0.5!(A73)$) -- ($(A66)! 0.5!(A73)$);
\draw[color10, thick] ($(A65)! 0.5!(A73)$) -- ($(A72)! 0.5!(A73)$);
\draw[color10, thick] ($(A66)! 0.5!(A73)$) -- ($(A82)! 0.5!(A73)$);
\draw[color10, thick] ($(A72)! 0.5!(A73)$) -- ($(A78)! 0.5!(A73)$);
\draw[color10, thick] ($(A78)! 0.5!(A73)$) -- ($(A82)! 0.5!(A73)$);
\draw[color9, thick] ($(A54)! 0.5!(A63)$) -- ($(A62)! 0.5!(A63)$);
\draw[color9, thick] ($(A54)! 0.5!(A63)$) -- ($(A64)! 0.5!(A63)$);
\draw[color9, thick] ($(A62)! 0.5!(A63)$) -- ($(A71)! 0.5!(A63)$);
\draw[color9, thick] ($(A64)! 0.5!(A63)$) -- ($(A72)! 0.5!(A63)$);
\draw[color9, thick] ($(A71)! 0.5!(A63)$) -- ($(A72)! 0.5!(A63)$);
\draw[color7, thick] ($(A28)! 0.5!(A53)$) -- ($(A54)! 0.5!(A53)$);
\draw[color7, thick] ($(A28)! 0.5!(A53)$) -- ($(A62)! 0.5!(A53)$);
\draw[color7, thick] ($(A54)! 0.5!(A53)$) -- ($(A62)! 0.5!(A53)$);
\draw[color5, thick] ($(A40)! 0.5!(A44)$) -- ($(A45)! 0.5!(A44)$);
\draw[color5, thick] ($(A40)! 0.5!(A44)$) -- ($(A46)! 0.5!(A44)$);
\draw[color5, thick] ($(A45)! 0.5!(A44)$) -- ($(A47)! 0.5!(A44)$);
\draw[color5, thick] ($(A46)! 0.5!(A44)$) -- ($(A47)! 0.5!(A44)$);
\draw[color1, thick] ($(A1)! 0.5!(A0)$) -- ($(A7)! 0.5!(A0)$);
\draw[color1, thick] ($(A0)! 0.5!(A1)$) -- ($(A54)! 0.5!(A1)$);
\draw[color1, thick] ($(A0)! 0.5!(A7)$) -- ($(A28)! 0.5!(A7)$);
\draw[color1, thick] ($(A2)! 0.5!(A1)$) -- ($(A9)! 0.5!(A1)$);
\draw[color1, thick] ($(A1)! 0.5!(A2)$) -- ($(A56)! 0.5!(A2)$);
\draw[color1, thick] ($(A1)! 0.5!(A15)$) -- ($(A8)! 0.5!(A15)$);
\draw[color1, thick] ($(A1)! 0.5!(A9)$) -- ($(A21)! 0.5!(A9)$);
\draw[color1, thick] ($(A1)! 0.5!(A15)$) -- ($(A21)! 0.5!(A15)$);
\draw[color1, thick] ($(A54)! 0.5!(A1)$) -- ($(A55)! 0.5!(A1)$);
\draw[color1, thick] ($(A55)! 0.5!(A1)$) -- ($(A66)! 0.5!(A1)$);
\draw[color1, thick] ($(A1)! 0.5!(A75)$) -- ($(A56)! 0.5!(A75)$);
\draw[color1, thick] ($(A66)! 0.5!(A1)$) -- ($(A75)! 0.5!(A1)$);
\draw[color1, thick] ($(A2)! 0.5!(A56)$) -- ($(A3)! 0.5!(A56)$);
\draw[color1, thick] ($(A3)! 0.5!(A56)$) -- ($(A4)! 0.5!(A56)$);
\draw[color1, thick] ($(A4)! 0.5!(A56)$) -- ($(A5)! 0.5!(A56)$);
\draw[color1, thick] ($(A5)! 0.5!(A6)$) -- ($(A12)! 0.5!(A6)$);
\draw[color1, thick] ($(A6)! 0.5!(A5)$) -- ($(A59)! 0.5!(A5)$);
\draw[color1, thick] ($(A56)! 0.5!(A5)$) -- ($(A57)! 0.5!(A5)$);
\draw[color1, thick] ($(A57)! 0.5!(A5)$) -- ($(A58)! 0.5!(A5)$);
\draw[color1, thick] ($(A58)! 0.5!(A5)$) -- ($(A59)! 0.5!(A5)$);
\draw[color1, thick] ($(A7)! 0.5!(A15)$) -- ($(A8)! 0.5!(A15)$);
\draw[color1, thick] ($(A13)! 0.5!(A7)$) -- ($(A14)! 0.5!(A7)$);
\draw[color1, thick] ($(A7)! 0.5!(A13)$) -- ($(A39)! 0.5!(A13)$);
\draw[color1, thick] ($(A14)! 0.5!(A7)$) -- ($(A15)! 0.5!(A7)$);
\draw[color1, thick] ($(A28)! 0.5!(A7)$) -- ($(A34)! 0.5!(A7)$);
\draw[color1, thick] ($(A34)! 0.5!(A7)$) -- ($(A40)! 0.5!(A7)$);
\draw[color1, thick] ($(A7)! 0.5!(A40)$) -- ($(A39)! 0.5!(A40)$);
\draw[color1, thick] ($(A9)! 0.5!(A17)$) -- ($(A10)! 0.5!(A17)$);
\draw[color1, thick] ($(A16)! 0.5!(A9)$) -- ($(A17)! 0.5!(A9)$);
\draw[color1, thick] ($(A16)! 0.5!(A9)$) -- ($(A21)! 0.5!(A9)$);
\draw[color1, thick] ($(A10)! 0.5!(A17)$) -- ($(A11)! 0.5!(A17)$);
\draw[color1, thick] ($(A11)! 0.5!(A17)$) -- ($(A12)! 0.5!(A17)$);
\draw[color1, thick] ($(A13)! 0.5!(A18)$) -- ($(A14)! 0.5!(A18)$);
\draw[color1, thick] ($(A18)! 0.5!(A13)$) -- ($(A39)! 0.5!(A13)$);
\draw[color1, thick] ($(A14)! 0.5!(A18)$) -- ($(A19)! 0.5!(A18)$);
\draw[color1, thick] ($(A15)! 0.5!(A25)$) -- ($(A19)! 0.5!(A25)$);
\draw[color1, thick] ($(A20)! 0.5!(A15)$) -- ($(A21)! 0.5!(A15)$);
\draw[color1, thick] ($(A20)! 0.5!(A15)$) -- ($(A23)! 0.5!(A15)$);
\draw[color1, thick] ($(A23)! 0.5!(A15)$) -- ($(A25)! 0.5!(A15)$);
\draw[color1, thick] ($(A19)! 0.5!(A18)$) -- ($(A22)! 0.5!(A18)$);
\draw[color1, thick] ($(A18)! 0.5!(A22)$) -- ($(A43)! 0.5!(A22)$);
\draw[color1, thick] ($(A19)! 0.5!(A25)$) -- ($(A22)! 0.5!(A25)$);
\draw[color1, thick] ($(A20)! 0.5!(A24)$) -- ($(A21)! 0.5!(A24)$);
\draw[color1, thick] ($(A20)! 0.5!(A24)$) -- ($(A23)! 0.5!(A24)$);
\draw[color1, thick] ($(A25)! 0.5!(A22)$) -- ($(A43)! 0.5!(A22)$);
\draw[color1, thick] ($(A23)! 0.5!(A24)$) -- ($(A26)! 0.5!(A24)$);
\draw[color1, thick] ($(A25)! 0.5!(A27)$) -- ($(A26)! 0.5!(A27)$);
\draw[color1, thick] ($(A27)! 0.5!(A25)$) -- ($(A48)! 0.5!(A25)$);
\draw[color1, thick] ($(A25)! 0.5!(A40)$) -- ($(A43)! 0.5!(A40)$);
\draw[color1, thick] ($(A40)! 0.5!(A25)$) -- ($(A46)! 0.5!(A25)$);
\draw[color1, thick] ($(A46)! 0.5!(A25)$) -- ($(A48)! 0.5!(A25)$);
\draw[color1, thick] ($(A28)! 0.5!(A52)$) -- ($(A29)! 0.5!(A52)$);
\draw[color1, thick] ($(A52)! 0.5!(A28)$) -- ($(A69)! 0.5!(A28)$);
\draw[color1, thick] ($(A28)! 0.5!(A69)$) -- ($(A62)! 0.5!(A69)$);
\draw[color1, thick] ($(A29)! 0.5!(A30)$) -- ($(A35)! 0.5!(A30)$);
\draw[color1, thick] ($(A30)! 0.5!(A29)$) -- ($(A52)! 0.5!(A29)$);
\draw[color1, thick] ($(A31)! 0.5!(A30)$) -- ($(A35)! 0.5!(A30)$);
\draw[color1, thick] ($(A30)! 0.5!(A31)$) -- ($(A52)! 0.5!(A31)$);
\draw[color1, thick] ($(A31)! 0.5!(A32)$) -- ($(A35)! 0.5!(A32)$);
\draw[color1, thick] ($(A32)! 0.5!(A31)$) -- ($(A52)! 0.5!(A31)$);
\draw[color1, thick] ($(A32)! 0.5!(A49)$) -- ($(A33)! 0.5!(A49)$);
\draw[color1, thick] ($(A32)! 0.5!(A35)$) -- ($(A36)! 0.5!(A35)$);
\draw[color1, thick] ($(A32)! 0.5!(A37)$) -- ($(A36)! 0.5!(A37)$);
\draw[color1, thick] ($(A32)! 0.5!(A37)$) -- ($(A38)! 0.5!(A37)$);
\draw[color1, thick] ($(A32)! 0.5!(A49)$) -- ($(A50)! 0.5!(A49)$);
\draw[color1, thick] ($(A32)! 0.5!(A51)$) -- ($(A50)! 0.5!(A51)$);
\draw[color1, thick] ($(A32)! 0.5!(A51)$) -- ($(A52)! 0.5!(A51)$);
\draw[color1, thick] ($(A36)! 0.5!(A35)$) -- ($(A42)! 0.5!(A35)$);
\draw[color1, thick] ($(A35)! 0.5!(A42)$) -- ($(A41)! 0.5!(A42)$);
\draw[color1, thick] ($(A36)! 0.5!(A37)$) -- ($(A42)! 0.5!(A37)$);
\draw[color1, thick] ($(A38)! 0.5!(A37)$) -- ($(A42)! 0.5!(A37)$);
\draw[color1, thick] ($(A39)! 0.5!(A40)$) -- ($(A43)! 0.5!(A40)$);
\draw[color1, thick] ($(A41)! 0.5!(A42)$) -- ($(A45)! 0.5!(A42)$);
\draw[color1, thick] ($(A50)! 0.5!(A49)$) -- ($(A60)! 0.5!(A49)$);
\draw[color1, thick] ($(A50)! 0.5!(A51)$) -- ($(A60)! 0.5!(A51)$);
\draw[color1, thick] ($(A52)! 0.5!(A51)$) -- ($(A60)! 0.5!(A51)$);
\draw[color1, thick] ($(A52)! 0.5!(A61)$) -- ($(A60)! 0.5!(A61)$);
\draw[color1, thick] ($(A52)! 0.5!(A61)$) -- ($(A69)! 0.5!(A61)$);
\draw[color1, thick] ($(A56)! 0.5!(A75)$) -- ($(A67)! 0.5!(A75)$);
\draw[color1, thick] ($(A60)! 0.5!(A61)$) -- ($(A69)! 0.5!(A61)$);
\draw[color1, thick] ($(A62)! 0.5!(A69)$) -- ($(A70)! 0.5!(A69)$);
\draw[color1, thick] ($(A67)! 0.5!(A75)$) -- ($(A68)! 0.5!(A75)$);
\draw[color1, thick] ($(A70)! 0.5!(A69)$) -- ($(A76)! 0.5!(A69)$);
\draw[color1, thick] ($(A76)! 0.5!(A81)$) -- ($(A77)! 0.5!(A81)$);
\draw[color1, thick] ($(A77)! 0.5!(A81)$) -- ($(A82)! 0.5!(A81)$);
\draw[color1, thick] ($(A82)! 0.5!(A81)$) -- ($(A84)! 0.5!(A81)$);
\draw[color8, thick] ($(A52)! 0.5!(A61)$) -- ($(A60)! 0.5!(A61)$);
\draw[color8, thick] ($(A52)! 0.5!(A61)$) -- ($(A69)! 0.5!(A61)$);
\draw[color8, thick] ($(A60)! 0.5!(A61)$) -- ($(A69)! 0.5!(A61)$);
\draw[color6, thick] ($(A32)! 0.5!(A51)$) -- ($(A50)! 0.5!(A51)$);
\draw[color6, thick] ($(A32)! 0.5!(A51)$) -- ($(A52)! 0.5!(A51)$);
\draw[color6, thick] ($(A50)! 0.5!(A51)$) -- ($(A60)! 0.5!(A51)$);
\draw[color6, thick] ($(A52)! 0.5!(A51)$) -- ($(A60)! 0.5!(A51)$);
\draw[color4, thick] ($(A32)! 0.5!(A37)$) -- ($(A36)! 0.5!(A37)$);
\draw[color4, thick] ($(A32)! 0.5!(A37)$) -- ($(A38)! 0.5!(A37)$);
\draw[color4, thick] ($(A36)! 0.5!(A37)$) -- ($(A42)! 0.5!(A37)$);
\draw[color4, thick] ($(A38)! 0.5!(A37)$) -- ($(A42)! 0.5!(A37)$);
\draw[color3, thick] ($(A1)! 0.5!(A15)$) -- ($(A8)! 0.5!(A15)$);
\draw[color3, thick] ($(A1)! 0.5!(A15)$) -- ($(A21)! 0.5!(A15)$);
\draw[color3, thick] ($(A7)! 0.5!(A15)$) -- ($(A8)! 0.5!(A15)$);
\draw[color3, thick] ($(A13)! 0.5!(A7)$) -- ($(A14)! 0.5!(A7)$);
\draw[color3, thick] ($(A7)! 0.5!(A13)$) -- ($(A39)! 0.5!(A13)$);
\draw[color3, thick] ($(A14)! 0.5!(A7)$) -- ($(A15)! 0.5!(A7)$);
\draw[color3, thick] ($(A13)! 0.5!(A18)$) -- ($(A14)! 0.5!(A18)$);
\draw[color3, thick] ($(A18)! 0.5!(A13)$) -- ($(A39)! 0.5!(A13)$);
\draw[color3, thick] ($(A14)! 0.5!(A18)$) -- ($(A19)! 0.5!(A18)$);
\draw[color3, thick] ($(A15)! 0.5!(A25)$) -- ($(A19)! 0.5!(A25)$);
\draw[color3, thick] ($(A20)! 0.5!(A15)$) -- ($(A21)! 0.5!(A15)$);
\draw[color3, thick] ($(A20)! 0.5!(A15)$) -- ($(A23)! 0.5!(A15)$);
\draw[color3, thick] ($(A23)! 0.5!(A15)$) -- ($(A25)! 0.5!(A15)$);
\draw[color3, thick] ($(A19)! 0.5!(A18)$) -- ($(A22)! 0.5!(A18)$);
\draw[color3, thick] ($(A18)! 0.5!(A22)$) -- ($(A43)! 0.5!(A22)$);
\draw[color3, thick] ($(A19)! 0.5!(A25)$) -- ($(A22)! 0.5!(A25)$);
\draw[color3, thick] ($(A25)! 0.5!(A22)$) -- ($(A43)! 0.5!(A22)$);
\draw[color2, thick] ($(A2)! 0.5!(A1)$) -- ($(A9)! 0.5!(A1)$);
\draw[color2, thick] ($(A1)! 0.5!(A2)$) -- ($(A56)! 0.5!(A2)$);
\draw[color2, thick] ($(A1)! 0.5!(A9)$) -- ($(A21)! 0.5!(A9)$);
\draw[color2, thick] ($(A2)! 0.5!(A56)$) -- ($(A3)! 0.5!(A56)$);
\draw[color2, thick] ($(A3)! 0.5!(A56)$) -- ($(A4)! 0.5!(A56)$);
\draw[color2, thick] ($(A4)! 0.5!(A56)$) -- ($(A5)! 0.5!(A56)$);
\draw[color2, thick] ($(A5)! 0.5!(A6)$) -- ($(A12)! 0.5!(A6)$);
\draw[color2, thick] ($(A6)! 0.5!(A5)$) -- ($(A59)! 0.5!(A5)$);
\draw[color2, thick] ($(A56)! 0.5!(A5)$) -- ($(A57)! 0.5!(A5)$);
\draw[color2, thick] ($(A57)! 0.5!(A5)$) -- ($(A58)! 0.5!(A5)$);
\draw[color2, thick] ($(A58)! 0.5!(A5)$) -- ($(A59)! 0.5!(A5)$);
\draw[color2, thick] ($(A9)! 0.5!(A17)$) -- ($(A10)! 0.5!(A17)$);
\draw[color2, thick] ($(A16)! 0.5!(A9)$) -- ($(A17)! 0.5!(A9)$);
\draw[color2, thick] ($(A16)! 0.5!(A9)$) -- ($(A21)! 0.5!(A9)$);
\draw[color2, thick] ($(A10)! 0.5!(A17)$) -- ($(A11)! 0.5!(A17)$);
\draw[color2, thick] ($(A11)! 0.5!(A17)$) -- ($(A12)! 0.5!(A17)$);
\draw[color2, thick] ($(A32)! 0.5!(A49)$) -- ($(A33)! 0.5!(A49)$);
\draw[color2, thick] ($(A32)! 0.5!(A49)$) -- ($(A50)! 0.5!(A49)$);
\draw[color2, thick] ($(A50)! 0.5!(A49)$) -- ($(A60)! 0.5!(A49)$);

\fill[black, opacity=.2] (A7) -- (A15) -- (A25) -- (A40);
\fill[black, opacity=.2] (A54) -- (A62) -- (A82)  -- (A66);

\foreach \a/\b/\c in {0/1/7,0/1/54,0/7/28,0/28/54,1/2/9,1/2/56,1/7/8,1/8/15,1/9/21,1/15/21,1/54/55,1/55/66,1/56/75,1/66/75,2/3/9,2/3/56,3/4/9,3/4/56,4/5/9,4/5/56,5/6/12,5/6/59,5/9/10,5/10/11,5/11/12,5/56/57,5/57/58,5/58/59,7/8/15,7/13/14,7/13/39,7/14/15,7/28/34,7/34/40,7/39/40,9/10/17,9/16/17,9/16/21,10/11/17,11/12/17,13/14/18,13/18/39,14/15/19,14/18/19,15/19/25,15/20/21,15/20/23,15/23/25,16/17/21,18/19/22,18/22/43,18/39/43,19/22/25,20/21/24,20/23/24,22/25/43,23/24/26,23/25/26,25/26/27,25/27/48,25/40/43,25/40/46,25/46/48,28/29/35,28/29/52,28/34/40,28/35/45,28/40/45,28/52/69,28/53/54,28/53/62,28/62/69,29/30/35,29/30/52,30/31/35,30/31/52,31/32/35,31/32/52,32/33/38,32/33/49,32/35/36,32/36/37,32/37/38,32/49/50,32/50/51,32/51/52,35/36/42,35/41/42,35/41/45,36/37/42,37/38/42,39/40/43,40/44/45,40/44/46,41/42/45,44/45/47,44/46/47,46/47/48,49/50/60,50/51/60,51/52/60,52/60/61,52/61/69,53/54/62,54/55/66,54/62/63,54/63/64,54/64/65,54/65/66,56/57/68,56/67/68,56/67/75,57/58/68,58/59/68,60/61/69,62/63/71,62/69/70,62/70/77,62/71/82,62/77/82,63/64/72,63/71/72,64/65/72,65/66/73,65/72/73,66/73/82,66/74/75,66/74/79,66/79/82,67/68/75,69/70/76,70/76/77,71/72/78,71/78/82,72/73/78,73/78/82,74/75/80,74/79/80,76/77/81,77/81/82,79/80/83,79/82/83,81/82/84,82/83/84}{
  \draw[black!30] (A\a) -- (A\b) -- (A\c) -- cycle;
}

\foreach \a/\b in {7/40,25/40,28/45,28/35,35/32,29/35,35/42,54/62,62/82,28/69,28/52,52/32,29/52,52/60,54/66,66/82,1/75,1/56,56/5,2/56,56/68}{
	\draw[black, thick] (A\a) -- (A\b);
}
\foreach \a/\b in {1/9,9/5}{
	\draw[black, loosely dashed, thick] (A\a) -- (A\b);
}
\foreach \a/\b in {2/9,9/17}{
	\draw[black, densely dotted, thick] (A\a) -- (A\b);
}
\foreach \a/\b in {7/15,15/25}{
	\draw[black, densely dashed, thick] (A\a) -- (A\b);
}
\foreach \a/\b in {1/21}{
	\draw[black, densely dashdotted, thick] (A\a) -- (A\b);
}

\foreach \i in {1,6,7,8,16,17,18,20,21,23,25,26,30,32,33,36,38,39,42,43,44,50,52,53,56,57,58,59,60,63,67,68,69,73,79,81}{
  \fill[myred] (A\i) circle (3pt);
}
\foreach \i in {0,2,3,4,5,9,10,11,12,13,14,15,19,22,24,27,28,29,31,34,35,37,40,41,45,46,47,48,49,51,54,55,61,62,64,65,66,70,71,72,74,75,76,77,78,80,82,83,84}{
  \fill[myblue] (A\i) circle (3pt);
}

\foreach \i in {-6,-4,-2,0,2,4,6}{
	\foreach \j in {-6,-4,-2,0,2,4,6}{
		\node[] at (\i,\j) {\input{tikz/even_even.tikz}};
	}
	\foreach \j in {-5,-3,-1,1,3,5}{
		\node[] at (\i,\j) {\input{tikz/even_odd.tikz}};
	}
}
\foreach \i in {-5,-3,-1,1,3,5}{
	\foreach \j in {-6,-4,-2,0,2,4,6}{
		\node[rotate=90] at (\i,\j) {\input{tikz/even_odd.tikz}};
	}
	\foreach \j in {-5,-3,-1,1,3,5}{
		\node[] at (\i,\j) {\input{tikz/odd_odd.tikz}};
	}
}

\iffalse
\foreach \i in {0,1,...,84}{
  \node[anchor=north] at (A\i) {\tiny{$\i$}};
}
\fi

%% file: tikz/simple_split_4.tikz
\coordinate (A0) at (0, 0);
\coordinate (A1) at (0, 1);
\coordinate (A2) at (0, 2);
\coordinate (A3) at (0, 3);
\coordinate (A4) at (0, 4);
\coordinate (A5) at (0, 5);
\coordinate (A6) at (0, 6);
\coordinate (A7) at (0, 7);
\coordinate (A8) at (0, 8);
\coordinate (A9) at (1, 0);
\coordinate (A10) at (1, 1);
\coordinate (A11) at (1, 2);
\coordinate (A12) at (1, 3);
\coordinate (A13) at (1, 4);
\coordinate (A14) at (1, 5);
\coordinate (A15) at (1, 6);
\coordinate (A16) at (1, 7);
\coordinate (A17) at (2, 0);
\coordinate (A18) at (2, 1);
\coordinate (A19) at (2, 2);
\coordinate (A20) at (2, 3);
\coordinate (A21) at (2, 4);
\coordinate (A22) at (2, 5);
\coordinate (A23) at (2, 6);
\coordinate (A24) at (3, 0);
\coordinate (A25) at (3, 1);
\coordinate (A26) at (3, 2);
\coordinate (A27) at (3, 3);
\coordinate (A28) at (3, 4);
\coordinate (A29) at (3, 5);
\coordinate (A30) at (4, 0);
\coordinate (A31) at (4, 1);
\coordinate (A32) at (4, 2);
\coordinate (A33) at (4, 3);
\coordinate (A34) at (4, 4);
\coordinate (A35) at (5, 0);
\coordinate (A36) at (5, 1);
\coordinate (A37) at (5, 2);
\coordinate (A38) at (5, 3);
\coordinate (A39) at (6, 0);
\coordinate (A40) at (6, 1);
\coordinate (A41) at (6, 2);
\coordinate (A42) at (7, 0);
\coordinate (A43) at (7, 1);
\coordinate (A44) at (8, 0);
\coordinate (A45) at (0, -1);
\coordinate (A46) at (0, -2);
\coordinate (A47) at (0, -3);
\coordinate (A48) at (0, -4);
\coordinate (A49) at (0, -5);
\coordinate (A50) at (0, -6);
\coordinate (A51) at (0, -7);
\coordinate (A52) at (0, -8);
\coordinate (A53) at (1, -1);
\coordinate (A54) at (1, -2);
\coordinate (A55) at (1, -3);
\coordinate (A56) at (1, -4);
\coordinate (A57) at (1, -5);
\coordinate (A58) at (1, -6);
\coordinate (A59) at (1, -7);
\coordinate (A60) at (2, -1);
\coordinate (A61) at (2, -2);
\coordinate (A62) at (2, -3);
\coordinate (A63) at (2, -4);
\coordinate (A64) at (2, -5);
\coordinate (A65) at (2, -6);
\coordinate (A66) at (3, -1);
\coordinate (A67) at (3, -2);
\coordinate (A68) at (3, -3);
\coordinate (A69) at (3, -4);
\coordinate (A70) at (3, -5);
\coordinate (A71) at (4, -1);
\coordinate (A72) at (4, -2);
\coordinate (A73) at (4, -3);
\coordinate (A74) at (4, -4);
\coordinate (A75) at (5, -1);
\coordinate (A76) at (5, -2);
\coordinate (A77) at (5, -3);
\coordinate (A78) at (6, -1);
\coordinate (A79) at (6, -2);
\coordinate (A80) at (7, -1);
\coordinate (A81) at (-1, -7);
\coordinate (A82) at (-1, -6);
\coordinate (A83) at (-1, -5);
\coordinate (A84) at (-1, -4);
\coordinate (A85) at (-1, -3);
\coordinate (A86) at (-1, -2);
\coordinate (A87) at (-1, -1);
\coordinate (A88) at (-1, 0);
\coordinate (A89) at (-1, 1);
\coordinate (A90) at (-1, 2);
\coordinate (A91) at (-1, 3);
\coordinate (A92) at (-1, 4);
\coordinate (A93) at (-1, 5);
\coordinate (A94) at (-1, 6);
\coordinate (A95) at (-1, 7);
\coordinate (A96) at (-2, -6);
\coordinate (A97) at (-2, -5);
\coordinate (A98) at (-2, -4);
\coordinate (A99) at (-2, -3);
\coordinate (A100) at (-2, -2);
\coordinate (A101) at (-2, -1);
\coordinate (A102) at (-2, 0);
\coordinate (A103) at (-2, 1);
\coordinate (A104) at (-2, 2);
\coordinate (A105) at (-2, 3);
\coordinate (A106) at (-2, 4);
\coordinate (A107) at (-2, 5);
\coordinate (A108) at (-2, 6);
\coordinate (A109) at (-3, -5);
\coordinate (A110) at (-3, -4);
\coordinate (A111) at (-3, -3);
\coordinate (A112) at (-3, -2);
\coordinate (A113) at (-3, -1);
\coordinate (A114) at (-3, 0);
\coordinate (A115) at (-3, 1);
\coordinate (A116) at (-3, 2);
\coordinate (A117) at (-3, 3);
\coordinate (A118) at (-3, 4);
\coordinate (A119) at (-3, 5);
\coordinate (A120) at (-4, -4);
\coordinate (A121) at (-4, -3);
\coordinate (A122) at (-4, -2);
\coordinate (A123) at (-4, -1);
\coordinate (A124) at (-4, 0);
\coordinate (A125) at (-4, 1);
\coordinate (A126) at (-4, 2);
\coordinate (A127) at (-4, 3);
\coordinate (A128) at (-4, 4);
\coordinate (A129) at (-5, -3);
\coordinate (A130) at (-5, -2);
\coordinate (A131) at (-5, -1);
\coordinate (A132) at (-5, 0);
\coordinate (A133) at (-5, 1);
\coordinate (A134) at (-5, 2);
\coordinate (A135) at (-5, 3);
\coordinate (A136) at (-6, -2);
\coordinate (A137) at (-6, -1);
\coordinate (A138) at (-6, 0);
\coordinate (A139) at (-6, 1);
\coordinate (A140) at (-6, 2);
\coordinate (A141) at (-7, -1);
\coordinate (A142) at (-7, 0);
\coordinate (A143) at (-7, 1);
\coordinate (A144) at (-8, 0);
%\definecolor{color1}{rgb}{0.489408871018358,0.655346115477001,0.731772744643193}
\colorlet{color1}{mycolor4}
\fill[color1!50] (A0) -- (A1) -- (A9) -- cycle;
\fill[color1!50] (A0) -- ($(A0)! 0.5!(A88)$) -- ($(A1)! 0.5!(A88)$) -- (A1) -- cycle;
\fill[color1!50] (A0) -- ($(A0)! 0.5!(A45)$) -- ($(A9)! 0.5!(A45)$) -- (A9) -- cycle;
\fill[color1!50] ($(A45)! 0.5!(A0)$) -- ($(A88)! 0.5!(A0)$) -- (A0) -- cycle;
\fill[color1!50] (A1) -- ($(A1)! 0.5!(A10)$) -- ($(A2)! 0.5!(A10)$) -- (A2) -- cycle;
\fill[color1!50] (A1) -- (A2) -- (A89) -- cycle;
\fill[color1!50] (A1) -- ($(A1)! 0.5!(A10)$) -- ($(A9)! 0.5!(A10)$) -- (A9) -- cycle;
\fill[color1!50] (A1) -- ($(A1)! 0.5!(A88)$) -- ($(A89)! 0.5!(A88)$) -- (A89) -- cycle;
\fill[color1!50] (A2) -- (A3) -- (A11) -- cycle;
\fill[color1!50] (A2) -- ($(A2)! 0.5!(A90)$) -- ($(A3)! 0.5!(A90)$) -- (A3) -- cycle;
\fill[color1!50] (A2) -- ($(A2)! 0.5!(A10)$) -- ($(A11)! 0.5!(A10)$) -- (A11) -- cycle;
\fill[color1!50] (A2) -- ($(A2)! 0.5!(A90)$) -- ($(A89)! 0.5!(A90)$) -- (A89) -- cycle;
\fill[color1!50] (A3) -- ($(A3)! 0.5!(A4)$) -- ($(A12)! 0.5!(A4)$) -- (A12) -- cycle;
\fill[color1!50] ($(A4)! 0.5!(A3)$) -- ($(A91)! 0.5!(A3)$) -- (A3) -- cycle;
\fill[color1!50] (A3) -- ($(A3)! 0.5!(A25)$) -- ($(A11)! 0.5!(A25)$) -- (A11) -- cycle;
\fill[color1!50] (A3) -- ($(A3)! 0.5!(A19)$) -- ($(A12)! 0.5!(A19)$) -- (A12) -- cycle;
\fill[color1!50] (A3) -- ($(A3)! 0.5!(A19)$) -- ($(A35)! 0.5!(A19)$) -- (A35) -- cycle;
\fill[color1!50] (A3) -- ($(A3)! 0.5!(A25)$) -- ($(A35)! 0.5!(A25)$) -- (A35) -- cycle;
\fill[color1!50] (A3) -- ($(A3)! 0.5!(A90)$) -- ($(A115)! 0.5!(A90)$) -- (A115) -- cycle;
\fill[color1!50] ($(A91)! 0.5!(A3)$) -- ($(A104)! 0.5!(A3)$) -- (A3) -- cycle;
\fill[color1!50] ($(A104)! 0.5!(A3)$) -- ($(A132)! 0.5!(A3)$) -- (A3) -- cycle;
\fill[color1!50] (A3) -- ($(A3)! 0.5!(A132)$) -- ($(A115)! 0.5!(A132)$) -- (A115) -- cycle;
\fill[color1!50] (A5) -- ($(A5)! 0.5!(A4)$) -- ($(A13)! 0.5!(A4)$) -- (A13) -- cycle;
\fill[color1!50] ($(A4)! 0.5!(A5)$) -- ($(A92)! 0.5!(A5)$) -- (A5) -- cycle;
\fill[color1!50] (A12) -- ($(A12)! 0.5!(A4)$) -- ($(A13)! 0.5!(A4)$) -- (A13) -- cycle;
\fill[color1!50] (A5) -- ($(A5)! 0.5!(A6)$) -- ($(A14)! 0.5!(A6)$) -- (A14) -- cycle;
\fill[color1!50] ($(A6)! 0.5!(A5)$) -- ($(A93)! 0.5!(A5)$) -- (A5) -- cycle;
\fill[color1!50] (A5) -- (A13) -- (A14) -- cycle;
\fill[color1!50] ($(A92)! 0.5!(A5)$) -- ($(A93)! 0.5!(A5)$) -- (A5) -- cycle;
\fill[color1!50] (A7) -- ($(A7)! 0.5!(A6)$) -- ($(A15)! 0.5!(A6)$) -- (A15) -- cycle;
\fill[color1!50] ($(A6)! 0.5!(A7)$) -- ($(A94)! 0.5!(A7)$) -- (A7) -- cycle;
\fill[color1!50] (A14) -- ($(A14)! 0.5!(A6)$) -- ($(A15)! 0.5!(A6)$) -- (A15) -- cycle;
\fill[color1!50] (A7) -- ($(A7)! 0.5!(A8)$) -- ($(A16)! 0.5!(A8)$) -- (A16) -- cycle;
\fill[color1!50] ($(A8)! 0.5!(A7)$) -- ($(A95)! 0.5!(A7)$) -- (A7) -- cycle;
\fill[color1!50] (A7) -- (A15) -- (A16) -- cycle;
\fill[color1!50] ($(A94)! 0.5!(A7)$) -- ($(A95)! 0.5!(A7)$) -- (A7) -- cycle;
\fill[color1!50] (A9) -- ($(A9)! 0.5!(A10)$) -- ($(A17)! 0.5!(A10)$) -- (A17) -- cycle;
\fill[color1!50] (A9) -- (A17) -- (A53) -- cycle;
\fill[color1!50] (A9) -- ($(A9)! 0.5!(A45)$) -- ($(A53)! 0.5!(A45)$) -- (A53) -- cycle;
\fill[color1!50] (A11) -- ($(A11)! 0.5!(A10)$) -- ($(A18)! 0.5!(A10)$) -- (A18) -- cycle;
\fill[color1!50] (A17) -- ($(A17)! 0.5!(A10)$) -- ($(A18)! 0.5!(A10)$) -- (A18) -- cycle;
\fill[color1!50] (A11) -- ($(A11)! 0.5!(A25)$) -- ($(A18)! 0.5!(A25)$) -- (A18) -- cycle;
\fill[color1!50] (A12) -- (A13) -- (A20) -- cycle;
\fill[color1!50] (A12) -- ($(A12)! 0.5!(A19)$) -- ($(A20)! 0.5!(A19)$) -- (A20) -- cycle;
\fill[color1!50] (A13) -- ($(A13)! 0.5!(A21)$) -- ($(A14)! 0.5!(A21)$) -- (A14) -- cycle;
\fill[color1!50] (A13) -- ($(A13)! 0.5!(A21)$) -- ($(A20)! 0.5!(A21)$) -- (A20) -- cycle;
\fill[color1!50] (A14) -- (A15) -- (A22) -- cycle;
\fill[color1!50] (A14) -- ($(A14)! 0.5!(A21)$) -- ($(A22)! 0.5!(A21)$) -- (A22) -- cycle;
\fill[color1!50] (A15) -- ($(A15)! 0.5!(A23)$) -- ($(A16)! 0.5!(A23)$) -- (A16) -- cycle;
\fill[color1!50] (A15) -- ($(A15)! 0.5!(A23)$) -- ($(A22)! 0.5!(A23)$) -- (A22) -- cycle;
\fill[color1!50] (A17) -- (A18) -- (A24) -- cycle;
\fill[color1!50] (A17) -- ($(A17)! 0.5!(A60)$) -- ($(A24)! 0.5!(A60)$) -- (A24) -- cycle;
\fill[color1!50] (A17) -- ($(A17)! 0.5!(A60)$) -- ($(A53)! 0.5!(A60)$) -- (A53) -- cycle;
\fill[color1!50] (A18) -- ($(A18)! 0.5!(A25)$) -- ($(A24)! 0.5!(A25)$) -- (A24) -- cycle;
\fill[color1!50] (A20) -- ($(A20)! 0.5!(A19)$) -- ($(A26)! 0.5!(A19)$) -- (A26) -- cycle;
\fill[color1!50] (A26) -- ($(A26)! 0.5!(A19)$) -- ($(A31)! 0.5!(A19)$) -- (A31) -- cycle;
\fill[color1!50] (A31) -- ($(A31)! 0.5!(A19)$) -- ($(A35)! 0.5!(A19)$) -- (A35) -- cycle;
\fill[color1!50] (A20) -- ($(A20)! 0.5!(A21)$) -- ($(A27)! 0.5!(A21)$) -- (A27) -- cycle;
\fill[color1!50] (A20) -- (A26) -- (A27) -- cycle;
\fill[color1!50] (A22) -- ($(A22)! 0.5!(A21)$) -- ($(A28)! 0.5!(A21)$) -- (A28) -- cycle;
\fill[color1!50] (A27) -- ($(A27)! 0.5!(A21)$) -- ($(A28)! 0.5!(A21)$) -- (A28) -- cycle;
\fill[color1!50] (A22) -- ($(A22)! 0.5!(A23)$) -- ($(A29)! 0.5!(A23)$) -- (A29) -- cycle;
\fill[color1!50] (A22) -- (A28) -- (A29) -- cycle;
\fill[color1!50] (A24) -- ($(A24)! 0.5!(A25)$) -- ($(A30)! 0.5!(A25)$) -- (A30) -- cycle;
\fill[color1!50] (A24) -- (A30) -- (A66) -- cycle;
\fill[color1!50] (A24) -- ($(A24)! 0.5!(A60)$) -- ($(A66)! 0.5!(A60)$) -- (A66) -- cycle;
\fill[color1!50] (A30) -- ($(A30)! 0.5!(A25)$) -- ($(A35)! 0.5!(A25)$) -- (A35) -- cycle;
\fill[color1!50] (A26) -- ($(A26)! 0.5!(A32)$) -- ($(A27)! 0.5!(A32)$) -- (A27) -- cycle;
\fill[color1!50] (A26) -- ($(A26)! 0.5!(A32)$) -- ($(A31)! 0.5!(A32)$) -- (A31) -- cycle;
\fill[color1!50] (A27) -- (A28) -- (A33) -- cycle;
\fill[color1!50] (A27) -- ($(A27)! 0.5!(A32)$) -- ($(A33)! 0.5!(A32)$) -- (A33) -- cycle;
\fill[color1!50] (A28) -- ($(A28)! 0.5!(A34)$) -- ($(A29)! 0.5!(A34)$) -- (A29) -- cycle;
\fill[color1!50] (A28) -- ($(A28)! 0.5!(A34)$) -- ($(A33)! 0.5!(A34)$) -- (A33) -- cycle;
\fill[color1!50] (A30) -- (A35) -- (A66) -- cycle;
\fill[color1!50] (A31) -- ($(A31)! 0.5!(A32)$) -- ($(A36)! 0.5!(A32)$) -- (A36) -- cycle;
\fill[color1!50] (A31) -- (A35) -- (A36) -- cycle;
\fill[color1!50] (A33) -- ($(A33)! 0.5!(A32)$) -- ($(A37)! 0.5!(A32)$) -- (A37) -- cycle;
\fill[color1!50] (A36) -- ($(A36)! 0.5!(A32)$) -- ($(A37)! 0.5!(A32)$) -- (A37) -- cycle;
\fill[color1!50] (A33) -- ($(A33)! 0.5!(A34)$) -- ($(A38)! 0.5!(A34)$) -- (A38) -- cycle;
\fill[color1!50] (A33) -- (A37) -- (A38) -- cycle;
\fill[color1!50] (A35) -- ($(A35)! 0.5!(A39)$) -- ($(A36)! 0.5!(A39)$) -- (A36) -- cycle;
\fill[color1!50] ($(A39)! 0.5!(A35)$) -- ($(A75)! 0.5!(A35)$) -- (A35) -- cycle;
\fill[color1!50] ($(A47)! 0.5!(A35)$) -- ($(A61)! 0.5!(A35)$) -- (A35) -- cycle;
\fill[color1!50] (A35) -- ($(A35)! 0.5!(A47)$) -- ($(A66)! 0.5!(A47)$) -- (A66) -- cycle;
\fill[color1!50] ($(A61)! 0.5!(A35)$) -- ($(A71)! 0.5!(A35)$) -- (A35) -- cycle;
\fill[color1!50] ($(A71)! 0.5!(A35)$) -- ($(A75)! 0.5!(A35)$) -- (A35) -- cycle;
\fill[color1!50] (A36) -- (A37) -- (A40) -- cycle;
\fill[color1!50] (A36) -- ($(A36)! 0.5!(A39)$) -- ($(A40)! 0.5!(A39)$) -- (A40) -- cycle;
\fill[color1!50] (A37) -- ($(A37)! 0.5!(A41)$) -- ($(A38)! 0.5!(A41)$) -- (A38) -- cycle;
\fill[color1!50] (A37) -- ($(A37)! 0.5!(A41)$) -- ($(A40)! 0.5!(A41)$) -- (A40) -- cycle;
\fill[color1!50] (A40) -- ($(A40)! 0.5!(A39)$) -- ($(A42)! 0.5!(A39)$) -- (A42) -- cycle;
\fill[color1!50] ($(A39)! 0.5!(A42)$) -- ($(A78)! 0.5!(A42)$) -- (A42) -- cycle;
\fill[color1!50] (A40) -- ($(A40)! 0.5!(A41)$) -- ($(A43)! 0.5!(A41)$) -- (A43) -- cycle;
\fill[color1!50] (A40) -- (A42) -- (A43) -- cycle;
\fill[color1!50] (A42) -- ($(A42)! 0.5!(A44)$) -- ($(A43)! 0.5!(A44)$) -- (A43) -- cycle;
\fill[color1!50] ($(A44)! 0.5!(A42)$) -- ($(A80)! 0.5!(A42)$) -- (A42) -- cycle;
\fill[color1!50] ($(A78)! 0.5!(A42)$) -- ($(A80)! 0.5!(A42)$) -- (A42) -- cycle;
\fill[color1!50] (A46) -- ($(A46)! 0.5!(A45)$) -- ($(A53)! 0.5!(A45)$) -- (A53) -- cycle;
\fill[color1!50] ($(A45)! 0.5!(A46)$) -- ($(A87)! 0.5!(A46)$) -- (A46) -- cycle;
\fill[color1!50] (A46) -- ($(A46)! 0.5!(A47)$) -- ($(A54)! 0.5!(A47)$) -- (A54) -- cycle;
\fill[color1!50] ($(A47)! 0.5!(A46)$) -- ($(A86)! 0.5!(A46)$) -- (A46) -- cycle;
\fill[color1!50] (A46) -- (A53) -- (A54) -- cycle;
\fill[color1!50] ($(A86)! 0.5!(A46)$) -- ($(A87)! 0.5!(A46)$) -- (A46) -- cycle;
\fill[color1!50] (A54) -- ($(A54)! 0.5!(A47)$) -- ($(A66)! 0.5!(A47)$) -- (A66) -- cycle;
\fill[color1!50] ($(A51)! 0.5!(A81)$) -- ($(A52)! 0.5!(A81)$) -- (A81) -- cycle;
\fill[color1!50] ($(A51)! 0.5!(A81)$) -- ($(A82)! 0.5!(A81)$) -- (A81) -- cycle;
\fill[color1!50] (A53) -- ($(A53)! 0.5!(A60)$) -- ($(A54)! 0.5!(A60)$) -- (A54) -- cycle;
\fill[color1!50] (A54) -- ($(A54)! 0.5!(A60)$) -- ($(A66)! 0.5!(A60)$) -- (A66) -- cycle;
\fill[color1!50] ($(A82)! 0.5!(A81)$) -- ($(A96)! 0.5!(A81)$) -- (A81) -- cycle;
\fill[color1!50] ($(A87)! 0.5!(A102)$) -- ($(A88)! 0.5!(A102)$) -- (A102) -- cycle;
\fill[color1!50] ($(A87)! 0.5!(A102)$) -- ($(A101)! 0.5!(A102)$) -- (A102) -- cycle;
\fill[color1!50] (A89) -- ($(A89)! 0.5!(A88)$) -- ($(A102)! 0.5!(A88)$) -- (A102) -- cycle;
\fill[color1!50] (A89) -- ($(A89)! 0.5!(A90)$) -- ($(A103)! 0.5!(A90)$) -- (A103) -- cycle;
\fill[color1!50] (A89) -- (A102) -- (A103) -- cycle;
\fill[color1!50] (A103) -- ($(A103)! 0.5!(A90)$) -- ($(A115)! 0.5!(A90)$) -- (A115) -- cycle;
\fill[color1!50] ($(A96)! 0.5!(A109)$) -- ($(A97)! 0.5!(A109)$) -- (A109) -- cycle;
\fill[color1!50] ($(A97)! 0.5!(A109)$) -- ($(A110)! 0.5!(A109)$) -- (A109) -- cycle;
\fill[color1!50] ($(A101)! 0.5!(A102)$) -- ($(A114)! 0.5!(A102)$) -- (A102) -- cycle;
\fill[color1!50] (A102) -- ($(A102)! 0.5!(A114)$) -- ($(A103)! 0.5!(A114)$) -- (A103) -- cycle;
\fill[color1!50] (A103) -- ($(A103)! 0.5!(A114)$) -- ($(A115)! 0.5!(A114)$) -- (A115) -- cycle;
\fill[color1!50] ($(A110)! 0.5!(A109)$) -- ($(A120)! 0.5!(A109)$) -- (A109) -- cycle;
\fill[color1!50] ($(A113)! 0.5!(A124)$) -- ($(A114)! 0.5!(A124)$) -- (A124) -- cycle;
\fill[color1!50] ($(A113)! 0.5!(A124)$) -- ($(A132)! 0.5!(A124)$) -- (A124) -- cycle;
\fill[color1!50] (A115) -- ($(A115)! 0.5!(A114)$) -- ($(A124)! 0.5!(A114)$) -- (A124) -- cycle;
\fill[color1!50] (A115) -- ($(A115)! 0.5!(A132)$) -- ($(A124)! 0.5!(A132)$) -- (A124) -- cycle;
\fill[color1!50] ($(A120)! 0.5!(A129)$) -- ($(A121)! 0.5!(A129)$) -- (A129) -- cycle;
\fill[color1!50] ($(A121)! 0.5!(A129)$) -- ($(A130)! 0.5!(A129)$) -- (A129) -- cycle;
\fill[color1!50] ($(A130)! 0.5!(A129)$) -- ($(A136)! 0.5!(A129)$) -- (A129) -- cycle;
\fill[color1!50] ($(A136)! 0.5!(A141)$) -- ($(A137)! 0.5!(A141)$) -- (A141) -- cycle;
\fill[color1!50] ($(A137)! 0.5!(A141)$) -- ($(A142)! 0.5!(A141)$) -- (A141) -- cycle;
\fill[color1!50] ($(A142)! 0.5!(A141)$) -- ($(A144)! 0.5!(A141)$) -- (A141) -- cycle;
%\definecolor{color2}{rgb}{0.858072647839766,0.685489635426028,0.566646115675435}
\colorlet{color2}{mycolor5}
\fill[color2!50] ($(A1)! 0.5!(A10)$) -- ($(A2)! 0.5!(A10)$) -- (A10) -- cycle;
\fill[color2!50] ($(A1)! 0.5!(A10)$) -- ($(A9)! 0.5!(A10)$) -- (A10) -- cycle;
\fill[color2!50] ($(A2)! 0.5!(A10)$) -- ($(A11)! 0.5!(A10)$) -- (A10) -- cycle;
\fill[color2!50] ($(A9)! 0.5!(A10)$) -- ($(A17)! 0.5!(A10)$) -- (A10) -- cycle;
\fill[color2!50] ($(A11)! 0.5!(A10)$) -- ($(A18)! 0.5!(A10)$) -- (A10) -- cycle;
\fill[color2!50] ($(A17)! 0.5!(A10)$) -- ($(A18)! 0.5!(A10)$) -- (A10) -- cycle;
%\definecolor{color3}{rgb}{0.760481017501623,0.605809316017773,0.611644572875402}
\colorlet{color3}{mycolor3}
\fill[color3!50] ($(A3)! 0.5!(A19)$) -- ($(A12)! 0.5!(A19)$) -- (A19) -- cycle;
\fill[color3!50] ($(A3)! 0.5!(A19)$) -- ($(A35)! 0.5!(A19)$) -- (A19) -- cycle;
\fill[color3!50] ($(A12)! 0.5!(A19)$) -- ($(A20)! 0.5!(A19)$) -- (A19) -- cycle;
\fill[color3!50] ($(A20)! 0.5!(A19)$) -- ($(A26)! 0.5!(A19)$) -- (A19) -- cycle;
\fill[color3!50] ($(A26)! 0.5!(A19)$) -- ($(A31)! 0.5!(A19)$) -- (A19) -- cycle;
\fill[color3!50] ($(A31)! 0.5!(A19)$) -- ($(A35)! 0.5!(A19)$) -- (A19) -- cycle;
%\definecolor{color4}{rgb}{0.217008843600708,0.531125846740666,0.60115988252581}
\colorlet{color4}{mycolor3}
\fill[color4!50] ($(A13)! 0.5!(A21)$) -- ($(A14)! 0.5!(A21)$) -- (A21) -- cycle;
\fill[color4!50] ($(A13)! 0.5!(A21)$) -- ($(A20)! 0.5!(A21)$) -- (A21) -- cycle;
\fill[color4!50] ($(A14)! 0.5!(A21)$) -- ($(A22)! 0.5!(A21)$) -- (A21) -- cycle;
\fill[color4!50] ($(A20)! 0.5!(A21)$) -- ($(A27)! 0.5!(A21)$) -- (A21) -- cycle;
\fill[color4!50] ($(A22)! 0.5!(A21)$) -- ($(A28)! 0.5!(A21)$) -- (A21) -- cycle;
\fill[color4!50] ($(A27)! 0.5!(A21)$) -- ($(A28)! 0.5!(A21)$) -- (A21) -- cycle;
%\definecolor{color5}{rgb}{0.562910234995884,0.517110881032526,0.470373263113858}
\colorlet{color5}{mycolor5}
\fill[color5!50] ($(A3)! 0.5!(A25)$) -- ($(A11)! 0.5!(A25)$) -- (A25) -- cycle;
\fill[color5!50] ($(A3)! 0.5!(A25)$) -- ($(A35)! 0.5!(A25)$) -- (A25) -- cycle;
\fill[color5!50] ($(A11)! 0.5!(A25)$) -- ($(A18)! 0.5!(A25)$) -- (A25) -- cycle;
\fill[color5!50] ($(A18)! 0.5!(A25)$) -- ($(A24)! 0.5!(A25)$) -- (A25) -- cycle;
\fill[color5!50] ($(A24)! 0.5!(A25)$) -- ($(A30)! 0.5!(A25)$) -- (A25) -- cycle;
\fill[color5!50] ($(A30)! 0.5!(A25)$) -- ($(A35)! 0.5!(A25)$) -- (A25) -- cycle;
%\definecolor{color6}{rgb}{0.450185019490017,0.513555248544456,0.519520541941275}
\colorlet{color6}{mycolor3}
\fill[color6!50] ($(A26)! 0.5!(A32)$) -- ($(A27)! 0.5!(A32)$) -- (A32) -- cycle;
\fill[color6!50] ($(A26)! 0.5!(A32)$) -- ($(A31)! 0.5!(A32)$) -- (A32) -- cycle;
\fill[color6!50] ($(A27)! 0.5!(A32)$) -- ($(A33)! 0.5!(A32)$) -- (A32) -- cycle;
\fill[color6!50] ($(A31)! 0.5!(A32)$) -- ($(A36)! 0.5!(A32)$) -- (A32) -- cycle;
\fill[color6!50] ($(A33)! 0.5!(A32)$) -- ($(A37)! 0.5!(A32)$) -- (A32) -- cycle;
\fill[color6!50] ($(A36)! 0.5!(A32)$) -- ($(A37)! 0.5!(A32)$) -- (A32) -- cycle;
%\definecolor{color7}{rgb}{0.747371988358532,0.517621386192285,0.422826115223542}
\colorlet{color7}{mycolor3}
\fill[color7!50] ($(A48)! 0.5!(A56)$) -- ($(A49)! 0.5!(A56)$) -- (A56) -- cycle;
\fill[color7!50] ($(A48)! 0.5!(A56)$) -- ($(A55)! 0.5!(A56)$) -- (A56) -- cycle;
\fill[color7!50] ($(A49)! 0.5!(A56)$) -- ($(A57)! 0.5!(A56)$) -- (A56) -- cycle;
\fill[color7!50] ($(A55)! 0.5!(A56)$) -- ($(A62)! 0.5!(A56)$) -- (A56) -- cycle;
\fill[color7!50] ($(A57)! 0.5!(A56)$) -- ($(A63)! 0.5!(A56)$) -- (A56) -- cycle;
\fill[color7!50] ($(A62)! 0.5!(A56)$) -- ($(A63)! 0.5!(A56)$) -- (A56) -- cycle;
%\definecolor{color8}{rgb}{0.649370816742935,0.476912453709824,0.713929698573435}
\colorlet{color8}{mycolor3}
\fill[color8!50] ($(A50)! 0.5!(A58)$) -- ($(A51)! 0.5!(A58)$) -- (A58) -- cycle;
\fill[color8!50] ($(A50)! 0.5!(A58)$) -- ($(A57)! 0.5!(A58)$) -- (A58) -- cycle;
\fill[color8!50] ($(A51)! 0.5!(A58)$) -- ($(A59)! 0.5!(A58)$) -- (A58) -- cycle;
\fill[color8!50] ($(A57)! 0.5!(A58)$) -- ($(A64)! 0.5!(A58)$) -- (A58) -- cycle;
\fill[color8!50] ($(A59)! 0.5!(A58)$) -- ($(A65)! 0.5!(A58)$) -- (A58) -- cycle;
\fill[color8!50] ($(A64)! 0.5!(A58)$) -- ($(A65)! 0.5!(A58)$) -- (A58) -- cycle;
%\definecolor{color9}{rgb}{0.66261567391136,0.373138227220015,0.712393243318094}
\colorlet{color9}{mycolor5}
\fill[color9!50] ($(A17)! 0.5!(A60)$) -- ($(A24)! 0.5!(A60)$) -- (A60) -- cycle;
\fill[color9!50] ($(A17)! 0.5!(A60)$) -- ($(A53)! 0.5!(A60)$) -- (A60) -- cycle;
\fill[color9!50] ($(A24)! 0.5!(A60)$) -- ($(A66)! 0.5!(A60)$) -- (A60) -- cycle;
\fill[color9!50] ($(A53)! 0.5!(A60)$) -- ($(A54)! 0.5!(A60)$) -- (A60) -- cycle;
\fill[color9!50] ($(A54)! 0.5!(A60)$) -- ($(A66)! 0.5!(A60)$) -- (A60) -- cycle;
%\definecolor{color10}{rgb}{0.724724563859856,0.615991217463074,0.64726433050794}
\colorlet{color10}{mycolor3}
\fill[color10!50] ($(A61)! 0.5!(A67)$) -- ($(A62)! 0.5!(A67)$) -- (A67) -- cycle;
\fill[color10!50] ($(A61)! 0.5!(A67)$) -- ($(A71)! 0.5!(A67)$) -- (A67) -- cycle;
\fill[color10!50] ($(A62)! 0.5!(A67)$) -- ($(A68)! 0.5!(A67)$) -- (A67) -- cycle;
\fill[color10!50] ($(A68)! 0.5!(A67)$) -- ($(A72)! 0.5!(A67)$) -- (A67) -- cycle;
\fill[color10!50] ($(A71)! 0.5!(A67)$) -- ($(A72)! 0.5!(A67)$) -- (A67) -- cycle;
%\definecolor{color11}{rgb}{0.751151720870593,0.64119344339794,0.689445159222198}
\colorlet{color11}{mycolor3}
\fill[color11!50] ($(A63)! 0.5!(A69)$) -- ($(A64)! 0.5!(A69)$) -- (A69) -- cycle;
\fill[color11!50] ($(A63)! 0.5!(A69)$) -- ($(A68)! 0.5!(A69)$) -- (A69) -- cycle;
\fill[color11!50] ($(A64)! 0.5!(A69)$) -- ($(A70)! 0.5!(A69)$) -- (A69) -- cycle;
\fill[color11!50] ($(A68)! 0.5!(A69)$) -- ($(A73)! 0.5!(A69)$) -- (A69) -- cycle;
\fill[color11!50] ($(A70)! 0.5!(A69)$) -- ($(A74)! 0.5!(A69)$) -- (A69) -- cycle;
\fill[color11!50] ($(A73)! 0.5!(A69)$) -- ($(A74)! 0.5!(A69)$) -- (A69) -- cycle;
%\definecolor{color12}{rgb}{0.646953469689286,0.58513004595989,0.828410294165587}
\colorlet{color12}{mycolor3}
\fill[color12!50] ($(A72)! 0.5!(A76)$) -- ($(A73)! 0.5!(A76)$) -- (A76) -- cycle;
\fill[color12!50] ($(A72)! 0.5!(A76)$) -- ($(A75)! 0.5!(A76)$) -- (A76) -- cycle;
\fill[color12!50] ($(A73)! 0.5!(A76)$) -- ($(A77)! 0.5!(A76)$) -- (A76) -- cycle;
\fill[color12!50] ($(A75)! 0.5!(A76)$) -- ($(A78)! 0.5!(A76)$) -- (A76) -- cycle;
\fill[color12!50] ($(A77)! 0.5!(A76)$) -- ($(A79)! 0.5!(A76)$) -- (A76) -- cycle;
\fill[color12!50] ($(A78)! 0.5!(A76)$) -- ($(A79)! 0.5!(A76)$) -- (A76) -- cycle;
\colorlet{color13}{mycolor3}
%\definecolor{color13}{rgb}{0.834259557073338,0.643195071638504,0.65280734899823}
\fill[color13!50] ($(A49)! 0.5!(A83)$) -- ($(A50)! 0.5!(A83)$) -- (A83) -- cycle;
\fill[color13!50] ($(A49)! 0.5!(A83)$) -- ($(A84)! 0.5!(A83)$) -- (A83) -- cycle;
\fill[color13!50] ($(A50)! 0.5!(A83)$) -- ($(A82)! 0.5!(A83)$) -- (A83) -- cycle;
\fill[color13!50] ($(A82)! 0.5!(A83)$) -- ($(A97)! 0.5!(A83)$) -- (A83) -- cycle;
\fill[color13!50] ($(A84)! 0.5!(A83)$) -- ($(A98)! 0.5!(A83)$) -- (A83) -- cycle;
\fill[color13!50] ($(A97)! 0.5!(A83)$) -- ($(A98)! 0.5!(A83)$) -- (A83) -- cycle;
%\definecolor{color14}{rgb}{0.567314771446259,0.512228331215659,0.96809306106272}
\colorlet{color14}{mycolor3}
\fill[color14!50] ($(A47)! 0.5!(A85)$) -- ($(A48)! 0.5!(A85)$) -- (A85) -- cycle;
\fill[color14!50] ($(A47)! 0.5!(A85)$) -- ($(A100)! 0.5!(A85)$) -- (A85) -- cycle;
\fill[color14!50] ($(A48)! 0.5!(A85)$) -- ($(A84)! 0.5!(A85)$) -- (A85) -- cycle;
\fill[color14!50] ($(A84)! 0.5!(A85)$) -- ($(A99)! 0.5!(A85)$) -- (A85) -- cycle;
\fill[color14!50] ($(A99)! 0.5!(A85)$) -- ($(A100)! 0.5!(A85)$) -- (A85) -- cycle;
%\definecolor{color15}{rgb}{0.627716956535899,0.524355254666345,0.385959122765975}
\colorlet{color15}{mycolor5}
\fill[color15!50] ($(A2)! 0.5!(A90)$) -- ($(A3)! 0.5!(A90)$) -- (A90) -- cycle;
\fill[color15!50] ($(A2)! 0.5!(A90)$) -- ($(A89)! 0.5!(A90)$) -- (A90) -- cycle;
\fill[color15!50] ($(A3)! 0.5!(A90)$) -- ($(A115)! 0.5!(A90)$) -- (A90) -- cycle;
\fill[color15!50] ($(A89)! 0.5!(A90)$) -- ($(A103)! 0.5!(A90)$) -- (A90) -- cycle;
\fill[color15!50] ($(A103)! 0.5!(A90)$) -- ($(A115)! 0.5!(A90)$) -- (A90) -- cycle;
%\definecolor{color16}{rgb}{0.593885436469775,0.490976094830578,0.696906969034984}
\colorlet{color16}{mycolor3}
\fill[color16!50] ($(A91)! 0.5!(A105)$) -- ($(A92)! 0.5!(A105)$) -- (A105) -- cycle;
\fill[color16!50] ($(A91)! 0.5!(A105)$) -- ($(A104)! 0.5!(A105)$) -- (A105) -- cycle;
\fill[color16!50] ($(A92)! 0.5!(A105)$) -- ($(A106)! 0.5!(A105)$) -- (A105) -- cycle;
\fill[color16!50] ($(A104)! 0.5!(A105)$) -- ($(A116)! 0.5!(A105)$) -- (A105) -- cycle;
\fill[color16!50] ($(A106)! 0.5!(A105)$) -- ($(A117)! 0.5!(A105)$) -- (A105) -- cycle;
\fill[color16!50] ($(A116)! 0.5!(A105)$) -- ($(A117)! 0.5!(A105)$) -- (A105) -- cycle;
%\definecolor{color17}{rgb}{0.74822561818672,0.616072526251147,0.726131008645574}
\colorlet{color17}{mycolor3}
\fill[color17!50] ($(A93)! 0.5!(A107)$) -- ($(A94)! 0.5!(A107)$) -- (A107) -- cycle;
\fill[color17!50] ($(A93)! 0.5!(A107)$) -- ($(A106)! 0.5!(A107)$) -- (A107) -- cycle;
\fill[color17!50] ($(A94)! 0.5!(A107)$) -- ($(A108)! 0.5!(A107)$) -- (A107) -- cycle;
\fill[color17!50] ($(A106)! 0.5!(A107)$) -- ($(A118)! 0.5!(A107)$) -- (A107) -- cycle;
\fill[color17!50] ($(A108)! 0.5!(A107)$) -- ($(A119)! 0.5!(A107)$) -- (A107) -- cycle;
\fill[color17!50] ($(A118)! 0.5!(A107)$) -- ($(A119)! 0.5!(A107)$) -- (A107) -- cycle;
%\definecolor{color18}{rgb}{0.554257003659849,0.607412158271001,0.820895776953617}
\colorlet{color18}{mycolor3}
\fill[color18!50] ($(A98)! 0.5!(A111)$) -- ($(A99)! 0.5!(A111)$) -- (A111) -- cycle;
\fill[color18!50] ($(A98)! 0.5!(A111)$) -- ($(A110)! 0.5!(A111)$) -- (A111) -- cycle;
\fill[color18!50] ($(A99)! 0.5!(A111)$) -- ($(A112)! 0.5!(A111)$) -- (A111) -- cycle;
\fill[color18!50] ($(A110)! 0.5!(A111)$) -- ($(A121)! 0.5!(A111)$) -- (A111) -- cycle;
\fill[color18!50] ($(A112)! 0.5!(A111)$) -- ($(A122)! 0.5!(A111)$) -- (A111) -- cycle;
\fill[color18!50] ($(A121)! 0.5!(A111)$) -- ($(A122)! 0.5!(A111)$) -- (A111) -- cycle;
%\definecolor{color19}{rgb}{0.66250035958106,0.614370291942401,0.586876119499875}
\colorlet{color19}{mycolor3}
\fill[color19!50] ($(A104)! 0.5!(A125)$) -- ($(A116)! 0.5!(A125)$) -- (A125) -- cycle;
\fill[color19!50] ($(A104)! 0.5!(A125)$) -- ($(A132)! 0.5!(A125)$) -- (A125) -- cycle;
\fill[color19!50] ($(A116)! 0.5!(A125)$) -- ($(A126)! 0.5!(A125)$) -- (A125) -- cycle;
\fill[color19!50] ($(A126)! 0.5!(A125)$) -- ($(A133)! 0.5!(A125)$) -- (A125) -- cycle;
\fill[color19!50] ($(A132)! 0.5!(A125)$) -- ($(A133)! 0.5!(A125)$) -- (A125) -- cycle;
%\definecolor{color20}{rgb}{0.638691059742165,0.645588909263044,0.548475035032595}
\colorlet{color20}{mycolor3}
\fill[color20!50] ($(A117)! 0.5!(A127)$) -- ($(A118)! 0.5!(A127)$) -- (A127) -- cycle;
\fill[color20!50] ($(A117)! 0.5!(A127)$) -- ($(A126)! 0.5!(A127)$) -- (A127) -- cycle;
\fill[color20!50] ($(A118)! 0.5!(A127)$) -- ($(A128)! 0.5!(A127)$) -- (A127) -- cycle;
\fill[color20!50] ($(A126)! 0.5!(A127)$) -- ($(A134)! 0.5!(A127)$) -- (A127) -- cycle;
\fill[color20!50] ($(A128)! 0.5!(A127)$) -- ($(A135)! 0.5!(A127)$) -- (A127) -- cycle;
\fill[color20!50] ($(A134)! 0.5!(A127)$) -- ($(A135)! 0.5!(A127)$) -- (A127) -- cycle;
%\definecolor{color21}{rgb}{0.510967642876474,0.826269352867521,0.916303467855384}
\colorlet{color21}{mycolor3}
\fill[color21!50] ($(A122)! 0.5!(A131)$) -- ($(A123)! 0.5!(A131)$) -- (A131) -- cycle;
\fill[color21!50] ($(A122)! 0.5!(A131)$) -- ($(A130)! 0.5!(A131)$) -- (A131) -- cycle;
\fill[color21!50] ($(A123)! 0.5!(A131)$) -- ($(A132)! 0.5!(A131)$) -- (A131) -- cycle;
\fill[color21!50] ($(A130)! 0.5!(A131)$) -- ($(A137)! 0.5!(A131)$) -- (A131) -- cycle;
\fill[color21!50] ($(A132)! 0.5!(A131)$) -- ($(A138)! 0.5!(A131)$) -- (A131) -- cycle;
\fill[color21!50] ($(A137)! 0.5!(A131)$) -- ($(A138)! 0.5!(A131)$) -- (A131) -- cycle;
%\definecolor{color22}{rgb}{0.770916622867055,0.810927562842745,0.926787501119097}
\colorlet{color22}{mycolor3}
\fill[color22!50] ($(A133)! 0.5!(A139)$) -- ($(A134)! 0.5!(A139)$) -- (A139) -- cycle;
\fill[color22!50] ($(A133)! 0.5!(A139)$) -- ($(A138)! 0.5!(A139)$) -- (A139) -- cycle;
\fill[color22!50] ($(A134)! 0.5!(A139)$) -- ($(A140)! 0.5!(A139)$) -- (A139) -- cycle;
\fill[color22!50] ($(A138)! 0.5!(A139)$) -- ($(A142)! 0.5!(A139)$) -- (A139) -- cycle;
\fill[color22!50] ($(A140)! 0.5!(A139)$) -- ($(A143)! 0.5!(A139)$) -- (A139) -- cycle;
\fill[color22!50] ($(A142)! 0.5!(A139)$) -- ($(A143)! 0.5!(A139)$) -- (A139) -- cycle;
\draw[color22, thick] ($(A133)! 0.5!(A139)$) -- ($(A134)! 0.5!(A139)$);
\draw[color22, thick] ($(A133)! 0.5!(A139)$) -- ($(A138)! 0.5!(A139)$);
\draw[color22, thick] ($(A134)! 0.5!(A139)$) -- ($(A140)! 0.5!(A139)$);
\draw[color22, thick] ($(A138)! 0.5!(A139)$) -- ($(A142)! 0.5!(A139)$);
\draw[color22, thick] ($(A140)! 0.5!(A139)$) -- ($(A143)! 0.5!(A139)$);
\draw[color22, thick] ($(A142)! 0.5!(A139)$) -- ($(A143)! 0.5!(A139)$);
\draw[color21, thick] ($(A122)! 0.5!(A131)$) -- ($(A123)! 0.5!(A131)$);
\draw[color21, thick] ($(A122)! 0.5!(A131)$) -- ($(A130)! 0.5!(A131)$);
\draw[color21, thick] ($(A123)! 0.5!(A131)$) -- ($(A132)! 0.5!(A131)$);
\draw[color21, thick] ($(A130)! 0.5!(A131)$) -- ($(A137)! 0.5!(A131)$);
\draw[color21, thick] ($(A132)! 0.5!(A131)$) -- ($(A138)! 0.5!(A131)$);
\draw[color21, thick] ($(A137)! 0.5!(A131)$) -- ($(A138)! 0.5!(A131)$);
\draw[color20, thick] ($(A117)! 0.5!(A127)$) -- ($(A118)! 0.5!(A127)$);
\draw[color20, thick] ($(A117)! 0.5!(A127)$) -- ($(A126)! 0.5!(A127)$);
\draw[color20, thick] ($(A118)! 0.5!(A127)$) -- ($(A128)! 0.5!(A127)$);
\draw[color20, thick] ($(A126)! 0.5!(A127)$) -- ($(A134)! 0.5!(A127)$);
\draw[color20, thick] ($(A128)! 0.5!(A127)$) -- ($(A135)! 0.5!(A127)$);
\draw[color20, thick] ($(A134)! 0.5!(A127)$) -- ($(A135)! 0.5!(A127)$);
\draw[color19, thick] ($(A104)! 0.5!(A125)$) -- ($(A116)! 0.5!(A125)$);
\draw[color19, thick] ($(A104)! 0.5!(A125)$) -- ($(A132)! 0.5!(A125)$);
\draw[color19, thick] ($(A116)! 0.5!(A125)$) -- ($(A126)! 0.5!(A125)$);
\draw[color19, thick] ($(A126)! 0.5!(A125)$) -- ($(A133)! 0.5!(A125)$);
\draw[color19, thick] ($(A132)! 0.5!(A125)$) -- ($(A133)! 0.5!(A125)$);
\draw[color18, thick] ($(A98)! 0.5!(A111)$) -- ($(A99)! 0.5!(A111)$);
\draw[color18, thick] ($(A98)! 0.5!(A111)$) -- ($(A110)! 0.5!(A111)$);
\draw[color18, thick] ($(A99)! 0.5!(A111)$) -- ($(A112)! 0.5!(A111)$);
\draw[color18, thick] ($(A110)! 0.5!(A111)$) -- ($(A121)! 0.5!(A111)$);
\draw[color18, thick] ($(A112)! 0.5!(A111)$) -- ($(A122)! 0.5!(A111)$);
\draw[color18, thick] ($(A121)! 0.5!(A111)$) -- ($(A122)! 0.5!(A111)$);
\draw[color17, thick] ($(A93)! 0.5!(A107)$) -- ($(A94)! 0.5!(A107)$);
\draw[color17, thick] ($(A93)! 0.5!(A107)$) -- ($(A106)! 0.5!(A107)$);
\draw[color17, thick] ($(A94)! 0.5!(A107)$) -- ($(A108)! 0.5!(A107)$);
\draw[color17, thick] ($(A106)! 0.5!(A107)$) -- ($(A118)! 0.5!(A107)$);
\draw[color17, thick] ($(A108)! 0.5!(A107)$) -- ($(A119)! 0.5!(A107)$);
\draw[color17, thick] ($(A118)! 0.5!(A107)$) -- ($(A119)! 0.5!(A107)$);
\draw[color16, thick] ($(A91)! 0.5!(A105)$) -- ($(A92)! 0.5!(A105)$);
\draw[color16, thick] ($(A91)! 0.5!(A105)$) -- ($(A104)! 0.5!(A105)$);
\draw[color16, thick] ($(A92)! 0.5!(A105)$) -- ($(A106)! 0.5!(A105)$);
\draw[color16, thick] ($(A104)! 0.5!(A105)$) -- ($(A116)! 0.5!(A105)$);
\draw[color16, thick] ($(A106)! 0.5!(A105)$) -- ($(A117)! 0.5!(A105)$);
\draw[color16, thick] ($(A116)! 0.5!(A105)$) -- ($(A117)! 0.5!(A105)$);
\draw[color14, thick] ($(A47)! 0.5!(A85)$) -- ($(A48)! 0.5!(A85)$);
\draw[color14, thick] ($(A47)! 0.5!(A85)$) -- ($(A100)! 0.5!(A85)$);
\draw[color14, thick] ($(A48)! 0.5!(A85)$) -- ($(A84)! 0.5!(A85)$);
\draw[color14, thick] ($(A84)! 0.5!(A85)$) -- ($(A99)! 0.5!(A85)$);
\draw[color14, thick] ($(A99)! 0.5!(A85)$) -- ($(A100)! 0.5!(A85)$);
\draw[color13, thick] ($(A49)! 0.5!(A83)$) -- ($(A50)! 0.5!(A83)$);
\draw[color13, thick] ($(A49)! 0.5!(A83)$) -- ($(A84)! 0.5!(A83)$);
\draw[color13, thick] ($(A50)! 0.5!(A83)$) -- ($(A82)! 0.5!(A83)$);
\draw[color13, thick] ($(A82)! 0.5!(A83)$) -- ($(A97)! 0.5!(A83)$);
\draw[color13, thick] ($(A84)! 0.5!(A83)$) -- ($(A98)! 0.5!(A83)$);
\draw[color13, thick] ($(A97)! 0.5!(A83)$) -- ($(A98)! 0.5!(A83)$);
\draw[color12, thick] ($(A72)! 0.5!(A76)$) -- ($(A73)! 0.5!(A76)$);
\draw[color12, thick] ($(A72)! 0.5!(A76)$) -- ($(A75)! 0.5!(A76)$);
\draw[color12, thick] ($(A73)! 0.5!(A76)$) -- ($(A77)! 0.5!(A76)$);
\draw[color12, thick] ($(A75)! 0.5!(A76)$) -- ($(A78)! 0.5!(A76)$);
\draw[color12, thick] ($(A77)! 0.5!(A76)$) -- ($(A79)! 0.5!(A76)$);
\draw[color12, thick] ($(A78)! 0.5!(A76)$) -- ($(A79)! 0.5!(A76)$);
\draw[color11, thick] ($(A63)! 0.5!(A69)$) -- ($(A64)! 0.5!(A69)$);
\draw[color11, thick] ($(A63)! 0.5!(A69)$) -- ($(A68)! 0.5!(A69)$);
\draw[color11, thick] ($(A64)! 0.5!(A69)$) -- ($(A70)! 0.5!(A69)$);
\draw[color11, thick] ($(A68)! 0.5!(A69)$) -- ($(A73)! 0.5!(A69)$);
\draw[color11, thick] ($(A70)! 0.5!(A69)$) -- ($(A74)! 0.5!(A69)$);
\draw[color11, thick] ($(A73)! 0.5!(A69)$) -- ($(A74)! 0.5!(A69)$);
\draw[color10, thick] ($(A61)! 0.5!(A67)$) -- ($(A62)! 0.5!(A67)$);
\draw[color10, thick] ($(A61)! 0.5!(A67)$) -- ($(A71)! 0.5!(A67)$);
\draw[color10, thick] ($(A62)! 0.5!(A67)$) -- ($(A68)! 0.5!(A67)$);
\draw[color10, thick] ($(A68)! 0.5!(A67)$) -- ($(A72)! 0.5!(A67)$);
\draw[color10, thick] ($(A71)! 0.5!(A67)$) -- ($(A72)! 0.5!(A67)$);
\draw[color8, thick] ($(A50)! 0.5!(A58)$) -- ($(A51)! 0.5!(A58)$);
\draw[color8, thick] ($(A50)! 0.5!(A58)$) -- ($(A57)! 0.5!(A58)$);
\draw[color8, thick] ($(A51)! 0.5!(A58)$) -- ($(A59)! 0.5!(A58)$);
\draw[color8, thick] ($(A57)! 0.5!(A58)$) -- ($(A64)! 0.5!(A58)$);
\draw[color8, thick] ($(A59)! 0.5!(A58)$) -- ($(A65)! 0.5!(A58)$);
\draw[color8, thick] ($(A64)! 0.5!(A58)$) -- ($(A65)! 0.5!(A58)$);
\draw[color7, thick] ($(A48)! 0.5!(A56)$) -- ($(A49)! 0.5!(A56)$);
\draw[color7, thick] ($(A48)! 0.5!(A56)$) -- ($(A55)! 0.5!(A56)$);
\draw[color7, thick] ($(A49)! 0.5!(A56)$) -- ($(A57)! 0.5!(A56)$);
\draw[color7, thick] ($(A55)! 0.5!(A56)$) -- ($(A62)! 0.5!(A56)$);
\draw[color7, thick] ($(A57)! 0.5!(A56)$) -- ($(A63)! 0.5!(A56)$);
\draw[color7, thick] ($(A62)! 0.5!(A56)$) -- ($(A63)! 0.5!(A56)$);
\draw[color1, thick] ($(A0)! 0.5!(A88)$) -- ($(A1)! 0.5!(A88)$);
\draw[color1, thick] ($(A0)! 0.5!(A45)$) -- ($(A9)! 0.5!(A45)$);
\draw[color1, thick] ($(A45)! 0.5!(A0)$) -- ($(A88)! 0.5!(A0)$);
\draw[color1, thick] ($(A1)! 0.5!(A10)$) -- ($(A2)! 0.5!(A10)$);
\draw[color1, thick] ($(A1)! 0.5!(A10)$) -- ($(A9)! 0.5!(A10)$);
\draw[color1, thick] ($(A1)! 0.5!(A88)$) -- ($(A89)! 0.5!(A88)$);
\draw[color1, thick] ($(A2)! 0.5!(A90)$) -- ($(A3)! 0.5!(A90)$);
\draw[color1, thick] ($(A2)! 0.5!(A10)$) -- ($(A11)! 0.5!(A10)$);
\draw[color1, thick] ($(A2)! 0.5!(A90)$) -- ($(A89)! 0.5!(A90)$);
\draw[color1, thick] ($(A3)! 0.5!(A4)$) -- ($(A12)! 0.5!(A4)$);
\draw[color1, thick] ($(A4)! 0.5!(A3)$) -- ($(A91)! 0.5!(A3)$);
\draw[color1, thick] ($(A3)! 0.5!(A25)$) -- ($(A11)! 0.5!(A25)$);
\draw[color1, thick] ($(A3)! 0.5!(A19)$) -- ($(A12)! 0.5!(A19)$);
\draw[color1, thick] ($(A3)! 0.5!(A19)$) -- ($(A35)! 0.5!(A19)$);
\draw[color1, thick] ($(A3)! 0.5!(A25)$) -- ($(A35)! 0.5!(A25)$);
\draw[color1, thick] ($(A3)! 0.5!(A90)$) -- ($(A115)! 0.5!(A90)$);
\draw[color1, thick] ($(A91)! 0.5!(A3)$) -- ($(A104)! 0.5!(A3)$);
\draw[color1, thick] ($(A104)! 0.5!(A3)$) -- ($(A132)! 0.5!(A3)$);
\draw[color1, thick] ($(A3)! 0.5!(A132)$) -- ($(A115)! 0.5!(A132)$);
\draw[color1, thick] ($(A5)! 0.5!(A4)$) -- ($(A13)! 0.5!(A4)$);
\draw[color1, thick] ($(A4)! 0.5!(A5)$) -- ($(A92)! 0.5!(A5)$);
\draw[color1, thick] ($(A12)! 0.5!(A4)$) -- ($(A13)! 0.5!(A4)$);
\draw[color1, thick] ($(A5)! 0.5!(A6)$) -- ($(A14)! 0.5!(A6)$);
\draw[color1, thick] ($(A6)! 0.5!(A5)$) -- ($(A93)! 0.5!(A5)$);
\draw[color1, thick] ($(A92)! 0.5!(A5)$) -- ($(A93)! 0.5!(A5)$);
\draw[color1, thick] ($(A7)! 0.5!(A6)$) -- ($(A15)! 0.5!(A6)$);
\draw[color1, thick] ($(A6)! 0.5!(A7)$) -- ($(A94)! 0.5!(A7)$);
\draw[color1, thick] ($(A14)! 0.5!(A6)$) -- ($(A15)! 0.5!(A6)$);
\draw[color1, thick] ($(A7)! 0.5!(A8)$) -- ($(A16)! 0.5!(A8)$);
\draw[color1, thick] ($(A8)! 0.5!(A7)$) -- ($(A95)! 0.5!(A7)$);
\draw[color1, thick] ($(A94)! 0.5!(A7)$) -- ($(A95)! 0.5!(A7)$);
\draw[color1, thick] ($(A9)! 0.5!(A10)$) -- ($(A17)! 0.5!(A10)$);
\draw[color1, thick] ($(A9)! 0.5!(A45)$) -- ($(A53)! 0.5!(A45)$);
\draw[color1, thick] ($(A11)! 0.5!(A10)$) -- ($(A18)! 0.5!(A10)$);
\draw[color1, thick] ($(A17)! 0.5!(A10)$) -- ($(A18)! 0.5!(A10)$);
\draw[color1, thick] ($(A11)! 0.5!(A25)$) -- ($(A18)! 0.5!(A25)$);
\draw[color1, thick] ($(A12)! 0.5!(A19)$) -- ($(A20)! 0.5!(A19)$);
\draw[color1, thick] ($(A13)! 0.5!(A21)$) -- ($(A14)! 0.5!(A21)$);
\draw[color1, thick] ($(A13)! 0.5!(A21)$) -- ($(A20)! 0.5!(A21)$);
\draw[color1, thick] ($(A14)! 0.5!(A21)$) -- ($(A22)! 0.5!(A21)$);
\draw[color1, thick] ($(A15)! 0.5!(A23)$) -- ($(A16)! 0.5!(A23)$);
\draw[color1, thick] ($(A15)! 0.5!(A23)$) -- ($(A22)! 0.5!(A23)$);
\draw[color1, thick] ($(A17)! 0.5!(A60)$) -- ($(A24)! 0.5!(A60)$);
\draw[color1, thick] ($(A17)! 0.5!(A60)$) -- ($(A53)! 0.5!(A60)$);
\draw[color1, thick] ($(A18)! 0.5!(A25)$) -- ($(A24)! 0.5!(A25)$);
\draw[color1, thick] ($(A20)! 0.5!(A19)$) -- ($(A26)! 0.5!(A19)$);
\draw[color1, thick] ($(A26)! 0.5!(A19)$) -- ($(A31)! 0.5!(A19)$);
\draw[color1, thick] ($(A31)! 0.5!(A19)$) -- ($(A35)! 0.5!(A19)$);
\draw[color1, thick] ($(A20)! 0.5!(A21)$) -- ($(A27)! 0.5!(A21)$);
\draw[color1, thick] ($(A22)! 0.5!(A21)$) -- ($(A28)! 0.5!(A21)$);
\draw[color1, thick] ($(A27)! 0.5!(A21)$) -- ($(A28)! 0.5!(A21)$);
\draw[color1, thick] ($(A22)! 0.5!(A23)$) -- ($(A29)! 0.5!(A23)$);
\draw[color1, thick] ($(A24)! 0.5!(A25)$) -- ($(A30)! 0.5!(A25)$);
\draw[color1, thick] ($(A24)! 0.5!(A60)$) -- ($(A66)! 0.5!(A60)$);
\draw[color1, thick] ($(A30)! 0.5!(A25)$) -- ($(A35)! 0.5!(A25)$);
\draw[color1, thick] ($(A26)! 0.5!(A32)$) -- ($(A27)! 0.5!(A32)$);
\draw[color1, thick] ($(A26)! 0.5!(A32)$) -- ($(A31)! 0.5!(A32)$);
\draw[color1, thick] ($(A27)! 0.5!(A32)$) -- ($(A33)! 0.5!(A32)$);
\draw[color1, thick] ($(A28)! 0.5!(A34)$) -- ($(A29)! 0.5!(A34)$);
\draw[color1, thick] ($(A28)! 0.5!(A34)$) -- ($(A33)! 0.5!(A34)$);
\draw[color1, thick] ($(A31)! 0.5!(A32)$) -- ($(A36)! 0.5!(A32)$);
\draw[color1, thick] ($(A33)! 0.5!(A32)$) -- ($(A37)! 0.5!(A32)$);
\draw[color1, thick] ($(A36)! 0.5!(A32)$) -- ($(A37)! 0.5!(A32)$);
\draw[color1, thick] ($(A33)! 0.5!(A34)$) -- ($(A38)! 0.5!(A34)$);
\draw[color1, thick] ($(A35)! 0.5!(A39)$) -- ($(A36)! 0.5!(A39)$);
\draw[color1, thick] ($(A39)! 0.5!(A35)$) -- ($(A75)! 0.5!(A35)$);
\draw[color1, thick] ($(A47)! 0.5!(A35)$) -- ($(A61)! 0.5!(A35)$);
\draw[color1, thick] ($(A35)! 0.5!(A47)$) -- ($(A66)! 0.5!(A47)$);
\draw[color1, thick] ($(A61)! 0.5!(A35)$) -- ($(A71)! 0.5!(A35)$);
\draw[color1, thick] ($(A71)! 0.5!(A35)$) -- ($(A75)! 0.5!(A35)$);
\draw[color1, thick] ($(A36)! 0.5!(A39)$) -- ($(A40)! 0.5!(A39)$);
\draw[color1, thick] ($(A37)! 0.5!(A41)$) -- ($(A38)! 0.5!(A41)$);
\draw[color1, thick] ($(A37)! 0.5!(A41)$) -- ($(A40)! 0.5!(A41)$);
\draw[color1, thick] ($(A40)! 0.5!(A39)$) -- ($(A42)! 0.5!(A39)$);
\draw[color1, thick] ($(A39)! 0.5!(A42)$) -- ($(A78)! 0.5!(A42)$);
\draw[color1, thick] ($(A40)! 0.5!(A41)$) -- ($(A43)! 0.5!(A41)$);
\draw[color1, thick] ($(A42)! 0.5!(A44)$) -- ($(A43)! 0.5!(A44)$);
\draw[color1, thick] ($(A44)! 0.5!(A42)$) -- ($(A80)! 0.5!(A42)$);
\draw[color1, thick] ($(A78)! 0.5!(A42)$) -- ($(A80)! 0.5!(A42)$);
\draw[color1, thick] ($(A46)! 0.5!(A45)$) -- ($(A53)! 0.5!(A45)$);
\draw[color1, thick] ($(A45)! 0.5!(A46)$) -- ($(A87)! 0.5!(A46)$);
\draw[color1, thick] ($(A46)! 0.5!(A47)$) -- ($(A54)! 0.5!(A47)$);
\draw[color1, thick] ($(A47)! 0.5!(A46)$) -- ($(A86)! 0.5!(A46)$);
\draw[color1, thick] ($(A86)! 0.5!(A46)$) -- ($(A87)! 0.5!(A46)$);
\draw[color1, thick] ($(A54)! 0.5!(A47)$) -- ($(A66)! 0.5!(A47)$);
\draw[color1, thick] ($(A51)! 0.5!(A81)$) -- ($(A52)! 0.5!(A81)$);
\draw[color1, thick] ($(A51)! 0.5!(A81)$) -- ($(A82)! 0.5!(A81)$);
\draw[color1, thick] ($(A53)! 0.5!(A60)$) -- ($(A54)! 0.5!(A60)$);
\draw[color1, thick] ($(A54)! 0.5!(A60)$) -- ($(A66)! 0.5!(A60)$);
\draw[color1, thick] ($(A82)! 0.5!(A81)$) -- ($(A96)! 0.5!(A81)$);
\draw[color1, thick] ($(A87)! 0.5!(A102)$) -- ($(A88)! 0.5!(A102)$);
\draw[color1, thick] ($(A87)! 0.5!(A102)$) -- ($(A101)! 0.5!(A102)$);
\draw[color1, thick] ($(A89)! 0.5!(A88)$) -- ($(A102)! 0.5!(A88)$);
\draw[color1, thick] ($(A89)! 0.5!(A90)$) -- ($(A103)! 0.5!(A90)$);
\draw[color1, thick] ($(A103)! 0.5!(A90)$) -- ($(A115)! 0.5!(A90)$);
\draw[color1, thick] ($(A96)! 0.5!(A109)$) -- ($(A97)! 0.5!(A109)$);
\draw[color1, thick] ($(A97)! 0.5!(A109)$) -- ($(A110)! 0.5!(A109)$);
\draw[color1, thick] ($(A101)! 0.5!(A102)$) -- ($(A114)! 0.5!(A102)$);
\draw[color1, thick] ($(A102)! 0.5!(A114)$) -- ($(A103)! 0.5!(A114)$);
\draw[color1, thick] ($(A103)! 0.5!(A114)$) -- ($(A115)! 0.5!(A114)$);
\draw[color1, thick] ($(A110)! 0.5!(A109)$) -- ($(A120)! 0.5!(A109)$);
\draw[color1, thick] ($(A113)! 0.5!(A124)$) -- ($(A114)! 0.5!(A124)$);
\draw[color1, thick] ($(A113)! 0.5!(A124)$) -- ($(A132)! 0.5!(A124)$);
\draw[color1, thick] ($(A115)! 0.5!(A114)$) -- ($(A124)! 0.5!(A114)$);
\draw[color1, thick] ($(A115)! 0.5!(A132)$) -- ($(A124)! 0.5!(A132)$);
\draw[color1, thick] ($(A120)! 0.5!(A129)$) -- ($(A121)! 0.5!(A129)$);
\draw[color1, thick] ($(A121)! 0.5!(A129)$) -- ($(A130)! 0.5!(A129)$);
\draw[color1, thick] ($(A130)! 0.5!(A129)$) -- ($(A136)! 0.5!(A129)$);
\draw[color1, thick] ($(A136)! 0.5!(A141)$) -- ($(A137)! 0.5!(A141)$);
\draw[color1, thick] ($(A137)! 0.5!(A141)$) -- ($(A142)! 0.5!(A141)$);
\draw[color1, thick] ($(A142)! 0.5!(A141)$) -- ($(A144)! 0.5!(A141)$);
\draw[color15, thick] ($(A2)! 0.5!(A90)$) -- ($(A3)! 0.5!(A90)$);
\draw[color15, thick] ($(A2)! 0.5!(A90)$) -- ($(A89)! 0.5!(A90)$);
\draw[color15, thick] ($(A3)! 0.5!(A90)$) -- ($(A115)! 0.5!(A90)$);
\draw[color15, thick] ($(A89)! 0.5!(A90)$) -- ($(A103)! 0.5!(A90)$);
\draw[color15, thick] ($(A103)! 0.5!(A90)$) -- ($(A115)! 0.5!(A90)$);
\draw[color9, thick] ($(A17)! 0.5!(A60)$) -- ($(A24)! 0.5!(A60)$);
\draw[color9, thick] ($(A17)! 0.5!(A60)$) -- ($(A53)! 0.5!(A60)$);
\draw[color9, thick] ($(A24)! 0.5!(A60)$) -- ($(A66)! 0.5!(A60)$);
\draw[color9, thick] ($(A53)! 0.5!(A60)$) -- ($(A54)! 0.5!(A60)$);
\draw[color9, thick] ($(A54)! 0.5!(A60)$) -- ($(A66)! 0.5!(A60)$);
\draw[color6, thick] ($(A26)! 0.5!(A32)$) -- ($(A27)! 0.5!(A32)$);
\draw[color6, thick] ($(A26)! 0.5!(A32)$) -- ($(A31)! 0.5!(A32)$);
\draw[color6, thick] ($(A27)! 0.5!(A32)$) -- ($(A33)! 0.5!(A32)$);
\draw[color6, thick] ($(A31)! 0.5!(A32)$) -- ($(A36)! 0.5!(A32)$);
\draw[color6, thick] ($(A33)! 0.5!(A32)$) -- ($(A37)! 0.5!(A32)$);
\draw[color6, thick] ($(A36)! 0.5!(A32)$) -- ($(A37)! 0.5!(A32)$);
\draw[color5, thick] ($(A3)! 0.5!(A25)$) -- ($(A11)! 0.5!(A25)$);
\draw[color5, thick] ($(A3)! 0.5!(A25)$) -- ($(A35)! 0.5!(A25)$);
\draw[color5, thick] ($(A11)! 0.5!(A25)$) -- ($(A18)! 0.5!(A25)$);
\draw[color5, thick] ($(A18)! 0.5!(A25)$) -- ($(A24)! 0.5!(A25)$);
\draw[color5, thick] ($(A24)! 0.5!(A25)$) -- ($(A30)! 0.5!(A25)$);
\draw[color5, thick] ($(A30)! 0.5!(A25)$) -- ($(A35)! 0.5!(A25)$);
\draw[color4, thick] ($(A13)! 0.5!(A21)$) -- ($(A14)! 0.5!(A21)$);
\draw[color4, thick] ($(A13)! 0.5!(A21)$) -- ($(A20)! 0.5!(A21)$);
\draw[color4, thick] ($(A14)! 0.5!(A21)$) -- ($(A22)! 0.5!(A21)$);
\draw[color4, thick] ($(A20)! 0.5!(A21)$) -- ($(A27)! 0.5!(A21)$);
\draw[color4, thick] ($(A22)! 0.5!(A21)$) -- ($(A28)! 0.5!(A21)$);
\draw[color4, thick] ($(A27)! 0.5!(A21)$) -- ($(A28)! 0.5!(A21)$);
\draw[color3, thick] ($(A3)! 0.5!(A19)$) -- ($(A12)! 0.5!(A19)$);
\draw[color3, thick] ($(A3)! 0.5!(A19)$) -- ($(A35)! 0.5!(A19)$);
\draw[color3, thick] ($(A12)! 0.5!(A19)$) -- ($(A20)! 0.5!(A19)$);
\draw[color3, thick] ($(A20)! 0.5!(A19)$) -- ($(A26)! 0.5!(A19)$);
\draw[color3, thick] ($(A26)! 0.5!(A19)$) -- ($(A31)! 0.5!(A19)$);
\draw[color3, thick] ($(A31)! 0.5!(A19)$) -- ($(A35)! 0.5!(A19)$);
\draw[color2, thick] ($(A1)! 0.5!(A10)$) -- ($(A2)! 0.5!(A10)$);
\draw[color2, thick] ($(A1)! 0.5!(A10)$) -- ($(A9)! 0.5!(A10)$);
\draw[color2, thick] ($(A2)! 0.5!(A10)$) -- ($(A11)! 0.5!(A10)$);
\draw[color2, thick] ($(A9)! 0.5!(A10)$) -- ($(A17)! 0.5!(A10)$);
\draw[color2, thick] ($(A11)! 0.5!(A10)$) -- ($(A18)! 0.5!(A10)$);
\draw[color2, thick] ($(A17)! 0.5!(A10)$) -- ($(A18)! 0.5!(A10)$);
\foreach \a/\b/\c in {0/1/9,0/1/88,0/9/45,0/45/88,1/2/10,1/2/89,1/9/10,1/88/89,2/3/11,2/3/90,2/10/11,2/89/90,3/4/12,3/4/91,3/11/25,3/12/19,3/19/35,3/25/35,3/90/115,3/91/104,3/104/132,3/115/132,4/5/13,4/5/92,4/12/13,4/91/92,5/6/14,5/6/93,5/13/14,5/92/93,6/7/15,6/7/94,6/14/15,6/93/94,7/8/16,7/8/95,7/15/16,7/94/95,9/10/17,9/17/53,9/45/53,10/11/18,10/17/18,11/18/25,12/13/20,12/19/20,13/14/21,13/20/21,14/15/22,14/21/22,15/16/23,15/22/23,17/18/24,17/24/60,17/53/60,18/24/25,19/20/26,19/26/31,19/31/35,20/21/27,20/26/27,21/22/28,21/27/28,22/23/29,22/28/29,24/25/30,24/30/66,24/60/66,25/30/35,26/27/32,26/31/32,27/28/33,27/32/33,28/29/34,28/33/34,30/35/66,31/32/36,31/35/36,32/33/37,32/36/37,33/34/38,33/37/38,35/36/39,35/39/75,35/47/61,35/47/66,35/61/71,35/71/75,36/37/40,36/39/40,37/38/41,37/40/41,39/40/42,39/42/78,39/75/78,40/41/43,40/42/43,42/43/44,42/44/80,42/78/80,45/46/53,45/46/87,45/87/88,46/47/54,46/47/86,46/53/54,46/86/87,47/48/55,47/48/85,47/54/66,47/55/61,47/85/100,47/86/113,47/100/132,47/113/132,48/49/56,48/49/84,48/55/56,48/84/85,49/50/57,49/50/83,49/56/57,49/83/84,50/51/58,50/51/82,50/57/58,50/82/83,51/52/59,51/52/81,51/58/59,51/81/82,53/54/60,54/60/66,55/56/62,55/61/62,56/57/63,56/62/63,57/58/64,57/63/64,58/59/65,58/64/65,61/62/67,61/67/71,62/63/68,62/67/68,63/64/69,63/68/69,64/65/70,64/69/70,67/68/72,67/71/72,68/69/73,68/72/73,69/70/74,69/73/74,71/72/75,72/73/76,72/75/76,73/74/77,73/76/77,75/76/78,76/77/79,76/78/79,78/79/80,81/82/96,82/83/97,82/96/97,83/84/98,83/97/98,84/85/99,84/98/99,85/99/100,86/87/101,86/101/113,87/88/102,87/101/102,88/89/102,89/90/103,89/102/103,90/103/115,91/92/105,91/104/105,92/93/106,92/105/106,93/94/107,93/106/107,94/95/108,94/107/108,96/97/109,97/98/110,97/109/110,98/99/111,98/110/111,99/100/112,99/111/112,100/112/123,100/123/132,101/102/114,101/113/114,102/103/114,103/114/115,104/105/116,104/116/125,104/125/132,105/106/117,105/116/117,106/107/118,106/117/118,107/108/119,107/118/119,109/110/120,110/111/121,110/120/121,111/112/122,111/121/122,112/122/123,113/114/124,113/124/132,114/115/124,115/124/132,116/117/126,116/125/126,117/118/127,117/126/127,118/119/128,118/127/128,120/121/129,121/122/130,121/129/130,122/123/131,122/130/131,123/131/132,125/126/133,125/132/133,126/127/134,126/133/134,127/128/135,127/134/135,129/130/136,130/131/137,130/136/137,131/132/138,131/137/138,132/133/138,133/134/139,133/138/139,134/135/140,134/139/140,136/137/141,137/138/142,137/141/142,138/139/142,139/140/143,139/142/143,141/142/144,142/143/144}{
  \draw[black!30] (A\a) -- (A\b) -- (A\c) -- cycle;
}

\draw[black, thick] (A3) -- (A35) -- (A47) -- (A132) --cycle;

\foreach \i in {0,1,2,3,5,7,9,11,12,13,14,15,16,17,18,20,22,24,26,27,28,29,30,31,33,35,36,37,38,40,42,43,46,53,54,56,58,66,67,69,76,81,83,85,89,102,103,105,107,109,111,115,124,125,127,129,131,139,141}{
  \fill[myred] (A\i) circle (3pt);
}
\foreach \i in {4,6,8,10,19,21,23,25,32,34,39,41,44,45,47,48,49,50,51,52,55,57,59,60,61,62,63,64,65,68,70,71,72,73,74,75,77,78,79,80,82,84,86,87,88,90,91,92,93,94,95,96,97,98,99,100,101,104,106,108,110,112,113,114,116,117,118,119,120,121,122,123,126,128,130,132,133,134,135,136,137,138,140,142,143,144}{
  \fill[myblue] (A\i) circle (3pt);
}

%\foreach \i in {-8,-6,-4,-2,0,2,4,6,8}{
%	\foreach \j in {-8,-6,-4,-2,0,2,4,6,8}{
%		\node[] at (\i,\j) {\input{tikz/even_even.tikz}};
%	}
%	\foreach \j in {-7,-5,-3,-1,1,3,5,7}{
%		\node[] at (\i,\j) {\input{tikz/even_odd.tikz}};
%	}
%}
%\foreach \i in {-7,-5,-3,-1,1,3,5,7}{
%	\foreach \j in {-8,-6,-4,-2,0,2,4,6,8}{
%		\node[rotate=90] at (\i,\j) {\input{tikz/even_odd.tikz}};
%	}
%	\foreach \j in {-7,-5,-3,-1,1,3,5,7}{
%		\node[] at (\i,\j) {\input{tikz/odd_odd.tikz}};
%	}
%}

%% file: tikz/simple_split_6.tikz
\coordinate (A0) at (0, 0);
\coordinate (A1) at (0, 1);
\coordinate (A2) at (0, 2);
\coordinate (A3) at (0, 3);
\coordinate (A4) at (0, 4);
\coordinate (A5) at (0, 5);
\coordinate (A6) at (0, 6);
\coordinate (A7) at (0, 7);
\coordinate (A8) at (0, 8);
\coordinate (A9) at (1, 0);
\coordinate (A10) at (1, 1);
\coordinate (A11) at (1, 2);
\coordinate (A12) at (1, 3);
\coordinate (A13) at (1, 4);
\coordinate (A14) at (1, 5);
\coordinate (A15) at (1, 6);
\coordinate (A16) at (1, 7);
\coordinate (A17) at (2, 0);
\coordinate (A18) at (2, 1);
\coordinate (A19) at (2, 2);
\coordinate (A20) at (2, 3);
\coordinate (A21) at (2, 4);
\coordinate (A22) at (2, 5);
\coordinate (A23) at (2, 6);
\coordinate (A24) at (3, 0);
\coordinate (A25) at (3, 1);
\coordinate (A26) at (3, 2);
\coordinate (A27) at (3, 3);
\coordinate (A28) at (3, 4);
\coordinate (A29) at (3, 5);
\coordinate (A30) at (4, 0);
\coordinate (A31) at (4, 1);
\coordinate (A32) at (4, 2);
\coordinate (A33) at (4, 3);
\coordinate (A34) at (4, 4);
\coordinate (A35) at (5, 0);
\coordinate (A36) at (5, 1);
\coordinate (A37) at (5, 2);
\coordinate (A38) at (5, 3);
\coordinate (A39) at (6, 0);
\coordinate (A40) at (6, 1);
\coordinate (A41) at (6, 2);
\coordinate (A42) at (7, 0);
\coordinate (A43) at (7, 1);
\coordinate (A44) at (8, 0);
\coordinate (A45) at (0, -1);
\coordinate (A46) at (0, -2);
\coordinate (A47) at (0, -3);
\coordinate (A48) at (0, -4);
\coordinate (A49) at (0, -5);
\coordinate (A50) at (0, -6);
\coordinate (A51) at (0, -7);
\coordinate (A52) at (0, -8);
\coordinate (A53) at (1, -1);
\coordinate (A54) at (1, -2);
\coordinate (A55) at (1, -3);
\coordinate (A56) at (1, -4);
\coordinate (A57) at (1, -5);
\coordinate (A58) at (1, -6);
\coordinate (A59) at (1, -7);
\coordinate (A60) at (2, -1);
\coordinate (A61) at (2, -2);
\coordinate (A62) at (2, -3);
\coordinate (A63) at (2, -4);
\coordinate (A64) at (2, -5);
\coordinate (A65) at (2, -6);
\coordinate (A66) at (3, -1);
\coordinate (A67) at (3, -2);
\coordinate (A68) at (3, -3);
\coordinate (A69) at (3, -4);
\coordinate (A70) at (3, -5);
\coordinate (A71) at (4, -1);
\coordinate (A72) at (4, -2);
\coordinate (A73) at (4, -3);
\coordinate (A74) at (4, -4);
\coordinate (A75) at (5, -1);
\coordinate (A76) at (5, -2);
\coordinate (A77) at (5, -3);
\coordinate (A78) at (6, -1);
\coordinate (A79) at (6, -2);
\coordinate (A80) at (7, -1);
\coordinate (A81) at (-1, -7);
\coordinate (A82) at (-1, -6);
\coordinate (A83) at (-1, -5);
\coordinate (A84) at (-1, -4);
\coordinate (A85) at (-1, -3);
\coordinate (A86) at (-1, -2);
\coordinate (A87) at (-1, -1);
\coordinate (A88) at (-1, 0);
\coordinate (A89) at (-1, 1);
\coordinate (A90) at (-1, 2);
\coordinate (A91) at (-1, 3);
\coordinate (A92) at (-1, 4);
\coordinate (A93) at (-1, 5);
\coordinate (A94) at (-1, 6);
\coordinate (A95) at (-1, 7);
\coordinate (A96) at (-2, -6);
\coordinate (A97) at (-2, -5);
\coordinate (A98) at (-2, -4);
\coordinate (A99) at (-2, -3);
\coordinate (A100) at (-2, -2);
\coordinate (A101) at (-2, -1);
\coordinate (A102) at (-2, 0);
\coordinate (A103) at (-2, 1);
\coordinate (A104) at (-2, 2);
\coordinate (A105) at (-2, 3);
\coordinate (A106) at (-2, 4);
\coordinate (A107) at (-2, 5);
\coordinate (A108) at (-2, 6);
\coordinate (A109) at (-3, -5);
\coordinate (A110) at (-3, -4);
\coordinate (A111) at (-3, -3);
\coordinate (A112) at (-3, -2);
\coordinate (A113) at (-3, -1);
\coordinate (A114) at (-3, 0);
\coordinate (A115) at (-3, 1);
\coordinate (A116) at (-3, 2);
\coordinate (A117) at (-3, 3);
\coordinate (A118) at (-3, 4);
\coordinate (A119) at (-3, 5);
\coordinate (A120) at (-4, -4);
\coordinate (A121) at (-4, -3);
\coordinate (A122) at (-4, -2);
\coordinate (A123) at (-4, -1);
\coordinate (A124) at (-4, 0);
\coordinate (A125) at (-4, 1);
\coordinate (A126) at (-4, 2);
\coordinate (A127) at (-4, 3);
\coordinate (A128) at (-4, 4);
\coordinate (A129) at (-5, -3);
\coordinate (A130) at (-5, -2);
\coordinate (A131) at (-5, -1);
\coordinate (A132) at (-5, 0);
\coordinate (A133) at (-5, 1);
\coordinate (A134) at (-5, 2);
\coordinate (A135) at (-5, 3);
\coordinate (A136) at (-6, -2);
\coordinate (A137) at (-6, -1);
\coordinate (A138) at (-6, 0);
\coordinate (A139) at (-6, 1);
\coordinate (A140) at (-6, 2);
\coordinate (A141) at (-7, -1);
\coordinate (A142) at (-7, 0);
\coordinate (A143) at (-7, 1);
\coordinate (A144) at (-8, 0);
%\definecolor{color1}{rgb}{0.492091553781074,0.557025354868187,0.534160755376992}
\colorlet{color1}{mycolor4}
\fill[color1!50] (A0) -- (A1) -- (A9) -- cycle;
\fill[color1!50] (A0) -- ($(A0)! 0.5!(A88)$) -- ($(A1)! 0.5!(A88)$) -- (A1) -- cycle;
\fill[color1!50] (A0) -- ($(A0)! 0.5!(A45)$) -- ($(A9)! 0.5!(A45)$) -- (A9) -- cycle;
\fill[color1!50] ($(A45)! 0.5!(A0)$) -- ($(A88)! 0.5!(A0)$) -- (A0) -- cycle;
\fill[color1!50] (A1) -- ($(A1)! 0.5!(A10)$) -- ($(A2)! 0.5!(A10)$) -- (A2) -- cycle;
\fill[color1!50] (A1) -- (A2) -- (A89) -- cycle;
\fill[color1!50] (A1) -- ($(A1)! 0.5!(A10)$) -- ($(A9)! 0.5!(A10)$) -- (A9) -- cycle;
\fill[color1!50] (A1) -- ($(A1)! 0.5!(A88)$) -- ($(A89)! 0.5!(A88)$) -- (A89) -- cycle;
\fill[color1!50] (A2) -- (A3) -- (A11) -- cycle;
\fill[color1!50] (A2) -- ($(A2)! 0.5!(A90)$) -- ($(A3)! 0.5!(A90)$) -- (A3) -- cycle;
\fill[color1!50] (A2) -- ($(A2)! 0.5!(A10)$) -- ($(A11)! 0.5!(A10)$) -- (A11) -- cycle;
\fill[color1!50] (A2) -- ($(A2)! 0.5!(A90)$) -- ($(A89)! 0.5!(A90)$) -- (A89) -- cycle;
\fill[color1!50] (A3) -- ($(A3)! 0.5!(A4)$) -- ($(A12)! 0.5!(A4)$) -- (A12) -- cycle;
\fill[color1!50] ($(A4)! 0.5!(A3)$) -- ($(A91)! 0.5!(A3)$) -- (A3) -- cycle;
\fill[color1!50] (A3) -- (A11) -- (A19) -- cycle;
\fill[color1!50] (A3) -- (A12) -- (A26) -- cycle;
\fill[color1!50] (A3) -- (A19) -- (A42) -- cycle;
\fill[color1!50] (A3) -- (A26) -- (A36) -- cycle;
\fill[color1!50] (A3) -- (A36) -- (A42) -- cycle;
\fill[color1!50] (A3) -- ($(A3)! 0.5!(A90)$) -- ($(A104)! 0.5!(A90)$) -- (A104) -- cycle;
\fill[color1!50] ($(A91)! 0.5!(A3)$) -- ($(A116)! 0.5!(A3)$) -- (A3) -- cycle;
\fill[color1!50] (A3) -- ($(A3)! 0.5!(A142)$) -- ($(A104)! 0.5!(A142)$) -- (A104) -- cycle;
\fill[color1!50] ($(A116)! 0.5!(A3)$) -- ($(A133)! 0.5!(A3)$) -- (A3) -- cycle;
\fill[color1!50] ($(A133)! 0.5!(A3)$) -- ($(A142)! 0.5!(A3)$) -- (A3) -- cycle;
\fill[color1!50] (A5) -- ($(A5)! 0.5!(A4)$) -- ($(A13)! 0.5!(A4)$) -- (A13) -- cycle;
\fill[color1!50] ($(A4)! 0.5!(A5)$) -- ($(A92)! 0.5!(A5)$) -- (A5) -- cycle;
\fill[color1!50] (A12) -- ($(A12)! 0.5!(A4)$) -- ($(A13)! 0.5!(A4)$) -- (A13) -- cycle;
\fill[color1!50] (A5) -- ($(A5)! 0.5!(A6)$) -- ($(A14)! 0.5!(A6)$) -- (A14) -- cycle;
\fill[color1!50] ($(A6)! 0.5!(A5)$) -- ($(A93)! 0.5!(A5)$) -- (A5) -- cycle;
\fill[color1!50] (A5) -- (A13) -- (A14) -- cycle;
\fill[color1!50] ($(A92)! 0.5!(A5)$) -- ($(A93)! 0.5!(A5)$) -- (A5) -- cycle;
\fill[color1!50] (A7) -- ($(A7)! 0.5!(A6)$) -- ($(A15)! 0.5!(A6)$) -- (A15) -- cycle;
\fill[color1!50] ($(A6)! 0.5!(A7)$) -- ($(A94)! 0.5!(A7)$) -- (A7) -- cycle;
\fill[color1!50] (A14) -- ($(A14)! 0.5!(A6)$) -- ($(A15)! 0.5!(A6)$) -- (A15) -- cycle;
\fill[color1!50] (A7) -- ($(A7)! 0.5!(A8)$) -- ($(A16)! 0.5!(A8)$) -- (A16) -- cycle;
\fill[color1!50] ($(A8)! 0.5!(A7)$) -- ($(A95)! 0.5!(A7)$) -- (A7) -- cycle;
\fill[color1!50] (A7) -- (A15) -- (A16) -- cycle;
\fill[color1!50] ($(A94)! 0.5!(A7)$) -- ($(A95)! 0.5!(A7)$) -- (A7) -- cycle;
\fill[color1!50] (A9) -- ($(A9)! 0.5!(A10)$) -- ($(A17)! 0.5!(A10)$) -- (A17) -- cycle;
\fill[color1!50] (A9) -- (A17) -- (A53) -- cycle;
\fill[color1!50] (A9) -- ($(A9)! 0.5!(A45)$) -- ($(A53)! 0.5!(A45)$) -- (A53) -- cycle;
\fill[color1!50] (A11) -- ($(A11)! 0.5!(A10)$) -- ($(A18)! 0.5!(A10)$) -- (A18) -- cycle;
\fill[color1!50] (A17) -- ($(A17)! 0.5!(A10)$) -- ($(A18)! 0.5!(A10)$) -- (A18) -- cycle;
\fill[color1!50] (A11) -- (A18) -- (A19) -- cycle;
\fill[color1!50] (A12) -- (A13) -- (A20) -- cycle;
\fill[color1!50] (A12) -- (A20) -- (A26) -- cycle;
\fill[color1!50] (A13) -- ($(A13)! 0.5!(A21)$) -- ($(A14)! 0.5!(A21)$) -- (A14) -- cycle;
\fill[color1!50] (A13) -- ($(A13)! 0.5!(A21)$) -- ($(A20)! 0.5!(A21)$) -- (A20) -- cycle;
\fill[color1!50] (A14) -- (A15) -- (A22) -- cycle;
\fill[color1!50] (A14) -- ($(A14)! 0.5!(A21)$) -- ($(A22)! 0.5!(A21)$) -- (A22) -- cycle;
\fill[color1!50] (A15) -- ($(A15)! 0.5!(A23)$) -- ($(A16)! 0.5!(A23)$) -- (A16) -- cycle;
\fill[color1!50] (A15) -- ($(A15)! 0.5!(A23)$) -- ($(A22)! 0.5!(A23)$) -- (A22) -- cycle;
\fill[color1!50] (A17) -- (A18) -- (A24) -- cycle;
\fill[color1!50] (A17) -- ($(A17)! 0.5!(A60)$) -- ($(A24)! 0.5!(A60)$) -- (A24) -- cycle;
\fill[color1!50] (A17) -- ($(A17)! 0.5!(A60)$) -- ($(A53)! 0.5!(A60)$) -- (A53) -- cycle;
\fill[color1!50] (A18) -- ($(A18)! 0.5!(A25)$) -- ($(A19)! 0.5!(A25)$) -- (A19) -- cycle;
\fill[color1!50] (A18) -- ($(A18)! 0.5!(A25)$) -- ($(A24)! 0.5!(A25)$) -- (A24) -- cycle;
\fill[color1!50] (A19) -- ($(A19)! 0.5!(A25)$) -- ($(A31)! 0.5!(A25)$) -- (A31) -- cycle;
\fill[color1!50] (A19) -- (A31) -- (A42) -- cycle;
\fill[color1!50] (A20) -- ($(A20)! 0.5!(A21)$) -- ($(A27)! 0.5!(A21)$) -- (A27) -- cycle;
\fill[color1!50] (A20) -- (A26) -- (A27) -- cycle;
\fill[color1!50] (A22) -- ($(A22)! 0.5!(A21)$) -- ($(A28)! 0.5!(A21)$) -- (A28) -- cycle;
\fill[color1!50] (A27) -- ($(A27)! 0.5!(A21)$) -- ($(A28)! 0.5!(A21)$) -- (A28) -- cycle;
\fill[color1!50] (A22) -- ($(A22)! 0.5!(A23)$) -- ($(A29)! 0.5!(A23)$) -- (A29) -- cycle;
\fill[color1!50] (A22) -- (A28) -- (A29) -- cycle;
\fill[color1!50] (A24) -- ($(A24)! 0.5!(A25)$) -- ($(A30)! 0.5!(A25)$) -- (A30) -- cycle;
\fill[color1!50] (A24) -- (A30) -- (A66) -- cycle;
\fill[color1!50] (A24) -- ($(A24)! 0.5!(A60)$) -- ($(A66)! 0.5!(A60)$) -- (A66) -- cycle;
\fill[color1!50] (A30) -- ($(A30)! 0.5!(A25)$) -- ($(A31)! 0.5!(A25)$) -- (A31) -- cycle;
\fill[color1!50] (A26) -- ($(A26)! 0.5!(A32)$) -- ($(A27)! 0.5!(A32)$) -- (A27) -- cycle;
\fill[color1!50] (A26) -- ($(A26)! 0.5!(A32)$) -- ($(A36)! 0.5!(A32)$) -- (A36) -- cycle;
\fill[color1!50] (A27) -- (A28) -- (A33) -- cycle;
\fill[color1!50] (A27) -- ($(A27)! 0.5!(A32)$) -- ($(A33)! 0.5!(A32)$) -- (A33) -- cycle;
\fill[color1!50] (A28) -- ($(A28)! 0.5!(A34)$) -- ($(A29)! 0.5!(A34)$) -- (A29) -- cycle;
\fill[color1!50] (A28) -- ($(A28)! 0.5!(A34)$) -- ($(A33)! 0.5!(A34)$) -- (A33) -- cycle;
\fill[color1!50] (A30) -- (A31) -- (A35) -- cycle;
\fill[color1!50] (A30) -- ($(A30)! 0.5!(A71)$) -- ($(A35)! 0.5!(A71)$) -- (A35) -- cycle;
\fill[color1!50] (A30) -- ($(A30)! 0.5!(A71)$) -- ($(A66)! 0.5!(A71)$) -- (A66) -- cycle;
\fill[color1!50] (A31) -- (A35) -- (A39) -- cycle;
\fill[color1!50] (A31) -- (A39) -- (A42) -- cycle;
\fill[color1!50] (A33) -- ($(A33)! 0.5!(A32)$) -- ($(A37)! 0.5!(A32)$) -- (A37) -- cycle;
\fill[color1!50] (A36) -- ($(A36)! 0.5!(A32)$) -- ($(A37)! 0.5!(A32)$) -- (A37) -- cycle;
\fill[color1!50] (A33) -- ($(A33)! 0.5!(A34)$) -- ($(A38)! 0.5!(A34)$) -- (A38) -- cycle;
\fill[color1!50] (A33) -- (A37) -- (A38) -- cycle;
\fill[color1!50] (A35) -- ($(A35)! 0.5!(A71)$) -- ($(A39)! 0.5!(A71)$) -- (A39) -- cycle;
\fill[color1!50] (A36) -- (A37) -- (A40) -- cycle;
\fill[color1!50] (A36) -- (A40) -- (A42) -- cycle;
\fill[color1!50] (A37) -- ($(A37)! 0.5!(A41)$) -- ($(A38)! 0.5!(A41)$) -- (A38) -- cycle;
\fill[color1!50] (A37) -- ($(A37)! 0.5!(A41)$) -- ($(A40)! 0.5!(A41)$) -- (A40) -- cycle;
\fill[color1!50] (A39) -- ($(A39)! 0.5!(A71)$) -- ($(A42)! 0.5!(A71)$) -- (A42) -- cycle;
\fill[color1!50] (A40) -- ($(A40)! 0.5!(A41)$) -- ($(A43)! 0.5!(A41)$) -- (A43) -- cycle;
\fill[color1!50] (A40) -- (A42) -- (A43) -- cycle;
\fill[color1!50] (A42) -- ($(A42)! 0.5!(A44)$) -- ($(A43)! 0.5!(A44)$) -- (A43) -- cycle;
\fill[color1!50] ($(A44)! 0.5!(A42)$) -- ($(A80)! 0.5!(A42)$) -- (A42) -- cycle;
\fill[color1!50] (A42) -- ($(A42)! 0.5!(A47)$) -- ($(A61)! 0.5!(A47)$) -- (A61) -- cycle;
\fill[color1!50] ($(A47)! 0.5!(A42)$) -- ($(A75)! 0.5!(A42)$) -- (A42) -- cycle;
\fill[color1!50] (A42) -- ($(A42)! 0.5!(A71)$) -- ($(A61)! 0.5!(A71)$) -- (A61) -- cycle;
\fill[color1!50] ($(A75)! 0.5!(A42)$) -- ($(A78)! 0.5!(A42)$) -- (A42) -- cycle;
\fill[color1!50] ($(A78)! 0.5!(A42)$) -- ($(A80)! 0.5!(A42)$) -- (A42) -- cycle;
\fill[color1!50] (A46) -- ($(A46)! 0.5!(A45)$) -- ($(A53)! 0.5!(A45)$) -- (A53) -- cycle;
\fill[color1!50] ($(A45)! 0.5!(A46)$) -- ($(A87)! 0.5!(A46)$) -- (A46) -- cycle;
\fill[color1!50] (A46) -- ($(A46)! 0.5!(A47)$) -- ($(A54)! 0.5!(A47)$) -- (A54) -- cycle;
\fill[color1!50] ($(A47)! 0.5!(A46)$) -- ($(A86)! 0.5!(A46)$) -- (A46) -- cycle;
\fill[color1!50] (A46) -- (A53) -- (A54) -- cycle;
\fill[color1!50] ($(A86)! 0.5!(A46)$) -- ($(A87)! 0.5!(A46)$) -- (A46) -- cycle;
\fill[color1!50] (A54) -- ($(A54)! 0.5!(A47)$) -- ($(A61)! 0.5!(A47)$) -- (A61) -- cycle;
\fill[color1!50] ($(A51)! 0.5!(A81)$) -- ($(A52)! 0.5!(A81)$) -- (A81) -- cycle;
\fill[color1!50] ($(A51)! 0.5!(A81)$) -- ($(A82)! 0.5!(A81)$) -- (A81) -- cycle;
\fill[color1!50] (A53) -- ($(A53)! 0.5!(A60)$) -- ($(A54)! 0.5!(A60)$) -- (A54) -- cycle;
\fill[color1!50] (A54) -- ($(A54)! 0.5!(A60)$) -- ($(A61)! 0.5!(A60)$) -- (A61) -- cycle;
\fill[color1!50] (A61) -- ($(A61)! 0.5!(A60)$) -- ($(A66)! 0.5!(A60)$) -- (A66) -- cycle;
\fill[color1!50] (A61) -- ($(A61)! 0.5!(A71)$) -- ($(A66)! 0.5!(A71)$) -- (A66) -- cycle;
\fill[color1!50] ($(A82)! 0.5!(A81)$) -- ($(A96)! 0.5!(A81)$) -- (A81) -- cycle;
\fill[color1!50] ($(A87)! 0.5!(A102)$) -- ($(A88)! 0.5!(A102)$) -- (A102) -- cycle;
\fill[color1!50] ($(A87)! 0.5!(A102)$) -- ($(A101)! 0.5!(A102)$) -- (A102) -- cycle;
\fill[color1!50] (A89) -- ($(A89)! 0.5!(A88)$) -- ($(A102)! 0.5!(A88)$) -- (A102) -- cycle;
\fill[color1!50] (A89) -- ($(A89)! 0.5!(A90)$) -- ($(A103)! 0.5!(A90)$) -- (A103) -- cycle;
\fill[color1!50] (A89) -- (A102) -- (A103) -- cycle;
\fill[color1!50] (A103) -- ($(A103)! 0.5!(A90)$) -- ($(A104)! 0.5!(A90)$) -- (A104) -- cycle;
\fill[color1!50] ($(A96)! 0.5!(A109)$) -- ($(A97)! 0.5!(A109)$) -- (A109) -- cycle;
\fill[color1!50] ($(A97)! 0.5!(A109)$) -- ($(A110)! 0.5!(A109)$) -- (A109) -- cycle;
\fill[color1!50] ($(A101)! 0.5!(A102)$) -- ($(A114)! 0.5!(A102)$) -- (A102) -- cycle;
\fill[color1!50] (A102) -- ($(A102)! 0.5!(A114)$) -- ($(A103)! 0.5!(A114)$) -- (A103) -- cycle;
\fill[color1!50] (A103) -- (A104) -- (A115) -- cycle;
\fill[color1!50] (A103) -- ($(A103)! 0.5!(A114)$) -- ($(A115)! 0.5!(A114)$) -- (A115) -- cycle;
\fill[color1!50] (A104) -- (A115) -- (A125) -- cycle;
\fill[color1!50] (A104) -- ($(A104)! 0.5!(A142)$) -- ($(A125)! 0.5!(A142)$) -- (A125) -- cycle;
\fill[color1!50] ($(A110)! 0.5!(A109)$) -- ($(A120)! 0.5!(A109)$) -- (A109) -- cycle;
\fill[color1!50] ($(A113)! 0.5!(A124)$) -- ($(A114)! 0.5!(A124)$) -- (A124) -- cycle;
\fill[color1!50] ($(A113)! 0.5!(A124)$) -- ($(A123)! 0.5!(A124)$) -- (A124) -- cycle;
\fill[color1!50] (A115) -- ($(A115)! 0.5!(A114)$) -- ($(A124)! 0.5!(A114)$) -- (A124) -- cycle;
\fill[color1!50] (A115) -- (A124) -- (A125) -- cycle;
\fill[color1!50] ($(A120)! 0.5!(A129)$) -- ($(A121)! 0.5!(A129)$) -- (A129) -- cycle;
\fill[color1!50] ($(A121)! 0.5!(A129)$) -- ($(A130)! 0.5!(A129)$) -- (A129) -- cycle;
\fill[color1!50] ($(A123)! 0.5!(A124)$) -- ($(A132)! 0.5!(A124)$) -- (A124) -- cycle;
\fill[color1!50] ($(A123)! 0.5!(A138)$) -- ($(A132)! 0.5!(A138)$) -- (A138) -- cycle;
\fill[color1!50] ($(A123)! 0.5!(A138)$) -- ($(A142)! 0.5!(A138)$) -- (A138) -- cycle;
\fill[color1!50] (A124) -- ($(A124)! 0.5!(A132)$) -- ($(A125)! 0.5!(A132)$) -- (A125) -- cycle;
\fill[color1!50] (A125) -- ($(A125)! 0.5!(A132)$) -- ($(A138)! 0.5!(A132)$) -- (A138) -- cycle;
\fill[color1!50] (A125) -- ($(A125)! 0.5!(A142)$) -- ($(A138)! 0.5!(A142)$) -- (A138) -- cycle;
\fill[color1!50] ($(A130)! 0.5!(A129)$) -- ($(A136)! 0.5!(A129)$) -- (A129) -- cycle;
\fill[color1!50] ($(A136)! 0.5!(A141)$) -- ($(A137)! 0.5!(A141)$) -- (A141) -- cycle;
\fill[color1!50] ($(A137)! 0.5!(A141)$) -- ($(A142)! 0.5!(A141)$) -- (A141) -- cycle;
\fill[color1!50] ($(A142)! 0.5!(A141)$) -- ($(A144)! 0.5!(A141)$) -- (A141) -- cycle;
\colorlet{color2}{mycolor5}
\fill[color2!50] ($(A1)! 0.5!(A10)$) -- ($(A2)! 0.5!(A10)$) -- (A10) -- cycle;
\fill[color2!50] ($(A1)! 0.5!(A10)$) -- ($(A9)! 0.5!(A10)$) -- (A10) -- cycle;
\fill[color2!50] ($(A2)! 0.5!(A10)$) -- ($(A11)! 0.5!(A10)$) -- (A10) -- cycle;
\fill[color2!50] ($(A9)! 0.5!(A10)$) -- ($(A17)! 0.5!(A10)$) -- (A10) -- cycle;
\fill[color2!50] ($(A11)! 0.5!(A10)$) -- ($(A18)! 0.5!(A10)$) -- (A10) -- cycle;
\fill[color2!50] ($(A17)! 0.5!(A10)$) -- ($(A18)! 0.5!(A10)$) -- (A10) -- cycle;
\colorlet{color3}{mycolor3}
\fill[color3!50] ($(A13)! 0.5!(A21)$) -- ($(A14)! 0.5!(A21)$) -- (A21) -- cycle;
\fill[color3!50] ($(A13)! 0.5!(A21)$) -- ($(A20)! 0.5!(A21)$) -- (A21) -- cycle;
\fill[color3!50] ($(A14)! 0.5!(A21)$) -- ($(A22)! 0.5!(A21)$) -- (A21) -- cycle;
\fill[color3!50] ($(A20)! 0.5!(A21)$) -- ($(A27)! 0.5!(A21)$) -- (A21) -- cycle;
\fill[color3!50] ($(A22)! 0.5!(A21)$) -- ($(A28)! 0.5!(A21)$) -- (A21) -- cycle;
\fill[color3!50] ($(A27)! 0.5!(A21)$) -- ($(A28)! 0.5!(A21)$) -- (A21) -- cycle;
\colorlet{color4}{mycolor5}
\fill[color4!50] ($(A18)! 0.5!(A25)$) -- ($(A19)! 0.5!(A25)$) -- (A25) -- cycle;
\fill[color4!50] ($(A18)! 0.5!(A25)$) -- ($(A24)! 0.5!(A25)$) -- (A25) -- cycle;
\fill[color4!50] ($(A19)! 0.5!(A25)$) -- ($(A31)! 0.5!(A25)$) -- (A25) -- cycle;
\fill[color4!50] ($(A24)! 0.5!(A25)$) -- ($(A30)! 0.5!(A25)$) -- (A25) -- cycle;
\fill[color4!50] ($(A30)! 0.5!(A25)$) -- ($(A31)! 0.5!(A25)$) -- (A25) -- cycle;
\colorlet{color5}{mycolor3}
\fill[color5!50] ($(A26)! 0.5!(A32)$) -- ($(A27)! 0.5!(A32)$) -- (A32) -- cycle;
\fill[color5!50] ($(A26)! 0.5!(A32)$) -- ($(A36)! 0.5!(A32)$) -- (A32) -- cycle;
\fill[color5!50] ($(A27)! 0.5!(A32)$) -- ($(A33)! 0.5!(A32)$) -- (A32) -- cycle;
\fill[color5!50] ($(A33)! 0.5!(A32)$) -- ($(A37)! 0.5!(A32)$) -- (A32) -- cycle;
\fill[color5!50] ($(A36)! 0.5!(A32)$) -- ($(A37)! 0.5!(A32)$) -- (A32) -- cycle;
\colorlet{color6}{mycolor3}
\fill[color6!50] ($(A48)! 0.5!(A56)$) -- ($(A49)! 0.5!(A56)$) -- (A56) -- cycle;
\fill[color6!50] ($(A48)! 0.5!(A56)$) -- ($(A55)! 0.5!(A56)$) -- (A56) -- cycle;
\fill[color6!50] ($(A49)! 0.5!(A56)$) -- ($(A57)! 0.5!(A56)$) -- (A56) -- cycle;
\fill[color6!50] ($(A55)! 0.5!(A56)$) -- ($(A62)! 0.5!(A56)$) -- (A56) -- cycle;
\fill[color6!50] ($(A57)! 0.5!(A56)$) -- ($(A63)! 0.5!(A56)$) -- (A56) -- cycle;
\fill[color6!50] ($(A62)! 0.5!(A56)$) -- ($(A63)! 0.5!(A56)$) -- (A56) -- cycle;
\colorlet{color7}{mycolor3}
\fill[color7!50] ($(A50)! 0.5!(A58)$) -- ($(A51)! 0.5!(A58)$) -- (A58) -- cycle;
\fill[color7!50] ($(A50)! 0.5!(A58)$) -- ($(A57)! 0.5!(A58)$) -- (A58) -- cycle;
\fill[color7!50] ($(A51)! 0.5!(A58)$) -- ($(A59)! 0.5!(A58)$) -- (A58) -- cycle;
\fill[color7!50] ($(A57)! 0.5!(A58)$) -- ($(A64)! 0.5!(A58)$) -- (A58) -- cycle;
\fill[color7!50] ($(A59)! 0.5!(A58)$) -- ($(A65)! 0.5!(A58)$) -- (A58) -- cycle;
\fill[color7!50] ($(A64)! 0.5!(A58)$) -- ($(A65)! 0.5!(A58)$) -- (A58) -- cycle;
\colorlet{color8}{mycolor5}
\fill[color8!50] ($(A17)! 0.5!(A60)$) -- ($(A24)! 0.5!(A60)$) -- (A60) -- cycle;
\fill[color8!50] ($(A17)! 0.5!(A60)$) -- ($(A53)! 0.5!(A60)$) -- (A60) -- cycle;
\fill[color8!50] ($(A24)! 0.5!(A60)$) -- ($(A66)! 0.5!(A60)$) -- (A60) -- cycle;
\fill[color8!50] ($(A53)! 0.5!(A60)$) -- ($(A54)! 0.5!(A60)$) -- (A60) -- cycle;
\fill[color8!50] ($(A54)! 0.5!(A60)$) -- ($(A61)! 0.5!(A60)$) -- (A60) -- cycle;
\fill[color8!50] ($(A61)! 0.5!(A60)$) -- ($(A66)! 0.5!(A60)$) -- (A60) -- cycle;
\colorlet{color9}{mycolor3}
\fill[color9!50] ($(A47)! 0.5!(A67)$) -- ($(A55)! 0.5!(A67)$) -- (A67) -- cycle;
\fill[color9!50] ($(A47)! 0.5!(A67)$) -- ($(A75)! 0.5!(A67)$) -- (A67) -- cycle;
\fill[color9!50] ($(A55)! 0.5!(A67)$) -- ($(A62)! 0.5!(A67)$) -- (A67) -- cycle;
\fill[color9!50] ($(A62)! 0.5!(A67)$) -- ($(A68)! 0.5!(A67)$) -- (A67) -- cycle;
\fill[color9!50] ($(A68)! 0.5!(A67)$) -- ($(A72)! 0.5!(A67)$) -- (A67) -- cycle;
\fill[color9!50] ($(A72)! 0.5!(A67)$) -- ($(A75)! 0.5!(A67)$) -- (A67) -- cycle;
\colorlet{color10}{mycolor3}
\fill[color10!50] ($(A63)! 0.5!(A69)$) -- ($(A64)! 0.5!(A69)$) -- (A69) -- cycle;
\fill[color10!50] ($(A63)! 0.5!(A69)$) -- ($(A68)! 0.5!(A69)$) -- (A69) -- cycle;
\fill[color10!50] ($(A64)! 0.5!(A69)$) -- ($(A70)! 0.5!(A69)$) -- (A69) -- cycle;
\fill[color10!50] ($(A68)! 0.5!(A69)$) -- ($(A73)! 0.5!(A69)$) -- (A69) -- cycle;
\fill[color10!50] ($(A70)! 0.5!(A69)$) -- ($(A74)! 0.5!(A69)$) -- (A69) -- cycle;
\fill[color10!50] ($(A73)! 0.5!(A69)$) -- ($(A74)! 0.5!(A69)$) -- (A69) -- cycle;
\colorlet{color11}{mycolor5}
\fill[color11!50] ($(A30)! 0.5!(A71)$) -- ($(A35)! 0.5!(A71)$) -- (A71) -- cycle;
\fill[color11!50] ($(A30)! 0.5!(A71)$) -- ($(A66)! 0.5!(A71)$) -- (A71) -- cycle;
\fill[color11!50] ($(A35)! 0.5!(A71)$) -- ($(A39)! 0.5!(A71)$) -- (A71) -- cycle;
\fill[color11!50] ($(A39)! 0.5!(A71)$) -- ($(A42)! 0.5!(A71)$) -- (A71) -- cycle;
\fill[color11!50] ($(A42)! 0.5!(A71)$) -- ($(A61)! 0.5!(A71)$) -- (A71) -- cycle;
\fill[color11!50] ($(A61)! 0.5!(A71)$) -- ($(A66)! 0.5!(A71)$) -- (A71) -- cycle;
\colorlet{color12}{mycolor3}
\fill[color12!50] ($(A72)! 0.5!(A76)$) -- ($(A73)! 0.5!(A76)$) -- (A76) -- cycle;
\fill[color12!50] ($(A72)! 0.5!(A76)$) -- ($(A75)! 0.5!(A76)$) -- (A76) -- cycle;
\fill[color12!50] ($(A73)! 0.5!(A76)$) -- ($(A77)! 0.5!(A76)$) -- (A76) -- cycle;
\fill[color12!50] ($(A75)! 0.5!(A76)$) -- ($(A78)! 0.5!(A76)$) -- (A76) -- cycle;
\fill[color12!50] ($(A77)! 0.5!(A76)$) -- ($(A79)! 0.5!(A76)$) -- (A76) -- cycle;
\fill[color12!50] ($(A78)! 0.5!(A76)$) -- ($(A79)! 0.5!(A76)$) -- (A76) -- cycle;
\colorlet{color13}{mycolor3}
\fill[color13!50] ($(A49)! 0.5!(A83)$) -- ($(A50)! 0.5!(A83)$) -- (A83) -- cycle;
\fill[color13!50] ($(A49)! 0.5!(A83)$) -- ($(A84)! 0.5!(A83)$) -- (A83) -- cycle;
\fill[color13!50] ($(A50)! 0.5!(A83)$) -- ($(A82)! 0.5!(A83)$) -- (A83) -- cycle;
\fill[color13!50] ($(A82)! 0.5!(A83)$) -- ($(A97)! 0.5!(A83)$) -- (A83) -- cycle;
\fill[color13!50] ($(A84)! 0.5!(A83)$) -- ($(A98)! 0.5!(A83)$) -- (A83) -- cycle;
\fill[color13!50] ($(A97)! 0.5!(A83)$) -- ($(A98)! 0.5!(A83)$) -- (A83) -- cycle;
\colorlet{color14}{mycolor3}
\fill[color14!50] ($(A47)! 0.5!(A85)$) -- ($(A48)! 0.5!(A85)$) -- (A85) -- cycle;
\fill[color14!50] ($(A47)! 0.5!(A85)$) -- ($(A112)! 0.5!(A85)$) -- (A85) -- cycle;
\fill[color14!50] ($(A48)! 0.5!(A85)$) -- ($(A84)! 0.5!(A85)$) -- (A85) -- cycle;
\fill[color14!50] ($(A84)! 0.5!(A85)$) -- ($(A99)! 0.5!(A85)$) -- (A85) -- cycle;
\fill[color14!50] ($(A99)! 0.5!(A85)$) -- ($(A112)! 0.5!(A85)$) -- (A85) -- cycle;
\colorlet{color15}{mycolor5}
\fill[color15!50] ($(A2)! 0.5!(A90)$) -- ($(A3)! 0.5!(A90)$) -- (A90) -- cycle;
\fill[color15!50] ($(A2)! 0.5!(A90)$) -- ($(A89)! 0.5!(A90)$) -- (A90) -- cycle;
\fill[color15!50] ($(A3)! 0.5!(A90)$) -- ($(A104)! 0.5!(A90)$) -- (A90) -- cycle;
\fill[color15!50] ($(A89)! 0.5!(A90)$) -- ($(A103)! 0.5!(A90)$) -- (A90) -- cycle;
\fill[color15!50] ($(A103)! 0.5!(A90)$) -- ($(A104)! 0.5!(A90)$) -- (A90) -- cycle;
\colorlet{color16}{mycolor5}
\fill[color16!50] ($(A47)! 0.5!(A100)$) -- ($(A86)! 0.5!(A100)$) -- (A100) -- cycle;
\fill[color16!50] ($(A47)! 0.5!(A100)$) -- ($(A142)! 0.5!(A100)$) -- (A100) -- cycle;
\fill[color16!50] ($(A86)! 0.5!(A100)$) -- ($(A101)! 0.5!(A100)$) -- (A100) -- cycle;
\fill[color16!50] ($(A101)! 0.5!(A100)$) -- ($(A113)! 0.5!(A100)$) -- (A100) -- cycle;
\fill[color16!50] ($(A113)! 0.5!(A100)$) -- ($(A123)! 0.5!(A100)$) -- (A100) -- cycle;
\fill[color16!50] ($(A123)! 0.5!(A100)$) -- ($(A142)! 0.5!(A100)$) -- (A100) -- cycle;
\colorlet{color17}{mycolor3}
\fill[color17!50] ($(A91)! 0.5!(A105)$) -- ($(A92)! 0.5!(A105)$) -- (A105) -- cycle;
\fill[color17!50] ($(A91)! 0.5!(A105)$) -- ($(A116)! 0.5!(A105)$) -- (A105) -- cycle;
\fill[color17!50] ($(A92)! 0.5!(A105)$) -- ($(A106)! 0.5!(A105)$) -- (A105) -- cycle;
\fill[color17!50] ($(A106)! 0.5!(A105)$) -- ($(A117)! 0.5!(A105)$) -- (A105) -- cycle;
\fill[color17!50] ($(A116)! 0.5!(A105)$) -- ($(A117)! 0.5!(A105)$) -- (A105) -- cycle;
\colorlet{color18}{mycolor3}
\fill[color18!50] ($(A93)! 0.5!(A107)$) -- ($(A94)! 0.5!(A107)$) -- (A107) -- cycle;
\fill[color18!50] ($(A93)! 0.5!(A107)$) -- ($(A106)! 0.5!(A107)$) -- (A107) -- cycle;
\fill[color18!50] ($(A94)! 0.5!(A107)$) -- ($(A108)! 0.5!(A107)$) -- (A107) -- cycle;
\fill[color18!50] ($(A106)! 0.5!(A107)$) -- ($(A118)! 0.5!(A107)$) -- (A107) -- cycle;
\fill[color18!50] ($(A108)! 0.5!(A107)$) -- ($(A119)! 0.5!(A107)$) -- (A107) -- cycle;
\fill[color18!50] ($(A118)! 0.5!(A107)$) -- ($(A119)! 0.5!(A107)$) -- (A107) -- cycle;
\colorlet{color19}{mycolor3}
\fill[color19!50] ($(A98)! 0.5!(A111)$) -- ($(A99)! 0.5!(A111)$) -- (A111) -- cycle;
\fill[color19!50] ($(A98)! 0.5!(A111)$) -- ($(A110)! 0.5!(A111)$) -- (A111) -- cycle;
\fill[color19!50] ($(A99)! 0.5!(A111)$) -- ($(A112)! 0.5!(A111)$) -- (A111) -- cycle;
\fill[color19!50] ($(A110)! 0.5!(A111)$) -- ($(A121)! 0.5!(A111)$) -- (A111) -- cycle;
\fill[color19!50] ($(A112)! 0.5!(A111)$) -- ($(A122)! 0.5!(A111)$) -- (A111) -- cycle;
\fill[color19!50] ($(A121)! 0.5!(A111)$) -- ($(A122)! 0.5!(A111)$) -- (A111) -- cycle;
\colorlet{color20}{mycolor3}
\fill[color20!50] ($(A117)! 0.5!(A127)$) -- ($(A118)! 0.5!(A127)$) -- (A127) -- cycle;
\fill[color20!50] ($(A117)! 0.5!(A127)$) -- ($(A126)! 0.5!(A127)$) -- (A127) -- cycle;
\fill[color20!50] ($(A118)! 0.5!(A127)$) -- ($(A128)! 0.5!(A127)$) -- (A127) -- cycle;
\fill[color20!50] ($(A126)! 0.5!(A127)$) -- ($(A134)! 0.5!(A127)$) -- (A127) -- cycle;
\fill[color20!50] ($(A128)! 0.5!(A127)$) -- ($(A135)! 0.5!(A127)$) -- (A127) -- cycle;
\fill[color20!50] ($(A134)! 0.5!(A127)$) -- ($(A135)! 0.5!(A127)$) -- (A127) -- cycle;
\colorlet{color21}{mycolor3}
\fill[color21!50] ($(A47)! 0.5!(A131)$) -- ($(A112)! 0.5!(A131)$) -- (A131) -- cycle;
\fill[color21!50] ($(A47)! 0.5!(A131)$) -- ($(A142)! 0.5!(A131)$) -- (A131) -- cycle;
\fill[color21!50] ($(A112)! 0.5!(A131)$) -- ($(A122)! 0.5!(A131)$) -- (A131) -- cycle;
\fill[color21!50] ($(A122)! 0.5!(A131)$) -- ($(A130)! 0.5!(A131)$) -- (A131) -- cycle;
\fill[color21!50] ($(A130)! 0.5!(A131)$) -- ($(A137)! 0.5!(A131)$) -- (A131) -- cycle;
\fill[color21!50] ($(A137)! 0.5!(A131)$) -- ($(A142)! 0.5!(A131)$) -- (A131) -- cycle;
\colorlet{color22}{mycolor3}
\fill[color22!50] ($(A133)! 0.5!(A139)$) -- ($(A134)! 0.5!(A139)$) -- (A139) -- cycle;
\fill[color22!50] ($(A133)! 0.5!(A139)$) -- ($(A142)! 0.5!(A139)$) -- (A139) -- cycle;
\fill[color22!50] ($(A134)! 0.5!(A139)$) -- ($(A140)! 0.5!(A139)$) -- (A139) -- cycle;
\fill[color22!50] ($(A140)! 0.5!(A139)$) -- ($(A143)! 0.5!(A139)$) -- (A139) -- cycle;
\fill[color22!50] ($(A142)! 0.5!(A139)$) -- ($(A143)! 0.5!(A139)$) -- (A139) -- cycle;
\draw[color22, thick] ($(A133)! 0.5!(A139)$) -- ($(A134)! 0.5!(A139)$);
\draw[color22, thick] ($(A133)! 0.5!(A139)$) -- ($(A142)! 0.5!(A139)$);
\draw[color22, thick] ($(A134)! 0.5!(A139)$) -- ($(A140)! 0.5!(A139)$);
\draw[color22, thick] ($(A140)! 0.5!(A139)$) -- ($(A143)! 0.5!(A139)$);
\draw[color22, thick] ($(A142)! 0.5!(A139)$) -- ($(A143)! 0.5!(A139)$);
\draw[color21, thick] ($(A47)! 0.5!(A131)$) -- ($(A112)! 0.5!(A131)$);
\draw[color21, thick] ($(A47)! 0.5!(A131)$) -- ($(A142)! 0.5!(A131)$);
\draw[color21, thick] ($(A112)! 0.5!(A131)$) -- ($(A122)! 0.5!(A131)$);
\draw[color21, thick] ($(A122)! 0.5!(A131)$) -- ($(A130)! 0.5!(A131)$);
\draw[color21, thick] ($(A130)! 0.5!(A131)$) -- ($(A137)! 0.5!(A131)$);
\draw[color21, thick] ($(A137)! 0.5!(A131)$) -- ($(A142)! 0.5!(A131)$);
\draw[color20, thick] ($(A117)! 0.5!(A127)$) -- ($(A118)! 0.5!(A127)$);
\draw[color20, thick] ($(A117)! 0.5!(A127)$) -- ($(A126)! 0.5!(A127)$);
\draw[color20, thick] ($(A118)! 0.5!(A127)$) -- ($(A128)! 0.5!(A127)$);
\draw[color20, thick] ($(A126)! 0.5!(A127)$) -- ($(A134)! 0.5!(A127)$);
\draw[color20, thick] ($(A128)! 0.5!(A127)$) -- ($(A135)! 0.5!(A127)$);
\draw[color20, thick] ($(A134)! 0.5!(A127)$) -- ($(A135)! 0.5!(A127)$);
\draw[color19, thick] ($(A98)! 0.5!(A111)$) -- ($(A99)! 0.5!(A111)$);
\draw[color19, thick] ($(A98)! 0.5!(A111)$) -- ($(A110)! 0.5!(A111)$);
\draw[color19, thick] ($(A99)! 0.5!(A111)$) -- ($(A112)! 0.5!(A111)$);
\draw[color19, thick] ($(A110)! 0.5!(A111)$) -- ($(A121)! 0.5!(A111)$);
\draw[color19, thick] ($(A112)! 0.5!(A111)$) -- ($(A122)! 0.5!(A111)$);
\draw[color19, thick] ($(A121)! 0.5!(A111)$) -- ($(A122)! 0.5!(A111)$);
\draw[color18, thick] ($(A93)! 0.5!(A107)$) -- ($(A94)! 0.5!(A107)$);
\draw[color18, thick] ($(A93)! 0.5!(A107)$) -- ($(A106)! 0.5!(A107)$);
\draw[color18, thick] ($(A94)! 0.5!(A107)$) -- ($(A108)! 0.5!(A107)$);
\draw[color18, thick] ($(A106)! 0.5!(A107)$) -- ($(A118)! 0.5!(A107)$);
\draw[color18, thick] ($(A108)! 0.5!(A107)$) -- ($(A119)! 0.5!(A107)$);
\draw[color18, thick] ($(A118)! 0.5!(A107)$) -- ($(A119)! 0.5!(A107)$);
\draw[color17, thick] ($(A91)! 0.5!(A105)$) -- ($(A92)! 0.5!(A105)$);
\draw[color17, thick] ($(A91)! 0.5!(A105)$) -- ($(A116)! 0.5!(A105)$);
\draw[color17, thick] ($(A92)! 0.5!(A105)$) -- ($(A106)! 0.5!(A105)$);
\draw[color17, thick] ($(A106)! 0.5!(A105)$) -- ($(A117)! 0.5!(A105)$);
\draw[color17, thick] ($(A116)! 0.5!(A105)$) -- ($(A117)! 0.5!(A105)$);
\draw[color16, thick] ($(A47)! 0.5!(A100)$) -- ($(A86)! 0.5!(A100)$);
\draw[color16, thick] ($(A47)! 0.5!(A100)$) -- ($(A142)! 0.5!(A100)$);
\draw[color16, thick] ($(A86)! 0.5!(A100)$) -- ($(A101)! 0.5!(A100)$);
\draw[color16, thick] ($(A101)! 0.5!(A100)$) -- ($(A113)! 0.5!(A100)$);
\draw[color16, thick] ($(A113)! 0.5!(A100)$) -- ($(A123)! 0.5!(A100)$);
\draw[color16, thick] ($(A123)! 0.5!(A100)$) -- ($(A142)! 0.5!(A100)$);
\draw[color14, thick] ($(A47)! 0.5!(A85)$) -- ($(A48)! 0.5!(A85)$);
\draw[color14, thick] ($(A47)! 0.5!(A85)$) -- ($(A112)! 0.5!(A85)$);
\draw[color14, thick] ($(A48)! 0.5!(A85)$) -- ($(A84)! 0.5!(A85)$);
\draw[color14, thick] ($(A84)! 0.5!(A85)$) -- ($(A99)! 0.5!(A85)$);
\draw[color14, thick] ($(A99)! 0.5!(A85)$) -- ($(A112)! 0.5!(A85)$);
\draw[color13, thick] ($(A49)! 0.5!(A83)$) -- ($(A50)! 0.5!(A83)$);
\draw[color13, thick] ($(A49)! 0.5!(A83)$) -- ($(A84)! 0.5!(A83)$);
\draw[color13, thick] ($(A50)! 0.5!(A83)$) -- ($(A82)! 0.5!(A83)$);
\draw[color13, thick] ($(A82)! 0.5!(A83)$) -- ($(A97)! 0.5!(A83)$);
\draw[color13, thick] ($(A84)! 0.5!(A83)$) -- ($(A98)! 0.5!(A83)$);
\draw[color13, thick] ($(A97)! 0.5!(A83)$) -- ($(A98)! 0.5!(A83)$);
\draw[color12, thick] ($(A72)! 0.5!(A76)$) -- ($(A73)! 0.5!(A76)$);
\draw[color12, thick] ($(A72)! 0.5!(A76)$) -- ($(A75)! 0.5!(A76)$);
\draw[color12, thick] ($(A73)! 0.5!(A76)$) -- ($(A77)! 0.5!(A76)$);
\draw[color12, thick] ($(A75)! 0.5!(A76)$) -- ($(A78)! 0.5!(A76)$);
\draw[color12, thick] ($(A77)! 0.5!(A76)$) -- ($(A79)! 0.5!(A76)$);
\draw[color12, thick] ($(A78)! 0.5!(A76)$) -- ($(A79)! 0.5!(A76)$);
\draw[color10, thick] ($(A63)! 0.5!(A69)$) -- ($(A64)! 0.5!(A69)$);
\draw[color10, thick] ($(A63)! 0.5!(A69)$) -- ($(A68)! 0.5!(A69)$);
\draw[color10, thick] ($(A64)! 0.5!(A69)$) -- ($(A70)! 0.5!(A69)$);
\draw[color10, thick] ($(A68)! 0.5!(A69)$) -- ($(A73)! 0.5!(A69)$);
\draw[color10, thick] ($(A70)! 0.5!(A69)$) -- ($(A74)! 0.5!(A69)$);
\draw[color10, thick] ($(A73)! 0.5!(A69)$) -- ($(A74)! 0.5!(A69)$);
\draw[color9, thick] ($(A47)! 0.5!(A67)$) -- ($(A55)! 0.5!(A67)$);
\draw[color9, thick] ($(A47)! 0.5!(A67)$) -- ($(A75)! 0.5!(A67)$);
\draw[color9, thick] ($(A55)! 0.5!(A67)$) -- ($(A62)! 0.5!(A67)$);
\draw[color9, thick] ($(A62)! 0.5!(A67)$) -- ($(A68)! 0.5!(A67)$);
\draw[color9, thick] ($(A68)! 0.5!(A67)$) -- ($(A72)! 0.5!(A67)$);
\draw[color9, thick] ($(A72)! 0.5!(A67)$) -- ($(A75)! 0.5!(A67)$);
\draw[color7, thick] ($(A50)! 0.5!(A58)$) -- ($(A51)! 0.5!(A58)$);
\draw[color7, thick] ($(A50)! 0.5!(A58)$) -- ($(A57)! 0.5!(A58)$);
\draw[color7, thick] ($(A51)! 0.5!(A58)$) -- ($(A59)! 0.5!(A58)$);
\draw[color7, thick] ($(A57)! 0.5!(A58)$) -- ($(A64)! 0.5!(A58)$);
\draw[color7, thick] ($(A59)! 0.5!(A58)$) -- ($(A65)! 0.5!(A58)$);
\draw[color7, thick] ($(A64)! 0.5!(A58)$) -- ($(A65)! 0.5!(A58)$);
\draw[color6, thick] ($(A48)! 0.5!(A56)$) -- ($(A49)! 0.5!(A56)$);
\draw[color6, thick] ($(A48)! 0.5!(A56)$) -- ($(A55)! 0.5!(A56)$);
\draw[color6, thick] ($(A49)! 0.5!(A56)$) -- ($(A57)! 0.5!(A56)$);
\draw[color6, thick] ($(A55)! 0.5!(A56)$) -- ($(A62)! 0.5!(A56)$);
\draw[color6, thick] ($(A57)! 0.5!(A56)$) -- ($(A63)! 0.5!(A56)$);
\draw[color6, thick] ($(A62)! 0.5!(A56)$) -- ($(A63)! 0.5!(A56)$);
\draw[color1, thick] ($(A0)! 0.5!(A88)$) -- ($(A1)! 0.5!(A88)$);
\draw[color1, thick] ($(A0)! 0.5!(A45)$) -- ($(A9)! 0.5!(A45)$);
\draw[color1, thick] ($(A45)! 0.5!(A0)$) -- ($(A88)! 0.5!(A0)$);
\draw[color1, thick] ($(A1)! 0.5!(A10)$) -- ($(A2)! 0.5!(A10)$);
\draw[color1, thick] ($(A1)! 0.5!(A10)$) -- ($(A9)! 0.5!(A10)$);
\draw[color1, thick] ($(A1)! 0.5!(A88)$) -- ($(A89)! 0.5!(A88)$);
\draw[color1, thick] ($(A2)! 0.5!(A90)$) -- ($(A3)! 0.5!(A90)$);
\draw[color1, thick] ($(A2)! 0.5!(A10)$) -- ($(A11)! 0.5!(A10)$);
\draw[color1, thick] ($(A2)! 0.5!(A90)$) -- ($(A89)! 0.5!(A90)$);
\draw[color1, thick] ($(A3)! 0.5!(A4)$) -- ($(A12)! 0.5!(A4)$);
\draw[color1, thick] ($(A4)! 0.5!(A3)$) -- ($(A91)! 0.5!(A3)$);
\draw[color1, thick] ($(A3)! 0.5!(A90)$) -- ($(A104)! 0.5!(A90)$);
\draw[color1, thick] ($(A91)! 0.5!(A3)$) -- ($(A116)! 0.5!(A3)$);
\draw[color1, thick] ($(A3)! 0.5!(A142)$) -- ($(A104)! 0.5!(A142)$);
\draw[color1, thick] ($(A116)! 0.5!(A3)$) -- ($(A133)! 0.5!(A3)$);
\draw[color1, thick] ($(A133)! 0.5!(A3)$) -- ($(A142)! 0.5!(A3)$);
\draw[color1, thick] ($(A5)! 0.5!(A4)$) -- ($(A13)! 0.5!(A4)$);
\draw[color1, thick] ($(A4)! 0.5!(A5)$) -- ($(A92)! 0.5!(A5)$);
\draw[color1, thick] ($(A12)! 0.5!(A4)$) -- ($(A13)! 0.5!(A4)$);
\draw[color1, thick] ($(A5)! 0.5!(A6)$) -- ($(A14)! 0.5!(A6)$);
\draw[color1, thick] ($(A6)! 0.5!(A5)$) -- ($(A93)! 0.5!(A5)$);
\draw[color1, thick] ($(A92)! 0.5!(A5)$) -- ($(A93)! 0.5!(A5)$);
\draw[color1, thick] ($(A7)! 0.5!(A6)$) -- ($(A15)! 0.5!(A6)$);
\draw[color1, thick] ($(A6)! 0.5!(A7)$) -- ($(A94)! 0.5!(A7)$);
\draw[color1, thick] ($(A14)! 0.5!(A6)$) -- ($(A15)! 0.5!(A6)$);
\draw[color1, thick] ($(A7)! 0.5!(A8)$) -- ($(A16)! 0.5!(A8)$);
\draw[color1, thick] ($(A8)! 0.5!(A7)$) -- ($(A95)! 0.5!(A7)$);
\draw[color1, thick] ($(A94)! 0.5!(A7)$) -- ($(A95)! 0.5!(A7)$);
\draw[color1, thick] ($(A9)! 0.5!(A10)$) -- ($(A17)! 0.5!(A10)$);
\draw[color1, thick] ($(A9)! 0.5!(A45)$) -- ($(A53)! 0.5!(A45)$);
\draw[color1, thick] ($(A11)! 0.5!(A10)$) -- ($(A18)! 0.5!(A10)$);
\draw[color1, thick] ($(A17)! 0.5!(A10)$) -- ($(A18)! 0.5!(A10)$);
\draw[color1, thick] ($(A13)! 0.5!(A21)$) -- ($(A14)! 0.5!(A21)$);
\draw[color1, thick] ($(A13)! 0.5!(A21)$) -- ($(A20)! 0.5!(A21)$);
\draw[color1, thick] ($(A14)! 0.5!(A21)$) -- ($(A22)! 0.5!(A21)$);
\draw[color1, thick] ($(A15)! 0.5!(A23)$) -- ($(A16)! 0.5!(A23)$);
\draw[color1, thick] ($(A15)! 0.5!(A23)$) -- ($(A22)! 0.5!(A23)$);
\draw[color1, thick] ($(A17)! 0.5!(A60)$) -- ($(A24)! 0.5!(A60)$);
\draw[color1, thick] ($(A17)! 0.5!(A60)$) -- ($(A53)! 0.5!(A60)$);
\draw[color1, thick] ($(A18)! 0.5!(A25)$) -- ($(A19)! 0.5!(A25)$);
\draw[color1, thick] ($(A18)! 0.5!(A25)$) -- ($(A24)! 0.5!(A25)$);
\draw[color1, thick] ($(A19)! 0.5!(A25)$) -- ($(A31)! 0.5!(A25)$);
\draw[color1, thick] ($(A20)! 0.5!(A21)$) -- ($(A27)! 0.5!(A21)$);
\draw[color1, thick] ($(A22)! 0.5!(A21)$) -- ($(A28)! 0.5!(A21)$);
\draw[color1, thick] ($(A27)! 0.5!(A21)$) -- ($(A28)! 0.5!(A21)$);
\draw[color1, thick] ($(A22)! 0.5!(A23)$) -- ($(A29)! 0.5!(A23)$);
\draw[color1, thick] ($(A24)! 0.5!(A25)$) -- ($(A30)! 0.5!(A25)$);
\draw[color1, thick] ($(A24)! 0.5!(A60)$) -- ($(A66)! 0.5!(A60)$);
\draw[color1, thick] ($(A30)! 0.5!(A25)$) -- ($(A31)! 0.5!(A25)$);
\draw[color1, thick] ($(A26)! 0.5!(A32)$) -- ($(A27)! 0.5!(A32)$);
\draw[color1, thick] ($(A26)! 0.5!(A32)$) -- ($(A36)! 0.5!(A32)$);
\draw[color1, thick] ($(A27)! 0.5!(A32)$) -- ($(A33)! 0.5!(A32)$);
\draw[color1, thick] ($(A28)! 0.5!(A34)$) -- ($(A29)! 0.5!(A34)$);
\draw[color1, thick] ($(A28)! 0.5!(A34)$) -- ($(A33)! 0.5!(A34)$);
\draw[color1, thick] ($(A30)! 0.5!(A71)$) -- ($(A35)! 0.5!(A71)$);
\draw[color1, thick] ($(A30)! 0.5!(A71)$) -- ($(A66)! 0.5!(A71)$);
\draw[color1, thick] ($(A33)! 0.5!(A32)$) -- ($(A37)! 0.5!(A32)$);
\draw[color1, thick] ($(A36)! 0.5!(A32)$) -- ($(A37)! 0.5!(A32)$);
\draw[color1, thick] ($(A33)! 0.5!(A34)$) -- ($(A38)! 0.5!(A34)$);
\draw[color1, thick] ($(A35)! 0.5!(A71)$) -- ($(A39)! 0.5!(A71)$);
\draw[color1, thick] ($(A37)! 0.5!(A41)$) -- ($(A38)! 0.5!(A41)$);
\draw[color1, thick] ($(A37)! 0.5!(A41)$) -- ($(A40)! 0.5!(A41)$);
\draw[color1, thick] ($(A39)! 0.5!(A71)$) -- ($(A42)! 0.5!(A71)$);
\draw[color1, thick] ($(A40)! 0.5!(A41)$) -- ($(A43)! 0.5!(A41)$);
\draw[color1, thick] ($(A42)! 0.5!(A44)$) -- ($(A43)! 0.5!(A44)$);
\draw[color1, thick] ($(A44)! 0.5!(A42)$) -- ($(A80)! 0.5!(A42)$);
\draw[color1, thick] ($(A42)! 0.5!(A47)$) -- ($(A61)! 0.5!(A47)$);
\draw[color1, thick] ($(A47)! 0.5!(A42)$) -- ($(A75)! 0.5!(A42)$);
\draw[color1, thick] ($(A42)! 0.5!(A71)$) -- ($(A61)! 0.5!(A71)$);
\draw[color1, thick] ($(A75)! 0.5!(A42)$) -- ($(A78)! 0.5!(A42)$);
\draw[color1, thick] ($(A78)! 0.5!(A42)$) -- ($(A80)! 0.5!(A42)$);
\draw[color1, thick] ($(A46)! 0.5!(A45)$) -- ($(A53)! 0.5!(A45)$);
\draw[color1, thick] ($(A45)! 0.5!(A46)$) -- ($(A87)! 0.5!(A46)$);
\draw[color1, thick] ($(A46)! 0.5!(A47)$) -- ($(A54)! 0.5!(A47)$);
\draw[color1, thick] ($(A47)! 0.5!(A46)$) -- ($(A86)! 0.5!(A46)$);
\draw[color1, thick] ($(A86)! 0.5!(A46)$) -- ($(A87)! 0.5!(A46)$);
\draw[color1, thick] ($(A54)! 0.5!(A47)$) -- ($(A61)! 0.5!(A47)$);
\draw[color1, thick] ($(A51)! 0.5!(A81)$) -- ($(A52)! 0.5!(A81)$);
\draw[color1, thick] ($(A51)! 0.5!(A81)$) -- ($(A82)! 0.5!(A81)$);
\draw[color1, thick] ($(A53)! 0.5!(A60)$) -- ($(A54)! 0.5!(A60)$);
\draw[color1, thick] ($(A54)! 0.5!(A60)$) -- ($(A61)! 0.5!(A60)$);
\draw[color1, thick] ($(A61)! 0.5!(A60)$) -- ($(A66)! 0.5!(A60)$);
\draw[color1, thick] ($(A61)! 0.5!(A71)$) -- ($(A66)! 0.5!(A71)$);
\draw[color1, thick] ($(A82)! 0.5!(A81)$) -- ($(A96)! 0.5!(A81)$);
\draw[color1, thick] ($(A87)! 0.5!(A102)$) -- ($(A88)! 0.5!(A102)$);
\draw[color1, thick] ($(A87)! 0.5!(A102)$) -- ($(A101)! 0.5!(A102)$);
\draw[color1, thick] ($(A89)! 0.5!(A88)$) -- ($(A102)! 0.5!(A88)$);
\draw[color1, thick] ($(A89)! 0.5!(A90)$) -- ($(A103)! 0.5!(A90)$);
\draw[color1, thick] ($(A103)! 0.5!(A90)$) -- ($(A104)! 0.5!(A90)$);
\draw[color1, thick] ($(A96)! 0.5!(A109)$) -- ($(A97)! 0.5!(A109)$);
\draw[color1, thick] ($(A97)! 0.5!(A109)$) -- ($(A110)! 0.5!(A109)$);
\draw[color1, thick] ($(A101)! 0.5!(A102)$) -- ($(A114)! 0.5!(A102)$);
\draw[color1, thick] ($(A102)! 0.5!(A114)$) -- ($(A103)! 0.5!(A114)$);
\draw[color1, thick] ($(A103)! 0.5!(A114)$) -- ($(A115)! 0.5!(A114)$);
\draw[color1, thick] ($(A104)! 0.5!(A142)$) -- ($(A125)! 0.5!(A142)$);
\draw[color1, thick] ($(A110)! 0.5!(A109)$) -- ($(A120)! 0.5!(A109)$);
\draw[color1, thick] ($(A113)! 0.5!(A124)$) -- ($(A114)! 0.5!(A124)$);
\draw[color1, thick] ($(A113)! 0.5!(A124)$) -- ($(A123)! 0.5!(A124)$);
\draw[color1, thick] ($(A115)! 0.5!(A114)$) -- ($(A124)! 0.5!(A114)$);
\draw[color1, thick] ($(A120)! 0.5!(A129)$) -- ($(A121)! 0.5!(A129)$);
\draw[color1, thick] ($(A121)! 0.5!(A129)$) -- ($(A130)! 0.5!(A129)$);
\draw[color1, thick] ($(A123)! 0.5!(A124)$) -- ($(A132)! 0.5!(A124)$);
\draw[color1, thick] ($(A123)! 0.5!(A138)$) -- ($(A132)! 0.5!(A138)$);
\draw[color1, thick] ($(A123)! 0.5!(A138)$) -- ($(A142)! 0.5!(A138)$);
\draw[color1, thick] ($(A124)! 0.5!(A132)$) -- ($(A125)! 0.5!(A132)$);
\draw[color1, thick] ($(A125)! 0.5!(A132)$) -- ($(A138)! 0.5!(A132)$);
\draw[color1, thick] ($(A125)! 0.5!(A142)$) -- ($(A138)! 0.5!(A142)$);
\draw[color1, thick] ($(A130)! 0.5!(A129)$) -- ($(A136)! 0.5!(A129)$);
\draw[color1, thick] ($(A136)! 0.5!(A141)$) -- ($(A137)! 0.5!(A141)$);
\draw[color1, thick] ($(A137)! 0.5!(A141)$) -- ($(A142)! 0.5!(A141)$);
\draw[color1, thick] ($(A142)! 0.5!(A141)$) -- ($(A144)! 0.5!(A141)$);
\draw[color15, thick] ($(A2)! 0.5!(A90)$) -- ($(A3)! 0.5!(A90)$);
\draw[color15, thick] ($(A2)! 0.5!(A90)$) -- ($(A89)! 0.5!(A90)$);
\draw[color15, thick] ($(A3)! 0.5!(A90)$) -- ($(A104)! 0.5!(A90)$);
\draw[color15, thick] ($(A89)! 0.5!(A90)$) -- ($(A103)! 0.5!(A90)$);
\draw[color15, thick] ($(A103)! 0.5!(A90)$) -- ($(A104)! 0.5!(A90)$);
\draw[color11, thick] ($(A30)! 0.5!(A71)$) -- ($(A35)! 0.5!(A71)$);
\draw[color11, thick] ($(A30)! 0.5!(A71)$) -- ($(A66)! 0.5!(A71)$);
\draw[color11, thick] ($(A35)! 0.5!(A71)$) -- ($(A39)! 0.5!(A71)$);
\draw[color11, thick] ($(A39)! 0.5!(A71)$) -- ($(A42)! 0.5!(A71)$);
\draw[color11, thick] ($(A42)! 0.5!(A71)$) -- ($(A61)! 0.5!(A71)$);
\draw[color11, thick] ($(A61)! 0.5!(A71)$) -- ($(A66)! 0.5!(A71)$);
\draw[color8, thick] ($(A17)! 0.5!(A60)$) -- ($(A24)! 0.5!(A60)$);
\draw[color8, thick] ($(A17)! 0.5!(A60)$) -- ($(A53)! 0.5!(A60)$);
\draw[color8, thick] ($(A24)! 0.5!(A60)$) -- ($(A66)! 0.5!(A60)$);
\draw[color8, thick] ($(A53)! 0.5!(A60)$) -- ($(A54)! 0.5!(A60)$);
\draw[color8, thick] ($(A54)! 0.5!(A60)$) -- ($(A61)! 0.5!(A60)$);
\draw[color8, thick] ($(A61)! 0.5!(A60)$) -- ($(A66)! 0.5!(A60)$);
\draw[color5, thick] ($(A26)! 0.5!(A32)$) -- ($(A27)! 0.5!(A32)$);
\draw[color5, thick] ($(A26)! 0.5!(A32)$) -- ($(A36)! 0.5!(A32)$);
\draw[color5, thick] ($(A27)! 0.5!(A32)$) -- ($(A33)! 0.5!(A32)$);
\draw[color5, thick] ($(A33)! 0.5!(A32)$) -- ($(A37)! 0.5!(A32)$);
\draw[color5, thick] ($(A36)! 0.5!(A32)$) -- ($(A37)! 0.5!(A32)$);
\draw[color4, thick] ($(A18)! 0.5!(A25)$) -- ($(A19)! 0.5!(A25)$);
\draw[color4, thick] ($(A18)! 0.5!(A25)$) -- ($(A24)! 0.5!(A25)$);
\draw[color4, thick] ($(A19)! 0.5!(A25)$) -- ($(A31)! 0.5!(A25)$);
\draw[color4, thick] ($(A24)! 0.5!(A25)$) -- ($(A30)! 0.5!(A25)$);
\draw[color4, thick] ($(A30)! 0.5!(A25)$) -- ($(A31)! 0.5!(A25)$);
\draw[color3, thick] ($(A13)! 0.5!(A21)$) -- ($(A14)! 0.5!(A21)$);
\draw[color3, thick] ($(A13)! 0.5!(A21)$) -- ($(A20)! 0.5!(A21)$);
\draw[color3, thick] ($(A14)! 0.5!(A21)$) -- ($(A22)! 0.5!(A21)$);
\draw[color3, thick] ($(A20)! 0.5!(A21)$) -- ($(A27)! 0.5!(A21)$);
\draw[color3, thick] ($(A22)! 0.5!(A21)$) -- ($(A28)! 0.5!(A21)$);
\draw[color3, thick] ($(A27)! 0.5!(A21)$) -- ($(A28)! 0.5!(A21)$);
\draw[color2, thick] ($(A1)! 0.5!(A10)$) -- ($(A2)! 0.5!(A10)$);
\draw[color2, thick] ($(A1)! 0.5!(A10)$) -- ($(A9)! 0.5!(A10)$);
\draw[color2, thick] ($(A2)! 0.5!(A10)$) -- ($(A11)! 0.5!(A10)$);
\draw[color2, thick] ($(A9)! 0.5!(A10)$) -- ($(A17)! 0.5!(A10)$);
\draw[color2, thick] ($(A11)! 0.5!(A10)$) -- ($(A18)! 0.5!(A10)$);
\draw[color2, thick] ($(A17)! 0.5!(A10)$) -- ($(A18)! 0.5!(A10)$);
\foreach \a/\b/\c in {0/1/9,0/1/88,0/9/45,0/45/88,1/2/10,1/2/89,1/9/10,1/88/89,2/3/11,2/3/90,2/10/11,2/89/90,3/4/12,3/4/91,3/11/19,3/12/26,3/19/42,3/26/36,3/36/42,3/90/104,3/91/116,3/104/142,3/116/133,3/133/142,4/5/13,4/5/92,4/12/13,4/91/92,5/6/14,5/6/93,5/13/14,5/92/93,6/7/15,6/7/94,6/14/15,6/93/94,7/8/16,7/8/95,7/15/16,7/94/95,9/10/17,9/17/53,9/45/53,10/11/18,10/17/18,11/18/19,12/13/20,12/20/26,13/14/21,13/20/21,14/15/22,14/21/22,15/16/23,15/22/23,17/18/24,17/24/60,17/53/60,18/19/25,18/24/25,19/25/31,19/31/42,20/21/27,20/26/27,21/22/28,21/27/28,22/23/29,22/28/29,24/25/30,24/30/66,24/60/66,25/30/31,26/27/32,26/32/36,27/28/33,27/32/33,28/29/34,28/33/34,30/31/35,30/35/71,30/66/71,31/35/39,31/39/42,32/33/37,32/36/37,33/34/38,33/37/38,35/39/71,36/37/40,36/40/42,37/38/41,37/40/41,39/42/71,40/41/43,40/42/43,42/43/44,42/44/80,42/47/61,42/47/75,42/61/71,42/75/78,42/78/80,45/46/53,45/46/87,45/87/88,46/47/54,46/47/86,46/53/54,46/86/87,47/48/55,47/48/85,47/54/61,47/55/67,47/67/75,47/85/112,47/86/100,47/100/142,47/112/131,47/131/142,48/49/56,48/49/84,48/55/56,48/84/85,49/50/57,49/50/83,49/56/57,49/83/84,50/51/58,50/51/82,50/57/58,50/82/83,51/52/59,51/52/81,51/58/59,51/81/82,53/54/60,54/60/61,55/56/62,55/62/67,56/57/63,56/62/63,57/58/64,57/63/64,58/59/65,58/64/65,60/61/66,61/66/71,62/63/68,62/67/68,63/64/69,63/68/69,64/65/70,64/69/70,67/68/72,67/72/75,68/69/73,68/72/73,69/70/74,69/73/74,72/73/76,72/75/76,73/74/77,73/76/77,75/76/78,76/77/79,76/78/79,78/79/80,81/82/96,82/83/97,82/96/97,83/84/98,83/97/98,84/85/99,84/98/99,85/99/112,86/87/101,86/100/101,87/88/102,87/101/102,88/89/102,89/90/103,89/102/103,90/103/104,91/92/105,91/105/116,92/93/106,92/105/106,93/94/107,93/106/107,94/95/108,94/107/108,96/97/109,97/98/110,97/109/110,98/99/111,98/110/111,99/111/112,100/101/113,100/113/123,100/123/142,101/102/114,101/113/114,102/103/114,103/104/115,103/114/115,104/115/125,104/125/142,105/106/117,105/116/117,106/107/118,106/117/118,107/108/119,107/118/119,109/110/120,110/111/121,110/120/121,111/112/122,111/121/122,112/122/131,113/114/124,113/123/124,114/115/124,115/124/125,116/117/126,116/126/133,117/118/127,117/126/127,118/119/128,118/127/128,120/121/129,121/122/130,121/129/130,122/130/131,123/124/132,123/132/138,123/138/142,124/125/132,125/132/138,125/138/142,126/127/134,126/133/134,127/128/135,127/134/135,129/130/136,130/131/137,130/136/137,131/137/142,133/134/139,133/139/142,134/135/140,134/139/140,136/137/141,137/141/142,139/140/143,139/142/143,141/142/144,142/143/144}{
  \draw[black!30] (A\a) -- (A\b) -- (A\c) -- cycle;
}

\draw[black, thick] (A3) -- (A42) -- (A47) -- (A142) --cycle;

\foreach \i in {0,1,2,3,5,7,9,11,12,13,14,15,16,17,18,19,20,22,24,26,27,28,29,30,31,33,35,36,37,38,39,40,42,43,46,53,54,56,58,61,66,67,69,76,81,83,85,89,100,102,103,104,105,107,109,111,115,124,125,127,129,131,138,139,141}{
  \fill[myred] (A\i) circle (3pt);
}
\foreach \i in {4,6,8,10,21,23,25,32,34,41,44,45,47,48,49,50,51,52,55,57,59,60,62,63,64,65,68,70,71,72,73,74,75,77,78,79,80,82,84,86,87,88,90,91,92,93,94,95,96,97,98,99,101,106,108,110,112,113,114,116,117,118,119,120,121,122,123,126,128,130,132,133,134,135,136,137,140,142,143,144}{
  \fill[myblue] (A\i) circle (3pt);
}

%\foreach \i in {-8,-6,-4,-2,0,2,4,6,8}{
%	\foreach \j in {-8,-6,-4,-2,0,2,4,6,8}{
%		\node[] at (\i,\j) {\input{tikz/even_even.tikz}};
%	}
%	\foreach \j in {-7,-5,-3,-1,1,3,5,7}{
%		\node[] at (\i,\j) {\input{tikz/even_odd.tikz}};
%	}
%}
%\foreach \i in {-7,-5,-3,-1,1,3,5,7}{
%	\foreach \j in {-8,-6,-4,-2,0,2,4,6,8}{
%		\node[rotate=90] at (\i,\j) {\input{tikz/even_odd.tikz}};
%	}
%	\foreach \j in {-7,-5,-3,-1,1,3,5,7}{
%		\node[] at (\i,\j) {\input{tikz/odd_odd.tikz}};
%	}
%}

%% file: tikz/simple_split_12.tikz
\coordinate (A0) at (0, 0);
\coordinate (A1) at (0, 1);
\coordinate (A2) at (0, 2);
\coordinate (A3) at (0, 3);
\coordinate (A4) at (0, 4);
\coordinate (A5) at (0, 5);
\coordinate (A6) at (0, 6);
\coordinate (A7) at (0, 7);
\coordinate (A8) at (0, 8);
\coordinate (A9) at (1, 0);
\coordinate (A10) at (1, 1);
\coordinate (A11) at (1, 2);
\coordinate (A12) at (1, 3);
\coordinate (A13) at (1, 4);
\coordinate (A14) at (1, 5);
\coordinate (A15) at (1, 6);
\coordinate (A16) at (1, 7);
\coordinate (A17) at (2, 0);
\coordinate (A18) at (2, 1);
\coordinate (A19) at (2, 2);
\coordinate (A20) at (2, 3);
\coordinate (A21) at (2, 4);
\coordinate (A22) at (2, 5);
\coordinate (A23) at (2, 6);
\coordinate (A24) at (3, 0);
\coordinate (A25) at (3, 1);
\coordinate (A26) at (3, 2);
\coordinate (A27) at (3, 3);
\coordinate (A28) at (3, 4);
\coordinate (A29) at (3, 5);
\coordinate (A30) at (4, 0);
\coordinate (A31) at (4, 1);
\coordinate (A32) at (4, 2);
\coordinate (A33) at (4, 3);
\coordinate (A34) at (4, 4);
\coordinate (A35) at (5, 0);
\coordinate (A36) at (5, 1);
\coordinate (A37) at (5, 2);
\coordinate (A38) at (5, 3);
\coordinate (A39) at (6, 0);
\coordinate (A40) at (6, 1);
\coordinate (A41) at (6, 2);
\coordinate (A42) at (7, 0);
\coordinate (A43) at (7, 1);
\coordinate (A44) at (8, 0);
\coordinate (A45) at (0, -1);
\coordinate (A46) at (0, -2);
\coordinate (A47) at (0, -3);
\coordinate (A48) at (0, -4);
\coordinate (A49) at (0, -5);
\coordinate (A50) at (0, -6);
\coordinate (A51) at (0, -7);
\coordinate (A52) at (0, -8);
\coordinate (A53) at (1, -1);
\coordinate (A54) at (1, -2);
\coordinate (A55) at (1, -3);
\coordinate (A56) at (1, -4);
\coordinate (A57) at (1, -5);
\coordinate (A58) at (1, -6);
\coordinate (A59) at (1, -7);
\coordinate (A60) at (2, -1);
\coordinate (A61) at (2, -2);
\coordinate (A62) at (2, -3);
\coordinate (A63) at (2, -4);
\coordinate (A64) at (2, -5);
\coordinate (A65) at (2, -6);
\coordinate (A66) at (3, -1);
\coordinate (A67) at (3, -2);
\coordinate (A68) at (3, -3);
\coordinate (A69) at (3, -4);
\coordinate (A70) at (3, -5);
\coordinate (A71) at (4, -1);
\coordinate (A72) at (4, -2);
\coordinate (A73) at (4, -3);
\coordinate (A74) at (4, -4);
\coordinate (A75) at (5, -1);
\coordinate (A76) at (5, -2);
\coordinate (A77) at (5, -3);
\coordinate (A78) at (6, -1);
\coordinate (A79) at (6, -2);
\coordinate (A80) at (7, -1);
\coordinate (A81) at (-1, -7);
\coordinate (A82) at (-1, -6);
\coordinate (A83) at (-1, -5);
\coordinate (A84) at (-1, -4);
\coordinate (A85) at (-1, -3);
\coordinate (A86) at (-1, -2);
\coordinate (A87) at (-1, -1);
\coordinate (A88) at (-1, 0);
\coordinate (A89) at (-1, 1);
\coordinate (A90) at (-1, 2);
\coordinate (A91) at (-1, 3);
\coordinate (A92) at (-1, 4);
\coordinate (A93) at (-1, 5);
\coordinate (A94) at (-1, 6);
\coordinate (A95) at (-1, 7);
\coordinate (A96) at (-2, -6);
\coordinate (A97) at (-2, -5);
\coordinate (A98) at (-2, -4);
\coordinate (A99) at (-2, -3);
\coordinate (A100) at (-2, -2);
\coordinate (A101) at (-2, -1);
\coordinate (A102) at (-2, 0);
\coordinate (A103) at (-2, 1);
\coordinate (A104) at (-2, 2);
\coordinate (A105) at (-2, 3);
\coordinate (A106) at (-2, 4);
\coordinate (A107) at (-2, 5);
\coordinate (A108) at (-2, 6);
\coordinate (A109) at (-3, -5);
\coordinate (A110) at (-3, -4);
\coordinate (A111) at (-3, -3);
\coordinate (A112) at (-3, -2);
\coordinate (A113) at (-3, -1);
\coordinate (A114) at (-3, 0);
\coordinate (A115) at (-3, 1);
\coordinate (A116) at (-3, 2);
\coordinate (A117) at (-3, 3);
\coordinate (A118) at (-3, 4);
\coordinate (A119) at (-3, 5);
\coordinate (A120) at (-4, -4);
\coordinate (A121) at (-4, -3);
\coordinate (A122) at (-4, -2);
\coordinate (A123) at (-4, -1);
\coordinate (A124) at (-4, 0);
\coordinate (A125) at (-4, 1);
\coordinate (A126) at (-4, 2);
\coordinate (A127) at (-4, 3);
\coordinate (A128) at (-4, 4);
\coordinate (A129) at (-5, -3);
\coordinate (A130) at (-5, -2);
\coordinate (A131) at (-5, -1);
\coordinate (A132) at (-5, 0);
\coordinate (A133) at (-5, 1);
\coordinate (A134) at (-5, 2);
\coordinate (A135) at (-5, 3);
\coordinate (A136) at (-6, -2);
\coordinate (A137) at (-6, -1);
\coordinate (A138) at (-6, 0);
\coordinate (A139) at (-6, 1);
\coordinate (A140) at (-6, 2);
\coordinate (A141) at (-7, -1);
\coordinate (A142) at (-7, 0);
\coordinate (A143) at (-7, 1);
\coordinate (A144) at (-8, 0);
%\definecolor{color1}{rgb}{0.488270119170684,0.630944969743452,0.795622715898143}
\colorlet{color1}{mycolor4}
\fill[color1!50] (A0) -- (A1) -- (A9) -- cycle;
\fill[color1!50] (A0) -- ($(A0)! 0.5!(A88)$) -- ($(A1)! 0.5!(A88)$) -- (A1) -- cycle;
\fill[color1!50] (A0) -- ($(A0)! 0.5!(A45)$) -- ($(A9)! 0.5!(A45)$) -- (A9) -- cycle;
\fill[color1!50] ($(A45)! 0.5!(A0)$) -- ($(A88)! 0.5!(A0)$) -- (A0) -- cycle;
\fill[color1!50] (A1) -- ($(A1)! 0.5!(A10)$) -- ($(A2)! 0.5!(A10)$) -- (A2) -- cycle;
\fill[color1!50] (A1) -- (A2) -- (A89) -- cycle;
\fill[color1!50] (A1) -- ($(A1)! 0.5!(A10)$) -- ($(A9)! 0.5!(A10)$) -- (A9) -- cycle;
\fill[color1!50] (A1) -- ($(A1)! 0.5!(A88)$) -- ($(A89)! 0.5!(A88)$) -- (A89) -- cycle;
\fill[color1!50] (A2) -- (A3) -- (A11) -- cycle;
\fill[color1!50] (A2) -- ($(A2)! 0.5!(A90)$) -- ($(A3)! 0.5!(A90)$) -- (A3) -- cycle;
\fill[color1!50] (A2) -- ($(A2)! 0.5!(A10)$) -- ($(A11)! 0.5!(A10)$) -- (A11) -- cycle;
\fill[color1!50] (A2) -- ($(A2)! 0.5!(A90)$) -- ($(A89)! 0.5!(A90)$) -- (A89) -- cycle;
\fill[color1!50] (A3) -- ($(A3)! 0.5!(A12)$) -- ($(A4)! 0.5!(A12)$) -- (A4) -- cycle;
\fill[color1!50] (A3) -- (A4) -- (A91) -- cycle;
\fill[color1!50] (A3) -- ($(A3)! 0.5!(A12)$) -- ($(A11)! 0.5!(A12)$) -- (A11) -- cycle;
\fill[color1!50] (A3) -- ($(A3)! 0.5!(A90)$) -- ($(A91)! 0.5!(A90)$) -- (A91) -- cycle;
\fill[color1!50] (A4) -- (A5) -- (A13) -- cycle;
\fill[color1!50] (A4) -- ($(A4)! 0.5!(A92)$) -- ($(A5)! 0.5!(A92)$) -- (A5) -- cycle;
\fill[color1!50] (A4) -- ($(A4)! 0.5!(A12)$) -- ($(A13)! 0.5!(A12)$) -- (A13) -- cycle;
\fill[color1!50] (A4) -- ($(A4)! 0.5!(A92)$) -- ($(A91)! 0.5!(A92)$) -- (A91) -- cycle;
\fill[color1!50] (A5) -- ($(A5)! 0.5!(A6)$) -- ($(A14)! 0.5!(A6)$) -- (A14) -- cycle;
\fill[color1!50] ($(A6)! 0.5!(A5)$) -- ($(A93)! 0.5!(A5)$) -- (A5) -- cycle;
\fill[color1!50] (A5) -- (A13) -- (A32) -- cycle;
\fill[color1!50] (A5) -- ($(A5)! 0.5!(A21)$) -- ($(A14)! 0.5!(A21)$) -- (A14) -- cycle;
\fill[color1!50] (A5) -- ($(A5)! 0.5!(A21)$) -- ($(A27)! 0.5!(A21)$) -- (A27) -- cycle;
\fill[color1!50] (A5) -- (A27) -- (A42) -- cycle;
\fill[color1!50] (A5) -- (A32) -- (A42) -- cycle;
\fill[color1!50] (A5) -- ($(A5)! 0.5!(A92)$) -- ($(A126)! 0.5!(A92)$) -- (A126) -- cycle;
\fill[color1!50] ($(A93)! 0.5!(A5)$) -- ($(A106)! 0.5!(A5)$) -- (A5) -- cycle;
\fill[color1!50] ($(A106)! 0.5!(A5)$) -- ($(A117)! 0.5!(A5)$) -- (A5) -- cycle;
\fill[color1!50] ($(A117)! 0.5!(A5)$) -- ($(A142)! 0.5!(A5)$) -- (A5) -- cycle;
\fill[color1!50] (A5) -- ($(A5)! 0.5!(A142)$) -- ($(A126)! 0.5!(A142)$) -- (A126) -- cycle;
\fill[color1!50] (A7) -- ($(A7)! 0.5!(A6)$) -- ($(A15)! 0.5!(A6)$) -- (A15) -- cycle;
\fill[color1!50] ($(A6)! 0.5!(A7)$) -- ($(A94)! 0.5!(A7)$) -- (A7) -- cycle;
\fill[color1!50] (A14) -- ($(A14)! 0.5!(A6)$) -- ($(A15)! 0.5!(A6)$) -- (A15) -- cycle;
\fill[color1!50] (A7) -- ($(A7)! 0.5!(A8)$) -- ($(A16)! 0.5!(A8)$) -- (A16) -- cycle;
\fill[color1!50] ($(A8)! 0.5!(A7)$) -- ($(A95)! 0.5!(A7)$) -- (A7) -- cycle;
\fill[color1!50] (A7) -- (A15) -- (A16) -- cycle;
\fill[color1!50] ($(A94)! 0.5!(A7)$) -- ($(A95)! 0.5!(A7)$) -- (A7) -- cycle;
\fill[color1!50] (A9) -- ($(A9)! 0.5!(A10)$) -- ($(A17)! 0.5!(A10)$) -- (A17) -- cycle;
\fill[color1!50] (A9) -- (A17) -- (A53) -- cycle;
\fill[color1!50] (A9) -- ($(A9)! 0.5!(A45)$) -- ($(A53)! 0.5!(A45)$) -- (A53) -- cycle;
\fill[color1!50] (A11) -- ($(A11)! 0.5!(A10)$) -- ($(A18)! 0.5!(A10)$) -- (A18) -- cycle;
\fill[color1!50] (A17) -- ($(A17)! 0.5!(A10)$) -- ($(A18)! 0.5!(A10)$) -- (A18) -- cycle;
\fill[color1!50] (A11) -- ($(A11)! 0.5!(A12)$) -- ($(A19)! 0.5!(A12)$) -- (A19) -- cycle;
\fill[color1!50] (A11) -- (A18) -- (A19) -- cycle;
\fill[color1!50] (A13) -- ($(A13)! 0.5!(A12)$) -- ($(A20)! 0.5!(A12)$) -- (A20) -- cycle;
\fill[color1!50] (A19) -- ($(A19)! 0.5!(A12)$) -- ($(A20)! 0.5!(A12)$) -- (A20) -- cycle;
\fill[color1!50] (A13) -- (A20) -- (A32) -- cycle;
\fill[color1!50] (A14) -- (A15) -- (A22) -- cycle;
\fill[color1!50] (A14) -- ($(A14)! 0.5!(A21)$) -- ($(A22)! 0.5!(A21)$) -- (A22) -- cycle;
\fill[color1!50] (A15) -- ($(A15)! 0.5!(A23)$) -- ($(A16)! 0.5!(A23)$) -- (A16) -- cycle;
\fill[color1!50] (A15) -- ($(A15)! 0.5!(A23)$) -- ($(A22)! 0.5!(A23)$) -- (A22) -- cycle;
\fill[color1!50] (A17) -- (A18) -- (A24) -- cycle;
\fill[color1!50] (A17) -- ($(A17)! 0.5!(A60)$) -- ($(A24)! 0.5!(A60)$) -- (A24) -- cycle;
\fill[color1!50] (A17) -- ($(A17)! 0.5!(A60)$) -- ($(A53)! 0.5!(A60)$) -- (A53) -- cycle;
\fill[color1!50] (A18) -- ($(A18)! 0.5!(A25)$) -- ($(A19)! 0.5!(A25)$) -- (A19) -- cycle;
\fill[color1!50] (A18) -- ($(A18)! 0.5!(A25)$) -- ($(A24)! 0.5!(A25)$) -- (A24) -- cycle;
\fill[color1!50] (A19) -- (A20) -- (A26) -- cycle;
\fill[color1!50] (A19) -- ($(A19)! 0.5!(A25)$) -- ($(A26)! 0.5!(A25)$) -- (A26) -- cycle;
\fill[color1!50] (A20) -- (A26) -- (A32) -- cycle;
\fill[color1!50] (A22) -- ($(A22)! 0.5!(A21)$) -- ($(A28)! 0.5!(A21)$) -- (A28) -- cycle;
\fill[color1!50] (A27) -- ($(A27)! 0.5!(A21)$) -- ($(A28)! 0.5!(A21)$) -- (A28) -- cycle;
\fill[color1!50] (A22) -- ($(A22)! 0.5!(A23)$) -- ($(A29)! 0.5!(A23)$) -- (A29) -- cycle;
\fill[color1!50] (A22) -- (A28) -- (A29) -- cycle;
\fill[color1!50] (A24) -- ($(A24)! 0.5!(A25)$) -- ($(A30)! 0.5!(A25)$) -- (A30) -- cycle;
\fill[color1!50] (A24) -- (A30) -- (A66) -- cycle;
\fill[color1!50] (A24) -- ($(A24)! 0.5!(A60)$) -- ($(A66)! 0.5!(A60)$) -- (A66) -- cycle;
\fill[color1!50] (A26) -- ($(A26)! 0.5!(A25)$) -- ($(A31)! 0.5!(A25)$) -- (A31) -- cycle;
\fill[color1!50] (A30) -- ($(A30)! 0.5!(A25)$) -- ($(A31)! 0.5!(A25)$) -- (A31) -- cycle;
\fill[color1!50] (A26) -- (A31) -- (A32) -- cycle;
\fill[color1!50] (A27) -- (A28) -- (A33) -- cycle;
\fill[color1!50] (A27) -- (A33) -- (A37) -- cycle;
\fill[color1!50] (A27) -- (A37) -- (A40) -- cycle;
\fill[color1!50] (A27) -- (A40) -- (A42) -- cycle;
\fill[color1!50] (A28) -- ($(A28)! 0.5!(A34)$) -- ($(A29)! 0.5!(A34)$) -- (A29) -- cycle;
\fill[color1!50] (A28) -- ($(A28)! 0.5!(A34)$) -- ($(A33)! 0.5!(A34)$) -- (A33) -- cycle;
\fill[color1!50] (A30) -- (A31) -- (A35) -- cycle;
\fill[color1!50] (A30) -- ($(A30)! 0.5!(A71)$) -- ($(A35)! 0.5!(A71)$) -- (A35) -- cycle;
\fill[color1!50] (A30) -- ($(A30)! 0.5!(A71)$) -- ($(A66)! 0.5!(A71)$) -- (A66) -- cycle;
\fill[color1!50] (A31) -- ($(A31)! 0.5!(A36)$) -- ($(A32)! 0.5!(A36)$) -- (A32) -- cycle;
\fill[color1!50] (A31) -- ($(A31)! 0.5!(A36)$) -- ($(A35)! 0.5!(A36)$) -- (A35) -- cycle;
\fill[color1!50] (A32) -- ($(A32)! 0.5!(A36)$) -- ($(A42)! 0.5!(A36)$) -- (A42) -- cycle;
\fill[color1!50] (A33) -- ($(A33)! 0.5!(A34)$) -- ($(A38)! 0.5!(A34)$) -- (A38) -- cycle;
\fill[color1!50] (A33) -- (A37) -- (A38) -- cycle;
\fill[color1!50] (A35) -- ($(A35)! 0.5!(A36)$) -- ($(A39)! 0.5!(A36)$) -- (A39) -- cycle;
\fill[color1!50] (A35) -- (A39) -- (A75) -- cycle;
\fill[color1!50] (A35) -- ($(A35)! 0.5!(A71)$) -- ($(A75)! 0.5!(A71)$) -- (A75) -- cycle;
\fill[color1!50] (A39) -- ($(A39)! 0.5!(A36)$) -- ($(A42)! 0.5!(A36)$) -- (A42) -- cycle;
\fill[color1!50] (A37) -- ($(A37)! 0.5!(A41)$) -- ($(A38)! 0.5!(A41)$) -- (A38) -- cycle;
\fill[color1!50] (A37) -- ($(A37)! 0.5!(A41)$) -- ($(A40)! 0.5!(A41)$) -- (A40) -- cycle;
\fill[color1!50] (A39) -- (A42) -- (A75) -- cycle;
\fill[color1!50] (A40) -- ($(A40)! 0.5!(A41)$) -- ($(A43)! 0.5!(A41)$) -- (A43) -- cycle;
\fill[color1!50] (A40) -- (A42) -- (A43) -- cycle;
\fill[color1!50] (A42) -- ($(A42)! 0.5!(A44)$) -- ($(A43)! 0.5!(A44)$) -- (A43) -- cycle;
\fill[color1!50] ($(A44)! 0.5!(A42)$) -- ($(A80)! 0.5!(A42)$) -- (A42) -- cycle;
\fill[color1!50] ($(A49)! 0.5!(A42)$) -- ($(A68)! 0.5!(A42)$) -- (A42) -- cycle;
\fill[color1!50] (A42) -- ($(A42)! 0.5!(A49)$) -- ($(A72)! 0.5!(A49)$) -- (A72) -- cycle;
\fill[color1!50] ($(A68)! 0.5!(A42)$) -- ($(A78)! 0.5!(A42)$) -- (A42) -- cycle;
\fill[color1!50] (A42) -- (A72) -- (A75) -- cycle;
\fill[color1!50] ($(A78)! 0.5!(A42)$) -- ($(A80)! 0.5!(A42)$) -- (A42) -- cycle;
\fill[color1!50] (A46) -- ($(A46)! 0.5!(A45)$) -- ($(A53)! 0.5!(A45)$) -- (A53) -- cycle;
\fill[color1!50] ($(A45)! 0.5!(A46)$) -- ($(A87)! 0.5!(A46)$) -- (A46) -- cycle;
\fill[color1!50] (A46) -- ($(A46)! 0.5!(A47)$) -- ($(A54)! 0.5!(A47)$) -- (A54) -- cycle;
\fill[color1!50] ($(A47)! 0.5!(A46)$) -- ($(A86)! 0.5!(A46)$) -- (A46) -- cycle;
\fill[color1!50] (A46) -- (A53) -- (A54) -- cycle;
\fill[color1!50] ($(A86)! 0.5!(A46)$) -- ($(A87)! 0.5!(A46)$) -- (A46) -- cycle;
\fill[color1!50] (A48) -- ($(A48)! 0.5!(A47)$) -- ($(A55)! 0.5!(A47)$) -- (A55) -- cycle;
\fill[color1!50] ($(A47)! 0.5!(A48)$) -- ($(A85)! 0.5!(A48)$) -- (A48) -- cycle;
\fill[color1!50] (A54) -- ($(A54)! 0.5!(A47)$) -- ($(A55)! 0.5!(A47)$) -- (A55) -- cycle;
\fill[color1!50] (A48) -- ($(A48)! 0.5!(A49)$) -- ($(A56)! 0.5!(A49)$) -- (A56) -- cycle;
\fill[color1!50] ($(A49)! 0.5!(A48)$) -- ($(A84)! 0.5!(A48)$) -- (A48) -- cycle;
\fill[color1!50] (A48) -- (A55) -- (A56) -- cycle;
\fill[color1!50] ($(A84)! 0.5!(A48)$) -- ($(A85)! 0.5!(A48)$) -- (A48) -- cycle;
\fill[color1!50] (A56) -- ($(A56)! 0.5!(A49)$) -- ($(A72)! 0.5!(A49)$) -- (A72) -- cycle;
\fill[color1!50] ($(A51)! 0.5!(A81)$) -- ($(A52)! 0.5!(A81)$) -- (A81) -- cycle;
\fill[color1!50] ($(A51)! 0.5!(A81)$) -- ($(A82)! 0.5!(A81)$) -- (A81) -- cycle;
\fill[color1!50] (A53) -- ($(A53)! 0.5!(A60)$) -- ($(A54)! 0.5!(A60)$) -- (A54) -- cycle;
\fill[color1!50] (A54) -- (A55) -- (A61) -- cycle;
\fill[color1!50] (A54) -- ($(A54)! 0.5!(A60)$) -- ($(A61)! 0.5!(A60)$) -- (A61) -- cycle;
\fill[color1!50] (A55) -- ($(A55)! 0.5!(A62)$) -- ($(A56)! 0.5!(A62)$) -- (A56) -- cycle;
\fill[color1!50] (A55) -- ($(A55)! 0.5!(A62)$) -- ($(A61)! 0.5!(A62)$) -- (A61) -- cycle;
\fill[color1!50] (A56) -- ($(A56)! 0.5!(A62)$) -- ($(A72)! 0.5!(A62)$) -- (A72) -- cycle;
\fill[color1!50] (A61) -- ($(A61)! 0.5!(A60)$) -- ($(A66)! 0.5!(A60)$) -- (A66) -- cycle;
\fill[color1!50] (A61) -- ($(A61)! 0.5!(A62)$) -- ($(A67)! 0.5!(A62)$) -- (A67) -- cycle;
\fill[color1!50] (A61) -- (A66) -- (A67) -- cycle;
\fill[color1!50] (A67) -- ($(A67)! 0.5!(A62)$) -- ($(A72)! 0.5!(A62)$) -- (A72) -- cycle;
\fill[color1!50] (A66) -- ($(A66)! 0.5!(A71)$) -- ($(A67)! 0.5!(A71)$) -- (A67) -- cycle;
\fill[color1!50] (A67) -- ($(A67)! 0.5!(A71)$) -- ($(A72)! 0.5!(A71)$) -- (A72) -- cycle;
\fill[color1!50] (A72) -- ($(A72)! 0.5!(A71)$) -- ($(A75)! 0.5!(A71)$) -- (A75) -- cycle;
\fill[color1!50] ($(A82)! 0.5!(A81)$) -- ($(A96)! 0.5!(A81)$) -- (A81) -- cycle;
\fill[color1!50] ($(A87)! 0.5!(A102)$) -- ($(A88)! 0.5!(A102)$) -- (A102) -- cycle;
\fill[color1!50] ($(A87)! 0.5!(A102)$) -- ($(A101)! 0.5!(A102)$) -- (A102) -- cycle;
\fill[color1!50] (A89) -- ($(A89)! 0.5!(A88)$) -- ($(A102)! 0.5!(A88)$) -- (A102) -- cycle;
\fill[color1!50] (A89) -- ($(A89)! 0.5!(A90)$) -- ($(A103)! 0.5!(A90)$) -- (A103) -- cycle;
\fill[color1!50] (A89) -- (A102) -- (A103) -- cycle;
\fill[color1!50] (A91) -- ($(A91)! 0.5!(A90)$) -- ($(A104)! 0.5!(A90)$) -- (A104) -- cycle;
\fill[color1!50] (A103) -- ($(A103)! 0.5!(A90)$) -- ($(A104)! 0.5!(A90)$) -- (A104) -- cycle;
\fill[color1!50] (A91) -- ($(A91)! 0.5!(A92)$) -- ($(A105)! 0.5!(A92)$) -- (A105) -- cycle;
\fill[color1!50] (A91) -- (A104) -- (A105) -- cycle;
\fill[color1!50] (A105) -- ($(A105)! 0.5!(A92)$) -- ($(A126)! 0.5!(A92)$) -- (A126) -- cycle;
\fill[color1!50] ($(A96)! 0.5!(A109)$) -- ($(A97)! 0.5!(A109)$) -- (A109) -- cycle;
\fill[color1!50] ($(A97)! 0.5!(A109)$) -- ($(A110)! 0.5!(A109)$) -- (A109) -- cycle;
\fill[color1!50] ($(A101)! 0.5!(A102)$) -- ($(A114)! 0.5!(A102)$) -- (A102) -- cycle;
\fill[color1!50] (A102) -- ($(A102)! 0.5!(A114)$) -- ($(A103)! 0.5!(A114)$) -- (A103) -- cycle;
\fill[color1!50] (A103) -- (A104) -- (A115) -- cycle;
\fill[color1!50] (A103) -- ($(A103)! 0.5!(A114)$) -- ($(A115)! 0.5!(A114)$) -- (A115) -- cycle;
\fill[color1!50] (A104) -- ($(A104)! 0.5!(A116)$) -- ($(A105)! 0.5!(A116)$) -- (A105) -- cycle;
\fill[color1!50] (A104) -- ($(A104)! 0.5!(A116)$) -- ($(A115)! 0.5!(A116)$) -- (A115) -- cycle;
\fill[color1!50] (A105) -- ($(A105)! 0.5!(A116)$) -- ($(A126)! 0.5!(A116)$) -- (A126) -- cycle;
\fill[color1!50] ($(A110)! 0.5!(A109)$) -- ($(A120)! 0.5!(A109)$) -- (A109) -- cycle;
\fill[color1!50] ($(A113)! 0.5!(A124)$) -- ($(A114)! 0.5!(A124)$) -- (A124) -- cycle;
\fill[color1!50] ($(A113)! 0.5!(A124)$) -- ($(A123)! 0.5!(A124)$) -- (A124) -- cycle;
\fill[color1!50] (A115) -- ($(A115)! 0.5!(A114)$) -- ($(A124)! 0.5!(A114)$) -- (A124) -- cycle;
\fill[color1!50] (A115) -- ($(A115)! 0.5!(A116)$) -- ($(A125)! 0.5!(A116)$) -- (A125) -- cycle;
\fill[color1!50] (A115) -- (A124) -- (A125) -- cycle;
\fill[color1!50] (A125) -- ($(A125)! 0.5!(A116)$) -- ($(A126)! 0.5!(A116)$) -- (A126) -- cycle;
\fill[color1!50] ($(A120)! 0.5!(A129)$) -- ($(A121)! 0.5!(A129)$) -- (A129) -- cycle;
\fill[color1!50] ($(A121)! 0.5!(A129)$) -- ($(A130)! 0.5!(A129)$) -- (A129) -- cycle;
\fill[color1!50] ($(A123)! 0.5!(A124)$) -- ($(A132)! 0.5!(A124)$) -- (A124) -- cycle;
\fill[color1!50] (A124) -- ($(A124)! 0.5!(A132)$) -- ($(A125)! 0.5!(A132)$) -- (A125) -- cycle;
\fill[color1!50] (A125) -- (A126) -- (A133) -- cycle;
\fill[color1!50] (A125) -- ($(A125)! 0.5!(A132)$) -- ($(A133)! 0.5!(A132)$) -- (A133) -- cycle;
\fill[color1!50] (A126) -- ($(A126)! 0.5!(A142)$) -- ($(A133)! 0.5!(A142)$) -- (A133) -- cycle;
\fill[color1!50] ($(A130)! 0.5!(A129)$) -- ($(A136)! 0.5!(A129)$) -- (A129) -- cycle;
\fill[color1!50] ($(A131)! 0.5!(A138)$) -- ($(A132)! 0.5!(A138)$) -- (A138) -- cycle;
\fill[color1!50] ($(A131)! 0.5!(A138)$) -- ($(A142)! 0.5!(A138)$) -- (A138) -- cycle;
\fill[color1!50] (A133) -- ($(A133)! 0.5!(A132)$) -- ($(A138)! 0.5!(A132)$) -- (A138) -- cycle;
\fill[color1!50] (A133) -- ($(A133)! 0.5!(A142)$) -- ($(A138)! 0.5!(A142)$) -- (A138) -- cycle;
\fill[color1!50] ($(A136)! 0.5!(A141)$) -- ($(A137)! 0.5!(A141)$) -- (A141) -- cycle;
\fill[color1!50] ($(A137)! 0.5!(A141)$) -- ($(A142)! 0.5!(A141)$) -- (A141) -- cycle;
\fill[color1!50] ($(A142)! 0.5!(A141)$) -- ($(A144)! 0.5!(A141)$) -- (A141) -- cycle;
\colorlet{color2}{mycolor5}
\fill[color2!50] ($(A1)! 0.5!(A10)$) -- ($(A2)! 0.5!(A10)$) -- (A10) -- cycle;
\fill[color2!50] ($(A1)! 0.5!(A10)$) -- ($(A9)! 0.5!(A10)$) -- (A10) -- cycle;
\fill[color2!50] ($(A2)! 0.5!(A10)$) -- ($(A11)! 0.5!(A10)$) -- (A10) -- cycle;
\fill[color2!50] ($(A9)! 0.5!(A10)$) -- ($(A17)! 0.5!(A10)$) -- (A10) -- cycle;
\fill[color2!50] ($(A11)! 0.5!(A10)$) -- ($(A18)! 0.5!(A10)$) -- (A10) -- cycle;
\fill[color2!50] ($(A17)! 0.5!(A10)$) -- ($(A18)! 0.5!(A10)$) -- (A10) -- cycle;
\colorlet{color3}{mycolor5}
\fill[color3!50] ($(A3)! 0.5!(A12)$) -- ($(A4)! 0.5!(A12)$) -- (A12) -- cycle;
\fill[color3!50] ($(A3)! 0.5!(A12)$) -- ($(A11)! 0.5!(A12)$) -- (A12) -- cycle;
\fill[color3!50] ($(A4)! 0.5!(A12)$) -- ($(A13)! 0.5!(A12)$) -- (A12) -- cycle;
\fill[color3!50] ($(A11)! 0.5!(A12)$) -- ($(A19)! 0.5!(A12)$) -- (A12) -- cycle;
\fill[color3!50] ($(A13)! 0.5!(A12)$) -- ($(A20)! 0.5!(A12)$) -- (A12) -- cycle;
\fill[color3!50] ($(A19)! 0.5!(A12)$) -- ($(A20)! 0.5!(A12)$) -- (A12) -- cycle;
\colorlet{color4}{mycolor3}
\fill[color4!50] ($(A5)! 0.5!(A21)$) -- ($(A14)! 0.5!(A21)$) -- (A21) -- cycle;
\fill[color4!50] ($(A5)! 0.5!(A21)$) -- ($(A27)! 0.5!(A21)$) -- (A21) -- cycle;
\fill[color4!50] ($(A14)! 0.5!(A21)$) -- ($(A22)! 0.5!(A21)$) -- (A21) -- cycle;
\fill[color4!50] ($(A22)! 0.5!(A21)$) -- ($(A28)! 0.5!(A21)$) -- (A21) -- cycle;
\fill[color4!50] ($(A27)! 0.5!(A21)$) -- ($(A28)! 0.5!(A21)$) -- (A21) -- cycle;
\colorlet{color5}{mycolor5}
\fill[color5!50] ($(A18)! 0.5!(A25)$) -- ($(A19)! 0.5!(A25)$) -- (A25) -- cycle;
\fill[color5!50] ($(A18)! 0.5!(A25)$) -- ($(A24)! 0.5!(A25)$) -- (A25) -- cycle;
\fill[color5!50] ($(A19)! 0.5!(A25)$) -- ($(A26)! 0.5!(A25)$) -- (A25) -- cycle;
\fill[color5!50] ($(A24)! 0.5!(A25)$) -- ($(A30)! 0.5!(A25)$) -- (A25) -- cycle;
\fill[color5!50] ($(A26)! 0.5!(A25)$) -- ($(A31)! 0.5!(A25)$) -- (A25) -- cycle;
\fill[color5!50] ($(A30)! 0.5!(A25)$) -- ($(A31)! 0.5!(A25)$) -- (A25) -- cycle;
\colorlet{color6}{mycolor5}
\fill[color6!50] ($(A31)! 0.5!(A36)$) -- ($(A32)! 0.5!(A36)$) -- (A36) -- cycle;
\fill[color6!50] ($(A31)! 0.5!(A36)$) -- ($(A35)! 0.5!(A36)$) -- (A36) -- cycle;
\fill[color6!50] ($(A32)! 0.5!(A36)$) -- ($(A42)! 0.5!(A36)$) -- (A36) -- cycle;
\fill[color6!50] ($(A35)! 0.5!(A36)$) -- ($(A39)! 0.5!(A36)$) -- (A36) -- cycle;
\fill[color6!50] ($(A39)! 0.5!(A36)$) -- ($(A42)! 0.5!(A36)$) -- (A36) -- cycle;
\colorlet{color7}{mycolor3}
\fill[color7!50] ($(A50)! 0.5!(A58)$) -- ($(A51)! 0.5!(A58)$) -- (A58) -- cycle;
\fill[color7!50] ($(A50)! 0.5!(A58)$) -- ($(A57)! 0.5!(A58)$) -- (A58) -- cycle;
\fill[color7!50] ($(A51)! 0.5!(A58)$) -- ($(A59)! 0.5!(A58)$) -- (A58) -- cycle;
\fill[color7!50] ($(A57)! 0.5!(A58)$) -- ($(A64)! 0.5!(A58)$) -- (A58) -- cycle;
\fill[color7!50] ($(A59)! 0.5!(A58)$) -- ($(A65)! 0.5!(A58)$) -- (A58) -- cycle;
\fill[color7!50] ($(A64)! 0.5!(A58)$) -- ($(A65)! 0.5!(A58)$) -- (A58) -- cycle;
\colorlet{color8}{mycolor5}
\fill[color8!50] ($(A17)! 0.5!(A60)$) -- ($(A24)! 0.5!(A60)$) -- (A60) -- cycle;
\fill[color8!50] ($(A17)! 0.5!(A60)$) -- ($(A53)! 0.5!(A60)$) -- (A60) -- cycle;
\fill[color8!50] ($(A24)! 0.5!(A60)$) -- ($(A66)! 0.5!(A60)$) -- (A60) -- cycle;
\fill[color8!50] ($(A53)! 0.5!(A60)$) -- ($(A54)! 0.5!(A60)$) -- (A60) -- cycle;
\fill[color8!50] ($(A54)! 0.5!(A60)$) -- ($(A61)! 0.5!(A60)$) -- (A60) -- cycle;
\fill[color8!50] ($(A61)! 0.5!(A60)$) -- ($(A66)! 0.5!(A60)$) -- (A60) -- cycle;
\colorlet{color9}{mycolor5}
\fill[color9!50] ($(A55)! 0.5!(A62)$) -- ($(A56)! 0.5!(A62)$) -- (A62) -- cycle;
\fill[color9!50] ($(A55)! 0.5!(A62)$) -- ($(A61)! 0.5!(A62)$) -- (A62) -- cycle;
\fill[color9!50] ($(A56)! 0.5!(A62)$) -- ($(A72)! 0.5!(A62)$) -- (A62) -- cycle;
\fill[color9!50] ($(A61)! 0.5!(A62)$) -- ($(A67)! 0.5!(A62)$) -- (A62) -- cycle;
\fill[color9!50] ($(A67)! 0.5!(A62)$) -- ($(A72)! 0.5!(A62)$) -- (A62) -- cycle;
\colorlet{color10}{mycolor3}
\fill[color10!50] ($(A63)! 0.5!(A69)$) -- ($(A64)! 0.5!(A69)$) -- (A69) -- cycle;
\fill[color10!50] ($(A63)! 0.5!(A69)$) -- ($(A68)! 0.5!(A69)$) -- (A69) -- cycle;
\fill[color10!50] ($(A64)! 0.5!(A69)$) -- ($(A70)! 0.5!(A69)$) -- (A69) -- cycle;
\fill[color10!50] ($(A68)! 0.5!(A69)$) -- ($(A73)! 0.5!(A69)$) -- (A69) -- cycle;
\fill[color10!50] ($(A70)! 0.5!(A69)$) -- ($(A74)! 0.5!(A69)$) -- (A69) -- cycle;
\fill[color10!50] ($(A73)! 0.5!(A69)$) -- ($(A74)! 0.5!(A69)$) -- (A69) -- cycle;
\colorlet{color11}{mycolor5}
\fill[color11!50] ($(A30)! 0.5!(A71)$) -- ($(A35)! 0.5!(A71)$) -- (A71) -- cycle;
\fill[color11!50] ($(A30)! 0.5!(A71)$) -- ($(A66)! 0.5!(A71)$) -- (A71) -- cycle;
\fill[color11!50] ($(A35)! 0.5!(A71)$) -- ($(A75)! 0.5!(A71)$) -- (A71) -- cycle;
\fill[color11!50] ($(A66)! 0.5!(A71)$) -- ($(A67)! 0.5!(A71)$) -- (A71) -- cycle;
\fill[color11!50] ($(A67)! 0.5!(A71)$) -- ($(A72)! 0.5!(A71)$) -- (A71) -- cycle;
\fill[color11!50] ($(A72)! 0.5!(A71)$) -- ($(A75)! 0.5!(A71)$) -- (A71) -- cycle;
\colorlet{color12}{mycolor3}
\fill[color12!50] ($(A68)! 0.5!(A76)$) -- ($(A73)! 0.5!(A76)$) -- (A76) -- cycle;
\fill[color12!50] ($(A68)! 0.5!(A76)$) -- ($(A78)! 0.5!(A76)$) -- (A76) -- cycle;
\fill[color12!50] ($(A73)! 0.5!(A76)$) -- ($(A77)! 0.5!(A76)$) -- (A76) -- cycle;
\fill[color12!50] ($(A77)! 0.5!(A76)$) -- ($(A79)! 0.5!(A76)$) -- (A76) -- cycle;
\fill[color12!50] ($(A78)! 0.5!(A76)$) -- ($(A79)! 0.5!(A76)$) -- (A76) -- cycle;
\colorlet{color13}{mycolor3}
\fill[color13!50] ($(A49)! 0.5!(A83)$) -- ($(A50)! 0.5!(A83)$) -- (A83) -- cycle;
\fill[color13!50] ($(A49)! 0.5!(A83)$) -- ($(A98)! 0.5!(A83)$) -- (A83) -- cycle;
\fill[color13!50] ($(A50)! 0.5!(A83)$) -- ($(A82)! 0.5!(A83)$) -- (A83) -- cycle;
\fill[color13!50] ($(A82)! 0.5!(A83)$) -- ($(A97)! 0.5!(A83)$) -- (A83) -- cycle;
\fill[color13!50] ($(A97)! 0.5!(A83)$) -- ($(A98)! 0.5!(A83)$) -- (A83) -- cycle;
\colorlet{color14}{mycolor5}
\fill[color14!50] ($(A2)! 0.5!(A90)$) -- ($(A3)! 0.5!(A90)$) -- (A90) -- cycle;
\fill[color14!50] ($(A2)! 0.5!(A90)$) -- ($(A89)! 0.5!(A90)$) -- (A90) -- cycle;
\fill[color14!50] ($(A3)! 0.5!(A90)$) -- ($(A91)! 0.5!(A90)$) -- (A90) -- cycle;
\fill[color14!50] ($(A89)! 0.5!(A90)$) -- ($(A103)! 0.5!(A90)$) -- (A90) -- cycle;
\fill[color14!50] ($(A91)! 0.5!(A90)$) -- ($(A104)! 0.5!(A90)$) -- (A90) -- cycle;
\fill[color14!50] ($(A103)! 0.5!(A90)$) -- ($(A104)! 0.5!(A90)$) -- (A90) -- cycle;
\colorlet{color15}{mycolor5}
\fill[color15!50] ($(A4)! 0.5!(A92)$) -- ($(A5)! 0.5!(A92)$) -- (A92) -- cycle;
\fill[color15!50] ($(A4)! 0.5!(A92)$) -- ($(A91)! 0.5!(A92)$) -- (A92) -- cycle;
\fill[color15!50] ($(A5)! 0.5!(A92)$) -- ($(A126)! 0.5!(A92)$) -- (A92) -- cycle;
\fill[color15!50] ($(A91)! 0.5!(A92)$) -- ($(A105)! 0.5!(A92)$) -- (A92) -- cycle;
\fill[color15!50] ($(A105)! 0.5!(A92)$) -- ($(A126)! 0.5!(A92)$) -- (A92) -- cycle;
\colorlet{color16}{mycolor5}
\fill[color16!50] ($(A85)! 0.5!(A100)$) -- ($(A86)! 0.5!(A100)$) -- (A100) -- cycle;
\fill[color16!50] ($(A85)! 0.5!(A100)$) -- ($(A99)! 0.5!(A100)$) -- (A100) -- cycle;
\fill[color16!50] ($(A86)! 0.5!(A100)$) -- ($(A101)! 0.5!(A100)$) -- (A100) -- cycle;
\fill[color16!50] ($(A99)! 0.5!(A100)$) -- ($(A112)! 0.5!(A100)$) -- (A100) -- cycle;
\fill[color16!50] ($(A101)! 0.5!(A100)$) -- ($(A113)! 0.5!(A100)$) -- (A100) -- cycle;
\fill[color16!50] ($(A112)! 0.5!(A100)$) -- ($(A113)! 0.5!(A100)$) -- (A100) -- cycle;
\colorlet{color17}{mycolor3}
\fill[color17!50] ($(A93)! 0.5!(A107)$) -- ($(A94)! 0.5!(A107)$) -- (A107) -- cycle;
\fill[color17!50] ($(A93)! 0.5!(A107)$) -- ($(A106)! 0.5!(A107)$) -- (A107) -- cycle;
\fill[color17!50] ($(A94)! 0.5!(A107)$) -- ($(A108)! 0.5!(A107)$) -- (A107) -- cycle;
\fill[color17!50] ($(A106)! 0.5!(A107)$) -- ($(A118)! 0.5!(A107)$) -- (A107) -- cycle;
\fill[color17!50] ($(A108)! 0.5!(A107)$) -- ($(A119)! 0.5!(A107)$) -- (A107) -- cycle;
\fill[color17!50] ($(A118)! 0.5!(A107)$) -- ($(A119)! 0.5!(A107)$) -- (A107) -- cycle;
\colorlet{color18}{mycolor3}
\fill[color18!50] ($(A49)! 0.5!(A111)$) -- ($(A98)! 0.5!(A111)$) -- (A111) -- cycle;
\fill[color18!50] ($(A49)! 0.5!(A111)$) -- ($(A142)! 0.5!(A111)$) -- (A111) -- cycle;
\fill[color18!50] ($(A98)! 0.5!(A111)$) -- ($(A110)! 0.5!(A111)$) -- (A111) -- cycle;
\fill[color18!50] ($(A110)! 0.5!(A111)$) -- ($(A121)! 0.5!(A111)$) -- (A111) -- cycle;
\fill[color18!50] ($(A121)! 0.5!(A111)$) -- ($(A130)! 0.5!(A111)$) -- (A111) -- cycle;
\fill[color18!50] ($(A130)! 0.5!(A111)$) -- ($(A137)! 0.5!(A111)$) -- (A111) -- cycle;
\fill[color18!50] ($(A137)! 0.5!(A111)$) -- ($(A142)! 0.5!(A111)$) -- (A111) -- cycle;
\colorlet{color19}{mycolor5}
\fill[color19!50] ($(A104)! 0.5!(A116)$) -- ($(A105)! 0.5!(A116)$) -- (A116) -- cycle;
\fill[color19!50] ($(A104)! 0.5!(A116)$) -- ($(A115)! 0.5!(A116)$) -- (A116) -- cycle;
\fill[color19!50] ($(A105)! 0.5!(A116)$) -- ($(A126)! 0.5!(A116)$) -- (A116) -- cycle;
\fill[color19!50] ($(A115)! 0.5!(A116)$) -- ($(A125)! 0.5!(A116)$) -- (A116) -- cycle;
\fill[color19!50] ($(A125)! 0.5!(A116)$) -- ($(A126)! 0.5!(A116)$) -- (A116) -- cycle;
\colorlet{color20}{mycolor5}
\fill[color20!50] ($(A49)! 0.5!(A122)$) -- ($(A84)! 0.5!(A122)$) -- (A122) -- cycle;
\fill[color20!50] ($(A49)! 0.5!(A122)$) -- ($(A142)! 0.5!(A122)$) -- (A122) -- cycle;
\fill[color20!50] ($(A84)! 0.5!(A122)$) -- ($(A99)! 0.5!(A122)$) -- (A122) -- cycle;
\fill[color20!50] ($(A99)! 0.5!(A122)$) -- ($(A112)! 0.5!(A122)$) -- (A122) -- cycle;
\fill[color20!50] ($(A112)! 0.5!(A122)$) -- ($(A123)! 0.5!(A122)$) -- (A122) -- cycle;
\fill[color20!50] ($(A123)! 0.5!(A122)$) -- ($(A131)! 0.5!(A122)$) -- (A122) -- cycle;
\fill[color20!50] ($(A131)! 0.5!(A122)$) -- ($(A142)! 0.5!(A122)$) -- (A122) -- cycle;
\colorlet{color21}{mycolor3}
\fill[color21!50] ($(A117)! 0.5!(A127)$) -- ($(A118)! 0.5!(A127)$) -- (A127) -- cycle;
\fill[color21!50] ($(A117)! 0.5!(A127)$) -- ($(A134)! 0.5!(A127)$) -- (A127) -- cycle;
\fill[color21!50] ($(A118)! 0.5!(A127)$) -- ($(A128)! 0.5!(A127)$) -- (A127) -- cycle;
\fill[color21!50] ($(A128)! 0.5!(A127)$) -- ($(A135)! 0.5!(A127)$) -- (A127) -- cycle;
\fill[color21!50] ($(A134)! 0.5!(A127)$) -- ($(A135)! 0.5!(A127)$) -- (A127) -- cycle;
\colorlet{color22}{mycolor3}
\fill[color22!50] ($(A117)! 0.5!(A139)$) -- ($(A134)! 0.5!(A139)$) -- (A139) -- cycle;
\fill[color22!50] ($(A117)! 0.5!(A139)$) -- ($(A142)! 0.5!(A139)$) -- (A139) -- cycle;
\fill[color22!50] ($(A134)! 0.5!(A139)$) -- ($(A140)! 0.5!(A139)$) -- (A139) -- cycle;
\fill[color22!50] ($(A140)! 0.5!(A139)$) -- ($(A143)! 0.5!(A139)$) -- (A139) -- cycle;
\fill[color22!50] ($(A142)! 0.5!(A139)$) -- ($(A143)! 0.5!(A139)$) -- (A139) -- cycle;
\draw[color22, thick] ($(A117)! 0.5!(A139)$) -- ($(A134)! 0.5!(A139)$);
\draw[color22, thick] ($(A117)! 0.5!(A139)$) -- ($(A142)! 0.5!(A139)$);
\draw[color22, thick] ($(A134)! 0.5!(A139)$) -- ($(A140)! 0.5!(A139)$);
\draw[color22, thick] ($(A140)! 0.5!(A139)$) -- ($(A143)! 0.5!(A139)$);
\draw[color22, thick] ($(A142)! 0.5!(A139)$) -- ($(A143)! 0.5!(A139)$);
\draw[color21, thick] ($(A117)! 0.5!(A127)$) -- ($(A118)! 0.5!(A127)$);
\draw[color21, thick] ($(A117)! 0.5!(A127)$) -- ($(A134)! 0.5!(A127)$);
\draw[color21, thick] ($(A118)! 0.5!(A127)$) -- ($(A128)! 0.5!(A127)$);
\draw[color21, thick] ($(A128)! 0.5!(A127)$) -- ($(A135)! 0.5!(A127)$);
\draw[color21, thick] ($(A134)! 0.5!(A127)$) -- ($(A135)! 0.5!(A127)$);
\draw[color20, thick] ($(A49)! 0.5!(A122)$) -- ($(A84)! 0.5!(A122)$);
\draw[color20, thick] ($(A49)! 0.5!(A122)$) -- ($(A142)! 0.5!(A122)$);
\draw[color20, thick] ($(A84)! 0.5!(A122)$) -- ($(A99)! 0.5!(A122)$);
\draw[color20, thick] ($(A99)! 0.5!(A122)$) -- ($(A112)! 0.5!(A122)$);
\draw[color20, thick] ($(A112)! 0.5!(A122)$) -- ($(A123)! 0.5!(A122)$);
\draw[color20, thick] ($(A123)! 0.5!(A122)$) -- ($(A131)! 0.5!(A122)$);
\draw[color20, thick] ($(A131)! 0.5!(A122)$) -- ($(A142)! 0.5!(A122)$);
\draw[color18, thick] ($(A49)! 0.5!(A111)$) -- ($(A98)! 0.5!(A111)$);
\draw[color18, thick] ($(A49)! 0.5!(A111)$) -- ($(A142)! 0.5!(A111)$);
\draw[color18, thick] ($(A98)! 0.5!(A111)$) -- ($(A110)! 0.5!(A111)$);
\draw[color18, thick] ($(A110)! 0.5!(A111)$) -- ($(A121)! 0.5!(A111)$);
\draw[color18, thick] ($(A121)! 0.5!(A111)$) -- ($(A130)! 0.5!(A111)$);
\draw[color18, thick] ($(A130)! 0.5!(A111)$) -- ($(A137)! 0.5!(A111)$);
\draw[color18, thick] ($(A137)! 0.5!(A111)$) -- ($(A142)! 0.5!(A111)$);
\draw[color17, thick] ($(A93)! 0.5!(A107)$) -- ($(A94)! 0.5!(A107)$);
\draw[color17, thick] ($(A93)! 0.5!(A107)$) -- ($(A106)! 0.5!(A107)$);
\draw[color17, thick] ($(A94)! 0.5!(A107)$) -- ($(A108)! 0.5!(A107)$);
\draw[color17, thick] ($(A106)! 0.5!(A107)$) -- ($(A118)! 0.5!(A107)$);
\draw[color17, thick] ($(A108)! 0.5!(A107)$) -- ($(A119)! 0.5!(A107)$);
\draw[color17, thick] ($(A118)! 0.5!(A107)$) -- ($(A119)! 0.5!(A107)$);
\draw[color16, thick] ($(A85)! 0.5!(A100)$) -- ($(A86)! 0.5!(A100)$);
\draw[color16, thick] ($(A85)! 0.5!(A100)$) -- ($(A99)! 0.5!(A100)$);
\draw[color16, thick] ($(A86)! 0.5!(A100)$) -- ($(A101)! 0.5!(A100)$);
\draw[color16, thick] ($(A99)! 0.5!(A100)$) -- ($(A112)! 0.5!(A100)$);
\draw[color16, thick] ($(A101)! 0.5!(A100)$) -- ($(A113)! 0.5!(A100)$);
\draw[color16, thick] ($(A112)! 0.5!(A100)$) -- ($(A113)! 0.5!(A100)$);
\draw[color13, thick] ($(A49)! 0.5!(A83)$) -- ($(A50)! 0.5!(A83)$);
\draw[color13, thick] ($(A49)! 0.5!(A83)$) -- ($(A98)! 0.5!(A83)$);
\draw[color13, thick] ($(A50)! 0.5!(A83)$) -- ($(A82)! 0.5!(A83)$);
\draw[color13, thick] ($(A82)! 0.5!(A83)$) -- ($(A97)! 0.5!(A83)$);
\draw[color13, thick] ($(A97)! 0.5!(A83)$) -- ($(A98)! 0.5!(A83)$);
\draw[color12, thick] ($(A68)! 0.5!(A76)$) -- ($(A73)! 0.5!(A76)$);
\draw[color12, thick] ($(A68)! 0.5!(A76)$) -- ($(A78)! 0.5!(A76)$);
\draw[color12, thick] ($(A73)! 0.5!(A76)$) -- ($(A77)! 0.5!(A76)$);
\draw[color12, thick] ($(A77)! 0.5!(A76)$) -- ($(A79)! 0.5!(A76)$);
\draw[color12, thick] ($(A78)! 0.5!(A76)$) -- ($(A79)! 0.5!(A76)$);
\draw[color10, thick] ($(A63)! 0.5!(A69)$) -- ($(A64)! 0.5!(A69)$);
\draw[color10, thick] ($(A63)! 0.5!(A69)$) -- ($(A68)! 0.5!(A69)$);
\draw[color10, thick] ($(A64)! 0.5!(A69)$) -- ($(A70)! 0.5!(A69)$);
\draw[color10, thick] ($(A68)! 0.5!(A69)$) -- ($(A73)! 0.5!(A69)$);
\draw[color10, thick] ($(A70)! 0.5!(A69)$) -- ($(A74)! 0.5!(A69)$);
\draw[color10, thick] ($(A73)! 0.5!(A69)$) -- ($(A74)! 0.5!(A69)$);
\draw[color7, thick] ($(A50)! 0.5!(A58)$) -- ($(A51)! 0.5!(A58)$);
\draw[color7, thick] ($(A50)! 0.5!(A58)$) -- ($(A57)! 0.5!(A58)$);
\draw[color7, thick] ($(A51)! 0.5!(A58)$) -- ($(A59)! 0.5!(A58)$);
\draw[color7, thick] ($(A57)! 0.5!(A58)$) -- ($(A64)! 0.5!(A58)$);
\draw[color7, thick] ($(A59)! 0.5!(A58)$) -- ($(A65)! 0.5!(A58)$);
\draw[color7, thick] ($(A64)! 0.5!(A58)$) -- ($(A65)! 0.5!(A58)$);
\draw[color1, thick] ($(A0)! 0.5!(A88)$) -- ($(A1)! 0.5!(A88)$);
\draw[color1, thick] ($(A0)! 0.5!(A45)$) -- ($(A9)! 0.5!(A45)$);
\draw[color1, thick] ($(A45)! 0.5!(A0)$) -- ($(A88)! 0.5!(A0)$);
\draw[color1, thick] ($(A1)! 0.5!(A10)$) -- ($(A2)! 0.5!(A10)$);
\draw[color1, thick] ($(A1)! 0.5!(A10)$) -- ($(A9)! 0.5!(A10)$);
\draw[color1, thick] ($(A1)! 0.5!(A88)$) -- ($(A89)! 0.5!(A88)$);
\draw[color1, thick] ($(A2)! 0.5!(A90)$) -- ($(A3)! 0.5!(A90)$);
\draw[color1, thick] ($(A2)! 0.5!(A10)$) -- ($(A11)! 0.5!(A10)$);
\draw[color1, thick] ($(A2)! 0.5!(A90)$) -- ($(A89)! 0.5!(A90)$);
\draw[color1, thick] ($(A3)! 0.5!(A12)$) -- ($(A4)! 0.5!(A12)$);
\draw[color1, thick] ($(A3)! 0.5!(A12)$) -- ($(A11)! 0.5!(A12)$);
\draw[color1, thick] ($(A3)! 0.5!(A90)$) -- ($(A91)! 0.5!(A90)$);
\draw[color1, thick] ($(A4)! 0.5!(A92)$) -- ($(A5)! 0.5!(A92)$);
\draw[color1, thick] ($(A4)! 0.5!(A12)$) -- ($(A13)! 0.5!(A12)$);
\draw[color1, thick] ($(A4)! 0.5!(A92)$) -- ($(A91)! 0.5!(A92)$);
\draw[color1, thick] ($(A5)! 0.5!(A6)$) -- ($(A14)! 0.5!(A6)$);
\draw[color1, thick] ($(A6)! 0.5!(A5)$) -- ($(A93)! 0.5!(A5)$);
\draw[color1, thick] ($(A5)! 0.5!(A21)$) -- ($(A14)! 0.5!(A21)$);
\draw[color1, thick] ($(A5)! 0.5!(A21)$) -- ($(A27)! 0.5!(A21)$);
\draw[color1, thick] ($(A5)! 0.5!(A92)$) -- ($(A126)! 0.5!(A92)$);
\draw[color1, thick] ($(A93)! 0.5!(A5)$) -- ($(A106)! 0.5!(A5)$);
\draw[color1, thick] ($(A106)! 0.5!(A5)$) -- ($(A117)! 0.5!(A5)$);
\draw[color1, thick] ($(A117)! 0.5!(A5)$) -- ($(A142)! 0.5!(A5)$);
\draw[color1, thick] ($(A5)! 0.5!(A142)$) -- ($(A126)! 0.5!(A142)$);
\draw[color1, thick] ($(A7)! 0.5!(A6)$) -- ($(A15)! 0.5!(A6)$);
\draw[color1, thick] ($(A6)! 0.5!(A7)$) -- ($(A94)! 0.5!(A7)$);
\draw[color1, thick] ($(A14)! 0.5!(A6)$) -- ($(A15)! 0.5!(A6)$);
\draw[color1, thick] ($(A7)! 0.5!(A8)$) -- ($(A16)! 0.5!(A8)$);
\draw[color1, thick] ($(A8)! 0.5!(A7)$) -- ($(A95)! 0.5!(A7)$);
\draw[color1, thick] ($(A94)! 0.5!(A7)$) -- ($(A95)! 0.5!(A7)$);
\draw[color1, thick] ($(A9)! 0.5!(A10)$) -- ($(A17)! 0.5!(A10)$);
\draw[color1, thick] ($(A9)! 0.5!(A45)$) -- ($(A53)! 0.5!(A45)$);
\draw[color1, thick] ($(A11)! 0.5!(A10)$) -- ($(A18)! 0.5!(A10)$);
\draw[color1, thick] ($(A17)! 0.5!(A10)$) -- ($(A18)! 0.5!(A10)$);
\draw[color1, thick] ($(A11)! 0.5!(A12)$) -- ($(A19)! 0.5!(A12)$);
\draw[color1, thick] ($(A13)! 0.5!(A12)$) -- ($(A20)! 0.5!(A12)$);
\draw[color1, thick] ($(A19)! 0.5!(A12)$) -- ($(A20)! 0.5!(A12)$);
\draw[color1, thick] ($(A14)! 0.5!(A21)$) -- ($(A22)! 0.5!(A21)$);
\draw[color1, thick] ($(A15)! 0.5!(A23)$) -- ($(A16)! 0.5!(A23)$);
\draw[color1, thick] ($(A15)! 0.5!(A23)$) -- ($(A22)! 0.5!(A23)$);
\draw[color1, thick] ($(A17)! 0.5!(A60)$) -- ($(A24)! 0.5!(A60)$);
\draw[color1, thick] ($(A17)! 0.5!(A60)$) -- ($(A53)! 0.5!(A60)$);
\draw[color1, thick] ($(A18)! 0.5!(A25)$) -- ($(A19)! 0.5!(A25)$);
\draw[color1, thick] ($(A18)! 0.5!(A25)$) -- ($(A24)! 0.5!(A25)$);
\draw[color1, thick] ($(A19)! 0.5!(A25)$) -- ($(A26)! 0.5!(A25)$);
\draw[color1, thick] ($(A22)! 0.5!(A21)$) -- ($(A28)! 0.5!(A21)$);
\draw[color1, thick] ($(A27)! 0.5!(A21)$) -- ($(A28)! 0.5!(A21)$);
\draw[color1, thick] ($(A22)! 0.5!(A23)$) -- ($(A29)! 0.5!(A23)$);
\draw[color1, thick] ($(A24)! 0.5!(A25)$) -- ($(A30)! 0.5!(A25)$);
\draw[color1, thick] ($(A24)! 0.5!(A60)$) -- ($(A66)! 0.5!(A60)$);
\draw[color1, thick] ($(A26)! 0.5!(A25)$) -- ($(A31)! 0.5!(A25)$);
\draw[color1, thick] ($(A30)! 0.5!(A25)$) -- ($(A31)! 0.5!(A25)$);
\draw[color1, thick] ($(A28)! 0.5!(A34)$) -- ($(A29)! 0.5!(A34)$);
\draw[color1, thick] ($(A28)! 0.5!(A34)$) -- ($(A33)! 0.5!(A34)$);
\draw[color1, thick] ($(A30)! 0.5!(A71)$) -- ($(A35)! 0.5!(A71)$);
\draw[color1, thick] ($(A30)! 0.5!(A71)$) -- ($(A66)! 0.5!(A71)$);
\draw[color1, thick] ($(A31)! 0.5!(A36)$) -- ($(A32)! 0.5!(A36)$);
\draw[color1, thick] ($(A31)! 0.5!(A36)$) -- ($(A35)! 0.5!(A36)$);
\draw[color1, thick] ($(A32)! 0.5!(A36)$) -- ($(A42)! 0.5!(A36)$);
\draw[color1, thick] ($(A33)! 0.5!(A34)$) -- ($(A38)! 0.5!(A34)$);
\draw[color1, thick] ($(A35)! 0.5!(A36)$) -- ($(A39)! 0.5!(A36)$);
\draw[color1, thick] ($(A35)! 0.5!(A71)$) -- ($(A75)! 0.5!(A71)$);
\draw[color1, thick] ($(A39)! 0.5!(A36)$) -- ($(A42)! 0.5!(A36)$);
\draw[color1, thick] ($(A37)! 0.5!(A41)$) -- ($(A38)! 0.5!(A41)$);
\draw[color1, thick] ($(A37)! 0.5!(A41)$) -- ($(A40)! 0.5!(A41)$);
\draw[color1, thick] ($(A40)! 0.5!(A41)$) -- ($(A43)! 0.5!(A41)$);
\draw[color1, thick] ($(A42)! 0.5!(A44)$) -- ($(A43)! 0.5!(A44)$);
\draw[color1, thick] ($(A44)! 0.5!(A42)$) -- ($(A80)! 0.5!(A42)$);
\draw[color1, thick] ($(A49)! 0.5!(A42)$) -- ($(A68)! 0.5!(A42)$);
\draw[color1, thick] ($(A42)! 0.5!(A49)$) -- ($(A72)! 0.5!(A49)$);
\draw[color1, thick] ($(A68)! 0.5!(A42)$) -- ($(A78)! 0.5!(A42)$);
\draw[color1, thick] ($(A78)! 0.5!(A42)$) -- ($(A80)! 0.5!(A42)$);
\draw[color1, thick] ($(A46)! 0.5!(A45)$) -- ($(A53)! 0.5!(A45)$);
\draw[color1, thick] ($(A45)! 0.5!(A46)$) -- ($(A87)! 0.5!(A46)$);
\draw[color1, thick] ($(A46)! 0.5!(A47)$) -- ($(A54)! 0.5!(A47)$);
\draw[color1, thick] ($(A47)! 0.5!(A46)$) -- ($(A86)! 0.5!(A46)$);
\draw[color1, thick] ($(A86)! 0.5!(A46)$) -- ($(A87)! 0.5!(A46)$);
\draw[color1, thick] ($(A48)! 0.5!(A47)$) -- ($(A55)! 0.5!(A47)$);
\draw[color1, thick] ($(A47)! 0.5!(A48)$) -- ($(A85)! 0.5!(A48)$);
\draw[color1, thick] ($(A54)! 0.5!(A47)$) -- ($(A55)! 0.5!(A47)$);
\draw[color1, thick] ($(A48)! 0.5!(A49)$) -- ($(A56)! 0.5!(A49)$);
\draw[color1, thick] ($(A49)! 0.5!(A48)$) -- ($(A84)! 0.5!(A48)$);
\draw[color1, thick] ($(A84)! 0.5!(A48)$) -- ($(A85)! 0.5!(A48)$);
\draw[color1, thick] ($(A56)! 0.5!(A49)$) -- ($(A72)! 0.5!(A49)$);
\draw[color1, thick] ($(A51)! 0.5!(A81)$) -- ($(A52)! 0.5!(A81)$);
\draw[color1, thick] ($(A51)! 0.5!(A81)$) -- ($(A82)! 0.5!(A81)$);
\draw[color1, thick] ($(A53)! 0.5!(A60)$) -- ($(A54)! 0.5!(A60)$);
\draw[color1, thick] ($(A54)! 0.5!(A60)$) -- ($(A61)! 0.5!(A60)$);
\draw[color1, thick] ($(A55)! 0.5!(A62)$) -- ($(A56)! 0.5!(A62)$);
\draw[color1, thick] ($(A55)! 0.5!(A62)$) -- ($(A61)! 0.5!(A62)$);
\draw[color1, thick] ($(A56)! 0.5!(A62)$) -- ($(A72)! 0.5!(A62)$);
\draw[color1, thick] ($(A61)! 0.5!(A60)$) -- ($(A66)! 0.5!(A60)$);
\draw[color1, thick] ($(A61)! 0.5!(A62)$) -- ($(A67)! 0.5!(A62)$);
\draw[color1, thick] ($(A67)! 0.5!(A62)$) -- ($(A72)! 0.5!(A62)$);
\draw[color1, thick] ($(A66)! 0.5!(A71)$) -- ($(A67)! 0.5!(A71)$);
\draw[color1, thick] ($(A67)! 0.5!(A71)$) -- ($(A72)! 0.5!(A71)$);
\draw[color1, thick] ($(A72)! 0.5!(A71)$) -- ($(A75)! 0.5!(A71)$);
\draw[color1, thick] ($(A82)! 0.5!(A81)$) -- ($(A96)! 0.5!(A81)$);
\draw[color1, thick] ($(A87)! 0.5!(A102)$) -- ($(A88)! 0.5!(A102)$);
\draw[color1, thick] ($(A87)! 0.5!(A102)$) -- ($(A101)! 0.5!(A102)$);
\draw[color1, thick] ($(A89)! 0.5!(A88)$) -- ($(A102)! 0.5!(A88)$);
\draw[color1, thick] ($(A89)! 0.5!(A90)$) -- ($(A103)! 0.5!(A90)$);
\draw[color1, thick] ($(A91)! 0.5!(A90)$) -- ($(A104)! 0.5!(A90)$);
\draw[color1, thick] ($(A103)! 0.5!(A90)$) -- ($(A104)! 0.5!(A90)$);
\draw[color1, thick] ($(A91)! 0.5!(A92)$) -- ($(A105)! 0.5!(A92)$);
\draw[color1, thick] ($(A105)! 0.5!(A92)$) -- ($(A126)! 0.5!(A92)$);
\draw[color1, thick] ($(A96)! 0.5!(A109)$) -- ($(A97)! 0.5!(A109)$);
\draw[color1, thick] ($(A97)! 0.5!(A109)$) -- ($(A110)! 0.5!(A109)$);
\draw[color1, thick] ($(A101)! 0.5!(A102)$) -- ($(A114)! 0.5!(A102)$);
\draw[color1, thick] ($(A102)! 0.5!(A114)$) -- ($(A103)! 0.5!(A114)$);
\draw[color1, thick] ($(A103)! 0.5!(A114)$) -- ($(A115)! 0.5!(A114)$);
\draw[color1, thick] ($(A104)! 0.5!(A116)$) -- ($(A105)! 0.5!(A116)$);
\draw[color1, thick] ($(A104)! 0.5!(A116)$) -- ($(A115)! 0.5!(A116)$);
\draw[color1, thick] ($(A105)! 0.5!(A116)$) -- ($(A126)! 0.5!(A116)$);
\draw[color1, thick] ($(A110)! 0.5!(A109)$) -- ($(A120)! 0.5!(A109)$);
\draw[color1, thick] ($(A113)! 0.5!(A124)$) -- ($(A114)! 0.5!(A124)$);
\draw[color1, thick] ($(A113)! 0.5!(A124)$) -- ($(A123)! 0.5!(A124)$);
\draw[color1, thick] ($(A115)! 0.5!(A114)$) -- ($(A124)! 0.5!(A114)$);
\draw[color1, thick] ($(A115)! 0.5!(A116)$) -- ($(A125)! 0.5!(A116)$);
\draw[color1, thick] ($(A125)! 0.5!(A116)$) -- ($(A126)! 0.5!(A116)$);
\draw[color1, thick] ($(A120)! 0.5!(A129)$) -- ($(A121)! 0.5!(A129)$);
\draw[color1, thick] ($(A121)! 0.5!(A129)$) -- ($(A130)! 0.5!(A129)$);
\draw[color1, thick] ($(A123)! 0.5!(A124)$) -- ($(A132)! 0.5!(A124)$);
\draw[color1, thick] ($(A124)! 0.5!(A132)$) -- ($(A125)! 0.5!(A132)$);
\draw[color1, thick] ($(A125)! 0.5!(A132)$) -- ($(A133)! 0.5!(A132)$);
\draw[color1, thick] ($(A126)! 0.5!(A142)$) -- ($(A133)! 0.5!(A142)$);
\draw[color1, thick] ($(A130)! 0.5!(A129)$) -- ($(A136)! 0.5!(A129)$);
\draw[color1, thick] ($(A131)! 0.5!(A138)$) -- ($(A132)! 0.5!(A138)$);
\draw[color1, thick] ($(A131)! 0.5!(A138)$) -- ($(A142)! 0.5!(A138)$);
\draw[color1, thick] ($(A133)! 0.5!(A132)$) -- ($(A138)! 0.5!(A132)$);
\draw[color1, thick] ($(A133)! 0.5!(A142)$) -- ($(A138)! 0.5!(A142)$);
\draw[color1, thick] ($(A136)! 0.5!(A141)$) -- ($(A137)! 0.5!(A141)$);
\draw[color1, thick] ($(A137)! 0.5!(A141)$) -- ($(A142)! 0.5!(A141)$);
\draw[color1, thick] ($(A142)! 0.5!(A141)$) -- ($(A144)! 0.5!(A141)$);
\draw[color19, thick] ($(A104)! 0.5!(A116)$) -- ($(A105)! 0.5!(A116)$);
\draw[color19, thick] ($(A104)! 0.5!(A116)$) -- ($(A115)! 0.5!(A116)$);
\draw[color19, thick] ($(A105)! 0.5!(A116)$) -- ($(A126)! 0.5!(A116)$);
\draw[color19, thick] ($(A115)! 0.5!(A116)$) -- ($(A125)! 0.5!(A116)$);
\draw[color19, thick] ($(A125)! 0.5!(A116)$) -- ($(A126)! 0.5!(A116)$);
\draw[color15, thick] ($(A4)! 0.5!(A92)$) -- ($(A5)! 0.5!(A92)$);
\draw[color15, thick] ($(A4)! 0.5!(A92)$) -- ($(A91)! 0.5!(A92)$);
\draw[color15, thick] ($(A5)! 0.5!(A92)$) -- ($(A126)! 0.5!(A92)$);
\draw[color15, thick] ($(A91)! 0.5!(A92)$) -- ($(A105)! 0.5!(A92)$);
\draw[color15, thick] ($(A105)! 0.5!(A92)$) -- ($(A126)! 0.5!(A92)$);
\draw[color14, thick] ($(A2)! 0.5!(A90)$) -- ($(A3)! 0.5!(A90)$);
\draw[color14, thick] ($(A2)! 0.5!(A90)$) -- ($(A89)! 0.5!(A90)$);
\draw[color14, thick] ($(A3)! 0.5!(A90)$) -- ($(A91)! 0.5!(A90)$);
\draw[color14, thick] ($(A89)! 0.5!(A90)$) -- ($(A103)! 0.5!(A90)$);
\draw[color14, thick] ($(A91)! 0.5!(A90)$) -- ($(A104)! 0.5!(A90)$);
\draw[color14, thick] ($(A103)! 0.5!(A90)$) -- ($(A104)! 0.5!(A90)$);
\draw[color11, thick] ($(A30)! 0.5!(A71)$) -- ($(A35)! 0.5!(A71)$);
\draw[color11, thick] ($(A30)! 0.5!(A71)$) -- ($(A66)! 0.5!(A71)$);
\draw[color11, thick] ($(A35)! 0.5!(A71)$) -- ($(A75)! 0.5!(A71)$);
\draw[color11, thick] ($(A66)! 0.5!(A71)$) -- ($(A67)! 0.5!(A71)$);
\draw[color11, thick] ($(A67)! 0.5!(A71)$) -- ($(A72)! 0.5!(A71)$);
\draw[color11, thick] ($(A72)! 0.5!(A71)$) -- ($(A75)! 0.5!(A71)$);
\draw[color9, thick] ($(A55)! 0.5!(A62)$) -- ($(A56)! 0.5!(A62)$);
\draw[color9, thick] ($(A55)! 0.5!(A62)$) -- ($(A61)! 0.5!(A62)$);
\draw[color9, thick] ($(A56)! 0.5!(A62)$) -- ($(A72)! 0.5!(A62)$);
\draw[color9, thick] ($(A61)! 0.5!(A62)$) -- ($(A67)! 0.5!(A62)$);
\draw[color9, thick] ($(A67)! 0.5!(A62)$) -- ($(A72)! 0.5!(A62)$);
\draw[color8, thick] ($(A17)! 0.5!(A60)$) -- ($(A24)! 0.5!(A60)$);
\draw[color8, thick] ($(A17)! 0.5!(A60)$) -- ($(A53)! 0.5!(A60)$);
\draw[color8, thick] ($(A24)! 0.5!(A60)$) -- ($(A66)! 0.5!(A60)$);
\draw[color8, thick] ($(A53)! 0.5!(A60)$) -- ($(A54)! 0.5!(A60)$);
\draw[color8, thick] ($(A54)! 0.5!(A60)$) -- ($(A61)! 0.5!(A60)$);
\draw[color8, thick] ($(A61)! 0.5!(A60)$) -- ($(A66)! 0.5!(A60)$);
\draw[color6, thick] ($(A31)! 0.5!(A36)$) -- ($(A32)! 0.5!(A36)$);
\draw[color6, thick] ($(A31)! 0.5!(A36)$) -- ($(A35)! 0.5!(A36)$);
\draw[color6, thick] ($(A32)! 0.5!(A36)$) -- ($(A42)! 0.5!(A36)$);
\draw[color6, thick] ($(A35)! 0.5!(A36)$) -- ($(A39)! 0.5!(A36)$);
\draw[color6, thick] ($(A39)! 0.5!(A36)$) -- ($(A42)! 0.5!(A36)$);
\draw[color5, thick] ($(A18)! 0.5!(A25)$) -- ($(A19)! 0.5!(A25)$);
\draw[color5, thick] ($(A18)! 0.5!(A25)$) -- ($(A24)! 0.5!(A25)$);
\draw[color5, thick] ($(A19)! 0.5!(A25)$) -- ($(A26)! 0.5!(A25)$);
\draw[color5, thick] ($(A24)! 0.5!(A25)$) -- ($(A30)! 0.5!(A25)$);
\draw[color5, thick] ($(A26)! 0.5!(A25)$) -- ($(A31)! 0.5!(A25)$);
\draw[color5, thick] ($(A30)! 0.5!(A25)$) -- ($(A31)! 0.5!(A25)$);
\draw[color4, thick] ($(A5)! 0.5!(A21)$) -- ($(A14)! 0.5!(A21)$);
\draw[color4, thick] ($(A5)! 0.5!(A21)$) -- ($(A27)! 0.5!(A21)$);
\draw[color4, thick] ($(A14)! 0.5!(A21)$) -- ($(A22)! 0.5!(A21)$);
\draw[color4, thick] ($(A22)! 0.5!(A21)$) -- ($(A28)! 0.5!(A21)$);
\draw[color4, thick] ($(A27)! 0.5!(A21)$) -- ($(A28)! 0.5!(A21)$);
\draw[color3, thick] ($(A3)! 0.5!(A12)$) -- ($(A4)! 0.5!(A12)$);
\draw[color3, thick] ($(A3)! 0.5!(A12)$) -- ($(A11)! 0.5!(A12)$);
\draw[color3, thick] ($(A4)! 0.5!(A12)$) -- ($(A13)! 0.5!(A12)$);
\draw[color3, thick] ($(A11)! 0.5!(A12)$) -- ($(A19)! 0.5!(A12)$);
\draw[color3, thick] ($(A13)! 0.5!(A12)$) -- ($(A20)! 0.5!(A12)$);
\draw[color3, thick] ($(A19)! 0.5!(A12)$) -- ($(A20)! 0.5!(A12)$);
\draw[color2, thick] ($(A1)! 0.5!(A10)$) -- ($(A2)! 0.5!(A10)$);
\draw[color2, thick] ($(A1)! 0.5!(A10)$) -- ($(A9)! 0.5!(A10)$);
\draw[color2, thick] ($(A2)! 0.5!(A10)$) -- ($(A11)! 0.5!(A10)$);
\draw[color2, thick] ($(A9)! 0.5!(A10)$) -- ($(A17)! 0.5!(A10)$);
\draw[color2, thick] ($(A11)! 0.5!(A10)$) -- ($(A18)! 0.5!(A10)$);
\draw[color2, thick] ($(A17)! 0.5!(A10)$) -- ($(A18)! 0.5!(A10)$);
\foreach \a/\b/\c in {0/1/9,0/1/88,0/9/45,0/45/88,1/2/10,1/2/89,1/9/10,1/88/89,2/3/11,2/3/90,2/10/11,2/89/90,3/4/12,3/4/91,3/11/12,3/90/91,4/5/13,4/5/92,4/12/13,4/91/92,5/6/14,5/6/93,5/13/32,5/14/21,5/21/27,5/27/42,5/32/42,5/92/126,5/93/106,5/106/117,5/117/142,5/126/142,6/7/15,6/7/94,6/14/15,6/93/94,7/8/16,7/8/95,7/15/16,7/94/95,9/10/17,9/17/53,9/45/53,10/11/18,10/17/18,11/12/19,11/18/19,12/13/20,12/19/20,13/20/32,14/15/22,14/21/22,15/16/23,15/22/23,17/18/24,17/24/60,17/53/60,18/19/25,18/24/25,19/20/26,19/25/26,20/26/32,21/22/28,21/27/28,22/23/29,22/28/29,24/25/30,24/30/66,24/60/66,25/26/31,25/30/31,26/31/32,27/28/33,27/33/37,27/37/40,27/40/42,28/29/34,28/33/34,30/31/35,30/35/71,30/66/71,31/32/36,31/35/36,32/36/42,33/34/38,33/37/38,35/36/39,35/39/75,35/71/75,36/39/42,37/38/41,37/40/41,39/42/75,40/41/43,40/42/43,42/43/44,42/44/80,42/49/68,42/49/72,42/68/78,42/72/75,42/78/80,45/46/53,45/46/87,45/87/88,46/47/54,46/47/86,46/53/54,46/86/87,47/48/55,47/48/85,47/54/55,47/85/86,48/49/56,48/49/84,48/55/56,48/84/85,49/50/57,49/50/83,49/56/72,49/57/63,49/63/68,49/83/98,49/84/122,49/98/111,49/111/142,49/122/142,50/51/58,50/51/82,50/57/58,50/82/83,51/52/59,51/52/81,51/58/59,51/81/82,53/54/60,54/55/61,54/60/61,55/56/62,55/61/62,56/62/72,57/58/64,57/63/64,58/59/65,58/64/65,60/61/66,61/62/67,61/66/67,62/67/72,63/64/69,63/68/69,64/65/70,64/69/70,66/67/71,67/71/72,68/69/73,68/73/76,68/76/78,69/70/74,69/73/74,71/72/75,73/74/77,73/76/77,76/77/79,76/78/79,78/79/80,81/82/96,82/83/97,82/96/97,83/97/98,84/85/99,84/99/122,85/86/100,85/99/100,86/87/101,86/100/101,87/88/102,87/101/102,88/89/102,89/90/103,89/102/103,90/91/104,90/103/104,91/92/105,91/104/105,92/105/126,93/94/107,93/106/107,94/95/108,94/107/108,96/97/109,97/98/110,97/109/110,98/110/111,99/100/112,99/112/122,100/101/113,100/112/113,101/102/114,101/113/114,102/103/114,103/104/115,103/114/115,104/105/116,104/115/116,105/116/126,106/107/118,106/117/118,107/108/119,107/118/119,109/110/120,110/111/121,110/120/121,111/121/130,111/130/137,111/137/142,112/113/123,112/122/123,113/114/124,113/123/124,114/115/124,115/116/125,115/124/125,116/125/126,117/118/127,117/127/134,117/134/139,117/139/142,118/119/128,118/127/128,120/121/129,121/129/130,122/123/131,122/131/142,123/124/132,123/131/132,124/125/132,125/126/133,125/132/133,126/133/142,127/128/135,127/134/135,129/130/136,130/136/137,131/132/138,131/138/142,132/133/138,133/138/142,134/135/140,134/139/140,136/137/141,137/141/142,139/140/143,139/142/143,141/142/144,142/143/144}{
  \draw[black!30] (A\a) -- (A\b) -- (A\c) -- cycle;
}

\draw[black, thick] (A5) -- (A42) -- (A49) -- (A142) --cycle;

\foreach \i in {0,1,2,3,4,5,7,9,11,13,14,15,16,17,18,19,20,22,24,26,27,28,29,30,31,32,33,35,37,38,39,40,42,43,46,48,53,54,55,56,58,61,66,67,69,72,75,76,81,83,89,91,100,102,103,104,105,107,109,111,115,122,124,125,126,127,129,133,138,139,141}{
  \fill[myred] (A\i) circle (3pt);
}
\foreach \i in {6,8,10,12,21,23,25,34,36,41,44,45,47,49,50,51,52,57,59,60,62,63,64,65,68,70,71,73,74,77,78,79,80,82,84,85,86,87,88,90,92,93,94,95,96,97,98,99,101,106,108,110,112,113,114,116,117,118,119,120,121,123,128,130,131,132,134,135,136,137,140,142,143,144}{
  \fill[myblue] (A\i) circle (3pt);
}

%\foreach \i in {-8,-6,-4,-2,0,2,4,6,8}{
%	\foreach \j in {-8,-6,-4,-2,0,2,4,6,8}{
%		\node[] at (\i,\j) {\input{tikz/even_even.tikz}};
%	}
%	\foreach \j in {-7,-5,-3,-1,1,3,5,7}{
%		\node[] at (\i,\j) {\input{tikz/even_odd.tikz}};
%	}
%}
%\foreach \i in {-7,-5,-3,-1,1,3,5,7}{
%	\foreach \j in {-8,-6,-4,-2,0,2,4,6,8}{
%		\node[rotate=90] at (\i,\j) {\input{tikz/even_odd.tikz}};
%	}
%	\foreach \j in {-7,-5,-3,-1,1,3,5,7}{
%		\node[] at (\i,\j) {\input{tikz/odd_odd.tikz}};
%	}
%}

%% file: tikz/double_split_4_var.tikz
\coordinate (A0) at (0, 0);
\coordinate (A1) at (0, 1);
\coordinate (A2) at (0, 2);
\coordinate (A3) at (0, 3);
\coordinate (A4) at (0, 4);
\coordinate (A5) at (0, 5);
\coordinate (A6) at (0, 6);
\coordinate (A7) at (0, 7);
\coordinate (A8) at (0, 8);
\coordinate (A9) at (1, 0);
\coordinate (A10) at (1, 1);
\coordinate (A11) at (1, 2);
\coordinate (A12) at (1, 3);
\coordinate (A13) at (1, 4);
\coordinate (A14) at (1, 5);
\coordinate (A15) at (1, 6);
\coordinate (A16) at (1, 7);
\coordinate (A17) at (2, 0);
\coordinate (A18) at (2, 1);
\coordinate (A19) at (2, 2);
\coordinate (A20) at (2, 3);
\coordinate (A21) at (2, 4);
\coordinate (A22) at (2, 5);
\coordinate (A23) at (2, 6);
\coordinate (A24) at (3, 0);
\coordinate (A25) at (3, 1);
\coordinate (A26) at (3, 2);
\coordinate (A27) at (3, 3);
\coordinate (A28) at (3, 4);
\coordinate (A29) at (3, 5);
\coordinate (A30) at (4, 0);
\coordinate (A31) at (4, 1);
\coordinate (A32) at (4, 2);
\coordinate (A33) at (4, 3);
\coordinate (A34) at (4, 4);
\coordinate (A35) at (5, 0);
\coordinate (A36) at (5, 1);
\coordinate (A37) at (5, 2);
\coordinate (A38) at (5, 3);
\coordinate (A39) at (6, 0);
\coordinate (A40) at (6, 1);
\coordinate (A41) at (6, 2);
\coordinate (A42) at (7, 0);
\coordinate (A43) at (7, 1);
\coordinate (A44) at (8, 0);
\coordinate (A45) at (0, -1);
\coordinate (A46) at (0, -2);
\coordinate (A47) at (0, -3);
\coordinate (A48) at (0, -4);
\coordinate (A49) at (0, -5);
\coordinate (A50) at (0, -6);
\coordinate (A51) at (0, -7);
\coordinate (A52) at (0, -8);
\coordinate (A53) at (1, -1);
\coordinate (A54) at (1, -2);
\coordinate (A55) at (1, -3);
\coordinate (A56) at (1, -4);
\coordinate (A57) at (1, -5);
\coordinate (A58) at (1, -6);
\coordinate (A59) at (1, -7);
\coordinate (A60) at (2, -1);
\coordinate (A61) at (2, -2);
\coordinate (A62) at (2, -3);
\coordinate (A63) at (2, -4);
\coordinate (A64) at (2, -5);
\coordinate (A65) at (2, -6);
\coordinate (A66) at (3, -1);
\coordinate (A67) at (3, -2);
\coordinate (A68) at (3, -3);
\coordinate (A69) at (3, -4);
\coordinate (A70) at (3, -5);
\coordinate (A71) at (4, -1);
\coordinate (A72) at (4, -2);
\coordinate (A73) at (4, -3);
\coordinate (A74) at (4, -4);
\coordinate (A75) at (5, -1);
\coordinate (A76) at (5, -2);
\coordinate (A77) at (5, -3);
\coordinate (A78) at (6, -1);
\coordinate (A79) at (6, -2);
\coordinate (A80) at (7, -1);
\coordinate (A81) at (-1, -7);
\coordinate (A82) at (-1, -6);
\coordinate (A83) at (-1, -5);
\coordinate (A84) at (-1, -4);
\coordinate (A85) at (-1, -3);
\coordinate (A86) at (-1, -2);
\coordinate (A87) at (-1, -1);
\coordinate (A88) at (-1, 0);
\coordinate (A89) at (-1, 1);
\coordinate (A90) at (-1, 2);
\coordinate (A91) at (-1, 3);
\coordinate (A92) at (-1, 4);
\coordinate (A93) at (-1, 5);
\coordinate (A94) at (-1, 6);
\coordinate (A95) at (-1, 7);
\coordinate (A96) at (-2, -6);
\coordinate (A97) at (-2, -5);
\coordinate (A98) at (-2, -4);
\coordinate (A99) at (-2, -3);
\coordinate (A100) at (-2, -2);
\coordinate (A101) at (-2, -1);
\coordinate (A102) at (-2, 0);
\coordinate (A103) at (-2, 1);
\coordinate (A104) at (-2, 2);
\coordinate (A105) at (-2, 3);
\coordinate (A106) at (-2, 4);
\coordinate (A107) at (-2, 5);
\coordinate (A108) at (-2, 6);
\coordinate (A109) at (-3, -5);
\coordinate (A110) at (-3, -4);
\coordinate (A111) at (-3, -3);
\coordinate (A112) at (-3, -2);
\coordinate (A113) at (-3, -1);
\coordinate (A114) at (-3, 0);
\coordinate (A115) at (-3, 1);
\coordinate (A116) at (-3, 2);
\coordinate (A117) at (-3, 3);
\coordinate (A118) at (-3, 4);
\coordinate (A119) at (-3, 5);
\coordinate (A120) at (-4, -4);
\coordinate (A121) at (-4, -3);
\coordinate (A122) at (-4, -2);
\coordinate (A123) at (-4, -1);
\coordinate (A124) at (-4, 0);
\coordinate (A125) at (-4, 1);
\coordinate (A126) at (-4, 2);
\coordinate (A127) at (-4, 3);
\coordinate (A128) at (-4, 4);
\coordinate (A129) at (-5, -3);
\coordinate (A130) at (-5, -2);
\coordinate (A131) at (-5, -1);
\coordinate (A132) at (-5, 0);
\coordinate (A133) at (-5, 1);
\coordinate (A134) at (-5, 2);
\coordinate (A135) at (-5, 3);
\coordinate (A136) at (-6, -2);
\coordinate (A137) at (-6, -1);
\coordinate (A138) at (-6, 0);
\coordinate (A139) at (-6, 1);
\coordinate (A140) at (-6, 2);
\coordinate (A141) at (-7, -1);
\coordinate (A142) at (-7, 0);
\coordinate (A143) at (-7, 1);
\coordinate (A144) at (-8, 0);
\colorlet{color1}{mycolor4}
\fill[color1!50] (A1) -- ($(A1)! 0.5!(A0)$) -- ($(A9)! 0.5!(A0)$) -- (A9) -- cycle;
\fill[color1!50] ($(A0)! 0.5!(A1)$) -- ($(A88)! 0.5!(A1)$) -- (A1) -- cycle;
\fill[color1!50] ($(A0)! 0.5!(A9)$) -- ($(A45)! 0.5!(A9)$) -- (A9) -- cycle;
\fill[color1!50] (A1) -- ($(A1)! 0.5!(A2)$) -- ($(A10)! 0.5!(A2)$) -- (A10) -- cycle;
\fill[color1!50] ($(A2)! 0.5!(A1)$) -- ($(A89)! 0.5!(A1)$) -- (A1) -- cycle;
\fill[color1!50] (A1) -- (A9) -- (A10) -- cycle;
\fill[color1!50] ($(A88)! 0.5!(A1)$) -- ($(A89)! 0.5!(A1)$) -- (A1) -- cycle;
\fill[color1!50] (A3) -- ($(A3)! 0.5!(A2)$) -- ($(A11)! 0.5!(A2)$) -- (A11) -- cycle;
\fill[color1!50] ($(A2)! 0.5!(A3)$) -- ($(A90)! 0.5!(A3)$) -- (A3) -- cycle;
\fill[color1!50] (A10) -- ($(A10)! 0.5!(A2)$) -- ($(A11)! 0.5!(A2)$) -- (A11) -- cycle;
\fill[color1!50] (A3) -- ($(A3)! 0.5!(A4)$) -- ($(A12)! 0.5!(A4)$) -- (A12) -- cycle;
\fill[color1!50] ($(A4)! 0.5!(A3)$) -- ($(A91)! 0.5!(A3)$) -- (A3) -- cycle;
\fill[color1!50] (A3) -- (A11) -- (A12) -- cycle;
\fill[color1!50] ($(A90)! 0.5!(A3)$) -- ($(A91)! 0.5!(A3)$) -- (A3) -- cycle;
\fill[color1!50] (A5) -- ($(A5)! 0.5!(A4)$) -- ($(A13)! 0.5!(A4)$) -- (A13) -- cycle;
\fill[color1!50] ($(A4)! 0.5!(A5)$) -- ($(A92)! 0.5!(A5)$) -- (A5) -- cycle;
\fill[color1!50] (A12) -- ($(A12)! 0.5!(A4)$) -- ($(A13)! 0.5!(A4)$) -- (A13) -- cycle;
\fill[color1!50] (A5) -- ($(A5)! 0.5!(A6)$) -- ($(A14)! 0.5!(A6)$) -- (A14) -- cycle;
\fill[color1!50] ($(A6)! 0.5!(A5)$) -- ($(A93)! 0.5!(A5)$) -- (A5) -- cycle;
\fill[color1!50] (A5) -- (A13) -- (A14) -- cycle;
\fill[color1!50] ($(A92)! 0.5!(A5)$) -- ($(A93)! 0.5!(A5)$) -- (A5) -- cycle;
\fill[color1!50] (A7) -- ($(A7)! 0.5!(A6)$) -- ($(A15)! 0.5!(A6)$) -- (A15) -- cycle;
\fill[color1!50] ($(A6)! 0.5!(A7)$) -- ($(A94)! 0.5!(A7)$) -- (A7) -- cycle;
\fill[color1!50] (A14) -- ($(A14)! 0.5!(A6)$) -- ($(A15)! 0.5!(A6)$) -- (A15) -- cycle;
\fill[color1!50] (A7) -- ($(A7)! 0.5!(A8)$) -- ($(A16)! 0.5!(A8)$) -- (A16) -- cycle;
\fill[color1!50] ($(A8)! 0.5!(A7)$) -- ($(A95)! 0.5!(A7)$) -- (A7) -- cycle;
\fill[color1!50] (A7) -- (A15) -- (A16) -- cycle;
\fill[color1!50] ($(A94)! 0.5!(A7)$) -- ($(A95)! 0.5!(A7)$) -- (A7) -- cycle;
\fill[color1!50] (A9) -- ($(A9)! 0.5!(A17)$) -- ($(A10)! 0.5!(A17)$) -- (A10) -- cycle;
\fill[color1!50] ($(A17)! 0.5!(A9)$) -- ($(A53)! 0.5!(A9)$) -- (A9) -- cycle;
\fill[color1!50] ($(A45)! 0.5!(A9)$) -- ($(A53)! 0.5!(A9)$) -- (A9) -- cycle;
\fill[color1!50] (A10) -- (A11) -- (A18) -- cycle;
\fill[color1!50] (A10) -- ($(A10)! 0.5!(A17)$) -- ($(A18)! 0.5!(A17)$) -- (A18) -- cycle;
\fill[color1!50] (A11) -- ($(A11)! 0.5!(A19)$) -- ($(A12)! 0.5!(A19)$) -- (A12) -- cycle;
\fill[color1!50] (A11) -- ($(A11)! 0.5!(A19)$) -- ($(A18)! 0.5!(A19)$) -- (A18) -- cycle;
\fill[color1!50] (A12) -- (A13) -- (A20) -- cycle;
\fill[color1!50] (A12) -- ($(A12)! 0.5!(A19)$) -- ($(A20)! 0.5!(A19)$) -- (A20) -- cycle;
\fill[color1!50] (A13) -- ($(A13)! 0.5!(A21)$) -- ($(A14)! 0.5!(A21)$) -- (A14) -- cycle;
\fill[color1!50] (A13) -- ($(A13)! 0.5!(A21)$) -- ($(A20)! 0.5!(A21)$) -- (A20) -- cycle;
\fill[color1!50] (A14) -- (A15) -- (A22) -- cycle;
\fill[color1!50] (A14) -- ($(A14)! 0.5!(A21)$) -- ($(A22)! 0.5!(A21)$) -- (A22) -- cycle;
\fill[color1!50] (A15) -- ($(A15)! 0.5!(A23)$) -- ($(A16)! 0.5!(A23)$) -- (A16) -- cycle;
\fill[color1!50] (A15) -- ($(A15)! 0.5!(A23)$) -- ($(A22)! 0.5!(A23)$) -- (A22) -- cycle;
\fill[color1!50] (A18) -- ($(A18)! 0.5!(A17)$) -- ($(A24)! 0.5!(A17)$) -- (A24) -- cycle;
\fill[color1!50] ($(A17)! 0.5!(A24)$) -- ($(A60)! 0.5!(A24)$) -- (A24) -- cycle;
\fill[color1!50] (A18) -- ($(A18)! 0.5!(A19)$) -- ($(A24)! 0.5!(A19)$) -- (A24) -- cycle;
\fill[color1!50] (A20) -- ($(A20)! 0.5!(A19)$) -- ($(A24)! 0.5!(A19)$) -- (A24) -- cycle;
\fill[color1!50] (A20) -- ($(A20)! 0.5!(A21)$) -- ($(A27)! 0.5!(A21)$) -- (A27) -- cycle;
\fill[color1!50] (A20) -- ($(A20)! 0.5!(A25)$) -- ($(A24)! 0.5!(A25)$) -- (A24) -- cycle;
\fill[color1!50] (A20) -- ($(A20)! 0.5!(A25)$) -- ($(A26)! 0.5!(A25)$) -- (A26) -- cycle;
\fill[color1!50] (A20) -- ($(A20)! 0.5!(A36)$) -- ($(A26)! 0.5!(A36)$) -- (A26) -- cycle;
\fill[color1!50] (A20) -- ($(A20)! 0.5!(A32)$) -- ($(A27)! 0.5!(A32)$) -- (A27) -- cycle;
\fill[color1!50] (A20) -- ($(A20)! 0.5!(A32)$) -- ($(A42)! 0.5!(A32)$) -- (A42) -- cycle;
\fill[color1!50] (A20) -- ($(A20)! 0.5!(A36)$) -- ($(A42)! 0.5!(A36)$) -- (A42) -- cycle;
\fill[color1!50] (A22) -- ($(A22)! 0.5!(A21)$) -- ($(A28)! 0.5!(A21)$) -- (A28) -- cycle;
\fill[color1!50] (A27) -- ($(A27)! 0.5!(A21)$) -- ($(A28)! 0.5!(A21)$) -- (A28) -- cycle;
\fill[color1!50] (A22) -- ($(A22)! 0.5!(A23)$) -- ($(A29)! 0.5!(A23)$) -- (A29) -- cycle;
\fill[color1!50] (A22) -- (A28) -- (A29) -- cycle;
\fill[color1!50] (A24) -- ($(A24)! 0.5!(A25)$) -- ($(A30)! 0.5!(A25)$) -- (A30) -- cycle;
\fill[color1!50] (A24) -- (A30) -- (A66) -- cycle;
\fill[color1!50] ($(A60)! 0.5!(A24)$) -- ($(A61)! 0.5!(A24)$) -- (A24) -- cycle;
\fill[color1!50] ($(A61)! 0.5!(A24)$) -- ($(A62)! 0.5!(A24)$) -- (A24) -- cycle;
\fill[color1!50] (A24) -- ($(A24)! 0.5!(A62)$) -- ($(A66)! 0.5!(A62)$) -- (A66) -- cycle;
\fill[color1!50] (A26) -- ($(A26)! 0.5!(A25)$) -- ($(A31)! 0.5!(A25)$) -- (A31) -- cycle;
\fill[color1!50] (A30) -- ($(A30)! 0.5!(A25)$) -- ($(A31)! 0.5!(A25)$) -- (A31) -- cycle;
\fill[color1!50] (A26) -- ($(A26)! 0.5!(A36)$) -- ($(A31)! 0.5!(A36)$) -- (A31) -- cycle;
\fill[color1!50] (A27) -- (A28) -- (A33) -- cycle;
\fill[color1!50] (A27) -- ($(A27)! 0.5!(A32)$) -- ($(A33)! 0.5!(A32)$) -- (A33) -- cycle;
\fill[color1!50] (A28) -- ($(A28)! 0.5!(A34)$) -- ($(A29)! 0.5!(A34)$) -- (A29) -- cycle;
\fill[color1!50] (A28) -- ($(A28)! 0.5!(A34)$) -- ($(A33)! 0.5!(A34)$) -- (A33) -- cycle;
\fill[color1!50] (A30) -- (A31) -- (A35) -- cycle;
\fill[color1!50] (A30) -- ($(A30)! 0.5!(A71)$) -- ($(A35)! 0.5!(A71)$) -- (A35) -- cycle;
\fill[color1!50] (A30) -- ($(A30)! 0.5!(A71)$) -- ($(A66)! 0.5!(A71)$) -- (A66) -- cycle;
\fill[color1!50] (A31) -- ($(A31)! 0.5!(A36)$) -- ($(A35)! 0.5!(A36)$) -- (A35) -- cycle;
\fill[color1!50] (A33) -- ($(A33)! 0.5!(A32)$) -- ($(A37)! 0.5!(A32)$) -- (A37) -- cycle;
\fill[color1!50] (A37) -- ($(A37)! 0.5!(A32)$) -- ($(A40)! 0.5!(A32)$) -- (A40) -- cycle;
\fill[color1!50] (A40) -- ($(A40)! 0.5!(A32)$) -- ($(A42)! 0.5!(A32)$) -- (A42) -- cycle;
\fill[color1!50] (A33) -- ($(A33)! 0.5!(A34)$) -- ($(A38)! 0.5!(A34)$) -- (A38) -- cycle;
\fill[color1!50] (A33) -- (A37) -- (A38) -- cycle;
\fill[color1!50] (A35) -- ($(A35)! 0.5!(A36)$) -- ($(A39)! 0.5!(A36)$) -- (A39) -- cycle;
\fill[color1!50] (A35) -- (A39) -- (A75) -- cycle;
\fill[color1!50] (A35) -- ($(A35)! 0.5!(A71)$) -- ($(A75)! 0.5!(A71)$) -- (A75) -- cycle;
\fill[color1!50] (A39) -- ($(A39)! 0.5!(A36)$) -- ($(A42)! 0.5!(A36)$) -- (A42) -- cycle;
\fill[color1!50] (A37) -- ($(A37)! 0.5!(A41)$) -- ($(A38)! 0.5!(A41)$) -- (A38) -- cycle;
\fill[color1!50] (A37) -- ($(A37)! 0.5!(A41)$) -- ($(A40)! 0.5!(A41)$) -- (A40) -- cycle;
\fill[color1!50] (A39) -- (A42) -- (A75) -- cycle;
\fill[color1!50] (A40) -- ($(A40)! 0.5!(A41)$) -- ($(A43)! 0.5!(A41)$) -- (A43) -- cycle;
\fill[color1!50] (A40) -- (A42) -- (A43) -- cycle;
\fill[color1!50] (A42) -- ($(A42)! 0.5!(A44)$) -- ($(A43)! 0.5!(A44)$) -- (A43) -- cycle;
\fill[color1!50] ($(A44)! 0.5!(A42)$) -- ($(A80)! 0.5!(A42)$) -- (A42) -- cycle;
\fill[color1!50] ($(A62)! 0.5!(A42)$) -- ($(A72)! 0.5!(A42)$) -- (A42) -- cycle;
\fill[color1!50] (A42) -- ($(A42)! 0.5!(A62)$) -- ($(A75)! 0.5!(A62)$) -- (A75) -- cycle;
\fill[color1!50] ($(A72)! 0.5!(A42)$) -- ($(A78)! 0.5!(A42)$) -- (A42) -- cycle;
\fill[color1!50] ($(A78)! 0.5!(A42)$) -- ($(A80)! 0.5!(A42)$) -- (A42) -- cycle;
\fill[color1!50] ($(A51)! 0.5!(A81)$) -- ($(A52)! 0.5!(A81)$) -- (A81) -- cycle;
\fill[color1!50] ($(A51)! 0.5!(A81)$) -- ($(A82)! 0.5!(A81)$) -- (A81) -- cycle;
\fill[color1!50] (A66) -- ($(A66)! 0.5!(A62)$) -- ($(A67)! 0.5!(A62)$) -- (A67) -- cycle;
\fill[color1!50] (A67) -- ($(A67)! 0.5!(A62)$) -- ($(A75)! 0.5!(A62)$) -- (A75) -- cycle;
\fill[color1!50] (A66) -- ($(A66)! 0.5!(A71)$) -- ($(A67)! 0.5!(A71)$) -- (A67) -- cycle;
\fill[color1!50] (A67) -- ($(A67)! 0.5!(A71)$) -- ($(A75)! 0.5!(A71)$) -- (A75) -- cycle;
\fill[color1!50] ($(A82)! 0.5!(A81)$) -- ($(A96)! 0.5!(A81)$) -- (A81) -- cycle;
\fill[color1!50] ($(A96)! 0.5!(A109)$) -- ($(A97)! 0.5!(A109)$) -- (A109) -- cycle;
\fill[color1!50] ($(A97)! 0.5!(A109)$) -- ($(A110)! 0.5!(A109)$) -- (A109) -- cycle;
\fill[color1!50] ($(A110)! 0.5!(A109)$) -- ($(A120)! 0.5!(A109)$) -- (A109) -- cycle;
\fill[color1!50] ($(A120)! 0.5!(A129)$) -- ($(A121)! 0.5!(A129)$) -- (A129) -- cycle;
\fill[color1!50] ($(A121)! 0.5!(A129)$) -- ($(A130)! 0.5!(A129)$) -- (A129) -- cycle;
\fill[color1!50] ($(A130)! 0.5!(A129)$) -- ($(A136)! 0.5!(A129)$) -- (A129) -- cycle;
\fill[color1!50] ($(A136)! 0.5!(A141)$) -- ($(A137)! 0.5!(A141)$) -- (A141) -- cycle;
\fill[color1!50] ($(A137)! 0.5!(A141)$) -- ($(A142)! 0.5!(A141)$) -- (A141) -- cycle;
\fill[color1!50] ($(A142)! 0.5!(A141)$) -- ($(A144)! 0.5!(A141)$) -- (A141) -- cycle;
\colorlet{color2}{mycolor3}
\fill[color2!50] ($(A11)! 0.5!(A19)$) -- ($(A12)! 0.5!(A19)$) -- (A19) -- cycle;
\fill[color2!50] ($(A11)! 0.5!(A19)$) -- ($(A18)! 0.5!(A19)$) -- (A19) -- cycle;
\fill[color2!50] ($(A12)! 0.5!(A19)$) -- ($(A20)! 0.5!(A19)$) -- (A19) -- cycle;
\fill[color2!50] ($(A18)! 0.5!(A19)$) -- ($(A24)! 0.5!(A19)$) -- (A19) -- cycle;
\fill[color2!50] ($(A20)! 0.5!(A19)$) -- ($(A24)! 0.5!(A19)$) -- (A19) -- cycle;
\colorlet{color3}{mycolor3}
\fill[color3!50] ($(A13)! 0.5!(A21)$) -- ($(A14)! 0.5!(A21)$) -- (A21) -- cycle;
\fill[color3!50] ($(A13)! 0.5!(A21)$) -- ($(A20)! 0.5!(A21)$) -- (A21) -- cycle;
\fill[color3!50] ($(A14)! 0.5!(A21)$) -- ($(A22)! 0.5!(A21)$) -- (A21) -- cycle;
\fill[color3!50] ($(A20)! 0.5!(A21)$) -- ($(A27)! 0.5!(A21)$) -- (A21) -- cycle;
\fill[color3!50] ($(A22)! 0.5!(A21)$) -- ($(A28)! 0.5!(A21)$) -- (A21) -- cycle;
\fill[color3!50] ($(A27)! 0.5!(A21)$) -- ($(A28)! 0.5!(A21)$) -- (A21) -- cycle;
\colorlet{color4}{mycolor5}
\fill[color4!50] ($(A20)! 0.5!(A25)$) -- ($(A24)! 0.5!(A25)$) -- (A25) -- cycle;
\fill[color4!50] ($(A20)! 0.5!(A25)$) -- ($(A26)! 0.5!(A25)$) -- (A25) -- cycle;
\fill[color4!50] ($(A24)! 0.5!(A25)$) -- ($(A30)! 0.5!(A25)$) -- (A25) -- cycle;
\fill[color4!50] ($(A26)! 0.5!(A25)$) -- ($(A31)! 0.5!(A25)$) -- (A25) -- cycle;
\fill[color4!50] ($(A30)! 0.5!(A25)$) -- ($(A31)! 0.5!(A25)$) -- (A25) -- cycle;
\colorlet{color5}{mycolor3}
\fill[color5!50] ($(A20)! 0.5!(A32)$) -- ($(A27)! 0.5!(A32)$) -- (A32) -- cycle;
\fill[color5!50] ($(A20)! 0.5!(A32)$) -- ($(A42)! 0.5!(A32)$) -- (A32) -- cycle;
\fill[color5!50] ($(A27)! 0.5!(A32)$) -- ($(A33)! 0.5!(A32)$) -- (A32) -- cycle;
\fill[color5!50] ($(A33)! 0.5!(A32)$) -- ($(A37)! 0.5!(A32)$) -- (A32) -- cycle;
\fill[color5!50] ($(A37)! 0.5!(A32)$) -- ($(A40)! 0.5!(A32)$) -- (A32) -- cycle;
\fill[color5!50] ($(A40)! 0.5!(A32)$) -- ($(A42)! 0.5!(A32)$) -- (A32) -- cycle;
\colorlet{color6}{mycolor5}
\fill[color6!50] ($(A20)! 0.5!(A36)$) -- ($(A26)! 0.5!(A36)$) -- (A36) -- cycle;
\fill[color6!50] ($(A20)! 0.5!(A36)$) -- ($(A42)! 0.5!(A36)$) -- (A36) -- cycle;
\fill[color6!50] ($(A26)! 0.5!(A36)$) -- ($(A31)! 0.5!(A36)$) -- (A36) -- cycle;
\fill[color6!50] ($(A31)! 0.5!(A36)$) -- ($(A35)! 0.5!(A36)$) -- (A36) -- cycle;
\fill[color6!50] ($(A35)! 0.5!(A36)$) -- ($(A39)! 0.5!(A36)$) -- (A36) -- cycle;
\fill[color6!50] ($(A39)! 0.5!(A36)$) -- ($(A42)! 0.5!(A36)$) -- (A36) -- cycle;
\colorlet{color7}{mycolor3}
\fill[color7!50] ($(A46)! 0.5!(A54)$) -- ($(A47)! 0.5!(A54)$) -- (A54) -- cycle;
\fill[color7!50] ($(A46)! 0.5!(A54)$) -- ($(A53)! 0.5!(A54)$) -- (A54) -- cycle;
\fill[color7!50] ($(A47)! 0.5!(A54)$) -- ($(A55)! 0.5!(A54)$) -- (A54) -- cycle;
\fill[color7!50] ($(A53)! 0.5!(A54)$) -- ($(A60)! 0.5!(A54)$) -- (A54) -- cycle;
\fill[color7!50] ($(A55)! 0.5!(A54)$) -- ($(A61)! 0.5!(A54)$) -- (A54) -- cycle;
\fill[color7!50] ($(A60)! 0.5!(A54)$) -- ($(A61)! 0.5!(A54)$) -- (A54) -- cycle;
\colorlet{color8}{mycolor3}
\fill[color8!50] ($(A48)! 0.5!(A56)$) -- ($(A49)! 0.5!(A56)$) -- (A56) -- cycle;
\fill[color8!50] ($(A48)! 0.5!(A56)$) -- ($(A55)! 0.5!(A56)$) -- (A56) -- cycle;
\fill[color8!50] ($(A49)! 0.5!(A56)$) -- ($(A57)! 0.5!(A56)$) -- (A56) -- cycle;
\fill[color8!50] ($(A55)! 0.5!(A56)$) -- ($(A62)! 0.5!(A56)$) -- (A56) -- cycle;
\fill[color8!50] ($(A57)! 0.5!(A56)$) -- ($(A63)! 0.5!(A56)$) -- (A56) -- cycle;
\fill[color8!50] ($(A62)! 0.5!(A56)$) -- ($(A63)! 0.5!(A56)$) -- (A56) -- cycle;
\colorlet{color9}{mycolor3}
\fill[color9!50] ($(A50)! 0.5!(A58)$) -- ($(A51)! 0.5!(A58)$) -- (A58) -- cycle;
\fill[color9!50] ($(A50)! 0.5!(A58)$) -- ($(A57)! 0.5!(A58)$) -- (A58) -- cycle;
\fill[color9!50] ($(A51)! 0.5!(A58)$) -- ($(A59)! 0.5!(A58)$) -- (A58) -- cycle;
\fill[color9!50] ($(A57)! 0.5!(A58)$) -- ($(A64)! 0.5!(A58)$) -- (A58) -- cycle;
\fill[color9!50] ($(A59)! 0.5!(A58)$) -- ($(A65)! 0.5!(A58)$) -- (A58) -- cycle;
\fill[color9!50] ($(A64)! 0.5!(A58)$) -- ($(A65)! 0.5!(A58)$) -- (A58) -- cycle;
\colorlet{color10}{mycolor3}
\fill[color10!50] ($(A63)! 0.5!(A69)$) -- ($(A64)! 0.5!(A69)$) -- (A69) -- cycle;
\fill[color10!50] ($(A63)! 0.5!(A69)$) -- ($(A68)! 0.5!(A69)$) -- (A69) -- cycle;
\fill[color10!50] ($(A64)! 0.5!(A69)$) -- ($(A70)! 0.5!(A69)$) -- (A69) -- cycle;
\fill[color10!50] ($(A68)! 0.5!(A69)$) -- ($(A73)! 0.5!(A69)$) -- (A69) -- cycle;
\fill[color10!50] ($(A70)! 0.5!(A69)$) -- ($(A74)! 0.5!(A69)$) -- (A69) -- cycle;
\fill[color10!50] ($(A73)! 0.5!(A69)$) -- ($(A74)! 0.5!(A69)$) -- (A69) -- cycle;
\colorlet{color11}{mycolor5}
\fill[color11!50] ($(A30)! 0.5!(A71)$) -- ($(A35)! 0.5!(A71)$) -- (A71) -- cycle;
\fill[color11!50] ($(A30)! 0.5!(A71)$) -- ($(A66)! 0.5!(A71)$) -- (A71) -- cycle;
\fill[color11!50] ($(A35)! 0.5!(A71)$) -- ($(A75)! 0.5!(A71)$) -- (A71) -- cycle;
\fill[color11!50] ($(A66)! 0.5!(A71)$) -- ($(A67)! 0.5!(A71)$) -- (A71) -- cycle;
\fill[color11!50] ($(A67)! 0.5!(A71)$) -- ($(A75)! 0.5!(A71)$) -- (A71) -- cycle;
\colorlet{color12}{mycolor3}
\fill[color12!50] ($(A72)! 0.5!(A76)$) -- ($(A73)! 0.5!(A76)$) -- (A76) -- cycle;
\fill[color12!50] ($(A72)! 0.5!(A76)$) -- ($(A78)! 0.5!(A76)$) -- (A76) -- cycle;
\fill[color12!50] ($(A73)! 0.5!(A76)$) -- ($(A77)! 0.5!(A76)$) -- (A76) -- cycle;
\fill[color12!50] ($(A77)! 0.5!(A76)$) -- ($(A79)! 0.5!(A76)$) -- (A76) -- cycle;
\fill[color12!50] ($(A78)! 0.5!(A76)$) -- ($(A79)! 0.5!(A76)$) -- (A76) -- cycle;
\colorlet{color13}{mycolor3}
\fill[color13!50] ($(A49)! 0.5!(A83)$) -- ($(A50)! 0.5!(A83)$) -- (A83) -- cycle;
\fill[color13!50] ($(A49)! 0.5!(A83)$) -- ($(A84)! 0.5!(A83)$) -- (A83) -- cycle;
\fill[color13!50] ($(A50)! 0.5!(A83)$) -- ($(A82)! 0.5!(A83)$) -- (A83) -- cycle;
\fill[color13!50] ($(A82)! 0.5!(A83)$) -- ($(A97)! 0.5!(A83)$) -- (A83) -- cycle;
\fill[color13!50] ($(A84)! 0.5!(A83)$) -- ($(A98)! 0.5!(A83)$) -- (A83) -- cycle;
\fill[color13!50] ($(A97)! 0.5!(A83)$) -- ($(A98)! 0.5!(A83)$) -- (A83) -- cycle;
\colorlet{color14}{mycolor3}
\fill[color14!50] ($(A47)! 0.5!(A85)$) -- ($(A48)! 0.5!(A85)$) -- (A85) -- cycle;
\fill[color14!50] ($(A47)! 0.5!(A85)$) -- ($(A86)! 0.5!(A85)$) -- (A85) -- cycle;
\fill[color14!50] ($(A48)! 0.5!(A85)$) -- ($(A84)! 0.5!(A85)$) -- (A85) -- cycle;
\fill[color14!50] ($(A84)! 0.5!(A85)$) -- ($(A99)! 0.5!(A85)$) -- (A85) -- cycle;
\fill[color14!50] ($(A86)! 0.5!(A85)$) -- ($(A100)! 0.5!(A85)$) -- (A85) -- cycle;
\fill[color14!50] ($(A99)! 0.5!(A85)$) -- ($(A100)! 0.5!(A85)$) -- (A85) -- cycle;
\colorlet{color15}{mycolor3}
\fill[color15!50] ($(A45)! 0.5!(A87)$) -- ($(A46)! 0.5!(A87)$) -- (A87) -- cycle;
\fill[color15!50] ($(A45)! 0.5!(A87)$) -- ($(A88)! 0.5!(A87)$) -- (A87) -- cycle;
\fill[color15!50] ($(A46)! 0.5!(A87)$) -- ($(A86)! 0.5!(A87)$) -- (A87) -- cycle;
\fill[color15!50] ($(A86)! 0.5!(A87)$) -- ($(A101)! 0.5!(A87)$) -- (A87) -- cycle;
\fill[color15!50] ($(A88)! 0.5!(A87)$) -- ($(A102)! 0.5!(A87)$) -- (A87) -- cycle;
\fill[color15!50] ($(A101)! 0.5!(A87)$) -- ($(A102)! 0.5!(A87)$) -- (A87) -- cycle;
\colorlet{color16}{mycolor3}
\fill[color16!50] ($(A89)! 0.5!(A103)$) -- ($(A90)! 0.5!(A103)$) -- (A103) -- cycle;
\fill[color16!50] ($(A89)! 0.5!(A103)$) -- ($(A102)! 0.5!(A103)$) -- (A103) -- cycle;
\fill[color16!50] ($(A90)! 0.5!(A103)$) -- ($(A104)! 0.5!(A103)$) -- (A103) -- cycle;
\fill[color16!50] ($(A102)! 0.5!(A103)$) -- ($(A114)! 0.5!(A103)$) -- (A103) -- cycle;
\fill[color16!50] ($(A104)! 0.5!(A103)$) -- ($(A114)! 0.5!(A103)$) -- (A103) -- cycle;
\colorlet{color17}{mycolor3}
\fill[color17!50] ($(A91)! 0.5!(A105)$) -- ($(A92)! 0.5!(A105)$) -- (A105) -- cycle;
\fill[color17!50] ($(A91)! 0.5!(A105)$) -- ($(A104)! 0.5!(A105)$) -- (A105) -- cycle;
\fill[color17!50] ($(A92)! 0.5!(A105)$) -- ($(A106)! 0.5!(A105)$) -- (A105) -- cycle;
\fill[color17!50] ($(A104)! 0.5!(A105)$) -- ($(A114)! 0.5!(A105)$) -- (A105) -- cycle;
\fill[color17!50] ($(A106)! 0.5!(A105)$) -- ($(A117)! 0.5!(A105)$) -- (A105) -- cycle;
\fill[color17!50] (A105) -- ($(A105)! 0.5!(A114)$) -- ($(A115)! 0.5!(A114)$) -- (A115) -- cycle;
\fill[color17!50] (A105) -- ($(A105)! 0.5!(A116)$) -- ($(A115)! 0.5!(A116)$) -- (A115) -- cycle;
\fill[color17!50] (A105) -- ($(A105)! 0.5!(A116)$) -- ($(A133)! 0.5!(A116)$) -- (A133) -- cycle;
\fill[color17!50] ($(A117)! 0.5!(A105)$) -- ($(A126)! 0.5!(A105)$) -- (A105) -- cycle;
\fill[color17!50] ($(A126)! 0.5!(A105)$) -- ($(A142)! 0.5!(A105)$) -- (A105) -- cycle;
\fill[color17!50] (A105) -- ($(A105)! 0.5!(A142)$) -- ($(A133)! 0.5!(A142)$) -- (A133) -- cycle;
\fill[color17!50] ($(A113)! 0.5!(A124)$) -- ($(A114)! 0.5!(A124)$) -- (A124) -- cycle;
\fill[color17!50] ($(A113)! 0.5!(A124)$) -- ($(A123)! 0.5!(A124)$) -- (A124) -- cycle;
\fill[color17!50] (A115) -- ($(A115)! 0.5!(A114)$) -- ($(A124)! 0.5!(A114)$) -- (A124) -- cycle;
\fill[color17!50] (A115) -- ($(A115)! 0.5!(A116)$) -- ($(A125)! 0.5!(A116)$) -- (A125) -- cycle;
\fill[color17!50] (A115) -- (A124) -- (A125) -- cycle;
\fill[color17!50] (A125) -- ($(A125)! 0.5!(A116)$) -- ($(A133)! 0.5!(A116)$) -- (A133) -- cycle;
\fill[color17!50] ($(A123)! 0.5!(A124)$) -- ($(A132)! 0.5!(A124)$) -- (A124) -- cycle;
\fill[color17!50] (A124) -- ($(A124)! 0.5!(A132)$) -- ($(A125)! 0.5!(A132)$) -- (A125) -- cycle;
\fill[color17!50] (A125) -- ($(A125)! 0.5!(A132)$) -- ($(A133)! 0.5!(A132)$) -- (A133) -- cycle;
\fill[color17!50] ($(A131)! 0.5!(A138)$) -- ($(A132)! 0.5!(A138)$) -- (A138) -- cycle;
\fill[color17!50] ($(A131)! 0.5!(A138)$) -- ($(A142)! 0.5!(A138)$) -- (A138) -- cycle;
\fill[color17!50] (A133) -- ($(A133)! 0.5!(A132)$) -- ($(A138)! 0.5!(A132)$) -- (A138) -- cycle;
\fill[color17!50] (A133) -- ($(A133)! 0.5!(A142)$) -- ($(A138)! 0.5!(A142)$) -- (A138) -- cycle;
\colorlet{color18}{mycolor3}
\fill[color18!50] ($(A93)! 0.5!(A107)$) -- ($(A94)! 0.5!(A107)$) -- (A107) -- cycle;
\fill[color18!50] ($(A93)! 0.5!(A107)$) -- ($(A106)! 0.5!(A107)$) -- (A107) -- cycle;
\fill[color18!50] ($(A94)! 0.5!(A107)$) -- ($(A108)! 0.5!(A107)$) -- (A107) -- cycle;
\fill[color18!50] ($(A106)! 0.5!(A107)$) -- ($(A118)! 0.5!(A107)$) -- (A107) -- cycle;
\fill[color18!50] ($(A108)! 0.5!(A107)$) -- ($(A119)! 0.5!(A107)$) -- (A107) -- cycle;
\fill[color18!50] ($(A118)! 0.5!(A107)$) -- ($(A119)! 0.5!(A107)$) -- (A107) -- cycle;
\colorlet{color19}{mycolor3}
\fill[color19!50] ($(A98)! 0.5!(A111)$) -- ($(A99)! 0.5!(A111)$) -- (A111) -- cycle;
\fill[color19!50] ($(A98)! 0.5!(A111)$) -- ($(A110)! 0.5!(A111)$) -- (A111) -- cycle;
\fill[color19!50] ($(A99)! 0.5!(A111)$) -- ($(A122)! 0.5!(A111)$) -- (A111) -- cycle;
\fill[color19!50] ($(A110)! 0.5!(A111)$) -- ($(A121)! 0.5!(A111)$) -- (A111) -- cycle;
\fill[color19!50] ($(A121)! 0.5!(A111)$) -- ($(A122)! 0.5!(A111)$) -- (A111) -- cycle;
\colorlet{color20}{mycolor5}
\fill[color20!50] ($(A105)! 0.5!(A116)$) -- ($(A115)! 0.5!(A116)$) -- (A116) -- cycle;
\fill[color20!50] ($(A105)! 0.5!(A116)$) -- ($(A133)! 0.5!(A116)$) -- (A116) -- cycle;
\fill[color20!50] ($(A115)! 0.5!(A116)$) -- ($(A125)! 0.5!(A116)$) -- (A116) -- cycle;
\fill[color20!50] ($(A125)! 0.5!(A116)$) -- ($(A133)! 0.5!(A116)$) -- (A116) -- cycle;
\colorlet{color21}{mycolor3}
\fill[color21!50] ($(A117)! 0.5!(A127)$) -- ($(A118)! 0.5!(A127)$) -- (A127) -- cycle;
\fill[color21!50] ($(A117)! 0.5!(A127)$) -- ($(A126)! 0.5!(A127)$) -- (A127) -- cycle;
\fill[color21!50] ($(A118)! 0.5!(A127)$) -- ($(A128)! 0.5!(A127)$) -- (A127) -- cycle;
\fill[color21!50] ($(A126)! 0.5!(A127)$) -- ($(A134)! 0.5!(A127)$) -- (A127) -- cycle;
\fill[color21!50] ($(A128)! 0.5!(A127)$) -- ($(A135)! 0.5!(A127)$) -- (A127) -- cycle;
\fill[color21!50] ($(A134)! 0.5!(A127)$) -- ($(A135)! 0.5!(A127)$) -- (A127) -- cycle;
\colorlet{color22}{mycolor3}
\fill[color22!50] ($(A126)! 0.5!(A139)$) -- ($(A134)! 0.5!(A139)$) -- (A139) -- cycle;
\fill[color22!50] ($(A126)! 0.5!(A139)$) -- ($(A142)! 0.5!(A139)$) -- (A139) -- cycle;
\fill[color22!50] ($(A134)! 0.5!(A139)$) -- ($(A140)! 0.5!(A139)$) -- (A139) -- cycle;
\fill[color22!50] ($(A140)! 0.5!(A139)$) -- ($(A143)! 0.5!(A139)$) -- (A139) -- cycle;
\fill[color22!50] ($(A142)! 0.5!(A139)$) -- ($(A143)! 0.5!(A139)$) -- (A139) -- cycle;
\draw[color22, thick] ($(A126)! 0.5!(A139)$) -- ($(A134)! 0.5!(A139)$);
\draw[color22, thick] ($(A126)! 0.5!(A139)$) -- ($(A142)! 0.5!(A139)$);
\draw[color22, thick] ($(A134)! 0.5!(A139)$) -- ($(A140)! 0.5!(A139)$);
\draw[color22, thick] ($(A140)! 0.5!(A139)$) -- ($(A143)! 0.5!(A139)$);
\draw[color22, thick] ($(A142)! 0.5!(A139)$) -- ($(A143)! 0.5!(A139)$);
\draw[color21, thick] ($(A117)! 0.5!(A127)$) -- ($(A118)! 0.5!(A127)$);
\draw[color21, thick] ($(A117)! 0.5!(A127)$) -- ($(A126)! 0.5!(A127)$);
\draw[color21, thick] ($(A118)! 0.5!(A127)$) -- ($(A128)! 0.5!(A127)$);
\draw[color21, thick] ($(A126)! 0.5!(A127)$) -- ($(A134)! 0.5!(A127)$);
\draw[color21, thick] ($(A128)! 0.5!(A127)$) -- ($(A135)! 0.5!(A127)$);
\draw[color21, thick] ($(A134)! 0.5!(A127)$) -- ($(A135)! 0.5!(A127)$);
\draw[color19, thick] ($(A98)! 0.5!(A111)$) -- ($(A99)! 0.5!(A111)$);
\draw[color19, thick] ($(A98)! 0.5!(A111)$) -- ($(A110)! 0.5!(A111)$);
\draw[color19, thick] ($(A99)! 0.5!(A111)$) -- ($(A122)! 0.5!(A111)$);
\draw[color19, thick] ($(A110)! 0.5!(A111)$) -- ($(A121)! 0.5!(A111)$);
\draw[color19, thick] ($(A121)! 0.5!(A111)$) -- ($(A122)! 0.5!(A111)$);
\draw[color18, thick] ($(A93)! 0.5!(A107)$) -- ($(A94)! 0.5!(A107)$);
\draw[color18, thick] ($(A93)! 0.5!(A107)$) -- ($(A106)! 0.5!(A107)$);
\draw[color18, thick] ($(A94)! 0.5!(A107)$) -- ($(A108)! 0.5!(A107)$);
\draw[color18, thick] ($(A106)! 0.5!(A107)$) -- ($(A118)! 0.5!(A107)$);
\draw[color18, thick] ($(A108)! 0.5!(A107)$) -- ($(A119)! 0.5!(A107)$);
\draw[color18, thick] ($(A118)! 0.5!(A107)$) -- ($(A119)! 0.5!(A107)$);
\draw[color17, thick] ($(A91)! 0.5!(A105)$) -- ($(A92)! 0.5!(A105)$);
\draw[color17, thick] ($(A91)! 0.5!(A105)$) -- ($(A104)! 0.5!(A105)$);
\draw[color17, thick] ($(A92)! 0.5!(A105)$) -- ($(A106)! 0.5!(A105)$);
\draw[color17, thick] ($(A104)! 0.5!(A105)$) -- ($(A114)! 0.5!(A105)$);
\draw[color17, thick] ($(A106)! 0.5!(A105)$) -- ($(A117)! 0.5!(A105)$);
\draw[color17, thick] ($(A105)! 0.5!(A114)$) -- ($(A115)! 0.5!(A114)$);
\draw[color17, thick] ($(A105)! 0.5!(A116)$) -- ($(A115)! 0.5!(A116)$);
\draw[color17, thick] ($(A105)! 0.5!(A116)$) -- ($(A133)! 0.5!(A116)$);
\draw[color17, thick] ($(A117)! 0.5!(A105)$) -- ($(A126)! 0.5!(A105)$);
\draw[color17, thick] ($(A126)! 0.5!(A105)$) -- ($(A142)! 0.5!(A105)$);
\draw[color17, thick] ($(A105)! 0.5!(A142)$) -- ($(A133)! 0.5!(A142)$);
\draw[color17, thick] ($(A113)! 0.5!(A124)$) -- ($(A114)! 0.5!(A124)$);
\draw[color17, thick] ($(A113)! 0.5!(A124)$) -- ($(A123)! 0.5!(A124)$);
\draw[color17, thick] ($(A115)! 0.5!(A114)$) -- ($(A124)! 0.5!(A114)$);
\draw[color17, thick] ($(A115)! 0.5!(A116)$) -- ($(A125)! 0.5!(A116)$);
\draw[color17, thick] ($(A125)! 0.5!(A116)$) -- ($(A133)! 0.5!(A116)$);
\draw[color17, thick] ($(A123)! 0.5!(A124)$) -- ($(A132)! 0.5!(A124)$);
\draw[color17, thick] ($(A124)! 0.5!(A132)$) -- ($(A125)! 0.5!(A132)$);
\draw[color17, thick] ($(A125)! 0.5!(A132)$) -- ($(A133)! 0.5!(A132)$);
\draw[color17, thick] ($(A131)! 0.5!(A138)$) -- ($(A132)! 0.5!(A138)$);
\draw[color17, thick] ($(A131)! 0.5!(A138)$) -- ($(A142)! 0.5!(A138)$);
\draw[color17, thick] ($(A133)! 0.5!(A132)$) -- ($(A138)! 0.5!(A132)$);
\draw[color17, thick] ($(A133)! 0.5!(A142)$) -- ($(A138)! 0.5!(A142)$);
\draw[color20, thick] ($(A105)! 0.5!(A116)$) -- ($(A115)! 0.5!(A116)$);
\draw[color20, thick] ($(A105)! 0.5!(A116)$) -- ($(A133)! 0.5!(A116)$);
\draw[color20, thick] ($(A115)! 0.5!(A116)$) -- ($(A125)! 0.5!(A116)$);
\draw[color20, thick] ($(A125)! 0.5!(A116)$) -- ($(A133)! 0.5!(A116)$);
\draw[color16, thick] ($(A89)! 0.5!(A103)$) -- ($(A90)! 0.5!(A103)$);
\draw[color16, thick] ($(A89)! 0.5!(A103)$) -- ($(A102)! 0.5!(A103)$);
\draw[color16, thick] ($(A90)! 0.5!(A103)$) -- ($(A104)! 0.5!(A103)$);
\draw[color16, thick] ($(A102)! 0.5!(A103)$) -- ($(A114)! 0.5!(A103)$);
\draw[color16, thick] ($(A104)! 0.5!(A103)$) -- ($(A114)! 0.5!(A103)$);
\draw[color15, thick] ($(A45)! 0.5!(A87)$) -- ($(A46)! 0.5!(A87)$);
\draw[color15, thick] ($(A45)! 0.5!(A87)$) -- ($(A88)! 0.5!(A87)$);
\draw[color15, thick] ($(A46)! 0.5!(A87)$) -- ($(A86)! 0.5!(A87)$);
\draw[color15, thick] ($(A86)! 0.5!(A87)$) -- ($(A101)! 0.5!(A87)$);
\draw[color15, thick] ($(A88)! 0.5!(A87)$) -- ($(A102)! 0.5!(A87)$);
\draw[color15, thick] ($(A101)! 0.5!(A87)$) -- ($(A102)! 0.5!(A87)$);
\draw[color14, thick] ($(A47)! 0.5!(A85)$) -- ($(A48)! 0.5!(A85)$);
\draw[color14, thick] ($(A47)! 0.5!(A85)$) -- ($(A86)! 0.5!(A85)$);
\draw[color14, thick] ($(A48)! 0.5!(A85)$) -- ($(A84)! 0.5!(A85)$);
\draw[color14, thick] ($(A84)! 0.5!(A85)$) -- ($(A99)! 0.5!(A85)$);
\draw[color14, thick] ($(A86)! 0.5!(A85)$) -- ($(A100)! 0.5!(A85)$);
\draw[color14, thick] ($(A99)! 0.5!(A85)$) -- ($(A100)! 0.5!(A85)$);
\draw[color13, thick] ($(A49)! 0.5!(A83)$) -- ($(A50)! 0.5!(A83)$);
\draw[color13, thick] ($(A49)! 0.5!(A83)$) -- ($(A84)! 0.5!(A83)$);
\draw[color13, thick] ($(A50)! 0.5!(A83)$) -- ($(A82)! 0.5!(A83)$);
\draw[color13, thick] ($(A82)! 0.5!(A83)$) -- ($(A97)! 0.5!(A83)$);
\draw[color13, thick] ($(A84)! 0.5!(A83)$) -- ($(A98)! 0.5!(A83)$);
\draw[color13, thick] ($(A97)! 0.5!(A83)$) -- ($(A98)! 0.5!(A83)$);
\draw[color12, thick] ($(A72)! 0.5!(A76)$) -- ($(A73)! 0.5!(A76)$);
\draw[color12, thick] ($(A72)! 0.5!(A76)$) -- ($(A78)! 0.5!(A76)$);
\draw[color12, thick] ($(A73)! 0.5!(A76)$) -- ($(A77)! 0.5!(A76)$);
\draw[color12, thick] ($(A77)! 0.5!(A76)$) -- ($(A79)! 0.5!(A76)$);
\draw[color12, thick] ($(A78)! 0.5!(A76)$) -- ($(A79)! 0.5!(A76)$);
\draw[color10, thick] ($(A63)! 0.5!(A69)$) -- ($(A64)! 0.5!(A69)$);
\draw[color10, thick] ($(A63)! 0.5!(A69)$) -- ($(A68)! 0.5!(A69)$);
\draw[color10, thick] ($(A64)! 0.5!(A69)$) -- ($(A70)! 0.5!(A69)$);
\draw[color10, thick] ($(A68)! 0.5!(A69)$) -- ($(A73)! 0.5!(A69)$);
\draw[color10, thick] ($(A70)! 0.5!(A69)$) -- ($(A74)! 0.5!(A69)$);
\draw[color10, thick] ($(A73)! 0.5!(A69)$) -- ($(A74)! 0.5!(A69)$);
\draw[color9, thick] ($(A50)! 0.5!(A58)$) -- ($(A51)! 0.5!(A58)$);
\draw[color9, thick] ($(A50)! 0.5!(A58)$) -- ($(A57)! 0.5!(A58)$);
\draw[color9, thick] ($(A51)! 0.5!(A58)$) -- ($(A59)! 0.5!(A58)$);
\draw[color9, thick] ($(A57)! 0.5!(A58)$) -- ($(A64)! 0.5!(A58)$);
\draw[color9, thick] ($(A59)! 0.5!(A58)$) -- ($(A65)! 0.5!(A58)$);
\draw[color9, thick] ($(A64)! 0.5!(A58)$) -- ($(A65)! 0.5!(A58)$);
\draw[color8, thick] ($(A48)! 0.5!(A56)$) -- ($(A49)! 0.5!(A56)$);
\draw[color8, thick] ($(A48)! 0.5!(A56)$) -- ($(A55)! 0.5!(A56)$);
\draw[color8, thick] ($(A49)! 0.5!(A56)$) -- ($(A57)! 0.5!(A56)$);
\draw[color8, thick] ($(A55)! 0.5!(A56)$) -- ($(A62)! 0.5!(A56)$);
\draw[color8, thick] ($(A57)! 0.5!(A56)$) -- ($(A63)! 0.5!(A56)$);
\draw[color8, thick] ($(A62)! 0.5!(A56)$) -- ($(A63)! 0.5!(A56)$);
\draw[color7, thick] ($(A46)! 0.5!(A54)$) -- ($(A47)! 0.5!(A54)$);
\draw[color7, thick] ($(A46)! 0.5!(A54)$) -- ($(A53)! 0.5!(A54)$);
\draw[color7, thick] ($(A47)! 0.5!(A54)$) -- ($(A55)! 0.5!(A54)$);
\draw[color7, thick] ($(A53)! 0.5!(A54)$) -- ($(A60)! 0.5!(A54)$);
\draw[color7, thick] ($(A55)! 0.5!(A54)$) -- ($(A61)! 0.5!(A54)$);
\draw[color7, thick] ($(A60)! 0.5!(A54)$) -- ($(A61)! 0.5!(A54)$);
\draw[color1, thick] ($(A1)! 0.5!(A0)$) -- ($(A9)! 0.5!(A0)$);
\draw[color1, thick] ($(A0)! 0.5!(A1)$) -- ($(A88)! 0.5!(A1)$);
\draw[color1, thick] ($(A0)! 0.5!(A9)$) -- ($(A45)! 0.5!(A9)$);
\draw[color1, thick] ($(A1)! 0.5!(A2)$) -- ($(A10)! 0.5!(A2)$);
\draw[color1, thick] ($(A2)! 0.5!(A1)$) -- ($(A89)! 0.5!(A1)$);
\draw[color1, thick] ($(A88)! 0.5!(A1)$) -- ($(A89)! 0.5!(A1)$);
\draw[color1, thick] ($(A3)! 0.5!(A2)$) -- ($(A11)! 0.5!(A2)$);
\draw[color1, thick] ($(A2)! 0.5!(A3)$) -- ($(A90)! 0.5!(A3)$);
\draw[color1, thick] ($(A10)! 0.5!(A2)$) -- ($(A11)! 0.5!(A2)$);
\draw[color1, thick] ($(A3)! 0.5!(A4)$) -- ($(A12)! 0.5!(A4)$);
\draw[color1, thick] ($(A4)! 0.5!(A3)$) -- ($(A91)! 0.5!(A3)$);
\draw[color1, thick] ($(A90)! 0.5!(A3)$) -- ($(A91)! 0.5!(A3)$);
\draw[color1, thick] ($(A5)! 0.5!(A4)$) -- ($(A13)! 0.5!(A4)$);
\draw[color1, thick] ($(A4)! 0.5!(A5)$) -- ($(A92)! 0.5!(A5)$);
\draw[color1, thick] ($(A12)! 0.5!(A4)$) -- ($(A13)! 0.5!(A4)$);
\draw[color1, thick] ($(A5)! 0.5!(A6)$) -- ($(A14)! 0.5!(A6)$);
\draw[color1, thick] ($(A6)! 0.5!(A5)$) -- ($(A93)! 0.5!(A5)$);
\draw[color1, thick] ($(A92)! 0.5!(A5)$) -- ($(A93)! 0.5!(A5)$);
\draw[color1, thick] ($(A7)! 0.5!(A6)$) -- ($(A15)! 0.5!(A6)$);
\draw[color1, thick] ($(A6)! 0.5!(A7)$) -- ($(A94)! 0.5!(A7)$);
\draw[color1, thick] ($(A14)! 0.5!(A6)$) -- ($(A15)! 0.5!(A6)$);
\draw[color1, thick] ($(A7)! 0.5!(A8)$) -- ($(A16)! 0.5!(A8)$);
\draw[color1, thick] ($(A8)! 0.5!(A7)$) -- ($(A95)! 0.5!(A7)$);
\draw[color1, thick] ($(A94)! 0.5!(A7)$) -- ($(A95)! 0.5!(A7)$);
\draw[color1, thick] ($(A9)! 0.5!(A17)$) -- ($(A10)! 0.5!(A17)$);
\draw[color1, thick] ($(A17)! 0.5!(A9)$) -- ($(A53)! 0.5!(A9)$);
\draw[color1, thick] ($(A45)! 0.5!(A9)$) -- ($(A53)! 0.5!(A9)$);
\draw[color1, thick] ($(A10)! 0.5!(A17)$) -- ($(A18)! 0.5!(A17)$);
\draw[color1, thick] ($(A11)! 0.5!(A19)$) -- ($(A12)! 0.5!(A19)$);
\draw[color1, thick] ($(A11)! 0.5!(A19)$) -- ($(A18)! 0.5!(A19)$);
\draw[color1, thick] ($(A12)! 0.5!(A19)$) -- ($(A20)! 0.5!(A19)$);
\draw[color1, thick] ($(A13)! 0.5!(A21)$) -- ($(A14)! 0.5!(A21)$);
\draw[color1, thick] ($(A13)! 0.5!(A21)$) -- ($(A20)! 0.5!(A21)$);
\draw[color1, thick] ($(A14)! 0.5!(A21)$) -- ($(A22)! 0.5!(A21)$);
\draw[color1, thick] ($(A15)! 0.5!(A23)$) -- ($(A16)! 0.5!(A23)$);
\draw[color1, thick] ($(A15)! 0.5!(A23)$) -- ($(A22)! 0.5!(A23)$);
\draw[color1, thick] ($(A18)! 0.5!(A17)$) -- ($(A24)! 0.5!(A17)$);
\draw[color1, thick] ($(A17)! 0.5!(A24)$) -- ($(A60)! 0.5!(A24)$);
\draw[color1, thick] ($(A18)! 0.5!(A19)$) -- ($(A24)! 0.5!(A19)$);
\draw[color1, thick] ($(A20)! 0.5!(A19)$) -- ($(A24)! 0.5!(A19)$);
\draw[color1, thick] ($(A20)! 0.5!(A21)$) -- ($(A27)! 0.5!(A21)$);
\draw[color1, thick] ($(A20)! 0.5!(A25)$) -- ($(A24)! 0.5!(A25)$);
\draw[color1, thick] ($(A20)! 0.5!(A25)$) -- ($(A26)! 0.5!(A25)$);
\draw[color1, thick] ($(A20)! 0.5!(A36)$) -- ($(A26)! 0.5!(A36)$);
\draw[color1, thick] ($(A20)! 0.5!(A32)$) -- ($(A27)! 0.5!(A32)$);
\draw[color1, thick] ($(A20)! 0.5!(A32)$) -- ($(A42)! 0.5!(A32)$);
\draw[color1, thick] ($(A20)! 0.5!(A36)$) -- ($(A42)! 0.5!(A36)$);
\draw[color1, thick] ($(A22)! 0.5!(A21)$) -- ($(A28)! 0.5!(A21)$);
\draw[color1, thick] ($(A27)! 0.5!(A21)$) -- ($(A28)! 0.5!(A21)$);
\draw[color1, thick] ($(A22)! 0.5!(A23)$) -- ($(A29)! 0.5!(A23)$);
\draw[color1, thick] ($(A24)! 0.5!(A25)$) -- ($(A30)! 0.5!(A25)$);
\draw[color1, thick] ($(A60)! 0.5!(A24)$) -- ($(A61)! 0.5!(A24)$);
\draw[color1, thick] ($(A61)! 0.5!(A24)$) -- ($(A62)! 0.5!(A24)$);
\draw[color1, thick] ($(A24)! 0.5!(A62)$) -- ($(A66)! 0.5!(A62)$);
\draw[color1, thick] ($(A26)! 0.5!(A25)$) -- ($(A31)! 0.5!(A25)$);
\draw[color1, thick] ($(A30)! 0.5!(A25)$) -- ($(A31)! 0.5!(A25)$);
\draw[color1, thick] ($(A26)! 0.5!(A36)$) -- ($(A31)! 0.5!(A36)$);
\draw[color1, thick] ($(A27)! 0.5!(A32)$) -- ($(A33)! 0.5!(A32)$);
\draw[color1, thick] ($(A28)! 0.5!(A34)$) -- ($(A29)! 0.5!(A34)$);
\draw[color1, thick] ($(A28)! 0.5!(A34)$) -- ($(A33)! 0.5!(A34)$);
\draw[color1, thick] ($(A30)! 0.5!(A71)$) -- ($(A35)! 0.5!(A71)$);
\draw[color1, thick] ($(A30)! 0.5!(A71)$) -- ($(A66)! 0.5!(A71)$);
\draw[color1, thick] ($(A31)! 0.5!(A36)$) -- ($(A35)! 0.5!(A36)$);
\draw[color1, thick] ($(A33)! 0.5!(A32)$) -- ($(A37)! 0.5!(A32)$);
\draw[color1, thick] ($(A37)! 0.5!(A32)$) -- ($(A40)! 0.5!(A32)$);
\draw[color1, thick] ($(A40)! 0.5!(A32)$) -- ($(A42)! 0.5!(A32)$);
\draw[color1, thick] ($(A33)! 0.5!(A34)$) -- ($(A38)! 0.5!(A34)$);
\draw[color1, thick] ($(A35)! 0.5!(A36)$) -- ($(A39)! 0.5!(A36)$);
\draw[color1, thick] ($(A35)! 0.5!(A71)$) -- ($(A75)! 0.5!(A71)$);
\draw[color1, thick] ($(A39)! 0.5!(A36)$) -- ($(A42)! 0.5!(A36)$);
\draw[color1, thick] ($(A37)! 0.5!(A41)$) -- ($(A38)! 0.5!(A41)$);
\draw[color1, thick] ($(A37)! 0.5!(A41)$) -- ($(A40)! 0.5!(A41)$);
\draw[color1, thick] ($(A40)! 0.5!(A41)$) -- ($(A43)! 0.5!(A41)$);
\draw[color1, thick] ($(A42)! 0.5!(A44)$) -- ($(A43)! 0.5!(A44)$);
\draw[color1, thick] ($(A44)! 0.5!(A42)$) -- ($(A80)! 0.5!(A42)$);
\draw[color1, thick] ($(A62)! 0.5!(A42)$) -- ($(A72)! 0.5!(A42)$);
\draw[color1, thick] ($(A42)! 0.5!(A62)$) -- ($(A75)! 0.5!(A62)$);
\draw[color1, thick] ($(A72)! 0.5!(A42)$) -- ($(A78)! 0.5!(A42)$);
\draw[color1, thick] ($(A78)! 0.5!(A42)$) -- ($(A80)! 0.5!(A42)$);
\draw[color1, thick] ($(A51)! 0.5!(A81)$) -- ($(A52)! 0.5!(A81)$);
\draw[color1, thick] ($(A51)! 0.5!(A81)$) -- ($(A82)! 0.5!(A81)$);
\draw[color1, thick] ($(A66)! 0.5!(A62)$) -- ($(A67)! 0.5!(A62)$);
\draw[color1, thick] ($(A67)! 0.5!(A62)$) -- ($(A75)! 0.5!(A62)$);
\draw[color1, thick] ($(A66)! 0.5!(A71)$) -- ($(A67)! 0.5!(A71)$);
\draw[color1, thick] ($(A67)! 0.5!(A71)$) -- ($(A75)! 0.5!(A71)$);
\draw[color1, thick] ($(A82)! 0.5!(A81)$) -- ($(A96)! 0.5!(A81)$);
\draw[color1, thick] ($(A96)! 0.5!(A109)$) -- ($(A97)! 0.5!(A109)$);
\draw[color1, thick] ($(A97)! 0.5!(A109)$) -- ($(A110)! 0.5!(A109)$);
\draw[color1, thick] ($(A110)! 0.5!(A109)$) -- ($(A120)! 0.5!(A109)$);
\draw[color1, thick] ($(A120)! 0.5!(A129)$) -- ($(A121)! 0.5!(A129)$);
\draw[color1, thick] ($(A121)! 0.5!(A129)$) -- ($(A130)! 0.5!(A129)$);
\draw[color1, thick] ($(A130)! 0.5!(A129)$) -- ($(A136)! 0.5!(A129)$);
\draw[color1, thick] ($(A136)! 0.5!(A141)$) -- ($(A137)! 0.5!(A141)$);
\draw[color1, thick] ($(A137)! 0.5!(A141)$) -- ($(A142)! 0.5!(A141)$);
\draw[color1, thick] ($(A142)! 0.5!(A141)$) -- ($(A144)! 0.5!(A141)$);
\draw[color11, thick] ($(A30)! 0.5!(A71)$) -- ($(A35)! 0.5!(A71)$);
\draw[color11, thick] ($(A30)! 0.5!(A71)$) -- ($(A66)! 0.5!(A71)$);
\draw[color11, thick] ($(A35)! 0.5!(A71)$) -- ($(A75)! 0.5!(A71)$);
\draw[color11, thick] ($(A66)! 0.5!(A71)$) -- ($(A67)! 0.5!(A71)$);
\draw[color11, thick] ($(A67)! 0.5!(A71)$) -- ($(A75)! 0.5!(A71)$);
\draw[color6, thick] ($(A20)! 0.5!(A36)$) -- ($(A26)! 0.5!(A36)$);
\draw[color6, thick] ($(A20)! 0.5!(A36)$) -- ($(A42)! 0.5!(A36)$);
\draw[color6, thick] ($(A26)! 0.5!(A36)$) -- ($(A31)! 0.5!(A36)$);
\draw[color6, thick] ($(A31)! 0.5!(A36)$) -- ($(A35)! 0.5!(A36)$);
\draw[color6, thick] ($(A35)! 0.5!(A36)$) -- ($(A39)! 0.5!(A36)$);
\draw[color6, thick] ($(A39)! 0.5!(A36)$) -- ($(A42)! 0.5!(A36)$);
\draw[color5, thick] ($(A20)! 0.5!(A32)$) -- ($(A27)! 0.5!(A32)$);
\draw[color5, thick] ($(A20)! 0.5!(A32)$) -- ($(A42)! 0.5!(A32)$);
\draw[color5, thick] ($(A27)! 0.5!(A32)$) -- ($(A33)! 0.5!(A32)$);
\draw[color5, thick] ($(A33)! 0.5!(A32)$) -- ($(A37)! 0.5!(A32)$);
\draw[color5, thick] ($(A37)! 0.5!(A32)$) -- ($(A40)! 0.5!(A32)$);
\draw[color5, thick] ($(A40)! 0.5!(A32)$) -- ($(A42)! 0.5!(A32)$);
\draw[color4, thick] ($(A20)! 0.5!(A25)$) -- ($(A24)! 0.5!(A25)$);
\draw[color4, thick] ($(A20)! 0.5!(A25)$) -- ($(A26)! 0.5!(A25)$);
\draw[color4, thick] ($(A24)! 0.5!(A25)$) -- ($(A30)! 0.5!(A25)$);
\draw[color4, thick] ($(A26)! 0.5!(A25)$) -- ($(A31)! 0.5!(A25)$);
\draw[color4, thick] ($(A30)! 0.5!(A25)$) -- ($(A31)! 0.5!(A25)$);
\draw[color3, thick] ($(A13)! 0.5!(A21)$) -- ($(A14)! 0.5!(A21)$);
\draw[color3, thick] ($(A13)! 0.5!(A21)$) -- ($(A20)! 0.5!(A21)$);
\draw[color3, thick] ($(A14)! 0.5!(A21)$) -- ($(A22)! 0.5!(A21)$);
\draw[color3, thick] ($(A20)! 0.5!(A21)$) -- ($(A27)! 0.5!(A21)$);
\draw[color3, thick] ($(A22)! 0.5!(A21)$) -- ($(A28)! 0.5!(A21)$);
\draw[color3, thick] ($(A27)! 0.5!(A21)$) -- ($(A28)! 0.5!(A21)$);
\draw[color2, thick] ($(A11)! 0.5!(A19)$) -- ($(A12)! 0.5!(A19)$);
\draw[color2, thick] ($(A11)! 0.5!(A19)$) -- ($(A18)! 0.5!(A19)$);
\draw[color2, thick] ($(A12)! 0.5!(A19)$) -- ($(A20)! 0.5!(A19)$);
\draw[color2, thick] ($(A18)! 0.5!(A19)$) -- ($(A24)! 0.5!(A19)$);
\draw[color2, thick] ($(A20)! 0.5!(A19)$) -- ($(A24)! 0.5!(A19)$);
\foreach \a/\b/\c in {0/1/9,0/1/88,0/9/45,0/45/88,1/2/10,1/2/89,1/9/10,1/88/89,2/3/11,2/3/90,2/10/11,2/89/90,3/4/12,3/4/91,3/11/12,3/90/91,4/5/13,4/5/92,4/12/13,4/91/92,5/6/14,5/6/93,5/13/14,5/92/93,6/7/15,6/7/94,6/14/15,6/93/94,7/8/16,7/8/95,7/15/16,7/94/95,9/10/17,9/17/53,9/45/53,10/11/18,10/17/18,11/12/19,11/18/19,12/13/20,12/19/20,13/14/21,13/20/21,14/15/22,14/21/22,15/16/23,15/22/23,17/18/24,17/24/60,17/53/60,18/19/24,19/20/24,20/21/27,20/24/25,20/25/26,20/26/36,20/27/32,20/32/42,20/36/42,21/22/28,21/27/28,22/23/29,22/28/29,24/25/30,24/30/66,24/60/61,24/61/62,24/62/66,25/26/31,25/30/31,26/31/36,27/28/33,27/32/33,28/29/34,28/33/34,30/31/35,30/35/71,30/66/71,31/35/36,32/33/37,32/37/40,32/40/42,33/34/38,33/37/38,35/36/39,35/39/75,35/71/75,36/39/42,37/38/41,37/40/41,39/42/75,40/41/43,40/42/43,42/43/44,42/44/80,42/62/72,42/62/75,42/72/78,42/78/80,45/46/53,45/46/87,45/87/88,46/47/54,46/47/86,46/53/54,46/86/87,47/48/55,47/48/85,47/54/55,47/85/86,48/49/56,48/49/84,48/55/56,48/84/85,49/50/57,49/50/83,49/56/57,49/83/84,50/51/58,50/51/82,50/57/58,50/82/83,51/52/59,51/52/81,51/58/59,51/81/82,53/54/60,54/55/61,54/60/61,55/56/62,55/61/62,56/57/63,56/62/63,57/58/64,57/63/64,58/59/65,58/64/65,62/63/68,62/66/67,62/67/75,62/68/72,63/64/69,63/68/69,64/65/70,64/69/70,66/67/71,67/71/75,68/69/73,68/72/73,69/70/74,69/73/74,72/73/76,72/76/78,73/74/77,73/76/77,76/77/79,76/78/79,78/79/80,81/82/96,82/83/97,82/96/97,83/84/98,83/97/98,84/85/99,84/98/99,85/86/100,85/99/100,86/87/101,86/100/101,87/88/102,87/101/102,88/89/102,89/90/103,89/102/103,90/91/104,90/103/104,91/92/105,91/104/105,92/93/106,92/105/106,93/94/107,93/106/107,94/95/108,94/107/108,96/97/109,97/98/110,97/109/110,98/99/111,98/110/111,99/100/114,99/111/122,99/112/113,99/112/131,99/113/114,99/122/142,99/131/142,100/101/114,101/102/114,102/103/114,103/104/114,104/105/114,105/106/117,105/114/115,105/115/116,105/116/133,105/117/126,105/126/142,105/133/142,106/107/118,106/117/118,107/108/119,107/118/119,109/110/120,110/111/121,110/120/121,111/121/122,112/113/123,112/123/131,113/114/124,113/123/124,114/115/124,115/116/125,115/124/125,116/125/133,117/118/127,117/126/127,118/119/128,118/127/128,120/121/129,121/122/130,121/129/130,122/130/137,122/137/142,123/124/132,123/131/132,124/125/132,125/132/133,126/127/134,126/134/139,126/139/142,127/128/135,127/134/135,129/130/136,130/136/137,131/132/138,131/138/142,132/133/138,133/138/142,134/135/140,134/139/140,136/137/141,137/141/142,139/140/143,139/142/143,141/142/144,142/143/144}{
  \draw[black!30] (A\a) -- (A\b) -- (A\c) -- cycle;
}

\draw[black, thick] (A20) -- (A42) -- (A62) -- (A24) --cycle;
\draw[black, thick] (A105) -- (A142) -- (A99) -- (A114) --cycle;

\foreach \i in {1,3,5,7,9,10,11,12,13,14,15,16,18,20,22,24,26,27,28,29,30,31,33,35,37,38,39,40,42,43,54,56,58,66,67,69,75,76,81,83,85,87,103,105,107,109,111,115,124,125,127,129,133,138,139,141}{
  \fill[myred] (A\i) circle (3pt);
}
\foreach \i in {0,2,4,6,8,17,19,21,23,25,32,34,36,41,44,45,46,47,48,49,50,51,52,53,55,57,59,60,61,62,63,64,65,68,70,71,72,73,74,77,78,79,80,82,84,86,88,89,90,91,92,93,94,95,96,97,98,99,100,101,102,104,106,108,110,112,113,114,116,117,118,119,120,121,122,123,126,128,130,131,132,134,135,136,137,140,142,143,144}{
  \fill[myblue] (A\i) circle (3pt);
}

%\foreach \i in {-8,-6,-4,-2,0,2,4,6,8}{
%	\foreach \j in {-8,-6,-4,-2,0,2,4,6,8}{
%		\node[] at (\i,\j) {\input{tikz/even_even.tikz}};
%	}
%	\foreach \j in {-7,-5,-3,-1,1,3,5,7}{
%		\node[] at (\i,\j) {\input{tikz/even_odd.tikz}};
%	}
%}
%\foreach \i in {-7,-5,-3,-1,1,3,5,7}{
%	\foreach \j in {-8,-6,-4,-2,0,2,4,6,8}{
%		\node[rotate=90] at (\i,\j) {\input{tikz/even_odd.tikz}};
%	}
%	\foreach \j in {-7,-5,-3,-1,1,3,5,7}{
%		\node[] at (\i,\j) {\input{tikz/odd_odd.tikz}};
%	}
%}

%% file: tikz/double_split_6.tikz
\coordinate (A0) at (0, 0);
\coordinate (A1) at (0, 1);
\coordinate (A2) at (0, 2);
\coordinate (A3) at (0, 3);
\coordinate (A4) at (0, 4);
\coordinate (A5) at (0, 5);
\coordinate (A6) at (0, 6);
\coordinate (A7) at (0, 7);
\coordinate (A8) at (0, 8);
\coordinate (A9) at (1, 0);
\coordinate (A10) at (1, 1);
\coordinate (A11) at (1, 2);
\coordinate (A12) at (1, 3);
\coordinate (A13) at (1, 4);
\coordinate (A14) at (1, 5);
\coordinate (A15) at (1, 6);
\coordinate (A16) at (1, 7);
\coordinate (A17) at (2, 0);
\coordinate (A18) at (2, 1);
\coordinate (A19) at (2, 2);
\coordinate (A20) at (2, 3);
\coordinate (A21) at (2, 4);
\coordinate (A22) at (2, 5);
\coordinate (A23) at (2, 6);
\coordinate (A24) at (3, 0);
\coordinate (A25) at (3, 1);
\coordinate (A26) at (3, 2);
\coordinate (A27) at (3, 3);
\coordinate (A28) at (3, 4);
\coordinate (A29) at (3, 5);
\coordinate (A30) at (4, 0);
\coordinate (A31) at (4, 1);
\coordinate (A32) at (4, 2);
\coordinate (A33) at (4, 3);
\coordinate (A34) at (4, 4);
\coordinate (A35) at (5, 0);
\coordinate (A36) at (5, 1);
\coordinate (A37) at (5, 2);
\coordinate (A38) at (5, 3);
\coordinate (A39) at (6, 0);
\coordinate (A40) at (6, 1);
\coordinate (A41) at (6, 2);
\coordinate (A42) at (7, 0);
\coordinate (A43) at (7, 1);
\coordinate (A44) at (8, 0);
\coordinate (A45) at (0, -1);
\coordinate (A46) at (0, -2);
\coordinate (A47) at (0, -3);
\coordinate (A48) at (0, -4);
\coordinate (A49) at (0, -5);
\coordinate (A50) at (0, -6);
\coordinate (A51) at (0, -7);
\coordinate (A52) at (0, -8);
\coordinate (A53) at (1, -1);
\coordinate (A54) at (1, -2);
\coordinate (A55) at (1, -3);
\coordinate (A56) at (1, -4);
\coordinate (A57) at (1, -5);
\coordinate (A58) at (1, -6);
\coordinate (A59) at (1, -7);
\coordinate (A60) at (2, -1);
\coordinate (A61) at (2, -2);
\coordinate (A62) at (2, -3);
\coordinate (A63) at (2, -4);
\coordinate (A64) at (2, -5);
\coordinate (A65) at (2, -6);
\coordinate (A66) at (3, -1);
\coordinate (A67) at (3, -2);
\coordinate (A68) at (3, -3);
\coordinate (A69) at (3, -4);
\coordinate (A70) at (3, -5);
\coordinate (A71) at (4, -1);
\coordinate (A72) at (4, -2);
\coordinate (A73) at (4, -3);
\coordinate (A74) at (4, -4);
\coordinate (A75) at (5, -1);
\coordinate (A76) at (5, -2);
\coordinate (A77) at (5, -3);
\coordinate (A78) at (6, -1);
\coordinate (A79) at (6, -2);
\coordinate (A80) at (7, -1);
\coordinate (A81) at (-1, -7);
\coordinate (A82) at (-1, -6);
\coordinate (A83) at (-1, -5);
\coordinate (A84) at (-1, -4);
\coordinate (A85) at (-1, -3);
\coordinate (A86) at (-1, -2);
\coordinate (A87) at (-1, -1);
\coordinate (A88) at (-1, 0);
\coordinate (A89) at (-1, 1);
\coordinate (A90) at (-1, 2);
\coordinate (A91) at (-1, 3);
\coordinate (A92) at (-1, 4);
\coordinate (A93) at (-1, 5);
\coordinate (A94) at (-1, 6);
\coordinate (A95) at (-1, 7);
\coordinate (A96) at (-2, -6);
\coordinate (A97) at (-2, -5);
\coordinate (A98) at (-2, -4);
\coordinate (A99) at (-2, -3);
\coordinate (A100) at (-2, -2);
\coordinate (A101) at (-2, -1);
\coordinate (A102) at (-2, 0);
\coordinate (A103) at (-2, 1);
\coordinate (A104) at (-2, 2);
\coordinate (A105) at (-2, 3);
\coordinate (A106) at (-2, 4);
\coordinate (A107) at (-2, 5);
\coordinate (A108) at (-2, 6);
\coordinate (A109) at (-3, -5);
\coordinate (A110) at (-3, -4);
\coordinate (A111) at (-3, -3);
\coordinate (A112) at (-3, -2);
\coordinate (A113) at (-3, -1);
\coordinate (A114) at (-3, 0);
\coordinate (A115) at (-3, 1);
\coordinate (A116) at (-3, 2);
\coordinate (A117) at (-3, 3);
\coordinate (A118) at (-3, 4);
\coordinate (A119) at (-3, 5);
\coordinate (A120) at (-4, -4);
\coordinate (A121) at (-4, -3);
\coordinate (A122) at (-4, -2);
\coordinate (A123) at (-4, -1);
\coordinate (A124) at (-4, 0);
\coordinate (A125) at (-4, 1);
\coordinate (A126) at (-4, 2);
\coordinate (A127) at (-4, 3);
\coordinate (A128) at (-4, 4);
\coordinate (A129) at (-5, -3);
\coordinate (A130) at (-5, -2);
\coordinate (A131) at (-5, -1);
\coordinate (A132) at (-5, 0);
\coordinate (A133) at (-5, 1);
\coordinate (A134) at (-5, 2);
\coordinate (A135) at (-5, 3);
\coordinate (A136) at (-6, -2);
\coordinate (A137) at (-6, -1);
\coordinate (A138) at (-6, 0);
\coordinate (A139) at (-6, 1);
\coordinate (A140) at (-6, 2);
\coordinate (A141) at (-7, -1);
\coordinate (A142) at (-7, 0);
\coordinate (A143) at (-7, 1);
\coordinate (A144) at (-8, 0);
%\definecolor{color1}{rgb}{0.637767583106935,0.81868784811705,0.486911808319665}
\colorlet{color1}{mycolor4}
\fill[color1!50] (A1) -- ($(A1)! 0.5!(A0)$) -- ($(A9)! 0.5!(A0)$) -- (A9) -- cycle;
\fill[color1!50] ($(A0)! 0.5!(A1)$) -- ($(A88)! 0.5!(A1)$) -- (A1) -- cycle;
\fill[color1!50] ($(A0)! 0.5!(A9)$) -- ($(A45)! 0.5!(A9)$) -- (A9) -- cycle;
\fill[color1!50] (A1) -- ($(A1)! 0.5!(A2)$) -- ($(A10)! 0.5!(A2)$) -- (A10) -- cycle;
\fill[color1!50] ($(A2)! 0.5!(A1)$) -- ($(A89)! 0.5!(A1)$) -- (A1) -- cycle;
\fill[color1!50] (A1) -- (A9) -- (A10) -- cycle;
\fill[color1!50] ($(A88)! 0.5!(A1)$) -- ($(A89)! 0.5!(A1)$) -- (A1) -- cycle;
\fill[color1!50] (A3) -- ($(A3)! 0.5!(A2)$) -- ($(A11)! 0.5!(A2)$) -- (A11) -- cycle;
\fill[color1!50] ($(A2)! 0.5!(A3)$) -- ($(A90)! 0.5!(A3)$) -- (A3) -- cycle;
\fill[color1!50] (A10) -- ($(A10)! 0.5!(A2)$) -- ($(A11)! 0.5!(A2)$) -- (A11) -- cycle;
\fill[color1!50] (A3) -- ($(A3)! 0.5!(A4)$) -- ($(A12)! 0.5!(A4)$) -- (A12) -- cycle;
\fill[color1!50] ($(A4)! 0.5!(A3)$) -- ($(A91)! 0.5!(A3)$) -- (A3) -- cycle;
\fill[color1!50] (A3) -- (A11) -- (A12) -- cycle;
\fill[color1!50] ($(A90)! 0.5!(A3)$) -- ($(A91)! 0.5!(A3)$) -- (A3) -- cycle;
\fill[color1!50] (A5) -- ($(A5)! 0.5!(A4)$) -- ($(A13)! 0.5!(A4)$) -- (A13) -- cycle;
\fill[color1!50] ($(A4)! 0.5!(A5)$) -- ($(A92)! 0.5!(A5)$) -- (A5) -- cycle;
\fill[color1!50] (A12) -- ($(A12)! 0.5!(A4)$) -- ($(A13)! 0.5!(A4)$) -- (A13) -- cycle;
\fill[color1!50] (A5) -- ($(A5)! 0.5!(A6)$) -- ($(A14)! 0.5!(A6)$) -- (A14) -- cycle;
\fill[color1!50] ($(A6)! 0.5!(A5)$) -- ($(A93)! 0.5!(A5)$) -- (A5) -- cycle;
\fill[color1!50] (A5) -- (A13) -- (A14) -- cycle;
\fill[color1!50] ($(A92)! 0.5!(A5)$) -- ($(A93)! 0.5!(A5)$) -- (A5) -- cycle;
\fill[color1!50] (A7) -- ($(A7)! 0.5!(A6)$) -- ($(A15)! 0.5!(A6)$) -- (A15) -- cycle;
\fill[color1!50] ($(A6)! 0.5!(A7)$) -- ($(A94)! 0.5!(A7)$) -- (A7) -- cycle;
\fill[color1!50] (A14) -- ($(A14)! 0.5!(A6)$) -- ($(A15)! 0.5!(A6)$) -- (A15) -- cycle;
\fill[color1!50] (A7) -- ($(A7)! 0.5!(A8)$) -- ($(A16)! 0.5!(A8)$) -- (A16) -- cycle;
\fill[color1!50] ($(A8)! 0.5!(A7)$) -- ($(A95)! 0.5!(A7)$) -- (A7) -- cycle;
\fill[color1!50] (A7) -- (A15) -- (A16) -- cycle;
\fill[color1!50] ($(A94)! 0.5!(A7)$) -- ($(A95)! 0.5!(A7)$) -- (A7) -- cycle;
\fill[color1!50] (A9) -- (A10) -- (A20) -- cycle;
\fill[color1!50] (A9) -- (A17) -- (A18) -- cycle;
\fill[color1!50] (A9) -- ($(A9)! 0.5!(A60)$) -- ($(A17)! 0.5!(A60)$) -- (A17) -- cycle;
\fill[color1!50] (A9) -- (A18) -- (A19) -- cycle;
\fill[color1!50] (A9) -- (A19) -- (A20) -- cycle;
\fill[color1!50] ($(A45)! 0.5!(A9)$) -- ($(A53)! 0.5!(A9)$) -- (A9) -- cycle;
\fill[color1!50] ($(A53)! 0.5!(A9)$) -- ($(A62)! 0.5!(A9)$) -- (A9) -- cycle;
\fill[color1!50] (A9) -- ($(A9)! 0.5!(A60)$) -- ($(A61)! 0.5!(A60)$) -- (A61) -- cycle;
\fill[color1!50] (A9) -- ($(A9)! 0.5!(A62)$) -- ($(A61)! 0.5!(A62)$) -- (A61) -- cycle;
\fill[color1!50] (A10) -- (A11) -- (A20) -- cycle;
\fill[color1!50] (A11) -- (A12) -- (A20) -- cycle;
\fill[color1!50] (A12) -- (A13) -- (A20) -- cycle;
\fill[color1!50] (A13) -- ($(A13)! 0.5!(A21)$) -- ($(A14)! 0.5!(A21)$) -- (A14) -- cycle;
\fill[color1!50] (A13) -- ($(A13)! 0.5!(A21)$) -- ($(A20)! 0.5!(A21)$) -- (A20) -- cycle;
\fill[color1!50] (A14) -- (A15) -- (A22) -- cycle;
\fill[color1!50] (A14) -- ($(A14)! 0.5!(A21)$) -- ($(A22)! 0.5!(A21)$) -- (A22) -- cycle;
\fill[color1!50] (A15) -- ($(A15)! 0.5!(A23)$) -- ($(A16)! 0.5!(A23)$) -- (A16) -- cycle;
\fill[color1!50] (A15) -- ($(A15)! 0.5!(A23)$) -- ($(A22)! 0.5!(A23)$) -- (A22) -- cycle;
\fill[color1!50] (A17) -- (A18) -- (A24) -- cycle;
\fill[color1!50] (A17) -- ($(A17)! 0.5!(A60)$) -- ($(A24)! 0.5!(A60)$) -- (A24) -- cycle;
\fill[color1!50] (A18) -- (A19) -- (A26) -- cycle;
\fill[color1!50] (A18) -- ($(A18)! 0.5!(A25)$) -- ($(A24)! 0.5!(A25)$) -- (A24) -- cycle;
\fill[color1!50] (A18) -- ($(A18)! 0.5!(A25)$) -- ($(A26)! 0.5!(A25)$) -- (A26) -- cycle;
\fill[color1!50] (A19) -- (A20) -- (A26) -- cycle;
\fill[color1!50] (A20) -- ($(A20)! 0.5!(A21)$) -- ($(A27)! 0.5!(A21)$) -- (A27) -- cycle;
\fill[color1!50] (A20) -- ($(A20)! 0.5!(A36)$) -- ($(A26)! 0.5!(A36)$) -- (A26) -- cycle;
\fill[color1!50] (A20) -- ($(A20)! 0.5!(A32)$) -- ($(A27)! 0.5!(A32)$) -- (A27) -- cycle;
\fill[color1!50] (A20) -- ($(A20)! 0.5!(A32)$) -- ($(A42)! 0.5!(A32)$) -- (A42) -- cycle;
\fill[color1!50] (A20) -- ($(A20)! 0.5!(A36)$) -- ($(A42)! 0.5!(A36)$) -- (A42) -- cycle;
\fill[color1!50] (A22) -- ($(A22)! 0.5!(A21)$) -- ($(A28)! 0.5!(A21)$) -- (A28) -- cycle;
\fill[color1!50] (A27) -- ($(A27)! 0.5!(A21)$) -- ($(A28)! 0.5!(A21)$) -- (A28) -- cycle;
\fill[color1!50] (A22) -- ($(A22)! 0.5!(A23)$) -- ($(A29)! 0.5!(A23)$) -- (A29) -- cycle;
\fill[color1!50] (A22) -- (A28) -- (A29) -- cycle;
\fill[color1!50] (A24) -- ($(A24)! 0.5!(A25)$) -- ($(A30)! 0.5!(A25)$) -- (A30) -- cycle;
\fill[color1!50] (A24) -- (A30) -- (A66) -- cycle;
\fill[color1!50] (A24) -- ($(A24)! 0.5!(A60)$) -- ($(A66)! 0.5!(A60)$) -- (A66) -- cycle;
\fill[color1!50] (A26) -- ($(A26)! 0.5!(A25)$) -- ($(A31)! 0.5!(A25)$) -- (A31) -- cycle;
\fill[color1!50] (A30) -- ($(A30)! 0.5!(A25)$) -- ($(A31)! 0.5!(A25)$) -- (A31) -- cycle;
\fill[color1!50] (A26) -- ($(A26)! 0.5!(A36)$) -- ($(A31)! 0.5!(A36)$) -- (A31) -- cycle;
\fill[color1!50] (A27) -- (A28) -- (A33) -- cycle;
\fill[color1!50] (A27) -- ($(A27)! 0.5!(A32)$) -- ($(A33)! 0.5!(A32)$) -- (A33) -- cycle;
\fill[color1!50] (A28) -- ($(A28)! 0.5!(A34)$) -- ($(A29)! 0.5!(A34)$) -- (A29) -- cycle;
\fill[color1!50] (A28) -- ($(A28)! 0.5!(A34)$) -- ($(A33)! 0.5!(A34)$) -- (A33) -- cycle;
\fill[color1!50] (A30) -- (A31) -- (A35) -- cycle;
\fill[color1!50] (A30) -- ($(A30)! 0.5!(A71)$) -- ($(A35)! 0.5!(A71)$) -- (A35) -- cycle;
\fill[color1!50] (A30) -- ($(A30)! 0.5!(A71)$) -- ($(A66)! 0.5!(A71)$) -- (A66) -- cycle;
\fill[color1!50] (A31) -- ($(A31)! 0.5!(A36)$) -- ($(A35)! 0.5!(A36)$) -- (A35) -- cycle;
\fill[color1!50] (A33) -- ($(A33)! 0.5!(A32)$) -- ($(A37)! 0.5!(A32)$) -- (A37) -- cycle;
\fill[color1!50] (A37) -- ($(A37)! 0.5!(A32)$) -- ($(A40)! 0.5!(A32)$) -- (A40) -- cycle;
\fill[color1!50] (A40) -- ($(A40)! 0.5!(A32)$) -- ($(A42)! 0.5!(A32)$) -- (A42) -- cycle;
\fill[color1!50] (A33) -- ($(A33)! 0.5!(A34)$) -- ($(A38)! 0.5!(A34)$) -- (A38) -- cycle;
\fill[color1!50] (A33) -- (A37) -- (A38) -- cycle;
\fill[color1!50] (A35) -- ($(A35)! 0.5!(A36)$) -- ($(A39)! 0.5!(A36)$) -- (A39) -- cycle;
\fill[color1!50] (A35) -- (A39) -- (A75) -- cycle;
\fill[color1!50] (A35) -- ($(A35)! 0.5!(A71)$) -- ($(A75)! 0.5!(A71)$) -- (A75) -- cycle;
\fill[color1!50] (A39) -- ($(A39)! 0.5!(A36)$) -- ($(A42)! 0.5!(A36)$) -- (A42) -- cycle;
\fill[color1!50] (A37) -- ($(A37)! 0.5!(A41)$) -- ($(A38)! 0.5!(A41)$) -- (A38) -- cycle;
\fill[color1!50] (A37) -- ($(A37)! 0.5!(A41)$) -- ($(A40)! 0.5!(A41)$) -- (A40) -- cycle;
\fill[color1!50] (A39) -- (A42) -- (A75) -- cycle;
\fill[color1!50] (A40) -- ($(A40)! 0.5!(A41)$) -- ($(A43)! 0.5!(A41)$) -- (A43) -- cycle;
\fill[color1!50] (A40) -- (A42) -- (A43) -- cycle;
\fill[color1!50] (A42) -- ($(A42)! 0.5!(A44)$) -- ($(A43)! 0.5!(A44)$) -- (A43) -- cycle;
\fill[color1!50] ($(A44)! 0.5!(A42)$) -- ($(A80)! 0.5!(A42)$) -- (A42) -- cycle;
\fill[color1!50] ($(A62)! 0.5!(A42)$) -- ($(A72)! 0.5!(A42)$) -- (A42) -- cycle;
\fill[color1!50] (A42) -- ($(A42)! 0.5!(A62)$) -- ($(A75)! 0.5!(A62)$) -- (A75) -- cycle;
\fill[color1!50] ($(A72)! 0.5!(A42)$) -- ($(A78)! 0.5!(A42)$) -- (A42) -- cycle;
\fill[color1!50] ($(A78)! 0.5!(A42)$) -- ($(A80)! 0.5!(A42)$) -- (A42) -- cycle;
\fill[color1!50] ($(A51)! 0.5!(A81)$) -- ($(A52)! 0.5!(A81)$) -- (A81) -- cycle;
\fill[color1!50] ($(A51)! 0.5!(A81)$) -- ($(A82)! 0.5!(A81)$) -- (A81) -- cycle;
\fill[color1!50] (A61) -- ($(A61)! 0.5!(A60)$) -- ($(A67)! 0.5!(A60)$) -- (A67) -- cycle;
\fill[color1!50] (A66) -- ($(A66)! 0.5!(A60)$) -- ($(A67)! 0.5!(A60)$) -- (A67) -- cycle;
\fill[color1!50] (A61) -- ($(A61)! 0.5!(A62)$) -- ($(A67)! 0.5!(A62)$) -- (A67) -- cycle;
\fill[color1!50] (A67) -- ($(A67)! 0.5!(A62)$) -- ($(A75)! 0.5!(A62)$) -- (A75) -- cycle;
\fill[color1!50] (A66) -- ($(A66)! 0.5!(A71)$) -- ($(A67)! 0.5!(A71)$) -- (A67) -- cycle;
\fill[color1!50] (A67) -- ($(A67)! 0.5!(A71)$) -- ($(A75)! 0.5!(A71)$) -- (A75) -- cycle;
\fill[color1!50] ($(A82)! 0.5!(A81)$) -- ($(A96)! 0.5!(A81)$) -- (A81) -- cycle;
\fill[color1!50] ($(A96)! 0.5!(A109)$) -- ($(A97)! 0.5!(A109)$) -- (A109) -- cycle;
\fill[color1!50] ($(A97)! 0.5!(A109)$) -- ($(A110)! 0.5!(A109)$) -- (A109) -- cycle;
\fill[color1!50] ($(A110)! 0.5!(A109)$) -- ($(A120)! 0.5!(A109)$) -- (A109) -- cycle;
\fill[color1!50] ($(A120)! 0.5!(A129)$) -- ($(A121)! 0.5!(A129)$) -- (A129) -- cycle;
\fill[color1!50] ($(A121)! 0.5!(A129)$) -- ($(A130)! 0.5!(A129)$) -- (A129) -- cycle;
\fill[color1!50] ($(A130)! 0.5!(A129)$) -- ($(A136)! 0.5!(A129)$) -- (A129) -- cycle;
\fill[color1!50] ($(A136)! 0.5!(A141)$) -- ($(A137)! 0.5!(A141)$) -- (A141) -- cycle;
\fill[color1!50] ($(A137)! 0.5!(A141)$) -- ($(A142)! 0.5!(A141)$) -- (A141) -- cycle;
\fill[color1!50] ($(A142)! 0.5!(A141)$) -- ($(A144)! 0.5!(A141)$) -- (A141) -- cycle;
\colorlet{color2}{mycolor3}
\fill[color2!50] ($(A13)! 0.5!(A21)$) -- ($(A14)! 0.5!(A21)$) -- (A21) -- cycle;
\fill[color2!50] ($(A13)! 0.5!(A21)$) -- ($(A20)! 0.5!(A21)$) -- (A21) -- cycle;
\fill[color2!50] ($(A14)! 0.5!(A21)$) -- ($(A22)! 0.5!(A21)$) -- (A21) -- cycle;
\fill[color2!50] ($(A20)! 0.5!(A21)$) -- ($(A27)! 0.5!(A21)$) -- (A21) -- cycle;
\fill[color2!50] ($(A22)! 0.5!(A21)$) -- ($(A28)! 0.5!(A21)$) -- (A21) -- cycle;
\fill[color2!50] ($(A27)! 0.5!(A21)$) -- ($(A28)! 0.5!(A21)$) -- (A21) -- cycle;
\colorlet{color3}{mycolor5}
\fill[color3!50] ($(A18)! 0.5!(A25)$) -- ($(A24)! 0.5!(A25)$) -- (A25) -- cycle;
\fill[color3!50] ($(A18)! 0.5!(A25)$) -- ($(A26)! 0.5!(A25)$) -- (A25) -- cycle;
\fill[color3!50] ($(A24)! 0.5!(A25)$) -- ($(A30)! 0.5!(A25)$) -- (A25) -- cycle;
\fill[color3!50] ($(A26)! 0.5!(A25)$) -- ($(A31)! 0.5!(A25)$) -- (A25) -- cycle;
\fill[color3!50] ($(A30)! 0.5!(A25)$) -- ($(A31)! 0.5!(A25)$) -- (A25) -- cycle;
\colorlet{color4}{mycolor3}
\fill[color4!50] ($(A20)! 0.5!(A32)$) -- ($(A27)! 0.5!(A32)$) -- (A32) -- cycle;
\fill[color4!50] ($(A20)! 0.5!(A32)$) -- ($(A42)! 0.5!(A32)$) -- (A32) -- cycle;
\fill[color4!50] ($(A27)! 0.5!(A32)$) -- ($(A33)! 0.5!(A32)$) -- (A32) -- cycle;
\fill[color4!50] ($(A33)! 0.5!(A32)$) -- ($(A37)! 0.5!(A32)$) -- (A32) -- cycle;
\fill[color4!50] ($(A37)! 0.5!(A32)$) -- ($(A40)! 0.5!(A32)$) -- (A32) -- cycle;
\fill[color4!50] ($(A40)! 0.5!(A32)$) -- ($(A42)! 0.5!(A32)$) -- (A32) -- cycle;
\colorlet{color5}{mycolor5}
\fill[color5!50] ($(A20)! 0.5!(A36)$) -- ($(A26)! 0.5!(A36)$) -- (A36) -- cycle;
\fill[color5!50] ($(A20)! 0.5!(A36)$) -- ($(A42)! 0.5!(A36)$) -- (A36) -- cycle;
\fill[color5!50] ($(A26)! 0.5!(A36)$) -- ($(A31)! 0.5!(A36)$) -- (A36) -- cycle;
\fill[color5!50] ($(A31)! 0.5!(A36)$) -- ($(A35)! 0.5!(A36)$) -- (A36) -- cycle;
\fill[color5!50] ($(A35)! 0.5!(A36)$) -- ($(A39)! 0.5!(A36)$) -- (A36) -- cycle;
\fill[color5!50] ($(A39)! 0.5!(A36)$) -- ($(A42)! 0.5!(A36)$) -- (A36) -- cycle;
\colorlet{color6}{mycolor3}
\fill[color6!50] ($(A46)! 0.5!(A54)$) -- ($(A47)! 0.5!(A54)$) -- (A54) -- cycle;
\fill[color6!50] ($(A46)! 0.5!(A54)$) -- ($(A53)! 0.5!(A54)$) -- (A54) -- cycle;
\fill[color6!50] ($(A47)! 0.5!(A54)$) -- ($(A55)! 0.5!(A54)$) -- (A54) -- cycle;
\fill[color6!50] ($(A53)! 0.5!(A54)$) -- ($(A62)! 0.5!(A54)$) -- (A54) -- cycle;
\fill[color6!50] ($(A55)! 0.5!(A54)$) -- ($(A62)! 0.5!(A54)$) -- (A54) -- cycle;
\colorlet{color7}{mycolor3}
\fill[color7!50] ($(A48)! 0.5!(A56)$) -- ($(A49)! 0.5!(A56)$) -- (A56) -- cycle;
\fill[color7!50] ($(A48)! 0.5!(A56)$) -- ($(A55)! 0.5!(A56)$) -- (A56) -- cycle;
\fill[color7!50] ($(A49)! 0.5!(A56)$) -- ($(A57)! 0.5!(A56)$) -- (A56) -- cycle;
\fill[color7!50] ($(A55)! 0.5!(A56)$) -- ($(A62)! 0.5!(A56)$) -- (A56) -- cycle;
\fill[color7!50] ($(A57)! 0.5!(A56)$) -- ($(A63)! 0.5!(A56)$) -- (A56) -- cycle;
\fill[color7!50] ($(A62)! 0.5!(A56)$) -- ($(A63)! 0.5!(A56)$) -- (A56) -- cycle;
\colorlet{color8}{mycolor3}
\fill[color8!50] ($(A50)! 0.5!(A58)$) -- ($(A51)! 0.5!(A58)$) -- (A58) -- cycle;
\fill[color8!50] ($(A50)! 0.5!(A58)$) -- ($(A57)! 0.5!(A58)$) -- (A58) -- cycle;
\fill[color8!50] ($(A51)! 0.5!(A58)$) -- ($(A59)! 0.5!(A58)$) -- (A58) -- cycle;
\fill[color8!50] ($(A57)! 0.5!(A58)$) -- ($(A64)! 0.5!(A58)$) -- (A58) -- cycle;
\fill[color8!50] ($(A59)! 0.5!(A58)$) -- ($(A65)! 0.5!(A58)$) -- (A58) -- cycle;
\fill[color8!50] ($(A64)! 0.5!(A58)$) -- ($(A65)! 0.5!(A58)$) -- (A58) -- cycle;
\colorlet{color9}{mycolor5}
\fill[color9!50] ($(A9)! 0.5!(A60)$) -- ($(A17)! 0.5!(A60)$) -- (A60) -- cycle;
\fill[color9!50] ($(A9)! 0.5!(A60)$) -- ($(A61)! 0.5!(A60)$) -- (A60) -- cycle;
\fill[color9!50] ($(A17)! 0.5!(A60)$) -- ($(A24)! 0.5!(A60)$) -- (A60) -- cycle;
\fill[color9!50] ($(A24)! 0.5!(A60)$) -- ($(A66)! 0.5!(A60)$) -- (A60) -- cycle;
\fill[color9!50] ($(A61)! 0.5!(A60)$) -- ($(A67)! 0.5!(A60)$) -- (A60) -- cycle;
\fill[color9!50] ($(A66)! 0.5!(A60)$) -- ($(A67)! 0.5!(A60)$) -- (A60) -- cycle;
\colorlet{color10}{mycolor3}
\fill[color10!50] ($(A63)! 0.5!(A69)$) -- ($(A64)! 0.5!(A69)$) -- (A69) -- cycle;
\fill[color10!50] ($(A63)! 0.5!(A69)$) -- ($(A68)! 0.5!(A69)$) -- (A69) -- cycle;
\fill[color10!50] ($(A64)! 0.5!(A69)$) -- ($(A70)! 0.5!(A69)$) -- (A69) -- cycle;
\fill[color10!50] ($(A68)! 0.5!(A69)$) -- ($(A73)! 0.5!(A69)$) -- (A69) -- cycle;
\fill[color10!50] ($(A70)! 0.5!(A69)$) -- ($(A74)! 0.5!(A69)$) -- (A69) -- cycle;
\fill[color10!50] ($(A73)! 0.5!(A69)$) -- ($(A74)! 0.5!(A69)$) -- (A69) -- cycle;
\colorlet{color11}{mycolor5}
\fill[color11!50] ($(A30)! 0.5!(A71)$) -- ($(A35)! 0.5!(A71)$) -- (A71) -- cycle;
\fill[color11!50] ($(A30)! 0.5!(A71)$) -- ($(A66)! 0.5!(A71)$) -- (A71) -- cycle;
\fill[color11!50] ($(A35)! 0.5!(A71)$) -- ($(A75)! 0.5!(A71)$) -- (A71) -- cycle;
\fill[color11!50] ($(A66)! 0.5!(A71)$) -- ($(A67)! 0.5!(A71)$) -- (A71) -- cycle;
\fill[color11!50] ($(A67)! 0.5!(A71)$) -- ($(A75)! 0.5!(A71)$) -- (A71) -- cycle;
\colorlet{color12}{mycolor3}
\fill[color12!50] ($(A72)! 0.5!(A76)$) -- ($(A73)! 0.5!(A76)$) -- (A76) -- cycle;
\fill[color12!50] ($(A72)! 0.5!(A76)$) -- ($(A78)! 0.5!(A76)$) -- (A76) -- cycle;
\fill[color12!50] ($(A73)! 0.5!(A76)$) -- ($(A77)! 0.5!(A76)$) -- (A76) -- cycle;
\fill[color12!50] ($(A77)! 0.5!(A76)$) -- ($(A79)! 0.5!(A76)$) -- (A76) -- cycle;
\fill[color12!50] ($(A78)! 0.5!(A76)$) -- ($(A79)! 0.5!(A76)$) -- (A76) -- cycle;
\colorlet{color13}{mycolor3}
\fill[color13!50] ($(A49)! 0.5!(A83)$) -- ($(A50)! 0.5!(A83)$) -- (A83) -- cycle;
\fill[color13!50] ($(A49)! 0.5!(A83)$) -- ($(A84)! 0.5!(A83)$) -- (A83) -- cycle;
\fill[color13!50] ($(A50)! 0.5!(A83)$) -- ($(A82)! 0.5!(A83)$) -- (A83) -- cycle;
\fill[color13!50] ($(A82)! 0.5!(A83)$) -- ($(A97)! 0.5!(A83)$) -- (A83) -- cycle;
\fill[color13!50] ($(A84)! 0.5!(A83)$) -- ($(A98)! 0.5!(A83)$) -- (A83) -- cycle;
\fill[color13!50] ($(A97)! 0.5!(A83)$) -- ($(A98)! 0.5!(A83)$) -- (A83) -- cycle;
\colorlet{color14}{mycolor3}
\fill[color14!50] ($(A47)! 0.5!(A85)$) -- ($(A48)! 0.5!(A85)$) -- (A85) -- cycle;
\fill[color14!50] ($(A47)! 0.5!(A85)$) -- ($(A86)! 0.5!(A85)$) -- (A85) -- cycle;
\fill[color14!50] ($(A48)! 0.5!(A85)$) -- ($(A84)! 0.5!(A85)$) -- (A85) -- cycle;
\fill[color14!50] ($(A84)! 0.5!(A85)$) -- ($(A99)! 0.5!(A85)$) -- (A85) -- cycle;
\fill[color14!50] ($(A86)! 0.5!(A85)$) -- ($(A99)! 0.5!(A85)$) -- (A85) -- cycle;
\colorlet{color15}{mycolor3}
\fill[color15!50] ($(A45)! 0.5!(A87)$) -- ($(A46)! 0.5!(A87)$) -- (A87) -- cycle;
\fill[color15!50] ($(A45)! 0.5!(A87)$) -- ($(A88)! 0.5!(A87)$) -- (A87) -- cycle;
\fill[color15!50] ($(A46)! 0.5!(A87)$) -- ($(A86)! 0.5!(A87)$) -- (A87) -- cycle;
\fill[color15!50] ($(A86)! 0.5!(A87)$) -- ($(A99)! 0.5!(A87)$) -- (A87) -- cycle;
\fill[color15!50] ($(A88)! 0.5!(A87)$) -- ($(A99)! 0.5!(A87)$) -- (A87) -- cycle;
\colorlet{color16}{mycolor5}
\fill[color16!50] ($(A88)! 0.5!(A100)$) -- ($(A99)! 0.5!(A100)$) -- (A100) -- cycle;
\fill[color16!50] ($(A88)! 0.5!(A100)$) -- ($(A101)! 0.5!(A100)$) -- (A100) -- cycle;
\fill[color16!50] ($(A99)! 0.5!(A100)$) -- ($(A112)! 0.5!(A100)$) -- (A100) -- cycle;
\fill[color16!50] ($(A101)! 0.5!(A100)$) -- ($(A112)! 0.5!(A100)$) -- (A100) -- cycle;
\colorlet{color17}{mycolor3}
\fill[color17!50] ($(A88)! 0.5!(A105)$) -- ($(A89)! 0.5!(A105)$) -- (A105) -- cycle;
\fill[color17!50] ($(A88)! 0.5!(A102)$) -- ($(A101)! 0.5!(A102)$) -- (A102) -- cycle;
\fill[color17!50] (A102) -- ($(A102)! 0.5!(A88)$) -- ($(A103)! 0.5!(A88)$) -- (A103) -- cycle;
\fill[color17!50] (A103) -- ($(A103)! 0.5!(A88)$) -- ($(A104)! 0.5!(A88)$) -- (A104) -- cycle;
\fill[color17!50] (A104) -- ($(A104)! 0.5!(A88)$) -- ($(A105)! 0.5!(A88)$) -- (A105) -- cycle;
\fill[color17!50] ($(A89)! 0.5!(A105)$) -- ($(A90)! 0.5!(A105)$) -- (A105) -- cycle;
\fill[color17!50] ($(A90)! 0.5!(A105)$) -- ($(A91)! 0.5!(A105)$) -- (A105) -- cycle;
\fill[color17!50] ($(A91)! 0.5!(A105)$) -- ($(A92)! 0.5!(A105)$) -- (A105) -- cycle;
\fill[color17!50] ($(A92)! 0.5!(A105)$) -- ($(A106)! 0.5!(A105)$) -- (A105) -- cycle;
\fill[color17!50] ($(A101)! 0.5!(A102)$) -- ($(A114)! 0.5!(A102)$) -- (A102) -- cycle;
\fill[color17!50] (A102) -- ($(A102)! 0.5!(A114)$) -- ($(A103)! 0.5!(A114)$) -- (A103) -- cycle;
\fill[color17!50] (A103) -- ($(A103)! 0.5!(A116)$) -- ($(A104)! 0.5!(A116)$) -- (A104) -- cycle;
\fill[color17!50] (A103) -- ($(A103)! 0.5!(A114)$) -- ($(A115)! 0.5!(A114)$) -- (A115) -- cycle;
\fill[color17!50] (A103) -- ($(A103)! 0.5!(A116)$) -- ($(A115)! 0.5!(A116)$) -- (A115) -- cycle;
\fill[color17!50] (A104) -- ($(A104)! 0.5!(A116)$) -- ($(A105)! 0.5!(A116)$) -- (A105) -- cycle;
\fill[color17!50] ($(A106)! 0.5!(A105)$) -- ($(A117)! 0.5!(A105)$) -- (A105) -- cycle;
\fill[color17!50] (A105) -- ($(A105)! 0.5!(A116)$) -- ($(A133)! 0.5!(A116)$) -- (A133) -- cycle;
\fill[color17!50] ($(A117)! 0.5!(A105)$) -- ($(A126)! 0.5!(A105)$) -- (A105) -- cycle;
\fill[color17!50] ($(A126)! 0.5!(A105)$) -- ($(A142)! 0.5!(A105)$) -- (A105) -- cycle;
\fill[color17!50] (A105) -- ($(A105)! 0.5!(A142)$) -- ($(A133)! 0.5!(A142)$) -- (A133) -- cycle;
\fill[color17!50] ($(A113)! 0.5!(A124)$) -- ($(A114)! 0.5!(A124)$) -- (A124) -- cycle;
\fill[color17!50] ($(A113)! 0.5!(A124)$) -- ($(A123)! 0.5!(A124)$) -- (A124) -- cycle;
\fill[color17!50] (A115) -- ($(A115)! 0.5!(A114)$) -- ($(A124)! 0.5!(A114)$) -- (A124) -- cycle;
\fill[color17!50] (A115) -- ($(A115)! 0.5!(A116)$) -- ($(A125)! 0.5!(A116)$) -- (A125) -- cycle;
\fill[color17!50] (A115) -- (A124) -- (A125) -- cycle;
\fill[color17!50] (A125) -- ($(A125)! 0.5!(A116)$) -- ($(A133)! 0.5!(A116)$) -- (A133) -- cycle;
\fill[color17!50] ($(A123)! 0.5!(A124)$) -- ($(A132)! 0.5!(A124)$) -- (A124) -- cycle;
\fill[color17!50] (A124) -- ($(A124)! 0.5!(A132)$) -- ($(A125)! 0.5!(A132)$) -- (A125) -- cycle;
\fill[color17!50] (A125) -- ($(A125)! 0.5!(A132)$) -- ($(A133)! 0.5!(A132)$) -- (A133) -- cycle;
\fill[color17!50] ($(A131)! 0.5!(A138)$) -- ($(A132)! 0.5!(A138)$) -- (A138) -- cycle;
\fill[color17!50] ($(A131)! 0.5!(A138)$) -- ($(A142)! 0.5!(A138)$) -- (A138) -- cycle;
\fill[color17!50] (A133) -- ($(A133)! 0.5!(A132)$) -- ($(A138)! 0.5!(A132)$) -- (A138) -- cycle;
\fill[color17!50] (A133) -- ($(A133)! 0.5!(A142)$) -- ($(A138)! 0.5!(A142)$) -- (A138) -- cycle;
\colorlet{color18}{mycolor3}
\fill[color18!50] ($(A93)! 0.5!(A107)$) -- ($(A94)! 0.5!(A107)$) -- (A107) -- cycle;
\fill[color18!50] ($(A93)! 0.5!(A107)$) -- ($(A106)! 0.5!(A107)$) -- (A107) -- cycle;
\fill[color18!50] ($(A94)! 0.5!(A107)$) -- ($(A108)! 0.5!(A107)$) -- (A107) -- cycle;
\fill[color18!50] ($(A106)! 0.5!(A107)$) -- ($(A118)! 0.5!(A107)$) -- (A107) -- cycle;
\fill[color18!50] ($(A108)! 0.5!(A107)$) -- ($(A119)! 0.5!(A107)$) -- (A107) -- cycle;
\fill[color18!50] ($(A118)! 0.5!(A107)$) -- ($(A119)! 0.5!(A107)$) -- (A107) -- cycle;
\colorlet{color19}{mycolor3}
\fill[color19!50] ($(A98)! 0.5!(A111)$) -- ($(A99)! 0.5!(A111)$) -- (A111) -- cycle;
\fill[color19!50] ($(A98)! 0.5!(A111)$) -- ($(A110)! 0.5!(A111)$) -- (A111) -- cycle;
\fill[color19!50] ($(A99)! 0.5!(A111)$) -- ($(A122)! 0.5!(A111)$) -- (A111) -- cycle;
\fill[color19!50] ($(A110)! 0.5!(A111)$) -- ($(A121)! 0.5!(A111)$) -- (A111) -- cycle;
\fill[color19!50] ($(A121)! 0.5!(A111)$) -- ($(A122)! 0.5!(A111)$) -- (A111) -- cycle;
\colorlet{color20}{mycolor5}
\fill[color20!50] ($(A103)! 0.5!(A116)$) -- ($(A104)! 0.5!(A116)$) -- (A116) -- cycle;
\fill[color20!50] ($(A103)! 0.5!(A116)$) -- ($(A115)! 0.5!(A116)$) -- (A116) -- cycle;
\fill[color20!50] ($(A104)! 0.5!(A116)$) -- ($(A105)! 0.5!(A116)$) -- (A116) -- cycle;
\fill[color20!50] ($(A105)! 0.5!(A116)$) -- ($(A133)! 0.5!(A116)$) -- (A116) -- cycle;
\fill[color20!50] ($(A115)! 0.5!(A116)$) -- ($(A125)! 0.5!(A116)$) -- (A116) -- cycle;
\fill[color20!50] ($(A125)! 0.5!(A116)$) -- ($(A133)! 0.5!(A116)$) -- (A116) -- cycle;
\colorlet{color21}{mycolor3}
\fill[color21!50] ($(A117)! 0.5!(A127)$) -- ($(A118)! 0.5!(A127)$) -- (A127) -- cycle;
\fill[color21!50] ($(A117)! 0.5!(A127)$) -- ($(A126)! 0.5!(A127)$) -- (A127) -- cycle;
\fill[color21!50] ($(A118)! 0.5!(A127)$) -- ($(A128)! 0.5!(A127)$) -- (A127) -- cycle;
\fill[color21!50] ($(A126)! 0.5!(A127)$) -- ($(A134)! 0.5!(A127)$) -- (A127) -- cycle;
\fill[color21!50] ($(A128)! 0.5!(A127)$) -- ($(A135)! 0.5!(A127)$) -- (A127) -- cycle;
\fill[color21!50] ($(A134)! 0.5!(A127)$) -- ($(A135)! 0.5!(A127)$) -- (A127) -- cycle;
\colorlet{color22}{mycolor3}
\fill[color22!50] ($(A126)! 0.5!(A139)$) -- ($(A134)! 0.5!(A139)$) -- (A139) -- cycle;
\fill[color22!50] ($(A126)! 0.5!(A139)$) -- ($(A142)! 0.5!(A139)$) -- (A139) -- cycle;
\fill[color22!50] ($(A134)! 0.5!(A139)$) -- ($(A140)! 0.5!(A139)$) -- (A139) -- cycle;
\fill[color22!50] ($(A140)! 0.5!(A139)$) -- ($(A143)! 0.5!(A139)$) -- (A139) -- cycle;
\fill[color22!50] ($(A142)! 0.5!(A139)$) -- ($(A143)! 0.5!(A139)$) -- (A139) -- cycle;
\draw[color22, thick] ($(A126)! 0.5!(A139)$) -- ($(A134)! 0.5!(A139)$);
\draw[color22, thick] ($(A126)! 0.5!(A139)$) -- ($(A142)! 0.5!(A139)$);
\draw[color22, thick] ($(A134)! 0.5!(A139)$) -- ($(A140)! 0.5!(A139)$);
\draw[color22, thick] ($(A140)! 0.5!(A139)$) -- ($(A143)! 0.5!(A139)$);
\draw[color22, thick] ($(A142)! 0.5!(A139)$) -- ($(A143)! 0.5!(A139)$);
\draw[color21, thick] ($(A117)! 0.5!(A127)$) -- ($(A118)! 0.5!(A127)$);
\draw[color21, thick] ($(A117)! 0.5!(A127)$) -- ($(A126)! 0.5!(A127)$);
\draw[color21, thick] ($(A118)! 0.5!(A127)$) -- ($(A128)! 0.5!(A127)$);
\draw[color21, thick] ($(A126)! 0.5!(A127)$) -- ($(A134)! 0.5!(A127)$);
\draw[color21, thick] ($(A128)! 0.5!(A127)$) -- ($(A135)! 0.5!(A127)$);
\draw[color21, thick] ($(A134)! 0.5!(A127)$) -- ($(A135)! 0.5!(A127)$);
\draw[color19, thick] ($(A98)! 0.5!(A111)$) -- ($(A99)! 0.5!(A111)$);
\draw[color19, thick] ($(A98)! 0.5!(A111)$) -- ($(A110)! 0.5!(A111)$);
\draw[color19, thick] ($(A99)! 0.5!(A111)$) -- ($(A122)! 0.5!(A111)$);
\draw[color19, thick] ($(A110)! 0.5!(A111)$) -- ($(A121)! 0.5!(A111)$);
\draw[color19, thick] ($(A121)! 0.5!(A111)$) -- ($(A122)! 0.5!(A111)$);
\draw[color18, thick] ($(A93)! 0.5!(A107)$) -- ($(A94)! 0.5!(A107)$);
\draw[color18, thick] ($(A93)! 0.5!(A107)$) -- ($(A106)! 0.5!(A107)$);
\draw[color18, thick] ($(A94)! 0.5!(A107)$) -- ($(A108)! 0.5!(A107)$);
\draw[color18, thick] ($(A106)! 0.5!(A107)$) -- ($(A118)! 0.5!(A107)$);
\draw[color18, thick] ($(A108)! 0.5!(A107)$) -- ($(A119)! 0.5!(A107)$);
\draw[color18, thick] ($(A118)! 0.5!(A107)$) -- ($(A119)! 0.5!(A107)$);
\draw[color17, thick] ($(A88)! 0.5!(A105)$) -- ($(A89)! 0.5!(A105)$);
\draw[color17, thick] ($(A88)! 0.5!(A102)$) -- ($(A101)! 0.5!(A102)$);
\draw[color17, thick] ($(A102)! 0.5!(A88)$) -- ($(A103)! 0.5!(A88)$);
\draw[color17, thick] ($(A103)! 0.5!(A88)$) -- ($(A104)! 0.5!(A88)$);
\draw[color17, thick] ($(A104)! 0.5!(A88)$) -- ($(A105)! 0.5!(A88)$);
\draw[color17, thick] ($(A89)! 0.5!(A105)$) -- ($(A90)! 0.5!(A105)$);
\draw[color17, thick] ($(A90)! 0.5!(A105)$) -- ($(A91)! 0.5!(A105)$);
\draw[color17, thick] ($(A91)! 0.5!(A105)$) -- ($(A92)! 0.5!(A105)$);
\draw[color17, thick] ($(A92)! 0.5!(A105)$) -- ($(A106)! 0.5!(A105)$);
\draw[color17, thick] ($(A101)! 0.5!(A102)$) -- ($(A114)! 0.5!(A102)$);
\draw[color17, thick] ($(A102)! 0.5!(A114)$) -- ($(A103)! 0.5!(A114)$);
\draw[color17, thick] ($(A103)! 0.5!(A116)$) -- ($(A104)! 0.5!(A116)$);
\draw[color17, thick] ($(A103)! 0.5!(A114)$) -- ($(A115)! 0.5!(A114)$);
\draw[color17, thick] ($(A103)! 0.5!(A116)$) -- ($(A115)! 0.5!(A116)$);
\draw[color17, thick] ($(A104)! 0.5!(A116)$) -- ($(A105)! 0.5!(A116)$);
\draw[color17, thick] ($(A106)! 0.5!(A105)$) -- ($(A117)! 0.5!(A105)$);
\draw[color17, thick] ($(A105)! 0.5!(A116)$) -- ($(A133)! 0.5!(A116)$);
\draw[color17, thick] ($(A117)! 0.5!(A105)$) -- ($(A126)! 0.5!(A105)$);
\draw[color17, thick] ($(A126)! 0.5!(A105)$) -- ($(A142)! 0.5!(A105)$);
\draw[color17, thick] ($(A105)! 0.5!(A142)$) -- ($(A133)! 0.5!(A142)$);
\draw[color17, thick] ($(A113)! 0.5!(A124)$) -- ($(A114)! 0.5!(A124)$);
\draw[color17, thick] ($(A113)! 0.5!(A124)$) -- ($(A123)! 0.5!(A124)$);
\draw[color17, thick] ($(A115)! 0.5!(A114)$) -- ($(A124)! 0.5!(A114)$);
\draw[color17, thick] ($(A115)! 0.5!(A116)$) -- ($(A125)! 0.5!(A116)$);
\draw[color17, thick] ($(A125)! 0.5!(A116)$) -- ($(A133)! 0.5!(A116)$);
\draw[color17, thick] ($(A123)! 0.5!(A124)$) -- ($(A132)! 0.5!(A124)$);
\draw[color17, thick] ($(A124)! 0.5!(A132)$) -- ($(A125)! 0.5!(A132)$);
\draw[color17, thick] ($(A125)! 0.5!(A132)$) -- ($(A133)! 0.5!(A132)$);
\draw[color17, thick] ($(A131)! 0.5!(A138)$) -- ($(A132)! 0.5!(A138)$);
\draw[color17, thick] ($(A131)! 0.5!(A138)$) -- ($(A142)! 0.5!(A138)$);
\draw[color17, thick] ($(A133)! 0.5!(A132)$) -- ($(A138)! 0.5!(A132)$);
\draw[color17, thick] ($(A133)! 0.5!(A142)$) -- ($(A138)! 0.5!(A142)$);
\draw[color20, thick] ($(A103)! 0.5!(A116)$) -- ($(A104)! 0.5!(A116)$);
\draw[color20, thick] ($(A103)! 0.5!(A116)$) -- ($(A115)! 0.5!(A116)$);
\draw[color20, thick] ($(A104)! 0.5!(A116)$) -- ($(A105)! 0.5!(A116)$);
\draw[color20, thick] ($(A105)! 0.5!(A116)$) -- ($(A133)! 0.5!(A116)$);
\draw[color20, thick] ($(A115)! 0.5!(A116)$) -- ($(A125)! 0.5!(A116)$);
\draw[color20, thick] ($(A125)! 0.5!(A116)$) -- ($(A133)! 0.5!(A116)$);
\draw[color16, thick] ($(A88)! 0.5!(A100)$) -- ($(A99)! 0.5!(A100)$);
\draw[color16, thick] ($(A88)! 0.5!(A100)$) -- ($(A101)! 0.5!(A100)$);
\draw[color16, thick] ($(A99)! 0.5!(A100)$) -- ($(A112)! 0.5!(A100)$);
\draw[color16, thick] ($(A101)! 0.5!(A100)$) -- ($(A112)! 0.5!(A100)$);
\draw[color15, thick] ($(A45)! 0.5!(A87)$) -- ($(A46)! 0.5!(A87)$);
\draw[color15, thick] ($(A45)! 0.5!(A87)$) -- ($(A88)! 0.5!(A87)$);
\draw[color15, thick] ($(A46)! 0.5!(A87)$) -- ($(A86)! 0.5!(A87)$);
\draw[color15, thick] ($(A86)! 0.5!(A87)$) -- ($(A99)! 0.5!(A87)$);
\draw[color15, thick] ($(A88)! 0.5!(A87)$) -- ($(A99)! 0.5!(A87)$);
\draw[color14, thick] ($(A47)! 0.5!(A85)$) -- ($(A48)! 0.5!(A85)$);
\draw[color14, thick] ($(A47)! 0.5!(A85)$) -- ($(A86)! 0.5!(A85)$);
\draw[color14, thick] ($(A48)! 0.5!(A85)$) -- ($(A84)! 0.5!(A85)$);
\draw[color14, thick] ($(A84)! 0.5!(A85)$) -- ($(A99)! 0.5!(A85)$);
\draw[color14, thick] ($(A86)! 0.5!(A85)$) -- ($(A99)! 0.5!(A85)$);
\draw[color13, thick] ($(A49)! 0.5!(A83)$) -- ($(A50)! 0.5!(A83)$);
\draw[color13, thick] ($(A49)! 0.5!(A83)$) -- ($(A84)! 0.5!(A83)$);
\draw[color13, thick] ($(A50)! 0.5!(A83)$) -- ($(A82)! 0.5!(A83)$);
\draw[color13, thick] ($(A82)! 0.5!(A83)$) -- ($(A97)! 0.5!(A83)$);
\draw[color13, thick] ($(A84)! 0.5!(A83)$) -- ($(A98)! 0.5!(A83)$);
\draw[color13, thick] ($(A97)! 0.5!(A83)$) -- ($(A98)! 0.5!(A83)$);
\draw[color12, thick] ($(A72)! 0.5!(A76)$) -- ($(A73)! 0.5!(A76)$);
\draw[color12, thick] ($(A72)! 0.5!(A76)$) -- ($(A78)! 0.5!(A76)$);
\draw[color12, thick] ($(A73)! 0.5!(A76)$) -- ($(A77)! 0.5!(A76)$);
\draw[color12, thick] ($(A77)! 0.5!(A76)$) -- ($(A79)! 0.5!(A76)$);
\draw[color12, thick] ($(A78)! 0.5!(A76)$) -- ($(A79)! 0.5!(A76)$);
\draw[color10, thick] ($(A63)! 0.5!(A69)$) -- ($(A64)! 0.5!(A69)$);
\draw[color10, thick] ($(A63)! 0.5!(A69)$) -- ($(A68)! 0.5!(A69)$);
\draw[color10, thick] ($(A64)! 0.5!(A69)$) -- ($(A70)! 0.5!(A69)$);
\draw[color10, thick] ($(A68)! 0.5!(A69)$) -- ($(A73)! 0.5!(A69)$);
\draw[color10, thick] ($(A70)! 0.5!(A69)$) -- ($(A74)! 0.5!(A69)$);
\draw[color10, thick] ($(A73)! 0.5!(A69)$) -- ($(A74)! 0.5!(A69)$);
\draw[color8, thick] ($(A50)! 0.5!(A58)$) -- ($(A51)! 0.5!(A58)$);
\draw[color8, thick] ($(A50)! 0.5!(A58)$) -- ($(A57)! 0.5!(A58)$);
\draw[color8, thick] ($(A51)! 0.5!(A58)$) -- ($(A59)! 0.5!(A58)$);
\draw[color8, thick] ($(A57)! 0.5!(A58)$) -- ($(A64)! 0.5!(A58)$);
\draw[color8, thick] ($(A59)! 0.5!(A58)$) -- ($(A65)! 0.5!(A58)$);
\draw[color8, thick] ($(A64)! 0.5!(A58)$) -- ($(A65)! 0.5!(A58)$);
\draw[color7, thick] ($(A48)! 0.5!(A56)$) -- ($(A49)! 0.5!(A56)$);
\draw[color7, thick] ($(A48)! 0.5!(A56)$) -- ($(A55)! 0.5!(A56)$);
\draw[color7, thick] ($(A49)! 0.5!(A56)$) -- ($(A57)! 0.5!(A56)$);
\draw[color7, thick] ($(A55)! 0.5!(A56)$) -- ($(A62)! 0.5!(A56)$);
\draw[color7, thick] ($(A57)! 0.5!(A56)$) -- ($(A63)! 0.5!(A56)$);
\draw[color7, thick] ($(A62)! 0.5!(A56)$) -- ($(A63)! 0.5!(A56)$);
\draw[color6, thick] ($(A46)! 0.5!(A54)$) -- ($(A47)! 0.5!(A54)$);
\draw[color6, thick] ($(A46)! 0.5!(A54)$) -- ($(A53)! 0.5!(A54)$);
\draw[color6, thick] ($(A47)! 0.5!(A54)$) -- ($(A55)! 0.5!(A54)$);
\draw[color6, thick] ($(A53)! 0.5!(A54)$) -- ($(A62)! 0.5!(A54)$);
\draw[color6, thick] ($(A55)! 0.5!(A54)$) -- ($(A62)! 0.5!(A54)$);
\draw[color1, thick] ($(A1)! 0.5!(A0)$) -- ($(A9)! 0.5!(A0)$);
\draw[color1, thick] ($(A0)! 0.5!(A1)$) -- ($(A88)! 0.5!(A1)$);
\draw[color1, thick] ($(A0)! 0.5!(A9)$) -- ($(A45)! 0.5!(A9)$);
\draw[color1, thick] ($(A1)! 0.5!(A2)$) -- ($(A10)! 0.5!(A2)$);
\draw[color1, thick] ($(A2)! 0.5!(A1)$) -- ($(A89)! 0.5!(A1)$);
\draw[color1, thick] ($(A88)! 0.5!(A1)$) -- ($(A89)! 0.5!(A1)$);
\draw[color1, thick] ($(A3)! 0.5!(A2)$) -- ($(A11)! 0.5!(A2)$);
\draw[color1, thick] ($(A2)! 0.5!(A3)$) -- ($(A90)! 0.5!(A3)$);
\draw[color1, thick] ($(A10)! 0.5!(A2)$) -- ($(A11)! 0.5!(A2)$);
\draw[color1, thick] ($(A3)! 0.5!(A4)$) -- ($(A12)! 0.5!(A4)$);
\draw[color1, thick] ($(A4)! 0.5!(A3)$) -- ($(A91)! 0.5!(A3)$);
\draw[color1, thick] ($(A90)! 0.5!(A3)$) -- ($(A91)! 0.5!(A3)$);
\draw[color1, thick] ($(A5)! 0.5!(A4)$) -- ($(A13)! 0.5!(A4)$);
\draw[color1, thick] ($(A4)! 0.5!(A5)$) -- ($(A92)! 0.5!(A5)$);
\draw[color1, thick] ($(A12)! 0.5!(A4)$) -- ($(A13)! 0.5!(A4)$);
\draw[color1, thick] ($(A5)! 0.5!(A6)$) -- ($(A14)! 0.5!(A6)$);
\draw[color1, thick] ($(A6)! 0.5!(A5)$) -- ($(A93)! 0.5!(A5)$);
\draw[color1, thick] ($(A92)! 0.5!(A5)$) -- ($(A93)! 0.5!(A5)$);
\draw[color1, thick] ($(A7)! 0.5!(A6)$) -- ($(A15)! 0.5!(A6)$);
\draw[color1, thick] ($(A6)! 0.5!(A7)$) -- ($(A94)! 0.5!(A7)$);
\draw[color1, thick] ($(A14)! 0.5!(A6)$) -- ($(A15)! 0.5!(A6)$);
\draw[color1, thick] ($(A7)! 0.5!(A8)$) -- ($(A16)! 0.5!(A8)$);
\draw[color1, thick] ($(A8)! 0.5!(A7)$) -- ($(A95)! 0.5!(A7)$);
\draw[color1, thick] ($(A94)! 0.5!(A7)$) -- ($(A95)! 0.5!(A7)$);
\draw[color1, thick] ($(A9)! 0.5!(A60)$) -- ($(A17)! 0.5!(A60)$);
\draw[color1, thick] ($(A45)! 0.5!(A9)$) -- ($(A53)! 0.5!(A9)$);
\draw[color1, thick] ($(A53)! 0.5!(A9)$) -- ($(A62)! 0.5!(A9)$);
\draw[color1, thick] ($(A9)! 0.5!(A60)$) -- ($(A61)! 0.5!(A60)$);
\draw[color1, thick] ($(A9)! 0.5!(A62)$) -- ($(A61)! 0.5!(A62)$);
\draw[color1, thick] ($(A13)! 0.5!(A21)$) -- ($(A14)! 0.5!(A21)$);
\draw[color1, thick] ($(A13)! 0.5!(A21)$) -- ($(A20)! 0.5!(A21)$);
\draw[color1, thick] ($(A14)! 0.5!(A21)$) -- ($(A22)! 0.5!(A21)$);
\draw[color1, thick] ($(A15)! 0.5!(A23)$) -- ($(A16)! 0.5!(A23)$);
\draw[color1, thick] ($(A15)! 0.5!(A23)$) -- ($(A22)! 0.5!(A23)$);
\draw[color1, thick] ($(A17)! 0.5!(A60)$) -- ($(A24)! 0.5!(A60)$);
\draw[color1, thick] ($(A18)! 0.5!(A25)$) -- ($(A24)! 0.5!(A25)$);
\draw[color1, thick] ($(A18)! 0.5!(A25)$) -- ($(A26)! 0.5!(A25)$);
\draw[color1, thick] ($(A20)! 0.5!(A21)$) -- ($(A27)! 0.5!(A21)$);
\draw[color1, thick] ($(A20)! 0.5!(A36)$) -- ($(A26)! 0.5!(A36)$);
\draw[color1, thick] ($(A20)! 0.5!(A32)$) -- ($(A27)! 0.5!(A32)$);
\draw[color1, thick] ($(A20)! 0.5!(A32)$) -- ($(A42)! 0.5!(A32)$);
\draw[color1, thick] ($(A20)! 0.5!(A36)$) -- ($(A42)! 0.5!(A36)$);
\draw[color1, thick] ($(A22)! 0.5!(A21)$) -- ($(A28)! 0.5!(A21)$);
\draw[color1, thick] ($(A27)! 0.5!(A21)$) -- ($(A28)! 0.5!(A21)$);
\draw[color1, thick] ($(A22)! 0.5!(A23)$) -- ($(A29)! 0.5!(A23)$);
\draw[color1, thick] ($(A24)! 0.5!(A25)$) -- ($(A30)! 0.5!(A25)$);
\draw[color1, thick] ($(A24)! 0.5!(A60)$) -- ($(A66)! 0.5!(A60)$);
\draw[color1, thick] ($(A26)! 0.5!(A25)$) -- ($(A31)! 0.5!(A25)$);
\draw[color1, thick] ($(A30)! 0.5!(A25)$) -- ($(A31)! 0.5!(A25)$);
\draw[color1, thick] ($(A26)! 0.5!(A36)$) -- ($(A31)! 0.5!(A36)$);
\draw[color1, thick] ($(A27)! 0.5!(A32)$) -- ($(A33)! 0.5!(A32)$);
\draw[color1, thick] ($(A28)! 0.5!(A34)$) -- ($(A29)! 0.5!(A34)$);
\draw[color1, thick] ($(A28)! 0.5!(A34)$) -- ($(A33)! 0.5!(A34)$);
\draw[color1, thick] ($(A30)! 0.5!(A71)$) -- ($(A35)! 0.5!(A71)$);
\draw[color1, thick] ($(A30)! 0.5!(A71)$) -- ($(A66)! 0.5!(A71)$);
\draw[color1, thick] ($(A31)! 0.5!(A36)$) -- ($(A35)! 0.5!(A36)$);
\draw[color1, thick] ($(A33)! 0.5!(A32)$) -- ($(A37)! 0.5!(A32)$);
\draw[color1, thick] ($(A37)! 0.5!(A32)$) -- ($(A40)! 0.5!(A32)$);
\draw[color1, thick] ($(A40)! 0.5!(A32)$) -- ($(A42)! 0.5!(A32)$);
\draw[color1, thick] ($(A33)! 0.5!(A34)$) -- ($(A38)! 0.5!(A34)$);
\draw[color1, thick] ($(A35)! 0.5!(A36)$) -- ($(A39)! 0.5!(A36)$);
\draw[color1, thick] ($(A35)! 0.5!(A71)$) -- ($(A75)! 0.5!(A71)$);
\draw[color1, thick] ($(A39)! 0.5!(A36)$) -- ($(A42)! 0.5!(A36)$);
\draw[color1, thick] ($(A37)! 0.5!(A41)$) -- ($(A38)! 0.5!(A41)$);
\draw[color1, thick] ($(A37)! 0.5!(A41)$) -- ($(A40)! 0.5!(A41)$);
\draw[color1, thick] ($(A40)! 0.5!(A41)$) -- ($(A43)! 0.5!(A41)$);
\draw[color1, thick] ($(A42)! 0.5!(A44)$) -- ($(A43)! 0.5!(A44)$);
\draw[color1, thick] ($(A44)! 0.5!(A42)$) -- ($(A80)! 0.5!(A42)$);
\draw[color1, thick] ($(A62)! 0.5!(A42)$) -- ($(A72)! 0.5!(A42)$);
\draw[color1, thick] ($(A42)! 0.5!(A62)$) -- ($(A75)! 0.5!(A62)$);
\draw[color1, thick] ($(A72)! 0.5!(A42)$) -- ($(A78)! 0.5!(A42)$);
\draw[color1, thick] ($(A78)! 0.5!(A42)$) -- ($(A80)! 0.5!(A42)$);
\draw[color1, thick] ($(A51)! 0.5!(A81)$) -- ($(A52)! 0.5!(A81)$);
\draw[color1, thick] ($(A51)! 0.5!(A81)$) -- ($(A82)! 0.5!(A81)$);
\draw[color1, thick] ($(A61)! 0.5!(A60)$) -- ($(A67)! 0.5!(A60)$);
\draw[color1, thick] ($(A66)! 0.5!(A60)$) -- ($(A67)! 0.5!(A60)$);
\draw[color1, thick] ($(A61)! 0.5!(A62)$) -- ($(A67)! 0.5!(A62)$);
\draw[color1, thick] ($(A67)! 0.5!(A62)$) -- ($(A75)! 0.5!(A62)$);
\draw[color1, thick] ($(A66)! 0.5!(A71)$) -- ($(A67)! 0.5!(A71)$);
\draw[color1, thick] ($(A67)! 0.5!(A71)$) -- ($(A75)! 0.5!(A71)$);
\draw[color1, thick] ($(A82)! 0.5!(A81)$) -- ($(A96)! 0.5!(A81)$);
\draw[color1, thick] ($(A96)! 0.5!(A109)$) -- ($(A97)! 0.5!(A109)$);
\draw[color1, thick] ($(A97)! 0.5!(A109)$) -- ($(A110)! 0.5!(A109)$);
\draw[color1, thick] ($(A110)! 0.5!(A109)$) -- ($(A120)! 0.5!(A109)$);
\draw[color1, thick] ($(A120)! 0.5!(A129)$) -- ($(A121)! 0.5!(A129)$);
\draw[color1, thick] ($(A121)! 0.5!(A129)$) -- ($(A130)! 0.5!(A129)$);
\draw[color1, thick] ($(A130)! 0.5!(A129)$) -- ($(A136)! 0.5!(A129)$);
\draw[color1, thick] ($(A136)! 0.5!(A141)$) -- ($(A137)! 0.5!(A141)$);
\draw[color1, thick] ($(A137)! 0.5!(A141)$) -- ($(A142)! 0.5!(A141)$);
\draw[color1, thick] ($(A142)! 0.5!(A141)$) -- ($(A144)! 0.5!(A141)$);
\draw[color11, thick] ($(A30)! 0.5!(A71)$) -- ($(A35)! 0.5!(A71)$);
\draw[color11, thick] ($(A30)! 0.5!(A71)$) -- ($(A66)! 0.5!(A71)$);
\draw[color11, thick] ($(A35)! 0.5!(A71)$) -- ($(A75)! 0.5!(A71)$);
\draw[color11, thick] ($(A66)! 0.5!(A71)$) -- ($(A67)! 0.5!(A71)$);
\draw[color11, thick] ($(A67)! 0.5!(A71)$) -- ($(A75)! 0.5!(A71)$);
\draw[color9, thick] ($(A9)! 0.5!(A60)$) -- ($(A17)! 0.5!(A60)$);
\draw[color9, thick] ($(A9)! 0.5!(A60)$) -- ($(A61)! 0.5!(A60)$);
\draw[color9, thick] ($(A17)! 0.5!(A60)$) -- ($(A24)! 0.5!(A60)$);
\draw[color9, thick] ($(A24)! 0.5!(A60)$) -- ($(A66)! 0.5!(A60)$);
\draw[color9, thick] ($(A61)! 0.5!(A60)$) -- ($(A67)! 0.5!(A60)$);
\draw[color9, thick] ($(A66)! 0.5!(A60)$) -- ($(A67)! 0.5!(A60)$);
\draw[color5, thick] ($(A20)! 0.5!(A36)$) -- ($(A26)! 0.5!(A36)$);
\draw[color5, thick] ($(A20)! 0.5!(A36)$) -- ($(A42)! 0.5!(A36)$);
\draw[color5, thick] ($(A26)! 0.5!(A36)$) -- ($(A31)! 0.5!(A36)$);
\draw[color5, thick] ($(A31)! 0.5!(A36)$) -- ($(A35)! 0.5!(A36)$);
\draw[color5, thick] ($(A35)! 0.5!(A36)$) -- ($(A39)! 0.5!(A36)$);
\draw[color5, thick] ($(A39)! 0.5!(A36)$) -- ($(A42)! 0.5!(A36)$);
\draw[color4, thick] ($(A20)! 0.5!(A32)$) -- ($(A27)! 0.5!(A32)$);
\draw[color4, thick] ($(A20)! 0.5!(A32)$) -- ($(A42)! 0.5!(A32)$);
\draw[color4, thick] ($(A27)! 0.5!(A32)$) -- ($(A33)! 0.5!(A32)$);
\draw[color4, thick] ($(A33)! 0.5!(A32)$) -- ($(A37)! 0.5!(A32)$);
\draw[color4, thick] ($(A37)! 0.5!(A32)$) -- ($(A40)! 0.5!(A32)$);
\draw[color4, thick] ($(A40)! 0.5!(A32)$) -- ($(A42)! 0.5!(A32)$);
\draw[color3, thick] ($(A18)! 0.5!(A25)$) -- ($(A24)! 0.5!(A25)$);
\draw[color3, thick] ($(A18)! 0.5!(A25)$) -- ($(A26)! 0.5!(A25)$);
\draw[color3, thick] ($(A24)! 0.5!(A25)$) -- ($(A30)! 0.5!(A25)$);
\draw[color3, thick] ($(A26)! 0.5!(A25)$) -- ($(A31)! 0.5!(A25)$);
\draw[color3, thick] ($(A30)! 0.5!(A25)$) -- ($(A31)! 0.5!(A25)$);
\draw[color2, thick] ($(A13)! 0.5!(A21)$) -- ($(A14)! 0.5!(A21)$);
\draw[color2, thick] ($(A13)! 0.5!(A21)$) -- ($(A20)! 0.5!(A21)$);
\draw[color2, thick] ($(A14)! 0.5!(A21)$) -- ($(A22)! 0.5!(A21)$);
\draw[color2, thick] ($(A20)! 0.5!(A21)$) -- ($(A27)! 0.5!(A21)$);
\draw[color2, thick] ($(A22)! 0.5!(A21)$) -- ($(A28)! 0.5!(A21)$);
\draw[color2, thick] ($(A27)! 0.5!(A21)$) -- ($(A28)! 0.5!(A21)$);
\foreach \a/\b/\c in {0/1/9,0/1/88,0/9/45,0/45/88,1/2/10,1/2/89,1/9/10,1/88/89,2/3/11,2/3/90,2/10/11,2/89/90,3/4/12,3/4/91,3/11/12,3/90/91,4/5/13,4/5/92,4/12/13,4/91/92,5/6/14,5/6/93,5/13/14,5/92/93,6/7/15,6/7/94,6/14/15,6/93/94,7/8/16,7/8/95,7/15/16,7/94/95,9/10/20,9/17/18,9/17/60,9/18/19,9/19/20,9/45/53,9/53/62,9/60/61,9/61/62,10/11/20,11/12/20,12/13/20,13/14/21,13/20/21,14/15/22,14/21/22,15/16/23,15/22/23,17/18/24,17/24/60,18/19/26,18/24/25,18/25/26,19/20/26,20/21/27,20/26/36,20/27/32,20/32/42,20/36/42,21/22/28,21/27/28,22/23/29,22/28/29,24/25/30,24/30/66,24/60/66,25/26/31,25/30/31,26/31/36,27/28/33,27/32/33,28/29/34,28/33/34,30/31/35,30/35/71,30/66/71,31/35/36,32/33/37,32/37/40,32/40/42,33/34/38,33/37/38,35/36/39,35/39/75,35/71/75,36/39/42,37/38/41,37/40/41,39/42/75,40/41/43,40/42/43,42/43/44,42/44/80,42/62/72,42/62/75,42/72/78,42/78/80,45/46/53,45/46/87,45/87/88,46/47/54,46/47/86,46/53/54,46/86/87,47/48/55,47/48/85,47/54/55,47/85/86,48/49/56,48/49/84,48/55/56,48/84/85,49/50/57,49/50/83,49/56/57,49/83/84,50/51/58,50/51/82,50/57/58,50/82/83,51/52/59,51/52/81,51/58/59,51/81/82,53/54/62,54/55/62,55/56/62,56/57/63,56/62/63,57/58/64,57/63/64,58/59/65,58/64/65,60/61/67,60/66/67,61/62/67,62/63/68,62/67/75,62/68/72,63/64/69,63/68/69,64/65/70,64/69/70,66/67/71,67/71/75,68/69/73,68/72/73,69/70/74,69/73/74,72/73/76,72/76/78,73/74/77,73/76/77,76/77/79,76/78/79,78/79/80,81/82/96,82/83/97,82/96/97,83/84/98,83/97/98,84/85/99,84/98/99,85/86/99,86/87/99,87/88/99,88/89/105,88/99/100,88/100/101,88/101/102,88/102/103,88/103/104,88/104/105,89/90/105,90/91/105,91/92/105,92/93/106,92/105/106,93/94/107,93/106/107,94/95/108,94/107/108,96/97/109,97/98/110,97/109/110,98/99/111,98/110/111,99/100/112,99/111/122,99/112/131,99/122/142,99/131/142,100/101/112,101/102/114,101/112/113,101/113/114,102/103/114,103/104/116,103/114/115,103/115/116,104/105/116,105/106/117,105/116/133,105/117/126,105/126/142,105/133/142,106/107/118,106/117/118,107/108/119,107/118/119,109/110/120,110/111/121,110/120/121,111/121/122,112/113/123,112/123/131,113/114/124,113/123/124,114/115/124,115/116/125,115/124/125,116/125/133,117/118/127,117/126/127,118/119/128,118/127/128,120/121/129,121/122/130,121/129/130,122/130/137,122/137/142,123/124/132,123/131/132,124/125/132,125/132/133,126/127/134,126/134/139,126/139/142,127/128/135,127/134/135,129/130/136,130/136/137,131/132/138,131/138/142,132/133/138,133/138/142,134/135/140,134/139/140,136/137/141,137/141/142,139/140/143,139/142/143,141/142/144,142/143/144}{
  \draw[black!30] (A\a) -- (A\b) -- (A\c) -- cycle;
}

\draw[black, thick] (A20) -- (A42) -- (A62) -- (A9) --cycle;
\draw[black, thick] (A105) -- (A142) -- (A99) -- (A88) --cycle;

\foreach \i in {1,3,5,7,9,10,11,12,13,14,15,16,17,18,19,20,22,24,26,27,28,29,30,31,33,35,37,38,39,40,42,43,54,56,58,61,66,67,69,75,76,81,83,85,87,100,102,103,104,105,107,109,111,115,124,125,127,129,133,138,139,141}{
  \fill[myred] (A\i) circle (3pt);
}
\foreach \i in {0,2,4,6,8,21,23,25,32,34,36,41,44,45,46,47,48,49,50,51,52,53,55,57,59,60,62,63,64,65,68,70,71,72,73,74,77,78,79,80,82,84,86,88,89,90,91,92,93,94,95,96,97,98,99,101,106,108,110,112,113,114,116,117,118,119,120,121,122,123,126,128,130,131,132,134,135,136,137,140,142,143,144}{
  \fill[myblue] (A\i) circle (3pt);
}

%\foreach \i in {-8,-6,-4,-2,0,2,4,6,8}{
%	\foreach \j in {-8,-6,-4,-2,0,2,4,6,8}{
%		\node[] at (\i,\j) {\input{tikz/even_even.tikz}};
%	}
%	\foreach \j in {-7,-5,-3,-1,1,3,5,7}{
%		\node[] at (\i,\j) {\input{tikz/even_odd.tikz}};
%	}
%}
%\foreach \i in {-7,-5,-3,-1,1,3,5,7}{
%	\foreach \j in {-8,-6,-4,-2,0,2,4,6,8}{
%		\node[rotate=90] at (\i,\j) {\input{tikz/even_odd.tikz}};
%	}
%	\foreach \j in {-7,-5,-3,-1,1,3,5,7}{
%		\node[] at (\i,\j) {\input{tikz/odd_odd.tikz}};
%	}
%}

%% file: tikz/double_split_8.tikz
\coordinate (A0) at (0, 0);
\coordinate (A1) at (0, 1);
\coordinate (A2) at (0, 2);
\coordinate (A3) at (0, 3);
\coordinate (A4) at (0, 4);
\coordinate (A5) at (0, 5);
\coordinate (A6) at (0, 6);
\coordinate (A7) at (0, 7);
\coordinate (A8) at (0, 8);
\coordinate (A9) at (1, 0);
\coordinate (A10) at (1, 1);
\coordinate (A11) at (1, 2);
\coordinate (A12) at (1, 3);
\coordinate (A13) at (1, 4);
\coordinate (A14) at (1, 5);
\coordinate (A15) at (1, 6);
\coordinate (A16) at (1, 7);
\coordinate (A17) at (2, 0);
\coordinate (A18) at (2, 1);
\coordinate (A19) at (2, 2);
\coordinate (A20) at (2, 3);
\coordinate (A21) at (2, 4);
\coordinate (A22) at (2, 5);
\coordinate (A23) at (2, 6);
\coordinate (A24) at (3, 0);
\coordinate (A25) at (3, 1);
\coordinate (A26) at (3, 2);
\coordinate (A27) at (3, 3);
\coordinate (A28) at (3, 4);
\coordinate (A29) at (3, 5);
\coordinate (A30) at (4, 0);
\coordinate (A31) at (4, 1);
\coordinate (A32) at (4, 2);
\coordinate (A33) at (4, 3);
\coordinate (A34) at (4, 4);
\coordinate (A35) at (5, 0);
\coordinate (A36) at (5, 1);
\coordinate (A37) at (5, 2);
\coordinate (A38) at (5, 3);
\coordinate (A39) at (6, 0);
\coordinate (A40) at (6, 1);
\coordinate (A41) at (6, 2);
\coordinate (A42) at (7, 0);
\coordinate (A43) at (7, 1);
\coordinate (A44) at (8, 0);
\coordinate (A45) at (0, -1);
\coordinate (A46) at (0, -2);
\coordinate (A47) at (0, -3);
\coordinate (A48) at (0, -4);
\coordinate (A49) at (0, -5);
\coordinate (A50) at (0, -6);
\coordinate (A51) at (0, -7);
\coordinate (A52) at (0, -8);
\coordinate (A53) at (1, -1);
\coordinate (A54) at (1, -2);
\coordinate (A55) at (1, -3);
\coordinate (A56) at (1, -4);
\coordinate (A57) at (1, -5);
\coordinate (A58) at (1, -6);
\coordinate (A59) at (1, -7);
\coordinate (A60) at (2, -1);
\coordinate (A61) at (2, -2);
\coordinate (A62) at (2, -3);
\coordinate (A63) at (2, -4);
\coordinate (A64) at (2, -5);
\coordinate (A65) at (2, -6);
\coordinate (A66) at (3, -1);
\coordinate (A67) at (3, -2);
\coordinate (A68) at (3, -3);
\coordinate (A69) at (3, -4);
\coordinate (A70) at (3, -5);
\coordinate (A71) at (4, -1);
\coordinate (A72) at (4, -2);
\coordinate (A73) at (4, -3);
\coordinate (A74) at (4, -4);
\coordinate (A75) at (5, -1);
\coordinate (A76) at (5, -2);
\coordinate (A77) at (5, -3);
\coordinate (A78) at (6, -1);
\coordinate (A79) at (6, -2);
\coordinate (A80) at (7, -1);
\coordinate (A81) at (-1, -7);
\coordinate (A82) at (-1, -6);
\coordinate (A83) at (-1, -5);
\coordinate (A84) at (-1, -4);
\coordinate (A85) at (-1, -3);
\coordinate (A86) at (-1, -2);
\coordinate (A87) at (-1, -1);
\coordinate (A88) at (-1, 0);
\coordinate (A89) at (-1, 1);
\coordinate (A90) at (-1, 2);
\coordinate (A91) at (-1, 3);
\coordinate (A92) at (-1, 4);
\coordinate (A93) at (-1, 5);
\coordinate (A94) at (-1, 6);
\coordinate (A95) at (-1, 7);
\coordinate (A96) at (-2, -6);
\coordinate (A97) at (-2, -5);
\coordinate (A98) at (-2, -4);
\coordinate (A99) at (-2, -3);
\coordinate (A100) at (-2, -2);
\coordinate (A101) at (-2, -1);
\coordinate (A102) at (-2, 0);
\coordinate (A103) at (-2, 1);
\coordinate (A104) at (-2, 2);
\coordinate (A105) at (-2, 3);
\coordinate (A106) at (-2, 4);
\coordinate (A107) at (-2, 5);
\coordinate (A108) at (-2, 6);
\coordinate (A109) at (-3, -5);
\coordinate (A110) at (-3, -4);
\coordinate (A111) at (-3, -3);
\coordinate (A112) at (-3, -2);
\coordinate (A113) at (-3, -1);
\coordinate (A114) at (-3, 0);
\coordinate (A115) at (-3, 1);
\coordinate (A116) at (-3, 2);
\coordinate (A117) at (-3, 3);
\coordinate (A118) at (-3, 4);
\coordinate (A119) at (-3, 5);
\coordinate (A120) at (-4, -4);
\coordinate (A121) at (-4, -3);
\coordinate (A122) at (-4, -2);
\coordinate (A123) at (-4, -1);
\coordinate (A124) at (-4, 0);
\coordinate (A125) at (-4, 1);
\coordinate (A126) at (-4, 2);
\coordinate (A127) at (-4, 3);
\coordinate (A128) at (-4, 4);
\coordinate (A129) at (-5, -3);
\coordinate (A130) at (-5, -2);
\coordinate (A131) at (-5, -1);
\coordinate (A132) at (-5, 0);
\coordinate (A133) at (-5, 1);
\coordinate (A134) at (-5, 2);
\coordinate (A135) at (-5, 3);
\coordinate (A136) at (-6, -2);
\coordinate (A137) at (-6, -1);
\coordinate (A138) at (-6, 0);
\coordinate (A139) at (-6, 1);
\coordinate (A140) at (-6, 2);
\coordinate (A141) at (-7, -1);
\coordinate (A142) at (-7, 0);
\coordinate (A143) at (-7, 1);
\coordinate (A144) at (-8, 0);
\colorlet{color1}{mycolor4}
\fill[color1!50] (A1) -- ($(A1)! 0.5!(A0)$) -- ($(A9)! 0.5!(A0)$) -- (A9) -- cycle;
\fill[color1!50] ($(A0)! 0.5!(A1)$) -- ($(A88)! 0.5!(A1)$) -- (A1) -- cycle;
\fill[color1!50] ($(A0)! 0.5!(A9)$) -- ($(A45)! 0.5!(A9)$) -- (A9) -- cycle;
\fill[color1!50] (A1) -- ($(A1)! 0.5!(A2)$) -- ($(A10)! 0.5!(A2)$) -- (A10) -- cycle;
\fill[color1!50] ($(A2)! 0.5!(A1)$) -- ($(A89)! 0.5!(A1)$) -- (A1) -- cycle;
\fill[color1!50] (A1) -- (A9) -- (A10) -- cycle;
\fill[color1!50] ($(A88)! 0.5!(A1)$) -- ($(A89)! 0.5!(A1)$) -- (A1) -- cycle;
\fill[color1!50] (A3) -- ($(A3)! 0.5!(A2)$) -- ($(A11)! 0.5!(A2)$) -- (A11) -- cycle;
\fill[color1!50] ($(A2)! 0.5!(A3)$) -- ($(A90)! 0.5!(A3)$) -- (A3) -- cycle;
\fill[color1!50] (A10) -- ($(A10)! 0.5!(A2)$) -- ($(A11)! 0.5!(A2)$) -- (A11) -- cycle;
\fill[color1!50] (A3) -- ($(A3)! 0.5!(A4)$) -- ($(A12)! 0.5!(A4)$) -- (A12) -- cycle;
\fill[color1!50] ($(A4)! 0.5!(A3)$) -- ($(A91)! 0.5!(A3)$) -- (A3) -- cycle;
\fill[color1!50] (A3) -- (A11) -- (A12) -- cycle;
\fill[color1!50] ($(A90)! 0.5!(A3)$) -- ($(A91)! 0.5!(A3)$) -- (A3) -- cycle;
\fill[color1!50] (A5) -- ($(A5)! 0.5!(A4)$) -- ($(A13)! 0.5!(A4)$) -- (A13) -- cycle;
\fill[color1!50] ($(A4)! 0.5!(A5)$) -- ($(A92)! 0.5!(A5)$) -- (A5) -- cycle;
\fill[color1!50] (A12) -- ($(A12)! 0.5!(A4)$) -- ($(A13)! 0.5!(A4)$) -- (A13) -- cycle;
\fill[color1!50] (A5) -- ($(A5)! 0.5!(A6)$) -- ($(A14)! 0.5!(A6)$) -- (A14) -- cycle;
\fill[color1!50] ($(A6)! 0.5!(A5)$) -- ($(A93)! 0.5!(A5)$) -- (A5) -- cycle;
\fill[color1!50] (A5) -- (A13) -- (A14) -- cycle;
\fill[color1!50] ($(A92)! 0.5!(A5)$) -- ($(A93)! 0.5!(A5)$) -- (A5) -- cycle;
\fill[color1!50] (A7) -- ($(A7)! 0.5!(A6)$) -- ($(A15)! 0.5!(A6)$) -- (A15) -- cycle;
\fill[color1!50] ($(A6)! 0.5!(A7)$) -- ($(A94)! 0.5!(A7)$) -- (A7) -- cycle;
\fill[color1!50] (A14) -- ($(A14)! 0.5!(A6)$) -- ($(A15)! 0.5!(A6)$) -- (A15) -- cycle;
\fill[color1!50] (A7) -- ($(A7)! 0.5!(A8)$) -- ($(A16)! 0.5!(A8)$) -- (A16) -- cycle;
\fill[color1!50] ($(A8)! 0.5!(A7)$) -- ($(A95)! 0.5!(A7)$) -- (A7) -- cycle;
\fill[color1!50] (A7) -- (A15) -- (A16) -- cycle;
\fill[color1!50] ($(A94)! 0.5!(A7)$) -- ($(A95)! 0.5!(A7)$) -- (A7) -- cycle;
\fill[color1!50] (A9) -- (A10) -- (A22) -- cycle;
\fill[color1!50] (A9) -- (A17) -- (A18) -- cycle;
\fill[color1!50] (A9) -- ($(A9)! 0.5!(A60)$) -- ($(A17)! 0.5!(A60)$) -- (A17) -- cycle;
\fill[color1!50] (A9) -- (A18) -- (A19) -- cycle;
\fill[color1!50] (A9) -- (A19) -- (A20) -- cycle;
\fill[color1!50] (A9) -- (A20) -- (A21) -- cycle;
\fill[color1!50] (A9) -- (A21) -- (A22) -- cycle;
\fill[color1!50] ($(A45)! 0.5!(A9)$) -- ($(A53)! 0.5!(A9)$) -- (A9) -- cycle;
\fill[color1!50] ($(A53)! 0.5!(A9)$) -- ($(A64)! 0.5!(A9)$) -- (A9) -- cycle;
\fill[color1!50] (A9) -- ($(A9)! 0.5!(A60)$) -- ($(A61)! 0.5!(A60)$) -- (A61) -- cycle;
\fill[color1!50] (A9) -- ($(A9)! 0.5!(A62)$) -- ($(A61)! 0.5!(A62)$) -- (A61) -- cycle;
\fill[color1!50] (A9) -- ($(A9)! 0.5!(A62)$) -- ($(A63)! 0.5!(A62)$) -- (A63) -- cycle;
\fill[color1!50] (A9) -- ($(A9)! 0.5!(A64)$) -- ($(A63)! 0.5!(A64)$) -- (A63) -- cycle;
\fill[color1!50] (A10) -- (A11) -- (A22) -- cycle;
\fill[color1!50] (A11) -- (A12) -- (A22) -- cycle;
\fill[color1!50] (A12) -- (A13) -- (A22) -- cycle;
\fill[color1!50] (A13) -- (A14) -- (A22) -- cycle;
\fill[color1!50] (A14) -- (A15) -- (A22) -- cycle;
\fill[color1!50] (A15) -- ($(A15)! 0.5!(A23)$) -- ($(A16)! 0.5!(A23)$) -- (A16) -- cycle;
\fill[color1!50] (A15) -- ($(A15)! 0.5!(A23)$) -- ($(A22)! 0.5!(A23)$) -- (A22) -- cycle;
\fill[color1!50] (A17) -- (A18) -- (A24) -- cycle;
\fill[color1!50] (A17) -- ($(A17)! 0.5!(A60)$) -- ($(A24)! 0.5!(A60)$) -- (A24) -- cycle;
\fill[color1!50] (A18) -- (A19) -- (A26) -- cycle;
\fill[color1!50] (A18) -- ($(A18)! 0.5!(A25)$) -- ($(A24)! 0.5!(A25)$) -- (A24) -- cycle;
\fill[color1!50] (A18) -- ($(A18)! 0.5!(A25)$) -- ($(A26)! 0.5!(A25)$) -- (A26) -- cycle;
\fill[color1!50] (A19) -- (A20) -- (A26) -- cycle;
\fill[color1!50] (A20) -- ($(A20)! 0.5!(A27)$) -- ($(A21)! 0.5!(A27)$) -- (A21) -- cycle;
\fill[color1!50] (A20) -- ($(A20)! 0.5!(A27)$) -- ($(A26)! 0.5!(A27)$) -- (A26) -- cycle;
\fill[color1!50] (A21) -- ($(A21)! 0.5!(A27)$) -- ($(A22)! 0.5!(A27)$) -- (A22) -- cycle;
\fill[color1!50] (A22) -- ($(A22)! 0.5!(A23)$) -- ($(A29)! 0.5!(A23)$) -- (A29) -- cycle;
\fill[color1!50] (A22) -- ($(A22)! 0.5!(A27)$) -- ($(A35)! 0.5!(A27)$) -- (A35) -- cycle;
\fill[color1!50] (A22) -- (A28) -- (A29) -- cycle;
\fill[color1!50] (A22) -- ($(A22)! 0.5!(A32)$) -- ($(A28)! 0.5!(A32)$) -- (A28) -- cycle;
\fill[color1!50] (A22) -- ($(A22)! 0.5!(A32)$) -- ($(A35)! 0.5!(A32)$) -- (A35) -- cycle;
\fill[color1!50] (A24) -- ($(A24)! 0.5!(A25)$) -- ($(A30)! 0.5!(A25)$) -- (A30) -- cycle;
\fill[color1!50] (A24) -- (A30) -- (A66) -- cycle;
\fill[color1!50] (A24) -- ($(A24)! 0.5!(A60)$) -- ($(A66)! 0.5!(A60)$) -- (A66) -- cycle;
\fill[color1!50] (A26) -- ($(A26)! 0.5!(A25)$) -- ($(A31)! 0.5!(A25)$) -- (A31) -- cycle;
\fill[color1!50] (A30) -- ($(A30)! 0.5!(A25)$) -- ($(A31)! 0.5!(A25)$) -- (A31) -- cycle;
\fill[color1!50] (A26) -- ($(A26)! 0.5!(A27)$) -- ($(A31)! 0.5!(A27)$) -- (A31) -- cycle;
\fill[color1!50] (A31) -- ($(A31)! 0.5!(A27)$) -- ($(A35)! 0.5!(A27)$) -- (A35) -- cycle;
\fill[color1!50] (A28) -- ($(A28)! 0.5!(A34)$) -- ($(A29)! 0.5!(A34)$) -- (A29) -- cycle;
\fill[color1!50] (A28) -- ($(A28)! 0.5!(A32)$) -- ($(A33)! 0.5!(A32)$) -- (A33) -- cycle;
\fill[color1!50] (A28) -- ($(A28)! 0.5!(A34)$) -- ($(A33)! 0.5!(A34)$) -- (A33) -- cycle;
\fill[color1!50] (A30) -- (A31) -- (A35) -- cycle;
\fill[color1!50] (A30) -- ($(A30)! 0.5!(A71)$) -- ($(A35)! 0.5!(A71)$) -- (A35) -- cycle;
\fill[color1!50] (A30) -- ($(A30)! 0.5!(A71)$) -- ($(A66)! 0.5!(A71)$) -- (A66) -- cycle;
\fill[color1!50] (A33) -- ($(A33)! 0.5!(A32)$) -- ($(A37)! 0.5!(A32)$) -- (A37) -- cycle;
\fill[color1!50] (A35) -- ($(A35)! 0.5!(A32)$) -- ($(A36)! 0.5!(A32)$) -- (A36) -- cycle;
\fill[color1!50] (A36) -- ($(A36)! 0.5!(A32)$) -- ($(A37)! 0.5!(A32)$) -- (A37) -- cycle;
\fill[color1!50] (A33) -- ($(A33)! 0.5!(A34)$) -- ($(A38)! 0.5!(A34)$) -- (A38) -- cycle;
\fill[color1!50] (A33) -- (A37) -- (A38) -- cycle;
\fill[color1!50] (A35) -- ($(A35)! 0.5!(A39)$) -- ($(A36)! 0.5!(A39)$) -- (A36) -- cycle;
\fill[color1!50] ($(A39)! 0.5!(A35)$) -- ($(A75)! 0.5!(A35)$) -- (A35) -- cycle;
\fill[color1!50] (A35) -- ($(A35)! 0.5!(A64)$) -- ($(A68)! 0.5!(A64)$) -- (A68) -- cycle;
\fill[color1!50] ($(A64)! 0.5!(A35)$) -- ($(A72)! 0.5!(A35)$) -- (A35) -- cycle;
\fill[color1!50] (A35) -- ($(A35)! 0.5!(A71)$) -- ($(A68)! 0.5!(A71)$) -- (A68) -- cycle;
\fill[color1!50] ($(A72)! 0.5!(A35)$) -- ($(A75)! 0.5!(A35)$) -- (A35) -- cycle;
\fill[color1!50] (A36) -- (A37) -- (A40) -- cycle;
\fill[color1!50] (A36) -- ($(A36)! 0.5!(A39)$) -- ($(A40)! 0.5!(A39)$) -- (A40) -- cycle;
\fill[color1!50] (A37) -- ($(A37)! 0.5!(A41)$) -- ($(A38)! 0.5!(A41)$) -- (A38) -- cycle;
\fill[color1!50] (A37) -- ($(A37)! 0.5!(A41)$) -- ($(A40)! 0.5!(A41)$) -- (A40) -- cycle;
\fill[color1!50] (A40) -- ($(A40)! 0.5!(A39)$) -- ($(A42)! 0.5!(A39)$) -- (A42) -- cycle;
\fill[color1!50] ($(A39)! 0.5!(A42)$) -- ($(A78)! 0.5!(A42)$) -- (A42) -- cycle;
\fill[color1!50] (A40) -- ($(A40)! 0.5!(A41)$) -- ($(A43)! 0.5!(A41)$) -- (A43) -- cycle;
\fill[color1!50] (A40) -- (A42) -- (A43) -- cycle;
\fill[color1!50] (A42) -- ($(A42)! 0.5!(A44)$) -- ($(A43)! 0.5!(A44)$) -- (A43) -- cycle;
\fill[color1!50] ($(A44)! 0.5!(A42)$) -- ($(A80)! 0.5!(A42)$) -- (A42) -- cycle;
\fill[color1!50] ($(A78)! 0.5!(A42)$) -- ($(A80)! 0.5!(A42)$) -- (A42) -- cycle;
\fill[color1!50] ($(A51)! 0.5!(A81)$) -- ($(A52)! 0.5!(A81)$) -- (A81) -- cycle;
\fill[color1!50] ($(A51)! 0.5!(A81)$) -- ($(A82)! 0.5!(A81)$) -- (A81) -- cycle;
\fill[color1!50] (A61) -- ($(A61)! 0.5!(A60)$) -- ($(A67)! 0.5!(A60)$) -- (A67) -- cycle;
\fill[color1!50] (A66) -- ($(A66)! 0.5!(A60)$) -- ($(A67)! 0.5!(A60)$) -- (A67) -- cycle;
\fill[color1!50] (A61) -- ($(A61)! 0.5!(A62)$) -- ($(A67)! 0.5!(A62)$) -- (A67) -- cycle;
\fill[color1!50] (A63) -- ($(A63)! 0.5!(A62)$) -- ($(A68)! 0.5!(A62)$) -- (A68) -- cycle;
\fill[color1!50] (A67) -- ($(A67)! 0.5!(A62)$) -- ($(A68)! 0.5!(A62)$) -- (A68) -- cycle;
\fill[color1!50] (A63) -- ($(A63)! 0.5!(A64)$) -- ($(A68)! 0.5!(A64)$) -- (A68) -- cycle;
\fill[color1!50] (A66) -- ($(A66)! 0.5!(A71)$) -- ($(A67)! 0.5!(A71)$) -- (A67) -- cycle;
\fill[color1!50] (A67) -- ($(A67)! 0.5!(A71)$) -- ($(A68)! 0.5!(A71)$) -- (A68) -- cycle;
\fill[color1!50] ($(A82)! 0.5!(A81)$) -- ($(A96)! 0.5!(A81)$) -- (A81) -- cycle;
\fill[color1!50] ($(A96)! 0.5!(A109)$) -- ($(A97)! 0.5!(A109)$) -- (A109) -- cycle;
\fill[color1!50] ($(A97)! 0.5!(A109)$) -- ($(A110)! 0.5!(A109)$) -- (A109) -- cycle;
\fill[color1!50] ($(A110)! 0.5!(A109)$) -- ($(A120)! 0.5!(A109)$) -- (A109) -- cycle;
\fill[color1!50] ($(A120)! 0.5!(A129)$) -- ($(A121)! 0.5!(A129)$) -- (A129) -- cycle;
\fill[color1!50] ($(A121)! 0.5!(A129)$) -- ($(A130)! 0.5!(A129)$) -- (A129) -- cycle;
\fill[color1!50] ($(A130)! 0.5!(A129)$) -- ($(A136)! 0.5!(A129)$) -- (A129) -- cycle;
\fill[color1!50] ($(A136)! 0.5!(A141)$) -- ($(A137)! 0.5!(A141)$) -- (A141) -- cycle;
\fill[color1!50] ($(A137)! 0.5!(A141)$) -- ($(A142)! 0.5!(A141)$) -- (A141) -- cycle;
\fill[color1!50] ($(A142)! 0.5!(A141)$) -- ($(A144)! 0.5!(A141)$) -- (A141) -- cycle;
\colorlet{color2}{mycolor5}
\fill[color2!50] ($(A18)! 0.5!(A25)$) -- ($(A24)! 0.5!(A25)$) -- (A25) -- cycle;
\fill[color2!50] ($(A18)! 0.5!(A25)$) -- ($(A26)! 0.5!(A25)$) -- (A25) -- cycle;
\fill[color2!50] ($(A24)! 0.5!(A25)$) -- ($(A30)! 0.5!(A25)$) -- (A25) -- cycle;
\fill[color2!50] ($(A26)! 0.5!(A25)$) -- ($(A31)! 0.5!(A25)$) -- (A25) -- cycle;
\fill[color2!50] ($(A30)! 0.5!(A25)$) -- ($(A31)! 0.5!(A25)$) -- (A25) -- cycle;
\colorlet{color3}{mycolor5}
\fill[color3!50] ($(A20)! 0.5!(A27)$) -- ($(A21)! 0.5!(A27)$) -- (A27) -- cycle;
\fill[color3!50] ($(A20)! 0.5!(A27)$) -- ($(A26)! 0.5!(A27)$) -- (A27) -- cycle;
\fill[color3!50] ($(A21)! 0.5!(A27)$) -- ($(A22)! 0.5!(A27)$) -- (A27) -- cycle;
\fill[color3!50] ($(A22)! 0.5!(A27)$) -- ($(A35)! 0.5!(A27)$) -- (A27) -- cycle;
\fill[color3!50] ($(A26)! 0.5!(A27)$) -- ($(A31)! 0.5!(A27)$) -- (A27) -- cycle;
\fill[color3!50] ($(A31)! 0.5!(A27)$) -- ($(A35)! 0.5!(A27)$) -- (A27) -- cycle;
\colorlet{color4}{mycolor3}
\fill[color4!50] ($(A22)! 0.5!(A32)$) -- ($(A28)! 0.5!(A32)$) -- (A32) -- cycle;
\fill[color4!50] ($(A22)! 0.5!(A32)$) -- ($(A35)! 0.5!(A32)$) -- (A32) -- cycle;
\fill[color4!50] ($(A28)! 0.5!(A32)$) -- ($(A33)! 0.5!(A32)$) -- (A32) -- cycle;
\fill[color4!50] ($(A33)! 0.5!(A32)$) -- ($(A37)! 0.5!(A32)$) -- (A32) -- cycle;
\fill[color4!50] ($(A35)! 0.5!(A32)$) -- ($(A36)! 0.5!(A32)$) -- (A32) -- cycle;
\fill[color4!50] ($(A36)! 0.5!(A32)$) -- ($(A37)! 0.5!(A32)$) -- (A32) -- cycle;
\colorlet{color5}{mycolor3}
\fill[color5!50] ($(A46)! 0.5!(A54)$) -- ($(A47)! 0.5!(A54)$) -- (A54) -- cycle;
\fill[color5!50] ($(A46)! 0.5!(A54)$) -- ($(A53)! 0.5!(A54)$) -- (A54) -- cycle;
\fill[color5!50] ($(A47)! 0.5!(A54)$) -- ($(A55)! 0.5!(A54)$) -- (A54) -- cycle;
\fill[color5!50] ($(A53)! 0.5!(A54)$) -- ($(A64)! 0.5!(A54)$) -- (A54) -- cycle;
\fill[color5!50] ($(A55)! 0.5!(A54)$) -- ($(A64)! 0.5!(A54)$) -- (A54) -- cycle;
\colorlet{color6}{mycolor3}
\fill[color6!50] ($(A48)! 0.5!(A56)$) -- ($(A49)! 0.5!(A56)$) -- (A56) -- cycle;
\fill[color6!50] ($(A48)! 0.5!(A56)$) -- ($(A55)! 0.5!(A56)$) -- (A56) -- cycle;
\fill[color6!50] ($(A49)! 0.5!(A56)$) -- ($(A57)! 0.5!(A56)$) -- (A56) -- cycle;
\fill[color6!50] ($(A55)! 0.5!(A56)$) -- ($(A64)! 0.5!(A56)$) -- (A56) -- cycle;
\fill[color6!50] ($(A57)! 0.5!(A56)$) -- ($(A64)! 0.5!(A56)$) -- (A56) -- cycle;
\colorlet{color7}{mycolor3}
\fill[color7!50] ($(A50)! 0.5!(A58)$) -- ($(A51)! 0.5!(A58)$) -- (A58) -- cycle;
\fill[color7!50] ($(A50)! 0.5!(A58)$) -- ($(A57)! 0.5!(A58)$) -- (A58) -- cycle;
\fill[color7!50] ($(A51)! 0.5!(A58)$) -- ($(A59)! 0.5!(A58)$) -- (A58) -- cycle;
\fill[color7!50] ($(A57)! 0.5!(A58)$) -- ($(A64)! 0.5!(A58)$) -- (A58) -- cycle;
\fill[color7!50] ($(A59)! 0.5!(A58)$) -- ($(A65)! 0.5!(A58)$) -- (A58) -- cycle;
\fill[color7!50] ($(A64)! 0.5!(A58)$) -- ($(A65)! 0.5!(A58)$) -- (A58) -- cycle;
\colorlet{color8}{mycolor5}
\fill[color8!50] ($(A9)! 0.5!(A60)$) -- ($(A17)! 0.5!(A60)$) -- (A60) -- cycle;
\fill[color8!50] ($(A9)! 0.5!(A60)$) -- ($(A61)! 0.5!(A60)$) -- (A60) -- cycle;
\fill[color8!50] ($(A17)! 0.5!(A60)$) -- ($(A24)! 0.5!(A60)$) -- (A60) -- cycle;
\fill[color8!50] ($(A24)! 0.5!(A60)$) -- ($(A66)! 0.5!(A60)$) -- (A60) -- cycle;
\fill[color8!50] ($(A61)! 0.5!(A60)$) -- ($(A67)! 0.5!(A60)$) -- (A60) -- cycle;
\fill[color8!50] ($(A66)! 0.5!(A60)$) -- ($(A67)! 0.5!(A60)$) -- (A60) -- cycle;
\colorlet{color9}{mycolor5}
\fill[color9!50] ($(A9)! 0.5!(A62)$) -- ($(A61)! 0.5!(A62)$) -- (A62) -- cycle;
\fill[color9!50] ($(A9)! 0.5!(A62)$) -- ($(A63)! 0.5!(A62)$) -- (A62) -- cycle;
\fill[color9!50] ($(A61)! 0.5!(A62)$) -- ($(A67)! 0.5!(A62)$) -- (A62) -- cycle;
\fill[color9!50] ($(A63)! 0.5!(A62)$) -- ($(A68)! 0.5!(A62)$) -- (A62) -- cycle;
\fill[color9!50] ($(A67)! 0.5!(A62)$) -- ($(A68)! 0.5!(A62)$) -- (A62) -- cycle;
\colorlet{color10}{mycolor3}
\fill[color10!50] ($(A64)! 0.5!(A69)$) -- ($(A70)! 0.5!(A69)$) -- (A69) -- cycle;
\fill[color10!50] ($(A64)! 0.5!(A69)$) -- ($(A72)! 0.5!(A69)$) -- (A69) -- cycle;
\fill[color10!50] ($(A70)! 0.5!(A69)$) -- ($(A74)! 0.5!(A69)$) -- (A69) -- cycle;
\fill[color10!50] ($(A72)! 0.5!(A69)$) -- ($(A73)! 0.5!(A69)$) -- (A69) -- cycle;
\fill[color10!50] ($(A73)! 0.5!(A69)$) -- ($(A74)! 0.5!(A69)$) -- (A69) -- cycle;
\colorlet{color11}{mycolor5}
\fill[color11!50] ($(A30)! 0.5!(A71)$) -- ($(A35)! 0.5!(A71)$) -- (A71) -- cycle;
\fill[color11!50] ($(A30)! 0.5!(A71)$) -- ($(A66)! 0.5!(A71)$) -- (A71) -- cycle;
\fill[color11!50] ($(A35)! 0.5!(A71)$) -- ($(A68)! 0.5!(A71)$) -- (A71) -- cycle;
\fill[color11!50] ($(A66)! 0.5!(A71)$) -- ($(A67)! 0.5!(A71)$) -- (A71) -- cycle;
\fill[color11!50] ($(A67)! 0.5!(A71)$) -- ($(A68)! 0.5!(A71)$) -- (A71) -- cycle;
\colorlet{color12}{mycolor3}
\fill[color12!50] ($(A72)! 0.5!(A76)$) -- ($(A73)! 0.5!(A76)$) -- (A76) -- cycle;
\fill[color12!50] ($(A72)! 0.5!(A76)$) -- ($(A75)! 0.5!(A76)$) -- (A76) -- cycle;
\fill[color12!50] ($(A73)! 0.5!(A76)$) -- ($(A77)! 0.5!(A76)$) -- (A76) -- cycle;
\fill[color12!50] ($(A75)! 0.5!(A76)$) -- ($(A78)! 0.5!(A76)$) -- (A76) -- cycle;
\fill[color12!50] ($(A77)! 0.5!(A76)$) -- ($(A79)! 0.5!(A76)$) -- (A76) -- cycle;
\fill[color12!50] ($(A78)! 0.5!(A76)$) -- ($(A79)! 0.5!(A76)$) -- (A76) -- cycle;
\colorlet{color13}{mycolor3}
\fill[color13!50] ($(A49)! 0.5!(A83)$) -- ($(A50)! 0.5!(A83)$) -- (A83) -- cycle;
\fill[color13!50] ($(A49)! 0.5!(A83)$) -- ($(A84)! 0.5!(A83)$) -- (A83) -- cycle;
\fill[color13!50] ($(A50)! 0.5!(A83)$) -- ($(A82)! 0.5!(A83)$) -- (A83) -- cycle;
\fill[color13!50] ($(A82)! 0.5!(A83)$) -- ($(A97)! 0.5!(A83)$) -- (A83) -- cycle;
\fill[color13!50] ($(A84)! 0.5!(A83)$) -- ($(A97)! 0.5!(A83)$) -- (A83) -- cycle;
\colorlet{color14}{mycolor3}
\fill[color14!50] ($(A47)! 0.5!(A85)$) -- ($(A48)! 0.5!(A85)$) -- (A85) -- cycle;
\fill[color14!50] ($(A47)! 0.5!(A85)$) -- ($(A86)! 0.5!(A85)$) -- (A85) -- cycle;
\fill[color14!50] ($(A48)! 0.5!(A85)$) -- ($(A84)! 0.5!(A85)$) -- (A85) -- cycle;
\fill[color14!50] ($(A84)! 0.5!(A85)$) -- ($(A97)! 0.5!(A85)$) -- (A85) -- cycle;
\fill[color14!50] ($(A86)! 0.5!(A85)$) -- ($(A97)! 0.5!(A85)$) -- (A85) -- cycle;
\colorlet{color15}{mycolor3}
\fill[color15!50] ($(A45)! 0.5!(A87)$) -- ($(A46)! 0.5!(A87)$) -- (A87) -- cycle;
\fill[color15!50] ($(A45)! 0.5!(A87)$) -- ($(A88)! 0.5!(A87)$) -- (A87) -- cycle;
\fill[color15!50] ($(A46)! 0.5!(A87)$) -- ($(A86)! 0.5!(A87)$) -- (A87) -- cycle;
\fill[color15!50] ($(A86)! 0.5!(A87)$) -- ($(A97)! 0.5!(A87)$) -- (A87) -- cycle;
\fill[color15!50] ($(A88)! 0.5!(A87)$) -- ($(A97)! 0.5!(A87)$) -- (A87) -- cycle;
\colorlet{color16}{mycolor5}
\fill[color16!50] ($(A88)! 0.5!(A98)$) -- ($(A97)! 0.5!(A98)$) -- (A98) -- cycle;
\fill[color16!50] ($(A88)! 0.5!(A98)$) -- ($(A99)! 0.5!(A98)$) -- (A98) -- cycle;
\fill[color16!50] ($(A97)! 0.5!(A98)$) -- ($(A111)! 0.5!(A98)$) -- (A98) -- cycle;
\fill[color16!50] ($(A99)! 0.5!(A98)$) -- ($(A111)! 0.5!(A98)$) -- (A98) -- cycle;
\colorlet{color17}{mycolor5}
\fill[color17!50] ($(A88)! 0.5!(A100)$) -- ($(A99)! 0.5!(A100)$) -- (A100) -- cycle;
\fill[color17!50] ($(A88)! 0.5!(A100)$) -- ($(A101)! 0.5!(A100)$) -- (A100) -- cycle;
\fill[color17!50] ($(A99)! 0.5!(A100)$) -- ($(A112)! 0.5!(A100)$) -- (A100) -- cycle;
\fill[color17!50] ($(A101)! 0.5!(A100)$) -- ($(A112)! 0.5!(A100)$) -- (A100) -- cycle;
\colorlet{color18}{mycolor3}
\fill[color18!50] ($(A88)! 0.5!(A107)$) -- ($(A89)! 0.5!(A107)$) -- (A107) -- cycle;
\fill[color18!50] ($(A88)! 0.5!(A102)$) -- ($(A101)! 0.5!(A102)$) -- (A102) -- cycle;
\fill[color18!50] (A102) -- ($(A102)! 0.5!(A88)$) -- ($(A103)! 0.5!(A88)$) -- (A103) -- cycle;
\fill[color18!50] (A103) -- ($(A103)! 0.5!(A88)$) -- ($(A104)! 0.5!(A88)$) -- (A104) -- cycle;
\fill[color18!50] (A104) -- ($(A104)! 0.5!(A88)$) -- ($(A105)! 0.5!(A88)$) -- (A105) -- cycle;
\fill[color18!50] (A105) -- ($(A105)! 0.5!(A88)$) -- ($(A106)! 0.5!(A88)$) -- (A106) -- cycle;
\fill[color18!50] (A106) -- ($(A106)! 0.5!(A88)$) -- ($(A107)! 0.5!(A88)$) -- (A107) -- cycle;
\fill[color18!50] ($(A89)! 0.5!(A107)$) -- ($(A90)! 0.5!(A107)$) -- (A107) -- cycle;
\fill[color18!50] ($(A90)! 0.5!(A107)$) -- ($(A91)! 0.5!(A107)$) -- (A107) -- cycle;
\fill[color18!50] ($(A91)! 0.5!(A107)$) -- ($(A92)! 0.5!(A107)$) -- (A107) -- cycle;
\fill[color18!50] ($(A92)! 0.5!(A107)$) -- ($(A93)! 0.5!(A107)$) -- (A107) -- cycle;
\fill[color18!50] ($(A93)! 0.5!(A107)$) -- ($(A94)! 0.5!(A107)$) -- (A107) -- cycle;
\fill[color18!50] ($(A94)! 0.5!(A107)$) -- ($(A108)! 0.5!(A107)$) -- (A107) -- cycle;
\fill[color18!50] ($(A101)! 0.5!(A102)$) -- ($(A114)! 0.5!(A102)$) -- (A102) -- cycle;
\fill[color18!50] (A102) -- ($(A102)! 0.5!(A114)$) -- ($(A103)! 0.5!(A114)$) -- (A103) -- cycle;
\fill[color18!50] (A103) -- ($(A103)! 0.5!(A116)$) -- ($(A104)! 0.5!(A116)$) -- (A104) -- cycle;
\fill[color18!50] (A103) -- ($(A103)! 0.5!(A114)$) -- ($(A115)! 0.5!(A114)$) -- (A115) -- cycle;
\fill[color18!50] (A103) -- ($(A103)! 0.5!(A116)$) -- ($(A115)! 0.5!(A116)$) -- (A115) -- cycle;
\fill[color18!50] (A104) -- ($(A104)! 0.5!(A116)$) -- ($(A105)! 0.5!(A116)$) -- (A105) -- cycle;
\fill[color18!50] (A105) -- (A106) -- (A117) -- cycle;
\fill[color18!50] (A105) -- ($(A105)! 0.5!(A116)$) -- ($(A117)! 0.5!(A116)$) -- (A117) -- cycle;
\fill[color18!50] (A106) -- (A107) -- (A117) -- cycle;
\fill[color18!50] ($(A108)! 0.5!(A107)$) -- ($(A119)! 0.5!(A107)$) -- (A107) -- cycle;
\fill[color18!50] (A107) -- ($(A107)! 0.5!(A132)$) -- ($(A117)! 0.5!(A132)$) -- (A117) -- cycle;
\fill[color18!50] ($(A118)! 0.5!(A107)$) -- ($(A119)! 0.5!(A107)$) -- (A107) -- cycle;
\fill[color18!50] ($(A118)! 0.5!(A107)$) -- ($(A126)! 0.5!(A107)$) -- (A107) -- cycle;
\fill[color18!50] ($(A126)! 0.5!(A107)$) -- ($(A132)! 0.5!(A107)$) -- (A107) -- cycle;
\fill[color18!50] ($(A113)! 0.5!(A124)$) -- ($(A114)! 0.5!(A124)$) -- (A124) -- cycle;
\fill[color18!50] ($(A113)! 0.5!(A124)$) -- ($(A123)! 0.5!(A124)$) -- (A124) -- cycle;
\fill[color18!50] (A115) -- ($(A115)! 0.5!(A114)$) -- ($(A124)! 0.5!(A114)$) -- (A124) -- cycle;
\fill[color18!50] (A115) -- ($(A115)! 0.5!(A116)$) -- ($(A125)! 0.5!(A116)$) -- (A125) -- cycle;
\fill[color18!50] (A115) -- (A124) -- (A125) -- cycle;
\fill[color18!50] (A117) -- ($(A117)! 0.5!(A116)$) -- ($(A125)! 0.5!(A116)$) -- (A125) -- cycle;
\fill[color18!50] (A117) -- ($(A117)! 0.5!(A132)$) -- ($(A125)! 0.5!(A132)$) -- (A125) -- cycle;
\fill[color18!50] ($(A123)! 0.5!(A124)$) -- ($(A132)! 0.5!(A124)$) -- (A124) -- cycle;
\fill[color18!50] (A124) -- ($(A124)! 0.5!(A132)$) -- ($(A125)! 0.5!(A132)$) -- (A125) -- cycle;
\colorlet{color19}{mycolor5}
\fill[color19!50] ($(A103)! 0.5!(A116)$) -- ($(A104)! 0.5!(A116)$) -- (A116) -- cycle;
\fill[color19!50] ($(A103)! 0.5!(A116)$) -- ($(A115)! 0.5!(A116)$) -- (A116) -- cycle;
\fill[color19!50] ($(A104)! 0.5!(A116)$) -- ($(A105)! 0.5!(A116)$) -- (A116) -- cycle;
\fill[color19!50] ($(A105)! 0.5!(A116)$) -- ($(A117)! 0.5!(A116)$) -- (A116) -- cycle;
\fill[color19!50] ($(A115)! 0.5!(A116)$) -- ($(A125)! 0.5!(A116)$) -- (A116) -- cycle;
\fill[color19!50] ($(A117)! 0.5!(A116)$) -- ($(A125)! 0.5!(A116)$) -- (A116) -- cycle;
\colorlet{color20}{mycolor3}
\fill[color20!50] ($(A118)! 0.5!(A127)$) -- ($(A126)! 0.5!(A127)$) -- (A127) -- cycle;
\fill[color20!50] ($(A118)! 0.5!(A127)$) -- ($(A128)! 0.5!(A127)$) -- (A127) -- cycle;
\fill[color20!50] ($(A126)! 0.5!(A127)$) -- ($(A134)! 0.5!(A127)$) -- (A127) -- cycle;
\fill[color20!50] ($(A128)! 0.5!(A127)$) -- ($(A135)! 0.5!(A127)$) -- (A127) -- cycle;
\fill[color20!50] ($(A134)! 0.5!(A127)$) -- ($(A135)! 0.5!(A127)$) -- (A127) -- cycle;
\colorlet{color21}{mycolor3}
\fill[color21!50] ($(A122)! 0.5!(A131)$) -- ($(A130)! 0.5!(A131)$) -- (A131) -- cycle;
\fill[color21!50] ($(A122)! 0.5!(A131)$) -- ($(A132)! 0.5!(A131)$) -- (A131) -- cycle;
\fill[color21!50] ($(A130)! 0.5!(A131)$) -- ($(A137)! 0.5!(A131)$) -- (A131) -- cycle;
\fill[color21!50] ($(A132)! 0.5!(A131)$) -- ($(A138)! 0.5!(A131)$) -- (A131) -- cycle;
\fill[color21!50] ($(A137)! 0.5!(A131)$) -- ($(A138)! 0.5!(A131)$) -- (A131) -- cycle;
\colorlet{color22}{mycolor3}
\fill[color22!50] ($(A133)! 0.5!(A139)$) -- ($(A134)! 0.5!(A139)$) -- (A139) -- cycle;
\fill[color22!50] ($(A133)! 0.5!(A139)$) -- ($(A138)! 0.5!(A139)$) -- (A139) -- cycle;
\fill[color22!50] ($(A134)! 0.5!(A139)$) -- ($(A140)! 0.5!(A139)$) -- (A139) -- cycle;
\fill[color22!50] ($(A138)! 0.5!(A139)$) -- ($(A142)! 0.5!(A139)$) -- (A139) -- cycle;
\fill[color22!50] ($(A140)! 0.5!(A139)$) -- ($(A143)! 0.5!(A139)$) -- (A139) -- cycle;
\fill[color22!50] ($(A142)! 0.5!(A139)$) -- ($(A143)! 0.5!(A139)$) -- (A139) -- cycle;
\draw[color22, thick] ($(A133)! 0.5!(A139)$) -- ($(A134)! 0.5!(A139)$);
\draw[color22, thick] ($(A133)! 0.5!(A139)$) -- ($(A138)! 0.5!(A139)$);
\draw[color22, thick] ($(A134)! 0.5!(A139)$) -- ($(A140)! 0.5!(A139)$);
\draw[color22, thick] ($(A138)! 0.5!(A139)$) -- ($(A142)! 0.5!(A139)$);
\draw[color22, thick] ($(A140)! 0.5!(A139)$) -- ($(A143)! 0.5!(A139)$);
\draw[color22, thick] ($(A142)! 0.5!(A139)$) -- ($(A143)! 0.5!(A139)$);
\draw[color21, thick] ($(A122)! 0.5!(A131)$) -- ($(A130)! 0.5!(A131)$);
\draw[color21, thick] ($(A122)! 0.5!(A131)$) -- ($(A132)! 0.5!(A131)$);
\draw[color21, thick] ($(A130)! 0.5!(A131)$) -- ($(A137)! 0.5!(A131)$);
\draw[color21, thick] ($(A132)! 0.5!(A131)$) -- ($(A138)! 0.5!(A131)$);
\draw[color21, thick] ($(A137)! 0.5!(A131)$) -- ($(A138)! 0.5!(A131)$);
\draw[color20, thick] ($(A118)! 0.5!(A127)$) -- ($(A126)! 0.5!(A127)$);
\draw[color20, thick] ($(A118)! 0.5!(A127)$) -- ($(A128)! 0.5!(A127)$);
\draw[color20, thick] ($(A126)! 0.5!(A127)$) -- ($(A134)! 0.5!(A127)$);
\draw[color20, thick] ($(A128)! 0.5!(A127)$) -- ($(A135)! 0.5!(A127)$);
\draw[color20, thick] ($(A134)! 0.5!(A127)$) -- ($(A135)! 0.5!(A127)$);
\draw[color18, thick] ($(A88)! 0.5!(A107)$) -- ($(A89)! 0.5!(A107)$);
\draw[color18, thick] ($(A88)! 0.5!(A102)$) -- ($(A101)! 0.5!(A102)$);
\draw[color18, thick] ($(A102)! 0.5!(A88)$) -- ($(A103)! 0.5!(A88)$);
\draw[color18, thick] ($(A103)! 0.5!(A88)$) -- ($(A104)! 0.5!(A88)$);
\draw[color18, thick] ($(A104)! 0.5!(A88)$) -- ($(A105)! 0.5!(A88)$);
\draw[color18, thick] ($(A105)! 0.5!(A88)$) -- ($(A106)! 0.5!(A88)$);
\draw[color18, thick] ($(A106)! 0.5!(A88)$) -- ($(A107)! 0.5!(A88)$);
\draw[color18, thick] ($(A89)! 0.5!(A107)$) -- ($(A90)! 0.5!(A107)$);
\draw[color18, thick] ($(A90)! 0.5!(A107)$) -- ($(A91)! 0.5!(A107)$);
\draw[color18, thick] ($(A91)! 0.5!(A107)$) -- ($(A92)! 0.5!(A107)$);
\draw[color18, thick] ($(A92)! 0.5!(A107)$) -- ($(A93)! 0.5!(A107)$);
\draw[color18, thick] ($(A93)! 0.5!(A107)$) -- ($(A94)! 0.5!(A107)$);
\draw[color18, thick] ($(A94)! 0.5!(A107)$) -- ($(A108)! 0.5!(A107)$);
\draw[color18, thick] ($(A101)! 0.5!(A102)$) -- ($(A114)! 0.5!(A102)$);
\draw[color18, thick] ($(A102)! 0.5!(A114)$) -- ($(A103)! 0.5!(A114)$);
\draw[color18, thick] ($(A103)! 0.5!(A116)$) -- ($(A104)! 0.5!(A116)$);
\draw[color18, thick] ($(A103)! 0.5!(A114)$) -- ($(A115)! 0.5!(A114)$);
\draw[color18, thick] ($(A103)! 0.5!(A116)$) -- ($(A115)! 0.5!(A116)$);
\draw[color18, thick] ($(A104)! 0.5!(A116)$) -- ($(A105)! 0.5!(A116)$);
\draw[color18, thick] ($(A105)! 0.5!(A116)$) -- ($(A117)! 0.5!(A116)$);
\draw[color18, thick] ($(A108)! 0.5!(A107)$) -- ($(A119)! 0.5!(A107)$);
\draw[color18, thick] ($(A107)! 0.5!(A132)$) -- ($(A117)! 0.5!(A132)$);
\draw[color18, thick] ($(A118)! 0.5!(A107)$) -- ($(A119)! 0.5!(A107)$);
\draw[color18, thick] ($(A118)! 0.5!(A107)$) -- ($(A126)! 0.5!(A107)$);
\draw[color18, thick] ($(A126)! 0.5!(A107)$) -- ($(A132)! 0.5!(A107)$);
\draw[color18, thick] ($(A113)! 0.5!(A124)$) -- ($(A114)! 0.5!(A124)$);
\draw[color18, thick] ($(A113)! 0.5!(A124)$) -- ($(A123)! 0.5!(A124)$);
\draw[color18, thick] ($(A115)! 0.5!(A114)$) -- ($(A124)! 0.5!(A114)$);
\draw[color18, thick] ($(A115)! 0.5!(A116)$) -- ($(A125)! 0.5!(A116)$);
\draw[color18, thick] ($(A117)! 0.5!(A116)$) -- ($(A125)! 0.5!(A116)$);
\draw[color18, thick] ($(A117)! 0.5!(A132)$) -- ($(A125)! 0.5!(A132)$);
\draw[color18, thick] ($(A123)! 0.5!(A124)$) -- ($(A132)! 0.5!(A124)$);
\draw[color18, thick] ($(A124)! 0.5!(A132)$) -- ($(A125)! 0.5!(A132)$);
\draw[color19, thick] ($(A103)! 0.5!(A116)$) -- ($(A104)! 0.5!(A116)$);
\draw[color19, thick] ($(A103)! 0.5!(A116)$) -- ($(A115)! 0.5!(A116)$);
\draw[color19, thick] ($(A104)! 0.5!(A116)$) -- ($(A105)! 0.5!(A116)$);
\draw[color19, thick] ($(A105)! 0.5!(A116)$) -- ($(A117)! 0.5!(A116)$);
\draw[color19, thick] ($(A115)! 0.5!(A116)$) -- ($(A125)! 0.5!(A116)$);
\draw[color19, thick] ($(A117)! 0.5!(A116)$) -- ($(A125)! 0.5!(A116)$);
\draw[color17, thick] ($(A88)! 0.5!(A100)$) -- ($(A99)! 0.5!(A100)$);
\draw[color17, thick] ($(A88)! 0.5!(A100)$) -- ($(A101)! 0.5!(A100)$);
\draw[color17, thick] ($(A99)! 0.5!(A100)$) -- ($(A112)! 0.5!(A100)$);
\draw[color17, thick] ($(A101)! 0.5!(A100)$) -- ($(A112)! 0.5!(A100)$);
\draw[color16, thick] ($(A88)! 0.5!(A98)$) -- ($(A97)! 0.5!(A98)$);
\draw[color16, thick] ($(A88)! 0.5!(A98)$) -- ($(A99)! 0.5!(A98)$);
\draw[color16, thick] ($(A97)! 0.5!(A98)$) -- ($(A111)! 0.5!(A98)$);
\draw[color16, thick] ($(A99)! 0.5!(A98)$) -- ($(A111)! 0.5!(A98)$);
\draw[color15, thick] ($(A45)! 0.5!(A87)$) -- ($(A46)! 0.5!(A87)$);
\draw[color15, thick] ($(A45)! 0.5!(A87)$) -- ($(A88)! 0.5!(A87)$);
\draw[color15, thick] ($(A46)! 0.5!(A87)$) -- ($(A86)! 0.5!(A87)$);
\draw[color15, thick] ($(A86)! 0.5!(A87)$) -- ($(A97)! 0.5!(A87)$);
\draw[color15, thick] ($(A88)! 0.5!(A87)$) -- ($(A97)! 0.5!(A87)$);
\draw[color14, thick] ($(A47)! 0.5!(A85)$) -- ($(A48)! 0.5!(A85)$);
\draw[color14, thick] ($(A47)! 0.5!(A85)$) -- ($(A86)! 0.5!(A85)$);
\draw[color14, thick] ($(A48)! 0.5!(A85)$) -- ($(A84)! 0.5!(A85)$);
\draw[color14, thick] ($(A84)! 0.5!(A85)$) -- ($(A97)! 0.5!(A85)$);
\draw[color14, thick] ($(A86)! 0.5!(A85)$) -- ($(A97)! 0.5!(A85)$);
\draw[color13, thick] ($(A49)! 0.5!(A83)$) -- ($(A50)! 0.5!(A83)$);
\draw[color13, thick] ($(A49)! 0.5!(A83)$) -- ($(A84)! 0.5!(A83)$);
\draw[color13, thick] ($(A50)! 0.5!(A83)$) -- ($(A82)! 0.5!(A83)$);
\draw[color13, thick] ($(A82)! 0.5!(A83)$) -- ($(A97)! 0.5!(A83)$);
\draw[color13, thick] ($(A84)! 0.5!(A83)$) -- ($(A97)! 0.5!(A83)$);
\draw[color12, thick] ($(A72)! 0.5!(A76)$) -- ($(A73)! 0.5!(A76)$);
\draw[color12, thick] ($(A72)! 0.5!(A76)$) -- ($(A75)! 0.5!(A76)$);
\draw[color12, thick] ($(A73)! 0.5!(A76)$) -- ($(A77)! 0.5!(A76)$);
\draw[color12, thick] ($(A75)! 0.5!(A76)$) -- ($(A78)! 0.5!(A76)$);
\draw[color12, thick] ($(A77)! 0.5!(A76)$) -- ($(A79)! 0.5!(A76)$);
\draw[color12, thick] ($(A78)! 0.5!(A76)$) -- ($(A79)! 0.5!(A76)$);
\draw[color10, thick] ($(A64)! 0.5!(A69)$) -- ($(A70)! 0.5!(A69)$);
\draw[color10, thick] ($(A64)! 0.5!(A69)$) -- ($(A72)! 0.5!(A69)$);
\draw[color10, thick] ($(A70)! 0.5!(A69)$) -- ($(A74)! 0.5!(A69)$);
\draw[color10, thick] ($(A72)! 0.5!(A69)$) -- ($(A73)! 0.5!(A69)$);
\draw[color10, thick] ($(A73)! 0.5!(A69)$) -- ($(A74)! 0.5!(A69)$);
\draw[color7, thick] ($(A50)! 0.5!(A58)$) -- ($(A51)! 0.5!(A58)$);
\draw[color7, thick] ($(A50)! 0.5!(A58)$) -- ($(A57)! 0.5!(A58)$);
\draw[color7, thick] ($(A51)! 0.5!(A58)$) -- ($(A59)! 0.5!(A58)$);
\draw[color7, thick] ($(A57)! 0.5!(A58)$) -- ($(A64)! 0.5!(A58)$);
\draw[color7, thick] ($(A59)! 0.5!(A58)$) -- ($(A65)! 0.5!(A58)$);
\draw[color7, thick] ($(A64)! 0.5!(A58)$) -- ($(A65)! 0.5!(A58)$);
\draw[color6, thick] ($(A48)! 0.5!(A56)$) -- ($(A49)! 0.5!(A56)$);
\draw[color6, thick] ($(A48)! 0.5!(A56)$) -- ($(A55)! 0.5!(A56)$);
\draw[color6, thick] ($(A49)! 0.5!(A56)$) -- ($(A57)! 0.5!(A56)$);
\draw[color6, thick] ($(A55)! 0.5!(A56)$) -- ($(A64)! 0.5!(A56)$);
\draw[color6, thick] ($(A57)! 0.5!(A56)$) -- ($(A64)! 0.5!(A56)$);
\draw[color5, thick] ($(A46)! 0.5!(A54)$) -- ($(A47)! 0.5!(A54)$);
\draw[color5, thick] ($(A46)! 0.5!(A54)$) -- ($(A53)! 0.5!(A54)$);
\draw[color5, thick] ($(A47)! 0.5!(A54)$) -- ($(A55)! 0.5!(A54)$);
\draw[color5, thick] ($(A53)! 0.5!(A54)$) -- ($(A64)! 0.5!(A54)$);
\draw[color5, thick] ($(A55)! 0.5!(A54)$) -- ($(A64)! 0.5!(A54)$);
\draw[color1, thick] ($(A1)! 0.5!(A0)$) -- ($(A9)! 0.5!(A0)$);
\draw[color1, thick] ($(A0)! 0.5!(A1)$) -- ($(A88)! 0.5!(A1)$);
\draw[color1, thick] ($(A0)! 0.5!(A9)$) -- ($(A45)! 0.5!(A9)$);
\draw[color1, thick] ($(A1)! 0.5!(A2)$) -- ($(A10)! 0.5!(A2)$);
\draw[color1, thick] ($(A2)! 0.5!(A1)$) -- ($(A89)! 0.5!(A1)$);
\draw[color1, thick] ($(A88)! 0.5!(A1)$) -- ($(A89)! 0.5!(A1)$);
\draw[color1, thick] ($(A3)! 0.5!(A2)$) -- ($(A11)! 0.5!(A2)$);
\draw[color1, thick] ($(A2)! 0.5!(A3)$) -- ($(A90)! 0.5!(A3)$);
\draw[color1, thick] ($(A10)! 0.5!(A2)$) -- ($(A11)! 0.5!(A2)$);
\draw[color1, thick] ($(A3)! 0.5!(A4)$) -- ($(A12)! 0.5!(A4)$);
\draw[color1, thick] ($(A4)! 0.5!(A3)$) -- ($(A91)! 0.5!(A3)$);
\draw[color1, thick] ($(A90)! 0.5!(A3)$) -- ($(A91)! 0.5!(A3)$);
\draw[color1, thick] ($(A5)! 0.5!(A4)$) -- ($(A13)! 0.5!(A4)$);
\draw[color1, thick] ($(A4)! 0.5!(A5)$) -- ($(A92)! 0.5!(A5)$);
\draw[color1, thick] ($(A12)! 0.5!(A4)$) -- ($(A13)! 0.5!(A4)$);
\draw[color1, thick] ($(A5)! 0.5!(A6)$) -- ($(A14)! 0.5!(A6)$);
\draw[color1, thick] ($(A6)! 0.5!(A5)$) -- ($(A93)! 0.5!(A5)$);
\draw[color1, thick] ($(A92)! 0.5!(A5)$) -- ($(A93)! 0.5!(A5)$);
\draw[color1, thick] ($(A7)! 0.5!(A6)$) -- ($(A15)! 0.5!(A6)$);
\draw[color1, thick] ($(A6)! 0.5!(A7)$) -- ($(A94)! 0.5!(A7)$);
\draw[color1, thick] ($(A14)! 0.5!(A6)$) -- ($(A15)! 0.5!(A6)$);
\draw[color1, thick] ($(A7)! 0.5!(A8)$) -- ($(A16)! 0.5!(A8)$);
\draw[color1, thick] ($(A8)! 0.5!(A7)$) -- ($(A95)! 0.5!(A7)$);
\draw[color1, thick] ($(A94)! 0.5!(A7)$) -- ($(A95)! 0.5!(A7)$);
\draw[color1, thick] ($(A9)! 0.5!(A60)$) -- ($(A17)! 0.5!(A60)$);
\draw[color1, thick] ($(A45)! 0.5!(A9)$) -- ($(A53)! 0.5!(A9)$);
\draw[color1, thick] ($(A53)! 0.5!(A9)$) -- ($(A64)! 0.5!(A9)$);
\draw[color1, thick] ($(A9)! 0.5!(A60)$) -- ($(A61)! 0.5!(A60)$);
\draw[color1, thick] ($(A9)! 0.5!(A62)$) -- ($(A61)! 0.5!(A62)$);
\draw[color1, thick] ($(A9)! 0.5!(A62)$) -- ($(A63)! 0.5!(A62)$);
\draw[color1, thick] ($(A9)! 0.5!(A64)$) -- ($(A63)! 0.5!(A64)$);
\draw[color1, thick] ($(A15)! 0.5!(A23)$) -- ($(A16)! 0.5!(A23)$);
\draw[color1, thick] ($(A15)! 0.5!(A23)$) -- ($(A22)! 0.5!(A23)$);
\draw[color1, thick] ($(A17)! 0.5!(A60)$) -- ($(A24)! 0.5!(A60)$);
\draw[color1, thick] ($(A18)! 0.5!(A25)$) -- ($(A24)! 0.5!(A25)$);
\draw[color1, thick] ($(A18)! 0.5!(A25)$) -- ($(A26)! 0.5!(A25)$);
\draw[color1, thick] ($(A20)! 0.5!(A27)$) -- ($(A21)! 0.5!(A27)$);
\draw[color1, thick] ($(A20)! 0.5!(A27)$) -- ($(A26)! 0.5!(A27)$);
\draw[color1, thick] ($(A21)! 0.5!(A27)$) -- ($(A22)! 0.5!(A27)$);
\draw[color1, thick] ($(A22)! 0.5!(A23)$) -- ($(A29)! 0.5!(A23)$);
\draw[color1, thick] ($(A22)! 0.5!(A27)$) -- ($(A35)! 0.5!(A27)$);
\draw[color1, thick] ($(A22)! 0.5!(A32)$) -- ($(A28)! 0.5!(A32)$);
\draw[color1, thick] ($(A22)! 0.5!(A32)$) -- ($(A35)! 0.5!(A32)$);
\draw[color1, thick] ($(A24)! 0.5!(A25)$) -- ($(A30)! 0.5!(A25)$);
\draw[color1, thick] ($(A24)! 0.5!(A60)$) -- ($(A66)! 0.5!(A60)$);
\draw[color1, thick] ($(A26)! 0.5!(A25)$) -- ($(A31)! 0.5!(A25)$);
\draw[color1, thick] ($(A30)! 0.5!(A25)$) -- ($(A31)! 0.5!(A25)$);
\draw[color1, thick] ($(A26)! 0.5!(A27)$) -- ($(A31)! 0.5!(A27)$);
\draw[color1, thick] ($(A31)! 0.5!(A27)$) -- ($(A35)! 0.5!(A27)$);
\draw[color1, thick] ($(A28)! 0.5!(A34)$) -- ($(A29)! 0.5!(A34)$);
\draw[color1, thick] ($(A28)! 0.5!(A32)$) -- ($(A33)! 0.5!(A32)$);
\draw[color1, thick] ($(A28)! 0.5!(A34)$) -- ($(A33)! 0.5!(A34)$);
\draw[color1, thick] ($(A30)! 0.5!(A71)$) -- ($(A35)! 0.5!(A71)$);
\draw[color1, thick] ($(A30)! 0.5!(A71)$) -- ($(A66)! 0.5!(A71)$);
\draw[color1, thick] ($(A33)! 0.5!(A32)$) -- ($(A37)! 0.5!(A32)$);
\draw[color1, thick] ($(A35)! 0.5!(A32)$) -- ($(A36)! 0.5!(A32)$);
\draw[color1, thick] ($(A36)! 0.5!(A32)$) -- ($(A37)! 0.5!(A32)$);
\draw[color1, thick] ($(A33)! 0.5!(A34)$) -- ($(A38)! 0.5!(A34)$);
\draw[color1, thick] ($(A35)! 0.5!(A39)$) -- ($(A36)! 0.5!(A39)$);
\draw[color1, thick] ($(A39)! 0.5!(A35)$) -- ($(A75)! 0.5!(A35)$);
\draw[color1, thick] ($(A35)! 0.5!(A64)$) -- ($(A68)! 0.5!(A64)$);
\draw[color1, thick] ($(A64)! 0.5!(A35)$) -- ($(A72)! 0.5!(A35)$);
\draw[color1, thick] ($(A35)! 0.5!(A71)$) -- ($(A68)! 0.5!(A71)$);
\draw[color1, thick] ($(A72)! 0.5!(A35)$) -- ($(A75)! 0.5!(A35)$);
\draw[color1, thick] ($(A36)! 0.5!(A39)$) -- ($(A40)! 0.5!(A39)$);
\draw[color1, thick] ($(A37)! 0.5!(A41)$) -- ($(A38)! 0.5!(A41)$);
\draw[color1, thick] ($(A37)! 0.5!(A41)$) -- ($(A40)! 0.5!(A41)$);
\draw[color1, thick] ($(A40)! 0.5!(A39)$) -- ($(A42)! 0.5!(A39)$);
\draw[color1, thick] ($(A39)! 0.5!(A42)$) -- ($(A78)! 0.5!(A42)$);
\draw[color1, thick] ($(A40)! 0.5!(A41)$) -- ($(A43)! 0.5!(A41)$);
\draw[color1, thick] ($(A42)! 0.5!(A44)$) -- ($(A43)! 0.5!(A44)$);
\draw[color1, thick] ($(A44)! 0.5!(A42)$) -- ($(A80)! 0.5!(A42)$);
\draw[color1, thick] ($(A78)! 0.5!(A42)$) -- ($(A80)! 0.5!(A42)$);
\draw[color1, thick] ($(A51)! 0.5!(A81)$) -- ($(A52)! 0.5!(A81)$);
\draw[color1, thick] ($(A51)! 0.5!(A81)$) -- ($(A82)! 0.5!(A81)$);
\draw[color1, thick] ($(A61)! 0.5!(A60)$) -- ($(A67)! 0.5!(A60)$);
\draw[color1, thick] ($(A66)! 0.5!(A60)$) -- ($(A67)! 0.5!(A60)$);
\draw[color1, thick] ($(A61)! 0.5!(A62)$) -- ($(A67)! 0.5!(A62)$);
\draw[color1, thick] ($(A63)! 0.5!(A62)$) -- ($(A68)! 0.5!(A62)$);
\draw[color1, thick] ($(A67)! 0.5!(A62)$) -- ($(A68)! 0.5!(A62)$);
\draw[color1, thick] ($(A63)! 0.5!(A64)$) -- ($(A68)! 0.5!(A64)$);
\draw[color1, thick] ($(A66)! 0.5!(A71)$) -- ($(A67)! 0.5!(A71)$);
\draw[color1, thick] ($(A67)! 0.5!(A71)$) -- ($(A68)! 0.5!(A71)$);
\draw[color1, thick] ($(A82)! 0.5!(A81)$) -- ($(A96)! 0.5!(A81)$);
\draw[color1, thick] ($(A96)! 0.5!(A109)$) -- ($(A97)! 0.5!(A109)$);
\draw[color1, thick] ($(A97)! 0.5!(A109)$) -- ($(A110)! 0.5!(A109)$);
\draw[color1, thick] ($(A110)! 0.5!(A109)$) -- ($(A120)! 0.5!(A109)$);
\draw[color1, thick] ($(A120)! 0.5!(A129)$) -- ($(A121)! 0.5!(A129)$);
\draw[color1, thick] ($(A121)! 0.5!(A129)$) -- ($(A130)! 0.5!(A129)$);
\draw[color1, thick] ($(A130)! 0.5!(A129)$) -- ($(A136)! 0.5!(A129)$);
\draw[color1, thick] ($(A136)! 0.5!(A141)$) -- ($(A137)! 0.5!(A141)$);
\draw[color1, thick] ($(A137)! 0.5!(A141)$) -- ($(A142)! 0.5!(A141)$);
\draw[color1, thick] ($(A142)! 0.5!(A141)$) -- ($(A144)! 0.5!(A141)$);
\draw[color11, thick] ($(A30)! 0.5!(A71)$) -- ($(A35)! 0.5!(A71)$);
\draw[color11, thick] ($(A30)! 0.5!(A71)$) -- ($(A66)! 0.5!(A71)$);
\draw[color11, thick] ($(A35)! 0.5!(A71)$) -- ($(A68)! 0.5!(A71)$);
\draw[color11, thick] ($(A66)! 0.5!(A71)$) -- ($(A67)! 0.5!(A71)$);
\draw[color11, thick] ($(A67)! 0.5!(A71)$) -- ($(A68)! 0.5!(A71)$);
\draw[color9, thick] ($(A9)! 0.5!(A62)$) -- ($(A61)! 0.5!(A62)$);
\draw[color9, thick] ($(A9)! 0.5!(A62)$) -- ($(A63)! 0.5!(A62)$);
\draw[color9, thick] ($(A61)! 0.5!(A62)$) -- ($(A67)! 0.5!(A62)$);
\draw[color9, thick] ($(A63)! 0.5!(A62)$) -- ($(A68)! 0.5!(A62)$);
\draw[color9, thick] ($(A67)! 0.5!(A62)$) -- ($(A68)! 0.5!(A62)$);
\draw[color8, thick] ($(A9)! 0.5!(A60)$) -- ($(A17)! 0.5!(A60)$);
\draw[color8, thick] ($(A9)! 0.5!(A60)$) -- ($(A61)! 0.5!(A60)$);
\draw[color8, thick] ($(A17)! 0.5!(A60)$) -- ($(A24)! 0.5!(A60)$);
\draw[color8, thick] ($(A24)! 0.5!(A60)$) -- ($(A66)! 0.5!(A60)$);
\draw[color8, thick] ($(A61)! 0.5!(A60)$) -- ($(A67)! 0.5!(A60)$);
\draw[color8, thick] ($(A66)! 0.5!(A60)$) -- ($(A67)! 0.5!(A60)$);
\draw[color4, thick] ($(A22)! 0.5!(A32)$) -- ($(A28)! 0.5!(A32)$);
\draw[color4, thick] ($(A22)! 0.5!(A32)$) -- ($(A35)! 0.5!(A32)$);
\draw[color4, thick] ($(A28)! 0.5!(A32)$) -- ($(A33)! 0.5!(A32)$);
\draw[color4, thick] ($(A33)! 0.5!(A32)$) -- ($(A37)! 0.5!(A32)$);
\draw[color4, thick] ($(A35)! 0.5!(A32)$) -- ($(A36)! 0.5!(A32)$);
\draw[color4, thick] ($(A36)! 0.5!(A32)$) -- ($(A37)! 0.5!(A32)$);
\draw[color3, thick] ($(A20)! 0.5!(A27)$) -- ($(A21)! 0.5!(A27)$);
\draw[color3, thick] ($(A20)! 0.5!(A27)$) -- ($(A26)! 0.5!(A27)$);
\draw[color3, thick] ($(A21)! 0.5!(A27)$) -- ($(A22)! 0.5!(A27)$);
\draw[color3, thick] ($(A22)! 0.5!(A27)$) -- ($(A35)! 0.5!(A27)$);
\draw[color3, thick] ($(A26)! 0.5!(A27)$) -- ($(A31)! 0.5!(A27)$);
\draw[color3, thick] ($(A31)! 0.5!(A27)$) -- ($(A35)! 0.5!(A27)$);
\draw[color2, thick] ($(A18)! 0.5!(A25)$) -- ($(A24)! 0.5!(A25)$);
\draw[color2, thick] ($(A18)! 0.5!(A25)$) -- ($(A26)! 0.5!(A25)$);
\draw[color2, thick] ($(A24)! 0.5!(A25)$) -- ($(A30)! 0.5!(A25)$);
\draw[color2, thick] ($(A26)! 0.5!(A25)$) -- ($(A31)! 0.5!(A25)$);
\draw[color2, thick] ($(A30)! 0.5!(A25)$) -- ($(A31)! 0.5!(A25)$);
\foreach \a/\b/\c in {0/1/9,0/1/88,0/9/45,0/45/88,1/2/10,1/2/89,1/9/10,1/88/89,2/3/11,2/3/90,2/10/11,2/89/90,3/4/12,3/4/91,3/11/12,3/90/91,4/5/13,4/5/92,4/12/13,4/91/92,5/6/14,5/6/93,5/13/14,5/92/93,6/7/15,6/7/94,6/14/15,6/93/94,7/8/16,7/8/95,7/15/16,7/94/95,9/10/22,9/17/18,9/17/60,9/18/19,9/19/20,9/20/21,9/21/22,9/45/53,9/53/64,9/60/61,9/61/62,9/62/63,9/63/64,10/11/22,11/12/22,12/13/22,13/14/22,14/15/22,15/16/23,15/22/23,17/18/24,17/24/60,18/19/26,18/24/25,18/25/26,19/20/26,20/21/27,20/26/27,21/22/27,22/23/29,22/27/35,22/28/29,22/28/32,22/32/35,24/25/30,24/30/66,24/60/66,25/26/31,25/30/31,26/27/31,27/31/35,28/29/34,28/32/33,28/33/34,30/31/35,30/35/71,30/66/71,32/33/37,32/35/36,32/36/37,33/34/38,33/37/38,35/36/39,35/39/75,35/64/68,35/64/72,35/68/71,35/72/75,36/37/40,36/39/40,37/38/41,37/40/41,39/40/42,39/42/78,39/75/78,40/41/43,40/42/43,42/43/44,42/44/80,42/78/80,45/46/53,45/46/87,45/87/88,46/47/54,46/47/86,46/53/54,46/86/87,47/48/55,47/48/85,47/54/55,47/85/86,48/49/56,48/49/84,48/55/56,48/84/85,49/50/57,49/50/83,49/56/57,49/83/84,50/51/58,50/51/82,50/57/58,50/82/83,51/52/59,51/52/81,51/58/59,51/81/82,53/54/64,54/55/64,55/56/64,56/57/64,57/58/64,58/59/65,58/64/65,60/61/67,60/66/67,61/62/67,62/63/68,62/67/68,63/64/68,64/65/70,64/69/70,64/69/72,66/67/71,67/68/71,69/70/74,69/72/73,69/73/74,72/73/76,72/75/76,73/74/77,73/76/77,75/76/78,76/77/79,76/78/79,78/79/80,81/82/96,82/83/97,82/96/97,83/84/97,84/85/97,85/86/97,86/87/97,87/88/97,88/89/107,88/97/98,88/98/99,88/99/100,88/100/101,88/101/102,88/102/103,88/103/104,88/104/105,88/105/106,88/106/107,89/90/107,90/91/107,91/92/107,92/93/107,93/94/107,94/95/108,94/107/108,96/97/109,97/98/111,97/109/110,97/110/122,97/111/132,97/122/132,98/99/111,99/100/112,99/111/112,100/101/112,101/102/114,101/112/113,101/113/114,102/103/114,103/104/116,103/114/115,103/115/116,104/105/116,105/106/117,105/116/117,106/107/117,107/108/119,107/117/132,107/118/119,107/118/126,107/126/132,109/110/120,110/120/121,110/121/122,111/112/123,111/123/132,112/113/123,113/114/124,113/123/124,114/115/124,115/116/125,115/124/125,116/117/125,117/125/132,118/119/128,118/126/127,118/127/128,120/121/129,121/122/130,121/129/130,122/130/131,122/131/132,123/124/132,124/125/132,126/127/134,126/132/133,126/133/134,127/128/135,127/134/135,129/130/136,130/131/137,130/136/137,131/132/138,131/137/138,132/133/138,133/134/139,133/138/139,134/135/140,134/139/140,136/137/141,137/138/142,137/141/142,138/139/142,139/140/143,139/142/143,141/142/144,142/143/144}{
  \draw[black!30] (A\a) -- (A\b) -- (A\c) -- cycle;
}
\draw[black, thick] (A22) -- (A35) -- (A64) -- (A9) --cycle;
\draw[black, thick] (A107) -- (A132) -- (A97) -- (A88) --cycle;

\foreach \i in {1,3,5,7,9,10,11,12,13,14,15,16,17,18,19,20,21,22,24,26,28,29,30,31,33,35,36,37,38,40,42,43,54,56,58,61,63,66,67,68,69,76,81,83,85,87,98,100,102,103,104,105,106,107,109,115,117,124,125,127,129,131,139,141}{
  \fill[myred] (A\i) circle (3pt);
}
\foreach \i in {0,2,4,6,8,23,25,27,32,34,39,41,44,45,46,47,48,49,50,51,52,53,55,57,59,60,62,64,65,70,71,72,73,74,75,77,78,79,80,82,84,86,88,89,90,91,92,93,94,95,96,97,99,101,108,110,111,112,113,114,116,118,119,120,121,122,123,126,128,130,132,133,134,135,136,137,138,140,142,143,144}{
  \fill[myblue] (A\i) circle (3pt);
}

%\foreach \i in {-8,-6,-4,-2,0,2,4,6,8}{
%	\foreach \j in {-8,-6,-4,-2,0,2,4,6,8}{
%		\node[] at (\i,\j) {\input{tikz/even_even.tikz}};
%	}
%	\foreach \j in {-7,-5,-3,-1,1,3,5,7}{
%		\node[] at (\i,\j) {\input{tikz/even_odd.tikz}};
%	}
%}
%\foreach \i in {-7,-5,-3,-1,1,3,5,7}{
%	\foreach \j in {-8,-6,-4,-2,0,2,4,6,8}{
%		\node[rotate=90] at (\i,\j) {\input{tikz/even_odd.tikz}};
%	}
%	\foreach \j in {-7,-5,-3,-1,1,3,5,7}{
%		\node[] at (\i,\j) {\input{tikz/odd_odd.tikz}};
%	}
%}

%% file: tikz/double_split_4_narrow.tikz
\coordinate (A0) at (0, 0);
\coordinate (A1) at (0, 1);
\coordinate (A2) at (0, 2);
\coordinate (A3) at (0, 3);
\coordinate (A4) at (0, 4);
\coordinate (A5) at (0, 5);
\coordinate (A6) at (0, 6);
\coordinate (A7) at (0, 7);
\coordinate (A8) at (0, 8);
\coordinate (A9) at (1, 0);
\coordinate (A10) at (1, 1);
\coordinate (A11) at (1, 2);
\coordinate (A12) at (1, 3);
\coordinate (A13) at (1, 4);
\coordinate (A14) at (1, 5);
\coordinate (A15) at (1, 6);
\coordinate (A16) at (1, 7);
\coordinate (A17) at (2, 0);
\coordinate (A18) at (2, 1);
\coordinate (A19) at (2, 2);
\coordinate (A20) at (2, 3);
\coordinate (A21) at (2, 4);
\coordinate (A22) at (2, 5);
\coordinate (A23) at (2, 6);
\coordinate (A24) at (3, 0);
\coordinate (A25) at (3, 1);
\coordinate (A26) at (3, 2);
\coordinate (A27) at (3, 3);
\coordinate (A28) at (3, 4);
\coordinate (A29) at (3, 5);
\coordinate (A30) at (4, 0);
\coordinate (A31) at (4, 1);
\coordinate (A32) at (4, 2);
\coordinate (A33) at (4, 3);
\coordinate (A34) at (4, 4);
\coordinate (A35) at (5, 0);
\coordinate (A36) at (5, 1);
\coordinate (A37) at (5, 2);
\coordinate (A38) at (5, 3);
\coordinate (A39) at (6, 0);
\coordinate (A40) at (6, 1);
\coordinate (A41) at (6, 2);
\coordinate (A42) at (7, 0);
\coordinate (A43) at (7, 1);
\coordinate (A44) at (8, 0);
\coordinate (A45) at (0, -1);
\coordinate (A46) at (0, -2);
\coordinate (A47) at (0, -3);
\coordinate (A48) at (0, -4);
\coordinate (A49) at (0, -5);
\coordinate (A50) at (0, -6);
\coordinate (A51) at (0, -7);
\coordinate (A52) at (0, -8);
\coordinate (A53) at (1, -1);
\coordinate (A54) at (1, -2);
\coordinate (A55) at (1, -3);
\coordinate (A56) at (1, -4);
\coordinate (A57) at (1, -5);
\coordinate (A58) at (1, -6);
\coordinate (A59) at (1, -7);
\coordinate (A60) at (2, -1);
\coordinate (A61) at (2, -2);
\coordinate (A62) at (2, -3);
\coordinate (A63) at (2, -4);
\coordinate (A64) at (2, -5);
\coordinate (A65) at (2, -6);
\coordinate (A66) at (3, -1);
\coordinate (A67) at (3, -2);
\coordinate (A68) at (3, -3);
\coordinate (A69) at (3, -4);
\coordinate (A70) at (3, -5);
\coordinate (A71) at (4, -1);
\coordinate (A72) at (4, -2);
\coordinate (A73) at (4, -3);
\coordinate (A74) at (4, -4);
\coordinate (A75) at (5, -1);
\coordinate (A76) at (5, -2);
\coordinate (A77) at (5, -3);
\coordinate (A78) at (6, -1);
\coordinate (A79) at (6, -2);
\coordinate (A80) at (7, -1);
\coordinate (A81) at (-1, -7);
\coordinate (A82) at (-1, -6);
\coordinate (A83) at (-1, -5);
\coordinate (A84) at (-1, -4);
\coordinate (A85) at (-1, -3);
\coordinate (A86) at (-1, -2);
\coordinate (A87) at (-1, -1);
\coordinate (A88) at (-1, 0);
\coordinate (A89) at (-1, 1);
\coordinate (A90) at (-1, 2);
\coordinate (A91) at (-1, 3);
\coordinate (A92) at (-1, 4);
\coordinate (A93) at (-1, 5);
\coordinate (A94) at (-1, 6);
\coordinate (A95) at (-1, 7);
\coordinate (A96) at (-2, -6);
\coordinate (A97) at (-2, -5);
\coordinate (A98) at (-2, -4);
\coordinate (A99) at (-2, -3);
\coordinate (A100) at (-2, -2);
\coordinate (A101) at (-2, -1);
\coordinate (A102) at (-2, 0);
\coordinate (A103) at (-2, 1);
\coordinate (A104) at (-2, 2);
\coordinate (A105) at (-2, 3);
\coordinate (A106) at (-2, 4);
\coordinate (A107) at (-2, 5);
\coordinate (A108) at (-2, 6);
\coordinate (A109) at (-3, -5);
\coordinate (A110) at (-3, -4);
\coordinate (A111) at (-3, -3);
\coordinate (A112) at (-3, -2);
\coordinate (A113) at (-3, -1);
\coordinate (A114) at (-3, 0);
\coordinate (A115) at (-3, 1);
\coordinate (A116) at (-3, 2);
\coordinate (A117) at (-3, 3);
\coordinate (A118) at (-3, 4);
\coordinate (A119) at (-3, 5);
\coordinate (A120) at (-4, -4);
\coordinate (A121) at (-4, -3);
\coordinate (A122) at (-4, -2);
\coordinate (A123) at (-4, -1);
\coordinate (A124) at (-4, 0);
\coordinate (A125) at (-4, 1);
\coordinate (A126) at (-4, 2);
\coordinate (A127) at (-4, 3);
\coordinate (A128) at (-4, 4);
\coordinate (A129) at (-5, -3);
\coordinate (A130) at (-5, -2);
\coordinate (A131) at (-5, -1);
\coordinate (A132) at (-5, 0);
\coordinate (A133) at (-5, 1);
\coordinate (A134) at (-5, 2);
\coordinate (A135) at (-5, 3);
\coordinate (A136) at (-6, -2);
\coordinate (A137) at (-6, -1);
\coordinate (A138) at (-6, 0);
\coordinate (A139) at (-6, 1);
\coordinate (A140) at (-6, 2);
\coordinate (A141) at (-7, -1);
\coordinate (A142) at (-7, 0);
\coordinate (A143) at (-7, 1);
\coordinate (A144) at (-8, 0);
\colorlet{color1}{mycolor4}
\fill[color1!50] (A1) -- ($(A1)! 0.5!(A0)$) -- ($(A9)! 0.5!(A0)$) -- (A9) -- cycle;
\fill[color1!50] ($(A0)! 0.5!(A1)$) -- ($(A88)! 0.5!(A1)$) -- (A1) -- cycle;
\fill[color1!50] ($(A0)! 0.5!(A9)$) -- ($(A45)! 0.5!(A9)$) -- (A9) -- cycle;
\fill[color1!50] (A1) -- ($(A1)! 0.5!(A2)$) -- ($(A10)! 0.5!(A2)$) -- (A10) -- cycle;
\fill[color1!50] ($(A2)! 0.5!(A1)$) -- ($(A89)! 0.5!(A1)$) -- (A1) -- cycle;
\fill[color1!50] (A1) -- (A9) -- (A10) -- cycle;
\fill[color1!50] ($(A88)! 0.5!(A1)$) -- ($(A89)! 0.5!(A1)$) -- (A1) -- cycle;
\fill[color1!50] (A3) -- ($(A3)! 0.5!(A2)$) -- ($(A11)! 0.5!(A2)$) -- (A11) -- cycle;
\fill[color1!50] ($(A2)! 0.5!(A3)$) -- ($(A90)! 0.5!(A3)$) -- (A3) -- cycle;
\fill[color1!50] (A10) -- ($(A10)! 0.5!(A2)$) -- ($(A11)! 0.5!(A2)$) -- (A11) -- cycle;
\fill[color1!50] (A3) -- ($(A3)! 0.5!(A4)$) -- ($(A12)! 0.5!(A4)$) -- (A12) -- cycle;
\fill[color1!50] ($(A4)! 0.5!(A3)$) -- ($(A91)! 0.5!(A3)$) -- (A3) -- cycle;
\fill[color1!50] (A3) -- (A11) -- (A12) -- cycle;
\fill[color1!50] ($(A90)! 0.5!(A3)$) -- ($(A91)! 0.5!(A3)$) -- (A3) -- cycle;
\fill[color1!50] (A5) -- ($(A5)! 0.5!(A4)$) -- ($(A13)! 0.5!(A4)$) -- (A13) -- cycle;
\fill[color1!50] ($(A4)! 0.5!(A5)$) -- ($(A92)! 0.5!(A5)$) -- (A5) -- cycle;
\fill[color1!50] (A12) -- ($(A12)! 0.5!(A4)$) -- ($(A13)! 0.5!(A4)$) -- (A13) -- cycle;
\fill[color1!50] (A5) -- ($(A5)! 0.5!(A6)$) -- ($(A14)! 0.5!(A6)$) -- (A14) -- cycle;
\fill[color1!50] ($(A6)! 0.5!(A5)$) -- ($(A93)! 0.5!(A5)$) -- (A5) -- cycle;
\fill[color1!50] (A5) -- (A13) -- (A14) -- cycle;
\fill[color1!50] ($(A92)! 0.5!(A5)$) -- ($(A93)! 0.5!(A5)$) -- (A5) -- cycle;
\fill[color1!50] (A7) -- ($(A7)! 0.5!(A6)$) -- ($(A15)! 0.5!(A6)$) -- (A15) -- cycle;
\fill[color1!50] ($(A6)! 0.5!(A7)$) -- ($(A94)! 0.5!(A7)$) -- (A7) -- cycle;
\fill[color1!50] (A14) -- ($(A14)! 0.5!(A6)$) -- ($(A15)! 0.5!(A6)$) -- (A15) -- cycle;
\fill[color1!50] (A7) -- ($(A7)! 0.5!(A8)$) -- ($(A16)! 0.5!(A8)$) -- (A16) -- cycle;
\fill[color1!50] ($(A8)! 0.5!(A7)$) -- ($(A95)! 0.5!(A7)$) -- (A7) -- cycle;
\fill[color1!50] (A7) -- (A15) -- (A16) -- cycle;
\fill[color1!50] ($(A94)! 0.5!(A7)$) -- ($(A95)! 0.5!(A7)$) -- (A7) -- cycle;
\fill[color1!50] (A9) -- ($(A9)! 0.5!(A17)$) -- ($(A10)! 0.5!(A17)$) -- (A10) -- cycle;
\fill[color1!50] ($(A17)! 0.5!(A9)$) -- ($(A53)! 0.5!(A9)$) -- (A9) -- cycle;
\fill[color1!50] ($(A45)! 0.5!(A9)$) -- ($(A53)! 0.5!(A9)$) -- (A9) -- cycle;
\fill[color1!50] (A10) -- (A11) -- (A18) -- cycle;
\fill[color1!50] (A10) -- ($(A10)! 0.5!(A17)$) -- ($(A18)! 0.5!(A17)$) -- (A18) -- cycle;
\fill[color1!50] (A11) -- ($(A11)! 0.5!(A19)$) -- ($(A12)! 0.5!(A19)$) -- (A12) -- cycle;
\fill[color1!50] (A11) -- ($(A11)! 0.5!(A19)$) -- ($(A18)! 0.5!(A19)$) -- (A18) -- cycle;
\fill[color1!50] (A12) -- (A13) -- (A20) -- cycle;
\fill[color1!50] (A12) -- ($(A12)! 0.5!(A19)$) -- ($(A20)! 0.5!(A19)$) -- (A20) -- cycle;
\fill[color1!50] (A13) -- ($(A13)! 0.5!(A21)$) -- ($(A14)! 0.5!(A21)$) -- (A14) -- cycle;
\fill[color1!50] (A13) -- ($(A13)! 0.5!(A21)$) -- ($(A20)! 0.5!(A21)$) -- (A20) -- cycle;
\fill[color1!50] (A14) -- (A15) -- (A22) -- cycle;
\fill[color1!50] (A14) -- ($(A14)! 0.5!(A21)$) -- ($(A22)! 0.5!(A21)$) -- (A22) -- cycle;
\fill[color1!50] (A15) -- ($(A15)! 0.5!(A23)$) -- ($(A16)! 0.5!(A23)$) -- (A16) -- cycle;
\fill[color1!50] (A15) -- ($(A15)! 0.5!(A23)$) -- ($(A22)! 0.5!(A23)$) -- (A22) -- cycle;
\fill[color1!50] (A18) -- ($(A18)! 0.5!(A17)$) -- ($(A24)! 0.5!(A17)$) -- (A24) -- cycle;
\fill[color1!50] ($(A17)! 0.5!(A24)$) -- ($(A60)! 0.5!(A24)$) -- (A24) -- cycle;
\fill[color1!50] (A18) -- ($(A18)! 0.5!(A19)$) -- ($(A24)! 0.5!(A19)$) -- (A24) -- cycle;
\fill[color1!50] (A20) -- ($(A20)! 0.5!(A19)$) -- ($(A24)! 0.5!(A19)$) -- (A24) -- cycle;
\fill[color1!50] (A20) -- ($(A20)! 0.5!(A21)$) -- ($(A24)! 0.5!(A21)$) -- (A24) -- cycle;
\fill[color1!50] (A22) -- ($(A22)! 0.5!(A21)$) -- ($(A24)! 0.5!(A21)$) -- (A24) -- cycle;
\fill[color1!50] (A22) -- ($(A22)! 0.5!(A23)$) -- ($(A29)! 0.5!(A23)$) -- (A29) -- cycle;
\fill[color1!50] (A22) -- ($(A22)! 0.5!(A25)$) -- ($(A24)! 0.5!(A25)$) -- (A24) -- cycle;
\fill[color1!50] (A22) -- ($(A22)! 0.5!(A25)$) -- ($(A26)! 0.5!(A25)$) -- (A26) -- cycle;
\fill[color1!50] (A22) -- ($(A22)! 0.5!(A27)$) -- ($(A26)! 0.5!(A27)$) -- (A26) -- cycle;
\fill[color1!50] (A22) -- ($(A22)! 0.5!(A27)$) -- ($(A35)! 0.5!(A27)$) -- (A35) -- cycle;
\fill[color1!50] (A22) -- (A28) -- (A29) -- cycle;
\fill[color1!50] (A22) -- ($(A22)! 0.5!(A32)$) -- ($(A28)! 0.5!(A32)$) -- (A28) -- cycle;
\fill[color1!50] (A22) -- ($(A22)! 0.5!(A32)$) -- ($(A35)! 0.5!(A32)$) -- (A35) -- cycle;
\fill[color1!50] (A24) -- ($(A24)! 0.5!(A25)$) -- ($(A30)! 0.5!(A25)$) -- (A30) -- cycle;
\fill[color1!50] (A24) -- (A30) -- (A66) -- cycle;
\fill[color1!50] ($(A60)! 0.5!(A24)$) -- ($(A61)! 0.5!(A24)$) -- (A24) -- cycle;
\fill[color1!50] ($(A61)! 0.5!(A24)$) -- ($(A62)! 0.5!(A24)$) -- (A24) -- cycle;
\fill[color1!50] ($(A62)! 0.5!(A24)$) -- ($(A63)! 0.5!(A24)$) -- (A24) -- cycle;
\fill[color1!50] ($(A63)! 0.5!(A24)$) -- ($(A64)! 0.5!(A24)$) -- (A24) -- cycle;
\fill[color1!50] (A24) -- ($(A24)! 0.5!(A64)$) -- ($(A66)! 0.5!(A64)$) -- (A66) -- cycle;
\fill[color1!50] (A26) -- ($(A26)! 0.5!(A25)$) -- ($(A31)! 0.5!(A25)$) -- (A31) -- cycle;
\fill[color1!50] (A30) -- ($(A30)! 0.5!(A25)$) -- ($(A31)! 0.5!(A25)$) -- (A31) -- cycle;
\fill[color1!50] (A26) -- ($(A26)! 0.5!(A27)$) -- ($(A31)! 0.5!(A27)$) -- (A31) -- cycle;
\fill[color1!50] (A31) -- ($(A31)! 0.5!(A27)$) -- ($(A35)! 0.5!(A27)$) -- (A35) -- cycle;
\fill[color1!50] (A28) -- ($(A28)! 0.5!(A34)$) -- ($(A29)! 0.5!(A34)$) -- (A29) -- cycle;
\fill[color1!50] (A28) -- ($(A28)! 0.5!(A32)$) -- ($(A33)! 0.5!(A32)$) -- (A33) -- cycle;
\fill[color1!50] (A28) -- ($(A28)! 0.5!(A34)$) -- ($(A33)! 0.5!(A34)$) -- (A33) -- cycle;
\fill[color1!50] (A30) -- (A31) -- (A35) -- cycle;
\fill[color1!50] (A30) -- ($(A30)! 0.5!(A71)$) -- ($(A35)! 0.5!(A71)$) -- (A35) -- cycle;
\fill[color1!50] (A30) -- ($(A30)! 0.5!(A71)$) -- ($(A66)! 0.5!(A71)$) -- (A66) -- cycle;
\fill[color1!50] (A33) -- ($(A33)! 0.5!(A32)$) -- ($(A37)! 0.5!(A32)$) -- (A37) -- cycle;
\fill[color1!50] (A35) -- ($(A35)! 0.5!(A32)$) -- ($(A36)! 0.5!(A32)$) -- (A36) -- cycle;
\fill[color1!50] (A36) -- ($(A36)! 0.5!(A32)$) -- ($(A37)! 0.5!(A32)$) -- (A37) -- cycle;
\fill[color1!50] (A33) -- ($(A33)! 0.5!(A34)$) -- ($(A38)! 0.5!(A34)$) -- (A38) -- cycle;
\fill[color1!50] (A33) -- (A37) -- (A38) -- cycle;
\fill[color1!50] (A35) -- ($(A35)! 0.5!(A39)$) -- ($(A36)! 0.5!(A39)$) -- (A36) -- cycle;
\fill[color1!50] ($(A39)! 0.5!(A35)$) -- ($(A75)! 0.5!(A35)$) -- (A35) -- cycle;
\fill[color1!50] (A35) -- ($(A35)! 0.5!(A64)$) -- ($(A68)! 0.5!(A64)$) -- (A68) -- cycle;
\fill[color1!50] ($(A64)! 0.5!(A35)$) -- ($(A72)! 0.5!(A35)$) -- (A35) -- cycle;
\fill[color1!50] (A35) -- ($(A35)! 0.5!(A71)$) -- ($(A68)! 0.5!(A71)$) -- (A68) -- cycle;
\fill[color1!50] ($(A72)! 0.5!(A35)$) -- ($(A75)! 0.5!(A35)$) -- (A35) -- cycle;
\fill[color1!50] (A36) -- (A37) -- (A40) -- cycle;
\fill[color1!50] (A36) -- ($(A36)! 0.5!(A39)$) -- ($(A40)! 0.5!(A39)$) -- (A40) -- cycle;
\fill[color1!50] (A37) -- ($(A37)! 0.5!(A41)$) -- ($(A38)! 0.5!(A41)$) -- (A38) -- cycle;
\fill[color1!50] (A37) -- ($(A37)! 0.5!(A41)$) -- ($(A40)! 0.5!(A41)$) -- (A40) -- cycle;
\fill[color1!50] (A40) -- ($(A40)! 0.5!(A39)$) -- ($(A42)! 0.5!(A39)$) -- (A42) -- cycle;
\fill[color1!50] ($(A39)! 0.5!(A42)$) -- ($(A78)! 0.5!(A42)$) -- (A42) -- cycle;
\fill[color1!50] (A40) -- ($(A40)! 0.5!(A41)$) -- ($(A43)! 0.5!(A41)$) -- (A43) -- cycle;
\fill[color1!50] (A40) -- (A42) -- (A43) -- cycle;
\fill[color1!50] (A42) -- ($(A42)! 0.5!(A44)$) -- ($(A43)! 0.5!(A44)$) -- (A43) -- cycle;
\fill[color1!50] ($(A44)! 0.5!(A42)$) -- ($(A80)! 0.5!(A42)$) -- (A42) -- cycle;
\fill[color1!50] ($(A78)! 0.5!(A42)$) -- ($(A80)! 0.5!(A42)$) -- (A42) -- cycle;
\fill[color1!50] ($(A51)! 0.5!(A81)$) -- ($(A52)! 0.5!(A81)$) -- (A81) -- cycle;
\fill[color1!50] ($(A51)! 0.5!(A81)$) -- ($(A82)! 0.5!(A81)$) -- (A81) -- cycle;
\fill[color1!50] (A66) -- ($(A66)! 0.5!(A64)$) -- ($(A67)! 0.5!(A64)$) -- (A67) -- cycle;
\fill[color1!50] (A67) -- ($(A67)! 0.5!(A64)$) -- ($(A68)! 0.5!(A64)$) -- (A68) -- cycle;
\fill[color1!50] (A66) -- ($(A66)! 0.5!(A71)$) -- ($(A67)! 0.5!(A71)$) -- (A67) -- cycle;
\fill[color1!50] (A67) -- ($(A67)! 0.5!(A71)$) -- ($(A68)! 0.5!(A71)$) -- (A68) -- cycle;
\fill[color1!50] ($(A82)! 0.5!(A81)$) -- ($(A96)! 0.5!(A81)$) -- (A81) -- cycle;
\fill[color1!50] ($(A96)! 0.5!(A109)$) -- ($(A97)! 0.5!(A109)$) -- (A109) -- cycle;
\fill[color1!50] ($(A97)! 0.5!(A109)$) -- ($(A110)! 0.5!(A109)$) -- (A109) -- cycle;
\fill[color1!50] ($(A110)! 0.5!(A109)$) -- ($(A120)! 0.5!(A109)$) -- (A109) -- cycle;
\fill[color1!50] ($(A120)! 0.5!(A129)$) -- ($(A121)! 0.5!(A129)$) -- (A129) -- cycle;
\fill[color1!50] ($(A121)! 0.5!(A129)$) -- ($(A130)! 0.5!(A129)$) -- (A129) -- cycle;
\fill[color1!50] ($(A130)! 0.5!(A129)$) -- ($(A136)! 0.5!(A129)$) -- (A129) -- cycle;
\fill[color1!50] ($(A136)! 0.5!(A141)$) -- ($(A137)! 0.5!(A141)$) -- (A141) -- cycle;
\fill[color1!50] ($(A137)! 0.5!(A141)$) -- ($(A142)! 0.5!(A141)$) -- (A141) -- cycle;
\fill[color1!50] ($(A142)! 0.5!(A141)$) -- ($(A144)! 0.5!(A141)$) -- (A141) -- cycle;
\colorlet{color2}{mycolor3}
\fill[color2!50] ($(A11)! 0.5!(A19)$) -- ($(A12)! 0.5!(A19)$) -- (A19) -- cycle;
\fill[color2!50] ($(A11)! 0.5!(A19)$) -- ($(A18)! 0.5!(A19)$) -- (A19) -- cycle;
\fill[color2!50] ($(A12)! 0.5!(A19)$) -- ($(A20)! 0.5!(A19)$) -- (A19) -- cycle;
\fill[color2!50] ($(A18)! 0.5!(A19)$) -- ($(A24)! 0.5!(A19)$) -- (A19) -- cycle;
\fill[color2!50] ($(A20)! 0.5!(A19)$) -- ($(A24)! 0.5!(A19)$) -- (A19) -- cycle;
\colorlet{color3}{mycolor3}
\fill[color3!50] ($(A13)! 0.5!(A21)$) -- ($(A14)! 0.5!(A21)$) -- (A21) -- cycle;
\fill[color3!50] ($(A13)! 0.5!(A21)$) -- ($(A20)! 0.5!(A21)$) -- (A21) -- cycle;
\fill[color3!50] ($(A14)! 0.5!(A21)$) -- ($(A22)! 0.5!(A21)$) -- (A21) -- cycle;
\fill[color3!50] ($(A20)! 0.5!(A21)$) -- ($(A24)! 0.5!(A21)$) -- (A21) -- cycle;
\fill[color3!50] ($(A22)! 0.5!(A21)$) -- ($(A24)! 0.5!(A21)$) -- (A21) -- cycle;
\colorlet{color4}{mycolor5}
\fill[color4!50] ($(A22)! 0.5!(A25)$) -- ($(A24)! 0.5!(A25)$) -- (A25) -- cycle;
\fill[color4!50] ($(A22)! 0.5!(A25)$) -- ($(A26)! 0.5!(A25)$) -- (A25) -- cycle;
\fill[color4!50] ($(A24)! 0.5!(A25)$) -- ($(A30)! 0.5!(A25)$) -- (A25) -- cycle;
\fill[color4!50] ($(A26)! 0.5!(A25)$) -- ($(A31)! 0.5!(A25)$) -- (A25) -- cycle;
\fill[color4!50] ($(A30)! 0.5!(A25)$) -- ($(A31)! 0.5!(A25)$) -- (A25) -- cycle;
\colorlet{color5}{mycolor5}
\fill[color5!50] ($(A22)! 0.5!(A27)$) -- ($(A26)! 0.5!(A27)$) -- (A27) -- cycle;
\fill[color5!50] ($(A22)! 0.5!(A27)$) -- ($(A35)! 0.5!(A27)$) -- (A27) -- cycle;
\fill[color5!50] ($(A26)! 0.5!(A27)$) -- ($(A31)! 0.5!(A27)$) -- (A27) -- cycle;
\fill[color5!50] ($(A31)! 0.5!(A27)$) -- ($(A35)! 0.5!(A27)$) -- (A27) -- cycle;
\colorlet{color6}{mycolor3}
\fill[color6!50] ($(A22)! 0.5!(A32)$) -- ($(A28)! 0.5!(A32)$) -- (A32) -- cycle;
\fill[color6!50] ($(A22)! 0.5!(A32)$) -- ($(A35)! 0.5!(A32)$) -- (A32) -- cycle;
\fill[color6!50] ($(A28)! 0.5!(A32)$) -- ($(A33)! 0.5!(A32)$) -- (A32) -- cycle;
\fill[color6!50] ($(A33)! 0.5!(A32)$) -- ($(A37)! 0.5!(A32)$) -- (A32) -- cycle;
\fill[color6!50] ($(A35)! 0.5!(A32)$) -- ($(A36)! 0.5!(A32)$) -- (A32) -- cycle;
\fill[color6!50] ($(A36)! 0.5!(A32)$) -- ($(A37)! 0.5!(A32)$) -- (A32) -- cycle;
\colorlet{color7}{mycolor3}
\fill[color7!50] ($(A46)! 0.5!(A54)$) -- ($(A47)! 0.5!(A54)$) -- (A54) -- cycle;
\fill[color7!50] ($(A46)! 0.5!(A54)$) -- ($(A53)! 0.5!(A54)$) -- (A54) -- cycle;
\fill[color7!50] ($(A47)! 0.5!(A54)$) -- ($(A55)! 0.5!(A54)$) -- (A54) -- cycle;
\fill[color7!50] ($(A53)! 0.5!(A54)$) -- ($(A60)! 0.5!(A54)$) -- (A54) -- cycle;
\fill[color7!50] ($(A55)! 0.5!(A54)$) -- ($(A61)! 0.5!(A54)$) -- (A54) -- cycle;
\fill[color7!50] ($(A60)! 0.5!(A54)$) -- ($(A61)! 0.5!(A54)$) -- (A54) -- cycle;
\colorlet{color8}{mycolor3}
\fill[color8!50] ($(A48)! 0.5!(A56)$) -- ($(A49)! 0.5!(A56)$) -- (A56) -- cycle;
\fill[color8!50] ($(A48)! 0.5!(A56)$) -- ($(A55)! 0.5!(A56)$) -- (A56) -- cycle;
\fill[color8!50] ($(A49)! 0.5!(A56)$) -- ($(A57)! 0.5!(A56)$) -- (A56) -- cycle;
\fill[color8!50] ($(A55)! 0.5!(A56)$) -- ($(A62)! 0.5!(A56)$) -- (A56) -- cycle;
\fill[color8!50] ($(A57)! 0.5!(A56)$) -- ($(A63)! 0.5!(A56)$) -- (A56) -- cycle;
\fill[color8!50] ($(A62)! 0.5!(A56)$) -- ($(A63)! 0.5!(A56)$) -- (A56) -- cycle;
\colorlet{color9}{mycolor3}
\fill[color9!50] ($(A50)! 0.5!(A58)$) -- ($(A51)! 0.5!(A58)$) -- (A58) -- cycle;
\fill[color9!50] ($(A50)! 0.5!(A58)$) -- ($(A57)! 0.5!(A58)$) -- (A58) -- cycle;
\fill[color9!50] ($(A51)! 0.5!(A58)$) -- ($(A59)! 0.5!(A58)$) -- (A58) -- cycle;
\fill[color9!50] ($(A57)! 0.5!(A58)$) -- ($(A64)! 0.5!(A58)$) -- (A58) -- cycle;
\fill[color9!50] ($(A59)! 0.5!(A58)$) -- ($(A65)! 0.5!(A58)$) -- (A58) -- cycle;
\fill[color9!50] ($(A64)! 0.5!(A58)$) -- ($(A65)! 0.5!(A58)$) -- (A58) -- cycle;
\colorlet{color10}{mycolor3}
\fill[color10!50] ($(A64)! 0.5!(A69)$) -- ($(A70)! 0.5!(A69)$) -- (A69) -- cycle;
\fill[color10!50] ($(A64)! 0.5!(A69)$) -- ($(A72)! 0.5!(A69)$) -- (A69) -- cycle;
\fill[color10!50] ($(A70)! 0.5!(A69)$) -- ($(A74)! 0.5!(A69)$) -- (A69) -- cycle;
\fill[color10!50] ($(A72)! 0.5!(A69)$) -- ($(A73)! 0.5!(A69)$) -- (A69) -- cycle;
\fill[color10!50] ($(A73)! 0.5!(A69)$) -- ($(A74)! 0.5!(A69)$) -- (A69) -- cycle;
\colorlet{color11}{mycolor5}
\fill[color11!50] ($(A30)! 0.5!(A71)$) -- ($(A35)! 0.5!(A71)$) -- (A71) -- cycle;
\fill[color11!50] ($(A30)! 0.5!(A71)$) -- ($(A66)! 0.5!(A71)$) -- (A71) -- cycle;
\fill[color11!50] ($(A35)! 0.5!(A71)$) -- ($(A68)! 0.5!(A71)$) -- (A71) -- cycle;
\fill[color11!50] ($(A66)! 0.5!(A71)$) -- ($(A67)! 0.5!(A71)$) -- (A71) -- cycle;
\fill[color11!50] ($(A67)! 0.5!(A71)$) -- ($(A68)! 0.5!(A71)$) -- (A71) -- cycle;
\colorlet{color12}{mycolor3}
\fill[color12!50] ($(A72)! 0.5!(A76)$) -- ($(A73)! 0.5!(A76)$) -- (A76) -- cycle;
\fill[color12!50] ($(A72)! 0.5!(A76)$) -- ($(A75)! 0.5!(A76)$) -- (A76) -- cycle;
\fill[color12!50] ($(A73)! 0.5!(A76)$) -- ($(A77)! 0.5!(A76)$) -- (A76) -- cycle;
\fill[color12!50] ($(A75)! 0.5!(A76)$) -- ($(A78)! 0.5!(A76)$) -- (A76) -- cycle;
\fill[color12!50] ($(A77)! 0.5!(A76)$) -- ($(A79)! 0.5!(A76)$) -- (A76) -- cycle;
\fill[color12!50] ($(A78)! 0.5!(A76)$) -- ($(A79)! 0.5!(A76)$) -- (A76) -- cycle;
\colorlet{color13}{mycolor3}
\fill[color13!50] ($(A49)! 0.5!(A83)$) -- ($(A50)! 0.5!(A83)$) -- (A83) -- cycle;
\fill[color13!50] ($(A49)! 0.5!(A83)$) -- ($(A84)! 0.5!(A83)$) -- (A83) -- cycle;
\fill[color13!50] ($(A50)! 0.5!(A83)$) -- ($(A82)! 0.5!(A83)$) -- (A83) -- cycle;
\fill[color13!50] ($(A82)! 0.5!(A83)$) -- ($(A97)! 0.5!(A83)$) -- (A83) -- cycle;
\fill[color13!50] ($(A84)! 0.5!(A83)$) -- ($(A98)! 0.5!(A83)$) -- (A83) -- cycle;
\fill[color13!50] ($(A97)! 0.5!(A83)$) -- ($(A98)! 0.5!(A83)$) -- (A83) -- cycle;
\colorlet{color14}{mycolor3}
\fill[color14!50] ($(A47)! 0.5!(A85)$) -- ($(A48)! 0.5!(A85)$) -- (A85) -- cycle;
\fill[color14!50] ($(A47)! 0.5!(A85)$) -- ($(A86)! 0.5!(A85)$) -- (A85) -- cycle;
\fill[color14!50] ($(A48)! 0.5!(A85)$) -- ($(A84)! 0.5!(A85)$) -- (A85) -- cycle;
\fill[color14!50] ($(A84)! 0.5!(A85)$) -- ($(A99)! 0.5!(A85)$) -- (A85) -- cycle;
\fill[color14!50] ($(A86)! 0.5!(A85)$) -- ($(A100)! 0.5!(A85)$) -- (A85) -- cycle;
\fill[color14!50] ($(A99)! 0.5!(A85)$) -- ($(A100)! 0.5!(A85)$) -- (A85) -- cycle;
\colorlet{color15}{mycolor3}
\fill[color15!50] ($(A45)! 0.5!(A87)$) -- ($(A46)! 0.5!(A87)$) -- (A87) -- cycle;
\fill[color15!50] ($(A45)! 0.5!(A87)$) -- ($(A88)! 0.5!(A87)$) -- (A87) -- cycle;
\fill[color15!50] ($(A46)! 0.5!(A87)$) -- ($(A86)! 0.5!(A87)$) -- (A87) -- cycle;
\fill[color15!50] ($(A86)! 0.5!(A87)$) -- ($(A101)! 0.5!(A87)$) -- (A87) -- cycle;
\fill[color15!50] ($(A88)! 0.5!(A87)$) -- ($(A102)! 0.5!(A87)$) -- (A87) -- cycle;
\fill[color15!50] ($(A101)! 0.5!(A87)$) -- ($(A102)! 0.5!(A87)$) -- (A87) -- cycle;
\colorlet{color16}{mycolor3}
\fill[color16!50] ($(A89)! 0.5!(A103)$) -- ($(A90)! 0.5!(A103)$) -- (A103) -- cycle;
\fill[color16!50] ($(A89)! 0.5!(A103)$) -- ($(A102)! 0.5!(A103)$) -- (A103) -- cycle;
\fill[color16!50] ($(A90)! 0.5!(A103)$) -- ($(A104)! 0.5!(A103)$) -- (A103) -- cycle;
\fill[color16!50] ($(A102)! 0.5!(A103)$) -- ($(A114)! 0.5!(A103)$) -- (A103) -- cycle;
\fill[color16!50] ($(A104)! 0.5!(A103)$) -- ($(A114)! 0.5!(A103)$) -- (A103) -- cycle;
\colorlet{color17}{mycolor3}
\fill[color17!50] ($(A91)! 0.5!(A105)$) -- ($(A92)! 0.5!(A105)$) -- (A105) -- cycle;
\fill[color17!50] ($(A91)! 0.5!(A105)$) -- ($(A104)! 0.5!(A105)$) -- (A105) -- cycle;
\fill[color17!50] ($(A92)! 0.5!(A105)$) -- ($(A106)! 0.5!(A105)$) -- (A105) -- cycle;
\fill[color17!50] ($(A104)! 0.5!(A105)$) -- ($(A114)! 0.5!(A105)$) -- (A105) -- cycle;
\fill[color17!50] ($(A106)! 0.5!(A105)$) -- ($(A114)! 0.5!(A105)$) -- (A105) -- cycle;
\colorlet{color18}{mycolor3}
\fill[color18!50] ($(A93)! 0.5!(A107)$) -- ($(A94)! 0.5!(A107)$) -- (A107) -- cycle;
\fill[color18!50] ($(A93)! 0.5!(A107)$) -- ($(A106)! 0.5!(A107)$) -- (A107) -- cycle;
\fill[color18!50] ($(A94)! 0.5!(A107)$) -- ($(A108)! 0.5!(A107)$) -- (A107) -- cycle;
\fill[color18!50] ($(A106)! 0.5!(A107)$) -- ($(A114)! 0.5!(A107)$) -- (A107) -- cycle;
\fill[color18!50] ($(A108)! 0.5!(A107)$) -- ($(A119)! 0.5!(A107)$) -- (A107) -- cycle;
\fill[color18!50] (A107) -- ($(A107)! 0.5!(A114)$) -- ($(A115)! 0.5!(A114)$) -- (A115) -- cycle;
\fill[color18!50] (A107) -- ($(A107)! 0.5!(A116)$) -- ($(A115)! 0.5!(A116)$) -- (A115) -- cycle;
\fill[color18!50] (A107) -- ($(A107)! 0.5!(A116)$) -- ($(A117)! 0.5!(A116)$) -- (A117) -- cycle;
\fill[color18!50] (A107) -- ($(A107)! 0.5!(A132)$) -- ($(A117)! 0.5!(A132)$) -- (A117) -- cycle;
\fill[color18!50] ($(A118)! 0.5!(A107)$) -- ($(A119)! 0.5!(A107)$) -- (A107) -- cycle;
\fill[color18!50] ($(A118)! 0.5!(A107)$) -- ($(A126)! 0.5!(A107)$) -- (A107) -- cycle;
\fill[color18!50] ($(A126)! 0.5!(A107)$) -- ($(A132)! 0.5!(A107)$) -- (A107) -- cycle;
\fill[color18!50] ($(A113)! 0.5!(A124)$) -- ($(A114)! 0.5!(A124)$) -- (A124) -- cycle;
\fill[color18!50] ($(A113)! 0.5!(A124)$) -- ($(A123)! 0.5!(A124)$) -- (A124) -- cycle;
\fill[color18!50] (A115) -- ($(A115)! 0.5!(A114)$) -- ($(A124)! 0.5!(A114)$) -- (A124) -- cycle;
\fill[color18!50] (A115) -- ($(A115)! 0.5!(A116)$) -- ($(A125)! 0.5!(A116)$) -- (A125) -- cycle;
\fill[color18!50] (A115) -- (A124) -- (A125) -- cycle;
\fill[color18!50] (A117) -- ($(A117)! 0.5!(A116)$) -- ($(A125)! 0.5!(A116)$) -- (A125) -- cycle;
\fill[color18!50] (A117) -- ($(A117)! 0.5!(A132)$) -- ($(A125)! 0.5!(A132)$) -- (A125) -- cycle;
\fill[color18!50] ($(A123)! 0.5!(A124)$) -- ($(A132)! 0.5!(A124)$) -- (A124) -- cycle;
\fill[color18!50] (A124) -- ($(A124)! 0.5!(A132)$) -- ($(A125)! 0.5!(A132)$) -- (A125) -- cycle;
\colorlet{color19}{mycolor5}
\fill[color19!50] ($(A107)! 0.5!(A116)$) -- ($(A115)! 0.5!(A116)$) -- (A116) -- cycle;
\fill[color19!50] ($(A107)! 0.5!(A116)$) -- ($(A117)! 0.5!(A116)$) -- (A116) -- cycle;
\fill[color19!50] ($(A115)! 0.5!(A116)$) -- ($(A125)! 0.5!(A116)$) -- (A116) -- cycle;
\fill[color19!50] ($(A117)! 0.5!(A116)$) -- ($(A125)! 0.5!(A116)$) -- (A116) -- cycle;
\colorlet{color20}{mycolor3}
\fill[color20!50] ($(A118)! 0.5!(A127)$) -- ($(A126)! 0.5!(A127)$) -- (A127) -- cycle;
\fill[color20!50] ($(A118)! 0.5!(A127)$) -- ($(A128)! 0.5!(A127)$) -- (A127) -- cycle;
\fill[color20!50] ($(A126)! 0.5!(A127)$) -- ($(A134)! 0.5!(A127)$) -- (A127) -- cycle;
\fill[color20!50] ($(A128)! 0.5!(A127)$) -- ($(A135)! 0.5!(A127)$) -- (A127) -- cycle;
\fill[color20!50] ($(A134)! 0.5!(A127)$) -- ($(A135)! 0.5!(A127)$) -- (A127) -- cycle;
\colorlet{color21}{mycolor3}
\fill[color21!50] ($(A122)! 0.5!(A131)$) -- ($(A130)! 0.5!(A131)$) -- (A131) -- cycle;
\fill[color21!50] ($(A122)! 0.5!(A131)$) -- ($(A132)! 0.5!(A131)$) -- (A131) -- cycle;
\fill[color21!50] ($(A130)! 0.5!(A131)$) -- ($(A137)! 0.5!(A131)$) -- (A131) -- cycle;
\fill[color21!50] ($(A132)! 0.5!(A131)$) -- ($(A138)! 0.5!(A131)$) -- (A131) -- cycle;
\fill[color21!50] ($(A137)! 0.5!(A131)$) -- ($(A138)! 0.5!(A131)$) -- (A131) -- cycle;
\colorlet{color22}{mycolor3}
\fill[color22!50] ($(A133)! 0.5!(A139)$) -- ($(A134)! 0.5!(A139)$) -- (A139) -- cycle;
\fill[color22!50] ($(A133)! 0.5!(A139)$) -- ($(A138)! 0.5!(A139)$) -- (A139) -- cycle;
\fill[color22!50] ($(A134)! 0.5!(A139)$) -- ($(A140)! 0.5!(A139)$) -- (A139) -- cycle;
\fill[color22!50] ($(A138)! 0.5!(A139)$) -- ($(A142)! 0.5!(A139)$) -- (A139) -- cycle;
\fill[color22!50] ($(A140)! 0.5!(A139)$) -- ($(A143)! 0.5!(A139)$) -- (A139) -- cycle;
\fill[color22!50] ($(A142)! 0.5!(A139)$) -- ($(A143)! 0.5!(A139)$) -- (A139) -- cycle;
\draw[color22, thick] ($(A133)! 0.5!(A139)$) -- ($(A134)! 0.5!(A139)$);
\draw[color22, thick] ($(A133)! 0.5!(A139)$) -- ($(A138)! 0.5!(A139)$);
\draw[color22, thick] ($(A134)! 0.5!(A139)$) -- ($(A140)! 0.5!(A139)$);
\draw[color22, thick] ($(A138)! 0.5!(A139)$) -- ($(A142)! 0.5!(A139)$);
\draw[color22, thick] ($(A140)! 0.5!(A139)$) -- ($(A143)! 0.5!(A139)$);
\draw[color22, thick] ($(A142)! 0.5!(A139)$) -- ($(A143)! 0.5!(A139)$);
\draw[color21, thick] ($(A122)! 0.5!(A131)$) -- ($(A130)! 0.5!(A131)$);
\draw[color21, thick] ($(A122)! 0.5!(A131)$) -- ($(A132)! 0.5!(A131)$);
\draw[color21, thick] ($(A130)! 0.5!(A131)$) -- ($(A137)! 0.5!(A131)$);
\draw[color21, thick] ($(A132)! 0.5!(A131)$) -- ($(A138)! 0.5!(A131)$);
\draw[color21, thick] ($(A137)! 0.5!(A131)$) -- ($(A138)! 0.5!(A131)$);
\draw[color20, thick] ($(A118)! 0.5!(A127)$) -- ($(A126)! 0.5!(A127)$);
\draw[color20, thick] ($(A118)! 0.5!(A127)$) -- ($(A128)! 0.5!(A127)$);
\draw[color20, thick] ($(A126)! 0.5!(A127)$) -- ($(A134)! 0.5!(A127)$);
\draw[color20, thick] ($(A128)! 0.5!(A127)$) -- ($(A135)! 0.5!(A127)$);
\draw[color20, thick] ($(A134)! 0.5!(A127)$) -- ($(A135)! 0.5!(A127)$);
\draw[color18, thick] ($(A93)! 0.5!(A107)$) -- ($(A94)! 0.5!(A107)$);
\draw[color18, thick] ($(A93)! 0.5!(A107)$) -- ($(A106)! 0.5!(A107)$);
\draw[color18, thick] ($(A94)! 0.5!(A107)$) -- ($(A108)! 0.5!(A107)$);
\draw[color18, thick] ($(A106)! 0.5!(A107)$) -- ($(A114)! 0.5!(A107)$);
\draw[color18, thick] ($(A108)! 0.5!(A107)$) -- ($(A119)! 0.5!(A107)$);
\draw[color18, thick] ($(A107)! 0.5!(A114)$) -- ($(A115)! 0.5!(A114)$);
\draw[color18, thick] ($(A107)! 0.5!(A116)$) -- ($(A115)! 0.5!(A116)$);
\draw[color18, thick] ($(A107)! 0.5!(A116)$) -- ($(A117)! 0.5!(A116)$);
\draw[color18, thick] ($(A107)! 0.5!(A132)$) -- ($(A117)! 0.5!(A132)$);
\draw[color18, thick] ($(A118)! 0.5!(A107)$) -- ($(A119)! 0.5!(A107)$);
\draw[color18, thick] ($(A118)! 0.5!(A107)$) -- ($(A126)! 0.5!(A107)$);
\draw[color18, thick] ($(A126)! 0.5!(A107)$) -- ($(A132)! 0.5!(A107)$);
\draw[color18, thick] ($(A113)! 0.5!(A124)$) -- ($(A114)! 0.5!(A124)$);
\draw[color18, thick] ($(A113)! 0.5!(A124)$) -- ($(A123)! 0.5!(A124)$);
\draw[color18, thick] ($(A115)! 0.5!(A114)$) -- ($(A124)! 0.5!(A114)$);
\draw[color18, thick] ($(A115)! 0.5!(A116)$) -- ($(A125)! 0.5!(A116)$);
\draw[color18, thick] ($(A117)! 0.5!(A116)$) -- ($(A125)! 0.5!(A116)$);
\draw[color18, thick] ($(A117)! 0.5!(A132)$) -- ($(A125)! 0.5!(A132)$);
\draw[color18, thick] ($(A123)! 0.5!(A124)$) -- ($(A132)! 0.5!(A124)$);
\draw[color18, thick] ($(A124)! 0.5!(A132)$) -- ($(A125)! 0.5!(A132)$);
\draw[color19, thick] ($(A107)! 0.5!(A116)$) -- ($(A115)! 0.5!(A116)$);
\draw[color19, thick] ($(A107)! 0.5!(A116)$) -- ($(A117)! 0.5!(A116)$);
\draw[color19, thick] ($(A115)! 0.5!(A116)$) -- ($(A125)! 0.5!(A116)$);
\draw[color19, thick] ($(A117)! 0.5!(A116)$) -- ($(A125)! 0.5!(A116)$);
\draw[color17, thick] ($(A91)! 0.5!(A105)$) -- ($(A92)! 0.5!(A105)$);
\draw[color17, thick] ($(A91)! 0.5!(A105)$) -- ($(A104)! 0.5!(A105)$);
\draw[color17, thick] ($(A92)! 0.5!(A105)$) -- ($(A106)! 0.5!(A105)$);
\draw[color17, thick] ($(A104)! 0.5!(A105)$) -- ($(A114)! 0.5!(A105)$);
\draw[color17, thick] ($(A106)! 0.5!(A105)$) -- ($(A114)! 0.5!(A105)$);
\draw[color16, thick] ($(A89)! 0.5!(A103)$) -- ($(A90)! 0.5!(A103)$);
\draw[color16, thick] ($(A89)! 0.5!(A103)$) -- ($(A102)! 0.5!(A103)$);
\draw[color16, thick] ($(A90)! 0.5!(A103)$) -- ($(A104)! 0.5!(A103)$);
\draw[color16, thick] ($(A102)! 0.5!(A103)$) -- ($(A114)! 0.5!(A103)$);
\draw[color16, thick] ($(A104)! 0.5!(A103)$) -- ($(A114)! 0.5!(A103)$);
\draw[color15, thick] ($(A45)! 0.5!(A87)$) -- ($(A46)! 0.5!(A87)$);
\draw[color15, thick] ($(A45)! 0.5!(A87)$) -- ($(A88)! 0.5!(A87)$);
\draw[color15, thick] ($(A46)! 0.5!(A87)$) -- ($(A86)! 0.5!(A87)$);
\draw[color15, thick] ($(A86)! 0.5!(A87)$) -- ($(A101)! 0.5!(A87)$);
\draw[color15, thick] ($(A88)! 0.5!(A87)$) -- ($(A102)! 0.5!(A87)$);
\draw[color15, thick] ($(A101)! 0.5!(A87)$) -- ($(A102)! 0.5!(A87)$);
\draw[color14, thick] ($(A47)! 0.5!(A85)$) -- ($(A48)! 0.5!(A85)$);
\draw[color14, thick] ($(A47)! 0.5!(A85)$) -- ($(A86)! 0.5!(A85)$);
\draw[color14, thick] ($(A48)! 0.5!(A85)$) -- ($(A84)! 0.5!(A85)$);
\draw[color14, thick] ($(A84)! 0.5!(A85)$) -- ($(A99)! 0.5!(A85)$);
\draw[color14, thick] ($(A86)! 0.5!(A85)$) -- ($(A100)! 0.5!(A85)$);
\draw[color14, thick] ($(A99)! 0.5!(A85)$) -- ($(A100)! 0.5!(A85)$);
\draw[color13, thick] ($(A49)! 0.5!(A83)$) -- ($(A50)! 0.5!(A83)$);
\draw[color13, thick] ($(A49)! 0.5!(A83)$) -- ($(A84)! 0.5!(A83)$);
\draw[color13, thick] ($(A50)! 0.5!(A83)$) -- ($(A82)! 0.5!(A83)$);
\draw[color13, thick] ($(A82)! 0.5!(A83)$) -- ($(A97)! 0.5!(A83)$);
\draw[color13, thick] ($(A84)! 0.5!(A83)$) -- ($(A98)! 0.5!(A83)$);
\draw[color13, thick] ($(A97)! 0.5!(A83)$) -- ($(A98)! 0.5!(A83)$);
\draw[color12, thick] ($(A72)! 0.5!(A76)$) -- ($(A73)! 0.5!(A76)$);
\draw[color12, thick] ($(A72)! 0.5!(A76)$) -- ($(A75)! 0.5!(A76)$);
\draw[color12, thick] ($(A73)! 0.5!(A76)$) -- ($(A77)! 0.5!(A76)$);
\draw[color12, thick] ($(A75)! 0.5!(A76)$) -- ($(A78)! 0.5!(A76)$);
\draw[color12, thick] ($(A77)! 0.5!(A76)$) -- ($(A79)! 0.5!(A76)$);
\draw[color12, thick] ($(A78)! 0.5!(A76)$) -- ($(A79)! 0.5!(A76)$);
\draw[color10, thick] ($(A64)! 0.5!(A69)$) -- ($(A70)! 0.5!(A69)$);
\draw[color10, thick] ($(A64)! 0.5!(A69)$) -- ($(A72)! 0.5!(A69)$);
\draw[color10, thick] ($(A70)! 0.5!(A69)$) -- ($(A74)! 0.5!(A69)$);
\draw[color10, thick] ($(A72)! 0.5!(A69)$) -- ($(A73)! 0.5!(A69)$);
\draw[color10, thick] ($(A73)! 0.5!(A69)$) -- ($(A74)! 0.5!(A69)$);
\draw[color9, thick] ($(A50)! 0.5!(A58)$) -- ($(A51)! 0.5!(A58)$);
\draw[color9, thick] ($(A50)! 0.5!(A58)$) -- ($(A57)! 0.5!(A58)$);
\draw[color9, thick] ($(A51)! 0.5!(A58)$) -- ($(A59)! 0.5!(A58)$);
\draw[color9, thick] ($(A57)! 0.5!(A58)$) -- ($(A64)! 0.5!(A58)$);
\draw[color9, thick] ($(A59)! 0.5!(A58)$) -- ($(A65)! 0.5!(A58)$);
\draw[color9, thick] ($(A64)! 0.5!(A58)$) -- ($(A65)! 0.5!(A58)$);
\draw[color8, thick] ($(A48)! 0.5!(A56)$) -- ($(A49)! 0.5!(A56)$);
\draw[color8, thick] ($(A48)! 0.5!(A56)$) -- ($(A55)! 0.5!(A56)$);
\draw[color8, thick] ($(A49)! 0.5!(A56)$) -- ($(A57)! 0.5!(A56)$);
\draw[color8, thick] ($(A55)! 0.5!(A56)$) -- ($(A62)! 0.5!(A56)$);
\draw[color8, thick] ($(A57)! 0.5!(A56)$) -- ($(A63)! 0.5!(A56)$);
\draw[color8, thick] ($(A62)! 0.5!(A56)$) -- ($(A63)! 0.5!(A56)$);
\draw[color7, thick] ($(A46)! 0.5!(A54)$) -- ($(A47)! 0.5!(A54)$);
\draw[color7, thick] ($(A46)! 0.5!(A54)$) -- ($(A53)! 0.5!(A54)$);
\draw[color7, thick] ($(A47)! 0.5!(A54)$) -- ($(A55)! 0.5!(A54)$);
\draw[color7, thick] ($(A53)! 0.5!(A54)$) -- ($(A60)! 0.5!(A54)$);
\draw[color7, thick] ($(A55)! 0.5!(A54)$) -- ($(A61)! 0.5!(A54)$);
\draw[color7, thick] ($(A60)! 0.5!(A54)$) -- ($(A61)! 0.5!(A54)$);
\draw[color1, thick] ($(A1)! 0.5!(A0)$) -- ($(A9)! 0.5!(A0)$);
\draw[color1, thick] ($(A0)! 0.5!(A1)$) -- ($(A88)! 0.5!(A1)$);
\draw[color1, thick] ($(A0)! 0.5!(A9)$) -- ($(A45)! 0.5!(A9)$);
\draw[color1, thick] ($(A1)! 0.5!(A2)$) -- ($(A10)! 0.5!(A2)$);
\draw[color1, thick] ($(A2)! 0.5!(A1)$) -- ($(A89)! 0.5!(A1)$);
\draw[color1, thick] ($(A88)! 0.5!(A1)$) -- ($(A89)! 0.5!(A1)$);
\draw[color1, thick] ($(A3)! 0.5!(A2)$) -- ($(A11)! 0.5!(A2)$);
\draw[color1, thick] ($(A2)! 0.5!(A3)$) -- ($(A90)! 0.5!(A3)$);
\draw[color1, thick] ($(A10)! 0.5!(A2)$) -- ($(A11)! 0.5!(A2)$);
\draw[color1, thick] ($(A3)! 0.5!(A4)$) -- ($(A12)! 0.5!(A4)$);
\draw[color1, thick] ($(A4)! 0.5!(A3)$) -- ($(A91)! 0.5!(A3)$);
\draw[color1, thick] ($(A90)! 0.5!(A3)$) -- ($(A91)! 0.5!(A3)$);
\draw[color1, thick] ($(A5)! 0.5!(A4)$) -- ($(A13)! 0.5!(A4)$);
\draw[color1, thick] ($(A4)! 0.5!(A5)$) -- ($(A92)! 0.5!(A5)$);
\draw[color1, thick] ($(A12)! 0.5!(A4)$) -- ($(A13)! 0.5!(A4)$);
\draw[color1, thick] ($(A5)! 0.5!(A6)$) -- ($(A14)! 0.5!(A6)$);
\draw[color1, thick] ($(A6)! 0.5!(A5)$) -- ($(A93)! 0.5!(A5)$);
\draw[color1, thick] ($(A92)! 0.5!(A5)$) -- ($(A93)! 0.5!(A5)$);
\draw[color1, thick] ($(A7)! 0.5!(A6)$) -- ($(A15)! 0.5!(A6)$);
\draw[color1, thick] ($(A6)! 0.5!(A7)$) -- ($(A94)! 0.5!(A7)$);
\draw[color1, thick] ($(A14)! 0.5!(A6)$) -- ($(A15)! 0.5!(A6)$);
\draw[color1, thick] ($(A7)! 0.5!(A8)$) -- ($(A16)! 0.5!(A8)$);
\draw[color1, thick] ($(A8)! 0.5!(A7)$) -- ($(A95)! 0.5!(A7)$);
\draw[color1, thick] ($(A94)! 0.5!(A7)$) -- ($(A95)! 0.5!(A7)$);
\draw[color1, thick] ($(A9)! 0.5!(A17)$) -- ($(A10)! 0.5!(A17)$);
\draw[color1, thick] ($(A17)! 0.5!(A9)$) -- ($(A53)! 0.5!(A9)$);
\draw[color1, thick] ($(A45)! 0.5!(A9)$) -- ($(A53)! 0.5!(A9)$);
\draw[color1, thick] ($(A10)! 0.5!(A17)$) -- ($(A18)! 0.5!(A17)$);
\draw[color1, thick] ($(A11)! 0.5!(A19)$) -- ($(A12)! 0.5!(A19)$);
\draw[color1, thick] ($(A11)! 0.5!(A19)$) -- ($(A18)! 0.5!(A19)$);
\draw[color1, thick] ($(A12)! 0.5!(A19)$) -- ($(A20)! 0.5!(A19)$);
\draw[color1, thick] ($(A13)! 0.5!(A21)$) -- ($(A14)! 0.5!(A21)$);
\draw[color1, thick] ($(A13)! 0.5!(A21)$) -- ($(A20)! 0.5!(A21)$);
\draw[color1, thick] ($(A14)! 0.5!(A21)$) -- ($(A22)! 0.5!(A21)$);
\draw[color1, thick] ($(A15)! 0.5!(A23)$) -- ($(A16)! 0.5!(A23)$);
\draw[color1, thick] ($(A15)! 0.5!(A23)$) -- ($(A22)! 0.5!(A23)$);
\draw[color1, thick] ($(A18)! 0.5!(A17)$) -- ($(A24)! 0.5!(A17)$);
\draw[color1, thick] ($(A17)! 0.5!(A24)$) -- ($(A60)! 0.5!(A24)$);
\draw[color1, thick] ($(A18)! 0.5!(A19)$) -- ($(A24)! 0.5!(A19)$);
\draw[color1, thick] ($(A20)! 0.5!(A19)$) -- ($(A24)! 0.5!(A19)$);
\draw[color1, thick] ($(A20)! 0.5!(A21)$) -- ($(A24)! 0.5!(A21)$);
\draw[color1, thick] ($(A22)! 0.5!(A21)$) -- ($(A24)! 0.5!(A21)$);
\draw[color1, thick] ($(A22)! 0.5!(A23)$) -- ($(A29)! 0.5!(A23)$);
\draw[color1, thick] ($(A22)! 0.5!(A25)$) -- ($(A24)! 0.5!(A25)$);
\draw[color1, thick] ($(A22)! 0.5!(A25)$) -- ($(A26)! 0.5!(A25)$);
\draw[color1, thick] ($(A22)! 0.5!(A27)$) -- ($(A26)! 0.5!(A27)$);
\draw[color1, thick] ($(A22)! 0.5!(A27)$) -- ($(A35)! 0.5!(A27)$);
\draw[color1, thick] ($(A22)! 0.5!(A32)$) -- ($(A28)! 0.5!(A32)$);
\draw[color1, thick] ($(A22)! 0.5!(A32)$) -- ($(A35)! 0.5!(A32)$);
\draw[color1, thick] ($(A24)! 0.5!(A25)$) -- ($(A30)! 0.5!(A25)$);
\draw[color1, thick] ($(A60)! 0.5!(A24)$) -- ($(A61)! 0.5!(A24)$);
\draw[color1, thick] ($(A61)! 0.5!(A24)$) -- ($(A62)! 0.5!(A24)$);
\draw[color1, thick] ($(A62)! 0.5!(A24)$) -- ($(A63)! 0.5!(A24)$);
\draw[color1, thick] ($(A63)! 0.5!(A24)$) -- ($(A64)! 0.5!(A24)$);
\draw[color1, thick] ($(A24)! 0.5!(A64)$) -- ($(A66)! 0.5!(A64)$);
\draw[color1, thick] ($(A26)! 0.5!(A25)$) -- ($(A31)! 0.5!(A25)$);
\draw[color1, thick] ($(A30)! 0.5!(A25)$) -- ($(A31)! 0.5!(A25)$);
\draw[color1, thick] ($(A26)! 0.5!(A27)$) -- ($(A31)! 0.5!(A27)$);
\draw[color1, thick] ($(A31)! 0.5!(A27)$) -- ($(A35)! 0.5!(A27)$);
\draw[color1, thick] ($(A28)! 0.5!(A34)$) -- ($(A29)! 0.5!(A34)$);
\draw[color1, thick] ($(A28)! 0.5!(A32)$) -- ($(A33)! 0.5!(A32)$);
\draw[color1, thick] ($(A28)! 0.5!(A34)$) -- ($(A33)! 0.5!(A34)$);
\draw[color1, thick] ($(A30)! 0.5!(A71)$) -- ($(A35)! 0.5!(A71)$);
\draw[color1, thick] ($(A30)! 0.5!(A71)$) -- ($(A66)! 0.5!(A71)$);
\draw[color1, thick] ($(A33)! 0.5!(A32)$) -- ($(A37)! 0.5!(A32)$);
\draw[color1, thick] ($(A35)! 0.5!(A32)$) -- ($(A36)! 0.5!(A32)$);
\draw[color1, thick] ($(A36)! 0.5!(A32)$) -- ($(A37)! 0.5!(A32)$);
\draw[color1, thick] ($(A33)! 0.5!(A34)$) -- ($(A38)! 0.5!(A34)$);
\draw[color1, thick] ($(A35)! 0.5!(A39)$) -- ($(A36)! 0.5!(A39)$);
\draw[color1, thick] ($(A39)! 0.5!(A35)$) -- ($(A75)! 0.5!(A35)$);
\draw[color1, thick] ($(A35)! 0.5!(A64)$) -- ($(A68)! 0.5!(A64)$);
\draw[color1, thick] ($(A64)! 0.5!(A35)$) -- ($(A72)! 0.5!(A35)$);
\draw[color1, thick] ($(A35)! 0.5!(A71)$) -- ($(A68)! 0.5!(A71)$);
\draw[color1, thick] ($(A72)! 0.5!(A35)$) -- ($(A75)! 0.5!(A35)$);
\draw[color1, thick] ($(A36)! 0.5!(A39)$) -- ($(A40)! 0.5!(A39)$);
\draw[color1, thick] ($(A37)! 0.5!(A41)$) -- ($(A38)! 0.5!(A41)$);
\draw[color1, thick] ($(A37)! 0.5!(A41)$) -- ($(A40)! 0.5!(A41)$);
\draw[color1, thick] ($(A40)! 0.5!(A39)$) -- ($(A42)! 0.5!(A39)$);
\draw[color1, thick] ($(A39)! 0.5!(A42)$) -- ($(A78)! 0.5!(A42)$);
\draw[color1, thick] ($(A40)! 0.5!(A41)$) -- ($(A43)! 0.5!(A41)$);
\draw[color1, thick] ($(A42)! 0.5!(A44)$) -- ($(A43)! 0.5!(A44)$);
\draw[color1, thick] ($(A44)! 0.5!(A42)$) -- ($(A80)! 0.5!(A42)$);
\draw[color1, thick] ($(A78)! 0.5!(A42)$) -- ($(A80)! 0.5!(A42)$);
\draw[color1, thick] ($(A51)! 0.5!(A81)$) -- ($(A52)! 0.5!(A81)$);
\draw[color1, thick] ($(A51)! 0.5!(A81)$) -- ($(A82)! 0.5!(A81)$);
\draw[color1, thick] ($(A66)! 0.5!(A64)$) -- ($(A67)! 0.5!(A64)$);
\draw[color1, thick] ($(A67)! 0.5!(A64)$) -- ($(A68)! 0.5!(A64)$);
\draw[color1, thick] ($(A66)! 0.5!(A71)$) -- ($(A67)! 0.5!(A71)$);
\draw[color1, thick] ($(A67)! 0.5!(A71)$) -- ($(A68)! 0.5!(A71)$);
\draw[color1, thick] ($(A82)! 0.5!(A81)$) -- ($(A96)! 0.5!(A81)$);
\draw[color1, thick] ($(A96)! 0.5!(A109)$) -- ($(A97)! 0.5!(A109)$);
\draw[color1, thick] ($(A97)! 0.5!(A109)$) -- ($(A110)! 0.5!(A109)$);
\draw[color1, thick] ($(A110)! 0.5!(A109)$) -- ($(A120)! 0.5!(A109)$);
\draw[color1, thick] ($(A120)! 0.5!(A129)$) -- ($(A121)! 0.5!(A129)$);
\draw[color1, thick] ($(A121)! 0.5!(A129)$) -- ($(A130)! 0.5!(A129)$);
\draw[color1, thick] ($(A130)! 0.5!(A129)$) -- ($(A136)! 0.5!(A129)$);
\draw[color1, thick] ($(A136)! 0.5!(A141)$) -- ($(A137)! 0.5!(A141)$);
\draw[color1, thick] ($(A137)! 0.5!(A141)$) -- ($(A142)! 0.5!(A141)$);
\draw[color1, thick] ($(A142)! 0.5!(A141)$) -- ($(A144)! 0.5!(A141)$);
\draw[color11, thick] ($(A30)! 0.5!(A71)$) -- ($(A35)! 0.5!(A71)$);
\draw[color11, thick] ($(A30)! 0.5!(A71)$) -- ($(A66)! 0.5!(A71)$);
\draw[color11, thick] ($(A35)! 0.5!(A71)$) -- ($(A68)! 0.5!(A71)$);
\draw[color11, thick] ($(A66)! 0.5!(A71)$) -- ($(A67)! 0.5!(A71)$);
\draw[color11, thick] ($(A67)! 0.5!(A71)$) -- ($(A68)! 0.5!(A71)$);
\draw[color6, thick] ($(A22)! 0.5!(A32)$) -- ($(A28)! 0.5!(A32)$);
\draw[color6, thick] ($(A22)! 0.5!(A32)$) -- ($(A35)! 0.5!(A32)$);
\draw[color6, thick] ($(A28)! 0.5!(A32)$) -- ($(A33)! 0.5!(A32)$);
\draw[color6, thick] ($(A33)! 0.5!(A32)$) -- ($(A37)! 0.5!(A32)$);
\draw[color6, thick] ($(A35)! 0.5!(A32)$) -- ($(A36)! 0.5!(A32)$);
\draw[color6, thick] ($(A36)! 0.5!(A32)$) -- ($(A37)! 0.5!(A32)$);
\draw[color5, thick] ($(A22)! 0.5!(A27)$) -- ($(A26)! 0.5!(A27)$);
\draw[color5, thick] ($(A22)! 0.5!(A27)$) -- ($(A35)! 0.5!(A27)$);
\draw[color5, thick] ($(A26)! 0.5!(A27)$) -- ($(A31)! 0.5!(A27)$);
\draw[color5, thick] ($(A31)! 0.5!(A27)$) -- ($(A35)! 0.5!(A27)$);
\draw[color4, thick] ($(A22)! 0.5!(A25)$) -- ($(A24)! 0.5!(A25)$);
\draw[color4, thick] ($(A22)! 0.5!(A25)$) -- ($(A26)! 0.5!(A25)$);
\draw[color4, thick] ($(A24)! 0.5!(A25)$) -- ($(A30)! 0.5!(A25)$);
\draw[color4, thick] ($(A26)! 0.5!(A25)$) -- ($(A31)! 0.5!(A25)$);
\draw[color4, thick] ($(A30)! 0.5!(A25)$) -- ($(A31)! 0.5!(A25)$);
\draw[color3, thick] ($(A13)! 0.5!(A21)$) -- ($(A14)! 0.5!(A21)$);
\draw[color3, thick] ($(A13)! 0.5!(A21)$) -- ($(A20)! 0.5!(A21)$);
\draw[color3, thick] ($(A14)! 0.5!(A21)$) -- ($(A22)! 0.5!(A21)$);
\draw[color3, thick] ($(A20)! 0.5!(A21)$) -- ($(A24)! 0.5!(A21)$);
\draw[color3, thick] ($(A22)! 0.5!(A21)$) -- ($(A24)! 0.5!(A21)$);
\draw[color2, thick] ($(A11)! 0.5!(A19)$) -- ($(A12)! 0.5!(A19)$);
\draw[color2, thick] ($(A11)! 0.5!(A19)$) -- ($(A18)! 0.5!(A19)$);
\draw[color2, thick] ($(A12)! 0.5!(A19)$) -- ($(A20)! 0.5!(A19)$);
\draw[color2, thick] ($(A18)! 0.5!(A19)$) -- ($(A24)! 0.5!(A19)$);
\draw[color2, thick] ($(A20)! 0.5!(A19)$) -- ($(A24)! 0.5!(A19)$);
\foreach \a/\b/\c in {0/1/9,0/1/88,0/9/45,0/45/88,1/2/10,1/2/89,1/9/10,1/88/89,2/3/11,2/3/90,2/10/11,2/89/90,3/4/12,3/4/91,3/11/12,3/90/91,4/5/13,4/5/92,4/12/13,4/91/92,5/6/14,5/6/93,5/13/14,5/92/93,6/7/15,6/7/94,6/14/15,6/93/94,7/8/16,7/8/95,7/15/16,7/94/95,9/10/17,9/17/53,9/45/53,10/11/18,10/17/18,11/12/19,11/18/19,12/13/20,12/19/20,13/14/21,13/20/21,14/15/22,14/21/22,15/16/23,15/22/23,17/18/24,17/24/60,17/53/60,18/19/24,19/20/24,20/21/24,21/22/24,22/23/29,22/24/25,22/25/26,22/26/27,22/27/35,22/28/29,22/28/32,22/32/35,24/25/30,24/30/66,24/60/61,24/61/62,24/62/63,24/63/64,24/64/66,25/26/31,25/30/31,26/27/31,27/31/35,28/29/34,28/32/33,28/33/34,30/31/35,30/35/71,30/66/71,32/33/37,32/35/36,32/36/37,33/34/38,33/37/38,35/36/39,35/39/75,35/64/68,35/64/72,35/68/71,35/72/75,36/37/40,36/39/40,37/38/41,37/40/41,39/40/42,39/42/78,39/75/78,40/41/43,40/42/43,42/43/44,42/44/80,42/78/80,45/46/53,45/46/87,45/87/88,46/47/54,46/47/86,46/53/54,46/86/87,47/48/55,47/48/85,47/54/55,47/85/86,48/49/56,48/49/84,48/55/56,48/84/85,49/50/57,49/50/83,49/56/57,49/83/84,50/51/58,50/51/82,50/57/58,50/82/83,51/52/59,51/52/81,51/58/59,51/81/82,53/54/60,54/55/61,54/60/61,55/56/62,55/61/62,56/57/63,56/62/63,57/58/64,57/63/64,58/59/65,58/64/65,64/65/70,64/66/67,64/67/68,64/69/70,64/69/72,66/67/71,67/68/71,69/70/74,69/72/73,69/73/74,72/73/76,72/75/76,73/74/77,73/76/77,75/76/78,76/77/79,76/78/79,78/79/80,81/82/96,82/83/97,82/96/97,83/84/98,83/97/98,84/85/99,84/98/99,85/86/100,85/99/100,86/87/101,86/100/101,87/88/102,87/101/102,88/89/102,89/90/103,89/102/103,90/91/104,90/103/104,91/92/105,91/104/105,92/93/106,92/105/106,93/94/107,93/106/107,94/95/108,94/107/108,96/97/109,97/98/114,97/109/110,97/110/122,97/111/112,97/111/132,97/112/113,97/113/114,97/122/132,98/99/114,99/100/114,100/101/114,101/102/114,102/103/114,103/104/114,104/105/114,105/106/114,106/107/114,107/108/119,107/114/115,107/115/116,107/116/117,107/117/132,107/118/119,107/118/126,107/126/132,109/110/120,110/120/121,110/121/122,111/112/123,111/123/132,112/113/123,113/114/124,113/123/124,114/115/124,115/116/125,115/124/125,116/117/125,117/125/132,118/119/128,118/126/127,118/127/128,120/121/129,121/122/130,121/129/130,122/130/131,122/131/132,123/124/132,124/125/132,126/127/134,126/132/133,126/133/134,127/128/135,127/134/135,129/130/136,130/131/137,130/136/137,131/132/138,131/137/138,132/133/138,133/134/139,133/138/139,134/135/140,134/139/140,136/137/141,137/138/142,137/141/142,138/139/142,139/140/143,139/142/143,141/142/144,142/143/144}{
  \draw[black!30] (A\a) -- (A\b) -- (A\c) -- cycle;
}

\draw[black, thick] (A22) -- (A35) -- (A64);
\draw[black, thick] (A22) -- (A24) -- (A64);
\draw[black, thick] (A107) -- (A132) -- (A97);
\draw[black, thick] (A107) -- (A114) -- (A97);

\foreach \i in {1,3,5,7,9,10,11,12,13,14,15,16,18,20,22,24,26,28,29,30,31,33,35,36,37,38,40,42,43,54,56,58,66,67,68,69,76,81,83,85,87,103,105,107,109,115,117,124,125,127,129,131,139,141}{
  \fill[myred] (A\i) circle (3pt);
}
\foreach \i in {0,2,4,6,8,17,19,21,23,25,27,32,34,39,41,44,45,46,47,48,49,50,51,52,53,55,57,59,60,61,62,63,64,65,70,71,72,73,74,75,77,78,79,80,82,84,86,88,89,90,91,92,93,94,95,96,97,98,99,100,101,102,104,106,108,110,111,112,113,114,116,118,119,120,121,122,123,126,128,130,132,133,134,135,136,137,138,140,142,143,144}{
  \fill[myblue] (A\i) circle (3pt);
}

%% file: tikz/deg6_triangulations.tex
\begin{figure}[!t]
  \centering
   \begin{subfigure}{0.48\textwidth}
    \centering
  \begin{tikzpicture}[scale=\scalefactor, every node/.style={font=\tiny}]
    \draw (0,5) -- (0,6) -- (1,5) -- cycle;
    \draw (0,5) -- (1,4) -- (1,5) -- cycle;
    \draw (0,4) -- (0,5) -- (1,4) -- cycle;
    \draw (1,4) -- (1,5) -- (2,4) -- cycle;
    \draw (0,4) -- (1,3) -- (1,4) -- cycle;
    \draw (1,4) -- (2,3) -- (2,4) -- cycle;
    \draw (0,3) -- (0,4) -- (1,3) -- cycle;
    \draw (1,3) -- (1,4) -- (2,3) -- cycle;
    \draw (2,3) -- (2,4) -- (3,3) -- cycle;
    \draw (0,2) -- (0,3) -- (1,3) -- cycle;
    \draw (0,2) -- (1,2) -- (1,3) -- cycle;
    \draw (1,2) -- (1,3) -- (2,3) -- cycle;
    \draw (0,1) -- (0,2) -- (1,2) -- cycle;
    \draw (1,2) -- (2,3) -- (3,3) -- cycle;
    \draw (0,1) -- (1,2) -- (3,3) -- cycle;
    \draw (0,1) -- (2,2) -- (3,3) -- cycle;
    \draw (0,1) -- (1,1) -- (2,2) -- cycle;
    \draw (0,1) -- (1,0) -- (1,1) -- cycle;
    \draw (1,0) -- (1,1) -- (2,2) -- cycle;
    \draw (0,0) -- (0,1) -- (1,0) -- cycle;
    \draw (1,0) -- (2,2) -- (3,3) -- cycle;
    \draw (1,0) -- (2,1) -- (3,3) -- cycle;
    \draw (2,1) -- (3,2) -- (3,3) -- cycle;
    \draw (2,1) -- (3,1) -- (3,2) -- cycle;
    \draw (1,0) -- (2,0) -- (2,1) -- cycle;
    \draw (2,0) -- (2,1) -- (3,1) -- cycle;
    \draw (3,2) -- (3,3) -- (4,2) -- cycle;
    \draw (3,2) -- (4,1) -- (4,2) -- cycle;
    \draw (3,1) -- (3,2) -- (4,1) -- cycle;
    \draw (4,1) -- (4,2) -- (5,1) -- cycle;
    \draw (2,0) -- (3,0) -- (3,1) -- cycle;
    \draw (3,1) -- (4,0) -- (4,1) -- cycle;
    \draw (3,0) -- (3,1) -- (4,0) -- cycle;
    \draw (4,1) -- (5,0) -- (5,1) -- cycle;
    \draw (4,0) -- (4,1) -- (5,0) -- cycle;
    \draw (5,0) -- (5,1) -- (6,0) -- cycle;
    \node[circle, fill=white, draw=black, minimum size=3mm, inner sep=0pt, outer sep=0pt, font=\fontsize{4}{5}\selectfont, text width=2.5mm, align=center] at (0,0) {4};
    \node[circle, fill=white, draw=black, minimum size=3mm, inner sep=0pt, outer sep=0pt, font=\fontsize{4}{5}\selectfont, text width=2.5mm, align=center] at (0,1) {2};
    \node[circle, fill=white, draw=black, minimum size=3mm, inner sep=0pt, outer sep=0pt, font=\fontsize{4}{5}\selectfont, text width=2.5mm, align=center] at (0,2) {12};
    \node[circle, fill=white, draw=black, minimum size=3mm, inner sep=0pt, outer sep=0pt, font=\fontsize{4}{5}\selectfont, text width=2.5mm, align=center] at (0,3) {24};
    \node[circle, fill=white, draw=black, minimum size=3mm, inner sep=0pt, outer sep=0pt, font=\fontsize{4}{5}\selectfont, text width=2.5mm, align=center] at (0,4) {37};
    \node[circle, fill=white, draw=black, minimum size=3mm, inner sep=0pt, outer sep=0pt, font=\fontsize{4}{5}\selectfont, text width=2.5mm, align=center] at (0,5) {52};
    \node[circle, fill=white, draw=black, minimum size=3mm, inner sep=0pt, outer sep=0pt, font=\fontsize{4}{5}\selectfont, text width=2.5mm, align=center] at (0,6) {69};
    \node[circle, fill=white, draw=black, minimum size=3mm, inner sep=0pt, outer sep=0pt, font=\fontsize{4}{5}\selectfont, text width=2.5mm, align=center] at (1,0) {2};
    \node[circle, fill=white, draw=black, minimum size=3mm, inner sep=0pt, outer sep=0pt, font=\fontsize{4}{5}\selectfont, text width=2.5mm, align=center] at (1,1) {1};
    \node[circle, fill=white, draw=black, minimum size=3mm, inner sep=0pt, outer sep=0pt, font=\fontsize{4}{5}\selectfont, text width=2.5mm, align=center] at (1,2) {3};
    \node[circle, fill=white, draw=black, minimum size=3mm, inner sep=0pt, outer sep=0pt, font=\fontsize{4}{5}\selectfont, text width=2.5mm, align=center] at (1,3) {14};
    \node[circle, fill=white, draw=black, minimum size=3mm, inner sep=0pt, outer sep=0pt, font=\fontsize{4}{5}\selectfont, text width=2.5mm, align=center] at (1,4) {28};
    \node[circle, fill=white, draw=black, minimum size=3mm, inner sep=0pt, outer sep=0pt, font=\fontsize{4}{5}\selectfont, text width=2.5mm, align=center] at (1,5) {44};
    \node[circle, fill=white, draw=black, minimum size=3mm, inner sep=0pt, outer sep=0pt, font=\fontsize{4}{5}\selectfont, text width=2.5mm, align=center] at (2,0) {12};
    \node[circle, fill=white, draw=black, minimum size=3mm, inner sep=0pt, outer sep=0pt, font=\fontsize{4}{5}\selectfont, text width=2.5mm, align=center] at (2,1) {3};
    \node[circle, fill=white, draw=black, minimum size=3mm, inner sep=0pt, outer sep=0pt, font=\fontsize{4}{5}\selectfont, text width=2.5mm, align=center] at (2,2) {0};
    \node[circle, fill=white, draw=black, minimum size=3mm, inner sep=0pt, outer sep=0pt, font=\fontsize{4}{5}\selectfont, text width=2.5mm, align=center] at (2,3) {6};
    \node[circle, fill=white, draw=black, minimum size=3mm, inner sep=0pt, outer sep=0pt, font=\fontsize{4}{5}\selectfont, text width=2.5mm, align=center] at (2,4) {21};
    \node[circle, fill=white, draw=black, minimum size=3mm, inner sep=0pt, outer sep=0pt, font=\fontsize{4}{5}\selectfont, text width=2.5mm, align=center] at (3,0) {24};
    \node[circle, fill=white, draw=black, minimum size=3mm, inner sep=0pt, outer sep=0pt, font=\fontsize{4}{5}\selectfont, text width=2.5mm, align=center] at (3,1) {14};
    \node[circle, fill=white, draw=black, minimum size=3mm, inner sep=0pt, outer sep=0pt, font=\fontsize{4}{5}\selectfont, text width=2.5mm, align=center] at (3,2) {6};
    \node[circle, fill=white, draw=black, minimum size=3mm, inner sep=0pt, outer sep=0pt, font=\fontsize{4}{5}\selectfont, text width=2.5mm, align=center] at (3,3) {0};
    \node[circle, fill=white, draw=black, minimum size=3mm, inner sep=0pt, outer sep=0pt, font=\fontsize{4}{5}\selectfont, text width=2.5mm, align=center] at (4,0) {37};
    \node[circle, fill=white, draw=black, minimum size=3mm, inner sep=0pt, outer sep=0pt, font=\fontsize{4}{5}\selectfont, text width=2.5mm, align=center] at (4,1) {28};
    \node[circle, fill=white, draw=black, minimum size=3mm, inner sep=0pt, outer sep=0pt, font=\fontsize{4}{5}\selectfont, text width=2.5mm, align=center] at (4,2) {21};
    \node[circle, fill=white, draw=black, minimum size=3mm, inner sep=0pt, outer sep=0pt, font=\fontsize{4}{5}\selectfont, text width=2.5mm, align=center] at (5,0) {52};
    \node[circle, fill=white, draw=black, minimum size=3mm, inner sep=0pt, outer sep=0pt, font=\fontsize{4}{5}\selectfont, text width=2.5mm, align=center] at (5,1) {44};
    \node[circle, fill=white, draw=black, minimum size=3mm, inner sep=0pt, outer sep=0pt, font=\fontsize{4}{5}\selectfont, text width=2.5mm, align=center] at (6,0) {69};
  \end{tikzpicture}
  \caption{\degsixA{} (realizes \degsixAtypes{} types)}
  \label{fig:bat}
  \end{subfigure}
  \begin{subfigure}{0.48\textwidth}
    \centering
  \begin{tikzpicture}[scale=\scalefactor, every node/.style={font=\tiny}]
    \draw (0,0) -- (1,0) -- (0,1) -- cycle;
    \draw (0,1) -- (1,0) -- (1,1) -- cycle;
    \draw (0,1) -- (1,1) -- (1,2) -- cycle;
    \draw (0,1) -- (1,2) -- (0,2) -- cycle;
    \draw (0,2) -- (1,2) -- (1,3) -- cycle;
    \draw (0,2) -- (1,3) -- (0,3) -- cycle;
    \draw (0,3) -- (1,3) -- (0,4) -- cycle;
    \draw (0,4) -- (1,3) -- (0,5) -- cycle;
    \draw (0,5) -- (1,3) -- (0,6) -- cycle;
    \draw (0,6) -- (1,3) -- (1,4) -- cycle;
    \draw (0,6) -- (1,4) -- (1,5) -- cycle;
    \draw (1,0) -- (2,0) -- (2,1) -- cycle;
    \draw (1,0) -- (2,1) -- (1,1) -- cycle;
    \draw (1,1) -- (2,1) -- (3,2) -- cycle;
    \draw (1,1) -- (2,2) -- (2,3) -- cycle;
    \draw (1,1) -- (2,3) -- (1,2) -- cycle;
    \draw (1,1) -- (3,2) -- (2,2) -- cycle;
    \draw (1,2) -- (2,3) -- (1,3) -- cycle;
    \draw (1,3) -- (2,3) -- (1,4) -- cycle;
    \draw (1,4) -- (2,3) -- (1,5) -- cycle;
    \draw (1,5) -- (2,3) -- (2,4) -- cycle;
    \draw (2,0) -- (3,0) -- (3,1) -- cycle;
    \draw (2,0) -- (3,1) -- (2,1) -- cycle;
    \draw (2,1) -- (3,1) -- (3,2) -- cycle;
    \draw (2,2) -- (3,2) -- (2,3) -- cycle;
    \draw (2,3) -- (3,2) -- (3,3) -- cycle;
    \draw (2,3) -- (3,3) -- (2,4) -- cycle;
    \draw (3,0) -- (4,0) -- (3,1) -- cycle;
    \draw (3,1) -- (4,0) -- (5,0) -- cycle;
    \draw (3,1) -- (4,1) -- (3,2) -- cycle;
    \draw (3,1) -- (5,0) -- (6,0) -- cycle;
    \draw (3,1) -- (6,0) -- (4,1) -- cycle;
    \draw (3,2) -- (4,1) -- (5,1) -- cycle;
    \draw (3,2) -- (4,2) -- (3,3) -- cycle;
    \draw (3,2) -- (5,1) -- (4,2) -- cycle;
    \draw (4,1) -- (6,0) -- (5,1) -- cycle;
    \node[circle, fill=white, draw=black, minimum size=3mm, inner sep=0pt, outer sep=0pt, font=\fontsize{4}{5}\selectfont, text width=2.5mm, align=center] at (0,0) {24};
    \node[circle, fill=white, draw=black, minimum size=3mm, inner sep=0pt, outer sep=0pt, font=\fontsize{4}{5}\selectfont, text width=2.5mm, align=center] at (0,1) {12};
    \node[circle, fill=white, draw=black, minimum size=3mm, inner sep=0pt, outer sep=0pt, font=\fontsize{4}{5}\selectfont, text width=2.5mm, align=center] at (0,2) {14};
    \node[circle, fill=white, draw=black, minimum size=3mm, inner sep=0pt, outer sep=0pt, font=\fontsize{4}{5}\selectfont, text width=2.5mm, align=center] at (0,3) {18};
    \node[circle, fill=white, draw=black, minimum size=3mm, inner sep=0pt, outer sep=0pt, font=\fontsize{4}{5}\selectfont, text width=2.5mm, align=center] at (0,4) {23};
    \node[circle, fill=white, draw=black, minimum size=3mm, inner sep=0pt, outer sep=0pt, font=\fontsize{4}{5}\selectfont, text width=2.5mm, align=center] at (0,5) {29};
    \node[circle, fill=white, draw=black, minimum size=3mm, inner sep=0pt, outer sep=0pt, font=\fontsize{4}{5}\selectfont, text width=2.5mm, align=center] at (0,6) {36};
    \node[circle, fill=white, draw=black, minimum size=3mm, inner sep=0pt, outer sep=0pt, font=\fontsize{4}{5}\selectfont, text width=2.5mm, align=center] at (1,0) {12};
    \node[circle, fill=white, draw=black, minimum size=3mm, inner sep=0pt, outer sep=0pt, font=\fontsize{4}{5}\selectfont, text width=2.5mm, align=center] at (1,1) {1};
    \node[circle, fill=white, draw=black, minimum size=3mm, inner sep=0pt, outer sep=0pt, font=\fontsize{4}{5}\selectfont, text width=2.5mm, align=center] at (1,2) {2};
    \node[circle, fill=white, draw=black, minimum size=3mm, inner sep=0pt, outer sep=0pt, font=\fontsize{4}{5}\selectfont, text width=2.5mm, align=center] at (1,3) {5};
    \node[circle, fill=white, draw=black, minimum size=3mm, inner sep=0pt, outer sep=0pt, font=\fontsize{4}{5}\selectfont, text width=2.5mm, align=center] at (1,4) {13};
    \node[circle, fill=white, draw=black, minimum size=3mm, inner sep=0pt, outer sep=0pt, font=\fontsize{4}{5}\selectfont, text width=2.5mm, align=center] at (1,5) {22};
    \node[circle, fill=white, draw=black, minimum size=3mm, inner sep=0pt, outer sep=0pt, font=\fontsize{4}{5}\selectfont, text width=2.5mm, align=center] at (2,0) {14};
    \node[circle, fill=white, draw=black, minimum size=3mm, inner sep=0pt, outer sep=0pt, font=\fontsize{4}{5}\selectfont, text width=2.5mm, align=center] at (2,1) {2};
    \node[circle, fill=white, draw=black, minimum size=3mm, inner sep=0pt, outer sep=0pt, font=\fontsize{4}{5}\selectfont, text width=2.5mm, align=center] at (2,2) {0};
    \node[circle, fill=white, draw=black, minimum size=3mm, inner sep=0pt, outer sep=0pt, font=\fontsize{4}{5}\selectfont, text width=2.5mm, align=center] at (2,3) {0};
    \node[circle, fill=white, draw=black, minimum size=3mm, inner sep=0pt, outer sep=0pt, font=\fontsize{4}{5}\selectfont, text width=2.5mm, align=center] at (2,4) {10};
    \node[circle, fill=white, draw=black, minimum size=3mm, inner sep=0pt, outer sep=0pt, font=\fontsize{4}{5}\selectfont, text width=2.5mm, align=center] at (3,0) {18};
    \node[circle, fill=white, draw=black, minimum size=3mm, inner sep=0pt, outer sep=0pt, font=\fontsize{4}{5}\selectfont, text width=2.5mm, align=center] at (3,1) {5};
    \node[circle, fill=white, draw=black, minimum size=3mm, inner sep=0pt, outer sep=0pt, font=\fontsize{4}{5}\selectfont, text width=2.5mm, align=center] at (3,2) {0};
    \node[circle, fill=white, draw=black, minimum size=3mm, inner sep=0pt, outer sep=0pt, font=\fontsize{4}{5}\selectfont, text width=2.5mm, align=center] at (3,3) {1};
    \node[circle, fill=white, draw=black, minimum size=3mm, inner sep=0pt, outer sep=0pt, font=\fontsize{4}{5}\selectfont, text width=2.5mm, align=center] at (4,0) {23};
    \node[circle, fill=white, draw=black, minimum size=3mm, inner sep=0pt, outer sep=0pt, font=\fontsize{4}{5}\selectfont, text width=2.5mm, align=center] at (4,1) {13};
    \node[circle, fill=white, draw=black, minimum size=3mm, inner sep=0pt, outer sep=0pt, font=\fontsize{4}{5}\selectfont, text width=2.5mm, align=center] at (4,2) {10};
    \node[circle, fill=white, draw=black, minimum size=3mm, inner sep=0pt, outer sep=0pt, font=\fontsize{4}{5}\selectfont, text width=2.5mm, align=center] at (5,0) {29};
    \node[circle, fill=white, draw=black, minimum size=3mm, inner sep=0pt, outer sep=0pt, font=\fontsize{4}{5}\selectfont, text width=2.5mm, align=center] at (5,1) {22};
    \node[circle, fill=white, draw=black, minimum size=3mm, inner sep=0pt, outer sep=0pt, font=\fontsize{4}{5}\selectfont, text width=2.5mm, align=center] at (6,0) {36};
  \end{tikzpicture}
  \caption{\degsixB{} (realizes \degsixBtypes{} types)}
  \label{fig:moth}
  \end{subfigure}
  \caption{Two regular triangulations of $6\cdot\Delta_2$ realizing all nonempty real schemes types of degree six.  Values at the vertices indicate lifting functions.}
  \label{fig:deg6-triangulations}
\end{figure}

%% file: tables/proof_deg6.tex
\small
\begin{longtable}{lrrrlr}
\caption{Proof of support for degree 6}
\label{tab:deg-6} \\
\toprule
Real scheme & \multicolumn{1}{c}{$p$} & \multicolumn{1}{c}{$n$} & \multicolumn{1}{c}{$p+n$} & \multicolumn{1}{c}{$\cT$} & \multicolumn{1}{c}{$\sigma$} \\
\midrule
\endfirsthead
\toprule
Real scheme & \multicolumn{1}{c}{$p$} & \multicolumn{1}{c}{$n$} & \multicolumn{1}{c}{$p+n$} & \multicolumn{1}{c}{$\cT$} & \multicolumn{1}{c}{$\sigma$} \\
\midrule
\endhead
\midrule
\endfoot
\bottomrule
\endlastfoot
$\langle 1 \sqcup 1\langle 9 \rangle \rangle$ & 2 & 9 & 11 & \degsixA & 1110\,1001\,1010\,0100\,1101\,1110\,1000 \\
$\langle 5 \sqcup 1\langle 5 \rangle \rangle$ & 6 & 5 & 11 & \degsixA & 1110\,1001\,1010\,0010\,1101\,1101\,0000 \\
$\langle 9 \sqcup 1\langle 1 \rangle \rangle$ & 10 & 1 & 11 & \degsixA & 1100\,0001\,1111\,0010\,1001\,1101\,0000 \\
$\langle 1\langle 9 \rangle \rangle$ & 1 & 9 & 10 & \degsixA & 1110\,1001\,0010\,0100\,1101\,1110\,1000 \\
$\langle 1 \sqcup 1\langle 8 \rangle \rangle$ & 2 & 8 & 10 & \degsixA & 1100\,0101\,0000\,0001\,1000\,1000\,0000 \\
$\langle 4 \sqcup 1\langle 5 \rangle \rangle$ & 5 & 5 & 10 & \degsixA & 1110\,1001\,1010\,0010\,1100\,0101\,0000 \\
$\langle 5 \sqcup 1\langle 4 \rangle \rangle$ & 6 & 4 & 10 & \degsixA & 1110\,1101\,0101\,0001\,1100\,1000\,0000 \\
$\langle 8 \sqcup 1\langle 1 \rangle \rangle$ & 9 & 1 & 10 & \degsixA & 1100\,0001\,1111\,0010\,1000\,0101\,0000 \\
$\langle 10 \rangle$ & 10 & 0 & 10 & \degsixB & 1100\,0001\,1111\,1011\,1001\,1101\,0010 \\
$\langle 1\langle 8 \rangle \rangle$ & 1 & 8 & 9 & \degsixA & 1100\,0101\,1000\,0000\,1000\,1000\,0000 \\
$\langle 1 \sqcup 1\langle 7 \rangle \rangle$ & 2 & 7 & 9 & \degsixA & 1100\,0101\,0000\,0001\,1001\,0000\,0000 \\
$\langle 2 \sqcup 1\langle 6 \rangle \rangle$ & 3 & 6 & 9 & \degsixA & 1110\,1001\,1010\,0000\,1100\,0100\,0000 \\
$\langle 3 \sqcup 1\langle 5 \rangle \rangle$ & 4 & 5 & 9 & \degsixA & 1110\,1001\,1010\,0000\,1101\,0100\,0000 \\
$\langle 4 \sqcup 1\langle 4 \rangle \rangle$ & 5 & 4 & 9 & \degsixA & 1110\,1101\,0101\,0001\,1000\,1000\,0000 \\
$\langle 5 \sqcup 1\langle 3 \rangle \rangle$ & 6 & 3 & 9 & \degsixA & 1110\,1101\,0101\,0001\,1101\,0000\,0000 \\
$\langle 6 \sqcup 1\langle 2 \rangle \rangle$ & 7 & 2 & 9 & \degsixA & 1100\,0001\,1111\,0000\,1000\,0100\,0000 \\
$\langle 7 \sqcup 1\langle 1 \rangle \rangle$ & 8 & 1 & 9 & \degsixA & 1100\,0001\,1111\,0000\,1001\,0100\,0000 \\
$\langle 9 \rangle$ & 9 & 0 & 9 & \degsixA & 1100\,0001\,1111\,0011\,1001\,1101\,0000 \\
$\langle 1\langle 7 \rangle \rangle$ & 1 & 7 & 8 & \degsixA & 1100\,0101\,1000\,0000\,1001\,0000\,0000 \\
$\langle 1 \sqcup 1\langle 6 \rangle \rangle$ & 2 & 6 & 8 & \degsixA & 1100\,0101\,0000\,0011\,1000\,0000\,0000 \\
$\langle 2 \sqcup 1\langle 5 \rangle \rangle$ & 3 & 5 & 8 & \degsixA & 1100\,0101\,0000\,0101\,1000\,0000\,0000 \\
$\langle 3 \sqcup 1\langle 4 \rangle \rangle$ & 4 & 4 & 8 & \degsixA & 1110\,1001\,0000\,0001\,0000\,1000\,0000 \\
$\langle 4 \sqcup 1\langle 3 \rangle \rangle$ & 5 & 3 & 8 & \degsixA & 1110\,1101\,0101\,0001\,1001\,0000\,0000 \\
$\langle 5 \sqcup 1\langle 2 \rangle \rangle$ & 6 & 2 & 8 & \degsixA & 1110\,1101\,0101\,0011\,1100\,0000\,0000 \\
$\langle 6 \sqcup 1\langle 1 \rangle \rangle$ & 7 & 1 & 8 & \degsixA & 1110\,1101\,0101\,0101\,1100\,0000\,0000 \\
$\langle 8 \rangle$ & 8 & 0 & 8 & \degsixA & 1100\,0001\,1111\,0011\,1000\,0101\,0000 \\
$\langle 1\langle 6 \rangle \rangle$ & 1 & 6 & 7 & \degsixA & 1100\,0101\,1000\,0010\,1000\,0000\,0000 \\
$\langle 1 \sqcup 1\langle 5 \rangle \rangle$ & 2 & 5 & 7 & \degsixA & 1100\,0101\,0010\,0011\,0000\,0000\,0000 \\
$\langle 2 \sqcup 1\langle 4 \rangle \rangle$ & 3 & 4 & 7 & \degsixA & 1100\,0101\,0010\,0101\,0000\,0000\,0000 \\
$\langle 3 \sqcup 1\langle 3 \rangle \rangle$ & 4 & 3 & 7 & \degsixA & 1110\,1001\,0000\,0001\,0001\,0000\,0000 \\
$\langle 4 \sqcup 1\langle 2 \rangle \rangle$ & 5 & 2 & 7 & \degsixA & 1110\,1101\,0101\,0011\,1000\,0000\,0000 \\
$\langle 5 \sqcup 1\langle 1 \rangle \rangle$ & 6 & 1 & 7 & \degsixA & 1110\,1101\,0101\,0101\,1000\,0000\,0000 \\
$\langle 7 \rangle$ & 7 & 0 & 7 & \degsixA & 1100\,0001\,1111\,0001\,1001\,0100\,0000 \\
$\langle 1\langle 5 \rangle \rangle$ & 1 & 5 & 6 & \degsixA & 1100\,0101\,1010\,0010\,0000\,0000\,0000 \\
$\langle 1 \sqcup 1\langle 4 \rangle \rangle$ & 2 & 4 & 6 & \degsixA & 1100\,0101\,1010\,0100\,0000\,0000\,0000 \\
$\langle 2 \sqcup 1\langle 3 \rangle \rangle$ & 3 & 3 & 6 & \degsixA & 1100\,0101\,0100\,0101\,0000\,0000\,0000 \\
$\langle 3 \sqcup 1\langle 2 \rangle \rangle$ & 4 & 2 & 6 & \degsixA & 1100\,1101\,0001\,0101\,0000\,0000\,0000 \\
$\langle 4 \sqcup 1\langle 1 \rangle \rangle$ & 5 & 1 & 6 & \degsixA & 1110\,1001\,0000\,0101\,0000\,0000\,0000 \\
$\langle 6 \rangle$ & 6 & 0 & 6 & \degsixA & 1110\,1101\,0101\,0100\,1100\,0000\,0000 \\
$\langle 1\langle 4 \rangle \rangle$ & 1 & 4 & 5 & \degsixA & 1100\,0101\,1010\,0000\,0000\,0000\,0000 \\
$\langle 1 \sqcup 1\langle 3 \rangle \rangle$ & 2 & 3 & 5 & \degsixA & 1100\,0101\,1100\,0100\,0000\,0000\,0000 \\
$\langle 2 \sqcup 1\langle 2 \rangle \rangle$ & 3 & 2 & 5 & \degsixA & 1100\,1101\,1001\,0100\,0000\,0000\,0000 \\
$\langle 3 \sqcup 1\langle 1 \rangle \rangle$ & 4 & 1 & 5 & \degsixA & 1110\,1001\,1000\,0100\,0000\,0000\,0000 \\
$\langle 5 \rangle$ & 5 & 0 & 5 & \degsixA & 1110\,1101\,0101\,0100\,1000\,0000\,0000 \\
$\langle 1\langle 3 \rangle \rangle$ & 1 & 3 & 4 & \degsixA & 1100\,0101\,1100\,0000\,0000\,0000\,0000 \\
$\langle 1 \sqcup 1\langle 2 \rangle \rangle$ & 2 & 2 & 4 & \degsixA & 1100\,1101\,1001\,0000\,0000\,0000\,0000 \\
$\langle 2 \sqcup 1\langle 1 \rangle \rangle$ & 3 & 1 & 4 & \degsixA & 1110\,1001\,1000\,0000\,0000\,0000\,0000 \\
$\langle 4 \rangle$ & 4 & 0 & 4 & \degsixA & 1110\,1001\,0000\,0100\,0000\,0000\,0000 \\
$\langle 1\langle 2 \rangle \rangle$ & 1 & 2 & 3 & \degsixA & 1100\,0101\,1000\,0000\,0000\,0000\,0000 \\
$\langle 1 \sqcup 1\langle 1 \rangle \rangle$ & 2 & 1 & 3 & \degsixA & 1110\,0001\,1000\,0000\,0000\,0000\,0000 \\
$\langle 3 \rangle$ & 3 & 0 & 3 & \degsixA & 1110\,1001\,0000\,0000\,0000\,0000\,0000 \\
$\langle 1\langle 1\langle 1 \rangle \rangle \rangle$ & 2 & 1 & 3 & \degsixB & 1100\,0001\,1000\,0000\,0000\,0000\,0000 \\
$\langle 1\langle 1 \rangle \rangle$ & 1 & 1 & 2 & \degsixA & 1100\,0101\,0000\,0000\,0000\,0000\,0000 \\
$\langle 2 \rangle$ & 2 & 0 & 2 & \degsixA & 1110\,0001\,0000\,0000\,0000\,0000\,0000 \\
$\langle 1 \rangle$ & 1 & 0 & 1 & \degsixA & 1100\,0001\,0000\,0000\,0000\,0000\,0000 \\
\end{longtable}
\normalsize

%% file: tikz/deg7_triangulations_A.tex
\begin{figure}[th]
  \centering
  \begin{subfigure}{0.48\textwidth}
    \centering
  \begin{tikzpicture}[scale=\scalefactor, every node/.style={font=\tiny}]
    \draw (0,0) -- (0,1) -- (1,0) -- cycle;
    \draw (0,1) -- (0,2) -- (1,2) -- cycle;
    \draw (0,1) -- (1,0) -- (1,1) -- cycle;
    \draw (0,1) -- (1,1) -- (2,2) -- cycle;
    \draw (0,1) -- (1,2) -- (3,3) -- cycle;
    \draw (0,1) -- (2,2) -- (3,3) -- cycle;
    \draw (0,2) -- (0,3) -- (1,3) -- cycle;
    \draw (0,2) -- (1,2) -- (3,3) -- cycle;
    \draw (0,2) -- (1,3) -- (3,4) -- cycle;
    \draw (0,2) -- (2,3) -- (3,3) -- cycle;
    \draw (0,2) -- (2,3) -- (3,4) -- cycle;
    \draw (0,3) -- (0,4) -- (1,4) -- cycle;
    \draw (0,3) -- (1,3) -- (3,4) -- cycle;
    \draw (0,3) -- (1,4) -- (2,4) -- cycle;
    \draw (0,3) -- (2,4) -- (3,4) -- cycle;
    \draw (0,4) -- (0,5) -- (1,4) -- cycle;
    \draw (0,5) -- (0,6) -- (1,5) -- cycle;
    \draw (0,5) -- (1,4) -- (2,4) -- cycle;
    \draw (0,5) -- (1,5) -- (3,4) -- cycle;
    \draw (0,5) -- (2,4) -- (3,4) -- cycle;
    \draw (0,6) -- (0,7) -- (1,6) -- cycle;
    \draw (0,6) -- (1,5) -- (3,4) -- cycle;
    \draw (0,6) -- (1,6) -- (2,5) -- cycle;
    \draw (0,6) -- (2,5) -- (3,4) -- cycle;
    \draw (1,0) -- (1,1) -- (2,2) -- cycle;
    \draw (1,0) -- (2,0) -- (2,1) -- cycle;
    \draw (1,0) -- (2,1) -- (3,3) -- cycle;
    \draw (1,0) -- (2,2) -- (3,3) -- cycle;
    \draw (2,0) -- (2,1) -- (3,3) -- cycle;
    \draw (2,0) -- (3,0) -- (3,1) -- cycle;
    \draw (2,0) -- (3,1) -- (4,3) -- cycle;
    \draw (2,0) -- (3,2) -- (3,3) -- cycle;
    \draw (2,0) -- (3,2) -- (4,3) -- cycle;
    \draw (2,3) -- (3,3) -- (3,4) -- cycle;
    \draw (3,0) -- (3,1) -- (4,3) -- cycle;
    \draw (3,0) -- (4,0) -- (4,1) -- cycle;
    \draw (3,0) -- (4,1) -- (4,2) -- cycle;
    \draw (3,0) -- (4,2) -- (4,3) -- cycle;
    \draw (3,2) -- (3,3) -- (4,3) -- cycle;
    \draw (3,3) -- (3,4) -- (4,3) -- cycle;
    \draw (4,0) -- (4,1) -- (5,0) -- cycle;
    \draw (4,1) -- (4,2) -- (5,0) -- cycle;
    \draw (4,2) -- (4,3) -- (5,0) -- cycle;
    \draw (4,3) -- (5,0) -- (5,1) -- cycle;
    \draw (4,3) -- (5,1) -- (6,0) -- cycle;
    \draw (4,3) -- (5,2) -- (6,0) -- cycle;
    \draw (5,0) -- (5,1) -- (6,0) -- cycle;
    \draw (5,2) -- (6,0) -- (6,1) -- cycle;
    \draw (6,0) -- (6,1) -- (7,0) -- cycle;
    \node[circle, fill=white, draw=black, minimum size=4mm, inner sep=0pt, outer sep=0pt, font=\fontsize{4}{5}\selectfont, align=center] at (0,0) {4};
    \node[circle, fill=white, draw=black, minimum size=4mm, inner sep=0pt, outer sep=0pt, font=\fontsize{4}{5}\selectfont, align=center] at (0,1) {2};
    \node[circle, fill=white, draw=black, minimum size=4mm, inner sep=0pt, outer sep=0pt, font=\fontsize{4}{5}\selectfont, align=center] at (0,2) {8};
    \node[circle, fill=white, draw=black, minimum size=4mm, inner sep=0pt, outer sep=0pt, font=\fontsize{4}{5}\selectfont, align=center] at (0,3) {24};
    \node[circle, fill=white, draw=black, minimum size=4mm, inner sep=0pt, outer sep=0pt, font=\fontsize{4}{5}\selectfont, align=center] at (0,4) {47};
    \node[circle, fill=white, draw=black, minimum size=4mm, inner sep=0pt, outer sep=0pt, font=\fontsize{4}{5}\selectfont, align=center] at (0,5) {71};
    \node[circle, fill=white, draw=black, minimum size=4mm, inner sep=0pt, outer sep=0pt, font=\fontsize{4}{5}\selectfont, align=center] at (0,6) {102};
    \node[circle, fill=white, draw=black, minimum size=4mm, inner sep=0pt, outer sep=0pt, font=\fontsize{4}{5}\selectfont, align=center] at (0,7) {140};
    \node[circle, fill=white, draw=black, minimum size=4mm, inner sep=0pt, outer sep=0pt, font=\fontsize{4}{5}\selectfont, align=center] at (1,0) {2};
    \node[circle, fill=white, draw=black, minimum size=4mm, inner sep=0pt, outer sep=0pt, font=\fontsize{4}{5}\selectfont, align=center] at (1,1) {1};
    \node[circle, fill=white, draw=black, minimum size=4mm, inner sep=0pt, outer sep=0pt, font=\fontsize{4}{5}\selectfont, align=center] at (1,2) {3};
    \node[circle, fill=white, draw=black, minimum size=4mm, inner sep=0pt, outer sep=0pt, font=\fontsize{4}{5}\selectfont, align=center] at (1,3) {14};
    \node[circle, fill=white, draw=black, minimum size=4mm, inner sep=0pt, outer sep=0pt, font=\fontsize{4}{5}\selectfont, align=center] at (1,4) {34};
    \node[circle, fill=white, draw=black, minimum size=4mm, inner sep=0pt, outer sep=0pt, font=\fontsize{4}{5}\selectfont, align=center] at (1,5) {61};
    \node[circle, fill=white, draw=black, minimum size=4mm, inner sep=0pt, outer sep=0pt, font=\fontsize{4}{5}\selectfont, align=center] at (1,6) {96};
    \node[circle, fill=white, draw=black, minimum size=4mm, inner sep=0pt, outer sep=0pt, font=\fontsize{4}{5}\selectfont, align=center] at (2,0) {8};
    \node[circle, fill=white, draw=black, minimum size=4mm, inner sep=0pt, outer sep=0pt, font=\fontsize{4}{5}\selectfont, align=center] at (2,1) {3};
    \node[circle, fill=white, draw=black, minimum size=4mm, inner sep=0pt, outer sep=0pt, font=\fontsize{4}{5}\selectfont, align=center] at (2,2) {0};
    \node[circle, fill=white, draw=black, minimum size=4mm, inner sep=0pt, outer sep=0pt, font=\fontsize{4}{5}\selectfont, align=center] at (2,3) {6};
    \node[circle, fill=white, draw=black, minimum size=4mm, inner sep=0pt, outer sep=0pt, font=\fontsize{4}{5}\selectfont, align=center] at (2,4) {22};
    \node[circle, fill=white, draw=black, minimum size=4mm, inner sep=0pt, outer sep=0pt, font=\fontsize{4}{5}\selectfont, align=center] at (2,5) {53};
    \node[circle, fill=white, draw=black, minimum size=4mm, inner sep=0pt, outer sep=0pt, font=\fontsize{4}{5}\selectfont, align=center] at (3,0) {24};
    \node[circle, fill=white, draw=black, minimum size=4mm, inner sep=0pt, outer sep=0pt, font=\fontsize{4}{5}\selectfont, align=center] at (3,1) {14};
    \node[circle, fill=white, draw=black, minimum size=4mm, inner sep=0pt, outer sep=0pt, font=\fontsize{4}{5}\selectfont, align=center] at (3,2) {6};
    \node[circle, fill=white, draw=black, minimum size=4mm, inner sep=0pt, outer sep=0pt, font=\fontsize{4}{5}\selectfont, align=center] at (3,3) {0};
    \node[circle, fill=white, draw=black, minimum size=4mm, inner sep=0pt, outer sep=0pt, font=\fontsize{4}{5}\selectfont, align=center] at (3,4) {11};
    \node[circle, fill=white, draw=black, minimum size=4mm, inner sep=0pt, outer sep=0pt, font=\fontsize{4}{5}\selectfont, align=center] at (4,0) {47};
    \node[circle, fill=white, draw=black, minimum size=4mm, inner sep=0pt, outer sep=0pt, font=\fontsize{4}{5}\selectfont, align=center] at (4,1) {34};
    \node[circle, fill=white, draw=black, minimum size=4mm, inner sep=0pt, outer sep=0pt, font=\fontsize{4}{5}\selectfont, align=center] at (4,2) {22};
    \node[circle, fill=white, draw=black, minimum size=4mm, inner sep=0pt, outer sep=0pt, font=\fontsize{4}{5}\selectfont, align=center] at (4,3) {11};
    \node[circle, fill=white, draw=black, minimum size=4mm, inner sep=0pt, outer sep=0pt, font=\fontsize{4}{5}\selectfont, align=center] at (5,0) {71};
    \node[circle, fill=white, draw=black, minimum size=4mm, inner sep=0pt, outer sep=0pt, font=\fontsize{4}{5}\selectfont, align=center] at (5,1) {61};
    \node[circle, fill=white, draw=black, minimum size=4mm, inner sep=0pt, outer sep=0pt, font=\fontsize{4}{5}\selectfont, align=center] at (5,2) {53};
    \node[circle, fill=white, draw=black, minimum size=4mm, inner sep=0pt, outer sep=0pt, font=\fontsize{4}{5}\selectfont, align=center] at (6,0) {102};
    \node[circle, fill=white, draw=black, minimum size=4mm, inner sep=0pt, outer sep=0pt, font=\fontsize{4}{5}\selectfont, align=center] at (6,1) {96};
    \node[circle, fill=white, draw=black, minimum size=4mm, inner sep=0pt, outer sep=0pt, font=\fontsize{4}{5}\selectfont, align=center] at (7,0) {140};
  \end{tikzpicture}
    \caption{\degsevenA{} (realizes \degsevenAtypes{} types)}
    \label{fig:deg7-A}
  \end{subfigure}
  \begin{subfigure}{0.48\textwidth}
    \centering
  \begin{tikzpicture}[scale=\scalefactor, every node/.style={font=\tiny}]
    \draw (0,0) -- (0,1) -- (1,0) -- cycle;
    \draw (0,1) -- (0,2) -- (1,1) -- cycle;
    \draw (0,1) -- (1,0) -- (1,1) -- cycle;
    \draw (0,2) -- (0,3) -- (1,3) -- cycle;
    \draw (0,2) -- (1,1) -- (1,2) -- cycle;
    \draw (0,2) -- (1,2) -- (2,3) -- cycle;
    \draw (0,2) -- (1,3) -- (3,4) -- cycle;
    \draw (0,2) -- (2,3) -- (3,4) -- cycle;
    \draw (0,3) -- (0,4) -- (1,4) -- cycle;
    \draw (0,3) -- (1,3) -- (3,4) -- cycle;
    \draw (0,3) -- (1,4) -- (2,4) -- cycle;
    \draw (0,3) -- (2,4) -- (3,4) -- cycle;
    \draw (0,4) -- (0,5) -- (1,4) -- cycle;
    \draw (0,5) -- (0,6) -- (1,5) -- cycle;
    \draw (0,5) -- (1,4) -- (2,4) -- cycle;
    \draw (0,5) -- (1,5) -- (3,4) -- cycle;
    \draw (0,5) -- (2,4) -- (3,4) -- cycle;
    \draw (0,6) -- (0,7) -- (1,6) -- cycle;
    \draw (0,6) -- (1,5) -- (3,4) -- cycle;
    \draw (0,6) -- (1,6) -- (2,5) -- cycle;
    \draw (0,6) -- (2,5) -- (3,4) -- cycle;
    \draw (1,0) -- (1,1) -- (2,0) -- cycle;
    \draw (1,1) -- (1,2) -- (2,3) -- cycle;
    \draw (1,1) -- (2,0) -- (2,1) -- cycle;
    \draw (1,1) -- (2,1) -- (3,2) -- cycle;
    \draw (1,1) -- (2,2) -- (3,4) -- cycle;
    \draw (1,1) -- (2,2) -- (4,3) -- cycle;
    \draw (1,1) -- (2,3) -- (3,4) -- cycle;
    \draw (1,1) -- (3,2) -- (4,3) -- cycle;
    \draw (2,0) -- (2,1) -- (3,2) -- cycle;
    \draw (2,0) -- (3,0) -- (3,1) -- cycle;
    \draw (2,0) -- (3,1) -- (4,3) -- cycle;
    \draw (2,0) -- (3,2) -- (4,3) -- cycle;
    \draw (2,2) -- (3,3) -- (3,4) -- cycle;
    \draw (2,2) -- (3,3) -- (4,3) -- cycle;
    \draw (3,0) -- (3,1) -- (4,3) -- cycle;
    \draw (3,0) -- (4,0) -- (4,1) -- cycle;
    \draw (3,0) -- (4,1) -- (4,2) -- cycle;
    \draw (3,0) -- (4,2) -- (4,3) -- cycle;
    \draw (3,3) -- (3,4) -- (4,3) -- cycle;
    \draw (4,0) -- (4,1) -- (5,0) -- cycle;
    \draw (4,1) -- (4,2) -- (5,0) -- cycle;
    \draw (4,2) -- (4,3) -- (5,0) -- cycle;
    \draw (4,3) -- (5,0) -- (5,1) -- cycle;
    \draw (4,3) -- (5,1) -- (6,0) -- cycle;
    \draw (4,3) -- (5,2) -- (6,0) -- cycle;
    \draw (5,0) -- (5,1) -- (6,0) -- cycle;
    \draw (5,2) -- (6,0) -- (6,1) -- cycle;
    \draw (6,0) -- (6,1) -- (7,0) -- cycle;
    \node[circle, fill=white, draw=black, minimum size=4mm, inner sep=0pt, outer sep=0pt, font=\fontsize{4}{5}\selectfont, align=center] at (0,0) {19};
    \node[circle, fill=white, draw=black, minimum size=4mm, inner sep=0pt, outer sep=0pt, font=\fontsize{4}{5}\selectfont, align=center] at (0,1) {9};
    \node[circle, fill=white, draw=black, minimum size=4mm, inner sep=0pt, outer sep=0pt, font=\fontsize{4}{5}\selectfont, align=center] at (0,2) {13};
    \node[circle, fill=white, draw=black, minimum size=4mm, inner sep=0pt, outer sep=0pt, font=\fontsize{4}{5}\selectfont, align=center] at (0,3) {25};
    \node[circle, fill=white, draw=black, minimum size=4mm, inner sep=0pt, outer sep=0pt, font=\fontsize{4}{5}\selectfont, align=center] at (0,4) {44};
    \node[circle, fill=white, draw=black, minimum size=4mm, inner sep=0pt, outer sep=0pt, font=\fontsize{4}{5}\selectfont, align=center] at (0,5) {64};
    \node[circle, fill=white, draw=black, minimum size=4mm, inner sep=0pt, outer sep=0pt, font=\fontsize{4}{5}\selectfont, align=center] at (0,6) {91};
    \node[circle, fill=white, draw=black, minimum size=4mm, inner sep=0pt, outer sep=0pt, font=\fontsize{4}{5}\selectfont, align=center] at (0,7) {125};
    \node[circle, fill=white, draw=black, minimum size=4mm, inner sep=0pt, outer sep=0pt, font=\fontsize{4}{5}\selectfont, align=center] at (1,0) {9};
    \node[circle, fill=white, draw=black, minimum size=4mm, inner sep=0pt, outer sep=0pt, font=\fontsize{4}{5}\selectfont, align=center] at (1,1) {0};
    \node[circle, fill=white, draw=black, minimum size=4mm, inner sep=0pt, outer sep=0pt, font=\fontsize{4}{5}\selectfont, align=center] at (1,2) {5};
    \node[circle, fill=white, draw=black, minimum size=4mm, inner sep=0pt, outer sep=0pt, font=\fontsize{4}{5}\selectfont, align=center] at (1,3) {13};
    \node[circle, fill=white, draw=black, minimum size=4mm, inner sep=0pt, outer sep=0pt, font=\fontsize{4}{5}\selectfont, align=center] at (1,4) {29};
    \node[circle, fill=white, draw=black, minimum size=4mm, inner sep=0pt, outer sep=0pt, font=\fontsize{4}{5}\selectfont, align=center] at (1,5) {52};
    \node[circle, fill=white, draw=black, minimum size=4mm, inner sep=0pt, outer sep=0pt, font=\fontsize{4}{5}\selectfont, align=center] at (1,6) {83};
    \node[circle, fill=white, draw=black, minimum size=4mm, inner sep=0pt, outer sep=0pt, font=\fontsize{4}{5}\selectfont, align=center] at (2,0) {13};
    \node[circle, fill=white, draw=black, minimum size=4mm, inner sep=0pt, outer sep=0pt, font=\fontsize{4}{5}\selectfont, align=center] at (2,1) {5};
    \node[circle, fill=white, draw=black, minimum size=4mm, inner sep=0pt, outer sep=0pt, font=\fontsize{4}{5}\selectfont, align=center] at (2,2) {0};
    \node[circle, fill=white, draw=black, minimum size=4mm, inner sep=0pt, outer sep=0pt, font=\fontsize{4}{5}\selectfont, align=center] at (2,3) {3};
    \node[circle, fill=white, draw=black, minimum size=4mm, inner sep=0pt, outer sep=0pt, font=\fontsize{4}{5}\selectfont, align=center] at (2,4) {15};
    \node[circle, fill=white, draw=black, minimum size=4mm, inner sep=0pt, outer sep=0pt, font=\fontsize{4}{5}\selectfont, align=center] at (2,5) {42};
    \node[circle, fill=white, draw=black, minimum size=4mm, inner sep=0pt, outer sep=0pt, font=\fontsize{4}{5}\selectfont, align=center] at (3,0) {25};
    \node[circle, fill=white, draw=black, minimum size=4mm, inner sep=0pt, outer sep=0pt, font=\fontsize{4}{5}\selectfont, align=center] at (3,1) {13};
    \node[circle, fill=white, draw=black, minimum size=4mm, inner sep=0pt, outer sep=0pt, font=\fontsize{4}{5}\selectfont, align=center] at (3,2) {3};
    \node[circle, fill=white, draw=black, minimum size=4mm, inner sep=0pt, outer sep=0pt, font=\fontsize{4}{5}\selectfont, align=center] at (3,3) {1};
    \node[circle, fill=white, draw=black, minimum size=4mm, inner sep=0pt, outer sep=0pt, font=\fontsize{4}{5}\selectfont, align=center] at (3,4) {2};
    \node[circle, fill=white, draw=black, minimum size=4mm, inner sep=0pt, outer sep=0pt, font=\fontsize{4}{5}\selectfont, align=center] at (4,0) {44};
    \node[circle, fill=white, draw=black, minimum size=4mm, inner sep=0pt, outer sep=0pt, font=\fontsize{4}{5}\selectfont, align=center] at (4,1) {29};
    \node[circle, fill=white, draw=black, minimum size=4mm, inner sep=0pt, outer sep=0pt, font=\fontsize{4}{5}\selectfont, align=center] at (4,2) {15};
    \node[circle, fill=white, draw=black, minimum size=4mm, inner sep=0pt, outer sep=0pt, font=\fontsize{4}{5}\selectfont, align=center] at (4,3) {2};
    \node[circle, fill=white, draw=black, minimum size=4mm, inner sep=0pt, outer sep=0pt, font=\fontsize{4}{5}\selectfont, align=center] at (5,0) {64};
    \node[circle, fill=white, draw=black, minimum size=4mm, inner sep=0pt, outer sep=0pt, font=\fontsize{4}{5}\selectfont, align=center] at (5,1) {52};
    \node[circle, fill=white, draw=black, minimum size=4mm, inner sep=0pt, outer sep=0pt, font=\fontsize{4}{5}\selectfont, align=center] at (5,2) {42};
    \node[circle, fill=white, draw=black, minimum size=4mm, inner sep=0pt, outer sep=0pt, font=\fontsize{4}{5}\selectfont, align=center] at (6,0) {91};
    \node[circle, fill=white, draw=black, minimum size=4mm, inner sep=0pt, outer sep=0pt, font=\fontsize{4}{5}\selectfont, align=center] at (6,1) {83};
    \node[circle, fill=white, draw=black, minimum size=4mm, inner sep=0pt, outer sep=0pt, font=\fontsize{4}{5}\selectfont, align=center] at (7,0) {125};
  \end{tikzpicture}
    \caption{\degsevenB{} (realizes \degsevenBtypes{} types)}
    \label{fig:deg7-B}
  \end{subfigure}
  \caption{Two out of four regular triangulations of $7\cdot\Delta_2$ realizing all real schemes types of degree seven; see also Figure~\ref{fig:deg7-B-triangulations}.  Values at the vertices indicate lifting functions.}
  \label{fig:deg7-A-triangulations}
\end{figure}

%% file: tikz/deg7_triangulations_B.tex
\begin{figure}[th]
  \centering
  \begin{subfigure}{0.48\textwidth}
    \centering
  \begin{tikzpicture}[scale=\scalefactor, every node/.style={font=\tiny}]
    \draw (0,0) -- (0,1) -- (1,0) -- cycle;
    \draw (0,1) -- (0,2) -- (1,2) -- cycle;
    \draw (0,1) -- (1,0) -- (1,1) -- cycle;
    \draw (0,1) -- (1,1) -- (2,2) -- cycle;
    \draw (0,1) -- (1,2) -- (2,2) -- cycle;
    \draw (0,2) -- (0,3) -- (1,3) -- cycle;
    \draw (0,2) -- (1,2) -- (2,3) -- cycle;
    \draw (0,2) -- (1,3) -- (3,4) -- cycle;
    \draw (0,2) -- (2,3) -- (3,4) -- cycle;
    \draw (0,3) -- (0,4) -- (1,4) -- cycle;
    \draw (0,3) -- (1,3) -- (3,4) -- cycle;
    \draw (0,3) -- (1,4) -- (2,4) -- cycle;
    \draw (0,3) -- (2,4) -- (3,4) -- cycle;
    \draw (0,4) -- (0,5) -- (1,4) -- cycle;
    \draw (0,5) -- (0,6) -- (1,5) -- cycle;
    \draw (0,5) -- (1,4) -- (2,4) -- cycle;
    \draw (0,5) -- (1,5) -- (3,4) -- cycle;
    \draw (0,5) -- (2,4) -- (3,4) -- cycle;
    \draw (0,6) -- (0,7) -- (1,6) -- cycle;
    \draw (0,6) -- (1,5) -- (3,4) -- cycle;
    \draw (0,6) -- (1,6) -- (2,5) -- cycle;
    \draw (0,6) -- (2,5) -- (3,4) -- cycle;
    \draw (1,0) -- (1,1) -- (2,2) -- cycle;
    \draw (1,0) -- (2,0) -- (2,1) -- cycle;
    \draw (1,0) -- (2,1) -- (2,2) -- cycle;
    \draw (1,2) -- (2,2) -- (2,3) -- cycle;
    \draw (2,0) -- (2,1) -- (3,2) -- cycle;
    \draw (2,0) -- (3,0) -- (3,1) -- cycle;
    \draw (2,0) -- (3,1) -- (4,3) -- cycle;
    \draw (2,0) -- (3,2) -- (4,3) -- cycle;
    \draw (2,1) -- (2,2) -- (3,2) -- cycle;
    \draw (2,2) -- (2,3) -- (3,4) -- cycle;
    \draw (2,2) -- (3,2) -- (4,3) -- cycle;
    \draw (2,2) -- (3,3) -- (3,4) -- cycle;
    \draw (2,2) -- (3,3) -- (4,3) -- cycle;
    \draw (3,0) -- (3,1) -- (4,3) -- cycle;
    \draw (3,0) -- (4,0) -- (4,1) -- cycle;
    \draw (3,0) -- (4,1) -- (4,2) -- cycle;
    \draw (3,0) -- (4,2) -- (4,3) -- cycle;
    \draw (3,3) -- (3,4) -- (4,3) -- cycle;
    \draw (4,0) -- (4,1) -- (5,0) -- cycle;
    \draw (4,1) -- (4,2) -- (5,0) -- cycle;
    \draw (4,2) -- (4,3) -- (5,0) -- cycle;
    \draw (4,3) -- (5,0) -- (5,1) -- cycle;
    \draw (4,3) -- (5,1) -- (6,0) -- cycle;
    \draw (4,3) -- (5,2) -- (6,0) -- cycle;
    \draw (5,0) -- (5,1) -- (6,0) -- cycle;
    \draw (5,2) -- (6,0) -- (6,1) -- cycle;
    \draw (6,0) -- (6,1) -- (7,0) -- cycle;
    \node[circle, fill=white, draw=black, minimum size=4mm, inner sep=0pt, outer sep=0pt, font=\fontsize{4}{5}\selectfont, align=center] at (0,0) {8};
    \node[circle, fill=white, draw=black, minimum size=4mm, inner sep=0pt, outer sep=0pt, font=\fontsize{4}{5}\selectfont, align=center] at (0,1) {5};
    \node[circle, fill=white, draw=black, minimum size=4mm, inner sep=0pt, outer sep=0pt, font=\fontsize{4}{5}\selectfont, align=center] at (0,2) {7};
    \node[circle, fill=white, draw=black, minimum size=4mm, inner sep=0pt, outer sep=0pt, font=\fontsize{4}{5}\selectfont, align=center] at (0,3) {16};
    \node[circle, fill=white, draw=black, minimum size=4mm, inner sep=0pt, outer sep=0pt, font=\fontsize{4}{5}\selectfont, align=center] at (0,4) {32};
    \node[circle, fill=white, draw=black, minimum size=4mm, inner sep=0pt, outer sep=0pt, font=\fontsize{4}{5}\selectfont, align=center] at (0,5) {49};
    \node[circle, fill=white, draw=black, minimum size=4mm, inner sep=0pt, outer sep=0pt, font=\fontsize{4}{5}\selectfont, align=center] at (0,6) {73};
    \node[circle, fill=white, draw=black, minimum size=4mm, inner sep=0pt, outer sep=0pt, font=\fontsize{4}{5}\selectfont, align=center] at (0,7) {104};
    \node[circle, fill=white, draw=black, minimum size=4mm, inner sep=0pt, outer sep=0pt, font=\fontsize{4}{5}\selectfont, align=center] at (1,0) {5};
    \node[circle, fill=white, draw=black, minimum size=4mm, inner sep=0pt, outer sep=0pt, font=\fontsize{4}{5}\selectfont, align=center] at (1,1) {3};
    \node[circle, fill=white, draw=black, minimum size=4mm, inner sep=0pt, outer sep=0pt, font=\fontsize{4}{5}\selectfont, align=center] at (1,2) {3};
    \node[circle, fill=white, draw=black, minimum size=4mm, inner sep=0pt, outer sep=0pt, font=\fontsize{4}{5}\selectfont, align=center] at (1,3) {8};
    \node[circle, fill=white, draw=black, minimum size=4mm, inner sep=0pt, outer sep=0pt, font=\fontsize{4}{5}\selectfont, align=center] at (1,4) {21};
    \node[circle, fill=white, draw=black, minimum size=4mm, inner sep=0pt, outer sep=0pt, font=\fontsize{4}{5}\selectfont, align=center] at (1,5) {41};
    \node[circle, fill=white, draw=black, minimum size=4mm, inner sep=0pt, outer sep=0pt, font=\fontsize{4}{5}\selectfont, align=center] at (1,6) {69};
    \node[circle, fill=white, draw=black, minimum size=4mm, inner sep=0pt, outer sep=0pt, font=\fontsize{4}{5}\selectfont, align=center] at (2,0) {7};
    \node[circle, fill=white, draw=black, minimum size=4mm, inner sep=0pt, outer sep=0pt, font=\fontsize{4}{5}\selectfont, align=center] at (2,1) {3};
    \node[circle, fill=white, draw=black, minimum size=4mm, inner sep=0pt, outer sep=0pt, font=\fontsize{4}{5}\selectfont, align=center] at (2,2) {0};
    \node[circle, fill=white, draw=black, minimum size=4mm, inner sep=0pt, outer sep=0pt, font=\fontsize{4}{5}\selectfont, align=center] at (2,3) {2};
    \node[circle, fill=white, draw=black, minimum size=4mm, inner sep=0pt, outer sep=0pt, font=\fontsize{4}{5}\selectfont, align=center] at (2,4) {11};
    \node[circle, fill=white, draw=black, minimum size=4mm, inner sep=0pt, outer sep=0pt, font=\fontsize{4}{5}\selectfont, align=center] at (2,5) {35};
    \node[circle, fill=white, draw=black, minimum size=4mm, inner sep=0pt, outer sep=0pt, font=\fontsize{4}{5}\selectfont, align=center] at (3,0) {16};
    \node[circle, fill=white, draw=black, minimum size=4mm, inner sep=0pt, outer sep=0pt, font=\fontsize{4}{5}\selectfont, align=center] at (3,1) {8};
    \node[circle, fill=white, draw=black, minimum size=4mm, inner sep=0pt, outer sep=0pt, font=\fontsize{4}{5}\selectfont, align=center] at (3,2) {2};
    \node[circle, fill=white, draw=black, minimum size=4mm, inner sep=0pt, outer sep=0pt, font=\fontsize{4}{5}\selectfont, align=center] at (3,3) {1};
    \node[circle, fill=white, draw=black, minimum size=4mm, inner sep=0pt, outer sep=0pt, font=\fontsize{4}{5}\selectfont, align=center] at (3,4) {2};
    \node[circle, fill=white, draw=black, minimum size=4mm, inner sep=0pt, outer sep=0pt, font=\fontsize{4}{5}\selectfont, align=center] at (4,0) {32};
    \node[circle, fill=white, draw=black, minimum size=4mm, inner sep=0pt, outer sep=0pt, font=\fontsize{4}{5}\selectfont, align=center] at (4,1) {21};
    \node[circle, fill=white, draw=black, minimum size=4mm, inner sep=0pt, outer sep=0pt, font=\fontsize{4}{5}\selectfont, align=center] at (4,2) {11};
    \node[circle, fill=white, draw=black, minimum size=4mm, inner sep=0pt, outer sep=0pt, font=\fontsize{4}{5}\selectfont, align=center] at (4,3) {2};
    \node[circle, fill=white, draw=black, minimum size=4mm, inner sep=0pt, outer sep=0pt, font=\fontsize{4}{5}\selectfont, align=center] at (5,0) {49};
    \node[circle, fill=white, draw=black, minimum size=4mm, inner sep=0pt, outer sep=0pt, font=\fontsize{4}{5}\selectfont, align=center] at (5,1) {41};
    \node[circle, fill=white, draw=black, minimum size=4mm, inner sep=0pt, outer sep=0pt, font=\fontsize{4}{5}\selectfont, align=center] at (5,2) {35};
    \node[circle, fill=white, draw=black, minimum size=4mm, inner sep=0pt, outer sep=0pt, font=\fontsize{4}{5}\selectfont, align=center] at (6,0) {73};
    \node[circle, fill=white, draw=black, minimum size=4mm, inner sep=0pt, outer sep=0pt, font=\fontsize{4}{5}\selectfont, align=center] at (6,1) {69};
    \node[circle, fill=white, draw=black, minimum size=4mm, inner sep=0pt, outer sep=0pt, font=\fontsize{4}{5}\selectfont, align=center] at (7,0) {104};
  \end{tikzpicture}
    \caption{\degsevenC{} (realizes \degsevenCtypes{} types)}
    \label{fig:deg7-C}
  \end{subfigure}
  \begin{subfigure}{0.48\textwidth}
    \centering
  \begin{tikzpicture}[scale=\scalefactor, every node/.style={font=\tiny}]
    \draw (0,0) -- (0,1) -- (1,0) -- cycle;
    \draw (0,1) -- (0,2) -- (1,1) -- cycle;
    \draw (0,1) -- (1,0) -- (1,1) -- cycle;
    \draw (0,2) -- (0,3) -- (1,2) -- cycle;
    \draw (0,2) -- (1,1) -- (1,2) -- cycle;
    \draw (0,3) -- (0,4) -- (1,3) -- cycle;
    \draw (0,3) -- (1,2) -- (1,3) -- cycle;
    \draw (0,4) -- (0,5) -- (1,4) -- cycle;
    \draw (0,4) -- (1,3) -- (1,4) -- cycle;
    \draw (0,5) -- (0,6) -- (1,5) -- cycle;
    \draw (0,5) -- (1,4) -- (1,5) -- cycle;
    \draw (0,6) -- (0,7) -- (1,6) -- cycle;
    \draw (0,6) -- (1,5) -- (1,6) -- cycle;
    \draw (1,0) -- (1,1) -- (2,0) -- cycle;
    \draw (1,1) -- (1,2) -- (2,1) -- cycle;
    \draw (1,1) -- (2,0) -- (2,1) -- cycle;
    \draw (1,2) -- (1,3) -- (2,2) -- cycle;
    \draw (1,2) -- (2,1) -- (2,2) -- cycle;
    \draw (1,3) -- (1,4) -- (2,3) -- cycle;
    \draw (1,3) -- (2,2) -- (2,3) -- cycle;
    \draw (1,4) -- (1,5) -- (2,4) -- cycle;
    \draw (1,4) -- (2,3) -- (2,4) -- cycle;
    \draw (1,5) -- (1,6) -- (2,5) -- cycle;
    \draw (1,5) -- (2,4) -- (2,5) -- cycle;
    \draw (2,0) -- (2,1) -- (3,0) -- cycle;
    \draw (2,1) -- (2,2) -- (3,1) -- cycle;
    \draw (2,1) -- (3,0) -- (3,1) -- cycle;
    \draw (2,2) -- (2,3) -- (3,2) -- cycle;
    \draw (2,2) -- (3,1) -- (3,2) -- cycle;
    \draw (2,3) -- (2,4) -- (3,3) -- cycle;
    \draw (2,3) -- (3,2) -- (3,3) -- cycle;
    \draw (2,4) -- (2,5) -- (3,4) -- cycle;
    \draw (2,4) -- (3,3) -- (3,4) -- cycle;
    \draw (3,0) -- (3,1) -- (4,0) -- cycle;
    \draw (3,1) -- (3,2) -- (4,1) -- cycle;
    \draw (3,1) -- (4,0) -- (4,1) -- cycle;
    \draw (3,2) -- (3,3) -- (4,2) -- cycle;
    \draw (3,2) -- (4,1) -- (4,2) -- cycle;
    \draw (3,3) -- (3,4) -- (4,3) -- cycle;
    \draw (3,3) -- (4,2) -- (4,3) -- cycle;
    \draw (4,0) -- (4,1) -- (5,0) -- cycle;
    \draw (4,1) -- (4,2) -- (5,1) -- cycle;
    \draw (4,1) -- (5,0) -- (5,1) -- cycle;
    \draw (4,2) -- (4,3) -- (5,2) -- cycle;
    \draw (4,2) -- (5,1) -- (5,2) -- cycle;
    \draw (5,0) -- (5,1) -- (6,0) -- cycle;
    \draw (5,1) -- (5,2) -- (6,1) -- cycle;
    \draw (5,1) -- (6,0) -- (6,1) -- cycle;
    \draw (6,0) -- (6,1) -- (7,0) -- cycle;
    \node[circle, fill=white, draw=black, minimum size=4mm, inner sep=0pt, outer sep=0pt, font=\fontsize{4}{5}\selectfont, align=center] at (0,0) {16};
    \node[circle, fill=white, draw=black, minimum size=4mm, inner sep=0pt, outer sep=0pt, font=\fontsize{4}{5}\selectfont, align=center] at (0,1) {10};
    \node[circle, fill=white, draw=black, minimum size=4mm, inner sep=0pt, outer sep=0pt, font=\fontsize{4}{5}\selectfont, align=center] at (0,2) {6};
    \node[circle, fill=white, draw=black, minimum size=4mm, inner sep=0pt, outer sep=0pt, font=\fontsize{4}{5}\selectfont, align=center] at (0,3) {4};
    \node[circle, fill=white, draw=black, minimum size=4mm, inner sep=0pt, outer sep=0pt, font=\fontsize{4}{5}\selectfont, align=center] at (0,4) {4};
    \node[circle, fill=white, draw=black, minimum size=4mm, inner sep=0pt, outer sep=0pt, font=\fontsize{4}{5}\selectfont, align=center] at (0,5) {6};
    \node[circle, fill=white, draw=black, minimum size=4mm, inner sep=0pt, outer sep=0pt, font=\fontsize{4}{5}\selectfont, align=center] at (0,6) {10};
    \node[circle, fill=white, draw=black, minimum size=4mm, inner sep=0pt, outer sep=0pt, font=\fontsize{4}{5}\selectfont, align=center] at (0,7) {16};
    \node[circle, fill=white, draw=black, minimum size=4mm, inner sep=0pt, outer sep=0pt, font=\fontsize{4}{5}\selectfont, align=center] at (1,0) {10};
    \node[circle, fill=white, draw=black, minimum size=4mm, inner sep=0pt, outer sep=0pt, font=\fontsize{4}{5}\selectfont, align=center] at (1,1) {5};
    \node[circle, fill=white, draw=black, minimum size=4mm, inner sep=0pt, outer sep=0pt, font=\fontsize{4}{5}\selectfont, align=center] at (1,2) {2};
    \node[circle, fill=white, draw=black, minimum size=4mm, inner sep=0pt, outer sep=0pt, font=\fontsize{4}{5}\selectfont, align=center] at (1,3) {1};
    \node[circle, fill=white, draw=black, minimum size=4mm, inner sep=0pt, outer sep=0pt, font=\fontsize{4}{5}\selectfont, align=center] at (1,4) {2};
    \node[circle, fill=white, draw=black, minimum size=4mm, inner sep=0pt, outer sep=0pt, font=\fontsize{4}{5}\selectfont, align=center] at (1,5) {5};
    \node[circle, fill=white, draw=black, minimum size=4mm, inner sep=0pt, outer sep=0pt, font=\fontsize{4}{5}\selectfont, align=center] at (1,6) {10};
    \node[circle, fill=white, draw=black, minimum size=4mm, inner sep=0pt, outer sep=0pt, font=\fontsize{4}{5}\selectfont, align=center] at (2,0) {6};
    \node[circle, fill=white, draw=black, minimum size=4mm, inner sep=0pt, outer sep=0pt, font=\fontsize{4}{5}\selectfont, align=center] at (2,1) {2};
    \node[circle, fill=white, draw=black, minimum size=4mm, inner sep=0pt, outer sep=0pt, font=\fontsize{4}{5}\selectfont, align=center] at (2,2) {0};
    \node[circle, fill=white, draw=black, minimum size=4mm, inner sep=0pt, outer sep=0pt, font=\fontsize{4}{5}\selectfont, align=center] at (2,3) {0};
    \node[circle, fill=white, draw=black, minimum size=4mm, inner sep=0pt, outer sep=0pt, font=\fontsize{4}{5}\selectfont, align=center] at (2,4) {2};
    \node[circle, fill=white, draw=black, minimum size=4mm, inner sep=0pt, outer sep=0pt, font=\fontsize{4}{5}\selectfont, align=center] at (2,5) {6};
    \node[circle, fill=white, draw=black, minimum size=4mm, inner sep=0pt, outer sep=0pt, font=\fontsize{4}{5}\selectfont, align=center] at (3,0) {4};
    \node[circle, fill=white, draw=black, minimum size=4mm, inner sep=0pt, outer sep=0pt, font=\fontsize{4}{5}\selectfont, align=center] at (3,1) {1};
    \node[circle, fill=white, draw=black, minimum size=4mm, inner sep=0pt, outer sep=0pt, font=\fontsize{4}{5}\selectfont, align=center] at (3,2) {0};
    \node[circle, fill=white, draw=black, minimum size=4mm, inner sep=0pt, outer sep=0pt, font=\fontsize{4}{5}\selectfont, align=center] at (3,3) {1};
    \node[circle, fill=white, draw=black, minimum size=4mm, inner sep=0pt, outer sep=0pt, font=\fontsize{4}{5}\selectfont, align=center] at (3,4) {4};
    \node[circle, fill=white, draw=black, minimum size=4mm, inner sep=0pt, outer sep=0pt, font=\fontsize{4}{5}\selectfont, align=center] at (4,0) {4};
    \node[circle, fill=white, draw=black, minimum size=4mm, inner sep=0pt, outer sep=0pt, font=\fontsize{4}{5}\selectfont, align=center] at (4,1) {2};
    \node[circle, fill=white, draw=black, minimum size=4mm, inner sep=0pt, outer sep=0pt, font=\fontsize{4}{5}\selectfont, align=center] at (4,2) {2};
    \node[circle, fill=white, draw=black, minimum size=4mm, inner sep=0pt, outer sep=0pt, font=\fontsize{4}{5}\selectfont, align=center] at (4,3) {4};
    \node[circle, fill=white, draw=black, minimum size=4mm, inner sep=0pt, outer sep=0pt, font=\fontsize{4}{5}\selectfont, align=center] at (5,0) {6};
    \node[circle, fill=white, draw=black, minimum size=4mm, inner sep=0pt, outer sep=0pt, font=\fontsize{4}{5}\selectfont, align=center] at (5,1) {5};
    \node[circle, fill=white, draw=black, minimum size=4mm, inner sep=0pt, outer sep=0pt, font=\fontsize{4}{5}\selectfont, align=center] at (5,2) {6};
    \node[circle, fill=white, draw=black, minimum size=4mm, inner sep=0pt, outer sep=0pt, font=\fontsize{4}{5}\selectfont, align=center] at (6,0) {10};
    \node[circle, fill=white, draw=black, minimum size=4mm, inner sep=0pt, outer sep=0pt, font=\fontsize{4}{5}\selectfont, align=center] at (6,1) {10};
    \node[circle, fill=white, draw=black, minimum size=4mm, inner sep=0pt, outer sep=0pt, font=\fontsize{4}{5}\selectfont, align=center] at (7,0) {16};
  \end{tikzpicture}
    \caption{\degsevenD{} (realizes \degsevenDtypes{} types)}
    \label{fig:deg7-D}
  \end{subfigure}
  \caption{Two out of four regular triangulations of $7\cdot\Delta_2$ realizing all real schemes types of degree seven; see also Figure~\ref{fig:deg7-A-triangulations}.  Values at the vertices indicate lifting functions.}
  \label{fig:deg7-B-triangulations}
\end{figure}

%% file: tables/proof_deg7.tex
\small
\begin{longtable}{lrrrlr}
\caption{Proof of support for degree 7}
\label{tab:deg-7} \\
\toprule
Real scheme & \multicolumn{1}{c}{$p$} & \multicolumn{1}{c}{$n$} & \multicolumn{1}{c}{$p+n$} & \multicolumn{1}{c}{$\cT$} & \multicolumn{1}{c}{$\sigma$} \\
\midrule
\endfirsthead
\toprule
Real scheme & \multicolumn{1}{c}{$p$} & \multicolumn{1}{c}{$n$} & \multicolumn{1}{c}{$p+n$} & \multicolumn{1}{c}{$\cT$} & \multicolumn{1}{c}{$\sigma$} \\
\midrule
\endhead
\midrule
\endfoot
\bottomrule
\endlastfoot
$\langle J \sqcup 1 \sqcup 1\langle 13 \rangle \rangle$ & 2 & 13 & 15 & \degsevenA & 1101\,0110\,1001\,0000\,0100\,0110\,0100\,0111\,0000 \\
$\langle J \sqcup 2 \sqcup 1\langle 12 \rangle \rangle$ & 3 & 12 & 15 & \degsevenB & 1100\,0010\,1001\,0000\,0110\,0011\,0000\,0001\,0000 \\
$\langle J \sqcup 3 \sqcup 1\langle 11 \rangle \rangle$ & 4 & 11 & 15 & \degsevenC & 1110\,0010\,1000\,0001\,0110\,0001\,0000\,0001\,0000 \\
$\langle J \sqcup 4 \sqcup 1\langle 10 \rangle \rangle$ & 5 & 10 & 15 & \degsevenA & 1100\,0010\,1001\,0000\,0110\,0011\,0000\,0001\,0000 \\
$\langle J \sqcup 5 \sqcup 1\langle 9 \rangle \rangle$ & 6 & 9 & 15 & \degsevenA & 1101\,0110\,1001\,0000\,0100\,0000\,0100\,0100\,0000 \\
$\langle J \sqcup 6 \sqcup 1\langle 8 \rangle \rangle$ & 7 & 8 & 15 & \degsevenB & 1101\,0110\,1011\,0000\,1100\,0000\,0100\,0100\,0000 \\
$\langle J \sqcup 7 \sqcup 1\langle 7 \rangle \rangle$ & 8 & 7 & 15 & \degsevenA & 1101\,0110\,1001\,0001\,1100\,0110\,0100\,0010\,0000 \\
$\langle J \sqcup 8 \sqcup 1\langle 6 \rangle \rangle$ & 9 & 6 & 15 & \degsevenA & 1100\,0010\,1001\,0000\,0110\,0101\,0000\,0010\,0000 \\
$\langle J \sqcup 9 \sqcup 1\langle 5 \rangle \rangle$ & 10 & 5 & 15 & \degsevenA & 1100\,0010\,1001\,0001\,1110\,0011\,0000\,0100\,0000 \\
$\langle J \sqcup 10 \sqcup 1\langle 4 \rangle \rangle$ & 11 & 4 & 15 & \degsevenB & 1100\,0000\,1010\,0000\,1100\,0000\,0100\,0100\,0000 \\
$\langle J \sqcup 11 \sqcup 1\langle 3 \rangle \rangle$ & 12 & 3 & 15 & \degsevenA & 1110\,0010\,1010\,0000\,0100\,0000\,0000\,0100\,0000 \\
$\langle J \sqcup 12 \sqcup 1\langle 2 \rangle \rangle$ & 13 & 2 & 15 & \degsevenA & 1101\,0100\,1000\,0000\,0110\,0101\,0000\,0010\,0000 \\
$\langle J \sqcup 13 \sqcup 1\langle 1 \rangle \rangle$ & 14 & 1 & 15 & \degsevenA & 1101\,0100\,1000\,0001\,1110\,0011\,0000\,0100\,0000 \\
$\langle J \sqcup 15 \rangle$ & 15 & 0 & 15 & \degsevenA & 1111\,0110\,1010\,0001\,1110\,0011\,0100\,0100\,0000 \\
$\langle J \sqcup 1\langle 13 \rangle \rangle$ & 1 & 13 & 14 & \degsevenA & 1101\,0110\,1101\,0000\,0000\,0110\,0100\,0111\,0000 \\
$\langle J \sqcup 1 \sqcup 1\langle 12 \rangle \rangle$ & 2 & 12 & 14 & \degsevenA & 1101\,0110\,1001\,0000\,0100\,0110\,0100\,0110\,0000 \\
$\langle J \sqcup 2 \sqcup 1\langle 11 \rangle \rangle$ & 3 & 11 & 14 & \degsevenB & 1100\,0010\,1001\,0000\,0110\,0011\,0000\,0000\,0000 \\
$\langle J \sqcup 3 \sqcup 1\langle 10 \rangle \rangle$ & 4 & 10 & 14 & \degsevenA & 1100\,0010\,1101\,0000\,0010\,0011\,0000\,0001\,0000 \\
$\langle J \sqcup 4 \sqcup 1\langle 9 \rangle \rangle$ & 5 & 9 & 14 & \degsevenA & 1100\,0010\,1001\,0000\,0110\,0011\,0000\,0000\,0000 \\
$\langle J \sqcup 5 \sqcup 1\langle 8 \rangle \rangle$ & 6 & 8 & 14 & \degsevenA & 1101\,0110\,1000\,0000\,0100\,0000\,0100\,0100\,0000 \\
$\langle J \sqcup 6 \sqcup 1\langle 7 \rangle \rangle$ & 7 & 7 & 14 & \degsevenA & 1101\,0110\,1001\,0001\,1100\,0000\,0100\,0000\,0000 \\
$\langle J \sqcup 7 \sqcup 1\langle 6 \rangle \rangle$ & 8 & 6 & 14 & \degsevenA & 1100\,0010\,1101\,0000\,0010\,0101\,0000\,0010\,0000 \\
$\langle J \sqcup 8 \sqcup 1\langle 5 \rangle \rangle$ & 9 & 5 & 14 & \degsevenA & 1101\,0100\,1000\,0000\,0110\,0011\,0000\,0000\,0000 \\
$\langle J \sqcup 9 \sqcup 1\langle 4 \rangle \rangle$ & 10 & 4 & 14 & \degsevenA & 1111\,0110\,1010\,0000\,0110\,0011\,0100\,0000\,0000 \\
$\langle J \sqcup 10 \sqcup 1\langle 3 \rangle \rangle$ & 11 & 3 & 14 & \degsevenA & 1100\,0000\,1000\,0001\,1100\,0000\,0100\,0000\,0000 \\
$\langle J \sqcup 11 \sqcup 1\langle 2 \rangle \rangle$ & 12 & 2 & 14 & \degsevenA & 1110\,0010\,1010\,0000\,1100\,0000\,0000\,0100\,0000 \\
$\langle J \sqcup 12 \sqcup 1\langle 1 \rangle \rangle$ & 13 & 1 & 14 & \degsevenA & 1110\,0010\,1010\,0001\,1100\,0000\,0000\,0000\,0000 \\
$\langle J \sqcup 14 \rangle$ & 14 & 0 & 14 & \degsevenA & 1101\,0100\,1000\,0000\,0110\,0000\,0000\,0100\,0000 \\
$\langle J \sqcup 1\langle 12 \rangle \rangle$ & 1 & 12 & 13 & \degsevenA & 1101\,0110\,1101\,0000\,0000\,0110\,0100\,0110\,0000 \\
$\langle J \sqcup 1 \sqcup 1\langle 11 \rangle \rangle$ & 2 & 11 & 13 & \degsevenA & 1101\,0110\,1001\,0000\,0100\,0100\,0100\,0110\,0000 \\
$\langle J \sqcup 2 \sqcup 1\langle 10 \rangle \rangle$ & 3 & 10 & 13 & \degsevenA & 1101\,0110\,1001\,0000\,0100\,0110\,0100\,0100\,0000 \\
$\langle J \sqcup 3 \sqcup 1\langle 9 \rangle \rangle$ & 4 & 9 & 13 & \degsevenA & 1100\,0010\,1101\,0000\,0010\,0011\,0000\,0000\,0000 \\
$\langle J \sqcup 4 \sqcup 1\langle 8 \rangle \rangle$ & 5 & 8 & 13 & \degsevenA & 1100\,0010\,1001\,0000\,0110\,0010\,0000\,0000\,0000 \\
$\langle J \sqcup 5 \sqcup 1\langle 7 \rangle \rangle$ & 6 & 7 & 13 & \degsevenA & 1101\,0110\,1101\,0001\,1000\,0000\,0100\,0000\,0000 \\
$\langle J \sqcup 6 \sqcup 1\langle 6 \rangle \rangle$ & 7 & 6 & 13 & \degsevenA & 1101\,0110\,1000\,0001\,1100\,0000\,0100\,0000\,0000 \\
$\langle J \sqcup 7 \sqcup 1\langle 5 \rangle \rangle$ & 8 & 5 & 13 & \degsevenA & 1101\,0100\,1100\,0000\,0010\,0011\,0000\,0000\,0000 \\
$\langle J \sqcup 8 \sqcup 1\langle 4 \rangle \rangle$ & 9 & 4 & 13 & \degsevenA & 1101\,0100\,1000\,0000\,0110\,0010\,0000\,0000\,0000 \\
$\langle J \sqcup 9 \sqcup 1\langle 3 \rangle \rangle$ & 10 & 3 & 13 & \degsevenA & 1100\,0000\,1100\,0001\,1000\,0000\,0100\,0000\,0000 \\
$\langle J \sqcup 10 \sqcup 1\langle 2 \rangle \rangle$ & 11 & 2 & 13 & \degsevenA & 1100\,0000\,1010\,0001\,1100\,0000\,0100\,0000\,0000 \\
$\langle J \sqcup 11 \sqcup 1\langle 1 \rangle \rangle$ & 12 & 1 & 13 & \degsevenA & 1110\,0010\,1110\,0001\,1000\,0000\,0000\,0000\,0000 \\
$\langle J \sqcup 13 \rangle$ & 13 & 0 & 13 & \degsevenA & 1101\,0100\,1000\,0001\,1110\,0000\,0000\,0000\,0000 \\
$\langle J \sqcup 1\langle 11 \rangle \rangle$ & 1 & 11 & 12 & \degsevenA & 1101\,0110\,1101\,0000\,0000\,0100\,0100\,0110\,0000 \\
$\langle J \sqcup 1 \sqcup 1\langle 10 \rangle \rangle$ & 2 & 10 & 12 & \degsevenA & 1101\,0110\,1101\,0000\,0000\,0110\,0100\,0100\,0000 \\
$\langle J \sqcup 2 \sqcup 1\langle 9 \rangle \rangle$ & 3 & 9 & 12 & \degsevenA & 1100\,0010\,1001\,0000\,0010\,0011\,0000\,0000\,0000 \\
$\langle J \sqcup 3 \sqcup 1\langle 8 \rangle \rangle$ & 4 & 8 & 12 & \degsevenA & 1100\,0010\,1101\,0000\,0010\,0010\,0000\,0000\,0000 \\
$\langle J \sqcup 4 \sqcup 1\langle 7 \rangle \rangle$ & 5 & 7 & 12 & \degsevenA & 1100\,0010\,1001\,0000\,0110\,0000\,0000\,0000\,0000 \\
$\langle J \sqcup 5 \sqcup 1\langle 6 \rangle \rangle$ & 6 & 6 & 12 & \degsevenA & 1100\,0010\,1001\,0000\,0110\,0101\,0000\,0000\,0000 \\
$\langle J \sqcup 6 \sqcup 1\langle 5 \rangle \rangle$ & 7 & 5 & 12 & \degsevenA & 1101\,0100\,1000\,0000\,0010\,0011\,0000\,0000\,0000 \\
$\langle J \sqcup 7 \sqcup 1\langle 4 \rangle \rangle$ & 8 & 4 & 12 & \degsevenA & 1101\,0100\,1100\,0000\,0010\,0010\,0000\,0000\,0000 \\
$\langle J \sqcup 8 \sqcup 1\langle 3 \rangle \rangle$ & 9 & 3 & 12 & \degsevenA & 1101\,0100\,1000\,0000\,0110\,0000\,0000\,0000\,0000 \\
$\langle J \sqcup 9 \sqcup 1\langle 2 \rangle \rangle$ & 10 & 2 & 12 & \degsevenA & 1101\,0100\,1000\,0000\,0110\,0101\,0000\,0000\,0000 \\
$\langle J \sqcup 10 \sqcup 1\langle 1 \rangle \rangle$ & 11 & 1 & 12 & \degsevenA & 1110\,0010\,1110\,0001\,0000\,0000\,0000\,0000\,0000 \\
$\langle J \sqcup 12 \rangle$ & 12 & 0 & 12 & \degsevenA & 1110\,0010\,1010\,0001\,1000\,0000\,0000\,0000\,0000 \\
$\langle J \sqcup 1\langle 10 \rangle \rangle$ & 1 & 10 & 11 & \degsevenA & 1101\,0110\,1100\,0000\,0000\,0100\,0100\,0110\,0000 \\
$\langle J \sqcup 1 \sqcup 1\langle 9 \rangle \rangle$ & 2 & 9 & 11 & \degsevenA & 1100\,0010\,1011\,0000\,0010\,0011\,0000\,0000\,0000 \\
$\langle J \sqcup 2 \sqcup 1\langle 8 \rangle \rangle$ & 3 & 8 & 11 & \degsevenA & 1100\,0010\,1001\,0000\,0010\,0010\,0000\,0000\,0000 \\
$\langle J \sqcup 3 \sqcup 1\langle 7 \rangle \rangle$ & 4 & 7 & 11 & \degsevenA & 1100\,0010\,1101\,0000\,0010\,0000\,0000\,0000\,0000 \\
$\langle J \sqcup 4 \sqcup 1\langle 6 \rangle \rangle$ & 5 & 6 & 11 & \degsevenA & 1100\,0010\,1001\,0000\,0100\,0000\,0000\,0000\,0000 \\
$\langle J \sqcup 5 \sqcup 1\langle 5 \rangle \rangle$ & 6 & 5 & 11 & \degsevenA & 1100\,0010\,1001\,0000\,0110\,0100\,0000\,0000\,0000 \\
$\langle J \sqcup 6 \sqcup 1\langle 4 \rangle \rangle$ & 7 & 4 & 11 & \degsevenA & 1101\,0100\,1000\,0000\,0010\,0010\,0000\,0000\,0000 \\
$\langle J \sqcup 7 \sqcup 1\langle 3 \rangle \rangle$ & 8 & 3 & 11 & \degsevenA & 1101\,0100\,1100\,0000\,0010\,0000\,0000\,0000\,0000 \\
$\langle J \sqcup 8 \sqcup 1\langle 2 \rangle \rangle$ & 9 & 2 & 11 & \degsevenA & 1101\,0100\,1000\,0000\,0100\,0000\,0000\,0000\,0000 \\
$\langle J \sqcup 9 \sqcup 1\langle 1 \rangle \rangle$ & 10 & 1 & 11 & \degsevenA & 1110\,0010\,1100\,0001\,0000\,0000\,0000\,0000\,0000 \\
$\langle J \sqcup 11 \rangle$ & 11 & 0 & 11 & \degsevenA & 1110\,0010\,1010\,0001\,0000\,0000\,0000\,0000\,0000 \\
$\langle J \sqcup 1\langle 9 \rangle \rangle$ & 1 & 9 & 10 & \degsevenA & 1100\,0010\,1011\,0000\,1010\,0011\,0000\,0000\,0000 \\
$\langle J \sqcup 1 \sqcup 1\langle 8 \rangle \rangle$ & 2 & 8 & 10 & \degsevenA & 1100\,0010\,1011\,0000\,0010\,0010\,0000\,0000\,0000 \\
$\langle J \sqcup 2 \sqcup 1\langle 7 \rangle \rangle$ & 3 & 7 & 10 & \degsevenA & 1100\,0010\,1001\,0000\,0010\,0000\,0000\,0000\,0000 \\
$\langle J \sqcup 3 \sqcup 1\langle 6 \rangle \rangle$ & 4 & 6 & 10 & \degsevenA & 1100\,0010\,1101\,0000\,0000\,0000\,0000\,0000\,0000 \\
$\langle J \sqcup 4 \sqcup 1\langle 5 \rangle \rangle$ & 5 & 5 & 10 & \degsevenA & 1100\,0010\,1000\,0000\,0100\,0000\,0000\,0000\,0000 \\
$\langle J \sqcup 5 \sqcup 1\langle 4 \rangle \rangle$ & 6 & 4 & 10 & \degsevenA & 1101\,0010\,1000\,0000\,0110\,0000\,0000\,0000\,0000 \\
$\langle J \sqcup 6 \sqcup 1\langle 3 \rangle \rangle$ & 7 & 3 & 10 & \degsevenA & 1101\,0100\,1000\,0000\,0010\,0000\,0000\,0000\,0000 \\
$\langle J \sqcup 7 \sqcup 1\langle 2 \rangle \rangle$ & 8 & 2 & 10 & \degsevenA & 1101\,0100\,1100\,0000\,0000\,0000\,0000\,0000\,0000 \\
$\langle J \sqcup 8 \sqcup 1\langle 1 \rangle \rangle$ & 9 & 1 & 10 & \degsevenA & 1110\,0000\,1100\,0001\,0000\,0000\,0000\,0000\,0000 \\
$\langle J \sqcup 10 \rangle$ & 10 & 0 & 10 & \degsevenA & 1110\,0010\,1110\,0000\,0000\,0000\,0000\,0000\,0000 \\
$\langle J \sqcup 1\langle 8 \rangle \rangle$ & 1 & 8 & 9 & \degsevenA & 1100\,0010\,1011\,0000\,1010\,0010\,0000\,0000\,0000 \\
$\langle J \sqcup 1 \sqcup 1\langle 7 \rangle \rangle$ & 2 & 7 & 9 & \degsevenA & 1100\,0010\,1011\,0000\,0010\,0000\,0000\,0000\,0000 \\
$\langle J \sqcup 2 \sqcup 1\langle 6 \rangle \rangle$ & 3 & 6 & 9 & \degsevenA & 1100\,0010\,1001\,0000\,0000\,0000\,0000\,0000\,0000 \\
$\langle J \sqcup 3 \sqcup 1\langle 5 \rangle \rangle$ & 4 & 5 & 9 & \degsevenA & 1100\,0010\,1100\,0000\,0000\,0000\,0000\,0000\,0000 \\
$\langle J \sqcup 4 \sqcup 1\langle 4 \rangle \rangle$ & 5 & 4 & 9 & \degsevenA & 1100\,0000\,1000\,0000\,0100\,0000\,0000\,0000\,0000 \\
$\langle J \sqcup 5 \sqcup 1\langle 3 \rangle \rangle$ & 6 & 3 & 9 & \degsevenA & 1101\,0010\,1000\,0000\,0100\,0000\,0000\,0000\,0000 \\
$\langle J \sqcup 6 \sqcup 1\langle 2 \rangle \rangle$ & 7 & 2 & 9 & \degsevenA & 1101\,0100\,1000\,0000\,0000\,0000\,0000\,0000\,0000 \\
$\langle J \sqcup 7 \sqcup 1\langle 1 \rangle \rangle$ & 8 & 1 & 9 & \degsevenA & 1110\,0100\,1100\,0001\,0000\,0000\,0000\,0000\,0000 \\
$\langle J \sqcup 9 \rangle$ & 9 & 0 & 9 & \degsevenA & 1110\,0010\,1100\,0000\,0000\,0000\,0000\,0000\,0000 \\
$\langle J \sqcup 1\langle 7 \rangle \rangle$ & 1 & 7 & 8 & \degsevenA & 1100\,0010\,1011\,0000\,1010\,0000\,0000\,0000\,0000 \\
$\langle J \sqcup 1 \sqcup 1\langle 6 \rangle \rangle$ & 2 & 6 & 8 & \degsevenA & 1100\,0010\,1011\,0000\,0000\,0000\,0000\,0000\,0000 \\
$\langle J \sqcup 2 \sqcup 1\langle 5 \rangle \rangle$ & 3 & 5 & 8 & \degsevenA & 1100\,0010\,1000\,0000\,0000\,0000\,0000\,0000\,0000 \\
$\langle J \sqcup 3 \sqcup 1\langle 4 \rangle \rangle$ & 4 & 4 & 8 & \degsevenA & 1100\,0000\,1100\,0000\,0000\,0000\,0000\,0000\,0000 \\
$\langle J \sqcup 4 \sqcup 1\langle 3 \rangle \rangle$ & 5 & 3 & 8 & \degsevenA & 1101\,0010\,1100\,0000\,0000\,0000\,0000\,0000\,0000 \\
$\langle J \sqcup 5 \sqcup 1\langle 2 \rangle \rangle$ & 6 & 2 & 8 & \degsevenA & 1101\,0110\,1000\,0000\,0000\,0000\,0000\,0000\,0000 \\
$\langle J \sqcup 6 \sqcup 1\langle 1 \rangle \rangle$ & 7 & 1 & 8 & \degsevenA & 1111\,0010\,1100\,0001\,0000\,0000\,0000\,0000\,0000 \\
$\langle J \sqcup 8 \rangle$ & 8 & 0 & 8 & \degsevenA & 1110\,0010\,1000\,0000\,0000\,0000\,0000\,0000\,0000 \\
$\langle J \sqcup 1\langle 6 \rangle \rangle$ & 1 & 6 & 7 & \degsevenA & 1100\,0010\,1011\,0000\,1000\,0000\,0000\,0000\,0000 \\
$\langle J \sqcup 1 \sqcup 1\langle 5 \rangle \rangle$ & 2 & 5 & 7 & \degsevenA & 1100\,0010\,1010\,0000\,0000\,0000\,0000\,0000\,0000 \\
$\langle J \sqcup 2 \sqcup 1\langle 4 \rangle \rangle$ & 3 & 4 & 7 & \degsevenA & 1100\,0000\,1000\,0000\,0000\,0000\,0000\,0000\,0000 \\
$\langle J \sqcup 3 \sqcup 1\langle 3 \rangle \rangle$ & 4 & 3 & 7 & \degsevenA & 1101\,0010\,1000\,0000\,0000\,0000\,0000\,0000\,0000 \\
$\langle J \sqcup 4 \sqcup 1\langle 2 \rangle \rangle$ & 5 & 2 & 7 & \degsevenA & 1101\,0000\,1100\,0000\,0000\,0000\,0000\,0000\,0000 \\
$\langle J \sqcup 5 \sqcup 1\langle 1 \rangle \rangle$ & 6 & 1 & 7 & \degsevenA & 1111\,0000\,1100\,0001\,0000\,0000\,0000\,0000\,0000 \\
$\langle J \sqcup 7 \rangle$ & 7 & 0 & 7 & \degsevenA & 1110\,0000\,1000\,0000\,0000\,0000\,0000\,0000\,0000 \\
$\langle J \sqcup 1\langle 5 \rangle \rangle$ & 1 & 5 & 6 & \degsevenA & 1100\,0010\,1010\,0000\,1000\,0000\,0000\,0000\,0000 \\
$\langle J \sqcup 1 \sqcup 1\langle 4 \rangle \rangle$ & 2 & 4 & 6 & \degsevenA & 1100\,0000\,1010\,0000\,0000\,0000\,0000\,0000\,0000 \\
$\langle J \sqcup 2 \sqcup 1\langle 3 \rangle \rangle$ & 3 & 3 & 6 & \degsevenA & 1101\,0010\,1010\,0000\,0000\,0000\,0000\,0000\,0000 \\
$\langle J \sqcup 3 \sqcup 1\langle 2 \rangle \rangle$ & 4 & 2 & 6 & \degsevenA & 1101\,0000\,1000\,0000\,0000\,0000\,0000\,0000\,0000 \\
$\langle J \sqcup 4 \sqcup 1\langle 1 \rangle \rangle$ & 5 & 1 & 6 & \degsevenA & 1111\,0000\,1100\,0001\,0010\,0000\,0000\,0000\,0000 \\
$\langle J \sqcup 6 \rangle$ & 6 & 0 & 6 & \degsevenA & 1110\,0100\,1000\,0000\,0000\,0000\,0000\,0000\,0000 \\
$\langle J \sqcup 1\langle 4 \rangle \rangle$ & 1 & 4 & 5 & \degsevenA & 1100\,0000\,1010\,0000\,1000\,0000\,0000\,0000\,0000 \\
$\langle J \sqcup 1 \sqcup 1\langle 3 \rangle \rangle$ & 2 & 3 & 5 & \degsevenA & 1101\,0010\,1011\,0000\,0000\,0000\,0000\,0000\,0000 \\
$\langle J \sqcup 2 \sqcup 1\langle 2 \rangle \rangle$ & 3 & 2 & 5 & \degsevenA & 1100\,0110\,1000\,0000\,0000\,0000\,0000\,0000\,0000 \\
$\langle J \sqcup 3 \sqcup 1\langle 1 \rangle \rangle$ & 4 & 1 & 5 & \degsevenA & 1101\,0010\,1010\,0000\,0000\,0100\,0000\,0000\,0000 \\
$\langle J \sqcup 5 \rangle$ & 5 & 0 & 5 & \degsevenA & 1111\,0010\,1000\,0000\,0000\,0000\,0000\,0000\,0000 \\
$\langle J \sqcup 1\langle 3 \rangle \rangle$ & 1 & 3 & 4 & \degsevenA & 1101\,0010\,1011\,0000\,1000\,0000\,0000\,0000\,0000 \\
$\langle J \sqcup 1 \sqcup 1\langle 2 \rangle \rangle$ & 2 & 2 & 4 & \degsevenA & 1100\,0110\,1010\,0000\,0000\,0000\,0000\,0000\,0000 \\
$\langle J \sqcup 2 \sqcup 1\langle 1 \rangle \rangle$ & 3 & 1 & 4 & \degsevenA & 1101\,0010\,1011\,0000\,0000\,0100\,0000\,0000\,0000 \\
$\langle J \sqcup 4 \rangle$ & 4 & 0 & 4 & \degsevenA & 1111\,0000\,1000\,0000\,0000\,0000\,0000\,0000\,0000 \\
$\langle J \sqcup 1\langle 2 \rangle \rangle$ & 1 & 2 & 3 & \degsevenA & 1100\,0110\,1010\,0000\,1000\,0000\,0000\,0000\,0000 \\
$\langle J \sqcup 1\langle 1\langle 1 \rangle \rangle \rangle$ & 2 & 1 & 3 & \degsevenD & 1100\,0000\,1000\,0000\,0000\,0000\,0000\,0000\,0000 \\
$\langle J \sqcup 1 \sqcup 1\langle 1 \rangle \rangle$ & 2 & 1 & 3 & \degsevenA & 1101\,0010\,1011\,0000\,1000\,0100\,0000\,0000\,0000 \\
$\langle J \sqcup 3 \rangle$ & 3 & 0 & 3 & \degsevenA & 1111\,0000\,1000\,0000\,1000\,0000\,0000\,0000\,0000 \\
$\langle J \sqcup 1\langle 1 \rangle \rangle$ & 1 & 1 & 2 & \degsevenA & 1101\,0010\,1011\,0000\,1000\,0110\,0000\,0000\,0000 \\
$\langle J \sqcup 2 \rangle$ & 2 & 0 & 2 & \degsevenA & 1111\,0000\,1000\,0000\,1010\,0000\,0000\,0000\,0000 \\
$\langle J \sqcup 1 \rangle$ & 1 & 0 & 1 & \degsevenA & 1100\,0110\,1010\,0001\,0000\,0100\,0000\,0000\,0000 \\
$\langle J \rangle$ & 0 & 0 & 0 & \degsevenA & 1111\,0000\,1000\,0000\,1010\,0110\,0000\,0000\,0000 \\
\end{longtable}
\normalsize

%% file: tikz/gudkov_53_3_density.tex
\begin{tikzpicture}
    \begin{axis}[
        xlabel={Number of components},
        ylabel={Frequency},
        title={\degsixA (mean: 4.41)},
        xtick={0,1,2,3,4,5,6,7,8,9,10,11},
        ymajorgrids=true,
        grid style={dashed,gray!30},
        width=12cm,
        height=8cm,
        ymin=0,
        ymax=0.3214599609375,
        yticklabel style={/pgf/number format/.cd, fixed, fixed zerofill, precision=2},
    ]
    \addplot[const plot, draw=black, thick, no markers] coordinates {
        (-0.5,0.000000)
        (0.5,0.011963)
        (1.5,0.071777)
        (2.5,0.185425)
        (3.5,0.270020)
        (4.5,0.244278)
        (5.5,0.143341)
        (6.5,0.055771)
        (7.5,0.014587)
        (8.5,0.002548)
        (9.5,0.000275)
        (10.5,0.000015)
        (11.5,0)
    };
    \node[font=\fontsize{5}{6}\selectfont] at (axis cs:0,0.000000) [above] {0.00};
    \node[font=\fontsize{5}{6}\selectfont] at (axis cs:1,0.011963) [above] {1.20};
    \node[font=\fontsize{5}{6}\selectfont] at (axis cs:2,0.071777) [above] {7.18};
    \node[font=\fontsize{5}{6}\selectfont] at (axis cs:3,0.185425) [above] {18.54};
    \node[font=\fontsize{5}{6}\selectfont] at (axis cs:4,0.270020) [above] {27.00};
    \node[font=\fontsize{5}{6}\selectfont] at (axis cs:5,0.244278) [above] {24.43};
    \node[font=\fontsize{5}{6}\selectfont] at (axis cs:6,0.143341) [above] {14.33};
    \node[font=\fontsize{5}{6}\selectfont] at (axis cs:7,0.055771) [above] {5.58};
    \node[font=\fontsize{5}{6}\selectfont] at (axis cs:8,0.014587) [above] {1.46};
    \node[font=\fontsize{5}{6}\selectfont] at (axis cs:9,0.002548) [above] {0.25};
    \node[font=\fontsize{5}{6}\selectfont] at (axis cs:10,0.000275) [above] {0.03};
    \node[font=\fontsize{5}{6}\selectfont] at (axis cs:11,0.000015) [above] {0.00};
    \end{axis}
    \end{tikzpicture}

%% file: tikz/symmetric-14822_density.tex
\begin{tikzpicture}
    \begin{axis}[
        xlabel={Number of components},
        ylabel={Frequency},
        title={\degsixB (mean: 3.70)},
        xtick={0,1,2,3,4,5,6,7,8,9,10,11},
        ymajorgrids=true,
        grid style={dashed,gray!30},
        width=12cm,
        height=8cm,
        ymin=0,
        ymax=0.3214599609375,
        yticklabel style={/pgf/number format/.cd, fixed, fixed zerofill, precision=2},
    ]
    \addplot[const plot, draw=black, thick, no markers] coordinates {
        (-0.5,0.000000)
        (0.5,0.042480)
        (1.5,0.161133)
        (2.5,0.266724)
        (3.5,0.258423)
        (4.5,0.165344)
        (5.5,0.074341)
        (6.5,0.024323)
        (7.5,0.005951)
        (8.5,0.001114)
        (9.5,0.000153)
        (10.5,0.000015)
        (11.5,0)
    };
    \node[font=\fontsize{5}{6}\selectfont] at (axis cs:0,0.000000) [above] {0.00};
    \node[font=\fontsize{5}{6}\selectfont] at (axis cs:1,0.042480) [above] {4.25};
    \node[font=\fontsize{5}{6}\selectfont] at (axis cs:2,0.161133) [above] {16.11};
    \node[font=\fontsize{5}{6}\selectfont] at (axis cs:3,0.266724) [above] {26.67};
    \node[font=\fontsize{5}{6}\selectfont] at (axis cs:4,0.258423) [above] {25.84};
    \node[font=\fontsize{5}{6}\selectfont] at (axis cs:5,0.165344) [above] {16.53};
    \node[font=\fontsize{5}{6}\selectfont] at (axis cs:6,0.074341) [above] {7.43};
    \node[font=\fontsize{5}{6}\selectfont] at (axis cs:7,0.024323) [above] {2.43};
    \node[font=\fontsize{5}{6}\selectfont] at (axis cs:8,0.005951) [above] {0.60};
    \node[font=\fontsize{5}{6}\selectfont] at (axis cs:9,0.001114) [above] {0.11};
    \node[font=\fontsize{5}{6}\selectfont] at (axis cs:10,0.000153) [above] {0.02};
    \node[font=\fontsize{5}{6}\selectfont] at (axis cs:11,0.000015) [above] {0.00};
    \end{axis}
    \end{tikzpicture}

%% file: tables/degree6_cover.tex
\small
\begin{longtable}{lrrrrr}
\caption{Minimal support for degree 6 curves} \\
\toprule
Real scheme & \multicolumn{1}{c}{$p$} & \multicolumn{1}{c}{$n$} & \multicolumn{1}{c}{$p+n$} & \multicolumn{1}{c}{\degsixA} & \multicolumn{1}{c}{\degsixB} \\
\midrule
\endfirsthead
\toprule
Real scheme & \multicolumn{1}{c}{$p$} & \multicolumn{1}{c}{$n$} & \multicolumn{1}{c}{$p+n$} & \multicolumn{1}{c}{\degsixA} & \multicolumn{1}{c}{\degsixB} \\
\midrule
\endhead
\midrule
\endfoot
\bottomrule
\endlastfoot
$\langle 1 \sqcup 1\langle 9 \rangle \rangle$ & 2 & 9 & 11 & 128 &  \\ %
$\langle 5 \sqcup 1\langle 5 \rangle \rangle$ & 6 & 5 & 11 & 256 &  \\ %
$\langle 9 \sqcup 1\langle 1 \rangle \rangle$ & 10 & 1 & 11 & 128 & 512 \\ %
$\langle 1\langle 9 \rangle \rangle$ & 1 & 9 & 10 & 256 &  \\ %
$\langle 1 \sqcup 1\langle 8 \rangle \rangle$ & 2 & 8 & 10 & 2,048 &  \\ %
$\langle 4 \sqcup 1\langle 5 \rangle \rangle$ & 5 & 5 & 10 & 2,560 &  \\ %
$\langle 5 \sqcup 1\langle 4 \rangle \rangle$ & 6 & 4 & 10 & 2,048 &  \\ %
$\langle 8 \sqcup 1\langle 1 \rangle \rangle$ & 9 & 1 & 10 & 2,304 & 4,608 \\ %
$\langle 10 \rangle$ & 10 & 0 & 10 &  & 512 \\ %
$\langle 1\langle 8 \rangle \rangle$ & 1 & 8 & 9 & 4,224 & 512 \\ %
$\langle 1 \sqcup 1\langle 7 \rangle \rangle$ & 2 & 7 & 9 & 15,360 &  \\ %
$\langle 2 \sqcup 1\langle 6 \rangle \rangle$ & 3 & 6 & 9 & 3,584 &  \\ %
$\langle 3 \sqcup 1\langle 5 \rangle \rangle$ & 4 & 5 & 9 & 11,264 &  \\ %
$\langle 4 \sqcup 1\langle 4 \rangle \rangle$ & 5 & 4 & 9 & 20,736 & 1,024 \\ %
$\langle 5 \sqcup 1\langle 3 \rangle \rangle$ & 6 & 3 & 9 & 7,168 & 2,048 \\ %
$\langle 6 \sqcup 1\langle 2 \rangle \rangle$ & 7 & 2 & 9 & 3,584 & 2,048 \\ %
$\langle 7 \sqcup 1\langle 1 \rangle \rangle$ & 8 & 1 & 9 & 19,456 & 24,576 \\ %
$\langle 9 \rangle$ & 9 & 0 & 9 & 128 & 7,168 \\ %
$\langle 1\langle 7 \rangle \rangle$ & 1 & 7 & 8 & 32,768 & 4,608 \\ %
$\langle 1 \sqcup 1\langle 6 \rangle \rangle$ & 2 & 6 & 8 & 78,848 &  \\ %
$\langle 2 \sqcup 1\langle 5 \rangle \rangle$ & 3 & 5 & 8 & 57,344 &  \\ %
$\langle 3 \sqcup 1\langle 4 \rangle \rangle$ & 4 & 4 & 8 & 92,160 & 7,168 \\ %
$\langle 4 \sqcup 1\langle 3 \rangle \rangle$ & 5 & 3 & 8 & 73,728 & 17,408 \\ %
$\langle 5 \sqcup 1\langle 2 \rangle \rangle$ & 6 & 2 & 8 & 50,176 & 20,480 \\ %
$\langle 6 \sqcup 1\langle 1 \rangle \rangle$ & 7 & 1 & 8 & 102,400 & 100,352 \\ %
$\langle 8 \rangle$ & 8 & 0 & 8 & 2,048 & 49,664 \\ %
$\langle 1\langle 6 \rangle \rangle$ & 1 & 6 & 7 & 158,720 & 22,528 \\ %
$\langle 1 \sqcup 1\langle 5 \rangle \rangle$ & 2 & 5 & 7 & 318,976 & 10,240 \\ %
$\langle 2 \sqcup 1\langle 4 \rangle \rangle$ & 3 & 4 & 7 & 336,896 & 30,720 \\ %
$\langle 3 \sqcup 1\langle 3 \rangle \rangle$ & 4 & 3 & 7 & 356,864 & 87,040 \\ %
$\langle 4 \sqcup 1\langle 2 \rangle \rangle$ & 5 & 2 & 7 & 308,224 & 116,736 \\ %
$\langle 5 \sqcup 1\langle 1 \rangle \rangle$ & 6 & 1 & 7 & 376,320 & 329,728 \\ %
$\langle 7 \rangle$ & 7 & 0 & 7 & 15,360 & 219,136 \\ %
$\langle 1\langle 5 \rangle \rangle$ & 1 & 5 & 6 & 530,432 & 86,016 \\ %
$\langle 1 \sqcup 1\langle 4 \rangle \rangle$ & 2 & 4 & 6 & 989,184 & 116,736 \\ %
$\langle 2 \sqcup 1\langle 3 \rangle \rangle$ & 3 & 3 & 6 & 1,111,040 & 270,336 \\ %
$\langle 3 \sqcup 1\langle 2 \rangle \rangle$ & 4 & 2 & 6 & 1,089,536 & 448,512 \\ %
$\langle 4 \sqcup 1\langle 1 \rangle \rangle$ & 5 & 1 & 6 & 1,017,856 & 888,832 \\ %
$\langle 6 \rangle$ & 6 & 0 & 6 & 71,680 & 684,032 \\ %
$\langle 1\langle 4 \rangle \rangle$ & 1 & 4 & 5 & 1,261,568 & 288,768 \\ %
$\langle 1 \sqcup 1\langle 3 \rangle \rangle$ & 2 & 3 & 5 & 2,222,080 & 608,256 \\ %
$\langle 2 \sqcup 1\langle 2 \rangle \rangle$ & 3 & 2 & 5 & 2,433,536 & 1,122,304 \\ %
$\langle 3 \sqcup 1\langle 1 \rangle \rangle$ & 4 & 1 & 5 & 2,050,048 & 1,935,360 \\ %
$\langle 5 \rangle$ & 5 & 0 & 5 & 229,376 & 1,593,344 \\ %
$\langle 1\langle 3 \rangle \rangle$ & 1 & 3 & 4 & 2,121,728 & 802,816 \\ %
$\langle 1 \sqcup 1\langle 2 \rangle \rangle$ & 2 & 2 & 4 & 3,411,968 & 1,865,728 \\ %
$\langle 2 \sqcup 1\langle 1 \rangle \rangle$ & 3 & 1 & 4 & 3,010,560 & 3,196,928 \\ %
$\langle 4 \rangle$ & 4 & 0 & 4 & 516,096 & 2,805,760 \\ %
$\langle 1\langle 2 \rangle \rangle$ & 1 & 2 & 3 & 2,408,448 & 1,646,592 \\ %
$\langle 1 \sqcup 1\langle 1 \rangle \rangle$ & 2 & 1 & 3 & 3,010,560 & 3,665,920 \\ %
$\langle 3 \rangle$ & 3 & 0 & 3 & 802,816 & 3,629,056 \\ %
$\langle 1\langle 1\langle 1 \rangle \rangle \rangle$ & 2 & 1 & 3 &  & 8,192 \\ %
$\langle 1\langle 1 \rangle \rangle$ & 1 & 1 & 2 & 1,605,632 & 2,236,416 \\ %
$\langle 2 \rangle$ & 2 & 0 & 2 & 802,816 & 3,170,304 \\ %
$\langle 1 \rangle$ & 1 & 0 & 1 & 401,408 & 1,425,408 \\ %
\midrule
\text{total} & & & & 33,554,432 & 33,554,432 \\ %
\end{longtable}
\normalsize

%% file: tikz/Jv1v1_13__centrifugal_density.tex
\begin{tikzpicture}
    \begin{axis}[
        xlabel={Number of components},
        ylabel={Frequency},
        title={\degsevenA (mean: 8.50)},
        xtick={0,1,2,3,4,5,6,7,8,9,10,11,12,13,14,15,16},
        ymajorgrids=true,
        grid style={dashed,gray!30},
        width=12cm,
        height=8cm,
        ymin=0,
        ymax=0.26646813154220583,
        yticklabel style={/pgf/number format/.cd, fixed, fixed zerofill, precision=2},
    ]
    \addplot[const plot, draw=black, thick, no markers] coordinates {
        (-0.5,0.000000)
        (0.5,0.000122)
        (1.5,0.001343)
        (2.5,0.006958)
        (3.5,0.022827)
        (4.5,0.053833)
        (5.5,0.098022)
        (6.5,0.143677)
        (7.5,0.173218)
        (8.5,0.173218)
        (9.5,0.143677)
        (10.5,0.098022)
        (11.5,0.053833)
        (12.5,0.022827)
        (13.5,0.006958)
        (14.5,0.001343)
        (15.5,0.000122)
        (16.5,0)
    };
    \node[font=\fontsize{5}{6}\selectfont] at (axis cs:0,0.000000) [above] {0.00};
    \node[font=\fontsize{5}{6}\selectfont] at (axis cs:1,0.000122) [above] {0.01};
    \node[font=\fontsize{5}{6}\selectfont] at (axis cs:2,0.001343) [above] {0.13};
    \node[font=\fontsize{5}{6}\selectfont] at (axis cs:3,0.006958) [above] {0.70};
    \node[font=\fontsize{5}{6}\selectfont] at (axis cs:4,0.022827) [above] {2.28};
    \node[font=\fontsize{5}{6}\selectfont] at (axis cs:5,0.053833) [above] {5.38};
    \node[font=\fontsize{5}{6}\selectfont] at (axis cs:6,0.098022) [above] {9.80};
    \node[font=\fontsize{5}{6}\selectfont] at (axis cs:7,0.143677) [above] {14.37};
    \node[font=\fontsize{5}{6}\selectfont] at (axis cs:8,0.173218) [above] {17.32};
    \node[font=\fontsize{5}{6}\selectfont] at (axis cs:9,0.173218) [above] {17.32};
    \node[font=\fontsize{5}{6}\selectfont] at (axis cs:10,0.143677) [above] {14.37};
    \node[font=\fontsize{5}{6}\selectfont] at (axis cs:11,0.098022) [above] {9.80};
    \node[font=\fontsize{5}{6}\selectfont] at (axis cs:12,0.053833) [above] {5.38};
    \node[font=\fontsize{5}{6}\selectfont] at (axis cs:13,0.022827) [above] {2.28};
    \node[font=\fontsize{5}{6}\selectfont] at (axis cs:14,0.006958) [above] {0.70};
    \node[font=\fontsize{5}{6}\selectfont] at (axis cs:15,0.001343) [above] {0.13};
    \node[font=\fontsize{5}{6}\selectfont] at (axis cs:16,0.000122) [above] {0.01};
    \end{axis}
    \end{tikzpicture}

%% file: tikz/Jv2v1_12__split_centrifugal_density.tex
\begin{tikzpicture}
    \begin{axis}[
        xlabel={Number of components},
        ylabel={Frequency},
        title={\degsevenB (mean: 8.50)},
        xtick={0,1,2,3,4,5,6,7,8,9,10,11,12,13,14,15,16},
        ymajorgrids=true,
        grid style={dashed,gray!30},
        width=12cm,
        height=8cm,
        ymin=0,
        ymax=0.26646813154220583,
        yticklabel style={/pgf/number format/.cd, fixed, fixed zerofill, precision=2},
    ]
    \addplot[const plot, draw=black, thick, no markers] coordinates {
        (-0.5,0.000000)
        (0.5,0.000122)
        (1.5,0.001343)
        (2.5,0.006958)
        (3.5,0.022827)
        (4.5,0.053833)
        (5.5,0.098022)
        (6.5,0.143677)
        (7.5,0.173218)
        (8.5,0.173218)
        (9.5,0.143677)
        (10.5,0.098022)
        (11.5,0.053833)
        (12.5,0.022827)
        (13.5,0.006958)
        (14.5,0.001343)
        (15.5,0.000122)
        (16.5,0)
    };
    \node[font=\fontsize{5}{6}\selectfont] at (axis cs:0,0.000000) [above] {0.00};
    \node[font=\fontsize{5}{6}\selectfont] at (axis cs:1,0.000122) [above] {0.01};
    \node[font=\fontsize{5}{6}\selectfont] at (axis cs:2,0.001343) [above] {0.13};
    \node[font=\fontsize{5}{6}\selectfont] at (axis cs:3,0.006958) [above] {0.70};
    \node[font=\fontsize{5}{6}\selectfont] at (axis cs:4,0.022827) [above] {2.28};
    \node[font=\fontsize{5}{6}\selectfont] at (axis cs:5,0.053833) [above] {5.38};
    \node[font=\fontsize{5}{6}\selectfont] at (axis cs:6,0.098022) [above] {9.80};
    \node[font=\fontsize{5}{6}\selectfont] at (axis cs:7,0.143677) [above] {14.37};
    \node[font=\fontsize{5}{6}\selectfont] at (axis cs:8,0.173218) [above] {17.32};
    \node[font=\fontsize{5}{6}\selectfont] at (axis cs:9,0.173218) [above] {17.32};
    \node[font=\fontsize{5}{6}\selectfont] at (axis cs:10,0.143677) [above] {14.37};
    \node[font=\fontsize{5}{6}\selectfont] at (axis cs:11,0.098022) [above] {9.80};
    \node[font=\fontsize{5}{6}\selectfont] at (axis cs:12,0.053833) [above] {5.38};
    \node[font=\fontsize{5}{6}\selectfont] at (axis cs:13,0.022827) [above] {2.28};
    \node[font=\fontsize{5}{6}\selectfont] at (axis cs:14,0.006958) [above] {0.70};
    \node[font=\fontsize{5}{6}\selectfont] at (axis cs:15,0.001343) [above] {0.13};
    \node[font=\fontsize{5}{6}\selectfont] at (axis cs:16,0.000122) [above] {0.01};
    \end{axis}
    \end{tikzpicture}

%% file: tikz/Jv3v1_11__frayed_centrifugal_density.tex
\begin{tikzpicture}
    \begin{axis}[
        xlabel={Number of components},
        ylabel={Frequency},
        title={\degsevenC (mean: 7.78)},
        xtick={0,1,2,3,4,5,6,7,8,9,10,11,12,13,14,15,16},
        ymajorgrids=true,
        grid style={dashed,gray!30},
        width=12cm,
        height=8cm,
        ymin=0,
        ymax=0.26646813154220583,
        yticklabel style={/pgf/number format/.cd, fixed, fixed zerofill, precision=2},
    ]
    \addplot[const plot, draw=black, thick, no markers] coordinates {
        (-0.5,0.000000)
        (0.5,0.000000)
        (1.5,0.002441)
        (2.5,0.014648)
        (3.5,0.042908)
        (4.5,0.086548)
        (5.5,0.136230)
        (6.5,0.172241)
        (7.5,0.177917)
        (8.5,0.153076)
        (9.5,0.108887)
        (10.5,0.062744)
        (11.5,0.029114)
        (12.5,0.010376)
        (13.5,0.002441)
        (14.5,0.000366)
        (15.5,0.000061)
        (16.5,0)
    };
    \node[font=\fontsize{5}{6}\selectfont] at (axis cs:0,0.000000) [above] {0.00};
    \node[font=\fontsize{5}{6}\selectfont] at (axis cs:1,0.000000) [above] {0.00};
    \node[font=\fontsize{5}{6}\selectfont] at (axis cs:2,0.002441) [above] {0.24};
    \node[font=\fontsize{5}{6}\selectfont] at (axis cs:3,0.014648) [above] {1.46};
    \node[font=\fontsize{5}{6}\selectfont] at (axis cs:4,0.042908) [above] {4.29};
    \node[font=\fontsize{5}{6}\selectfont] at (axis cs:5,0.086548) [above] {8.65};
    \node[font=\fontsize{5}{6}\selectfont] at (axis cs:6,0.136230) [above] {13.62};
    \node[font=\fontsize{5}{6}\selectfont] at (axis cs:7,0.172241) [above] {17.22};
    \node[font=\fontsize{5}{6}\selectfont] at (axis cs:8,0.177917) [above] {17.79};
    \node[font=\fontsize{5}{6}\selectfont] at (axis cs:9,0.153076) [above] {15.31};
    \node[font=\fontsize{5}{6}\selectfont] at (axis cs:10,0.108887) [above] {10.89};
    \node[font=\fontsize{5}{6}\selectfont] at (axis cs:11,0.062744) [above] {6.27};
    \node[font=\fontsize{5}{6}\selectfont] at (axis cs:12,0.029114) [above] {2.91};
    \node[font=\fontsize{5}{6}\selectfont] at (axis cs:13,0.010376) [above] {1.04};
    \node[font=\fontsize{5}{6}\selectfont] at (axis cs:14,0.002441) [above] {0.24};
    \node[font=\fontsize{5}{6}\selectfont] at (axis cs:15,0.000366) [above] {0.04};
    \node[font=\fontsize{5}{6}\selectfont] at (axis cs:16,0.000061) [above] {0.01};
    \end{axis}
    \end{tikzpicture}

%% file: tikz/honeycomb-7_density.tex
\begin{tikzpicture}
    \begin{axis}[
        xlabel={Number of components},
        ylabel={Frequency},
        title={\degsevenD (mean: 3.98)},
        xtick={0,1,2,3,4,5,6,7,8,9,10,11,12,13,14,15,16},
        ymajorgrids=true,
        grid style={dashed,gray!30},
        width=12cm,
        height=8cm,
        ymin=0,
        ymax=0.26646813154220583,
        yticklabel style={/pgf/number format/.cd, fixed, fixed zerofill, precision=2},
    ]
    \addplot[const plot, draw=black, thick, no markers] coordinates {
        (-0.5,0.000000)
        (0.5,0.033215)
        (1.5,0.133775)
        (2.5,0.238544)
        (3.5,0.253779)
        (4.5,0.183636)
        (5.5,0.097851)
        (6.5,0.040561)
        (7.5,0.013617)
        (8.5,0.003829)
        (9.5,0.000935)
        (10.5,0.000206)
        (11.5,0.000042)
        (12.5,0.000008)
        (13.5,0.000001)
        (14.5,0.000000)
        (15.5,0.000000)
        (16.5,0)
    };
    \node[font=\fontsize{5}{6}\selectfont] at (axis cs:0,0.000000) [above] {0.00};
    \node[font=\fontsize{5}{6}\selectfont] at (axis cs:1,0.033215) [above] {3.32};
    \node[font=\fontsize{5}{6}\selectfont] at (axis cs:2,0.133775) [above] {13.38};
    \node[font=\fontsize{5}{6}\selectfont] at (axis cs:3,0.238544) [above] {23.85};
    \node[font=\fontsize{5}{6}\selectfont] at (axis cs:4,0.253779) [above] {25.38};
    \node[font=\fontsize{5}{6}\selectfont] at (axis cs:5,0.183636) [above] {18.36};
    \node[font=\fontsize{5}{6}\selectfont] at (axis cs:6,0.097851) [above] {9.79};
    \node[font=\fontsize{5}{6}\selectfont] at (axis cs:7,0.040561) [above] {4.06};
    \node[font=\fontsize{5}{6}\selectfont] at (axis cs:8,0.013617) [above] {1.36};
    \node[font=\fontsize{5}{6}\selectfont] at (axis cs:9,0.003829) [above] {0.38};
    \node[font=\fontsize{5}{6}\selectfont] at (axis cs:10,0.000935) [above] {0.09};
    \node[font=\fontsize{5}{6}\selectfont] at (axis cs:11,0.000206) [above] {0.02};
    \node[font=\fontsize{5}{6}\selectfont] at (axis cs:12,0.000042) [above] {0.00};
    \node[font=\fontsize{5}{6}\selectfont] at (axis cs:13,0.000008) [above] {0.00};
    \node[font=\fontsize{5}{6}\selectfont] at (axis cs:14,0.000001) [above] {0.00};
    \node[font=\fontsize{5}{6}\selectfont] at (axis cs:15,0.000000) [above] {0.00};
    \node[font=\fontsize{5}{6}\selectfont] at (axis cs:16,0.000000) [above] {0.00};
    \end{axis}
    \end{tikzpicture}

%% file: tables/degree7_cover.tex
\small
\begin{longtable}{lrrrrrrr}
\caption{Minimal support for degree 7 curves} \\
\toprule
Real scheme & \multicolumn{1}{c}{$p$} & \multicolumn{1}{c}{$n$} & \multicolumn{1}{c}{$p+n$} & \multicolumn{1}{c}{\degsevenA} & \multicolumn{1}{c}{\degsevenB} & \multicolumn{1}{c}{\degsevenC} & \multicolumn{1}{c}{\degsevenD} \\
\midrule
\endfirsthead
\toprule
Real scheme & \multicolumn{1}{c}{$p$} & \multicolumn{1}{c}{$n$} & \multicolumn{1}{c}{$p+n$} & \multicolumn{1}{c}{\degsevenA} & \multicolumn{1}{c}{\degsevenB} & \multicolumn{1}{c}{\degsevenC} & \multicolumn{1}{c}{\degsevenD} \\
\midrule
\endhead
\midrule
\endfoot
\bottomrule
\endlastfoot
$\langle J \sqcup 1 \sqcup 1\langle 13 \rangle \rangle$ & 2 & 13 & 15 & 32,768 &  &  &  \\ %
$\langle J \sqcup 2 \sqcup 1\langle 12 \rangle \rangle$ & 3 & 12 & 15 &  & 65,536 &  &  \\ %
$\langle J \sqcup 3 \sqcup 1\langle 11 \rangle \rangle$ & 4 & 11 & 15 &  &  & 32,768 &  \\ %
$\langle J \sqcup 4 \sqcup 1\langle 10 \rangle \rangle$ & 5 & 10 & 15 & 32,768 &  &  &  \\ %
$\langle J \sqcup 5 \sqcup 1\langle 9 \rangle \rangle$ & 6 & 9 & 15 & 65,536 &  &  &  \\ %
$\langle J \sqcup 6 \sqcup 1\langle 8 \rangle \rangle$ & 7 & 8 & 15 &  & 131,072 &  &  \\ %
$\langle J \sqcup 7 \sqcup 1\langle 7 \rangle \rangle$ & 8 & 7 & 15 & 131,072 &  & 65,536 &  \\ %
$\langle J \sqcup 8 \sqcup 1\langle 6 \rangle \rangle$ & 9 & 6 & 15 & 65,536 &  &  &  \\ %
$\langle J \sqcup 9 \sqcup 1\langle 5 \rangle \rangle$ & 10 & 5 & 15 & 163,840 & 262,144 & 131,072 &  \\ %
$\langle J \sqcup 10 \sqcup 1\langle 4 \rangle \rangle$ & 11 & 4 & 15 &  & 65,536 &  &  \\ %
$\langle J \sqcup 11 \sqcup 1\langle 3 \rangle \rangle$ & 12 & 3 & 15 & 131,072 &  & 32,768 &  \\ %
$\langle J \sqcup 12 \sqcup 1\langle 2 \rangle \rangle$ & 13 & 2 & 15 & 32,768 &  &  &  \\ %
$\langle J \sqcup 13 \sqcup 1\langle 1 \rangle \rangle$ & 14 & 1 & 15 & 262,144 & 262,144 & 131,072 &  \\ %
$\langle J \sqcup 15 \rangle$ & 15 & 0 & 15 & 131,072 & 262,144 & 131,072 & 64 \\ %
$\langle J \sqcup 1\langle 13 \rangle \rangle$ & 1 & 13 & 14 & 65,536 &  &  &  \\ %
$\langle J \sqcup 1 \sqcup 1\langle 12 \rangle \rangle$ & 2 & 12 & 14 & 262,144 & 262,144 &  &  \\ %
$\langle J \sqcup 2 \sqcup 1\langle 11 \rangle \rangle$ & 3 & 11 & 14 &  & 393,216 & 32,768 &  \\ %
$\langle J \sqcup 3 \sqcup 1\langle 10 \rangle \rangle$ & 4 & 10 & 14 & 131,072 &  & 163,840 &  \\ %
$\langle J \sqcup 4 \sqcup 1\langle 9 \rangle \rangle$ & 5 & 9 & 14 & 458,752 &  &  &  \\ %
$\langle J \sqcup 5 \sqcup 1\langle 8 \rangle \rangle$ & 6 & 8 & 14 & 393,216 & 786,432 &  &  \\ %
$\langle J \sqcup 6 \sqcup 1\langle 7 \rangle \rangle$ & 7 & 7 & 14 & 786,432 & 524,288 & 196,608 &  \\ %
$\langle J \sqcup 7 \sqcup 1\langle 6 \rangle \rangle$ & 8 & 6 & 14 & 917,504 &  & 196,608 &  \\ %
$\langle J \sqcup 8 \sqcup 1\langle 5 \rangle \rangle$ & 9 & 5 & 14 & 1,376,256 & 2,097,152 & 524,288 &  \\ %
$\langle J \sqcup 9 \sqcup 1\langle 4 \rangle \rangle$ & 10 & 4 & 14 & 524,288 & 1,048,576 & 262,144 &  \\ %
$\langle J \sqcup 10 \sqcup 1\langle 3 \rangle \rangle$ & 11 & 3 & 14 & 1,048,576 & 131,072 & 163,840 &  \\ %
$\langle J \sqcup 11 \sqcup 1\langle 2 \rangle \rangle$ & 12 & 2 & 14 & 524,288 &  & 32,768 &  \\ %
$\langle J \sqcup 12 \sqcup 1\langle 1 \rangle \rangle$ & 13 & 1 & 14 & 2,555,904 & 2,621,440 & 786,432 &  \\ %
$\langle J \sqcup 14 \rangle$ & 14 & 0 & 14 & 2,490,368 & 3,670,016 & 786,432 & 1,152 \\ %
$\langle J \sqcup 1\langle 12 \rangle \rangle$ & 1 & 12 & 13 & 557,056 & 262,144 &  &  \\ %
$\langle J \sqcup 1 \sqcup 1\langle 11 \rangle \rangle$ & 2 & 11 & 13 & 917,504 & 1,703,936 &  &  \\ %
$\langle J \sqcup 2 \sqcup 1\langle 10 \rangle \rangle$ & 3 & 10 & 13 & 327,680 & 983,040 & 425,984 &  \\ %
$\langle J \sqcup 3 \sqcup 1\langle 9 \rangle \rangle$ & 4 & 9 & 13 & 1,114,112 & 262,144 & 327,680 &  \\ %
$\langle J \sqcup 4 \sqcup 1\langle 8 \rangle \rangle$ & 5 & 8 & 13 & 2,129,920 & 1,703,936 & 131,072 &  \\ %
$\langle J \sqcup 5 \sqcup 1\langle 7 \rangle \rangle$ & 6 & 7 & 13 & 3,080,192 & 3,407,872 & 196,608 &  \\ %
$\langle J \sqcup 6 \sqcup 1\langle 6 \rangle \rangle$ & 7 & 6 & 13 & 4,390,912 & 786,432 & 1,114,112 &  \\ %
$\langle J \sqcup 7 \sqcup 1\langle 5 \rangle \rangle$ & 8 & 5 & 13 & 5,701,632 & 7,602,176 & 1,638,400 &  \\ %
$\langle J \sqcup 8 \sqcup 1\langle 4 \rangle \rangle$ & 9 & 4 & 13 & 4,227,072 & 6,160,384 & 1,703,936 &  \\ %
$\langle J \sqcup 9 \sqcup 1\langle 3 \rangle \rangle$ & 10 & 3 & 13 & 4,390,912 & 1,441,792 & 458,752 &  \\ %
$\langle J \sqcup 10 \sqcup 1\langle 2 \rangle \rangle$ & 11 & 2 & 13 & 3,407,872 & 589,824 & 688,128 &  \\ %
$\langle J \sqcup 11 \sqcup 1\langle 1 \rangle \rangle$ & 12 & 1 & 13 & 11,534,336 & 11,796,480 & 2,621,440 &  \\ %
$\langle J \sqcup 13 \rangle$ & 13 & 0 & 13 & 17,989,632 & 23,068,672 & 11,665,408 & 10,560 \\ %
$\langle J \sqcup 1\langle 11 \rangle \rangle$ & 1 & 11 & 12 & 2,097,152 & 1,835,008 &  &  \\ %
$\langle J \sqcup 1 \sqcup 1\langle 10 \rangle \rangle$ & 2 & 10 & 12 & 2,359,296 & 4,718,592 & 262,144 &  \\ %
$\langle J \sqcup 2 \sqcup 1\langle 9 \rangle \rangle$ & 3 & 9 & 12 & 2,228,224 & 2,621,440 & 1,638,400 &  \\ %
$\langle J \sqcup 3 \sqcup 1\langle 8 \rangle \rangle$ & 4 & 8 & 12 & 4,063,232 & 2,883,584 & 589,824 &  \\ %
$\langle J \sqcup 4 \sqcup 1\langle 7 \rangle \rangle$ & 5 & 7 & 12 & 8,519,680 & 7,864,320 & 589,824 &  \\ %
$\langle J \sqcup 5 \sqcup 1\langle 6 \rangle \rangle$ & 6 & 6 & 12 & 12,320,768 & 5,767,168 & 2,162,688 &  \\ %
$\langle J \sqcup 6 \sqcup 1\langle 5 \rangle \rangle$ & 7 & 5 & 12 & 15,859,712 & 17,563,648 & 5,308,416 &  \\ %
$\langle J \sqcup 7 \sqcup 1\langle 4 \rangle \rangle$ & 8 & 4 & 12 & 15,335,424 & 19,922,944 & 5,570,560 &  \\ %
$\langle J \sqcup 8 \sqcup 1\langle 3 \rangle \rangle$ & 9 & 3 & 12 & 13,631,488 & 6,815,744 & 2,162,688 &  \\ %
$\langle J \sqcup 9 \sqcup 1\langle 2 \rangle \rangle$ & 10 & 2 & 12 & 12,713,984 & 5,505,024 & 2,949,120 &  \\ %
$\langle J \sqcup 10 \sqcup 1\langle 1 \rangle \rangle$ & 11 & 1 & 12 & 33,161,216 & 32,505,856 & 7,077,888 &  \\ %
$\langle J \sqcup 12 \rangle$ & 12 & 0 & 12 & 73,793,536 & 88,080,384 & 60,817,408 & 67,968 \\ %
$\langle J \sqcup 1\langle 10 \rangle \rangle$ & 1 & 10 & 11 & 4,882,432 & 5,570,560 &  &  \\ %
$\langle J \sqcup 1 \sqcup 1\langle 9 \rangle \rangle$ & 2 & 9 & 11 & 5,963,776 & 9,306,112 & 1,835,008 &  \\ %
$\langle J \sqcup 2 \sqcup 1\langle 8 \rangle \rangle$ & 3 & 8 & 11 & 6,553,600 & 8,585,216 & 3,080,192 &  \\ %
$\langle J \sqcup 3 \sqcup 1\langle 7 \rangle \rangle$ & 4 & 7 & 11 & 12,779,520 & 10,223,616 & 2,260,992 &  \\ %
$\langle J \sqcup 4 \sqcup 1\langle 6 \rangle \rangle$ & 5 & 6 & 11 & 24,608,768 & 14,811,136 & 2,555,904 &  \\ %
$\langle J \sqcup 5 \sqcup 1\langle 5 \rangle \rangle$ & 6 & 5 & 11 & 32,473,088 & 31,195,136 & 11,075,584 &  \\ %
$\langle J \sqcup 6 \sqcup 1\langle 4 \rangle \rangle$ & 7 & 4 & 11 & 35,586,048 & 43,450,368 & 14,614,528 &  \\ %
$\langle J \sqcup 7 \sqcup 1\langle 3 \rangle \rangle$ & 8 & 3 & 11 & 32,899,072 & 18,874,368 & 7,897,088 &  \\ %
$\langle J \sqcup 8 \sqcup 1\langle 2 \rangle \rangle$ & 9 & 2 & 11 & 32,702,464 & 22,478,848 & 7,536,640 &  \\ %
$\langle J \sqcup 9 \sqcup 1\langle 1 \rangle \rangle$ & 10 & 1 & 11 & 70,287,360 & 64,618,496 & 16,515,072 & 5,824 \\ %
$\langle J \sqcup 11 \rangle$ & 11 & 0 & 11 & 203,685,888 & 233,308,160 & 182,714,368 & 357,696 \\ %
$\langle J \sqcup 1\langle 9 \rangle \rangle$ & 1 & 9 & 10 & 8,585,216 & 10,616,832 & 524,288 &  \\ %
$\langle J \sqcup 1 \sqcup 1\langle 8 \rangle \rangle$ & 2 & 8 & 10 & 12,976,128 & 18,612,224 & 5,242,880 &  \\ %
$\langle J \sqcup 2 \sqcup 1\langle 7 \rangle \rangle$ & 3 & 7 & 10 & 15,073,280 & 18,743,296 & 5,406,720 &  \\ %
$\langle J \sqcup 3 \sqcup 1\langle 6 \rangle \rangle$ & 4 & 6 & 10 & 32,243,712 & 19,136,512 & 6,324,224 &  \\ %
$\langle J \sqcup 4 \sqcup 1\langle 5 \rangle \rangle$ & 5 & 5 & 10 & 50,593,792 & 43,778,048 & 13,697,024 &  \\ %
$\langle J \sqcup 5 \sqcup 1\langle 4 \rangle \rangle$ & 6 & 4 & 10 & 61,603,840 & 70,516,736 & 30,736,384 &  \\ %
$\langle J \sqcup 6 \sqcup 1\langle 3 \rangle \rangle$ & 7 & 3 & 10 & 61,210,624 & 38,141,952 & 21,004,288 &  \\ %
$\langle J \sqcup 7 \sqcup 1\langle 2 \rangle \rangle$ & 8 & 2 & 10 & 63,963,136 & 53,739,520 & 16,678,912 &  \\ %
$\langle J \sqcup 8 \sqcup 1\langle 1 \rangle \rangle$ & 9 & 1 & 10 & 118,816,768 & 104,202,240 & 31,457,280 & 99,072 \\ %
$\langle J \sqcup 10 \rangle$ & 10 & 0 & 10 & 416,940,032 & 464,519,168 & 407,896,064 & 1,673,472 \\ %
$\langle J \sqcup 1\langle 8 \rangle \rangle$ & 1 & 8 & 9 & 12,877,824 & 16,449,536 & 2,621,440 &  \\ %
$\langle J \sqcup 1 \sqcup 1\langle 7 \rangle \rangle$ & 2 & 7 & 9 & 22,544,384 & 32,636,928 & 8,912,896 &  \\ %
$\langle J \sqcup 2 \sqcup 1\langle 6 \rangle \rangle$ & 3 & 6 & 9 & 31,653,888 & 27,197,440 & 12,615,680 &  \\ %
$\langle J \sqcup 3 \sqcup 1\langle 5 \rangle \rangle$ & 4 & 5 & 9 & 57,737,216 & 45,350,912 & 16,252,928 &  \\ %
$\langle J \sqcup 4 \sqcup 1\langle 4 \rangle \rangle$ & 5 & 4 & 9 & 81,920,000 & 85,852,160 & 44,498,944 &  \\ %
$\langle J \sqcup 5 \sqcup 1\langle 3 \rangle \rangle$ & 6 & 3 & 9 & 90,832,896 & 62,259,200 & 42,205,184 & 960 \\ %
$\langle J \sqcup 6 \sqcup 1\langle 2 \rangle \rangle$ & 7 & 2 & 9 & 98,631,680 & 88,276,992 & 33,849,344 & 12,288 \\ %
$\langle J \sqcup 7 \sqcup 1\langle 1 \rangle \rangle$ & 8 & 1 & 9 & 166,002,688 & 143,392,768 & 51,642,368 & 786,624 \\ %
$\langle J \sqcup 9 \rangle$ & 9 & 0 & 9 & 671,973,376 & 732,758,016 & 722,731,008 & 7,232,896 \\ %
$\langle J \sqcup 1\langle 7 \rangle \rangle$ & 1 & 7 & 8 & 16,908,288 & 22,806,528 & 5,242,880 &  \\ %
$\langle J \sqcup 1 \sqcup 1\langle 6 \rangle \rangle$ & 2 & 6 & 8 & 33,292,288 & 42,205,184 & 14,942,208 &  \\ %
$\langle J \sqcup 2 \sqcup 1\langle 5 \rangle \rangle$ & 3 & 5 & 8 & 51,380,224 & 44,302,336 & 25,034,752 &  \\ %
$\langle J \sqcup 3 \sqcup 1\langle 4 \rangle \rangle$ & 4 & 4 & 8 & 81,526,784 & 78,905,344 & 46,268,416 & 768 \\ %
$\langle J \sqcup 4 \sqcup 1\langle 3 \rangle \rangle$ & 5 & 3 & 8 & 107,610,112 & 76,283,904 & 61,210,624 & 12,288 \\ %
$\langle J \sqcup 5 \sqcup 1\langle 2 \rangle \rangle$ & 6 & 2 & 8 & 122,945,536 & 112,459,776 & 54,132,736 & 175,104 \\ %
$\langle J \sqcup 6 \sqcup 1\langle 1 \rangle \rangle$ & 7 & 1 & 8 & 193,200,128 & 166,985,728 & 81,526,784 & 3,923,328 \\ %
$\langle J \sqcup 8 \rangle$ & 8 & 0 & 8 & 881,065,984 & 943,980,544 & 1,026,555,904 & 28,781,952 \\ %
$\langle J \sqcup 1\langle 6 \rangle \rangle$ & 1 & 6 & 7 & 19,234,816 & 26,148,864 & 7,340,032 & 2,496 \\ %
$\langle J \sqcup 1 \sqcup 1\langle 5 \rangle \rangle$ & 2 & 5 & 7 & 41,713,664 & 47,579,136 & 26,476,544 & 10,944 \\ %
$\langle J \sqcup 2 \sqcup 1\langle 4 \rangle \rangle$ & 3 & 4 & 7 & 64,749,568 & 64,749,568 & 47,185,920 & 36,288 \\ %
$\langle J \sqcup 3 \sqcup 1\langle 3 \rangle \rangle$ & 4 & 3 & 7 & 97,320,960 & 67,633,152 & 66,322,432 & 124,672 \\ %
$\langle J \sqcup 4 \sqcup 1\langle 2 \rangle \rangle$ & 5 & 2 & 7 & 125,272,064 & 112,918,528 & 63,176,704 & 1,160,064 \\ %
$\langle J \sqcup 5 \sqcup 1\langle 1 \rangle \rangle$ & 6 & 1 & 7 & 187,334,656 & 163,708,928 & 114,950,144 & 14,107,968 \\ %
$\langle J \sqcup 7 \rangle$ & 7 & 0 & 7 & 952,303,616 & 1,005,191,168 & 1,202,847,744 & 101,523,584 \\ %
$\langle J \sqcup 1\langle 5 \rangle \rangle$ & 1 & 5 & 6 & 18,874,368 & 24,772,608 & 11,010,048 & 54,528 \\ %
$\langle J \sqcup 1 \sqcup 1\langle 4 \rangle \rangle$ & 2 & 4 & 6 & 43,515,904 & 50,855,936 & 40,370,176 & 254,208 \\ %
$\langle J \sqcup 2 \sqcup 1\langle 3 \rangle \rangle$ & 3 & 3 & 6 & 68,943,872 & 53,739,520 & 62,259,200 & 850,176 \\ %
$\langle J \sqcup 3 \sqcup 1\langle 2 \rangle \rangle$ & 4 & 2 & 6 & 101,187,584 & 83,886,080 & 60,424,192 & 4,919,808 \\ %
$\langle J \sqcup 4 \sqcup 1\langle 1 \rangle \rangle$ & 5 & 1 & 6 & 151,257,088 & 134,873,088 & 125,566,976 & 39,674,112 \\ %
$\langle J \sqcup 6 \rangle$ & 6 & 0 & 6 & 850,395,136 & 886,046,720 & 1,179,910,144 & 302,662,144 \\ %
$\langle J \sqcup 1\langle 4 \rangle \rangle$ & 1 & 4 & 5 & 16,056,320 & 21,037,056 & 15,204,352 & 625,664 \\ %
$\langle J \sqcup 1 \sqcup 1\langle 3 \rangle \rangle$ & 2 & 3 & 5 & 38,797,312 & 39,845,888 & 44,564,480 & 3,031,296 \\ %
$\langle J \sqcup 2 \sqcup 1\langle 2 \rangle \rangle$ & 3 & 2 & 5 & 63,307,776 & 49,414,144 & 54,132,736 & 14,579,712 \\ %
$\langle J \sqcup 3 \sqcup 1\langle 1 \rangle \rangle$ & 4 & 1 & 5 & 99,614,720 & 88,604,672 & 107,741,184 & 89,668,864 \\ %
$\langle J \sqcup 5 \rangle$ & 5 & 0 & 5 & 624,230,400 & 643,104,768 & 948,568,064 & 732,625,920 \\ %
$\langle J \sqcup 1\langle 3 \rangle \rangle$ & 1 & 3 & 4 & 11,927,552 & 14,942,208 & 14,680,064 & 4,609,024 \\ %
$\langle J \sqcup 1 \sqcup 1\langle 2 \rangle \rangle$ & 2 & 2 & 4 & 29,491,200 & 26,476,544 & 35,651,584 & 27,589,632 \\ %
$\langle J \sqcup 2 \sqcup 1\langle 1 \rangle \rangle$ & 3 & 1 & 4 & 50,987,008 & 43,778,048 & 78,643,200 & 155,397,120 \\ %
$\langle J \sqcup 4 \rangle$ & 4 & 0 & 4 & 370,016,256 & 377,225,216 & 614,465,536 & 1,389,829,120 \\ %
$\langle J \sqcup 1\langle 2 \rangle \rangle$ & 1 & 2 & 3 & 7,471,104 & 8,388,608 & 10,485,760 & 25,202,688 \\ %
$\langle J \sqcup 1\langle 1\langle 1 \rangle \rangle \rangle$ & 2 & 1 & 3 &  &  &  & 90,112 \\ %
$\langle J \sqcup 1 \sqcup 1\langle 1 \rangle \rangle$ & 2 & 1 & 3 & 18,612,224 & 16,252,928 & 41,943,040 & 180,443,136 \\ %
$\langle J \sqcup 3 \rangle$ & 3 & 0 & 3 & 170,000,384 & 171,442,176 & 316,145,664 & 1,974,210,560 \\ %
$\langle J \sqcup 1\langle 1 \rangle \rangle$ & 1 & 1 & 2 & 3,670,016 & 3,670,016 & 10,485,760 & 103,927,808 \\ %
$\langle J \sqcup 2 \rangle$ & 2 & 0 & 2 & 56,098,816 & 56,098,816 & 115,343,360 & 1,945,153,536 \\ %
$\langle J \sqcup 1 \rangle$ & 1 & 0 & 1 & 11,534,336 & 11,534,336 & 20,971,520 & 1,149,116,416 \\ %
$\langle J \rangle$ & 0 & 0 & 0 & 1,048,576 & 1,048,576 &  & 285,310,976 \\ %
\midrule
\text{total} & & & & 8,589,934,592 & 8,589,934,592 & 8,589,934,592 & 8,589,934,592 \\ %
\end{longtable}
\normalsize